\documentclass[a4,10pt]{article}
\usepackage{amsmath,amscd,amssymb}
\usepackage[margin=.8in]{geometry}

\newcommand{\nbiga}{\mathcal{A}}
\newcommand{\nbigb}{\mathcal{B}}
\newcommand{\nbigc}{\mathcal{C}}
\newcommand{\nbigd}{\mathcal{D}}
\newcommand{\nbige}{\mathcal{E}}
\newcommand{\nbigf}{\mathcal{F}}
\newcommand{\nbigg}{\mathcal{G}}
\newcommand{\nbigh}{\mathcal{H}}

\newcommand{\nbigj}{\mathcal{J}}
\newcommand{\nbigk}{\mathcal{K}}
\newcommand{\nbigl}{\mathcal{L}}
\newcommand{\nbigm}{\mathcal{M}}
\newcommand{\nbign}{\mathcal{N}}
\newcommand{\nbigo}{\mathcal{O}}
\newcommand{\nbigp}{\mathcal{P}}
\newcommand{\nbigq}{\mathcal{Q}}
\newcommand{\nbigr}{\mathcal{R}}
\newcommand{\nbigs}{\mathcal{S}}
\newcommand{\nbigt}{\mathcal{T}}
\newcommand{\nbigu}{\mathcal{U}}
\newcommand{\nbigv}{\mathcal{V}}
\newcommand{\nbigw}{\mathcal{W}}
\newcommand{\nbigx}{\mathcal{X}}
\newcommand{\nbigy}{\mathcal{Y}}
\newcommand{\nbigz}{\mathcal{Z}}

\newcommand{\proj}{\mathbb{P}}
\newcommand{\seisuu}{{\mathbb Z}}

\newcommand{\cnum}{{\mathbb C}}
\newcommand{\real}{{\mathbb R}}
\newcommand{\hyperh}{\mathbb{H}}

\newcommand{\Tate}{\mathbb{T}}
\newcommand{\newTate}{\pmb{T}}

\newcommand{\DD}{\mathbb{D}}
\newcommand{\EE}{\mathbb{E}}

\newcommand{\gbigc}{\mathfrak C}
\newcommand{\gbigd}{\mathfrak D}
\newcommand{\gbige}{\mathfrak E}
\newcommand{\gbigf}{\mathfrak F}

\newcommand{\gbigs}{\mathfrak S}
\newcommand{\gbigt}{\mathfrak T}
\newcommand{\gbigu}{\mathfrak U}
\newcommand{\gbigv}{\mathfrak V}

\newcommand{\gminib}{\mathfrak b}

\newcommand{\gminie}{\mathfrak e}

\newcommand{\gminik}{\mathfrak k}

\newcommand{\gminip}{\mathfrak p}

\newcommand{\gminis}{\mathfrak s}

\newcommand{\gminiv}{\mathfrak v}

\newcommand{\vecv}{{\boldsymbol v}}
\newcommand{\vecu}{{\boldsymbol u}}
\newcommand{\vecw}{{\boldsymbol w}}

\newcommand{\veczero}{{\boldsymbol 0}}

\newcommand{\vecalpha}{{\boldsymbol \alpha}}
\newcommand{\veca}{{\boldsymbol a}}
\newcommand{\vecb}{{\boldsymbol b}}
\newcommand{\vecbeta}{{\boldsymbol \beta}}
\newcommand{\vecdelta}{{\boldsymbol \delta}}

\newcommand{\vect}{{\boldsymbol t}}
\newcommand{\vecc}{{\boldsymbol c}}
\newcommand{\vecd}{{\boldsymbol d}}

\newcommand{\veck}{{\boldsymbol k}}
\newcommand{\vecm}{{\boldsymbol m}}

\newcommand{\vecN}{{\boldsymbol N}}

\newcommand{\vecx}{{\boldsymbol x}}
\newcommand{\vecf}{{\boldsymbol f}}
\newcommand{\vecepsilon}{{\boldsymbol \epsilon}}

\newcommand{\vecn}{{\boldsymbol n}}
\newcommand{\vecp}{{\boldsymbol p}}

\newcommand{\vecC}{{\boldsymbol C}}

\newcommand{\vecS}{{\boldsymbol S}}

\newcommand{\larr}{\leftarrow}
\newcommand{\llarr}{\longleftarrow}

\newcommand{\lrarr}{\longrightarrow}

\newcommand{\pf}{{\bf Proof}\hspace{.1in}}
\newcommand{\qed}{\mbox{\rule{1.2mm}{3mm}}}

\def\Hom{\mathop{\rm Hom}\nolimits}

\def\End{\mathop{\rm End}\nolimits}

\def\Cok{\mathop{\rm Cok}\nolimits}

\def\Image{\mathop{\rm Im}\nolimits}

\def\Re{\mathop{\rm Re}\nolimits}

\def\Gr{\mathop{\rm Gr}\nolimits}

\def\SL{\mathop{\rm SL}\nolimits}
\def\Tot{\mathop{\rm Tot}\nolimits}
\def\Cone{\mathop{\rm Cone}\nolimits}
\def\rank{\mathop{\rm rank}\nolimits}
\def\Spec{\mathop{\rm Spec}\nolimits}

\def\Ker{\mathop{\rm Ker}\nolimits}

\def\Gr{\mathop{\rm Gr}\nolimits}

\def\Res{\mathop{\rm Res}\nolimits}

\def\tr{\mathop{\rm tr}\nolimits}

\def\dvol{\mathop{\rm dvol}\nolimits}
\def\Diff{\mathop{\rm Diff}\nolimits}

\def\can{\mathop{\rm can}\nolimits}
\def\var{\mathop{\rm var}\nolimits}
\def\id{\mathop{\rm id}\nolimits}

\def\filt{\mathop{\rm fil}\nolimits}

\def\Supp{\mathop{\rm Supp}\nolimits}

\def\st{\mathop{\rm st}\nolimits}

\newcommand{\del}{\partial}
\newcommand{\delbar}{\overline{\del}}

\newcommand{\nbar}{\underline{n}}

\newcommand{\jbar}{\underline{j}}
\newcommand{\mbar}{\underline{m}}

\newcommand{\ibar}{\underline{i}}
\newcommand{\isitabar}{\underline{i}}

\newcommand{\sankaku}{\triangle}

\newcommand{\barz}{\overline{z}}
\newcommand{\zbar}{\barz}
\newcommand{\zetabar}{\overline{\zeta}}

\newcommand{\barlambda}{\overline{\lambda}}
\newcommand{\lambdabar}{\barlambda}
\newcommand{\fbar}{\overline{f}}

\newcommand{\abar}{\overline{a}}

\newcommand{\DDlambda}{\DD^{\lambda}}

\newcommand{\nbigelambda}{\nbige^{\lambda}}

\newcommand{\KMS}{{\mathcal{KMS}}}

\newcommand{\kmsmap}{\gminik}
\newcommand{\paramap}{\gminip}
\newcommand{\eigenmap}{\gminie}

\newcommand{\lefttop}[1]{{}^{#1}\!}

\def\reg{\mathop{\rm reg}\nolimits}

\def\Harm{\mathop{\rm Harm}\nolimits}

\newcommand{\Fzero}{F^{(\lambda_0)}}
\newcommand{\EEzero}{\EE^{(\lambda_0)}}

\newcommand{\Vzero}{V^{(\lambda_0)}}

\newcommand{\psizero}{\psi^{(\lambda_0)}}

\newcommand{\tildepsi}{\widetilde{\psi}}
\newcommand{\psitilde}{\tildepsi}
\newcommand{\nbigvzero}{\nbigv^{(\lambda_0)}}

\newcommand{\nbigqzero}{\nbigq^{(\lambda_0)}}

\newcommand{\deldel}{\eth}

\newcommand{\lamda}{\lambda}

\newcommand{\distribution}{\gbigd\gminib}

\newcommand{\supp}{\gminis}

\newcommand{\rtriplecat}{\nbigr\textrm{-Tri}}
\newcommand{\closedopen}[2]{[#1,#2[}
\newcommand{\openclosed}[2]{]#1,#2]}
\newcommand{\openopen}[2]{]#1,#2[}
\newcommand{\closedclosed}[2]{[#1,#2]}

\newcommand{\omegatilde}{\widetilde{\omega}}

\newcommand{\rhotilde}{\widetilde{\rho}}

\newcommand{\Etilde}{\widetilde{E}}

\newcommand{\thetatilde}{\widetilde{\theta}}

\newcommand{\Vtilde}{\widetilde{V}}

\newcommand{\htilde}{\widetilde{h}}

\newcommand{\ftilde}{\widetilde{f}}

\newcommand{\stilde}{\widetilde{s}}

\newcommand{\Psitilde}{\widetilde{\Psi}}

\newcommand{\nbigxhat}{\widehat{\nbigx}}
\newcommand{\nbigkhat}{\widehat{\nbigk}}

\newcommand{\Itilde}{\widetilde{I}}

\newcommand{\nbigvtilde}{\widetilde{\nbigv}}

\newcommand{\ellsitabar}{\underline{\ell}}
\newcommand{\Utilde}{\widetilde{U}}

\newcommand{\Xtilde}{\widetilde{X}}

\newcommand{\Ltilde}{\widetilde{L}}
\newcommand{\nbigmtilde}{\widetilde{\nbigm}}

\def\MT{\mathop{\rm MT}\nolimits}

\def\pol{\mathop{\rm pol}\nolimits}

\def\moderate{\mathop{\rm mod}\nolimits}

\def\exchange{\mathop{\rm exchange}\nolimits}

\def\TNIL{\mathop{\rm TNIL}\nolimits}

\def\DR{\mathop{\rm DR}\nolimits}

\def\Ch{\mathop{\rm Ch}\nolimits}

\def\hol{\mathop{\rm hol}\nolimits}

\def\MTM{\mathop{\rm MTM}\nolimits}

\def\sm{\mathop{\rm sm}\nolimits}

\def\pt{\mathop{\rm pt}\nolimits}
\def\IntC{\mathop{\rm IC}\nolimits}
\def\punc{\mathop{\rm punc}\nolimits}
\def\tw{\mathop{\rm tw}\nolimits}
\def\cp{\mathop{\rm cp}\nolimits}
\def\wc{\mathop{\rm wc}\nolimits}
\def\MTS{\mathop{\rm MTS}\nolimits}
\def\TS{\mathop{\rm TS}\nolimits}
\def\PTS{\mathop{\rm PTS}\nolimits}
\def\adm{\mathop{\rm adm}\nolimits}
\def\coh{\mathop{\rm coh}\nolimits}

\newcommand{\Mod}{{\tt Mod}}

\newcommand{\Ptilde}{\widetilde{P}}

\newcommand{\nbigstilde}{\widetilde{\nbigs}}

\newcommand{\nbigmhat}{\widehat{\nbigm}}

\newcommand{\Ztilde}{\widetilde{Z}}
\newcommand{\nbigttilde}{\widetilde{\nbigt}}
\newcommand{\gbigstilde}{\widetilde{\gbigs}}

\newcommand{\Ctilde}{\widetilde{C}}

\newcommand{\nbigpzero}{\nbigp^{(\lambda_0)}}
\newcommand{\nbiguzero}{\nbigu^{(\lambda_0)}}
\newcommand{\Mtilde}{\widetilde{M}}
\newcommand{\Ytilde}{\widetilde{Y}}
\newcommand{\Ybar}{\overline{Y}}

\newcommand{\nbigsbar}{\overline{\nbigs}}
\newcommand{\betabar}{\overline{\beta}}

\newcommand{\nbigktilde}{\widetilde{\nbigk}}

\newcommand{\vectheta}{{\boldsymbol \theta}}

\newcommand{\Lambdatilde}{\widetilde{\Lambda}}

\newcommand{\DDD}{\boldsymbol D}

\newcommand{\nbiggzero}{\nbigg^{(\lambda_0)}}

\newcommand{\Hhat}{\widehat{H}}
\newcommand{\Htilde}{\widetilde{H}}
\newcommand{\Omegatilde}{\widetilde{\Omega}}

\newcommand{\nbigxzero}{\nbigx^{(\lambda_0)}}
\newcommand{\nbigyzero}{\nbigy^{(\lambda_0)}}

\newcommand{\nbighzero}{\nbigh^{(\lambda_0)}}
\newcommand{\nbigzzero}{\nbigz^{(\lambda_0)}}

\newcommand{\nbigmzero}{\nbigm^{(\lambda_0)}}

\newcommand{\nbightilde}{\widetilde{\nbigh}}

\newcommand{\Xbar}{\overline{X}}

\newcommand{\nbigctilde}{\widetilde{\nbigc}}

\newcommand{\Zbar}{\overline{Z}}
\newcommand{\Fbar}{\overline{F}}

\newcommand{\Ccat}{{\sf C}}

\newcommand{\II}{\mathbb{I}}

\newcommand{\Ccheck}{\check{C}}

\newcommand{\ttS}{{\tt S}}
\newcommand{\ttG}{{\tt G}}

\newcommand{\ttF}{{\tt F}}

\newcommand{\ttY}{{\tt Y}}
\newcommand{\ttD}{{\tt D}}
\newcommand{\ttN}{{\tt N}}
\newcommand{\ttM}{{\tt M}}
\newcommand{\ttX}{{\tt X}}

\newcommand{\ttT}{{\tt T}}

\newcommand{\ttU}{{\tt U}}
\newcommand{\ttO}{{\tt O}}
\newcommand{\tte}{{\tt e}}
\newcommand{\ttv}{{\tt v}}
\newcommand{\ttd}{{\tt d}}
\newcommand{\ttC}{{\tt C}}

\newcommand{\vecnbign}{\boldsymbol\nbign}

\newcommand{\vecJ}{\boldsymbol{J}}

\newcommand{\Sigmabar}{\overline{\Sigma}}
\newcommand{\sigmabar}{\overline{\sigma}}

\newcommand{\ttNbar}{\overline{\ttN}}
\newcommand{\ttXbar}{\overline{\ttX}}
\newcommand{\ttYbar}{\overline{\ttY}}
\newcommand{\Fun}{{\tt Fun}}
\newcommand{\vecK}{\boldsymbol{K}}
\newcommand{\veciti}{\boldsymbol 1}
\newcommand{\Vtildezero}{\Vtilde^{(\lambda_0)}}
\newcommand{\Lbar}{\overline{L}}

\newcommand{\IItilde}{\widetilde{\II}}

\newcommand{\Xhat}{\widehat{X}}
\newcommand{\nbigshat}{\widehat{\nbigs}}
\newcommand{\nbigthat}{\widehat{\nbigt}}
\newcommand{\nbigchat}{\widehat{\nbigc}}
\newcommand{\upsilonhat}{\widehat{\upsilon}}
\newcommand{\nbighhat}{\widehat{\nbigh}}

\newtheorem{thm}{Theorem}[section]
\newtheorem{cor}[thm]{Corollary}

\newtheorem{rem}[thm]{Remark}
\newtheorem{lem}[thm]{Lemma}
\newtheorem{prop}[thm]{Proposition}
\newtheorem{df}[thm]{Definition}

\newtheorem{example}[thm]{Example}

\newtheorem{assumption}[thm]{Assumption}

\begin{document}

\title{$L^2$-complexes and twistor complexes of
tame harmonic bundles}

\author{Takuro Mochizuki}
\date{}
\maketitle

\begin{abstract}
Let $f:X\to Y$ be a morphism of complex manifolds.
Suppose that $X$ is a K\"ahler manifold.
Let $(\mathcal{T},\mathcal{S})$
be a regular polarized pure twistor $\mathcal{D}$-module
of weight $w$ on $X$ whose support is proper over $Y$. 
We prove the Hard Lefschetz Theorem for the push-forward of
 $(\mathcal{T},\mathcal{S})$ by $f$.
As one of the key steps,
we obtain the twistor version of
a theorem of Kashiwara and Kawai
about the Hodge structure on the intersection complex of
polarized variation of Hodge structure.

\vspace{.1in}
 \end{abstract}

\section{Introduction}

\subsection{$L^2$-cohomology and the intersection cohomology of
variation of Hodge structure}
\label{subsection;22.4.2.10}

Let us recall some classical results
on the $L^2$-cohomology of variation of Hodge structure.
Let $X$ be a compact connected K\"ahler manifold.
Let $H$ be a simple normal crossing hypersurface of $X$.
As in \cite{s1},
a polarized variation of complex Hodge structure of weight $w$
on $X\setminus H$
is formulated to be a graded vector bundle
$V=\bigoplus_{p+q=w} V^{p,q}$ on $X\setminus H$
equipped with a flat connection $\nabla$
and a flat $(-1)^w$-symmetric sesqui-linear pairing
$\langle\cdot,\cdot\rangle$
satisfying
(i) the Griffiths transversality condition,
i.e.,
\[
\nabla(V^{p,q})\subset
(V^{p,q}\oplus V^{p-1,q+1})
\otimes \Omega^{1,0}
\oplus
(V^{p,q}\oplus V^{p+1,q-1})
\otimes \Omega^{0,1},
\]
(ii) the decomposition
$V=\bigoplus V^{p,q}$ is orthogonal with respect to
$\langle\cdot,\cdot\rangle$,
(iii) $(\sqrt{-1})^{p-q}\langle\cdot,\cdot\rangle$
is positive definite on $V^{p,q}$.

There are two naturally induced complexes of sheaves on $X$.
One is the intersection complex of $(V,\nabla)$.
We obtain the $\nbigd_X$-module $V_{\min}$
as the minimal extension of $(V,\nabla)$.
The intersection complex of $(V,\nabla)$
is isomorphic to the de Rham complex
$\DR(V_{\min})=\Omega^{\bullet}_X\otimes V_{\min}$ on $X$
in the derived category of sheaves in a canonical way.
(Here, we do not make a shift of the degree
for the de Rham complex.)
Its cohomology group
$H^{\bullet}(X,\DR(V_{\min}))$
is isomorphic to the intersection cohomology of $(V,\nabla)$.
The other is the $L^2$-complex of $(V,\nabla)$
with respect to the Hodge metric
$h=\bigoplus (\sqrt{-1})^{p-q}\langle\cdot,\cdot\rangle_{|V^{p,q}}$
and a Poincar\'e like K\"ahler metric $g_{X\setminus H}$ of $X\setminus H$.
Recall that a K\"ahler metric $g_{X\setminus H}$ of $X\setminus H$
is called Poincar\'{e} like 
if for any $P\in H$
there exists a holomorphic coordinate neighbourhood
$(X_P,z_1,\ldots,z_n)$ such that
(i) $H\cap X_P=\bigcup_{i=1}^{\ell}\{z_i=0\}$,
(ii) $g_{X\setminus H|X_P\setminus H}$
is mutually bounded with 
$\sum_{i=1}^{\ell}|z_i|^{-2}(-\log |z_i|^2)^{-2}dz_i\,d\zbar_i
+\sum_{i=\ell+1}^ndz_i\,d\zbar_i$.
For any open subset $U$ of $X$,
let $\nbigc^k_{L^2}(V,\nabla,h)(U)$
denote the space of
sections $\tau$ of
$\Tot^k(V\otimes\Omega^{\bullet,\bullet})
=\bigoplus_{r+s=k}
(V\otimes\Omega^{r,s})$ on $U\setminus H$
such that $\tau$ and $\nabla\tau$
are $L^2$ with respect to $h$ and $g_{X\setminus H}$
around any point of $U$.
Thus, we obtain the complex
$\nbigc^{\bullet}_{L^2}(V,\nabla,h)$,
called the $L^2$-complex of $(V,\nabla,h)$.
By the construction,
there exists the natural inclusion,
which is a quasi-isomorphism:
\begin{equation}
\label{eq;22.4.3.1}
 \DR(V_{\min})_{|X\setminus H}
\lrarr
 \nbigc^{\bullet}_{L^2}(V,\nabla,h)_{|X\setminus H}.
\end{equation}
According to
Cattani-Kaplan-Schmid \cite{cks2}
and Kashiwara-Kawai \cite{k3},
the quasi-isomorphism (\ref{eq;22.4.3.1})
extends to an isomorphism
in the derived category
$\ttD(\cnum_X)$ of $\cnum_X$-modules.
This is due to Zucker \cite{z}
in the case $\dim X=1$.
\begin{thm}[Cattani-Kaplan-Schmid, Kashiwara-Kawai, Zucker]
\label{thm;22.4.3.2}
There exists an isomorphism
\begin{equation}
\DR(V_{\min})
\simeq
\nbigc^{\bullet}_{L^2}(V,\nabla,h)
\end{equation}
in $\ttD(\cnum_X)$
whose restriction to $X\setminus H$
is equal to {\rm(\ref{eq;22.4.3.1})}.
\hfill\qed
\end{thm}

Theorem \ref{thm;22.4.3.2} implies that the intersection cohomology
$H^{\bullet}(X,\DR(V_{\min}))$
is equipped with a Hodge structure.
Namely,
the space of harmonic $k$-forms
$\Harm^k(V,\nabla,h,g_{X\setminus H})$
of $(V,\nabla)$ with respect to $h$ and $g_{X\setminus H}$
is finite dimensional,
and it has the Hodge decomposition
\begin{equation}
\label{eq;22.4.12.2}
\Harm^k(V,\nabla,h,g_{X\setminus H})=
\bigoplus_{p+q=k+w} \Harm^{p,q}(V,\nabla,h,g_{X\setminus H}).
\end{equation}
The isomorphism
$H^k(X,\DR(V_{\min}))
\simeq
H^k(X,\nbigc^{\bullet}_{L^2}(V,\nabla,h))
\simeq
\Harm^k(V,\nabla,h,g_{X\setminus H})$
induces
the Hodge decomposition
\[
 H^k(X,\DR(V_{\min}))
 =\bigoplus_{p+q=k+w}
\Harm^{p,q}(V,\nabla,h,g_{X\setminus H}).
\]
In particular,
we obtain the Hodge filtration
$F^pH^k(X,\DR(V_{\min}))=
\bigoplus_{p'\geq p}\Harm^{p',q'}(V,\nabla,h,g_{X\setminus H})$.
Kashiwara and Kawai proved that
the filtration is independent of the choice of
a polarization $\langle\cdot,\cdot\rangle$
and $g_{X\setminus H}$.

Kashiwara and Kawai \cite{Kashiwara-Kawai-Hodge-holonomic}
studied how the Hodge filtration is described.
There exists the associated Hodge module
$(V_{\min},F)$ of weight $w+\dim X$
\cite{saito1,saito2}.
The Hodge filtration $F$ on $V_{\min}$
induces a filtration $F$ on $\DR(V_{\min})$.
Kashiwara and Kawai announced the theorem
\cite[Theorem 1]{Kashiwara-Kawai-Hodge-holonomic}
that
$H^{k}(X,F_p\DR(V_{\min}))
\lrarr
H^k(X,\DR(V_{\min}))$
is injective,
and the image is equal to the Hodge filtration
$F^{-p}H^k(X,\DR(V_{\min}))$
induced by the $L^2$-theory.
(See \S\ref{subsection;22.3.16.30}
for more details.)
In the one dimensional case,
it is again due to Zucker \cite{z}.

\subsection{Main result and some consequences}

Roughly speaking,
the main result of this paper
is a generalization 
of the theorem of Kashiwara and Kawai
to the context of tame harmonic bundles.

\subsubsection{Harmonic bundles}

Let $(E,\delbar_E)$ be a holomorphic vector bundle on
a complex manifold $Y$.
Let $\theta$ be a Higgs field of $(E,\delbar_E)$,
which means, $\theta$ is a holomorphic section of
$\End(E)\otimes\Omega^1$ such that $\theta\wedge\theta=0$.
Let $h$ be a Hermitian metric of $E$.
Then, we obtain the Chern connection
$\nabla_h=\del_{E,h}+\delbar_E$,
which is a unitary connection
whose $(0,1)$-part is equal to
the holomorphic structure $\delbar_E$.
We also obtain the adjoint $\theta^{\dagger}_h$
of $\theta$ with respect to $h$.
Then, $h$ is called a pluri-harmonic metric
of the Higgs bundle $(E,\delbar_E,\theta)$
if the connection $\DD^1=\nabla_h+\theta+\theta^{\dagger}$ is flat.
In that case, $(E,\delbar_E,\theta,h)$ is called a harmonic bundle.

As introduced in \cite{s3},
we obtain the holomorphic vector bundle
$\nbige$ on $\nbigy=\cnum_{\lambda}\times Y$.
Let $p_{\lambda}:\nbigy\to Y$ denote the projection,
and we set
$\nbige=\bigl(
p_{\lambda}^{-1}(E),
p_{\lambda}^{\ast}(\delbar_E)+\lambda\theta^{\dagger}_h
\bigr)$.
It is equipped with
the family of flat $\lambda$-connections
\[
\DD:\nbige\to\nbige\otimes
p_{\lambda}^{-1}\bigl(\Omega_Y^{1,0}
\oplus
\Omega^{0,1}_Y
\bigr)
\]
given by
$\DD=
 \lambda \del_{E,h}+\theta
+\delbar_E+\lambda\theta^{\dagger}$.
It is compatible with the holomorphic structure,
and satisfies
$\DD(fs)=(\lambda\del_Y+\delbar_Y)f\cdot s
+f\DD(s)$
and $\DD\circ\DD=0$.
We obtain $\DD:
\nbige\otimes p_{\lambda}^{-1}(\Omega^i_Y)
\to
\nbige\otimes p_{\lambda}^{-1}(\Omega^{i+1}_Y)$
induced by $\DD$ and $\lambda\del_Y$
for $\nbige$ and $p_{\lambda}^{-1}(\Omega^{\bullet})$,
respectively.
Thus, we obtain a complex of
$(\nbigo_{\cnum})_{\nbigy}$-modules
$\nbige\otimes p_{\lambda}^{-1}(\Omega^{\bullet}_Y)$.
Here, $(\nbigo_{\cnum})_{\nbigy}$ denotes the pull back of
$\nbigo_{\cnum}$ by the projection $\cnum\times Y\to \cnum$.
Together with $\delbar_Y$,
it extends to
$\nbige\otimes p_{\lambda}^{-1}(\Tot \Omega^{\bullet,\bullet}_Y)$,
and there exists the natural quasi-isomorphism
of $(\nbigo_{\cnum})_{\nbigy}$-modules
$\nbige\otimes p_{\lambda}^{-1}(\Omega^{\bullet}_Y)
\simeq
\nbige\otimes p_{\lambda}^{-1}(\Tot\Omega^{\bullet,\bullet}_Y)$.
As the restriction to $\{\lambda\}\times Y$,
we obtain
the holomorphic vector bundle
$\nbigelambda=(E,\delbar_E+\lambda\theta^{\dagger})$
and
the flat $\lambda$-connection
$\DDlambda=\lambda\del_{E,h}+\theta+\delbar_E+\lambda\theta^{\dagger}$.
We obtain the complexes
$\nbigelambda\otimes \Omega_Y^{\bullet}
\simeq
\nbigelambda\otimes \Tot\Omega_Y^{\bullet,\bullet}$.

\vspace{.1in}

For a polarized variation of complex Hodge structure
$(V=\bigoplus V^{p,q},\nabla,\langle\cdot,\cdot,\rangle)$ on $Y$,
we obtain the decomposition
$\nabla=\nabla'+\theta+\theta^{\dagger}$
by the Griffiths transversality condition,
where $\nabla'$ is a connection preserving the Hodge decomposition,
$\theta$ is a section of
$\bigoplus \Hom(V^{p,q},V^{p-1,q+1})\otimes\Omega^{1,0}$
and
$\theta^{\dagger}$ is a section of
$\bigoplus \Hom(V^{p,q},V^{p+1,q-1})\otimes\Omega^{0,1}$.
We have the decomposition
$\nabla'=\del_V+\delbar_V$
into the $(1,0)$-part and the $(0,1)$-part.
Then, $(V,\delbar_V,\theta)$ is a Higgs bundle.
Together with the Hodge metric $h$,
$(V,\delbar_V,\theta,h)$ is a harmonic bundle.
In this case,
the holomorphic vector bundle
$\nbigv=
(p_{\lambda}^{-1}V,p_{\lambda}^{\ast}(\delbar_V)+\lamda\theta^{\dagger})$
is isomorphic to the Rees module
associated with the Hodge filtration $(V,F)$,
where $F^p=\bigoplus_{p'\geq p} V^{p',w-p'}$.
The complex
$\nbigv\otimes
p_{\lambda}^{\ast}\Omega^{\bullet}_Y$
is the Rees module of the filtered de Rham complex.

\subsubsection{Tame harmonic bundles}

Let $X$ and $H$ be as in \S\ref{subsection;22.4.2.10}.
Let $(E,\delbar_E,\theta,h)$ be a harmonic bundle
on $X\setminus H$.
By the Higgs field,
we may regard $E$ as a module over
the symmetric product algebra of the tangent bundle
$\Theta_{X\setminus H}$.
It induces a coherent sheaf
on the cotangent bundle $\Omega^1_{X\setminus H}$.
The support $\Sigma_{(E,\theta)}$ is called
the spectral variety of the Higgs bundle.
The harmonic bundle is called tame on $(X,H)$
if the closure of the spectral variety
$\overline{\Sigma_{(E,\theta)}}$
in the logarithmic cotangent bundle
$\Omega^1_X(\log H)$ is proper over $X$.
Note that any polarized variation of Hodge structure on $X\setminus H$
is a tame harmonic bundle
because its spectral variety is the $0$-section.
Tame harmonic bundles were studied in {\rm\cite{Simpson-non-compact}}
in the one dimensional case,
and in {\rm\cite{mochi2}} in the higher dimensional case. 

\subsubsection{Main result}

For a tame harmonic bundle
$(E,\delbar_E,\theta,h)$ on $(X,H)$,
we obtain the $(\nbigo_{\cnum})_{\nbigx\setminus \nbigh}$-complexes
$\nbige\otimes p_{\lambda}^{\ast}\Omega_{X\setminus H}^{\bullet}$
on $\cnum\times(X\setminus H)$.
As in the Hodge case,
we consider two $(\nbigo_{\cnum})_{\nbigx}$-complexes on $\nbigx$
as extensions of 
the $(\nbigo_{\cnum})_{\nbigx\setminus\nbigh}$ -complexes
$\nbige\otimes p_{\lambda}^{\ast}\Omega_{X\setminus H}^{\bullet}$.

One is obtained by the theory of twistor $\nbigd$-modules.
Let $\nbigr_X$ denote the sheaf of subalgebras of
$\nbigd_{\nbigx}$ on $\cnum\times X$, generated by
$\lambda p_{\lambda}^{\ast}\Theta_X$ over $\nbigo_{\nbigx}$.
By the differential operator $\DD$,
we may naturally regard $\nbige$
as an $\nbigr_{X\setminus H}$-module.
In \cite{mochi2},
we obtained the $\nbigr_X$-module $\gbige$ on $\cnum\times X$
as the underlying pure twistor $\nbigd$-module associated with
the harmonic bundle $(E,\delbar_E,\theta,h)$,
which is an extension of $(\nbige,\DD)$
and a twistor analogue of the minimal extension of a flat bundle.
We obtain the complex of $(\nbigo_{\cnum_{\lambda}})_{\nbigx}$-modules
$\gbige\otimes p_{\lambda}^{\ast}\Omega_X^{\bullet}$
on $\cnum\times X$,
which is an extension of
$\nbige\otimes p_{\lambda}^{\ast}\Omega_{X\setminus H}^{\bullet}$.
In the Hodge case,
$\gbige$ is the Rees module of the filtered $\nbigd$-module
underlying the Hodge module,
and 
$\gbige\otimes p_{\lambda}^{\ast}\Omega_X^{\bullet}$
is the Rees module
of the associated filtered de Rham complex.

The other is obtained as the $L^2$-complex.
Let $U$ be any open subset of $\nbigx$.
Let $\nbigc^{k}_{L^2}(\nbige,\DD,h)(U)$
denote the space of measurable sections $\tau$ of
$\nbige\otimes p_{\lambda}^{-1}\Tot^k(\Omega^{\bullet,\bullet}_{X\setminus H})$
on $U\setminus(\cnum\times H)$,
such that
(i) $\del_{\lambdabar}\tau=0$,
(ii) $\tau$ and $\DD\tau$ are $L^2$
with respect to $h$
and $g_{X\setminus H}+d\lambda\,d\lambdabar$
locally around any point of $U$.
In this way,
we obtain the complex of $(\nbigo_{\cnum_{\lambda}})_{\nbigx}$-modules
$\nbigc^{\bullet}_{L^2}(\nbige,\DD,h)$
on $\nbigx$.

By the construction,
there exists the natural quasi-isomorphism
of $(\nbigo_{\cnum})_{\nbigx\setminus\nbigh}$-complexes
\begin{equation}
\label{eq;22.4.2.11}
 (\gbige\otimes p_{\lambda}^{\ast}\Omega^{\bullet})
 _{|\nbigx\setminus\nbigh}
\lrarr
\nbigc^{\bullet}_{L^2}(\nbige,\DD,h)_{|\nbigx\setminus\nbigh}.
\end{equation}

Let $\ttD\bigl((\nbigo_{\cnum_{\lambda}})_{\nbigx}\bigr)$
denote the derived category of
$(\nbigo_{\cnum_{\lambda}})_{\nbigx}$-complexes.

\begin{thm}[Lemma \ref{lem;22.3.17.11},
Corollary \ref{cor;22.2.26.1},
Lemma \ref{lem;22.2.18.20},
Theorem \ref{thm;22.2.17.111},
 Theorem \ref{thm;21.12.21.11}]
\label{thm;22.4.2.50}
There exists a natural isomorphism
\[
\DR(\gbige)\simeq
\nbigc^{\bullet}_{L^2}(\nbige,\DD,h)
\]
in $\ttD((\nbigo_{\cnum_{\lambda}})_{\nbigx})$,
whose restriction to
$\nbigx\setminus\nbigh$ is
the above quasi-isomorphism
{\rm(\ref{eq;22.4.2.11})}.
\end{thm}

We have a similar statement for the specialization
at $\lambda$.
We have the holomorphic vector bundle
$\nbigelambda$ with a flat $\lambda$-connection $\DDlambda$
on $X\setminus H$.
It induces a $\cnum_{X\setminus H}$-complex
$\nbigelambda\otimes\Omega^{\bullet}_{X\setminus H}$.
The differential is induced by $\DDlambda$ and $\lambda\del_X$.
It extends to a $\cnum_X$-complex
$\nbigc^{\bullet}_{\tw}(\nbigelambda,\DDlambda)$
which is related with the specialization of $\gbige$
to $\{\lambda\}\times X$.
(See \S\ref{subsection;22.4.2.12}.)
We also obtain the $L^2$-complex
$\nbigc^{\bullet}_{L^2}(\nbigelambda,\DDlambda,h)$.
(See \S\ref{subsection;22.4.2.20}.)

\begin{thm}[Corollary
\ref{cor;22.2.18.10},
Lemma \ref{lem;22.2.18.20},
Theorem \ref{thm;22.2.18.11}]
\label{thm;22.4.2.40}
There exists an isomorphism
\[
 \nbigc^{\bullet}_{\tw}(\gbige,\lambda)
\simeq
\nbigc^{\bullet}_{L^2}(\nbigelambda,\DDlambda,h)
\]
in $\ttD(\cnum_X)$
whose restriction to $X\setminus H$
is the natural quasi-isomorphism.
\end{thm}

\begin{rem}
If $\dim X=1$,
Theorems {\rm\ref{thm;22.4.2.50}} and
{\rm\ref{thm;22.4.2.40}}
were established in
{\rm\cite{sabbah2}} and {\rm\cite{mochi2}}
on the basis of the monumental study of Zucker {\rm\cite{z}}
for $L^2$-cohomology of variation of Hodge structure
on a punctured curve.
It was generalized to the wild case {\rm\cite{Mochizuki-wild}}
in the one dimensional case.
A related issue 
was studied in {\rm\cite{Jost-Yang-Zuo2}}
for a special class of tame harmonic bundles.
\hfill\qed
\end{rem}

\begin{rem}
Let $(E,\delbar_E,\theta,h)$ be a tame harmonic bundle
on $(X,H)$ such that the residue of the Higgs fields are nilpotent
and that the local monodromy of the flat connection
$\DD^1$ are quasi-unipotent.
In this case, it is easy to prove
Theorem {\rm\ref{thm;22.4.2.40}}
for $\lambda=1$
as a direct consequence of Theorem {\rm\ref{thm;22.4.3.2}}.
It is enough to check the claim locally around any point of
$P\in H$.
There exist a neighbourhood $X_P$ of $P$,
a polarized variation of Hodge structure
$(\nbigv_P,\nabla,F_P)$ on $X_P\setminus H$
and an isomorphism 
$(\nbige^1,\DD^1)_{|X_P\setminus H}
\simeq
(\nbigv_P,\nabla)_{|X_P\setminus H}$
under which $h_{|X_P\setminus H}$
and the Hodge metric of $(\nbigv_P,\nabla,F_P)$
are mutually bounded.
Hence,
we can deduce the claim for 
$(\nbige^1,\DD^1,h)$ from
Theorem {\rm\ref{thm;22.4.3.2}}.
\hfill\qed
\end{rem}
 
\subsubsection{A direct consequence}

We explain some consequences
in the theory of pure twistor $\nbigd$-modules
\cite{sabbah2, mochi2, Mochizuki-wild, Mochizuki-MTM}.
See \S\ref{section;22.4.2.1} for some explanations
of terminology.
(See \cite{sabbah2} and \cite{Mochizuki-MTM} for more details.
Note that
in \S\ref{section;22.4.2.1},
which is based on \cite{Mochizuki-MTM},
we adopt a different but equivalent convention for signatures
from that in \cite{sabbah2}.)

For a tame harmonic bundle $(E,\delbar_E,\theta,h)$ on $(X,H)$,
we have the pure twistor $\nbigd$-module
$\gbigt(E)=(\lambda^{\dim X}\gbige,\gbige,C_h)$
of weight $\dim X$ on $X$.
It is equipped with the polarization
$\nbigs:\gbigt(E)\to\gbigt(E)^{\ast}\otimes\newTate(-\dim X)$
given as the pair of the identity morphisms.
Let $a_X$ denote the canonical map of $X$
to the one point set $\pt$.
As the push-forward by $a_X$,
we obtain the graded $\nbigr$-triple
$\bigoplus a^j_{X\dagger}\gbigt(E)$.
Let $\omega$ be a K\"ahler form of $X$.
We obtain the morphisms
$L_{\omega}:a_{X\dagger}^j(\gbigt(E))
\lrarr
a_{X\dagger}^{j+2}(\gbigt(E))\otimes\newTate(1)$,
induced by the multiplication of 
$\sqrt{-1}\omega$.
We have the Hermitian sesqui-linear duality
$a_{\dagger}(\nbigs)$ of
$\bigoplus a^j_{X\dagger}\gbigt(E)$
of weight $\dim X$.

We obtain a Poincar\'{e} like K\"ahler metric $g_{X\setminus H}$
by modifying the K\"ahler metric associated with $\omega$
as in \cite{k3}.
(See \S\ref{subsection;22.4.2.21}).
Theorem \ref{thm;22.4.2.40} implies that
the space $\Harm^j(E)$
of harmonic forms of 
the tame harmonic bundle $(E,\delbar_E,\theta,h)$
with respect to $g_{X\setminus H}$.
We obtain the graded $\nbigr$-triple
$\bigoplus_k\nbightilde^k$,
associated with $\Harm^j(E)$.
(See \S\ref{subsection;22.4.2.41}.)
It is equipped with a Hermitian sesqui-linear duality
$a_{\dagger}\nbigs$.
Let $\omegatilde$ denote the K\"ahler form of $g_{X\setminus H}$.
We obtain the induced Lefschetz operator
$L_{\omegatilde}$ of
$\bigoplus_k\nbightilde^k$.
As in \cite[Hodge-Simpson theorem]{sabbah2} in the smooth case,
$(\bigoplus_k\nbightilde^k,L_{\omegatilde},a_{\dagger}(\nbigs))$
is a polarized graded Lefschetz twistor structure of
weight $\dim X$ and type $1$.

By using Theorem \ref{thm;22.4.2.50},
we obtain an isomorphism between two graded $\nbigr$-triples
\[
 \Bigl(
 \bigoplus a_{X\dagger}^j\gbigt(E),
 L_{\omega},a_{X\dagger}(\nbigs)
 \Bigr)
 \simeq
  \Bigl(
 \bigoplus \nbightilde^j,
 L_{\omegatilde},a_{X\dagger}(\nbigs)
 \Bigr).
\]
Hence, we obtain the following.

\begin{thm}[Theorem
\ref{thm;22.3.15.20}]
\label{thm;22.4.3.10}
The tuple
$(\bigoplus a^j_{X\dagger}\gbigt(E),L_{\omega},
a_{X\dagger}(\nbigs))$
is a polarized graded Lefschetz twistor structure
of weight $\dim X$ and type $1$.
 \end{thm}
Namely,
the morphisms
$L_{\omega}^j:a_{X\dagger}^{-j}(\gbigt(E))
\lrarr a_{X\dagger}^{j}(\gbigt(E))\otimes\newTate(j)$
are isomorphisms for $j\geq 0$,
and
$(-1)^j
a_{X\dagger}(\nbigs)\circ L_{\omega}^j$
induces a polarization of
the primitive part
\[
Pa_{X\dagger}^{-j}(\gbigt(E))
=\Ker\Bigl(
L_{\omega}^{j+1}:
 a_{X\dagger}^{-j}(\gbigt(E))
 \lrarr
  a_{X\dagger}^{j+2}(\gbigt(E))\otimes\newTate(j+1)
\Bigr). 
\]
(We note that
the signature of a Lefschetz operator such as  $L_{\omega}$
is chosen to be opposite to that in \cite{sabbah2},
hence the above positivity condition is equivalent to
that in \cite{sabbah2}.)

We also revisit the theorem
of Kashiwara and Kawai
(Theorem \ref{thm;22.4.10.11}).
Let $(V=\bigoplus_{p+q=w}V^{p,q},\nabla,\langle\cdot,\cdot\rangle)$
be a polarized variation of Hodge structure.
Then, $(\nbigv,\DD)$ is naturally $\cnum^{\ast}$-equivariant.
(See \S\ref{subsection;22.4.10.2}.)
It is inherited to
the $(\nbigo_{\cnum_{\lambda}})_{\nbigx}$-complexes
$\gbigv\otimes p_{\lambda}^{\ast}\Omega^{\bullet}_{\nbigx}$
and
$\nbigc^{\bullet}_{L^2}(\nbigv,\DD,h)$.
The morphisms interpolating two complexes
are $\cnum^{\ast}$-equivariant by the construction.
There exists a natural
$\cnum^{\ast}$-action on
$\Harm^k(V,\nabla,h,g_{X\setminus H})$
by the Hodge decomposition (\ref{eq;22.4.12.2}).
By the construction, the isomorphisms
\begin{equation}
\label{eq;22.4.12.3}
 R^k(\id_X\times a_X)_{\ast}
 \Bigl(
 \gbigv\otimes p_{\lambda}^{\ast}\Omega^{\bullet}_X
 \Bigr)
 \simeq
 R^k(\id_X\times a_X)_{\ast}
 \nbigc^{\bullet}_{L^2}(\nbigv,\DD,h)
 \simeq
 \Harm^k(V,\nabla,h,g_{X\setminus H})
 \otimes\nbigo_{\cnum}
\end{equation}
are $\cnum^{\ast}$-equivariant.
We also note that
$\gbigv$ is naturally isomorphic to
$\gbigv\otimes p_{\lambda}^{\ast}\Omega^{\bullet}_X$
is naturally isomorphic
to the Rees module of the filtered de Rham complex
$(V_{\min}\otimes\Omega^{\bullet}_X,F)$.
Because 
$R^k(\id_X\times a_X)_{\ast}
 \Bigl(
 \gbigv\otimes p_{\lambda}^{\ast}\Omega^{\bullet}_X
 \Bigr)$ 
is locally free $\nbigo_{\cnum_{\lamda}}$-module,
we obtain that
$H^k(X,F_p(V_{\min}\otimes\Omega^{\bullet}_X))
\to
 H^k(X,V_{\min}\otimes\Omega^{\bullet}_X)$
are injective for any $k$ and $p$.
The left hand side of (\ref{eq;22.4.12.3})
is the Rees module of $H^k(X,V_{\min}\otimes\Omega^{\bullet}_X)$
with the induced filtration.
The right hand side is isomorphic to
the Rees module of $H^k(X,V_{\min}\otimes\Omega^{\bullet}_X)$
with the Hodge filtration induced by the $L^2$-theory of
Cattani-Kaplan-Schmid and Kashiwara-Kawai.
Because (\ref{eq;22.4.12.3}) is
a $\cnum^{\ast}$-equivariant isomorphism,
we obtain that the two filtrations are the same.
In this way, we revisit theorem
\cite[Theorem 1]{Kashiwara-Kawai-Hodge-holonomic}
of Kashiwara and Kawai.
(See \S\ref{subsection;22.4.10.2} for more detailed explanation.)

\subsubsection{Hard Lefschetz Theorem
for regular pure twistor $\nbigd$-modules}

We generalize the Hard Lefschetz theorem
for regular pure twistor $\nbigd$-modules.
Let $X$ be a K\"ahler manifold.
Let $f:X\lrarr Y$ be a morphism of complex manifolds.
Let $(\nbigt,\nbigs)$ be
a regular pure twistor $\nbigd$-module of weight $w$
on $X$ whose support is proper over $Y$.
Let $\omega$ be a K\"ahler class.
We obtain the $\nbigr_Y$-triples
$f_{\dagger}^j(\nbigt)$ $(j\in\seisuu)$.
We have the induced morphisms
$L_{\omega}:f_{\dagger}^j(\nbigt)
\lrarr
f_{\dagger}^{j+2}(\nbigt)\otimes\newTate(1)$.
We also obtain the morphisms
$f_{\dagger}(\nbigs):
f_{\dagger}^j(\nbigt)
\lrarr
\bigl(
f_{\dagger}^{-j}(\nbigt)
\bigr)^{\ast}\otimes\newTate(-w)$.
We shall prove
the following type of Hard Lefschetz theorem
for regular pure twistor $\nbigd$-modules.

\begin{thm}[Theorem
\ref{thm;22.3.5.1}]
\label{thm;22.3.10.20}
 The tuple
 $(\bigoplus f_{\dagger}^j\nbigt,L_{\omega},f_{\dagger}\nbigs)$
is a polarized graded Lefschetz twistor $\nbigd$-module
of weight $w$ on $Y$.
Namely, the following holds.
\begin{description}
 \item[(a)]  Each $f_{\dagger}^{j}(\nbigt)$ is
	a pure twistor $\nbigd$-module of weight $w+j$.
 \item[(b)] For any $j\geq 0$,
       the induced morphism
       $L_{\omega}:f_{\dagger}^{-j}(\nbigt)
       \lrarr f_{\dagger}^j(\nbigt)\otimes\newTate(j)$
       is an isomorphism.
       Moreover,
       $(-1)^jf_{\dagger}(\nbigs)\circ L_{\omega}^j$
       is a polarization of
\[
       Pf_{\dagger}^{-j}(\nbigt)
       =\Ker\Bigl(
       L_{\omega}^{j+1}:
       f_{\dagger}^{-j}(\nbigt)
       \lrarr
       f_{\dagger}^{j+2}(\nbigt)\otimes\newTate(j+1)       
       \Bigr).
\]       
\end{description}
\end{thm}

If $f$ is a projective morphism and
if $\omega$ is the first Chern class of a relatively ample line bundle,
the claim of Theorem \ref{thm;22.3.10.20}
was proved in
the theory of pure twistor $D$-modules \cite{sabbah2,mochi2,Mochizuki-wild},
which is a twistor version \cite{s3} of
the Hard Lefschetz Theorem for 
polarized pure Hodge modules \cite{saito1,saito2} due to Saito.
It was established by Simpson in \cite{Simpson-family}
for the push-forward of harmonic bundles
by smooth proper morphisms of K\"ahler manifolds
which are not necessarily projective.
Theorem \ref{thm;22.3.10.20} is a generalization
which allows us to study
the push-forward of regular polarized pure twistor $\nbigd$-modules
by non-projective but proper morphisms
in the context of K\"ahler manifolds.
Saito \cite{MSaito-proper} established
the Hard Lefschetz theorem
for constant Hodge modules
by proper morphisms of K\"ahler manifolds.
Saito \cite{MSaito-proper}
also explained, in the case of Hodge modules,
how to deduce (a) in Theorem \ref{thm;22.3.10.20}
from the theorem of Kashiwara and Kawai.
(See \cite{MSaito-proper} for more general and precise statements.)

\begin{rem}
After almost completing this manuscript,
the author received
several preliminary versions of {\rm\cite{MSaito-2022}}.
\hfill\qed
\end{rem}

\subsection{Related recent works}

Though there are many related important works,
we just mention two recent stimulating studies
\cite{Shentu-Zhao} and \cite{Wei-Yang},
where the previous works and the backgrounds are explained in detail.
(See also an excellent survey paper \cite{Williamson-Bourbaki}
for the decomposition theorem of perverse sheaves
of geometric origin,
in particular the proof due to de Cataldo and Migliorini
\cite{de-Cataldo-Migliorini-semismall, de-Cataldo-Migliorini-decomposition}.)

In \cite{Shentu-Zhao},
Shentu and Zhao developed the $L^2$-theory
for polarized variation of Hodge structure
with singularity in a general situation,
and they proved the Hard Lefschetz theorem
for the associated intersection cohomology group.
Let $X$ be a K\"ahler manifold.
Let $Z_1\subset Z_0$ be closed complex analytic subvarieties of $X$
such that $Z_0\setminus Z_1$ is a locally closed complex submanifold of $X$.
Let $(V=\bigoplus V^{p,q},\nabla,\langle\cdot,\cdot\rangle)$
be a polarized variation of pure Hodge structure
of weight $w$ on $Z_0\setminus Z_1$.
Let $h$ denote the Hodge metric.
Shentu and Zhao proved that there exists
a complete K\"ahler metric $g_{Z_0\setminus Z_1}$
of $Z_0\setminus Z_1$ satisfying some good properties,
called a distinguished metric,
such that the $L^2$-complex of $(V,\nabla)$
with respect to $h$ and $g_{Z_0\setminus Z_1}$
is isomorphic to the intersection complex of $(V,\nabla)$.
It particularly implies that
the Hard Lefschetz Theorem for the intersection cohomology group
of $(V,\nabla)$.
They also obtained a similar result for $0$-tame harmonic bundles.

Both of the result of {\rm\cite{Shentu-Zhao}} and
Theorem {\rm\ref{thm;22.3.10.20}}
imply the Hard Lefschetz theorem for the intersection cohomology group
of a polarizable variation of Hodge structure
or a $0$-tame harmonic bundle.
However,
there are important difference in the more precise statements,
depending on their purpose and ours.
On one hand,
Shentu and Zhao obtained the comparison with the $L^2$-complex
and the intersection complex for $0$-tame harmonic bundles
on general $(Z_0,Z_1)$,
which is expected to allow us to apply the $L^2$-theory
efficiently in a general situation.
On the other hand,
though we study the $L^2$-complex only in the normal crossing singular case,
we also show that
the induced polarized Hodge structure
or the induced polarized pure twistor structure on the cohomology groups
are the same as those constructed in terms of pure twistor $\nbigd$-modules.
It allows us to enhance the theory of pure twistor $\nbigd$-modules
in the non-algebraic but K\"ahler situation,
and to obtain the Hard Lefschetz Theorem in the relative case
not only in the absolute case, as in Theorem {\rm\ref{thm;22.3.10.20}}.

In \cite{Wei-Yang},
Wei and Yang obtained another proof of
the Hard Lefschetz Theorem
for semisimple local systems on projective manifolds
with respect to a proper morphism
between projective manifolds,
which was originally proved by Sabbah \cite{sabbah2}.
They generalized the argument of
de Cataldo and Migliorini in the Hodge case
\cite{de-Cataldo-Migliorini-semismall, de-Cataldo-Migliorini-decomposition}
to the twistor case,
by following the idea of Simpson's Meta Theorem \cite{s3},
instead using the theory of pure twistor $\nbigd$-modules.
Stimulated by their work,
the author reminded an old question
to ask whether we could generalize
the method in \cite{de-Cataldo-Migliorini-decomposition}
to a more singular situation.
It motivated the author to study
the $L^2$-complex of tame harmonic bundles with normal crossing singularity
though he chose in this paper
to apply it to enhance the theory of pure twistor $\nbigd$-modules.

\subsection{An outline of the proof in the two dimensional case}

\subsubsection{The issue to be considered}

The natural morphism 
$\bigl(
\gbige\otimes p_{\lambda}^{\ast}\Omega^{\bullet}
\bigr)_{|\nbigx\setminus\nbigh}
\lrarr
\nbigc^{\bullet}_{L^2}(\nbige,\DD,h)_{|\nbigx\setminus\nbigh}$
does not extend to a morphism
$\bigl(
\gbige\otimes p_{\lambda}^{\ast}\Omega^{\bullet}
\bigr)
\lrarr
\nbigc^{\bullet}_{L^2}(\nbige,\DD,h)$
unless $(E,\delbar_E,\theta,h)$ extends
to a harmonic bundle on $X$.
Therefore, we need to find
a complex of $(\nbigo_{\cnum_{\lambda}})_{\nbigx}$-modules
$\nbigc^{\bullet}$
with quasi-isomorphisms
$\nbigc^{\bullet}\to
\gbige\otimes\Omega^{\bullet}_X$
and
$\nbigc^{\bullet}\to
\nbigc^{\bullet}_{L^2}(\nbige,\DD,h)$
for the proof of Theorem \ref{thm;22.4.2.50}.
We need to consider a similar problem
for the specialization to $\lambda$
(Theorem \ref{thm;22.4.2.40}).

For a polarized variation of Hodge structure
on a punctured curve,
this is solved by Zucker \cite{z} in an ideal way.
It was generalized for tame harmonic bundles
on punctured curves in \cite{sabbah2,mochi2}
by generalizing and refining the arguments in \cite{z}.
Roughly,
if we consider the sheaf of $L^2$ and holomorphic sections
of $\nbige\otimes p_{\lambda}^{\ast}\Omega^{\bullet}$,
it is naturally quasi-isomorphic to the both complexes.
More precisely,
for an open subset $U\subset\nbigx$,
let $\nbigc^{\bullet}_{L^2,\hol}(\nbige,\DD,h)$
denote the space of
holomorphic sections $\tau$ of
$\nbige\otimes p_{\lambda}^{\ast}\Omega^{\bullet}$
on $U\setminus\nbigh$
such that
(i) $\delbar_{\nbige}\tau=\del_{\lambdabar}\tau=0$,
(ii) $\tau$ and $\DD\tau$ are $L^2$
with respect to $h$ and $g_{X\setminus H}+d\lambda\,d\lambdabar$
locally around any point of $U$.
It is clearly an $(\nbigo_{\cnum})_{\nbigx}$-subcomplex of
both
$\nbigc^{\bullet}_{L^2}(\nbige,\DD,h)$
and
$\gbige\otimes p_{\lambda}^{-1}(\Omega_X^{\bullet})$.
By using a generalization of Dolbeault's lemma
to the singular case in \cite{z},
and by using some additional argument,
we can prove that the inclusion
$\nbigc^{\bullet}_{L^2,\hol}(\nbige,\DD,h)
\lrarr
\nbigc^{\bullet}_{L^2}(\nbige,\DD,h)$
is a quasi-isomorphism.
Moreover, we can explicitly describe 
$\nbigc^{\bullet}_{L^2,\hol}(\nbige,\DD,h)$
in terms of
the parabolic structure and the weight filtration
at each point of $\nbigh$,
we can prove that 
$\nbigc^{\bullet}_{L^2,\hol}(\nbige,\DD,h)
\lrarr
\gbige\otimes p_{\lambda}^{-1}(\Omega_X^{\bullet})$
is a quasi-isomorphism.

However, in the higher dimensional case,
in general,
this argument cannot work
as pointed out by Cattani and Kaplan \cite{ck2}.
We can similarly define the complex
$\nbigc^{\bullet}_{L^2,\hol}(\nbige,\DD,h)$
even in the higher dimensional case,
but both the morphisms
\[
\nbigc^{\bullet}_{L^2,\hol}(\nbige,\DD,h)
\lrarr
\nbigc^{\bullet}_{L^2}(\nbige,\DD,h),
\quad
\nbigc^{\bullet}_{L^2,\hol}(\nbige,\DD,h)
\lrarr
\gbige\otimes p_{\lambda}^{-1}(\Omega_X^{\bullet})
\]
are not necessarily quasi-isomorphisms.
Specializations at $\lambda$ are not necessarily
quasi-isomorphisms, too.
(See \S\ref{subsection;22.3.18.30}.)

\vspace{.1in}

In the theorem of Cattani-Kaplan-Schmid
and Kashiwara-Kawai
(Theorem \ref{thm;22.4.3.2}),
we need to compare the complexes
$V_{\min}\otimes\Omega^{\bullet}_X$
and
$\nbigc^{\bullet}_{L^2}(V,\nabla,h)$.
In this case,
$V_{\min}\otimes\Omega^{\bullet}_X$
is the intersection complex.
By setting
$H_I^{\circ}=H_I\setminus \bigcup_{i\not\in I}H_i$,
we obtain the stratification
$X=\coprod H_I^{\circ}$.
The intersection complex
is characterized by the cohomological property
along the strata.
Hence, it is reduced to check whether
the complex $\nbigc^{\bullet}_{L^2}(V,\nabla,h)$ also
satisfies the conditions to characterize the intersection complex.
It is still a difficult problem,
but solved by Cattani-Kaplan-Schmid and Kashiwara-Kawai.

However, for the comparison of 
$\gbige\otimes p_{\lambda}^{\ast}\Omega^{\bullet}_X$
and 
$\nbigc^{\bullet}_{L^2}(\nbige,\DD,h)$,
such a helpful characterization
of the complex $\gbige\otimes p_{\lambda}^{\ast}\Omega^{\bullet}_X$
is not known around $\lambda=0$,
even in the Hodge case.
Therefore, we need to construct
a sequence of complexes and quasi-isomorphisms
to connect the complexes.

\subsubsection{The suggestion of Kashiwara and Kawai}

As mentioned,
Theorem \ref{thm;22.4.3.10} is
a twistor version of the theorem
\cite[Theorem 1]{Kashiwara-Kawai-Hodge-holonomic}
of Kashiwara and Kawai.
Here, let us mention what is suggested in the last paragraph
of \cite{Kashiwara-Kawai-Hodge-holonomic}.

Let $V$ be a vector space equipped with
two nilpotent endomorphisms $N_1$ and $N_2$
which are commutative.
We consider the two dimensional vector space
$\cnum e_1\oplus\cnum e_2$.
For $k\in\seisuu$,
we set
$\nbigc_1^{k}:=
V\otimes\bigwedge^{k+2}(\cnum e_1\oplus\cnum e_2)$.
The morphisms $N_i$ together with
the exterior product of $e_i$
induces the differential $N_1e_1+N_2e_2$
on $\nbigc_1^{\bullet}$:
\[
 \nbigc^{\bullet}_1:
 V\lrarr
 V\cdot e_1
 \oplus
 V\cdot e_2
 \lrarr
 V\cdot e_1\wedge e_2.
\]
We obtain the following subcomplex of $\nbigc^{\bullet}_1$:
\[
\nbigc^{\bullet}_0:
 V\lrarr
 \Image(N_1)\cdot e_1
 \oplus
 \Image(N_2)\cdot e_2
 \lrarr
 \Image (N_1N_2)\cdot e_1\wedge e_2.
\]
\begin{rem}
Let $X$ be a neighbourhood of $O=(0,0)$ in $\cnum^2$.
We set $H=X\cap\{z_1z_2=0\}$.
 We consider the regular singular meromorphic flat bundle
$(\nbigv,\nabla)$ on $(X,H)$,
where 
$\nbigv=V\otimes\nbigo_{X}(\ast H)$
and $\nabla=d+\sum_{i=1,2}(N_i\,dz_i/z_i)$.
Then, $\nbigc^{\bullet}_1$ is naturally quasi-isomorphic to
the stalk of 
$\nbigv\otimes\Omega^{\bullet}_{X}[2]$
at $O$,
and $\nbigc^{\bullet}_0$ is naturally quasi-isomorphic to
 the stalk of
$\nbigv_{\min}\otimes\Omega^{\bullet}_X[2]$ at $O$.
\hfill\qed
\end{rem}

We have the monodromy weight filtrations $W(N_i)$ of $V$.
Namely, $W(N_i)=(W_k(N_i)\,|\,k\in\seisuu)$
is an increasing filtration such that
(i) $W_{-m_0}(N_i)=0$ and $W_{m_0}(N_i)=V$ for a sufficiently large $m_0$,
(ii) $N_i\cdot W_k(N_i)\subset W_{k-2}(N_i)$,
(iii) $N_i^k:\Gr^{W(N_i)}_k\simeq \Gr^{W(N_i)}_{-k}$ for any $k\geq 0$.
We also have the monodromy weight filtration $W(N_1+N_2)$
of $N_1+N_2$.
If $(V,N_1,N_2)$ comes from a polarized mixed Hodge structure,
then $W(N_1+N_2)$ is equal to 
$W(a_1N_1+a_2N_2)$ for any $(a_1,a_2)\in\real_{>0}^2$,
and $W(N_1+N_2)$ is a relative monodromy filtration of
$(V,W(N_1),N_2)$
and $(V,W(N_2),N_1)$.
We assume the conditions.
In particular,
$N_iW_k(N_1+N_2)\subset W_{k-2}(N_1+N_2)$
for $i=1,2$.
We note
$N_1W_k(N_2)\subset W_k(N_2)$
and
$N_2W_k(N_1)\subset W_k(N_1)$.
To simplify the description,
we set $W_k'=W_k(N_1)$,
$W_k''=W_k(N_2)$
and $W_k(N_1+N_2)$.

We consider the following subcomplexes of $\nbigc^{\bullet}_1$:
{\normalsize
\[
W_1\nbigc^{\bullet}_1:
 W_1(V)
 \lrarr
 W_{-1}(V)e_1
 \oplus
 W_{-1}(V)e_2
 \lrarr
 W_{-3}(V)\,e_1\wedge e_2,
\]
\[
(W_1\cap W'_0)
\nbigc^{\bullet}_1:
\bigl(
 W_1
 \cap
 W'_{0}
\bigr)(V)
 \lrarr
 \bigl(
 W_{-1}
 \cap
 W'_{-2}
 \bigr)(V)e_1
 \oplus
  \bigl(
 W_{-1}
 \cap
 W'_{0}
 \bigr)(V)e_2
 \lrarr
 \bigl(
 W_{-3}
 \cap
 W'_{-2}
 \bigr)(V) e_1\wedge e_2,
\]
\[
(W_1\cap W''_0)
 \nbigc^{\bullet}_1:
 \bigl(
 W_1
 \cap
 W''_{0}
\bigr)(V)
 \lrarr
 \bigl(
 W_{-1}
 \cap
 W''_{0}
 \bigr)(V)e_1
 \oplus
  \bigl(
 W_{-1}
 \cap
 W''_{-2}
 \bigr)(V)e_2
 \lrarr
 \bigl(
 W_{-3}
 \cap
 W''_{-2}
 \bigr)(V) e_1\wedge e_2.
\]
}
There exist the inclusions
$(W_1\cap W_0')\nbigc^{\bullet}_1\to W_1\nbigc^{\bullet}_1$
and $(W_1\cap W_0'')\nbigc^{\bullet}_1\to W_1\nbigc^{\bullet}_1$.
We obtain the following complex:
\[
\Cone\Bigl(
(W_1\cap W_0')\nbigc^{\bullet}_1
\oplus
(W_1\cap W_0'')\nbigc^{\bullet}_1
 \lrarr
 W_1\nbigc^{\bullet}_1
 \Bigr)[1].
\]
(Here $\Cone(A\to B)$ denotes the mapping cone of $A\to B$.
See \cite[\S1.4]{Kashiwara-Schapira} for the mapping cone.)

In the last paragraph of \cite{Kashiwara-Kawai-Hodge-holonomic},
when $(V,N_1,N_2)$ is induced by a polarized mixed Hodge structure,
Kashiwara and Kawai suggested
that there exists an isomorphism
between the two complexes
\begin{equation}
\label{eq;22.4.3.20}
\nbigc_0^{\bullet}\simeq
\Cone\Bigl(
(W_1\cap W_0')\nbigc^{\bullet}_1
\oplus
(W_1\cap W_0'')\nbigc^{\bullet}_1
\to
W_1\nbigc^{\bullet}_1
\Bigr)[1]
\end{equation}
in the derived category of complexes of $\cnum$-vector spaces,
i.e.,
a sequence of quasi-isomorphisms
connecting the both sides of (\ref{eq;22.4.3.20}).
\begin{rem}
As mentioned,
$\nbigc^{\bullet}_0$ is closely related with
the intersection complex of
$(\nbigv,\nabla)$.
The right hand side of {\rm(\ref{eq;22.4.3.20})}
is related with the $L^2$-complex
of $(\nbigv,\nabla)$ with the Hodge metric
if $(\nbigv,\nabla)$ is induced by
a polarized variation of Hodge structure
with unipotent local monodromy.
Indeed,
when we regard $s\in V$ as a holomorphic section of $\nbigv$,
the condition $s\in W_1$
is related with
the condition that $s$ is $L^2$ on
 $\{-C_1\log|z_1|<-\log|z_2|<-C_2\log|z_1|\}$
for $0<C_1<C_2$.
The condition $s\in W_1\cap W'_0$
is related with
the condition that $s$ is $L^2$ on
$\{-\log|z_2|<-C \log|z_1|\}$ for $0<C$.
We have similar claims for
$s\,dz_i/z_i$ or $s\,dz_1\,dz_2/(z_1z_2)$.
(See {\rm\S\ref{subsection;22.3.17.51}}.)
\hfill\qed 
\end{rem}

It would be instructive
to recall a similar complexes in the one dimensional case.
For a vector space $V$ with a nilpotent endomorphism $N$,
the complex $V\to \Image (N)$
is related with the intersection complex,
and
the complex 
$W_0(N)\to W_{-2}(N)$
is related with the $L^2$-holomorphic complex.
In this case,
there exists the natural inclusion
\[
 \Bigl(
 W_0(N)\to W_{-2}(N)
 \Bigr)
 \lrarr
 \Bigl(
 V\to \Image (N)
 \Bigr),
\]
and it is a quasi-isomorphism.
Kashiwara and Kawai suggested that 
there exists
analogue but more complicated quasi-isomorphisms
even in the higher dimensional case.

\subsubsection{Some quasi-isomorphisms to obtain (\ref{eq;22.4.3.20})}
\label{subsection;22.4.4.2}

We assume that $(V,N_1,N_2)$ underlies a polarized mixed Hodge structure.
Let us explain how to obtain (\ref{eq;22.4.3.20}).
Though we shall study the issue in a general situation
(Theorem \ref{thm;22.2.8.7}),
it would be useful to give an explanation
in this simpler case.
By using the vanishing cycle theorem (see \S\ref{subsection;22.4.3.31}),
we obtain the following subcomplexes of $\nbigc^{\bullet}_0$:
{\normalsize
\[
 (W_1\cap W_0')\nbigc_0^{\bullet}:
 (W_1\cap W_0')(V)
 \lrarr
 (W_0\cap W'_{-1})(\Image N_1)e_1
 \oplus
 (W_0\cap W'_{0})(\Image N_2)e_2
 \lrarr
 (W_{-1}\cap W_{-1}')(\Image N_1N_2)e_1\wedge e_2,
\]
\[
 (W_1\cap W_0'')\nbigc_0^{\bullet}:
 (W_1\cap W_0'')(V)
 \lrarr
 (W_0\cap W''_0)(\Image N_1)e_1
 \oplus
 (W_0\cap W''_{-1})(\Image N_2)e_2
 \lrarr
 (W_{-1}\cap W''_{-1})(\Image N_1N_2)e_1\wedge e_2,
\]
\[
 W_1\nbigc_0^{\bullet}:
 W_1(V)
 \lrarr
 W_0(\Image N_1)e_1
 \oplus
 W_0(\Image N_2)e_2
 \lrarr
 W_{-1}(\Image N_1N_2)e_1\wedge e_2.
\]
}
Here,
$W_{\bullet}(\Image N_j)$
denote the monodromy weight filtration of
$N_1+N_2$ on $\Image (N_j)$.
The notation
$W'_{\bullet}(\Image N_j)$,
$W''_{\bullet}(\Image N_j)$, etc.,
are used in similar ways.

Let us observe that 
the natural inclusions
\[
 (W_1\cap W_0')\nbigc_0^{\bullet}
 \stackrel{a_1}{\lrarr}
  \nbigc_0^{\bullet},
  \quad
 (W_1\cap W_0'')\nbigc_0^{\bullet}
 \stackrel{a_2}{\lrarr}
  \nbigc_0^{\bullet},
  \quad
  W_1\nbigc_0^{\bullet}
  \stackrel{a_3}{\lrarr}
  \nbigc_0^{\bullet}
\]
are quasi-isomorphisms.
Indeed,
$a_3$ is a quasi-isomorphism
by the purity theorem
for polarized mixed Hodge structure
(see Proposition \ref{prop;22.2.6.3} and Lemma \ref{lem;22.2.7.1}).
Because the inclusions
{\normalsize
\[
\bigl(
W_0'(V)\to W_{-1}'(\Image N_1)
\bigr)
\to
\bigl(
V\to \Image N_1
\bigr),
\quad\quad
\bigl(
W_0'(\Image N_2)\to W_{-1}'(\Image N_1N_2)
\bigr)
\to
\bigl(
\Image N_2\to \Image N_1N_2
\bigr)
\]
}
are quasi-isomorphisms,
the natural inclusion
\[
\Bigl(
 W_0'(V)
 \to
 W_{-1}'(\Image N_1)
 \oplus
 W_0'(\Image N_2)
 \to
 W_{-1}'(\Image N_1N_2)
 \Bigr)
 \lrarr
 \Bigl(
 V
 \to
 \Image N_1
 \oplus
 \Image N_2
 \to
 \Image N_1N_2
 \Bigr)
\]
is a quasi-isomorphism.
By using the truncation by the weight filtration
of the mixed twistor structures,
we obtain that $a_1$ is a quasi-isomorphism.
(See Proposition \ref{prop;22.2.7.3}
for a more general case.)

\vspace{.1in}

There exist natural morphisms
$(W_1\cap W_0')\nbigc_0^{\bullet}
\lrarr
W_1\nbigc_0^{\bullet}$
and
$(W_1\cap W_0'')\nbigc_0^{\bullet}
\lrarr
W_1\nbigc_0^{\bullet}$.
We obtain the following complex:
\[
 \Cone\Bigl(
 (W_1\cap W_0')\nbigc_0^{\bullet}
 \oplus
 (W_1\cap W_0'')\nbigc_0^{\bullet}
 \lrarr
 W_1\nbigc_0^{\bullet}
 \Bigr)[1].
\]
We obtain the following natural quasi-isomorphisms:
\[
 \nbigc^{\bullet}_0\lrarr
 \Cone\Bigl(
 \nbigc^{\bullet}_0
 \oplus
 \nbigc^{\bullet}_0
 \lrarr
 \nbigc^{\bullet}_0
 \Bigr)[1]
 \llarr
 \Cone\Bigl(
 (W_1\cap W_0')\nbigc_0^{\bullet}
 \oplus
 (W_1\cap W_0'')\nbigc_0^{\bullet}
 \lrarr
 W_1\nbigc_0^{\bullet}
 \Bigr)[1]. 
\]
Thus, we may replace $\nbigc_0$
with the right hand side.

\vspace{.1in}

By the vanishing cycle theorem
(see \S\ref{subsection;22.4.3.31}),
there exist natural inclusions:
\begin{equation}
 (W_1\cap W_0')
 \nbigc_0^{\bullet}
 \lrarr
 (W_1\cap W_0')
 \nbigc_1^{\bullet},
 \quad
 (W_1\cap W_0'')
 \nbigc_0^{\bullet}
 \lrarr
 (W_1\cap W_0'')
 \nbigc_1^{\bullet},
 \quad
 W_1\nbigc_0^{\bullet}
 \lrarr
 W_1\nbigc_1^{\bullet}.
\end{equation}
We obtain the following morphism:
\begin{equation}
 \label{eq;22.4.3.32}
 \Cone\Bigl(
 (W_1\cap W_0')
 \nbigc_0^{\bullet}
 \oplus
  (W_1\cap W_0'')
  \nbigc_0^{\bullet}
  \lrarr
   W_1\nbigc_0^{\bullet}
   \Bigr)[1]
   \lrarr
   \Cone\Bigl(
   (W_1\cap W_0')\nbigc_1^{\bullet}
   \oplus
   (W_1\cap W_0'')\nbigc_1^{\bullet}
   \lrarr
   W_1\nbigc_1^{\bullet}
 \Bigr).
\end{equation}
\begin{description}
 \item[Claim]
	    The morphism (\ref{eq;22.4.3.32}) is a quasi-isomorphism.
\end{description}
To prove the claim,
it is convenient to introduce some more complexes.
We consider the following complexes of $\nbigc_1^{\bullet}$:
\[
 \nbigc_2^{\bullet}:
 V\lrarr
 \Image(N_1)e_1
  \oplus
  V\cdot e_2
  \lrarr
  \Image(N_1) e_1\wedge e_2,
\]
\[
 \nbigc_3^{\bullet}:
 V\lrarr
 V\cdot e_1
  \oplus
 \Image(N_2)\cdot e_2
  \lrarr
  \Image(N_2) e_1\wedge e_2.
\]
We obtain the following complex
\begin{equation}
\label{eq;22.4.3.33}
 \Cone\Bigl(
 \nbigc_2^{\bullet}
 \oplus
 \nbigc_3^{\bullet}
 \lrarr
 \nbigc_1^{\bullet}
 \Bigr)[1].
\end{equation}
\begin{rem}
Let $j:X\setminus\{(0,0)\}\lrarr X$.
Then, {\rm(\ref{eq;22.4.3.33})} represents
the stalk of
$Rj_{\ast}j^{-1}
(\nbigv_{\min}\otimes\Omega^{\bullet}_X)$
at $(0,0)$.
\hfill\qed
\end{rem}

There exists the following natural morphism,
which is not necessarily a quasi-isomorphism:
\[
\Cone\Bigl(
\nbigc_0^{\bullet}\oplus
\nbigc_0^{\bullet}
\lrarr
\nbigc_0^{\bullet}
\Bigr)[1]
 \lrarr
 \Cone\Bigl(
 \nbigc_2^{\bullet}
 \oplus
 \nbigc_3^{\bullet}
 \lrarr
 \nbigc_1^{\bullet}
 \Bigr)[1].
\]
We consider the following subcomplexes of
$\nbigc_2^{\bullet}$ and $\nbigc_3^{\bullet}$, respectively:
\[
 W_1\nbigc_2^{\bullet}:
 W_1(V)\lrarr
 W_0( \Image N_1)e_1
  \oplus
  W_{-1}(V)\cdot e_2
  \lrarr
  W_{-2}(\Image(N_1)) e_1\wedge e_2,
\]
\[
 W_1\nbigc_3^{\bullet}:
 W_1(V)\lrarr
 W_{-1}(V)e_1
  \oplus
  W_{0}(\Image N_2)\cdot e_2
  \lrarr
  W_{-2}(\Image(N_2)) e_1\wedge e_2.
\]
By using the purity theorem and its dual statement,
we can prove that the following morphism is a quasi-isomorphism
(see Corollary \ref{cor;22.2.6.10}):
\begin{equation}
\label{eq;22.4.3.40}
\Cone\Bigl(
W_1\nbigc_0^{\bullet}\oplus
W_1\nbigc_0^{\bullet}
\lrarr
W_1\nbigc_0^{\bullet}
\Bigr)[1]
 \lrarr
 \Cone\Bigl(
 W_1\nbigc_2^{\bullet}
 \oplus
 W_1\nbigc_3^{\bullet}
 \lrarr
 W_1\nbigc_1^{\bullet}
 \Bigr).
\end{equation}
Moreover,
we consider the following subcomplexes
of $W_1\nbigc_2^{\bullet}$
and $W_1\nbigc_3^{\bullet}$, respectively:
{\normalsize
\[
 (W_1\cap W_0')\nbigc_2^{\bullet}:
 (W_1\cap W_0')(V)
 \lrarr
 (W_0\cap W_{-1}')(\Image N_1)e_1
 \oplus
 (W_{-1}\cap W_{0}')(V)e_2
 \lrarr
 (W_{-2}\cap W_{-1}')(\Image N_1)e_1\wedge e_2,
\]
\[
 (W_1\cap W_0'')\nbigc_3^{\bullet}:
 (W_1\cap W_0'')(V)
 \lrarr
 (W_{-1}\cap W_{0}'')(V)e_1
 \oplus
 (W_0\cap W_{-1}'')(\Image N_2)e_2
 \lrarr
 (W_{-2}\cap W_{-1}'')(\Image N_2)e_1\wedge e_2.
\]
}
We can observe that
$(W_1\cap W_0')\nbigc_2^{\bullet}
 \lrarr
 W_1\nbigc_2^{\bullet}$
is a quasi-isomorphism.
Indeed, because the natural inclusion
\[
\Bigl(
W_0'(V)\to W_{-1}'(\Image N_1)
\Bigr)
\lrarr
\Bigl(
V\to \Image N_1
\Bigr)
\]
is a quasi-isomorphism of complexes of mixed twistor structures,
the inclusions
\[
 \Bigl(
(W_1\cap W_0')(V)\to (W_0\cap W_{-1}')(\Image N_1)
\Bigr)
\lrarr
\Bigl(
W_1(V)\to W_{0}(\Image N_1)
\Bigr),
\]
\[
 \Bigl(
(W_{-1}\cap W_0')(V)\to (W_0\cap W_{-1}')(\Image N_1)
\Bigr)
\lrarr
\Bigl(
W_{-1}(V)\to W_{0}(\Image N_1)
\Bigr)
\]
are quasi-isomorphisms
by the vanishing cycle theorem.
Then, we can easily observe
$(W_1\cap W_0')\nbigc_2^{\bullet}
 \lrarr
 W_1\nbigc_2^{\bullet}$
 is a quasi-isomorphism.
Similarly,
the inclusions
\[
(W_1\cap W_0'')\nbigc_3^{\bullet}
\lrarr
W_1\nbigc_3^{\bullet},
\quad
(W_1\cap W_0')\nbigc_0^{\bullet}
\lrarr
W_1\nbigc_0^{\bullet},
\quad
(W_1\cap W_0'')\nbigc_0^{\bullet}
\lrarr
W_1\nbigc_0^{\bullet}
\]
are quasi-isomorphisms.
Hence,
we obtain the following commutative diagram:
{\normalsize
\[
 \begin{CD}
\Cone\Bigl(
(W_1\cap W_0')\nbigc_0\oplus
(W_1\cap W_0'')\nbigc_0
\lrarr
W_1\nbigc_0
\Bigr)[1]
@>{a}>>
 \Cone\Bigl(
 (W_1\cap W_0')\nbigc_2^{\bullet}
 \oplus
 (W_1\cap W_0'')\nbigc_3^{\bullet}
 \lrarr
 W_1\nbigc_1^{\bullet}
  \Bigr)[1]
  \\
@V{\simeq}VV @V{\simeq}VV \\
\Cone\Bigl(
W_1\nbigc_0\oplus
W_1\nbigc_0
\lrarr
W_1\nbigc_0
\Bigr)[1]
@>{\simeq}>>
 \Cone\Bigl(
 W_1\nbigc_2^{\bullet}
 \oplus
 W_1\nbigc_3^{\bullet}
 \lrarr
 W_1\nbigc_1^{\bullet}
  \Bigr)[1]
 \end{CD}
\]
}
We obtain, the morphism $a$ is also a quasi-isomorphism.

In the two dimensional case,
we have already obtained (\ref{eq;22.4.3.20}).
Indeed, because
\[
 W'_{-1}(\Image N_1)
 \simeq
 W'_{-2}(V),
 \quad\quad
 W''_{-1}(\Image N_2)
 \simeq
 W''_{-2}(V),
\]
we obtain
\begin{multline}
\label{eq;22.4.1.1}
 \Cone\Bigl(
 (W_1\cap W_0')\nbigc_2^{\bullet}
 \oplus
 (W_1\cap W_0'')\nbigc_3^{\bullet}
 \lrarr
 W_1\nbigc_1^{\bullet}
  \Bigr)[1]
 =
 \Cone\Bigl(
 (W_1\cap W_0')\nbigc_1^{\bullet}
 \oplus
 (W_1\cap W_0'')\nbigc_1^{\bullet}
 \lrarr
 W_1\nbigc_1^{\bullet}
  \Bigr)[1].
\end{multline}
Hence, we obtain the desired sequence of quasi-isomorphisms
connecting the both sides of
(\ref{eq;22.4.3.20}).
In the higher dimensional case,
(\ref{eq;22.4.1.1}) is not an isomorphism.
But, by using the result in the lower dimensional case,
we can prove that it is a quasi-isomorphism.
(See Theorem \ref{thm;22.2.8.7}.)
In this way,
we understand what Kashiwara and Kawai suggested.

\subsubsection{Complexes of sheaves}

Let $X$ be a neighbourhood of $O=(0,0)$ in $\cnum^2$,
and we set $H=X\cap\{z_1z_2=0\}$ and $H_i=X\cap\{z_i=0\}$.
Let $(E,\delbar_E,\theta,h)$ be a tame harmonic bundle
on $(X,H)$.
We explain what complexes we consider
to prove Theorem \ref{thm;22.4.2.40} with $\lambda=0$.

We obtain the filtered Higgs bundle
$(\nbigp_{\ast}E,\theta)$.
We obtain the endomorphisms $f_i$ of $\nbigp_{\ast}E$
by $\theta=f_1\,dz_1/z_1+f_2\,dz_2/z_2$.
For simplicity,
we assume that $\Res_i(\theta)=f_{i|H_i}$ are nilpotent,
and that the parabolic structure of
$\nbigp_{0,0}E$ is trivial,
i.e.,
$\nbigp_{a_1,a_2}E_{|O}=\nbigp_{0,0}E_{|O}$
for $-1<a_i\leq 0$.

We have the logarithmic complex
$\nbigc^{\bullet}_{\tw}(\gbige[\ast H],0):=
\nbigp_{0,0}E\otimes\Omega^{\bullet}_X(\log H)$.
For $K=O,H_1,H_2$,
we obtain the following complex on $K$,
where the differential is induced by
$f_{1|K}\,dz_1/z_1+f_{2|K}\,dz_2/z_2$:
\begin{equation}
 \nbigc^{\bullet}_{\tw}(\gbige[\ast H],0)_{|K}:=
  \nbigp_{0,0}E_{|K}
  \otimes
  \bigwedge^{\bullet}\Bigl(
  \cnum\,\frac{dz_1}{z_1}
  \oplus
  \cnum\,\frac{dz_2}{z_2}
\Bigr).
\end{equation}
There exists the natural morphism
\begin{equation}
\label{eq;22.4.4.1}
 \nbigp_{0,0}E\otimes\Omega^{\bullet}_X(\log H)
 \lrarr
 \bigoplus_{K=O,H_1,H_2}
 \nbigc^{\bullet}_{\tw}(\gbige[\ast H],0)_{|K}.
\end{equation}
When we are given subcomplexes
$\nbigj^{\bullet}_K\subset
\nbigc^{\bullet}_{\tw}(\gbige[\ast H],0)_{|K}$
$(K=O,H_1,H_2)$,
we obtain the subcomplex
$\nbigc^{\bullet}(\{\nbigj^{\bullet}_K\})
\subset
\nbigc^{\bullet}_{\tw}(\gbige[\ast H],0)$
as the inverse image of
$\bigoplus_{K=O,H_1,H_2} \nbigj_{K}$
by (\ref{eq;22.4.4.1}).

We obtain the subcomplex
$\nbigc^{\bullet}_{\tw}(\gbige,0)
=\nbigc^{\bullet}(\{\nbigj^{\bullet}_{0,K}\})$:
\[
 \nbigj_{0,O}:
\nbigp_{0,0}E_{|O}
       \lrarr
       \Image f_{1|O}\cdot \frac{dz_1}{z_1}
       \oplus
       \Image f_{2|O}\cdot \frac{dz_2}{z_2}
       \lrarr
       \Image (f_1f_2)_{|O}\cdot
       \frac{dz_1}{z_1}\frac{dz_2}{z_2},
\]
\[
 \nbigj_{0,H_1}:
 \nbigp_{0,0}E_{|H_1}
 \lrarr
 \Image f_{1|H_1}\cdot \frac{dz_1}{z_1}
 \oplus
 \nbigp_{0,0}E_{|H_1}\cdot dz_2
 \lrarr
 \Image f_{1|H_1}\cdot
 \frac{dz_1}{z_1}
 \frac{dz_2}{z_2},
\]
\[
 \nbigj_{0,H_2}:
 \nbigp_{0,0}E_{|H_2}
 \lrarr
 \nbigp_{0,0}E_{|H_2}\cdot dz_1
 \oplus
 \Image f_{2|H_2}\cdot \frac{dz_2}{z_2}
 \lrarr
 \Image f_{2|H_2}\cdot
 \frac{dz_1}{z_1}
 \frac{dz_2}{z_2}.
\]
The complex $\nbigc^{\bullet}_{\tw}(\gbige,0)$
is closely related with
$\gbige\otimes p_{\lambda}^{\ast}\Omega^{\bullet}$,
and for Theorem \ref{thm;22.4.2.40} with $\lambda=0$,
we would like to construct a sequence of
quasi-isomorphisms between
$\nbigc^{\bullet}_{\tw}(\gbige,0)$
and the $L^2$-complex.

\paragraph{Auxiliary complexes}

We introduce some more auxiliary complexes.
We set
$W_1\nbigc^{\bullet}_{\tw}(\gbige,0):=
\nbigc^{\bullet}(\{\nbigj^{\bullet}_{1,K}\})$:
\[
 \nbigj_{1,O}:
 W_1\nbigp_{0,0}E_{|O}
 \lrarr
 W_0(\Image f_{1|O})\,\frac{dz_1}{z_1}
 \oplus
 W_0(\Image f_{2|O})\,\frac{dz_2}{z_2}
 \lrarr
 W_{-1}\bigl(
 \Image (f_1f_2)_{|O}
 \bigr)
 \frac{dz_1}{z_1}\frac{dz_2}{z_2},
\]
\[
 \nbigj_{1,H_i}=\nbigj_{0,H_i}
 \quad (i=1,2).
\]
Here, $W_{\bullet}(U)$ denotes
the monodromy weight filtration of $f_{1|O}+f_{2|O}$
on $U$.

We set
$(W_1\cap W_0')\nbigc^{\bullet}_{\tw}(\gbige,0):=
\nbigc^{\bullet}(\{\nbigj_{2,K}\})$:
{\small
\[
 \nbigj_{2,O}:
 (W_1\cap W_0')\nbigp_{0,0}E_{|O}
 \lrarr
 (W_0\cap W_{-1}')(\Image f_{1|O})\,\frac{dz_1}{z_1}
 \oplus
 (W_{0}\cap W_0')(\Image f_{2|O})\,\frac{dz_2}{z_2}
 \lrarr
 (W_{-1}\cap W_{-1}')
 \bigl(
 \Image (f_1f_2)_{|O}
 \bigr)
 \frac{dz_1}{z_1}\frac{dz_2}{z_2},
\]
}
\[
\nbigj_{2,H_1}:
       W_0'\bigl(
       \nbigp_{0,0}E_{|H_1}
       \bigr)
       \lrarr
       W'_{-1}\bigl(
       \Image(f_{1|H_1})
       \bigr)\cdot
       \frac{dz_1}{z_1}
       \oplus
       W_0'\bigl(
       \nbigp_{0,0}E_{|H_1}
       \bigr)
       \frac{dz_2}{z_2}
       \lrarr
       W'_{-1}\bigl(
       \Image(f_{1|H_1})
       \bigr)\cdot
       \frac{dz_1}{z_1}
       \frac{dz_2}{z_2}.       
\]
\[
 \nbigj_{2,H_2}:=\nbigj_{0,H_2}.
\]
Here, $W'_{\bullet}(U)$ denotes the monodromy weight filtration of
$f_{1|K}$ $(K=O,H_1)$ on a vector bundle $U$ on $K$.

We set
$(W_1\cap W_0'')\nbigc^{\bullet}_{\tw}(\gbige,0):=
\nbigc^{\bullet}(\{\nbigj_{3,K}\})$:
{\small
\[
 \nbigj_{3,O}:
 (W_1\cap W_0'')\nbigp_{0,0}E_{|O}
 \lrarr
(W_{0}\cap W_0'')(\Image f_{1|O})\,\frac{dz_1}{z_1}
 \oplus
 (W_0\cap W_{-1}'')(\Image f_{2|O})\,\frac{dz_2}{z_2}
 \lrarr
 (W_{-1}\cap W_{-1}'')
 \bigl(
 \Image (f_1f_2)_{|O}
 \bigr)
 \frac{dz_1}{z_1}\frac{dz_2}{z_2},
\]
}
\[
 \nbigj_{3,H_1}:=\nbigj_{0,H_1},
\]
\[
\nbigj_{3,H_2}:
       W_0''\bigl(
       \nbigp_{0,0}E_{|H_2}
       \bigr)
       \lrarr
       W_0''\bigl(
       \nbigp_{0,0}E_{|H_2}
       \bigr)
       \frac{dz_1}{z_1}
       \oplus
       W''_{-1}\bigl(
       \Image(f_{2|H_2})
       \bigr)\cdot
       \frac{dz_2}{z_2}
       \lrarr
       W''_{-1}\bigl(
       \Image(f_{2|H_2})
       \bigr)\cdot
       \frac{dz_1}{z_1}
       \frac{dz_2}{z_2}.
\]
Here $W''_{\bullet}(U)$ denotes the monodromy weight filtration of
$f_{2|K}$ $(K=O,H_2)$
on a vector bundle on $K$.

We set
$W_1\nbigc^{\bullet}_{\tw}(\gbige[\ast H],0):=
\nbigc^{\bullet}(\{\nbigj_{4,K}\})$:
\[
 \nbigj_{4,O}:
        W_1\nbigp_{0,0}E_{|O}
       \lrarr
       W_{-1}\nbigp_{0,0}E_{|O}\frac{dz_1}{z_1}
       \oplus
       W_{-1}\nbigp_{0,0}E_{|O}\frac{dz_2}{z_2}
       \lrarr
       W_{-3}\nbigp_{0,0}E_{|O}
       \frac{dz_1}{z_1}\frac{dz_2}{z_2},
\]
\[
 \nbigj_{4,H_i}=\nbigc^{\bullet}_{\tw}(\gbige[\ast H],0)_{|H_i}
 \quad(i=1,2).
\]

We set
$(W_1\cap W_0')\nbigc^{\bullet}_{\tw}(\gbige[\ast H],0):=
\nbigc^{\bullet}(\{\nbigj_{5,K}\})$:
{\small
\[
 \nbigj_{5,O}:
       (W_1\cap W_0')\nbigp_{0,0}E_{|O}
       \lrarr
       (W_{-1}\cap W_{-2}')\nbigp_{0,0}E_{|O}\,\frac{dz_1}{z_1}
       \oplus
       (W_{-1}\cap W_{0}')\nbigp_{0,0}E_{|O}\,\frac{dz_2}{z_2}
       \lrarr
       (W_{-3}\cap W_{-2}')\nbigp_{0,0}E_{|O}\,
       \frac{dz_1}{z_1}\frac{dz_2}{z_2},
\]
}
\[
 \nbigj_{5,H_1}:
       W_0'(\nbigp_{0,0}E_{|H_1})
       \lrarr
       W_{-2}'(\nbigp_{0,0}E_{|H_1})\,\frac{dz_1}{z_1}
       \oplus
       W_{0}'(\nbigp_{0,0}E_{|H_1})\,\frac{dz_2}{z_2}
       \lrarr
       W_{-2}'(\nbigp_{0,0}E_{|H_1})\,
       \frac{dz_1}{z_1}\frac{dz_2}{z_2},
\]
\[
 \nbigj_{5,H_2}=
 \nbigc^{\bullet}_{\tw}(\gbige[\ast H],0)_{|H_2}.
\]

We set 
$(W_1\cap W_0'')\nbigc^{\bullet}_{\tw}(\gbige[\ast H],0):=
\nbigc^{\bullet}(\{\nbigj_{6,K}\})$:
{\small
\[
 \nbigj_{6,O}:
       (W_1\cap W_0'')\nbigp_{0,0}E_{|O}
       \lrarr
       (W_{-1}\cap W_{0}'')\nbigp_{0,0}E_{|O}\,\frac{dz_1}{z_1}
       \oplus
       (W_{-1}\cap W_{-2}'')\nbigp_{0,0}E_{|O}\,\frac{dz_2}{z_2}
       \lrarr
       (W_{-3}\cap W_{-2}'')\nbigp_{0,0}E_{|O}\,
       \frac{dz_1}{z_1}\frac{dz_2}{z_2},
\]
}
\[
 \nbigj_{6,H_1}=
 \nbigc^{\bullet}_{\tw}(\gbige[\ast H],0)_{|H_1},
\]
\[
 \nbigj_{6,H_2}:
       W_0''(\nbigp_{0,0}E_{|H_2})
       \lrarr
       W_{0}''(\nbigp_{0,0}E_{|H_2})\,\frac{dz_1}{z_1}
       \oplus
       W_{-2}''(\nbigp_{0,0}E_{|H_2})\,\frac{dz_2}{z_2}
       \lrarr
       W_{-2}'(\nbigp_{0,0}E_{|H_2})\,
       \frac{dz_1}{z_1}\frac{dz_2}{z_2}.
\]

By the construction,
$W_1\nbigc^{\bullet}_{\tw}(\gbige,0)$,
$(W_1\cap W_0')\nbigc^{\bullet}_{\tw}(\gbige,0)$
and
$(W_1\cap W_0'')\nbigc^{\bullet}_{\tw}(\gbige,0)$
are subcomplexes of 
$\nbigc^{\bullet}_{\tw}(\gbige,0)$.
Similarly,
$W_1\nbigc^{\bullet}_{\tw}(\gbige[\ast H],0)$,
$(W_1\cap W_0')\nbigc^{\bullet}_{\tw}(\gbige[\ast H],0)$
and
$(W_1\cap W_0'')\nbigc^{\bullet}_{\tw}(\gbige[\ast H],0)$
are subcomplexes of 
$\nbigc^{\bullet}_{\tw}(\gbige[\ast H],0)$.

\begin{rem}
This type of subcomplexes are introduced in
{\rm\cite{Kashiwara-Kawai-Hodge-holonomic}}.
\hfill\qed 
\end{rem}

\paragraph{Construction of quasi-isomorphisms}

Let $\Psi:X\to \real_{\geq 0}^2$
be the continuous map defined by
\[
 \Psi(z_1,z_2)=\bigl(
(-\log|z_1|)^{-1},
(-\log|z_2|)^{-1}
\bigr).
\]
Let $\widetilde{(\real_{\geq 0}^2)}\to \real_{\geq 0}^2$
denote the oriented real blow up
along $(0,0)$.
Let $\Xtilde$ denote the fiber product of
$X$ and $\widetilde{(\real_{\geq 0}^2)}$
over $\real_{\geq 0}^2$.
Let $p_1:\Xtilde\to X$ denote the projection.
Let $\Htilde_i$ denote the closure of
$p_1^{-1}(H_i\setminus \{O\})$.

Let $\nbigc^{\bullet}$ denote 
$\nbigc^{\bullet}_{\tw}(\gbige,0)$ or
$\nbigc^{\bullet}_{\tw}(\gbige[\ast H],0)$.
The stalk $p_1^{-1}(\nbigc^{\bullet})_P$
of $\nbigc^{\bullet}$ at $P\in \Xtilde$
is naturally identified with
the stalk $\nbigc^{\bullet}_{p_1(P)}$
of $\nbigc^{\bullet}$ at $p_1(P)$.
We define the subcomplex
$\nbigw\nbigc^{\bullet}
\subset
p_1^{-1}\nbigc^{\bullet}$
as follows:
\begin{itemize}
 \item
      $p_1^{-1}\nbigc^{\bullet}=\nbigc^{\bullet}$
      on $X\setminus H=\Xtilde\setminus p^{-1}(H)$.
 \item For $P\in \Htilde_1$,
       we set
       $\nbigw\nbigc^{\bullet}_P=(W_1\cap W_0')(\nbigc^{\bullet})_{p_1(P)}$.
 \item For $P\in \Htilde_2$,
       we set
       $\nbigw\nbigc^{\bullet}_P=
       (W_1\cap W_0'')(\nbigc^{\bullet})_{p_1(P)}$.
 \item For $P\in p_1^{-1}(O)\setminus(\Htilde_1\cup\Htilde_2)$,
       $\nbigw\nbigc^{\bullet}_P=(W_1\nbigc^{\bullet})_{O}$.
\end{itemize}
Then, the following morphisms are naturally defined.
\[
\begin{CD}
 p_1^{-1}\nbigc^{\bullet}_{\tw}(\gbige,0)
 @<{a_1}<<
 \nbigw\nbigc^{\bullet}_{\tw}(\gbige,0)
 @>{a_2}>>
 \nbigw\nbigc^{\bullet}_{\tw}(\gbige[\ast H],0).
\end{CD}
\]
Here, $a_1$ is a quasi-isomorphism,
but $a_2$ is not, in general.
The following induced morphisms
in $\ttD(\cnum_X)$ are isomorphisms:
\[
\begin{CD}
 \nbigc^{\bullet}_{\tw}(\gbige,0)
 \simeq
 Rp_{1\ast}p_1^{-1}
 \nbigc^{\bullet}_{\tw}(\gbige,0)
 @<{b_1}<<
 Rp_{1\ast}\nbigw \nbigc^{\bullet}_{\tw}(\gbige,0)
 @>{b_2}>>
 Rp_{1\ast} \nbigw\nbigc^{\bullet}_{\tw}(\gbige[\ast H],0).
\end{CD}
\]
Indeed, $b_1$ is an isomorphism because $a_1$ is a quasi-isomorphism.
We use the consideration in \S\ref{subsection;22.4.4.2}
to prove that $b_2$ is a quasi-isomorphism.
(See Theorem \ref{thm;22.2.17.103},
Theorem \ref{thm;22.2.17.111} and Corollary \ref{cor;22.2.18.10}.)

We can define the $L^2$-complex
$\nbigc^{\bullet}_{L^2,\Xtilde}(E,\theta,h)$ on $\Xtilde$
as in the cases of
$\nbigc^{\bullet}_{L^2}(\nbige,\DD,h)$
and
$\nbigc^{\bullet}_{L^2}(\nbigelambda,\DDlambda,h)$.
There exist natural quasi-isomorphisms:
\[
\nbigc^{\bullet}_{L^2}(E,\theta,h)
\simeq
 p_{1\ast}
 \nbigc^{\bullet}_{L^2,\Xtilde}(E,\theta,h)
\lrarr
 Rp_{1\ast}
\nbigc^{\bullet}_{L^2,\Xtilde}(E,\theta,h). 
\]

As explained in \S\ref{subsection;22.3.17.51},
there exists a natural inclusion
$\nbigw\nbigc^{\bullet}_{\tw}(\gbige[\ast H],0)
\lrarr
\nbigc^{\bullet}_{L^2,\Xtilde}(E,\theta,h)$.
This is not a quasi-isomorphism,
but the induced morphism 
\[
Rp_{1\ast}
\nbigw\nbigc^{\bullet}_{\tw}(\gbige[\ast H],0)
\lrarr
Rp_{1\ast}
\nbigc^{\bullet}_{L^2,\Xtilde}(E,\theta,h)
\]
is a quasi-isomorphism
(Theorem \ref{thm;22.2.18.11}).
As a result,
we obtain the following isomorphisms
\begin{multline}
\nbigc_{\tw}^{\bullet}(\gbige,0)
\lrarr
Rp_{1\ast}p_1^{-1}\nbigc_{\tw}^{\bullet}(\gbige,0)
\llarr
Rp_{1\ast}\nbigw\nbigc_{\tw}^{\bullet}(\gbige,0)
\lrarr
Rp_{1\ast}\nbigw\nbigc_{\tw}^{\bullet}(\gbige[\ast H],0)
\lrarr \\
Rp_{1\ast} \nbigc^{\bullet}_{L^2,\Xtilde}(E,\theta,h)
\llarr
p_{1\ast}\nbigc^{\bullet}_{L^2,\Xtilde}(E,\theta,h)
 \simeq
 \nbigc^{\bullet}_{L^2}(E,\theta,h).
\end{multline}
in $\ttD(\cnum_X)$.
This construction can be globalized.

\begin{rem}
In {\rm\cite{Kashiwara-Kawai-Hodge-holonomic, Kashiwara-Kawai-partition}},
it is proposed to construct quasi-isomorphisms
by using a kind of partition of unity,
which we do not use in this paper.
See {\rm\cite{Kashiwara-Kawai-Hodge-holonomic, Kashiwara-Kawai-partition}}
for more detail. 
\hfill\qed
\end{rem}

\paragraph{Acknowledgements}

This research has grown from my effort
to understand the last paragraph of the paper
\cite{Kashiwara-Kawai-Hodge-holonomic}
due to Masaki Kashiwara and Takahiro Kawai.
I thank Claude Sabbah and Carlos Simpson
for helpful discussions
about twistor $D$-modules and harmonic bundles
on many occasions.
Needless to say, their works are most fundamental in this study.
I thank Ron Donagi whose question about the $L^2$-complexes of
tame harmonic bundles is one of the motivations of this study.
I thank Maxim Kontsevich for a related stimulating question.
My interest in Hard Lefschetz Theorem
for pure twistor $\nbigd$-modules has been renewed
by recent works due to
Junchao Shentu and Chen Zhao \cite{Shentu-Zhao}
and Chuanhao Wei and Ruijie Yang \cite{Wei-Yang}
though this study is mathematically independent from theirs.
I thank J\"{u}rgen Jost for answering my question.
I am grateful to Morihiko Saito for remarks and questions.
I thank Takahiro Saito for discussions,
which are helpful to keep my interest in
twistor $\nbigd$-modules and Hodge modules.
I thank Yoshifumi Tsuchimoto and Akira Ishii
for their constant encouragements.

I am partially supported by
the Grant-in-Aid for Scientific Research (S) (No. 17H06127),
the Grant-in-Aid for Scientific Research (S) (No. 16H06335),
the Grant-in-Aid for Scientific Research (A) (No. 21H04429),
and the Grant-in-Aid for Scientific Research (C) (No. 20K03609),
Japan Society for the Promotion of Science.

\section{Preliminary}

\subsection{Some notation}

\subsubsection{Vector spaces associated with a finite subset}
\label{subsection;22.2.4.30}

Let $\Gamma$ be any non-empty finite set.
We set 
$\cnum\cdot\Gamma
=\bigoplus_{i\in\Gamma}\cnum\cdot\ttv_i$,
i.e.,
$\cnum\cdot\Gamma$ denotes the vector space
equipped with a base $\ttv_i$ $(i\in\Gamma)$.
We set
$\cnum\langle\Gamma\rangle=
\bigwedge^{\bullet}(\cnum\cdot\Gamma)$.

Let $2^{\Gamma}$ denote the set of subsets of $\Gamma$.
For any
$I=\{i_1,\ldots,i_k\}\in 2^{\Gamma}\setminus\{\emptyset\}$,
let $\cnum(I)$ denote the subspace
generated by
$\ttv_{i_1}\wedge\cdots\wedge\ttv_{i_k}$.
We set
$\cnum(\emptyset)=\bigwedge^0\cnum\langle\Gamma\rangle=\cnum$.
The exterior product of $\ttv_i$ induces
a morphism
$\ttv_i\wedge:
\cnum\langle\Gamma\rangle
\lrarr
\cnum\langle\Gamma\rangle$
such that
$\ttv_i\wedge\Bigl(
\bigwedge^k(\cnum\cdot\Gamma)
\Bigr)
\subset
\bigwedge^{k+1}(\cnum\cdot\Gamma)$.

Let
$(\cnum\cdot\Gamma)^{\lor}$
and 
$\cnum\langle\Gamma\rangle^{\lor}$
denote the dual spaces of
$\cnum\cdot\Gamma$ and
$\cnum\langle\Gamma\rangle$,
respectively.
Let $\ttv_i^{\lor}$ $(i\in\Gamma)$
denote the base of $(\cnum\cdot\Gamma)^{\lor}$,
which is dual to the base $\ttv_i$ $(i\in\Gamma)$ of
$\cnum\cdot\Gamma$.
For any $I=\{i_1,\ldots,i_k\}\in 2^{\Gamma}\setminus\{\emptyset\}$,
let $\cnum(I)^{\lor}\subset\cnum\langle\Gamma\rangle^{\lor}$
denote the subspace generated by
$\ttv_{i_1}^{\lor}\wedge\cdots\wedge \ttv_{i_k}^{\lor}$.
The exterior product of $\ttv_i^{\lor}$
induces
a morphism
$\ttv_i^{\lor}\wedge:\cnum\langle\Gamma\rangle
\lrarr\cnum\langle\Gamma\rangle$.

The inner product of $\ttv_i^{\lor}$
induces a morphism
$\ttv_i^{\lor}\vdash:\cnum\langle\Gamma\rangle
\lrarr\cnum\langle\Gamma\rangle$
for which
$\ttv_i^{\lor}\vdash\bigl(
\bigwedge^k(\cnum\cdot\Gamma)
\bigr)
\subset
\bigwedge^{k-1}(\cnum\cdot\Gamma)$.
Similarly,
the inner product of $\ttv_i$
induces a morphism
$\ttv_i\vdash:\cnum\langle\Gamma\rangle^{\lor}
\lrarr\cnum\langle\Gamma\rangle^{\lor}$
for which
$\ttv_i\vdash\bigl(
\bigwedge^k(\cnum\cdot\Gamma)^{\lor}
\bigr)
\subset
\bigwedge^{k-1}(\cnum\cdot\Gamma)^{\lor}$.

\subsubsection{Some categories of sheaves and complexes}

Let $S$ be a topological space.
Let $A_S$ be a sheaf of algebras on $S$.
Let $\Mod(A_S)$ denote the category of
$A_S$-modules.
Let $T$ be a topological space.
Let $p_T:T\times S\to S$  denote the projection.

\begin{lem}
If $T$ is contractible,
$p_T^{-1}:\Mod(A_S)\to \Mod(p_T^{-1}(A_S))$
is fully faithful.
Let
$\Mod(p_T^{-1}(A_S)/p_T)$
denote the essential image.
For $\ttF\in \Mod(p_T^{-1}(A_S)/p_T)$,
the natural morphism
$p_T^{-1}p_{T\ast}\ttF\lrarr\ttF$
is an isomorphism.
\end{lem}
\pf
According to \cite[Corollary 2.7.7 (i)]{Kashiwara-Schapira},
for any $\ttG\in\Mod(A_S)$,
the natural morphism
$\ttG\to Rp_{T\ast}\circ p_T^{-1}(\ttG)$
in the derived category of $A_S$-complexes
is an isomorphism.
It implies that the natural morphism
$\ttG\to p_{T\ast}\circ p_T^{-1}(\ttG)$
in $\Mod(A_S)$ is an isomorphism.
In general, the composite of the natural morphisms
\[
p_T^{-1}(\ttG)
\lrarr
p_T^{-1}\bigl(p_{T\ast}\circ p_T^{-1}(\ttG)\bigr)
=
\bigl(p_T^{-1}\circ p_{T\ast}\bigr)
p_T^{-1}(\ttG)
\lrarr p_T^{-1}(\ttG)
\]
is the identity.
Then, we can easily obtain the claim of the lemma.
\hfill\qed

\vspace{.1in}

Let $\ttC(A_S)$ denote the category of
bounded $A_S$-complexes.
The functor
$p_T^{-1}:\ttC(A_S)\to \ttC(p_T^{-1}(A_S))$
is fully faithful.
Let $\ttC(p_T^{-1}(A_S)/p_T)$
denote the essential image.
For any
$\ttF^{\bullet}\in\ttC(p_T^{-1}(A_S)/p_T)$,
the natural morphism
$p_T^{-1}p_{T\ast}\ttF^{\bullet}\to \ttF$
is an isomorphism.

\subsection{Toric manifolds and weakly constructible sheaves}

\subsubsection{Affine toric manifold associated with a cone}

We refer \cite{Toric-Varieties-Book, Fulton-toric-varieties}
as useful references for toric varieties.

Let $\ttN$ be a free $\seisuu$-module of finite rank.
We set $\ttN_{\real}=\real\otimes_{\seisuu}\ttN$.
Let $\ttM=\Hom_{\seisuu}(\ttN,\seisuu)$
and $\ttM_{\real}=\real\otimes \ttM$.
For each finitely generated cone $\sigma\subset \ttN_{\real}$,
we set
\[
 \sigma^{\lor}=\bigl\{
 m\in \ttM_{\real}\,\big|\,
 \langle m,n\rangle\geq 0,\,\,\forall n\in\sigma
 \bigr\},
 \quad
 \sigma^{\bot}=\bigl\{
 m\in \ttM_{\real}\,\big|\,
 \langle m,n\rangle=0,\,\,\forall n\in\sigma
 \bigr\}.
\]
We obtain the semigroup
$\ttS_{\sigma}=\sigma^{\lor}\cap \ttM$.
Let $\ttU_{\sigma}:=\Spec\cnum[\ttS_{\sigma}]$.
We identify $U_{\sigma}$ as the associated complex variety.
It is also identified with
the set of homomorphisms of semigroups
$\Hom_{\seisuu}(\ttS_{\sigma},\cnum)$,
where we regard $\cnum$ as a semigroup
by the multiplication.
We set
$\ttT_{\ttN}=\ttN\otimes_{\seisuu}\cnum^{\ast}
=\Hom_{\seisuu}(\ttM,\cnum^{\ast})=\Hom_{\seisuu}(\ttM,\cnum)$,
which naturally acts on $\ttU_{\sigma}$.
Recall that
$\ttU_{\sigma}$ has the distinguished point
$\gamma_{\sigma}\in\Hom_{\seisuu}(\ttS_{\sigma},\cnum)$
defined by
\[
 \gamma_{\sigma}(m)=
 \left\{
 \begin{array}{ll}
  1 & (m\in \sigma^{\bot}\cap \ttS_{\sigma})\\
  0 & (\mbox{otherwise}).
 \end{array}
 \right.
\]
We obtain the orbit
$\ttO(\sigma):=\ttT_{\ttN}\gamma_{\sigma}\subset \ttU_{\sigma}$.

We set
$(\ttU_{\sigma})_{\geq 0}:=
\Hom_{\seisuu}(\ttS_{\sigma},\real_{\geq 0})$,
where we regard $\real_{\geq 0}$ as a semigroup
by the multiplication.
Let $\ttT_{\ttN}^{\cp}\subset \ttT_{\ttN}$
denote the maximal compact subgroup.
The projection
$\ttU_{\sigma}\lrarr \ttU_{\sigma}/\ttT^{\cp}_{\ttN}$
induces a homeomorphism
$(\ttU_{\sigma})_{\geq 0}\simeq (\ttU_{\sigma})/\ttT_{\ttN}^{\cp}$.
We set
$\ttO(\sigma)_{\geq 0}:=(\ttU_{\sigma})_{\geq 0}\cap\ttO(\sigma)$.
We set
$\ttT_{\ttN,\real}:=\ttN\otimes\real_{>0}
=\Hom_{\seisuu}(\ttM,\real_{\geq 0})
=\Hom_{\seisuu}(\ttM,\real_{>0})$,
which naturally acts on $(\ttU_{\sigma})_{\geq 0}$.
We have
$\ttO(\sigma)_{\geq 0}
=\ttT_{\ttN,\real}\cdot\gamma_{\sigma}$.

We set $\ttd(\ttN):=\rank\ttN$.
Let $\sigma_0$ denote the cone
generated by a frame $\tte_1,\ldots,\tte_{\ttd(\ttN)}$ of $\ttN$.
Let $\tte_1^{\lor},\ldots,\tte_{\ttd(\ttN)}^{\lor}$
denote the dual frame of $\ttM$.
Because
$\ttS_{\sigma_0}=\bigoplus_{i=1}^{\ttd(\ttN)}
\seisuu_{\geq 0}\tte_i^{\lor}$,
we obtain the isomorphism
$\ttU_{\sigma_0}=\Hom_{\seisuu}(\ttS_{\sigma_0},\cnum)
\simeq \cnum^{\ttd(\ttN)}$
induced by
$f\longmapsto (f(\tte_1^{\lor}),\ldots,f(\tte_{\ttd(\ttN)}^{\lor}))$.
We may regard $(\tte_1^{\lor},\ldots,\tte_{\ttd(\ttN)}^{\lor})$
as a coordinate system.

Recall that $\ttU_{\sigma}$ is smooth
if and only if
there exists a frame $\tte_1,\ldots,\tte_{\ttd(\ttN)}$
of $\ttN$
such that $\sigma$ is generated by
$\tte_1,\ldots,\tte_k$ for some $0\leq k\leq \ttd(\ttN)$,
where the case $k=0$ means $\sigma=\{0\}$.
The coordinate system
$(\tte_1^{\lor},\ldots,\tte_{\ttd(\ttN)}^{\lor})$
induces the following isomorphisms:
\[
\ttU_{\sigma}\simeq
\cnum^{k}\times(\cnum^{\ast})^{\ttd(\ttN)-k},
\quad
(\ttU_{\sigma})_{\geq 0}\simeq
\real_{\geq 0}^{k}\times\real_{>0}^{\ttd(\ttN)-k},
\quad
\ttO(\sigma)\simeq
\{0\}\times
(\cnum^{\ast})^{\ttd(\ttN)-k}
\quad
\ttO(\sigma)_{\geq 0}\simeq
\{0\}\times
\real_{>0}^{\ttd(\ttN)-k}.
\]

\subsubsection{Toric manifolds and the non-negative parts}

Let $\Sigma$ be a fan in $\ttN_{\real}$,
and let $\ttX(\Sigma)$ denote the associated toric manifold,
which is obtained as the gluing of
$\ttU_{\sigma}$ $(\sigma\in\Sigma)$.
It is naturally equipped with the action of $\ttT_{\ttN}$,
and there exists the orbit decomposition:
\[
 \ttX(\Sigma)=\bigsqcup_{\sigma\in\Sigma}
 \ttO(\sigma). 
\]
We also obtain the non-negative part
$\ttX(\Sigma)_{\geq 0}$ of $\ttX(\Sigma)$.
It is equipped with the natural action of
$\ttT_{\ttN,\real}$,
and there exists the orbit decomposition:
\begin{equation}
\label{eq;21.11.8.1}
\ttX(\Sigma)_{\geq 0}=
 \bigsqcup_{\sigma\in\Sigma}
 \ttO(\sigma)_{\geq 0}.
\end{equation}
We have $\ttO(\sigma_1)_{\geq 0}
\subset\overline{\ttO(\sigma_2)}_{\geq 0}$
if and only if
$\sigma_1\supset\sigma_2$.
It is equivalent to
$\ttO(\sigma_1)
\subset\overline{\ttO(\sigma_2)}$.
According to \cite[Theorem 3.2.6]{Toric-Varieties-Book},
there exists the decomposition
\[
 \overline{\ttO(\sigma)}
=\bigsqcup_{\sigma\subset\sigma'}
 \ttO(\sigma'),
 \quad\quad
  \overline{\ttO(\sigma)}_{\geq 0}
=\bigsqcup_{\sigma\subset\sigma'}
 \ttO(\sigma')_{\geq 0}.
\]
For each $\sigma\in \Sigma$,
let $\iota_{\sigma}$ denote the inclusion
$\ttO(\sigma)_{\geq 0}\to\ttX(\Sigma)_{\geq 0}$.

Let $\Sigma(k)$ denote the set of the $k$-dimensional cones
in $\Sigma$.
For each $\sigma\in\Sigma$,
we set
$\Psi_{\Sigma}(\sigma):=
\bigl\{
 \tau\in\Sigma(1)\,\big|\,
 \tau\subset\sigma
\bigr\}$.
We obtain the injection
$\Psi_{\Sigma}:\Sigma\to 2^{\Sigma(1)}$.
For $\tau\in \Sigma(1)$,
let $\tte_{\tau}\in \ttN\cap \tau$ denote the generator.
Each $\sigma\in\Sigma$ is generated by
$\tte_{\tau}$ $(\tau\in\Psi_{\Sigma}(\sigma))$.

We say that $\Sigma$ is smooth,
if each $\sigma\in\Sigma$
is generated by a part of a basis of $\ttN$.
It is equivalent to the condition that
$\ttX(\Sigma)$ is smooth.
In that case,
$\ttX(\Sigma)_{\geq 0}$
is a closed real analytic submanifold with corner
$\ttX(\Sigma)$.
If $\Sigma$ is smooth,
we have $|\Psi_{\Sigma}(\sigma)|=\dim\sigma$
for any $\sigma\in\Sigma$.

\subsubsection{Weakly constructible sheaves
in the $\ttX(\Sigma)_{\geq 0}$-direction}
\label{subsection;22.2.4.10}

Let $S$ be a topological space.
Let $\Sigma$ be a smooth fan in $\ttN_{\real}$.
Let $\pi:\ttX(\Sigma)_{\geq 0}\times S\lrarr S$ denote the projection.
The inclusion
$\ttO(\sigma)_{\geq 0}\times S\to \ttX(\Sigma)_{\geq 0}\times S$
is denoted by $\iota_{\sigma,S}$.
The projection
$\ttO(\sigma)_{\geq 0}\times S\to S$ is denoted by
$\pi_{\sigma,S}$.

Let $A_S$ be a sheaf of algebras on $S$.
Let $\ttC^{\wc}(\ttX(\Sigma)_{\geq 0};A_S)
\subset\ttC(\pi^{-1}(A_S))$
denote the full subcategory of
$\ttF\in\ttC(\pi^{-1}(A_S))$
satisfying the following condition for any $\sigma\in\Sigma$:
\[
\iota_{\sigma,S}^{-1}(\ttF)
\in
\ttC(\pi_{\sigma,S}^{-1}(A_S)/\pi_{\sigma,S}).
\]
We set
$\Mod^{\wc}(\ttX(\Sigma)_{\geq 0};A_S)=
\ttC^{\wc}(\ttX(\Sigma)_{\geq 0};A_S)
\cap
\Mod(\ttX(\Sigma)_{\geq 0};A_S)$.
For any $\ttF\in\ttC^{\wc}(\ttX(\Sigma)_{\geq 0};A_S)$,
we obtain
\[
 \nbigf_{\ttF}(\sigma):=
 (\pi_{\sigma,S})_{\ast}\circ
 \iota_{\sigma,S}^{-1}(\ttF)
 \in \ttC(A_S).
\]

The inclusion
$\overline{\ttO(\sigma)_{\geq 0}}\times S\to
X(\Sigma)_{\geq 0}\times S$
is denoted by $\overline{\iota}_{\sigma,S}$,
and the projection
$\overline{\ttO(\sigma)_{\geq 0}}\times S \to S$
is denoted by $\overline{\pi}_{\sigma,S}$.
Note that there exist the following natural isomorphisms
in $\ttC^{\wc}(\ttX(\Sigma)_{\geq 0},A_S)$
for any $\sigma\in\Sigma$:
\[
 (\iota_{\sigma,S})_{\ast}\circ
\iota_{\sigma,S}^{-1}(\ttF)
 \simeq
 (\overline{\iota}_{\sigma,S})_{\ast}\circ
 \overline{\pi}_{\sigma,S}^{-1}(\nbigf_{\ttF}(\sigma)).
\]

Let $\sigma_1\subset\sigma_2$,
i.e.,
$\ttO(\sigma_1)_{\geq 0}\supset
 \ttO(\sigma_2)_{\geq 0}$.
The natural morphism
$\ttF\to
(\iota_{\sigma_1,S})_{\ast}\circ
\iota_{\sigma_1,S}^{-1}(\ttF)$
induces
\[
\pi_{\sigma_2,S}^{-1}(\nbigf_{\ttF}(\sigma_2))
\simeq
\iota_{\sigma_2,S}^{-1}(\ttF)
\to
\iota_{\sigma_2,S}^{-1}\circ
 (\iota_{\sigma_1,S})_{\ast}\circ
 \iota_{\sigma_1,S}^{-1}(\ttF)
 \simeq
\iota_{\sigma_2,S}^{-1}\circ
\overline{\iota}_{\sigma_1,S\ast}
 \overline{\pi}_{\sigma_1,S}^{-1}(\nbigf_{\ttF}(\sigma_1))
\simeq
\pi_{\sigma_2,S}^{-1}(\nbigf_{\ttF}(\sigma_1)).
\]
Hence, we obtain
$\nbigf_{\ttF}(\sigma_2)\to \nbigf_{\ttF}(\sigma_1)$.

We regard $\Sigma$ as a category,
where a morphism $\sigma_2\to\sigma_1$
corresponds to $\sigma_1\subset\sigma_2$.
Then, $\nbigf_{\ttF}$ is a functor from
$\Sigma$ to $\ttC(A_S)$.
Let $\Fun(\Sigma,\ttC(A_S))$ denote the category
of functors from $\Sigma$ to $\ttC(A_S)$.
Similarly,
let $\Fun(\Sigma,\Mod(A_S))$ denote the category
of functors from $\Sigma$ to $\Mod(A_S)$.
The following lemma is clear.

\begin{lem}
The above construction induces
\[
\ttC^{\wc}(\ttX(\Sigma)_{\geq 0};A_S)\simeq
 \Fun(\Sigma,\ttC(A_S)),
 \quad\quad
\Mod^{\wc}(\ttX(\Sigma)_{\geq 0};A_S)\simeq
\Fun(\Sigma,\Mod(A_S)).
\]
\hfill\qed
\end{lem}

\subsubsection{Push-forward}
\label{subsection;22.2.4.11}

Let $\nbigf\in\Fun(\Sigma,\Mod(A_S))$.
Let us construct
$\nbigc^{\bullet}(\nbigf)\in\ttC(A_S)$.
By using the injection $\Psi_{\Sigma}:\Sigma\to 2^{\Sigma(1)}$,
we set
$\cnum(\sigma)^{\lor}:=\cnum(\Psi_{\Sigma}(\sigma))^{\lor}
\subset
\cnum\langle 2^{\Sigma(1)}\rangle$,
where we use the notation in \S\ref{subsection;22.2.4.30}.
For $k\leq 0$,
we set
\[
 \nbigc^k(\nbigf):=
 \bigoplus_{\sigma\in\Sigma(-k)}
 \nbigf(\sigma)\otimes\cnum(\sigma)^{\lor}.
\]
For any $\sigma\in\Sigma$ and any $\tau\in \Psi_{\Sigma}(\sigma)$,
let $\sigma\setminus\tau$ denote the cone in $\Sigma$
generated by $\Psi_{\Sigma}(\sigma)\setminus\{\tau\}$.
The inner product with $\ttv_{\tau}$
and the natural morphism
$\nbigf(\sigma)\to\nbigf(\sigma\setminus\tau)$
induce
\[
 \nbigf(\sigma)\otimes\cnum(\sigma)^{\lor}
 \lrarr
 \nbigf(\sigma\setminus\tau)\otimes
 \cnum(\sigma\setminus\tau)^{\lor}.
\]
The morphisms for $\sigma\in\Sigma(-k)$ with
$\tau\in\Psi_{\Sigma}(\sigma)$
induce
$\nbigc^{k}(\nbigf)
\to \nbigc^{k+1}(\nbigf)$,
with which we obtain
$\nbigc^{\bullet}(\nbigf)\in\ttC(A_S)$.

For any $\nbigf^{\bullet}\in\Fun(\Sigma,\ttC(A_S))$,
we obtain the double complex
$\nbigc^{\bullet}(\nbigf^{\bullet})$
in $\Mod(A_S)$.
Let $\Tot\nbigc^{\bullet}(\nbigf^{\bullet})\in\ttC(A_S)$
denote the total complex.

\begin{prop}
\label{prop;22.2.4.3}
Choose an orientation of $\ttN_{\real}$.
We assume that $\ttX(\Sigma)$ is proper.
Then, there exists the following natural isomorphism
for any $\ttF^{\bullet}\in
\ttC^{\wc}(\ttX(\Sigma)_{\geq 0};A_S)$:
\[
 R\pi_{\ast}(\ttF^{\bullet})
 \simeq
 \Tot\nbigc^{\bullet}(\nbigf_{\ttF^{\bullet}})[-\ttd(\ttN)].
\]
\end{prop}
We shall construct an explicit isomorphism
in \S\ref{subsection;22.2.4.1} by using
the resolution in \S\ref{subsection;22.2.4.4}.

\subsubsection{Resolution of weakly constructible sheaves
in the $\ttX(\Sigma)_{\geq 0}$-direction}

\label{subsection;22.2.4.4}

In this subsection, we do not have to assume that
$\ttX(\Sigma)$ is proper.
Let $\ttF\in\Mod^{\wc}(\ttX(\Sigma)_{\geq 0};A_S)$.
For $k\leq 0$ and $\ell\geq 0$
such that $0\leq -k\leq \ell$,
we set
\[
C^{k,\ell}(\ttF)
:=\bigoplus_{\substack{
\sigma_1\in\Sigma(-k)\\
\sigma_2\in\Sigma(\ell)\\
\sigma_1\subset\sigma_2
}}
(\iota_{\sigma_2,S})_{\ast}\circ
\iota_{\sigma_2,S}^{-1}\circ
(\iota_{\sigma_1,S})_{\ast}\circ
\iota_{\sigma_1,S}^{-1}
(\ttF)
\otimes
\cnum(\sigma_2)^{\lor}
\otimes
\cnum(\sigma_1)^{\lor}.
\]
Note that there exists the following natural isomorphism:
\[
 (\iota_{\sigma_2,S})_{\ast}\circ
 \iota_{\sigma_2,S}^{-1}\circ
 (\iota_{\sigma_1,S})_{\ast}\circ
 \iota_{\sigma_1,S}^{-1}
(\ttF)
\simeq
 (\iota_{\sigma_2,S})_{\ast}\circ
 \pi_{\sigma_2,S}^{-1}\bigl(
 \nbigf_{\ttF}(\sigma_1)
 \bigr)
 \simeq
 (\overline{\iota}_{\sigma_2,S})_{\ast}\circ
 \overline{\pi}_{\sigma_2,S}^{-1}\bigl(
 \nbigf_{\ttF}(\sigma_1)
 \bigr).
\]

Let $\sigma_1,\sigma_2\in\Sigma$
such that $\sigma_1\subset\sigma_2$.
Let $\tau\in \Psi_{\Sigma}(\sigma_1)$.
The natural morphism
$\ttF\to
(\iota_{\sigma_1\setminus\tau,S})_{\ast}\circ
 \iota_{\sigma_1\setminus\tau,S}^{-1}(\ttF)$
induces
\begin{multline}
 (\iota_{\sigma_2,S})_{\ast}\circ
 \iota_{\sigma_2,S}^{-1}\circ
 (\iota_{\sigma_1,S})_{\ast}\circ
 \iota_{\sigma_1,S}^{-1}(\ttF)
\lrarr
 (\iota_{\sigma_2,S})_{\ast}\circ
 \iota_{\sigma_2,S}^{-1}\circ
 (\iota_{\sigma_1,S})_{\ast}\circ
 \iota_{\sigma_1,S}^{-1}\circ
 (\iota_{\sigma_1\setminus\tau,S})_{\ast}\circ
 \iota_{\sigma_1\setminus\tau,S}^{-1}(\ttF)
 \simeq
 \\
 (\iota_{\sigma_2,S})_{\ast}\circ
 \pi_{\sigma_2,S}^{-1}\bigl(
 \nbigf_{\ttF}(\sigma_1\setminus\tau)
 \bigr)
\simeq
 (\iota_{\sigma_2,S})_{\ast}\circ\iota_{\sigma_2,S}^{-1}
 \circ
 (\iota_{\sigma_1\setminus\tau,S})_{\ast}\circ
 \iota_{\sigma_1\setminus\tau,S}^{-1}(\ttF).
\end{multline}
Together with the inner product of $\ttv_{\tau}$,
we obtain
\[
 (\iota_{\sigma_2,S})_{\ast}\circ\iota_{\sigma_2,S}^{-1}\circ
 (\iota_{\sigma_1,S})_{\ast}\circ
 \iota_{\sigma_1,S}^{-1}(\ttF)
 \otimes\cnum(\sigma_2)^{\lor}
 \otimes\cnum(\sigma_1)^{\lor}
\lrarr
(\iota_{\sigma_2,S})_{\ast}\circ\iota_{\sigma_2,S}^{-1}\circ
(\iota_{\sigma_1\setminus\tau,S})_{\ast}\circ
\iota_{\sigma_1\setminus\tau,S}^{-1}(\ttF)
  \otimes\cnum(\sigma_2)^{\lor}
 \otimes\cnum(\sigma_1\setminus\tau)^{\lor}.
\]
They induce a morphism
$d^{(1)}:C^{k,\ell}(\ttF)\lrarr C^{k+1,\ell}(\ttF)$
such that  $d^{(1)}\circ d^{(1)}=0$.

For two cones $\tau_i$ $(i=1,2)$ in $\ttN_{\real}$,
let $\tau_1\tau_2$ denote the cone generated by
$\tau_1$ and $\tau_2$.
For $\tau\in\Sigma(1)$
and $\sigma_2\in\Sigma(\ell)$ such that
$\sigma_2\tau\in\Sigma(\ell+1)$,
and for $\sigma_1$ such that $\sigma_2\supset\sigma_1$,
there exists the following natural morphism:
\begin{multline}
(\iota_{\sigma_2,S})_{\ast}\circ
\iota_{\sigma_2,S}^{-1}\circ
(\iota_{\sigma_1,S})_{\ast}\circ
\iota_{\sigma_1,S}^{-1}
(\ttF)
 \simeq
 (\overline{\iota}_{\sigma_2,S})_{\ast}\circ
 \overline{\pi}_{\sigma_2,S}^{-1}\bigl(
 \nbigf_{\ttF}(\sigma_1)
 \bigr)
 \lrarr \\
 (\iota_{\sigma_2\tau,S})_{\ast}\circ
 \pi_{\sigma_2\tau,S}^{-1}\bigl(
 \nbigf_{\ttF}(\sigma_1)
 \bigr)
 \simeq
 (\iota_{\sigma_2\tau,S})_{\ast}\circ
 \iota_{\sigma_2\tau}^{-1}\circ
 (\iota_{\sigma_1,S})_{\ast}\circ\iota_{\sigma_1,S}^{-1}
(\ttF).
\end{multline}
Together with the exterior product of $\ttv_{\tau}^{\lor}$,
we obtain the following morphism:
\[
 (\iota_{\sigma_2,S})_{\ast}\circ
 \iota_{\sigma_2,S}^{-1}\circ
 (\iota_{\sigma_1,S})_{\ast}\circ
 \iota_{\sigma_1,S}^{-1}(\ttF)
 \otimes
 \cnum(\sigma_2)^{\lor}
 \otimes
 \cnum(\sigma_1)^{\lor}
\lrarr
 (\iota_{\sigma_2\tau,S})_{\ast}\circ
 \iota_{\sigma_2\tau,S}^{-1}\circ
 (\iota_{\sigma_1,S})_{\ast}\circ
 \iota_{\sigma_1,S}^{-1}(\ttF)
 \otimes
  \cnum(\sigma_2\tau)^{\lor}
 \otimes
 \cnum(\sigma_1)^{\lor}.
\]
They induce
$d^{(2)}:C^{k,\ell}(\ttF)\lrarr C^{k,\ell+1}(\ttF)$
such that $d^{(2)}\circ d^{(2)}=0$.

By the construction,
we have
$d^{(1)}\circ d^{(2)}=d^{(2)}\circ d^{(1)}$.
We obtain
$\Tot C^{\bullet,\bullet}(\ttF)\in
\ttC^{\wc}(\ttX(\Sigma)_{\geq 0};A_S)$
as the total complex.

For any $\sigma\in\Sigma$,
there exists the canonical morphism
\[
 \ttF\lrarr
  (\iota_{\sigma,S})_{\ast}\circ
  \iota_{\sigma,S}^{-1}\circ
  (\iota_{\sigma,S})_{\ast}\circ
  \iota_{\sigma,S}^{-1}(\ttF)
 \otimes
 \cnum(\sigma)^{\lor}
 \otimes
 \cnum(\sigma)^{\lor}.
\]
They induce the following morphism of complexes
\begin{equation}
\label{eq;21.11.8.2}
 \ttF\lrarr
 \Tot C^{\bullet,\bullet}(\ttF).
\end{equation}

\begin{lem}
\label{lem;22.2.4.2}
The morphism {\rm(\ref{eq;21.11.8.2})}
is a quasi-isomorphism.
\end{lem}
\pf
Let $\sigma_1\in\Sigma$.
For $m\geq 0$,
we set
\[
 C_{\sigma_1}^{m}(\ttF)=
 \bigoplus_{\substack{
 \sigma_2\in\Sigma\\
 \sigma_1\subset\sigma_2\\
 \dim\sigma_2-\dim\sigma_1=m
 }}
 (\iota_{\sigma_2,S})_{\ast}\circ
 \iota_{\sigma_2,S}^{-1}\circ
 (\iota_{\sigma_1,S})_{\ast}\circ
 \iota_{\sigma_1,S}^{-1}(\ttF)
 \otimes
 \cnum(\sigma_2)^{\lor}.
\]
The morphism
$C_{\sigma_1}^m(\ttF)\to C_{\sigma_1}^{m+1}(\ttF)$
is defined as in the case of $d^{(2)}$,
and thus we obtain a complex
$C^{\bullet}_{\sigma_1}(\ttF)$.
The morphism
$(\iota_{\sigma_1,S})_{!}\circ
\iota_{\sigma_1,S}^{-1}(\ttF)
\lrarr
(\iota_{\sigma_1,S})_{\ast}\circ
\iota_{\sigma_1,S}^{-1}(\ttF)$
naturally induces the following morphism of complexes:
\begin{equation}
\label{eq;21.11.8.3}
 (\iota_{\sigma_1,S})_{!}\circ
 \iota_{\sigma_1,S}^{-1}(\ttF)
 \lrarr
 C^{\bullet}_{\sigma_1}(\ttF)\otimes\cnum(\sigma_1)^{\lor}.
\end{equation}
We can check that 
(\ref{eq;21.11.8.3}) is a quasi-isomorphism.

At each point of $\ttX(\Sigma)_{\geq 0}\times S$,
we obtain that
$\Tot C^{\bullet,\bullet}(\ttF)$
is quasi-isomorphic to $\ttF$,
by considering the cohomology with respect to $d^{(2)}$.
Then, it is easy to see that
(\ref{eq;21.11.8.2}) induces an isomorphism.
\hfill\qed

\vspace{.1in}

We note that
\[
 R\pi_{\ast}\circ
 (\iota_{\sigma_2,S})_{\ast}\circ
 \iota_{\sigma_2,S}^{-1}\circ
 (\iota_{\sigma_1,S})_{\ast}\circ
 \iota_{\sigma_1,S}^{-1}(\ttF)
\simeq
\pi_{\ast}\circ
(\iota_{\sigma_2,S})_{\ast}\circ
\iota_{\sigma_2,S}^{-1}\circ
(\iota_{\sigma_1,S})_{\ast}\circ
\iota_{\sigma_1,S}^{-1}(\ttF)
\simeq
\nbigf_{\ttF}(\sigma_1).
\]
By Lemma \ref{lem;22.2.4.2},
$R\pi_{\ast}(\ttF)$ is represented by
\begin{equation}
\label{eq;22.2.4.6}
 R\pi_{\ast}
 \Tot C^{\bullet,\bullet}(\ttF)
 \simeq
 \pi_{\ast}
  \Tot C^{\bullet,\bullet}(\ttF).
\end{equation}

\subsubsection{Proof of Proposition \ref{prop;22.2.4.3}}
\label{subsection;22.2.4.1}

We assume that $\ttX(\Sigma)$ is proper.
Let us look at the following complex
\begin{equation}
\label{eq;22.2.4.5}
 R\pi_{\ast}C^{\bullet}_{\sigma_1}(\ttF)
 \simeq
 \pi_{\ast}C^{\bullet}_{\sigma_1}(\ttF).
\end{equation}
Because (\ref{eq;21.11.8.3}) is a quasi-isomorphism,
the $j$-th cohomology sheaf of the complex (\ref{eq;22.2.4.5})
is $0$ unless
$j=\dim \ttO(\sigma)_{\geq 0}=\ttd(\ttN)-\dim\sigma$.
The $(\ttd(\ttN)-\dim\sigma)$-th cohomology sheaf
is the cokernel $\nbigc(\ttF,\sigma)$
of the following morphism:
\[
\bigoplus_{\substack{\sigma'_2\in\Sigma(\ttd(\ttN)-1)\\
\sigma_1\subset\sigma'_2
}}
 \nbigf_{\ttF}(\sigma_1)\otimes
 \cnum(\sigma'_2)^{\lor}
\lrarr
\bigoplus_{\substack{\sigma_2\in\Sigma(\ttd(\ttN))\\
\sigma_1\subset\sigma_2
}}
 \nbigf_{\ttF}(\sigma_1)\otimes
 \cnum(\sigma_2)^{\lor}.
\]
For $\sigma_1\supset\sigma_1'$,
we have the naturally defined morphism
$\nbigc(\ttF,\sigma_1)\to\nbigc(\ttF,\sigma_1')$.
For $k\leq 0$,
we set
\[
 \check{\nbigc}^k(\ttF):=
 \bigoplus_{\sigma_1\in\Sigma(-k)}
 \nbigc(\ttF,\sigma_1)
 \otimes\cnum(\sigma_1)^{\lor}.
\]
For $\sigma_1\in\Sigma(-k)$
and $\tau\in\Psi_{\Sigma}(\sigma)$,
we obtain the following morphism
by the inner product of $\ttv_{\tau}$
and the natural morphism
$\nbigc(\ttF,\sigma)\lrarr
\nbigc(\ttF,\sigma\setminus\tau)$:
\[
 \nbigc(\ttF,\sigma)\otimes\cnum(\sigma)^{\lor}
 \lrarr
 \nbigc(\ttF,\sigma\setminus\tau)
 \otimes\cnum(\sigma\setminus\tau)^{\lor}.
\]
The morphisms induce
$\check{\nbigc}^{k}(\ttF)\lrarr\check{\nbigc}^{k+1}(\ttF)$
with which
$\check{\nbigc}^{\bullet}(\ttF)$ is a complex.
By the construction,
there exists the natural quasi-isomorphism:
\[
 \pi_{\ast}
 \Tot C^{\bullet,\bullet}(\ttF)
 \simeq
 \check{\nbigc}^{\bullet}(\ttF)[-\ttd(\ttN)].
\]

An orientation of $\ttN$ induces the isomorphisms
$\cnum(\sigma_2)^{\lor}\simeq\cnum$
for any $\sigma_2\in\Sigma(\ttd(\ttN))$.
They induce the isomorphism
$\nbigc(\ttF,\sigma_1)\simeq \nbigf_{\ttF}(\sigma_1)$
for any $\sigma_1$.
We obtain
\[
 \check{\nbigc}^{\bullet}(\ttF)
 \simeq
 \nbigc^{\bullet}(\ttF).
\]
Thus, we obtain the claim of
Proposition \ref{prop;22.2.4.3}
for $\ttF\in \Mod^{\wc}(\ttX(\Sigma)_{\geq 0}\times S)$.
By the construction, it is easily generalized
to the case $\ttF^{\bullet}\in
\ttC^{\wc}(\ttX(\Sigma)_{\geq 0}\times S)$.
\hfill\qed

\subsection{Some toric manifolds}
\label{subsection;21.12.8.1}

\subsubsection{Increasing sequences}
\label{subsection;22.4.26.1}

Let $\Lambda$ be a non-empty finite set.
We set $\ttN(\Lambda):=\bigoplus_{i\in\Lambda}\seisuu \tte_i$.
The dual $\ttM(\Lambda)$ is equipped with the dual frame
$(\tte_i^{\lor})_{i\in\Lambda}$.
For a positive integer $m$,
an increasing sequence $\vecJ$ of length $m$
in $\Lambda$
is a subset
$\vecJ=(J_1,\ldots,J_m)\subset
2^{\Lambda}\setminus\{\emptyset\}$
such that
\[
J_1\subsetneq J_2\subsetneq\cdots
\subsetneq J_m\subset\Lambda.
\]
An increasing sequence of length $0$
means the empty subset of $2^{\Lambda}\setminus\{\emptyset\}$.
Let $\nbigs(\Lambda)$
denote the set of increasing sequences in $\Lambda$.
The length of $\vecJ\in\nbigs(\Lambda)$
is denoted by $|\vecJ|$.

For any non-empty subset $J\subset\Lambda$,
let $\tte_J:=\sum_{i\in J}\tte_i\in \ttN(\Lambda)$.
For each non-empty $\vecJ\in\nbigs(\Lambda)$,
let $\sigma(\vecJ)$ denote the cone in $\ttN(\Lambda)_{\real}$
generated by $\tte_J$ $(J\in\vecJ)$.
If $\vecJ$ is empty, we set $\sigma(\vecJ):=\{0\}$.

Let $\Sigma(\Lambda)$ denote the fan of
the cones $\sigma(\vecJ)$ $(\vecJ\in\nbigs(\Lambda))$.
We obtain the complex toric variety
$\ttX(\Lambda):=\ttX(\Sigma(\Lambda))$
and the non-negative part
$\ttX(\Lambda)_{\geq 0}:=\ttX(\Sigma(\Lambda))_{\geq 0}$.
For each $\vecJ\in\nbigs(\Lambda)$,
we set
$\ttU_{\vecJ}:=\ttU_{\sigma(\vecJ)}$,
$\ttO(\vecJ)=\ttO(\sigma(\vecJ))$
and
$(\ttU_{\vecJ})_{\geq 0}:=
(\ttU_{\sigma(\vecJ)})_{\geq 0}$.

Let $\sigma_0$ denote the cone in $\ttN(\Lambda)_{\real}$
generated by
$\tte_i$ $(i\in\Lambda)$.
Because $\sigma(\vecJ)\subset\sigma_0$ for any $\vecJ\in\nbigs(\Lambda)$,
there exists the natural toric morphism
$\pi_{\Lambda}:\ttX(\Sigma)\lrarr
\ttU_{\sigma_0}=\cnum^{\Lambda}$.
Let $(z_i)_{i\in\Lambda}$ be the coordinate system
of $\cnum^{\Lambda}$
induced by $\tte_i^{\lor}$ $(i\in\Lambda)$.
Each $J\subset\Lambda$
induces an increasing sequence of length one
$(J)\in\nbigs(\Lambda)$.
We have
\[
\pi_{\Lambda}^{-1}\left(
\bigcap_{i\in J}\{z_i=0\}
\right)
=\overline{\ttO(J)}.
\]
In particular,
$\pi_{\Lambda}^{-1}(0,\ldots,0)
=\overline{\ttO(\Lambda)}$.

\subsubsection{Orbits}
\label{subsection;22.4.26.30}

For $\vecJ,\vecJ'\in\nbigs(\Lambda)$,
we say that $\vecJ'$ is a refinement of $\vecJ$
if $\vecJ\subset\vecJ'$ as a subset of
$2^{\Lambda}\setminus\{\emptyset\}$.
It is equivalent to
$\sigma(\vecJ)\subset\sigma(\vecJ')$.
Let $\nbigs(\Lambda,\vecJ)$ denote the set of
$\vecJ'\in\nbigs(\Lambda)$
which is a refinement of $\vecJ$.
We have the following orbit decompositions
\begin{equation}
\label{eq;22.2.4.8}
 \overline{\ttO(\vecJ)}
 =\bigsqcup_{\vecJ'\in\nbigs(\Lambda,\vecJ)}
 \ttO(\vecJ'),
 \quad\quad
 \overline{\ttO(\vecJ)}_{\geq 0}
 =\bigsqcup_{\vecJ'\in\nbigs(\Lambda,\vecJ)}
 \ttO(\vecJ')_{\geq 0}.
\end{equation}
By using the decompositions (\ref{eq;22.2.4.8}),
we obtain the following lemma.
\begin{lem}
If $\vecJ=\bigsqcup \vecJ_i\in \nbigs(\Lambda)$,
then we have
$\overline{\ttO(\vecJ)}
=\bigcap \overline{\ttO(\vecJ_i)}$
and
$\overline{\ttO(\vecJ)}_{\geq 0}
=\bigcap \overline{\ttO(\vecJ_i)}_{\geq 0}$.
\hfill\qed
\end{lem}

We note that
$\dim_{\real}\ttO(\vecJ)_{\geq 0}
=|\Lambda|-|\vecJ|$.

\subsubsection{An affine chart}
\label{subsection;22.2.18.3}

We set $\isitabar=\{1,\ldots,i\}$ for any $i\in\seisuu_{\geq 1}$.
We choose a bijection $\Lambda\simeq \ellsitabar$
by setting $\ell=|\Lambda|$.
We set
$\vecJ_{\st}=(\isitabar\,|\,i=1,\ldots,\ell)\in\nbigs(\ellsitabar)$.
Then,
\[
 \sigma(\vecJ_{\st})^{\lor}
 =\left\{
 \sum_{i=1}^{\ell} u_i\tte_i^{\lor}\,\left|\,
 \sum_{j=1}^m u_j\geq 0\,\,
 (m=1,\ldots,\ell)
 \right.
 \right\}
 =\bigoplus_{i=1}^{\ell-1} \seisuu_{\geq 0}
 (\tte_i^{\lor}-\tte_{i+1}^{\lor})
\oplus
 \seisuu_{\geq 0}\tte_{\ell}^{\lor}.
\]
By setting
$t^{\st}_i=z_i/z_{i+1}$ $(i=1,\ldots,\ell-1)$
and $t^{\st}_{\ell}=z_{\ell}$,
we have the natural identification
\[
\ttU_{\vecJ_{\st}}=\{(t^{\st}_1,\ldots,t^{\st}_{\ell})\}
=\cnum^{\ell},
\]
and the restriction of $\pi_{\Lambda}$ to $\ttU_{\vecJ_{\st}}$
is described by
\[
\pi_{\Lambda}^{\ast}(z_{i})=\prod_{i\leq j\leq \ell}t^{\st}_j
\quad
(i=1,\ldots,\ell).
\]
We have
$(\ttU_{\vecJ_{\st}})_{\geq 0}
 =\bigl\{
 (t^{\st}_1,\ldots,t^{\st}_{\ell})\,\big|\,
 t^{\st}_{j}\geq 0\,\,\,(j=1,\ldots,\ell)
 \bigr\}$
and
$\ttO(\vecJ_{\st})=\{(0,\ldots,0)\}$
with respect to
the coordinate system
$(t^{\st}_1,\ldots,t^{\st}_{\ell})$.
For $I\subset\ellsitabar$,
we set $\vecJ_{\st,I}:=(\ibar\,|\,i\in I)\subset \vecJ_{\st}$,
and then
\[
 \ttU_{\vecJ_{\st,I}}
 =\bigl\{
 (t^{\st}_1,\ldots,t^{\st}_{\ell})\in \ttU_{\vecJ_{\st}}
 \,\big|\,
 t^{\st}_i\neq 0\,\,(i\not\in I)
 \bigr\}.
\]
We have
$(\ttU_{\vecJ_{\st,I}})_{\geq 0}
=\prod_{i\in I}\{t^{\st}_i\geq 0\}
 \times
 \prod_{i\not\in I}\{t^{\st}_i>0\}$,
$\ttO(\vecJ_{\st,I})=
\prod_{i\in I}\{t^{\st}_i=0\}
\times
\prod_{i\not\in I}\{t^{\st}_i\neq 0\}$,
and 
$\ttO(\vecJ_{\st,I})_{\geq 0}=
\prod_{i\in I}\{t^{\st}_i=0\}
\times
\prod_{i\not\in I}\{t^{\st}_i>0\}$.
We have
\[
 \overline{\ttO(\isitabar)}\cap \ttU_{\vecJ_{\st}}
=\prod_{j\neq i}
\{t^{\st}_j\in\cnum\}
\times
\bigl\{
 t^{\st}_{i}=0
\bigr\}\simeq \cnum^{\ell-1}.
\]
In particular,
$\overline{\ttO(\ellsitabar)}\cap \ttU_{\vecJ_{\st}}
=\prod_{j=1}^{\ell-1}\{t^{\st}_j\in\cnum\}
\times
\bigl\{
 t^{\st}_{\ell}=0
\bigr\}$.
We have
$\ttO(\vecJ_{\st,I})\subset
\overline{\ttO(\ellsitabar)}$
if and only if $\ell\in I$.
We have
\[
(\ttU_{\vecJ_{\st}})_{\geq 0}
=\coprod_{I\subset\ellsitabar}
\ttO(\vecJ_{\st,I})_{\geq 0},
\quad
\overline{\ttO(\ellsitabar)}\cap
(\ttU_{\vecJ_{\st}})_{\geq 0}
=\coprod_{\substack{
\ell\in I\subset\ellsitabar}}
\ttO(\vecJ_{\st,I})_{\geq 0}.
\]
If $J\neq \isitabar$ for any $i$,
then
$\ttO(J)\cap \ttU_{\vecJ_{\st}}=\emptyset$.

\subsubsection{The closure of $\ttO(\Lambda)$}
\label{subsection;22.2.7.20}

According to \cite[Proposition 3.2.7]{Toric-Varieties-Book},
the closure of any orbit is also a toric variety.
We recall the fan corresponding to
$\overline{\ttO(\Lambda)}$.
Let $\overline{\ttN}(\Lambda)=\ttN(\Lambda)/\seisuu \tte_{\Lambda}$
and $\overline{\ttN}(\Lambda)_{\real}
=\overline{\ttN}(\Lambda)\otimes\real
\simeq \ttN(\Lambda)_{\real}/\real e_{\Lambda}$.
Let $\bar{\tte}_{J}\in \ttNbar(\Lambda)$
denote the image of $\tte_J$.
Note that $\bar{\tte}_{\Lambda}=0$ in $\ttNbar(\Lambda)$.

We set $\nbigsbar(\Lambda)=\nbigs(\Lambda,\Lambda)$,
i.e.,
$\nbigsbar(\Lambda)\subset\nbigs(\Lambda)$
denotes the set of
increasing sequences $\vecJ=(J_i)$ in $\Lambda$
which contains $\Lambda$,
i.e.,
\[
 \emptyset\neq
 J_1\subsetneq J_2\subsetneq\cdots\subsetneq
 J_{m}=\Lambda.
\]
For each $\vecJ=(J_i)\in\nbigsbar(\Lambda)$,
let $\bar{\sigma}(\vecJ)$
denote the cone in $\ttNbar(\Lambda)_{\real}$
generated by $\bar{\tte}_{J_i}$.
It is equal to the image of $\sigma(\vecJ)$
by $\ttN(\Lambda)_{\real}\lrarr \ttNbar(\Lambda)_{\real}$.
Thus, we obtain a fan
$\Sigmabar(\Lambda)=
\{\bar{\sigma}(\vecJ)\,|\,\vecJ\in\nbigsbar(\Lambda)\}$.
Let $\ttXbar(\Lambda)$ denote the associated toric manifold,
which is projective
and isomorphic to $\overline{\ttO(\Lambda)}$.
We have
$\overline{\ttO(\Lambda)}_{\geq 0}
=\ttXbar(\Lambda)_{\geq 0}$.

Because the $\ttT_{\ttNbar(\Lambda)}$-action
on $\overline{\ttO(\Lambda)}$
is induced by 
the $\ttT_{\ttN(\Lambda)}$-action
on $\overline{\ttO(\Lambda)}$,
we have the same orbit decompositions:
\[
 \overline{\ttO(\Lambda)}
=\coprod_{\vecJ\in\nbigsbar(\Lambda)}
 \ttO(\vecJ),
 \quad\quad
 \overline{\ttO(\Lambda)}_{\geq 0}
=\coprod_{\vecJ\in\nbigsbar(\Lambda)}
 \ttO(\vecJ)_{\geq 0}.
\]

\subsubsection{The closure of $\ttO(J)$}

Let us also look at $\overline{\ttO(J)}$.
We set
$\ttNbar_J(\Lambda):=\ttN(\Lambda)/\seisuu \tte_{J}$
and
$\ttNbar_{J}(\Lambda)_{\real}:=
\ttNbar_{J}(\Lambda)\otimes\real$.
For $\tte_{J'}$,
let $\bar{\tte}_{J,J'}\in\ttNbar_{J}(\Lambda)$
denote the image of $\tte_{J'}$.
We set
$\nbigsbar_{J}(\Lambda)=\nbigs(\Lambda,J)$.
For each $\vecJ\in\nbigsbar_{J}(\Lambda)$,
let $\sigmabar_{J}(\vecJ)$
denote the image of $\sigma(\vecJ)$.
Thus, we obtain a fan
$\Sigmabar_{J}(\Lambda)$ in $\ttNbar_{J}(\Lambda)_{\real}$.
Let $\ttXbar_J(\Lambda)$ denote the associated toric manifold.
We have
$\overline{\ttO(J)}
=\ttXbar_J(\Lambda)$.

The $\ttT_{\ttNbar_J(\Lambda)}$-action
on $\ttXbar_J(\Lambda)$
is induced by the $\ttT_{\ttN(\Lambda)}$-action
on $\overline{\ttO(J)}$.
Hence, the orbit decompositions are the same:
\[
 \overline{\ttO(J)}
 =\coprod_{\vecJ\in\nbigsbar_J(\Lambda)}
 \ttO(\vecJ),
 \quad\quad
 \overline{\ttO(J)}_{\geq 0}
 =\coprod_{\vecJ\in\nbigsbar_J(\Lambda)}
 \ttO(\vecJ)_{\geq 0}.
\]

We set $I=\Lambda\setminus J$.
There exists a natural isomorphism
$\ttNbar_J(\Lambda)\simeq \ttNbar(J)\times \ttN(I)$.
There exists the natural bijection
$\nbigsbar_J(\Lambda)\simeq\nbigsbar(J)\times\nbigs(I)$,
and there exists the isomorphism of fans
$\Sigmabar_J(\Lambda)\simeq
\Sigmabar(J)\times\Sigma(I)$.
Hence, we obtain
\[
\overline{\ttO(J)}=\ttXbar_J(\Lambda)\simeq
 \ttXbar(J)\times \ttX(I),
 \quad\quad
\overline{\ttO(J)}_{\geq 0}
=\ttXbar_J(\Lambda)_{\geq 0}
 \simeq
 \ttXbar(J)_{\geq 0}\times \ttX(I)_{\geq 0}.
\]

\subsubsection{The closure of the orbits}

More generally,
let us look at $\overline{\ttO(\vecJ)}$
for $\vecJ\in\nbigs(\Lambda)$.
We set
$\ttNbar_{\vecJ}(\Lambda)=
\ttN(\Lambda)\big/\bigoplus_{J\in \vecJ}\seisuu \tte_{J}$
and 
$\ttNbar_{\vecJ}(\Lambda)_{\real}
=\ttNbar_{\vecJ}(\Lambda)\otimes\real$.
For $\tte_{J'}$,
let $\tte_{\vecJ,J'}\in \ttNbar_{\vecJ}(\Lambda)$
denote the image of $\tte_{J'}$.
Note that $\tte_{\vecJ,J'}=0$ if $J'\in\vecJ$.
We set
$\nbigsbar_{\vecJ}(\Lambda):=\nbigs(\Lambda,\vecJ)$.
For each $\vecJ'\in\nbigs_{\vecJ}(\Lambda)$,
let $\sigmabar_{\vecJ}(\vecJ')$
denote the image of $\sigma(\vecJ')$.
Thus, we obtain a fan
$\Sigmabar_{\vecJ}(\Lambda)$
in $\ttNbar_{\vecJ}(\Lambda)_{\real}$.
Let $\ttXbar_{\vecJ}(\Lambda)$
denote the associated toric manifold.
We have
$\overline{\ttO(\vecJ)}=\ttXbar_{\vecJ}(\Lambda)$.
If $\Lambda\in\vecJ$,
then $\ttXbar_{\vecJ}(\Lambda)$ is projective,
and in particular,
$\ttXbar_{\vecJ}(\Lambda)_{\geq 0}$ is compact.

The $\ttT_{\ttNbar_{\vecJ}(\Lambda)}$-action
on $\ttXbar_{\vecJ}(\Lambda)$
is induced by
the $\ttT_{\ttN(\Lambda)}$-action on $\overline{\ttO(\vecJ)}$.
The orbit decompositions are the same:
\[
 \overline{\ttO(\vecJ)}
 =\coprod_{\vecJ'\in\nbigsbar_{\vecJ}(\Lambda)}
 \ttO(\vecJ'),
 \quad\quad
 \overline{O(\vecJ)}_{\geq 0}
 =\coprod_{\vecJ'\in\nbigsbar_{\vecJ}(\Lambda)}
 \ttO(\vecJ')_{\geq 0}.
\]

We set
$I_i=J_{i}\setminus J_{i-1}$
for $i=1,\ldots,|\vecJ|+1$,
where $J_0:=\emptyset$
and $J_{|\vecJ|+1}:=\Lambda$.
Note that
$I_{|\vecJ|+1}=\emptyset$
if $J_{|\vecJ|}=\Lambda$.
There exist the following natural isomorphisms:
\[
 \ttNbar_{\vecJ}(\Lambda)\simeq
 \prod_{i=1}^{|\vecJ|}
 \ttNbar(I_i)
 \times
 \ttN(I_{|\vecJ|+1}),
\quad
 \nbigsbar_{\vecJ}(\Lambda)\simeq
 \prod_{i=1}^{|\vecJ|}
 \nbigsbar(I_i)\times
 \nbigs(I_{|\vecJ|+1}),
\quad
 \Sigmabar_{\vecJ}(\Lambda)\simeq
 \prod_{i=1}^{|\vecJ|}
 \Sigmabar(I_i)\times
 \Sigma(I_{|\vecJ|+1}).
\]
They induce
\[
 \ttXbar_{\vecJ}(\Lambda)\simeq
 \prod_{i=1}^{|\vecJ|}
 \ttXbar(I_i)\times
 \ttX(I_{|\vecJ|+1}),
 \quad\quad
 \ttXbar_{\vecJ}(\Lambda)_{\geq 0}
 \simeq
 \prod_{i=1}^{|\vecJ|}
 \ttXbar(I_i)_{\geq 0}
 \times
 \ttX(I_{|\vecJ|+1})_{\geq 0}.
\]

\subsubsection{Comparison with the blow up at the origin}

For $J\in 2^{\Lambda}\setminus\{\Lambda\}$,
let $\sigma_{\ttY}(J,\Lambda)$
denote the cone in $\ttN(\Lambda)$
generated by
$\tte_{i}$ $(i\in J)$ and $\tte_{\Lambda}$,
and let $\sigma_{\ttY}(J)$ denote the cone
generated by $\tte_i$ $(i\in J)$.
Let $\Sigma_{\ttY}(\Lambda)$ be the fan
of the cones
$\sigma_{\ttY}(J)$, $\sigma_{\ttY}(J,\Lambda)$
$(J\in 2^{\Lambda}\setminus\{\Lambda\})$.
Let $\ttY(\Lambda)$ denote the associated toric manifold.
It is well known to be isomorphic to
the complex blow up of $\cnum^{\Lambda}$
at $(0,\ldots,0)$.
For any $\sigma\in\Sigma_{\ttY}(\Lambda)$,
let $\ttO_{\ttY}(\sigma)$ denote the corresponding orbit.

For any $J\in 2^{\Lambda}\setminus\Lambda$,
let $\sigmabar_{\ttY}(J)$ be the cone
in $\ttNbar(\Lambda)$
generated by
$\bar{\tte}_i$ $(i\in J)$.
(See \S\ref{subsection;22.2.7.20}
for $\ttNbar(\Lambda)$.)
Let $\Sigmabar_{\ttY}(\Lambda)$
be the fan of the cones
$\sigmabar_{\ttY}(J)$ $(J\in 2^{\Lambda}\setminus\{\Lambda\})$.
Let $\ttYbar(\Lambda)$ denote the associated toric manifold,
which is well known to be isomorphic to
$\proj^{\Lambda}$.
It is the exceptional fiber of the blow up
$\ttY(\Lambda)\lrarr \cnum^{\Lambda}$ at $(0,\ldots,0)$.
For any $\sigmabar(J)\in\Sigmabar_{\ttY}(\Lambda)$,
let $\ttO_{\ttYbar}(J)$ denote the corresponding orbit.

Let $\vecJ=(J_i\,|\,i=1,\ldots,m)\in\nbigs(\Lambda)$.
If $J_m\neq\Lambda$,
$\sigma(\vecJ)$ is contained in the cone
generated by $\tte_i$ $(i\in J_m)$.
If $J_m=\Lambda$,
$\sigma(\vecJ)$ is contained in the cone
generated by $\tte_i$ $(i\in J_{m-1})$
and $\tte_{\Lambda}$.
Hence, there exist the natural morphisms
$\ttX(\Lambda)\lrarr \ttY(\Lambda)$
and
$\ttX(\Lambda)_{\geq 0}\lrarr \ttY(\Lambda)_{\geq 0}$.
(See \cite[\S3.3]{Toric-Varieties-Book}.)

Let $\vecJ=(J_1,\ldots,J_{m-1},\Lambda)\in\nbigsbar(\Lambda)$.
Note that $\sigmabar(\vecJ)$ is contained
in the cone generated by $\bar{\tte}_i$ $(i\in J_{m-1})$.
There exist the natural morphisms
$\pi_{\ttYbar,\ttXbar}:
\ttXbar(\Lambda)\lrarr \ttYbar(\Lambda)$
and
$\pi_{\ttYbar,\ttXbar,\geq 0}:
\ttXbar(\Lambda)_{\geq 0}\lrarr
\ttYbar(\Lambda)_{\geq 0}$.

There exists the decompositions
\[
\ttYbar(\Lambda)
=\coprod_{J\subsetneq\Lambda}\ttO_{\Ybar}(J),
\quad\quad
\ttYbar(\Lambda)_{\geq 0}
=\coprod_{J\subsetneq\Lambda}
 \ttO_{\Ybar}(J)_{\geq 0}.
\]
For $J\subsetneq\Lambda$,
let $\nbigsbar_{\ttX,J}(\Lambda)$
denote the set of
$\vecJ=(J_1,\ldots,J_{m-1},\Lambda)\in\nbigsbar(\Lambda)$
such that $J_{m-1}=J$.
We have
\[
\pi_{\ttYbar,\ttXbar}^{-1}(\ttO_{\ttYbar}(J))
=\coprod_{\vecJ\in\nbigsbar_{\ttX,J}(\Lambda)}
\ttO(\vecJ)
=\ttO_{\ttYbar}(J)\times
 \ttXbar(J),
\]
\[
 \pi_{\ttYbar,\ttXbar,\geq 0}^{-1}(\ttO_{\ttYbar}(J)_{\geq 0})
=\coprod_{\vecJ\in\nbigsbar_{X,J}(\Lambda)}
\ttO(\vecJ)_{\geq 0}
=\ttO_{\ttYbar}(J)_{\geq 0}
\times
 \ttXbar(J)_{\geq 0}.
\]

\subsubsection{Weakly constructible sheaves}
\label{subsection;22.2.7.21}

Let $S$ be a topological space
with a sheaf of algebras $A_S$.
We regard
$\nbigs_{\ttYbar}(\Lambda)=2^{\Lambda}\setminus\{\Lambda\}$
as a category
where a morphism $J_1\to J_2$ is given by
the relation $J_1\supset J_2$.
By the procedure in \S\ref{subsection;22.2.4.10},
we obtain the following lemma.
\begin{lem}
There exist natural equivalences
 $\Fun(\nbigs_{\ttYbar}(\Lambda),\Mod(A_S))
 \simeq
 \Mod^{\wc}(\ttYbar(\Lambda)_{\geq 0};A_S)$
and
 $\Fun(\nbigs_{\ttYbar}(\Lambda),\ttC(A_S))
 \simeq
 \ttC^{\wc}(\ttYbar(\Lambda)_{\geq 0};A_S)$.
\hfill\qed
\end{lem}

We regard $\nbigsbar(\Lambda)$ as a category
where a morphism
$\vecJ_1\to\vecJ_2$ corresponds to the relation
$\vecJ_1\in\nbigsbar(\Lambda,\vecJ_2)$,
which is equivalent to
$\ttO(\vecJ_1)\subset \overline{\ttO(\vecJ_2)}$.
By the procedure in \S\ref{subsection;22.2.4.10},
we obtain the following lemma.
\begin{lem}
There exist natural equivalences
$\Fun(\nbigsbar(\Lambda),\Mod(A_S))
\simeq
\Mod^{\wc}(\ttXbar(\Lambda)_{\geq 0};A_S)$ 
and
$\Fun(\nbigsbar(\Lambda),\ttC(A_S))
\simeq
\ttC^{\wc}(\ttXbar(\Lambda)_{\geq 0};A_S)$.
\hfill\qed
\end{lem}

Let $\pi_{\ttYbar,\ttXbar}:
\ttXbar(\Lambda)_{\geq 0}\times S
\lrarr \ttYbar(\Lambda)_{\geq 0}\times S$
denote the natural map.
The following lemma is obvious.

\begin{lem}
For any $\ttF^{\bullet}\in
\ttC^{\wc}(\Ybar(\Lambda)_{\geq 0};A_S)$,
$\pi_{\ttYbar,\ttXbar}^{-1}(\ttF^{\bullet})$
is an object of
$\ttC^{\wc}(\Xbar(\Lambda)_{\geq 0};A_S)$.
The natural morphism
$\ttF^{\bullet}\lrarr
R\pi_{\ttYbar,\ttXbar\ast}
\pi_{\ttYbar,\ttXbar}^{-1}(\ttF^{\bullet})$
in the derived category
is an isomorphism.
\hfill\qed
\end{lem}

Let $\nbigf\in\Fun(\nbigsbar(\Lambda),\Mod(A_S))$.
Let us recall the associated complex in \S\ref{subsection;22.2.4.11}
with the terminology in this situation.
For any $\vecJ=(J_1,\ldots,J_{m-1},\Lambda)\in \nbigsbar(\Lambda)$,
let $\cnum(\vecJ)^{\lor}$ denote the subspace of
$\cnum\langle 2^{\Lambda}\setminus\{\emptyset,\Lambda\}\rangle^{\lor}$
generated by
$\ttv_{J_1}^{\lor}\wedge\cdots\wedge \ttv_{J_{m-1}}^{\lor}$.
For $k\leq 0$, we set
\[
 \nbigc^k(\nbigf)=
 \bigoplus_{|\vecJ|=-k+1}
 \nbigf(\vecJ)\otimes\cnum(\vecJ)^{\lor}.
\]
If $J\in \vecJ$,
the morphism
$\nbigf(\vecJ)\to\nbigf(\vecJ\setminus J)$
and
the inner product of $\ttv_J$
induces
\[
 \nbigf(\vecJ)\otimes\cnum(\vecJ)^{\lor}
 \lrarr
 \nbigf(\vecJ\setminus J)\otimes\cnum(\vecJ\setminus J)^{\lor}.
\]
They induce
$d:\nbigc^{k}(\nbigf)\lrarr\nbigc^{k+1}(\nbigf)$.
Thus, we obtain a complex
$\nbigc^{\bullet}(\nbigf)$.

Let $\pi_{\ttXbar}:\ttXbar(\Lambda)_{\geq 0}\times S\to S$
denote the projection.
For $\ttF^{\bullet}\in
\ttC^{\wc}(\Xbar(\Lambda)_{\geq 0};A_S)$,
let $\nbigf_{\ttF^{\bullet}}$ be the corresponding
object in $\Fun(\nbigsbar(\Lambda),\ttC(A_S))$.
Then, there exists the following natural isomorphism
in the derived category of $A_S$-complexes:
\begin{equation}
\label{eq;22.4.25.22}
 R\pi_{\ttXbar\ast}(\ttF^{\bullet})
 \simeq
 \Tot
 \nbigc^{\bullet}(\nbigf_{\ttF^{\bullet}})[-|\Lambda|+1].
\end{equation}

\subsubsection{Contractibility of some subspaces}
\label{subsection;22.1.22.13}

Let $\emptyset\neq\Lambda_0\subsetneq\Lambda$.
For $\vecJ=(J_1\subsetneq\cdots\subsetneq J_{m-1}\subsetneq \Lambda_0)
\in\nbigsbar(\Lambda_0)$,
we set
$\vecJ^{\Lambda}=
(J_1\subsetneq\cdots\subsetneq J_{m-1}\subsetneq\Lambda)
\in\nbigsbar(\Lambda)$
and $\vecJ\cdot\Lambda:=
 (J_1\subsetneq \cdots\subsetneq
 J_{m-1}\subsetneq \Lambda_0\subsetneq\Lambda)$.

\begin{lem}
Let $\vecJ\in\nbigsbar(\Lambda_0)$.
Let $\vecJ_1,\ldots,\vecJ_k\in\nbigsbar(\Lambda_0)$
be refinements of $\vecJ$.
Then,
\begin{equation}
 \overline{\ttO(\vecJ\cdot\Lambda)}_{\geq 0}
 \cup
\bigcup_{j=1}^k
 \overline{\ttO(\vecJ_j^{\Lambda})}_{\geq 0}
\end{equation}
is contractible.
\end{lem}
\pf
If $k=0$, the claim is clear.
We use an induction on $k$.
Note that
\begin{equation}
\label{eq;22.1.22.10}
 \overline{\ttO(\vecJ_k^{\Lambda})}_{\geq 0}
 \cap
\left(
  \overline{\ttO(\vecJ\cdot\Lambda)}_{\geq 0}
 \cup
\bigcup_{j=1}^{k-1}
\overline{\ttO(\vecJ_j^{\Lambda})}_{\geq 0}
\right)
 =\overline{\ttO(\vecJ_k\cdot\Lambda)}_{\geq 0}
 \cup
 \bigcup_{j=1}^{k-1}
 \bigl(
 \overline{\ttO(\vecJ_k^{\Lambda})}_{\geq 0}
 \cap
 \overline{\ttO(\vecJ_j^{\Lambda})}_{\geq 0}
 \bigr).
\end{equation}
If
$\overline{\ttO(\vecJ_k^{\Lambda})}_{\geq 0}
 \cap
 \overline{\ttO(\vecJ_j^{\Lambda})}_{\geq 0}\neq\emptyset$,
there exists a refinement
$\vecJ_{j,k}$ of $\vecJ_k$
such that
\[
 \overline{\ttO(\vecJ_k^{\Lambda})}_{\geq 0}
 \cap
 \overline{\ttO(\vecJ_j^{\Lambda})}_{\geq 0}
=\overline{\ttO(\vecJ_{j,k}^{\Lambda})}_{\geq 0}.
\]
Hence, by the assumption of the induction,
(\ref{eq;22.1.22.10}) is contractible.
Thus, we obtain the claim of the lemma.
\hfill\qed

\section{Purity theorem and some consequences}

\subsection{Vector spaces
with a commuting tuple of nilpotent endomorphisms}

\subsubsection{Linearized intersection complexes}
\label{subsection;22.2.6.5}

Let $\Lambda$ be a non-empty finite set.
Let $V$ be a finite dimensional $\cnum$-vector space
equipped with a tuple of mutually commuting
nilpotent endomorphisms
$\vecN=(N_i\,|\,i\in \Lambda)$.
For any $J\in 2^{\Lambda}\setminus\{\emptyset\}$,
we set $N_J:=\prod_{j\in J}N_j$.
We set $N_{\emptyset}=\id_V$.

We use the notation in \S\ref{subsection;22.2.4.30}.
For $-|\Lambda|\leq k\leq 0$,
we set
\[
 \IntC^{k}(V,\vecN)
=\bigoplus_{\substack{J\subset\Lambda\\ |J|=|\Lambda|+k}}
\Image (N_J)\otimes\cnum(J)
\subset
 V\otimes \bigwedge^{|\Lambda|+k} (\cnum\cdot\Lambda).
\]
For $k<-|\Lambda|$ or $k>0$,
we set $\IntC^{k}(V,\vecN)=0$.
We obtain the morphism
\begin{equation}
\label{eq;21.11.2.1}
 d:\IntC^{k}(V,\vecN)\lrarr
 \IntC^{k+1}(V,\vecN),
 \quad
 d(v)=
 \sum_{i\in\Lambda}\ttv_i\wedge(N_i(v)).
\end{equation}
Thus, we obtain a complex
$\IntC^{\bullet}(V,\vecN)$.

For $0\leq k\leq |\Lambda|$,
we set
\[
 \IntC^k_{c}(V,\vecN)=\bigoplus
 _{\substack{J\subset \Lambda\\ |J|=|\Lambda|-k}}
 \Image (N_J)\otimes\cnum(J)
 \subset
 V\otimes \bigwedge^{|\Lambda|-k}(\cnum\cdot\Lambda).
\]
For $k<0$ or $k>|\Lambda|$,
we set $\IntC^k_{c}(V,\vecN)=0$.
We obtain the morphism
\begin{equation}
\label{eq;21.11.2.2}
 d:\IntC^{k}_{c}(V,\vecN)
 \lrarr
 \IntC^{k+1}_{c}(V,\vecN),
 \quad
 d(v)=\sum_{i\in\Lambda} \ttv_i^{\lor}\vdash v.
\end{equation}
Thus, we obtain a complex
$\IntC^{\bullet}_{c}(V,\vecN)$.

The identity map on $\Image N_{\Lambda}\otimes\cnum(\Lambda)$
induces the naturally defined morphism
$\IntC^{\bullet}_c(V,\vecN)
\lrarr
\IntC^{\bullet}(V,\vecN)$.

\subsubsection{Linearized partially intersection complexes}

We introduce a generalization.
Let $\Lambda_0\subset\Lambda$.
We set
$N^{(\Lambda_0)}_J:=
N_{J\setminus\Lambda_0}$.

For $-|\Lambda|\leq k\leq 0$,
we set
\[
 \IntC^k(V,\vecN;\ast\Lambda_0)
=\bigoplus_{\substack{J\subset\Lambda\\ |J|=|\Lambda|+k}}
\Image (N^{(\Lambda_0)}_J)\otimes\cnum(J).
\]
For $k<-|\Lambda|$ or $k>0$,
we set $\IntC^{k}(V,\vecN;\ast\Lambda_0)=0$.
We define
$d:\IntC^k(V,\vecN;\ast\Lambda_0)
\lrarr\IntC^{k+1}(V,\vecN;\ast\Lambda_0)$
by (\ref{eq;21.11.2.1}),
and we obtain a complex
$\IntC^{\bullet}(V,\vecN;\ast\Lambda_0)$.

For $0\leq k\leq |\Lambda|$,
we set
\[
 \IntC^k_{c}(V,\vecN;\ast\Lambda_0)=\bigoplus
 _{\substack{J\subset \Lambda\\ |J|=|\Lambda|-k}}
 \Image (N_J^{(\Lambda_0)})\otimes\cnum(J).
\]
For $k<0$ or $k>|\Lambda|$,
we set $\IntC^k_{c}(V,\vecN;\ast\Lambda_0)=0$.
We obtain the morphism
$d:\IntC^{k}_{c}(V,\vecN;\ast\Lambda_0)
 \lrarr
 \IntC^{k+1}_{c}(V,\vecN;\ast\Lambda_0)$
by (\ref{eq;21.11.2.2}).
Thus, we obtain a complex
$\IntC^{\bullet}_{c}(V,\vecN;\ast\Lambda_0)$.
The following lemma is clear.
\begin{lem}
If $\Lambda_0\neq\emptyset$,
$\IntC^{\bullet}_{c}(V,\vecN;\ast\Lambda_0)$
is acyclic.
\hfill\qed
\end{lem}

There exists the naturally defined morphism
$\IntC^{\bullet}_{c}(V,\vecN;\ast\Lambda_0)
\lrarr
\IntC^{\bullet}(V,\vecN;\ast\Lambda_0)$.
For $\Lambda_0\subset\Lambda_1$,
there exist the naturally defined commutative diagram:
\[
 \begin{CD}
  \IntC^{\bullet}_{c}(V,\vecN;\ast\Lambda_0)
  @>>>
  \IntC^{\bullet}(V,\vecN;\ast\Lambda_0)
  \\
  @VVV @VVV \\
  \IntC^{\bullet}_{c}(V,\vecN;\ast\Lambda_1)
  @>>>
  \IntC^{\bullet}(V,\vecN;\ast\Lambda_1).
 \end{CD}
\]

\subsubsection{Linearized punctured intersection complexes}

Let $\Ccat(\Lambda)$ denote the category of
non-empty subsets of $\Lambda$,
i.e.,
objects are non-empty subsets of $\Lambda$,
and morphisms
$\iota_{\Lambda_1,\Lambda_0}:
\Lambda_0\lrarr\Lambda_1$ are the inclusions.
Let $\nbigf^{\bullet}$ be a functor from
$\Ccat(\Lambda)$ to the category of bounded complexes
of finite dimensional $\cnum$-vector spaces.
For $k\in\seisuu_{\geq 0}$,
we define
\[
 \Ccheck^k(\nbigf^{\bullet}):=
 \bigoplus_{\substack{\emptyset\neq\Lambda_0\subset\Lambda\\
  |\Lambda_0|=k+1
  }}
  \nbigf^{\bullet}(\Lambda_0)\otimes \cnum(\Lambda_0).
\]
For $k<0$, we set
$\Ccheck^k(\nbigf^{\bullet}):=0$.
Let $d: \Ccheck^k(\nbigf^{\bullet})\lrarr\Ccheck^{k+1}(\nbigf^{\bullet})$
denote the morphism induced by
\[
 \sum_{\substack{\emptyset\neq\Lambda_0\subset \Lambda\\ |\Lambda_0|=k+1}}
 \sum_{i\in\Lambda\setminus \Lambda_0}
 \nbigf^{\bullet}(\iota_{\Lambda_0\cup\{i\},\Lambda_0})
 \ttv_i\wedge .
\]
Thus, we obtain a double complex
$\Ccheck^{\bullet}(\nbigf^{\bullet})$,
and the associated total complex
$\Tot\Ccheck^{\bullet}(\nbigf^{\bullet})$.

By using the above construction
to the functors
$\IntC^{\bullet}(V,\vecN;\ast\bullet)$
defined by
$\Lambda_0\longmapsto
\IntC^{\bullet}(V,\vecN;\ast\Lambda_0)$,
we obtain
\[
 \IntC^{\bullet}_{\punc}(V,\vecN):=
 \Tot\Ccheck^{\bullet}
 \IntC^{\bullet}(V,\vecN;\ast\bullet).
\]
Similarly, from the functor
$\IntC^{\bullet}_c(V,\vecN;\ast\bullet)$
defined by
$\Lambda_0\longmapsto
\IntC^{\bullet}_c(V,\vecN;\ast\Lambda_0)$,
we obtain
\[
 \IntC^{\bullet}_{\punc,c}(V,\vecN):=
 \Tot\Ccheck^{\bullet}
 \IntC_c^{\bullet}(V,\vecN;\ast\bullet).
\]
Note that $\IntC^{\bullet}_{\punc,c}(V,\vecN)$ is acyclic.

\vspace{.1in}
 
The natural transform
$\IntC^{\bullet}_c(V,\vecN;\ast\bullet)\lrarr
\IntC^{\bullet}(V,\vecN;\ast\bullet)$
induces
\[
 \IntC^{\bullet}_{\punc,c}(V,\vecN)\lrarr
 \IntC^{\bullet}_{\punc}(V,\vecN).
\]

\subsubsection{Some morphisms}

By the construction,
there exist the following naturally defined
commutative diagram:
\begin{equation}
\label{eq;22.2.4.47}
\begin{CD}
 \IntC^{\bullet}_c(V,\vecN)
 @>>>
 \IntC^{\bullet}_{\punc,c}(V,\vecN),
 \\ @VVV @VVV \\
  \IntC^{\bullet}(V,\vecN)
 @>>>
 \IntC^{\bullet}_{\punc}(V,\vecN).
\end{CD}
\end{equation}
We shall prove the following proposition
in \S\ref{subsection;22.2.4.20}.

\begin{prop}
\label{prop;21.11.3.1}
The natural morphism
\[
 \Cone\Bigl(
 \IntC^{\bullet}_c(V,\vecN)
 \lrarr
 \IntC^{\bullet}(V,\vecN)
 \Bigr)
 \lrarr
 \Cone\Bigl(
 \IntC^{\bullet}_{\punc,c}(V,\vecN)
 \lrarr
 \IntC^{\bullet}_{\punc}(V,\vecN)
 \Bigr)
\]
is a quasi-isomorphism. 
Here,
$\Cone(A^{\bullet}\to B^{\bullet})$
denotes the mapping cone
for a morphism of complexes
$A^{\bullet}\to B^{\bullet}$ 
(see {\rm\cite[\S1.4]{Kashiwara-Schapira}}).
\end{prop}

There exists the following natural commutative diagram:
\begin{equation}
\label{eq;21.11.3.2}
 \begin{CD}
  \IntC^{\bullet}(V,\vecN)@>>>
  \IntC^{\bullet}_{\punc}(V,\vecN)\\
@VVV @V{b}VV\\
 \Cone\Bigl(
 \IntC^{\bullet}_c(V,\vecN)
 \to
 \IntC^{\bullet}(V,\vecN)
 \Bigr)
 @>{a}>>
 \Cone\Bigl(
 \IntC^{\bullet}_{\punc,c}(V,\vecN)
 \to
 \IntC^{\bullet}_{\punc}(V,\vecN)
 \Bigr)
 \\
 @VVV \\
 \IntC^{\bullet}_{c}(V,\vecN)[1].
 \end{CD}
\end{equation}
Here, $a$ and $b$ are quasi-isomorphisms.
\begin{cor}
\label{cor;21.11.2.3}
We obtain
\[
 H^{k}(\IntC^{\bullet}(V,\vecN))
 \simeq
 \left\{
 \begin{array}{ll}
 H^k(\IntC^{\bullet}_{\punc}(V,\vecN))  & (k<0)\\
 0 & (k\geq 0)
 \end{array}
\right.
\]
\[
 H^k(\IntC_c^{\bullet}(V,\vecN))
 \simeq
 \left\{
\begin{array}{ll}
 H^{k-1}(\IntC^{\bullet}_{\punc}(V,\vecN))& (k\geq 1)\\
 0 & (k\leq 0).
\end{array}
 \right.
\]
\[
 H^k(\IntC_{\punc}^{\bullet}(V,\vecN))
 \simeq
 \left\{
\begin{array}{ll}
 H^k(\IntC^{\bullet}(V,\vecN))
  & (k<0)\\
 \\
 H^{k+1}(\IntC^{\bullet}_c(V,\vecN))
  & (k\geq 0)
\end{array}
 \right.
\]
\end{cor}
\pf
It is easy to see that
$H^k(\IntC^{\bullet}(V,\vecN))=0$ $(k\geq 0)$
and
$H^k(\IntC^{\bullet}_c(V,\vecN))=0$ $(k\leq 0)$
by the construction.
By the diagram (\ref{eq;21.11.3.2})
and Proposition \ref{prop;21.11.3.1},
we obtain the following long exact sequence:
\[
 \cdots\lrarr
 H^k\bigl(
 \IntC_c^{\bullet}(V,\vecN)
 \bigr) 
 \lrarr
 H^k\bigl(
 \IntC^{\bullet}(V,\vecN)
 \bigr)
 \lrarr
  H^k\bigl(
 \IntC_{\punc}^{\bullet}(V,\vecN)
 \bigr)
 \lrarr
 H^{k+1}\bigl(
 \IntC_c^{\bullet}(V,\vecN)
 \bigr) \lrarr\cdots.
\]
Then, we obtain the claim of the corollary.
\hfill\qed

\subsubsection{The associated $\nbigd$-modules}

We set $X=(\proj^1)^{\Lambda}$.
Let $p_i:X\lrarr\proj^1$ denote the projection onto
the $i$-th component $(i\in\Lambda)$.
We set
$H^{(\infty)}_{i}=p_i^{-1}(\infty)$
and
$H^{(0)}_{i}=p_i^{-1}(0)$.
We set
$H^{(\infty)}=\bigcup_{i\in\Lambda} H^{(\infty)}_{i}$,
$H^{(0)}=\bigcup_{i\in\Lambda}H^{(0)}_{i}$
and $H=H^{(0)}\cup H^{(\infty)}$.
For any $K\subset\Lambda$,
we set $H^{(0)}(K):=\bigcup_{i\in K}H^{(0)}_i$
and $H^{(\infty)}(K):=\bigcup_{i\in K}H^{(\infty)}_i$.
In general, for a regular holonomic $\nbigd_X$-module $\nbign$
and a hypersurface $H'$ of $X$,
we set
$\nbign(\ast H'):=\nbign\otimes_{\nbigo_X}\nbigo_X(\ast H')$
and
$\nbign(!H'):=\DDD_X(\DDD_X(\nbign)(\ast H'))$,
where $\DDD_X$ denote the duality functor for
coherent $\nbigd_X$-modules.

We consider $\nbigv=V\otimes\nbigo_{X}(\ast H)$
with the integrable connection $\nabla$
determined by
\[
\nabla(1\otimes v)=-\sum 1\otimes N_i(v)\,dz_i/z_i
\]
for $v\in V$.
We naturally regard $\nbigv$ as a $\nbigd_X$-module.
For $K^{(0)},K^{(\infty)}\subset\Lambda$,
we obtain
the $\nbigd_X$-module
\[
\nbigv\Bigl(
!\Bigl(
 H^{(0)}(K^{(0)})
 \cup
 H^{(\infty)}(K^{(\infty)})
\Bigr)
\Bigr).
\]

\begin{lem}
\label{lem;22.3.22.1}
For $K^{(0)}_1,K^{(\infty)}_1\subset\Lambda$,
there exists a natural isomorphism:
\begin{multline}
 \nbigv\Bigl(
!\Bigl(
 H^{(0)}(K^{(0)})
 \cup
 H^{(\infty)}(K^{(\infty)})
\Bigr)
 \Bigr)
 \Bigl(
 \ast \Bigl(
 H^{(0)}(K^{(0)}_1)
 \cup
 H^{(\infty)}(K^{(\infty)}_1)
 \Bigr)
 \Bigr)
 \simeq
 \\
 \nbigv\Bigl(
 !\Bigl(
 H^{(0)}(K^{(0)}\setminus K^{(0)}_1)
 \cup
 H^{(\infty)}(K^{(\infty)}\setminus K^{(\infty)}_1)
 \Bigr)
 \Bigr).
\end{multline}
\end{lem}
\pf
It follows from \cite[Lemma 3.1.11]{Mochizuki-holonomic-Betti}.
\hfill\qed

\vspace{.1in}

Let $\nbigm$ denote the image of
$\nbigv(!H^{(0)})\to \nbigv(\ast H^{(0)})$.
We obtain the $\nbigd_{X}$-modules
$\nbigm(\star H^{(\infty)})$ $(\star=!,\ast)$.
\begin{lem}
We have $\nbigm(\ast H^{(\infty)})=\nbigm$. 
The image of
$\nbigv\Bigl(
 !\bigl(
 H^{(0)}\cup H^{(\infty)}
 \bigr)
 \Bigr)
 \lrarr
 \nbigv(!H^{(\infty)})$
is equal to
$\nbigm(!H^{(\infty)})$.
\end{lem}
\pf
It follows from Lemma \ref{lem;22.3.22.1}.
\hfill\qed

\vspace{.1in}

We set $O=(0,\ldots,0)\in X$.
Let $j:X\setminus\{O\}\lrarr X$
denote the open embedding.
We have the natural functors of
the derived category of coherent algebraic $\nbigd$-modules
$j_{+}:D^b_{\coh}(\nbigd_{X\setminus O})
\lrarr D^b_{\coh}(\nbigd_{X})$.
There exists the following natural morphism
in $D^b_{\coh}(\nbigd_X)$:
\begin{equation}
\label{eq;22.2.4.31}
\Cone\bigl(
 \nbigm(!H^{(\infty)})
 \lrarr
 \nbigm(\ast H^{(\infty)})
 \bigr)
\lrarr
\Cone\bigl(
 j_{+}j^{\ast}\nbigm(!H^{(\infty)})
 \lrarr
 j_{+}j^{\ast}\nbigm(\ast H^{(\infty)})
 \bigr).
\end{equation}

\begin{lem}
\label{lem;22.2.4.32}
The morphism {\rm(\ref{eq;22.2.4.31})} is 
a quasi-isomorphism.
\end{lem}
\pf
At each point of $X\setminus H^{(\infty)}$,
both sides of (\ref{eq;22.2.4.31}) are acyclic.
At each point of $H^{(\infty)}$,
the morphism (\ref{eq;22.2.4.31}) is an isomorphism.
Hence, the claim is clear.
\hfill\qed

\subsubsection{Representatives of
$j_+j^{\ast}\nbigm(\star H^{(\infty)})$}

We set $\nbign=\nbigm(\star H^{(\infty)})$ $(\star=!,\ast)$.
For $k\in\seisuu_{\geq 0}$,
we set
\[
 \Ccheck^k_{\punc}(\nbign):=
 \bigoplus_{\substack{K\subset\Lambda\\ |K|=k+1}}
 \nbign(\ast H^{(0)}(K))
 \otimes\cnum(K).
\]
We set $\Ccheck^k_{\punc}(\nbign)=0$ if $k<0$ or $k+1>|\Lambda|$.
For $\emptyset\neq K\subset \Lambda$
and $i\in \Lambda\setminus K$,
the natural morphism
$\nbign(\ast (H^{(0)}(K)))
\lrarr
\nbign(\ast (H^{(0)}(K\cup\{i\})))$
and the exterior product of $\ttv_i$
induces
\[
\nbign(\ast (H^{(0)}(K)))
\otimes\cnum(K)
\lrarr
\nbign(\ast (H^{(0)}(K\cup\{i\})))
\otimes\cnum(K\cup\{i\}).
\]
They define a morphism
$\Ccheck_{\punc}^k(\nbign)
\to
\Ccheck^{k+1}_{\punc}(\nbign)$.
Thus, we obtain a complex
$\Ccheck^{\bullet}_{\punc}(\nbign)$,
which represents
$j_+j^{\ast}\nbign$.

There exists the following naturally defined commutative diagram:
\[
 \begin{CD}
  \nbigm(!H^{(\infty)})
  @>>>
  \Ccheck^{\bullet}_{\punc}(\nbigm(!H^{(\infty)}))
  \\
  @VVV @VVV \\
 \nbigm(\ast H^{(\infty)})
  @>>>
  \Ccheck^{\bullet}_{\punc}(\nbigm(\ast H^{(\infty)})).
 \end{CD}
\]

\subsubsection{Some complexes associated with
$\nbigm\bigl(\ast H^{(0)}(K)\bigr)$}
\label{subsection;22.4.28.10}

We set $\nbign_{K}=\nbigm\bigl(\ast H^{(0)}(K)\bigr)$
for a subset $K\subset\Lambda$.
For $-|\Lambda|\leq j\leq 0$,
we set
\[
 C^j(\nbign_K):=
 \bigoplus_{|J|=-j}
 \nbign_K(!H^{(0)}(J))
 \otimes
 \cnum(J)^{\lor}.
\]
For $j<-|\Lambda|$ or $j>0$,
we set
$C^j(\nbign_K)=0$.
For a non-empty subset $J\subset \Lambda$
and $i\in\Lambda$,
the natural morphism
$\nbign_K(!H^{(0)}(J))
\lrarr
\nbign_K(!H^{(0)}(J\setminus\{i\}))$
and the inner product of $\ttv_i$
induce
\[
 \nbign_K(!H^{(0)}(J))
 \otimes\cnum(J)^{\lor}
\lrarr
\nbign_K(!H^{(0)}(J\setminus\{i\}))
 \otimes\cnum(J\setminus\{i\})^{\lor}.
\]
They define a morphism
$C^j(\nbign_K)
\lrarr
C^{j+1}(\nbign_K)$.
Thus, we obtain a complex
$C^{\bullet}(\nbign_K)$.
There exists a natural morphism
$\nbign_K\to C^{\bullet}(\nbign_K)$.

Let $a_X$ denote the natural morphism
from $X$ to a one point set.
\begin{lem}
\label{lem;22.2.4.40}
The induced morphism
 $a_{X\dagger}(\nbign_K)
 \lrarr
 a_{X\dagger}(C^{\bullet}(\nbign_K))$ 
is an isomorphism
in the derived category of complexes of $\cnum$-vector spaces.
\end{lem}
\pf
Let $J\neq \emptyset$.
For any $j\in J$,
by \cite[Lemma 3.1.11]{Mochizuki-holonomic-Betti},
there exists a natural isomorphism
\[
\nbign_K(!H^{(0)}(J))
\simeq
\Bigl(
\bigl(
 \nbign_K(!H^{(0)}(J))
 \bigr)
 (!H^{(0)}_j)
  \Bigr)
 (\ast H^{(\infty)}_j).
\]
Moreover, the $\nbigd_X$-module is monodromic with respect to
with respect to the fibration
$X\to (\proj^1)^{\Lambda\setminus\{j\}}$.
By Lemma \ref{lem;22.2.4.49} below,
we obtain $a_{X\dagger}(\nbign_K(!H^{(0)}(J)))=0$.
It implies the claim of Lemma \ref{lem;22.2.4.40}.
\hfill\qed

\vspace{.1in}

For $K_1\subset K_2$,
we obtain $\nbign_{K_1}\lrarr \nbign_{K_2}$,
which induces
$C^{\bullet}(\nbign_{K_1})
\lrarr
 C^{\bullet}(\nbign_{K_2})$.
Note $\nbigm=\nbigm(\ast H^{(\infty)})$.
We obtain
a double complex
$C^{\bullet}\Ccheck^{\bullet}_{\punc}(\nbigm(\ast H^{(\infty)}))$.
There exists the following naturally defined commutative diagram:
\[
\begin{CD}
 \nbigm(\ast H^{(\infty)}) @>>>
 \Ccheck^{\bullet}_{\punc}(\nbigm(\ast H^{(\infty)}))\\
 @VVV @VVV \\
 C^{\bullet}(\nbigm(\ast H^{(\infty)})) @>>>
 \Tot C^{\bullet}\bigl(
 \Ccheck^{\bullet}_{\punc}(\nbigm(\ast H^{(\infty)}))
 \bigr)
\end{CD}
\]
We obtain the following commutative diagram:
\begin{equation}
\label{eq;22.2.4.43}
\begin{CD}
 a_{X\dagger}\nbigm(\ast H^{(\infty)}) @>>>
 a_{X\dagger}\Ccheck^{\bullet}_{\punc}(\nbigm(\ast H^{(\infty)}))\\
 @V{\simeq}VV @V{\simeq}VV \\
 a_{X\dagger}C^{\bullet}(\nbigm(\ast H^{(\infty)})) @>>>
 a_{X\dagger}\Tot C^{\bullet}\bigl(
 \Ccheck^{\bullet}_{\punc}(\nbigm(\ast H^{(\infty)}))
 \bigr).
\end{CD}
\end{equation}
The vertical arrows are isomorphisms
because of Lemma \ref{lem;22.2.4.40}.

Let $\phi^{(0)}_{z_i}$ denote the vanishing cycle functor
along $z_i$.
(See \cite{beilinson2}. See also \cite[\S4]{Mochizuki-MTM}
for the explanation in the context of $\nbigr$-modules.)

\begin{lem}
\label{lem;22.2.4.45}
We choose a bijection $\Lambda\simeq\ellsitabar=\{1,\ldots,\ell\}$.
The isomorphisms
$
 \cnum(\ellsitabar\setminus J)\simeq
\cnum(J)^{\lor}\otimes\cnum(\ellsitabar)$
induces the following isomorphism:
\[
 a_{X\dagger}
 \phi^{(0)}_{z_1}\circ\cdots\circ
 \phi^{(0)}_{z_{\ell}}
 C^{\bullet}(\nbign_K)
 \otimes\cnum(\ellsitabar)
 \simeq
 \IntC^{\bullet}(V,\vecN;\ast K).
\]
\end{lem}
\pf
The $\nbigd_X$-module
$\nbign_K(!H^{(0)}(J))$ is the image of
$\nbigv\Bigl(
 !\bigl(H^{(0)}(J\cup(\Lambda\setminus K))\bigr)
 \Bigr)
 \lrarr
 \nbigv\Bigl(
 !\bigl(H^{(0)}(J)\bigr)
 \Bigr)$.
The morphism
\[
 a_{X\dagger}
 \phi^{(0)}_{z_1}\circ\cdots\circ
 \phi^{(0)}_{z_{\ell}}
 \nbigv\Bigl(
 !\bigl(H^{(0)}(J\cup(\Lambda\setminus K))\bigr)
 \Bigr)
 \lrarr
  a_{X\dagger}
 \phi^{(0)}_{z_1}\circ\cdots\circ
 \phi^{(0)}_{z_{\ell}}
 \nbigv\Bigl(
 !\bigl(H^{(0)}(J)\bigr)
 \Bigr)
\]
is identified with
$\prod_{i\in \Lambda\setminus (J\cup K)}N_i:V\to V$.
We obtain
\[
 a_{X\dagger}
 \phi^{(0)}_{z_1}\circ\cdots\circ
 \phi^{(0)}_{z_{\ell}}
 \nbign_K(!H^{(0)}(J))
 \simeq
 \Image N^{(K)}_{\Lambda\setminus J}.
\]
Then, we obtain the claim of Lemma \ref{lem;22.2.4.45}
by the construction.
\hfill\qed

\subsubsection{Some complexes associated with
$\nbigm(!H^{(\infty)}\ast H^{(0)}(K))$}

Let $\nbign_K:=\nbigm(!H^{(\infty)}\ast H^{(0)}(K))$
for a subset $K\subset\Lambda$.
For $0\leq k\leq |\Lambda|$,
we set
\[
 C_c^{k}(\nbign_K):=
 \bigoplus_{|J|=k}
 \nbign_K(\ast H^{(0)}(J))
 \otimes
 \cnum(J)^{\lor}. 
\]
For $J\subset\Lambda$
and $i\in \Lambda\setminus J$,
the natural morphism
$\nbign_K(\ast H^{(0)}(J))
\lrarr
\nbign_K(\ast H^{(0)}(J\cup\{i\}))$
and the exterior product of $\ttv_i^{\lor}$
induce
\[
\nbign_K(\ast H^{(0)}(J))
\otimes\cnum(J)^{\lor}
\lrarr
\nbign_K(\ast H^{(0)}(J\cup\{i\}))
\otimes
\cnum(J\cup\{i\})^{\lor}.
\]
They induce a morphism
$C_c^{k}(\nbign_K)\lrarr
 C_c^{k+1}(\nbign_K)$.
Thus, we obtain a complex
$C_c^{\bullet}(\nbign_K)$.
There exists a naturally defined morphism
$C_c^{\bullet}(\nbign_K)\to \nbign_K$.

\begin{lem}
\label{lem;22.2.4.41}
 The induced morphism
$a_{X\dagger}C_c^{\bullet}(\nbign_K)\to
a_{X\dagger}\nbign_K$
is an isomorphism
in the derived category of $\cnum$-complexes. 
\end{lem}
\pf
If $J\neq\emptyset$,
we have $a_{X\dagger}\nbign_K(\ast H^{(0)}(J))=0$
by Lemma \ref{lem;22.2.4.49} below.
It implies the claim of the lemma.
\hfill\qed

\vspace{.1in}
For $K_1\subset K_2$,
we obtain
$\nbign_{K_1}\to \nbign_{K_2}$,
which induce
$C_c^{\bullet}(\nbign_{K_1})\to
C_c^{\bullet}(\nbign_{K_2})$.
Hence, we obtain a double complex
$C^{\bullet}_c\Ccheck^{\bullet}_{\punc}(\nbigm(!H^{(\infty)}))$.
There exists the following naturally defined commutative diagram:
\[
 \begin{CD}
  C^{\bullet}_c(\nbigm(!H^{(\infty)}))
  @>>>
  \Tot
  C^{\bullet}_c\Ccheck^{\bullet}_{\punc}(\nbigm(!H^{(\infty)}))
  \\
  @VVV @VVV \\
  \nbigm(!H^{(\infty)})
  @>>>
  \Ccheck^{\bullet}_{\punc}(\nbigm(!H^{(\infty)})).
 \end{CD}
\]
We obtain the following commutative diagram:
\begin{equation}
\label{eq;22.2.4.44}
 \begin{CD}
  a_{X\dagger}C^{\bullet}_c(\nbigm(!H^{(\infty)}))
  @>>>
  a_{X\dagger}\Tot
  C^{\bullet}_c\Ccheck^{\bullet}_{\punc}(\nbigm(!H^{(\infty)}))
  \\
  @V{\simeq}VV @V{\simeq}VV \\
  a_{X\dagger}\nbigm(!H^{(\infty)})
  @>>>
  a_{X\dagger}\Ccheck^{\bullet}_{\punc}(\nbigm(!H^{(\infty)})).
 \end{CD}
\end{equation}
The vertical morphisms are isomorphisms
because of Lemma \ref{lem;22.2.4.41}.

\begin{lem}
\label{lem;22.2.4.46}
We choose a bijection $\Lambda\simeq\ellsitabar=\{1,\ldots,\ell\}$.
The isomorphisms
$\cnum(\ellsitabar\setminus J)\simeq
\cnum(J)^{\lor}\otimes\cnum(\ellsitabar)$
induces the following isomorphism:
\[
 a_{X\dagger}
 \phi^{(0)}_{z_1}\circ\cdots\circ
 \phi^{(0)}_{z_{\ell}}
 C^{\bullet}_c(\nbign_K)
 \otimes\cnum(\ellsitabar)
 \simeq
 \IntC^{\bullet}_c(V,\vecN;\ast K).
\]
\end{lem}
\pf
The $\nbigd_X$-module $\nbign_K(\ast H^{(0)}(J))$
is isomorphic to the image of
$\nbigv\Bigl(!
 \bigl(
 H^{(\infty)}
 \cup
 H^{(0)}\bigl(
 \Lambda\setminus(K\cup J)
 \bigr)
 \bigr)
 \Bigr)
 \lrarr
 \nbigv(!H^{(\infty)})$.
The morphism
\[
 a_{X\dagger}
 \phi^{(0)}_{z_1}\circ\cdots\circ
 \phi^{(0)}_{z_{\ell}}
 \nbigv\Bigl(!\bigl(
 H^{(\infty)}
 \cup
 H^{(0)}\bigl(
\Lambda\setminus(K\cup J)
 \bigr)
 \bigr)
 \Bigr)
 \lrarr
  a_{X\dagger}
 \phi^{(0)}_{z_1}\circ\cdots\circ
 \phi^{(0)}_{z_{\ell}}
 \nbigv(!H^{(\infty)}) 
\]
is identified with
$N^{(K)}_{\Lambda\setminus J}$.
Hence, we obtain the following natural isomorphism:
\[
 a_{X\dagger}
 \phi^{(0)}_{z_1}\circ\cdots\circ
 \phi^{(0)}_{z_{\ell}}
 \nbign_K(\ast H^{(0)}(J))
 \simeq
 \Image N^{(K)}_{\Lambda\setminus J}.
\]
Then, we obtain the claim of Lemma \ref{lem;22.2.4.46}
by the construction.
\hfill\qed

\subsubsection{Proof of Proposition \ref{prop;21.11.3.1}}
\label{subsection;22.2.4.20}

There exists the following naturally defined commutative diagram:
\[
 \begin{CD}
  C^{\bullet}_c(\nbigm(!H^{(\infty)}))
  @>{f_1}>>
  \Tot
  C^{\bullet}_c\Ccheck^{\bullet}_{\punc}(\nbigm(!H^{(\infty)}))
  \\
  @V{f_2}VV @V{f_3}VV \\
   C^{\bullet}(\nbigm(\ast H^{(\infty)}))
  @>{f_4}>>
  \Tot
  C^{\bullet}\Ccheck^{\bullet}_{\punc}(\nbigm(\ast H^{(\infty)})).
 \end{CD}
\]
We consider the complex $\nbigc^{\bullet}$
obtained as the total complex of the following double complex:
\begin{multline}
 0\lrarr
 C^{\bullet}_c(\nbigm(!H^{(\infty)}))
 \stackrel{f_1-f_2}{\lrarr}
 \Tot
 C^{\bullet}_c\Ccheck^{\bullet}_{\punc}(\nbigm(!H^{(\infty)}))
 \oplus
 C^{\bullet}(\nbigm(\ast H^{(\infty)}))
  \\
 \stackrel{f_3+f_4}{\lrarr}
  \Tot
 C^{\bullet}\Ccheck^{\bullet}_{\punc}(\nbigm(\ast H^{(\infty)}))
 \lrarr 0.
\end{multline}
Because of Lemma \ref{lem;22.2.4.32}
and the commutative diagrams (\ref{eq;22.2.4.43})
and (\ref{eq;22.2.4.44}),
$a_{X\dagger}\nbigc^{\bullet}$
is acyclic.

We choose a bijection
$\Lambda\simeq\ellsitabar=\{1,\ldots,\ell\}$.
We note that the cohomological support of
$\nbigc^{\bullet}$ is contained in the origin $O$.
Hence, we obtain that
$a_{X\dagger}
\phi_{z_1}^{(0)}\circ
\cdots\circ
\phi_{z_{\ell}}^{(0)}
\nbigc^{\bullet}$
is acyclic.
By Lemma \ref{lem;22.2.4.45} and Lemma \ref{lem;22.2.4.46},
the induced commutative diagram
\[
 \begin{CD}
  a_{X\dagger}\phi^{(0)}_{z_1}\circ\cdots\circ\phi^{(0)}_{z_{\ell}}
  C^{\bullet}_c(\nbigm(!H^{(\infty)}))
  @>>>
    a_{X\dagger}\phi^{(0)}_{z_1}\circ\cdots\circ\phi^{(0)}_{z_{\ell}}
  \Tot
  C^{\bullet}_c\Ccheck^{\bullet}_{\punc}(\nbigm(!H^{(\infty)}))
  \\
  @VVV @VVV \\
    a_{X\dagger}\phi^{(0)}_{z_1}\circ\cdots\circ\phi^{(0)}_{z_{\ell}}
   C^{\bullet}(\nbigm(\ast H^{(\infty)}))
  @>>>
    a_{X\dagger}\phi^{(0)}_{z_1}\circ\cdots\circ\phi^{(0)}_{z_{\ell}}
  \Tot
  C^{\bullet}\Ccheck^{\bullet}_{\punc}(\nbigm(\ast H^{(\infty)})).
 \end{CD}
\]
is identified with (\ref{eq;22.2.4.47}).
Hence, we obtain the claim of Proposition \ref{prop;21.11.3.1}.
\hfill\qed

\subsubsection{Appendix: Some vanishing}

Let $L$ be an ample line bundle on an algebraic manifold $X$.
Let $\proj_L$ denote the projective completion of $L$.
Let $H_{\infty}=\proj_L\setminus L$.
Let $H_0\subset L$ denote the $0$-section of $L$.
Let $\pi:\proj_L\lrarr X$ denote the projection.

Let $M$ be a monodromic algebraic
$\nbigd_{\proj_L}(\ast (H_0\cup H_{\infty}))$-module.

\begin{lem}
\label{lem;22.2.4.49}
 We have
$\pi_{+}\bigl(
M(!H_{\infty})(\ast H_{0})
\bigr)=0$. 
\end{lem}
\pf
It is enough to consider the case
$L=\cnum\times X$.
Let $t$ denote the standard coordinate of $\cnum$.
Let $\gminiv=t\del_t$ denote the Euler vector field on $\proj_L$.
There exists a generalized eigen decomposition
$\pi_{\ast}M=\bigoplus _{\alpha\in\cnum}M_{\alpha}$
with respect to the action of $\gminiv$,
where $\gminiv-\alpha$ is nilpotent on $M_{\alpha}$.
We fix an appropriate total order on $\cnum$.
We have
\[
\pi_{\ast}\lefttop{0}V_{0}(M)(\ast H_{\infty})
=\bigoplus_{\alpha\geq -1}M_{\alpha},
\quad
\pi_{\ast}\lefttop{\infty}V_{-1}(M)(\ast H_{0})
=\bigoplus_{\alpha\leq -2}M_{\alpha}.
\]
The relative de Rham complex of
$M(!H_{\infty})(\ast H_{0})$
is represented by
\[
 \lefttop{\infty}V_{-1}
 \lefttop{0}V_0(M)
 \lrarr
 \lefttop{\infty}V_{-1}
 \lefttop{0}V_0(M)\otimes dt/t,
 \quad
 v\longmapsto \gminiv(v)\,dt/t.
\]
Because
$R\pi_{\ast}
 \lefttop{\infty}V_{-1}
 \lefttop{0}V_0(M)=0$,
we obtain the claim of the lemma.
\hfill\qed

\subsection{Purity theorem for polarized mixed twistor structures
and some consequences}

\subsubsection{Purity theorem}
\label{subsection;22.2.7.5}

Let $\Lambda$ be a finite subset.
Let $w\in\seisuu$.
Let $(\nbigv,W,\vecf)$ be a $(w,\Lambda)$-polarizable
mixed twistor structure.
(See \S\ref{subsection;22.3.24.30}.)
For $-|\Lambda|\leq k\leq 0$,
we set
\[
 \IntC^k(\nbigv,W,\vecf)
=\bigoplus_{\substack{J\subset|\Lambda|\\ |J|=|\Lambda|+k }}
 (\Image f_J,W)\otimes\cnum(J).
\]
For $k<-|\Lambda|$ or $k>0$,
we set $\IntC^k(\nbigv,W,\vecf)=0$.
The differential is given as in \S\ref{subsection;22.2.6.5}.
Thus, we obtain a complex
of mixed twistor structures
$\IntC^{\bullet}(\nbigv,W,\vecf)$.

\begin{prop}
\label{prop;22.2.6.3}
 $\Gr^W_{m}H^k\bigl(
  \IntC^{\bullet}(\nbigv,W,\vecf)
 \bigr)=0$ if $m>w+k+|\Lambda|$.
\end{prop}
\pf
We set $\nbigv^{(0)}:=\Gr^W_{\bullet}(\nbigv)$.
It is equipped with the weight filtration
$W_m(\nbigv^{(0)})=\bigoplus_{j\leq m}\Gr^W_j(\nbigv)$
and the induced commuting tuple of morphisms
$f^{(0)}_i:(\nbigv^{(0)},W)\to (\nbigv^{(0)},W)\otimes\Tate(-1)$,
and $(\nbigv^{(0)},W,\vecf^{(0)})$
is a $(w,\Lambda)$-polarizable mixed twistor structure.
Indeed, it underlies a polarizable mixed Hodge structure,
as explained in \cite[\S3.7.7]{mochi2}.
The claim of Proposition \ref{prop;22.2.6.3}
for $(\nbigv^{(0)},W,\vecf^{(0)})$
follows from the purity theorem for
polarized mixed Hodge structures,
due to Cattani-Kaplan-Schmid \cite{cks2}
and Kashiwara-Kawai \cite{k3},
which was originally conjectured by Deligne.
By general Lemma \ref{lem;22.3.24.50},
we obtain the claim of Proposition \ref{prop;22.2.6.3}
for $(\nbigv,W,\vecf)$.
\hfill\qed

\vspace{.1in}
We note that
$H^0\bigl(
\IntC^{\bullet}(\nbigv,W,\vecf)
\bigr)=0$
by the construction.
We also note that the differentials are strictly compatible
with the weight filtrations $W$.
\begin{cor}
\label{cor;22.4.3.30}
The natural morphism
\[
 W_{w+|\Lambda|-1}
 \IntC^{\bullet}(\nbigv,W,\vecf)
 \lrarr
  \IntC^{\bullet}(\nbigv,W,\vecf)
\]
is a quasi-isomorphism.
\hfill\qed
\end{cor}

We obtain the dual statement of Proposition \ref{prop;22.2.6.3}
from the dual statement of the purity theorem
in \cite{k3}.
For $0\leq k\leq|\Lambda|$,
we set
\[
 \IntC^k_c(\nbigv,W,\vecf)
 =\bigoplus_{\substack{J\subset\Lambda \\ |J|=|\Lambda|-k}}
 (\Image f_J,W)\otimes\Tate(-|\Lambda|+|J|)
 \otimes\cnum(J).
\]
If $J_1\subset J_2$,
there exists the following natural monomorphism
as subobjects of
$(\nbigv,W)\otimes\Tate(-|\Lambda|)$:
\[
 (\Image (f_{J_2}),W)\otimes\Tate(-|\Lambda|+|J_2|)
\lrarr
 (\Image (f_{J_1}),W)\otimes\Tate(-|\Lambda|+|J_1|).
\]
Hence, for $J\subset\Lambda$ and $i\in J$,
the inner product of $\ttv_i^{\lor}$ induces
\[
 (\Image (f_J),W)\otimes\Tate(-|\Lambda|+|J|)
 \otimes\cnum(J)
 \lrarr
 (\Image (f_{J\setminus\{i\}}),W)\otimes\Tate(-|\Lambda|+|J\setminus\{i\}|)
 \otimes\cnum(J\setminus\{i\}).
\]
They induce the morphism of mixed twistor structures
\[
 d:\IntC^k_c(\nbigv,W,\vecf)
 \lrarr
 \IntC^{k+1}_c(\nbigv,W,\vecf)
\]
such that $d\circ d=0$.
Thus, we obtain a complex of mixed twistor structure.

\begin{cor}
 $\Gr^W_{m}
 H^k\bigl(
 \IntC^{\bullet}_c(\nbigv,W,\vecf)
 \bigr)=0$
if $m<w+|\Lambda|+k$.
\end{cor}
\pf
It follows from
\cite[Theorem 4.0.2]{k3}
as in the case of Proposition \ref{prop;22.2.6.3}.
\hfill\qed

\vspace{.1in}

We note that
$H^0\bigl(
\IntC^{\bullet}_c(\nbigv,W,\vecf)
\bigr)=0$
by the construction.
\begin{cor}
 $W_{w+|\Lambda|}\IntC^{\bullet}_c(\nbigv,W,\vecf)$
is acyclic.
\hfill\qed
\end{cor}

\subsubsection{Linearized partially intersection complexes}
\label{subsection;22.2.7.4}

Let $\Lambda_0\subset\Lambda$.
For any $J\subset\Lambda$,
we set $f^{(\Lambda_0)}_J:=f_{J\setminus\Lambda_0}$.
For $-|\Lambda|\leq k\leq 0$,
we set
\[
 \IntC^k(\nbigv,W,\vecf;\ast \Lambda_0)
=\bigoplus_{\substack{J\subset\Lambda\\ |J|=|\Lambda|+k}}
\bigl(
\Image(f^{(\Lambda_0)}_J),W
\bigr)\otimes\Tate(-|J|+|J\setminus\Lambda_0|)
\otimes
\cnum(J).
\]
We obtain a complex of mixed twistor structures
$\IntC^{\bullet}(\nbigv,W,\vecf;\ast\Lambda_0)$.

For $0\leq k\leq |\Lambda|$,
we set
\[
 \IntC^k_c(\nbigv,W,\vecf;\ast\Lambda_0)
=\bigoplus_{\substack{J\subset\Lambda\\ |J|=|\Lambda|-k}}
\bigl(
\Image(f^{(\Lambda_0)}_J),W
\bigr)
\otimes\Tate(-|\Lambda|+|J\setminus\Lambda_0|)
\otimes
\cnum(J).
\]
We obtain a complex of mixed twistor structures
$\IntC^{\bullet}_c(\nbigv,W,\vecf;\ast\Lambda_0)$.

There exists the naturally defined morphism
$\IntC^{\bullet}_c(\nbigv,W,\vecf;\ast\Lambda_0)
\lrarr
\IntC^{\bullet}(\nbigv,W,\vecf;\ast\Lambda_0)$
of complex of mixed twistor structures.

\subsubsection{Linearized punctured intersection complexes}

Let $\nbigf^{\bullet}$ be a functor
from $\Ccat(\Lambda)$ to the category of mixed twistor structures.
For $k\in\seisuu_{\geq 0}$,
we set
\[
 \Ccheck^k(\nbigf^{\bullet})
 :=
 \bigoplus_{\substack{\Lambda_0\subset\Lambda\\
  |\Lambda_0|=k+1}}
  \nbigf^{\bullet}(\Lambda_0)\otimes\cnum(\Lambda_0).
\]
For $k<0$, we set 
$\Ccheck^k(\nbigf^{\bullet})=0$.
Thus, we obtain a double complex of mixed twistor structures
$\Ccheck^{\bullet}(\nbigf^{\bullet})$,
and the associated total complex
$\Tot\Ccheck^{\bullet}(\nbigf^{\bullet})$.

By applying the above construction to the functor
$\IntC^{\bullet}(\nbigv,W,\vecf;\ast\star)$
given by
$\Lambda_0\mapsto
\IntC^{\bullet}(\nbigv,W,\vecf;\ast \Lambda_0)$,
we obtain the following complex of mixed twistor structures:
\[
 \IntC^{\bullet}_{\punc}(\nbigv,W,\vecf)
:=\Tot \Ccheck^{\bullet}\IntC^{\bullet}(\nbigv,W,\vecf;\ast\star),
\]
Similarly, 
from the functor
$\IntC^{\bullet}_c(\nbigv,\vecf;\ast\star)$
given by
$\Lambda_0\mapsto
\IntC^{\bullet}_c(\nbigv,W,\vecf;\ast \Lambda_0)$,
we obtain the following complex of mixed twistor structures:
\[
 \IntC^{\bullet}_{\punc,c}(\nbigv,W,\vecf)
 :=\Tot \Ccheck^{\bullet}\IntC_c^{\bullet}(\nbigv,W,\vecf;\ast\star).
\]
The natural transform
$\IntC^{\bullet}(\nbigv,W,\vecf;\ast\star)
\lrarr
\IntC^{\bullet}_c(\nbigv,W,\vecf;\ast\star)$
induces
\[
 \IntC^{\bullet}_{\punc,c}(\nbigv,W,\vecf)
 \lrarr
 \IntC^{\bullet}_{\punc}(\nbigv,W,\vecf).
\]

\subsubsection{Some morphisms}

By the construction,
there exist the following naturally defined
commutative diagram:
\[
\begin{CD}
 \IntC^{\bullet}_c(\nbigv,W,\vecf)
 @>>>
 \IntC^{\bullet}_{\punc,c}(\nbigv,W,\vecf),
 \\ @VVV @VVV \\
  \IntC^{\bullet}(\nbigv,W,\vecf)
 @>>>
 \IntC^{\bullet}_{\punc}(\nbigv,W,\vecf).
\end{CD}
\]
We obtain the following lemma from
Proposition \ref{prop;21.11.3.1}
and Lemma \ref{lem;22.4.25.1}.
\begin{lem}
\label{lem;22.2.7.6}
The natural morphism
\[
\Cone\Bigl(
 \IntC^{\bullet}_c(\nbigv,W,\vecf)
 \lrarr
 \IntC^{\bullet}(\nbigv,W,\vecf)
 \Bigr)
 \lrarr \\
 \Cone\Bigl(
 \IntC^{\bullet}_{\punc,c}(\nbigv,W,\vecf)
 \lrarr
 \IntC^{\bullet}_{\punc}(\nbigv,W,\vecf)
 \Bigr)
\]
is a quasi-isomorphism.
\hfill\qed
\end{lem}

We note the following natural commutative diagram,
which is an enrichment of (\ref{eq;21.11.3.2}):
\[
 \begin{CD}
  \IntC^{\bullet}(\nbigv,W,\vecf)@>>>
  \IntC^{\bullet}_{\punc}(\nbigv,W,\vecf)\\
@VVV @V{b}VV\\
 \Cone\Bigl(
 \IntC^{\bullet}_c(\nbigv,W,\vecf)
 \to
 \IntC^{\bullet}(\nbigv,W,\vecf)
 \Bigr)
 @>{a}>>
 \Cone\Bigl(
 \IntC^{\bullet}_{\punc,c}(\nbigv,W,\vecf)
 \to
 \IntC^{\bullet}_{\punc}(\nbigv,W,\vecf)
 \Bigr)
 \\
 @VVV \\
 \IntC^{\bullet}_{c}(\nbigv,W,\vecf)[1].
 \end{CD}
\]
Here, $a$ and $b$ are quasi-isomorphisms.
We obtain the following proposition from
Corollary \ref{cor;21.11.2.3}.
\begin{prop}
There exists the following
isomorphisms of mixed twistor structures:
\[
 H^k\bigl(
 \IntC^{\bullet}_{\punc}(\nbigv,W,\vecf)
 \bigr)
 \simeq
 \left\{
 \begin{array}{ll}
  H^k\bigl(
   \IntC^{\bullet}(\nbigv,W,\vecf)
   \bigr)& (k<0)\\
 \mbox{{}} \\
 H^{k+1}\bigl(
 \IntC_c^{\bullet}(\nbigv,W,\vecf)
 \bigr)& (k\geq 0).
 \end{array}
 \right.
\] 
\hfill\qed
\end{prop}

\begin{cor}
For $k<0$,
we have the following commutative diagram of isomorphisms:
\[
\begin{CD}
 H^k W_{w+|\Lambda|-1}
 \IntC^{\bullet}
 (\nbigv,W,\vecf) 
 @>{\simeq}>>
 H^k W_{w+|\Lambda|-1}
 \IntC_{\punc}^{\bullet}
 (\nbigv,W,\vecf) \\
 @V{\simeq}VV @V{\simeq}VV \\
  H^k 
 \IntC^{\bullet}
 (\nbigv,W,\vecf) 
 @>{\simeq}>>
 H^k 
 \IntC_{\punc}^{\bullet}
 (\nbigv,W,\vecf).
\end{CD}
\]
For $k\geq 0$,
we have
\[
 H^kW_{w+|\Lambda|-1}
 \IntC_{\punc}^{\bullet}(\nbigv,W,\vecf)=0.
\]
\hfill\qed
\end{cor}

\begin{cor}
\label{cor;22.2.6.10}
The following morphisms are quasi-isomorphisms:
\[
\IntC^{\bullet}
 (\nbigv,W,\vecf) 
 \llarr
 W_{w+|\Lambda|-1}\IntC^{\bullet}
 (\nbigv,W,\vecf) 
\lrarr
 W_{w+|\Lambda|-1}\IntC_{\punc}^{\bullet}
 (\nbigv,W,\vecf).
\]
\hfill\qed
\end{cor}

\subsubsection{Tuple of weight filtrations}

Let $\Lambda_0,\Lambda_1\subset\Lambda$ be subsets.
Let $\vecf_{\Lambda_1}=(f_i\,|\,i\in\Lambda_1)$.
Let $\vecK=(K_1,\ldots,K_p)\in\nbigs(\Lambda)$.
Let $\vecm\in\seisuu^{p}$.
For $k\leq 0$,
we set
\[
 \lefttop{\vecK}W_{\vecm}\IntC^k(\nbigv,W,\vecf_{\Lambda_1};\ast\Lambda_0)
 =\\
 \bigoplus_{\substack{J\subset\Lambda_1 \\ |J|=|\Lambda_1|+k}}
 \bigcap_{i=1}^p
 \lefttop{K_i}
 W_{m_i-|K_i\cap J|
 -|K_i\cap\Lambda_0\cap J|}
 \bigl(
 \Image f^{(\Lambda_0)}_J,W
 \bigr)
  \otimes\Tate(-|J|+|J\setminus\Lambda_0|)
 \otimes\cnum(J).
 \]
Thus, we obtain a complex
$\lefttop{\vecK}W_{\vecm}
\IntC^{\bullet}(\nbigv,W,\vecf_{\Lambda_1};\ast\Lambda_0)$
of mixed twistor structures.
If $\Lambda_0=\emptyset$,
we set
\[
\lefttop{\vecK}W_{\vecm}
\IntC^{\bullet}(\nbigv,W,\vecf_{\Lambda_1}):=
\lefttop{\vecK}W_{\vecm}
\IntC^{\bullet}(\nbigv,W,\vecf_{\Lambda_1};\ast\emptyset).
\]

\begin{lem}
\label{lem;22.2.7.1}
The following morphism is a quasi-isomorphism:
\[
\lefttop{\Lambda}W_{|\Lambda|-1}
 \IntC^{\bullet}(\nbigv,W,\vecf)
 \lrarr
 \IntC^{\bullet}(\nbigv,W,\vecf).
\]
\end{lem}
\pf
Because
$\lefttop{\Lambda}W_m\IntC^{\bullet}(\nbigv,W,\vecf)
=W_{m+w}\IntC^{\bullet}(\nbigv,W,\vecf)$,
we obtain the claim of
Lemma \ref{lem;22.2.7.1}
from Corollary \ref{cor;22.4.3.30}.

\hfill\qed

\vspace{.1in}

Let $L\subset \Lambda$ be any subset.
We define the decreasing filtration $F_L$
on
$\lefttop{\vecK}W_{\vecm}
\IntC^{\bullet}(\nbigv,W,\vecf_{\Lambda_1};\ast\Lambda_0)$
as follows:
\begin{equation}
\label{eq;22.2.7.30}
 F_L^j
 \lefttop{\vecK}W_{\vecm}\IntC^k(\nbigv,W,\vecf_{\Lambda_1};\ast\Lambda_0)
 =\\
 \bigoplus_{
 \substack{
 J\subset\Lambda_1 \\
 |J|=|\Lambda_1|+k \\
 |J\setminus L|\geq j
 }}
 \bigcap_{i=1}^p
 \lefttop{K_i}
 W_{m_i-|K_i\cap J|
 -|K_i\cap\Lambda_0\cap J|}
 \bigl(
 \Image f^{(\Lambda_0)}_J,W
 \bigr)
  \otimes\Tate(-|J|+|J\setminus\Lambda_0|)
 \otimes\cnum(J).
\end{equation}
By construction,
the following lemma is clear.
\begin{lem}
\label{lem;22.2.8.1}
$F_L^j
 \lefttop{\vecK}W_{\vecm}
 \IntC^{\bullet}(\nbigv,W,\vecf_{\Lambda_1};\ast\Lambda_0)$
is a subcomplex.
There exist the following natural isomorphisms:
\begin{equation}
\label{eq;22.2.8.2}
 \Gr_{F_L}^j
  \lefttop{\vecK}W_{\vecm}
 \IntC^{\bullet}(\nbigv,W,\vecf_{\Lambda_1};\ast\Lambda_0)
 \simeq
 \bigoplus_{\substack{
 I\subset\Lambda_1\\
 I\cap L=\emptyset\\
 |I|=j
 }}
 \lefttop{\vecK}W_{\vecm}
 \IntC^{\bullet}
 \Bigl(
 \Image f_I^{(\Lambda_0)}\otimes\Tate(-|I|+|I\setminus\Lambda_0|),
 W,\vecf_{\Lambda_1\cap L};\ast\Lambda_0
 \Bigr)
 \Bigl[
 \bigl|\Lambda_1\setminus (I\cup L)\bigr|
 \Bigr].
\end{equation}
\hfill\qed
\end{lem}

We set
$\vecw(\vecK):=(|K_1|-1,|K_2|-1,\ldots,|K_p|-1)
\in\seisuu^p$.
\begin{prop}
\label{prop;22.2.7.3}
The following natural morphism is a quasi-isomorphism:
\[
 \lefttop{\vecK}
  W_{\vecw(\vecK)}
 \IntC^{\bullet}(\nbigv,W,\vecf)
 \lrarr
 \IntC^{\bullet}(\nbigv,W,\vecf).
\] 
As a result,
\[
 \lefttop{\Lambda}
 W_{|\Lambda|-1}
 \lefttop{\vecK}
 W_{\vecw(\vecK)}
 \IntC^{\bullet}(\nbigv,W,\vecf)
 \lrarr
 \lefttop{\Lambda}
 W_{|\Lambda|-1}
 \IntC^{\bullet}(\nbigv,W,\vecf)
 \lrarr
 \IntC^{\bullet}(\nbigv,W,\vecf)
\] 
are quasi-isomorphisms.
\end{prop}
\pf
We use an induction on $p$.
Let us study the case $p=1$,
i.e., $\vecK$ consists of $K\subset\Lambda$.
\begin{lem}
\label{lem;22.3.24.101}
We set $K':=\Lambda\setminus K$.
Let $J\subset K'$ be any subset.
The following morphism is a quasi-isomorphism:
\begin{equation}
\label{eq;22.2.7.31}
 \lefttop{K}W_{|K|-1}
 \IntC^{\bullet}(\Image f_J,W,\vecf_{K})
 \lrarr
 \IntC^{\bullet}(\Image f_J,W,\vecf_{K}).
\end{equation} 
\end{lem}
\pf
Recall that $(\Image f_J,\vecf)$ is
a $(w+|J|,\Lambda)$-polarizable mixed twistor structure.
For $i\in K'$,
we set $H_i:=\cnum^{K'\setminus\{i\}}
\subset\cnum^{K'}$
and $H(K')=\bigcup_{i\in K'}H_i$.
By Corollary \ref{cor;22.3.24.100},
there exists a neighbourhood $U$ of
$(0,\ldots,0)$ in $\cnum^{K'}$
such that
$(\TNIL(\Image f_J,\vecf_{K'}),\vecf_K)_{|U\setminus H(K')}$
is a $(w+|J|,K)$-polarizable mixed twistor structure on
$U\setminus H(K')$.
Let $P$ be any point of $U\setminus H(K')$.
By Lemma \ref{lem;22.2.7.1},
we obtain that the following is a quasi-isomorphism:
\[
 \lefttop{K}W_{|K|-1}
 \IntC^{\bullet}
 \Bigl(
\TNIL(\Image f_J,\vecf_{K'}),\vecf_K
 \Bigr)
 \lrarr
  \IntC^{\bullet}
 \Bigl(
\TNIL(\Image f_J,\vecf_{K'}),\vecf_K
 \Bigr).
\]
By the construction of twistor nilpotent orbit,
we obtain the claim of Lemma \ref{lem;22.3.24.101}.
\hfill\qed

\vspace{.1in}

We consider the filtration $F^{\bullet}=F_K^{\bullet}$
for
$\lefttop{K}W_{|K|-1}
\IntC^{\bullet}(\nbigv,W,\vecf)$
and 
$\IntC^{\bullet}(\nbigv,W,\vecf)$
defined by (\ref{eq;22.2.7.30}) with $L=K$.
Because (\ref{eq;22.2.7.31}) are quasi-isomorphisms
for any $J$ such that $J\cap K=\emptyset$,
the following induced morphisms are quasi-isomorphisms
for any $j$:
\[
 \Gr_F^j
\lefttop{K}W_{|K|-1}
\IntC^{\bullet}(\nbigv,W,\vecf)
\lrarr
 \Gr_F^j
\IntC^{\bullet}(\nbigv,W,\vecf).
\]
Hence, we obtain that
\[
  \lefttop{K}W_{|K|-1}
 \IntC^{\bullet}(\nbigv,W,\vecf)
 \lrarr
 \IntC^{\bullet}(\nbigv,W,\vecf)
\]
is a quasi-isomorphism.

Let us consider the general case.
For $\vecK=(K_1,\ldots,K_p)$,
we set $\vecK'=(K_1,\ldots,K_{p-1})$.
Suppose that we have already proved the claim for $\vecK'$.
Let $J\subset\Lambda\setminus K_{p-1}$ be any subset.
By using the twistor nilpotent orbit with respect to
$\vecf_{\Lambda\setminus K_{p-1}}$,
and by using the hypothesis of the induction,
we obtain that the following is a quasi-isomorphism
(see the proof of Lemma \ref{lem;22.3.24.101}):
\begin{equation}
\label{eq;22.2.7.32}
 \lefttop{\vecK'}W_{\vecw(\vecK')}
 \IntC^{\bullet}(\Image f_J,W,\vecf_{K_{p-1}})
 \lrarr
 \IntC^{\bullet}(\Image f_J,W,\vecf_{K_{p-1}}).
\end{equation}
Let $I$ be any subset of $\Lambda\setminus K_p$.
We consider the filtrations $F^{\bullet}=F^{\bullet}_{K_{p-1}}$
on the complexes
$\lefttop{\vecK'}W_{\vecw(\vecK')}
 \IntC^{\bullet}(\Image f_I,W,\vecf_{K_{p}})$
and 
$\IntC^{\bullet}(\Image f_I,W,\vecf_{K_{p}})$
defined by (\ref{eq;22.2.7.30}) with $L=K_{p-1}$.
For any $J\subset K_p\setminus K_{p-1}$,
\begin{equation}
 \lefttop{\vecK'}W_{\vecw(\vecK')}
 \IntC^{\bullet}(\Image f_{I\cup J},W,\vecf_{K_{p-1}})
 \lrarr
 \IntC^{\bullet}(\Image f_{I\cup J},W,\vecf_{K_{p-1}})
\end{equation}
is a quasi-isomorphism.
Hence, the induced morphisms
\[
 \Gr_F^j
\lefttop{\vecK'}W_{\vecw(\vecK')}
 \IntC^{\bullet}(\Image f_I,W,\vecf_{K_{p}})
\lrarr
  \Gr_F^j
 \IntC^{\bullet}(\Image f_I,W,\vecf_{K_{p}})
\]
are quasi-isomorphisms for any $j$.
We obtain that the following is a quasi-isomorphism:
\[
 \lefttop{\vecK'}W_{\vecw(\vecK')}
 \IntC^{\bullet}(\Image f_I,W,\vecf_{K_{p}})
 \lrarr
 \IntC^{\bullet}(\Image f_I,W,\vecf_{K_{p}}).
\]
By considering the twistor nilpotent orbit
with respect to $\vecf_{\Lambda\setminus K_p}$
(see the proof of Lemma \ref{lem;22.3.24.101}),
we obtain that
\[
 \lefttop{K_p}W_{|K_p|-1}
 \lefttop{\vecK'}W_{\vecw(\vecK')}
 \IntC^{\bullet}(\Image f_I,W,\vecf_{K_{p}})
 \lrarr
 \lefttop{K_p}W_{|K_p|-1}
 \IntC^{\bullet}(\Image f_I,W,\vecf_{K_{p}})
\]
is a quasi-isomorphism.
By using the above argument again,
we obtain that
\[
 \lefttop{K_p}W_{|K_p|-1}
 \lefttop{\vecK'}W_{\vecw(\vecK')}
 \IntC^{\bullet}(\nbigv,W,\vecf)
 \lrarr
 \lefttop{K_p}W_{|K_p|-1}
 \IntC^{\bullet}(\nbigv,W,\vecf)
\]
is a quasi-isomorphism.
Because
\[
  \lefttop{K_p}W_{|K_p|-1}
  \IntC^{\bullet}(\nbigv,W,\vecf)
  \lrarr
  \IntC^{\bullet}(\nbigv,W,\vecf)
\]
is a quasi-isomorphism,
we obtain that
\[
 \lefttop{\vecK}W_{\vecw(\vecK)}
 \IntC^{\bullet}(\nbigv,W,\vecf)
 \lrarr
 \IntC^{\bullet}(\nbigv,W,\vecf)
\]
is a quasi-isomorphism.
\hfill\qed

\vspace{.1in}
We have the following refinement of
Proposition \ref{prop;22.2.7.3}.

\begin{prop}
\label{prop;22.2.7.2}
If $\Lambda_0\cap K_p=\emptyset$
and $K_p\subset \Lambda_1$,
the following natural morphism is a quasi-isomorphism:
\begin{equation}
\label{eq;22.2.7.35}
 \lefttop{\vecK}W_{\vecw(\vecK)}
 \IntC^{\bullet}(\nbigv,W,\vecf_{\Lambda_1};\ast\Lambda_0)
 \lrarr
 \IntC^{\bullet}(\nbigv,W,\vecf_{\Lambda_1};\ast\Lambda_0).
\end{equation}
\end{prop}
\pf
Let $J\subset \Lambda_1\setminus K_p$.
By using the twistor nilpotent orbit
with respect to $\vecf_{\Lambda\setminus K_p}$,
and by using Proposition \ref{prop;22.2.7.3},
we obtain that the following is a quasi-isomorphism
(see the proof of Lemma \ref{lem;22.3.24.101}):
\[
 \lefttop{\vecK} W_{\vecw(\vecK)}
 \IntC^{\bullet}\Bigl(
 \Image f_J^{(\Lambda_0)}\otimes\Tate(-|J|+|J\setminus\Lambda_0|),
 W,\vecf_{K_p}\Bigr)
 \lrarr
 \IntC^{\bullet}\Bigl(
 \Image f_J^{(\Lambda_0)}\otimes\Tate(-|J|+|J\setminus\Lambda_0|),
  W,\vecf_{K_p}\Bigr).
\]
Then, by using the argument in the proof of Proposition \ref{prop;22.2.7.3},
we obtain that (\ref{eq;22.2.7.35}) is a quasi-isomorphism.
\hfill\qed

\subsection{Quasi-isomorphisms for some complexes of mixed twistor structures}

\subsubsection{Preliminary}
\label{subsection;22.2.9.1}

Let $\nbigsbar(\Lambda)$ be as in \S\ref{subsection;22.2.7.20}.
We regard $\nbigsbar(\Lambda)$ be a category
as in \S\ref{subsection;22.2.7.21}.
Let $\MTS$ denote the category of mixed twistor structures.
Let $\ttC(\MTS)$ denote the category of
bounded complexes of mixed twistor structures.

Let $\nbigf^{\bullet}$
be a functor $\nbigsbar(\Lambda)\to\ttC(\MTS)$.
For any $\vecJ=(J_1,\ldots,J_{m-1},\Lambda)\in \nbigsbar(\Lambda)$,
let $\cnum(\vecJ)^{\lor}$ denote the subspace of
$\cnum\langle 2^{\Lambda}\setminus\{\emptyset,\Lambda\}\rangle^{\lor}$
generated by
$\ttv_{J_1}^{\lor}\wedge\cdots\wedge \ttv_{J_{m-1}}^{\lor}$.
For $k\leq 0$, we set
\[
 \nbigc^k(\nbigf^{\bullet})=
 \bigoplus_{|\vecJ|=-k+1}
 \nbigf^{\bullet}(\vecJ)\otimes\cnum(\vecJ)^{\lor}.
\]
For $k>0$, we set $\nbigc^k(\nbigf^{\bullet})=0$.
We define
$d:\nbigc^k(\nbigf^{\bullet})\lrarr \nbigc^{k+1}(\nbigf^{\bullet})$
as in \S\ref{subsection;22.2.7.21}.
Thus, we obtain a double complex
$\nbigc^{\bullet}(\nbigf^{\bullet})$ in $\MTS$,
and the associated complex
$\Tot\nbigc^{\bullet}(\nbigf^{\bullet})
\in \ttC(\MTS)$.

\subsubsection{Some complexes associated with
a polarizable mixed twistor structure}

Let $(\nbigv,W,\vecf)$ be a polarizable mixed twistor structure
of weight $w$.
Let $\nbigf^{\bullet}_{10}$ be a functor
$\nbigsbar(\Lambda)\to \ttC(\MTS)$
defined as follows:
\[
 \nbigf^{\bullet}_{10}(\vecJ)
 =\lefttop{\vecJ}W_{\vecw(\vecJ)}
 \IntC^{\bullet}(\nbigv,W,\vecf;\ast\Lambda).
\]
For
$\vecJ=(J_1\subset\cdots \subset J_{m-1}\subset\Lambda)
\in\nbigsbar(\Lambda)$,
we set
\[
 \nbigf^{\bullet}_{11}(\vecJ)
 =\lefttop{\vecJ}W_{\vecw(\vecJ)}
 \IntC^{\bullet}(\nbigv,W,\vecf;\ast(\Lambda\setminus J_{m-1})).
\]
Thus, we obtain a functor
$\nbigf^{\bullet}_{11}:
\nbigsbar(\Lambda)\to \ttC(\MTS)$.
We also obtain the functor
$\nbigf^{\bullet}_{12}:\nbigsbar(\Lambda)\to\ttC(\MTS)$ as follows:
\[
 \nbigf^{\bullet}_{12}(\vecJ)
 =\lefttop{\vecJ}
 W_{\vecw(\vecJ)}\IntC^{\bullet}(\nbigv,W,\vecf).
\]
We have the following constant functor
$\nbigf^{\bullet}_{13}:\nbigsbar(\Lambda)\to\ttC(\MTS)$:
\[
 \nbigf^{\bullet}_{13}(\vecJ)=\IntC^{\bullet}(\nbigv,W,\vecf).
\]
There exists the following natural transformations:
\[
 \nbigf^{\bullet}_{10}
 \llarr
 \nbigf^{\bullet}_{11}
 \llarr
 \nbigf^{\bullet}_{12}
 \lrarr
 \nbigf^{\bullet}_{13}.
\]

The following theorem is suggested in the last paragraph of 
\cite{Kashiwara-Kawai-Hodge-holonomic}.
\begin{thm}
\label{thm;22.2.8.7}
The following induced morphisms are quasi-isomorphisms:
\begin{equation}
 \begin{CD}
\Tot\nbigc^{\bullet}(\nbigf^{\bullet}_{10})
  @<{g_1}<<
\Tot\nbigc^{\bullet}(\nbigf^{\bullet}_{11})
  @<{g_2}<<
\Tot\nbigc^{\bullet}(\nbigf^{\bullet}_{12})
  @>{g_3}>>
\Tot\nbigc^{\bullet}(\nbigf^{\bullet}_{13}).
 \end{CD}
\end{equation}
\end{thm}

By applying Proposition \ref{prop;22.2.7.2} with $\Lambda_0=\emptyset$,
we obtain that
$g_3$ is a quasi-isomorphism.
We shall prove the claim for $g_2$
in \S\ref{subsection;22.2.8.8}
and the claim for $g_1$ in
\S\ref{subsection;22.2.8.9}.

\subsubsection{Morphism $g_2$}
\label{subsection;22.2.8.8}

For $\vecJ=(J_1,\ldots,J_{m-1},\Lambda)\in\nbigsbar(\Lambda)$,
we set 
\[
 \nbigf_{11}^{\prime\bullet}(\vecJ)
 =\lefttop{\Lambda}W_{|\Lambda|-1}
 \IntC^{\bullet}(\nbigv,W,\vecf;\ast (\Lambda\setminus J_{m-1})),
 \quad
 \nbigf_{12}^{\prime\bullet}(\vecJ)
 =\lefttop{\Lambda}W_{|\Lambda|-1}\IntC^{\bullet}(\nbigv,W,\vecf).
\]
Thus, we obtain functors
$\nbigf_{11}^{\prime\bullet}$
and
$\nbigf_{12}^{\prime\bullet}$
from $\nbigsbar(\Lambda)$
to $\ttC(\MTS)$.
There exists following commutative diagram of
natural transforms:
\[
 \begin{CD}
  \nbigf_{11}^{\bullet}
 @<<<
  \nbigf_{12}^{\bullet}
  \\
  @VVV @VVV \\
 \nbigf_{11}^{\prime\bullet}
 @<<<
 \nbigf_{12}^{\prime\bullet}.
 \end{CD}
\]
Because
$\nbigf^{\bullet}_{i}(\vecJ)
\to
\nbigf^{\prime\bullet}_i(\vecJ)$ $(i=11,12)$
are quasi-isomorphisms
due to Proposition \ref{prop;22.2.7.2},
it is enough to prove that
the following morphism is a quasi-isomorphism:
\[
 g_2':
 \Tot\nbigc^{\bullet}(\nbigf^{\prime\bullet}_{12})
 \lrarr
 \Tot\nbigc^{\bullet}(\nbigf^{\prime\bullet}_{11}).
\]
Choose any point $P\in\proj^1$.
Let $i_P:\{P\}\lrarr\proj^1$ denote the inclusion.
Because $g_2'$ is a complex of morphism of mixed twistor structures,
it is enough to prove that
the induced morphism
\[
 i_P^{\ast}(g_2'):
 i_P^{\ast}\Tot\nbigc^{\bullet}(\nbigf^{\prime\bullet}_{12})
 \lrarr
 i_P^{\ast}\Tot\nbigc^{\bullet}(\nbigf^{\prime\bullet}_{11})
\]
is a quasi-isomorphism (Lemma \ref{lem;22.4.25.1}).

Let $\ttC(\cnum)$ denote the category of
bounded $\cnum$-complexes.
We define the functors
$\nbigf^{\prime\bullet}_{11,P}$
and
$\nbigf^{\prime\bullet}_{12,P}$
from
$\nbigsbar(\Lambda)$ to $\ttC(\cnum)$
by
\[
 \nbigf_{11,P}^{\prime\bullet}(\vecJ)
=i_P^{\ast}\lefttop{\Lambda}W_{|\Lambda|-1}
 \IntC^{\bullet}(\nbigv,W,\vecf;\ast (\Lambda\setminus J_{m-1})),
 \quad
 \nbigf_{12,P}^{\prime\bullet}(\vecJ)
 =i_P^{\ast}\lefttop{\Lambda}W_{|\Lambda|-1}\IntC^{\bullet}(\nbigv,W,\vecf).
\]
We obtain complexes 
$\ttF_{i,P}^{\prime\bullet}$ $(i=11,12)$
in $\ttC^{\wc}(\ttXbar(\Lambda)_{\geq 0})$.
By applying the construction in \S\ref{subsection;22.2.7.21}
to the functor $\nbigf^{\prime\bullet}_{i,P}$,
we obtain the complexes
$\Tot\nbigc^{\bullet}\bigl(
 \nbigf^{\prime\bullet}_{i,P}
 \bigr)$,
which are equal to
$i_P^{\ast}\Tot\nbigc^{\bullet}\bigl(
 \nbigf^{\prime\bullet}_{i}
 \bigr)$
by the construction.
Moreover, 
the natural morphism
\[
  \Tot\nbigc^{\bullet}\bigl(
 \nbigf^{\prime\bullet}_{12,P}
 \bigr)
 \lrarr
   \Tot\nbigc^{\bullet}\bigl(
 \nbigf^{\prime\bullet}_{11,P}
 \bigr)
\]
is identified with
$R\pi_{\ttXbar\ast}(\ttF_{12,P}^{\prime\bullet})
\lrarr
 R\pi_{\ttXbar\ast}(\ttF_{11,P}^{\prime\bullet})$.

Recall that
$\nbigsbar_{\ttY}(\Lambda):=2^{\Lambda}\setminus\{\Lambda\}$,
which we regard as a category as in \S\ref{subsection;22.2.7.21}.
We define functors
$\nbigf_{\ttY,i,P}^{\prime\bullet}$ $(i=11,12)$
from $\nbigsbar_{\ttY}(\Lambda)$
to $\ttC(\cnum)$ by 
\[
 \nbigf_{\ttY,11,P}^{\prime\bullet}(J)
=i_P^{\ast}\lefttop{\Lambda}W_{|\Lambda|-1}
 \IntC^{\bullet}(\nbigv,W,\vecf;\ast (\Lambda\setminus J)),
 \quad
 \nbigf_{\ttY,12,P}^{\prime\bullet}(J)
 =i_P^{\ast}\lefttop{\Lambda}W_{|\Lambda|-1}\IntC^{\bullet}(\nbigv,W,\vecf).
\]
We obtain
$\ttF_{\ttY,i,P}^{\prime\bullet}\in
\ttC^{\wc}(\ttYbar(\Lambda)_{\geq 0};\cnum)$
$(i=11,12)$.
By the construction, there exists the following commutative diagram
in $\ttC^{\wc}(\ttXbar(\Lambda)_{\geq 0};\cnum)$:
\[
 \begin{CD}
  \pi_{\ttYbar,\ttXbar}^{-1}
  \ttF_{\ttY,12,P}^{\prime\bullet}
  @>>>
   \pi_{\ttYbar,\ttXbar}^{-1}
  \ttF_{\ttY,11,P}^{\prime\bullet}
  \\
  @V{=}VV @V=VV \\
  \ttF_{12,P}^{\prime\bullet}
  @>>>
  \ttF_{11,P}^{\prime\bullet}.
 \end{CD}
\]
Therefore,
we obtain the following commutative diagram
in the derived category of $\cnum$-complexes:
\[
 \begin{CD}
  R\pi_{\ttYbar\ast}
  \ttF_{\ttY,12,P}^{\prime\bullet}
  @>>>
   R\pi_{\ttYbar\ast}
  \ttF_{\ttY,11,P}^{\prime\bullet}
  \\
  @V{\simeq}VV @V{\simeq}VV \\
  R\pi_{\ttXbar\ast}
  \ttF_{12,P}^{\prime\bullet}
  @>>>
  R\pi_{\ttXbar\ast}
  \ttF_{11,P}^{\prime\bullet}.
 \end{CD}
\]

By applying the construction in \S\ref{subsection;22.2.7.21}
to the functor $\nbigf^{\prime\bullet}_{\ttY,i,P}$,
we obtain the complexes
$\Tot\nbigc^{\bullet}\bigl(
 \nbigf^{\prime\bullet}_{\ttY,i,P}
 \bigr)$.
There exists the following natural commutative diagram
\[
  \begin{CD}
   R\pi_{\ttYbar\ast}
  \ttF_{\ttY,12,P}^{\prime\bullet}
  @>>>
   R\pi_{\ttYbar\ast}
   \ttF_{\ttY,11,P}^{\prime\bullet}\\
   @V{\simeq}VV @V{\simeq}VV \\
\Tot\nbigc^{\bullet}\bigl(
 \nbigf^{\prime\bullet}_{\ttY,12,P}
   \bigr)
   @>{a}>>
\Tot\nbigc^{\bullet}\bigl(
 \nbigf^{\prime\bullet}_{\ttY,11,P}
   \bigr).
   \end{CD}
\]
Moreover, the morphism $a$ is naturally identified with
\[
 i_P^{\ast}
 \lefttop{\Lambda}W_{|\Lambda|-1}
 \IntC^{\bullet}(\nbigv,W,\vecf)
 \lrarr
 i_P^{\ast}
 \lefttop{\Lambda}W_{|\Lambda|-1}
 \IntC^{\bullet}_{\punc}(\nbigv,W,\vecf).
\]
It is a quasi-isomorphism due to Corollary \ref{cor;22.2.6.10}.
Thus, we obtain that $g_2$ is a quasi-isomorphism.

\subsubsection{Filtration and the associated graded complex}

We make a preliminary to study $g_1$.
Let $\nbigf^{\bullet}$ be a functor from $\nbigsbar(\Lambda)$
to $\ttC(\MTS)$.
For $j\geq 0$ and $k\leq 0$,
we set
\[
 G_j\nbigc^k(\nbigf^{\bullet})=
 \bigoplus_{\substack{|\vecJ|=-k+1 \\ |J_{-k}|\leq j}}
 \nbigf^{\bullet}(\vecJ)\otimes\cnum(\vecJ)^{\lor}.
\]
Then,
$\Tot G_j\nbigc^{\bullet}(\nbigf^{\bullet})$
is a subcomplex of $\Tot\nbigc^{\bullet}(\nbigf^{\bullet})$.
The filtration $G$ is increasing.
Let $K\subsetneq\Lambda$.
For $\vecJ=(J_1,\ldots,J_{m-2},K,\Lambda)\in\nbigsbar_{\ttX,K}(\Lambda)$,
let $\cnum(\vecJ,K)^{\lor}$
denote the subspace of
$\cnum\langle 2^K\setminus\{\emptyset,K\}\rangle^{\lor}$
generated by
$\ttv_{J_1}^{\lor}\wedge\cdots\wedge\ttv_{J_{m-2}}^{\lor}$.
We set
\[
 \nbigc^k_K(\nbigf^{\bullet}):=
 \bigoplus_{\substack{\vecJ\in\nbigsbar_{\ttX,K}(\Lambda)\\
  |\vecJ|=-k+1
  }}
  \nbigf^{\bullet}(\vecJ)\otimes
  \cnum(\vecJ,K)^{\lor}.
\]
Then, for $j>0$,
there exists the following natural isomorphism
\[
 \Gr^G_j\nbigc^{\bullet}(\nbigf^{\bullet})
 \simeq
 \bigoplus_{|K|=j}
 \Tot\nbigc^{\bullet}_K(\nbigf^{\bullet}).
\]
We also have
$\Gr^G_0\Tot\nbigc^{\bullet}(\nbigf^{\bullet})
 \simeq
 \nbigf^{\bullet}(\Lambda)$.

Let $K\subsetneq\Lambda$.
We obtain
 $\nbigf^{\bullet}_{|K}:\nbigsbar(K)\to \ttC(\MTS)$
by
$\nbigf^{\bullet}_{|K}(J_1,\ldots,J_{m-1},K)
=\nbigf^{\bullet}(J_1,\ldots,J_{m-1},K,\Lambda)$.
We have
\[
  \Tot\nbigc^{\bullet}_K(\nbigf^{\bullet})
  \simeq
  \Tot\nbigc^{\bullet}(\nbigf^{\bullet}_{|K})[1].
\]

\subsubsection{Morphism $g_1$}
\label{subsection;22.2.8.9}

Let us study $g_1$.
Let $K\subsetneq \Lambda$.
Let $\nbigf_{i,K}^{\bullet}:\nbigsbar(K)\to\ttC(\MTS)$
$(i=20,21,22,23)$
be the functors defined as follows:
\[
 \nbigf^{\bullet}_{20,K}(\vecJ)=
 \lefttop{\vecJ}W_{\vecw(\vecJ)}
 \IntC^{\bullet}(\nbigv,\vecf_{K};\ast K),
 \quad
 \nbigf^{\bullet}_{21,K}(\vecJ)
 =\lefttop{\vecJ}W_{\vecw(\vecJ)}\IntC^{\bullet}(\nbigv,\vecf_K),
\]
\[
 \nbigf^{\bullet}_{22,K}(\vecJ)=
 \lefttop{\Lambda}W_{|\Lambda|-1}
 \lefttop{\vecJ}W_{\vecw(\vecJ)}
 \IntC^{\bullet}(\nbigv,\vecf_{K};\ast K),
 \quad
 \nbigf^{\bullet}_{23,K}(\vecJ)
 = \lefttop{\Lambda}W_{|\Lambda|-1}
 \lefttop{\vecJ}W_{\vecw(\vecJ)}\IntC^{\bullet}(\nbigv,\vecf_K).
\]
By using the hypothesis of the induction,
and by using the twistor nilpotent orbit
with respect to $\vecf_{\Lambda\setminus K}$
(see the proof of Lemma \ref{lem;22.3.24.101}),
we obtain that
\[
  \Tot\nbigc^{\bullet}(
 \nbigf^{\bullet}_{21,K})
\lrarr
  \Tot\nbigc^{\bullet}(
 \nbigf^{\bullet}_{20,K})
\]
is a quasi-isomorphism in $\ttC(\MTS)$.
Note that
\[
 W_{w+|\Lambda|-1}
\Tot\nbigc^{\bullet}(
 \nbigf^{\bullet}_{20,K})
=
 \Tot\nbigc^{\bullet}(
 \nbigf^{\bullet}_{22,K}),
 \quad\quad
  W_{w+|\Lambda|-1}
\Tot\nbigc^{\bullet}(
 \nbigf^{\bullet}_{21,K})
=\Tot\nbigc^{\bullet}(
 \nbigf^{\bullet}_{23,K}),
\]
where $W$ denotes the weight filtration of
mixed twistor structures.
Hence, we obtain that
\begin{equation}
\label{eq;22.2.8.3}
\Tot\nbigc^{\bullet}(
 \nbigf^{\bullet}_{23,K})
\lrarr
  \Tot\nbigc^{\bullet}(
 \nbigf^{\bullet}_{22,K})
\end{equation}
is a quasi-isomorphism.

\vspace{.1in}
Let us study
$\nbigf^{\bullet}_{i|K}$ $(i=10,11)$.
For $\vecJ=(J_1,\ldots,J_{p-1},K)\in\nbigsbar(K)$,
we have
\[
 \nbigf^{\bullet}_{10|K}(\vecJ)
 =\lefttop{\Lambda}W_{|\Lambda|-1}
 \lefttop{\vecJ}W_{\vecw(\vecJ)}
 \IntC^{\bullet}(\nbigv,W,\vecf;\ast\Lambda)
\]
\[
 \nbigf^{\bullet}_{11|K}(\vecJ)
 =\lefttop{\Lambda}W_{|\Lambda|-1}
 \lefttop{\vecJ}W_{\vecw(\vecJ)}
 \IntC^{\bullet}(\nbigv,W,\vecf;\ast(\Lambda\setminus K)).
\]
We consider the filtration $F_K$ defined by
(\ref{eq;22.2.7.30}) with $L=K$.
We note that if $J\cap K=\emptyset$,
both $f^{(\Lambda)}_J$
and $f^{(\Lambda\setminus K)}_J$ are the identity.
Hence,
each direct summand
in (\ref{eq;22.2.8.2})
for $\Gr_{F_K}^j \nbigf^{\bullet}_{10|K}(\vecJ)$
is the shift of
$\Tot\nbigc^{\bullet}(
 \nbigf^{\bullet}_{22,K})$,
and  each direct summand
in (\ref{eq;22.2.8.2})
for $\Gr_{F_K}^j \nbigf^{\bullet}_{11|K}(\vecJ)$
is the shift of
$\Tot\nbigc^{\bullet}(
 \nbigf^{\bullet}_{23,K})$,
 and the morphism
\begin{equation}
\label{eq;22.2.8.4}
\Gr_{F_K}^j \nbigf^{\bullet}_{11|K}(\vecJ)
\lrarr
\Gr_{F_K}^j \nbigf^{\bullet}_{10|K}(\vecJ)
\end{equation}
is the direct sum of the shifts of (\ref{eq;22.2.8.3}).
Hence, (\ref{eq;22.2.8.4})
is a quasi-isomorphism.
We obtain that
\[
 \nbigf^{\bullet}_{11|K}(\vecJ)
\lrarr
 \nbigf^{\bullet}_{10|K}(\vecJ)
\]
is a quasi-isomorphism for any $K\subsetneq \Lambda$.
Hence, we obtain that 
\begin{equation}
\label{eq;22.2.8.6}
 \Gr^G_j
\Tot\nbigc^{\bullet}(
 \nbigf^{\bullet}_{11})
 \lrarr
  \Gr^G_j
\Tot\nbigc^{\bullet}(
 \nbigf^{\bullet}_{10})
\end{equation}
are quasi-isomorphisms for any $j>0$.
It is easy to check that
(\ref{eq;22.2.8.6}) is an isomorphism if $j=0$.
Thus, we obtain that $g_1$ is a quasi-isomorphism,
and the proof of Theorem \ref{thm;22.2.8.7}
is completed.
\hfill\qed

\subsection{Refinements and some consequences}

\subsubsection{Complexes
associated with a regular admissible polarizable mixed twistor structure}
\label{subsection;22.2.7.13}

Let $X$ be a complex manifold with a simple normal crossing
hypersurface $H$.
Let $\MTS^{\adm}_{\reg}(X,H)$ denote the category of
regular admissible mixed twistor structures on $(X,H)$.
(See \S\ref{subsection;22.3.26.1}.)
Let $\ttC(\MTS^{\adm}_{\reg}(X,H))$
denote the category of complexes in
$\MTS^{\adm}_{\reg}(X,H)$.

Let $(\nbigt,W)\in\MTS^{\adm}_{\reg}(X,H)$.
Let $f_j:(\nbigt,W)\to (\nbigt,W)\otimes\newTate(-1)$ $(j\in\Lambda)$
be commuting tuple of morphisms
such that $(\nbigt,W,\vecf)_{|X\setminus H}$
is a $(w,\Lambda)$-polarizable mixed twistor structure
on $X\setminus H$.
(See \S\ref{subsection;22.3.24.30} and \S\ref{subsection;22.3.27.1}.)
For any non-empty $J\subset\Lambda$,
we set $f_J:=\prod_{j\in J}f_j$.
If $J=\emptyset$, $f_{\emptyset}:=\id$.
We obtain
$(\Image f_J,W,\vecf)
\in\MTS^{\adm}_{\reg}(X,H)$.

Let $\Lambda_0,\Lambda_1\subset\Lambda$.
We set $\vecf_{\Lambda_1}=(f_j\,|\,j\in\Lambda_1)$.
We set
$f^{(\Lambda_0)}_J:=f_{J\setminus\Lambda_0}$.
For $-|\Lambda_1|\leq k\leq 0$,
we set
\[
 \IntC^k(\nbigt,W,\vecf_{\Lambda_1};\ast\Lambda_0):=
 \bigoplus_{\substack{J\subset\Lambda_1\\ |J|=|\Lambda_1|+k }}
 (\Image f^{(\Lambda_0)}_J,W)
 \otimes\newTate(-|J|+|J\setminus\Lambda_0|)
 \otimes\cnum(J).
\]
For $k<-|\Lambda_1|$ or $k>0$,
we set $\IntC^k(\nbigt,W,\vecf_{\Lambda_1};\ast\Lambda_0)=0$.
The differential is given as in \S\ref{subsection;22.2.6.5}.
Thus, we obtain
\[
\IntC^{\bullet}(\nbigt,W,\vecf_{\Lambda_1};\ast\Lambda_0)
\in\ttC(\MTS^{\adm}_{\reg}(X,H)).
\]
If $\Lambda_0=\emptyset$,
it is denoted by 
$\IntC^{\bullet}(\nbigt,W,\vecf_{\Lambda_1})$.

As in \S\ref{subsection;22.2.7.4},
we obtain the following complex in $\MTS^{\adm}_{\reg}(X,H)$:
\[
 \IntC^{\bullet}_{\punc}(\nbigt,W,\vecf):=
 \Tot\Ccheck^{\bullet}
 \IntC^{\bullet}(\nbigt,W,\vecf;\ast\star).
\]

By Corollary \ref{cor;22.2.6.10},
we obtain the following.
\begin{prop}
\label{prop;22.2.7.7}
The following morphisms are quasi-isomorphisms
in $\ttC(\MTS^{\adm}_{\reg}(X,H))$:
 \[
 \IntC^{\bullet}
 (\nbigt,W,\vecf) 
 \llarr
 W_{w+|\Lambda|-1}\IntC^{\bullet}
 (\nbigt,W,\vecf) 
\lrarr
 W_{w+|\Lambda|-1}\IntC_{\punc}^{\bullet}
 (\nbigt,W,\vecf).
\]
\hfill\qed
\end{prop}

Let $\Lambda_0,\Lambda_1\subset\Lambda$ be subsets.
Let $\vecK=(K_1,\ldots,K_p)\in\nbigs(\Lambda)$.
Let $\vecm\in\seisuu^{p}$.
We set
\[
 \lefttop{\vecK}W_{\vecm}
 \IntC^k(\nbigt,W,\vecf_{\Lambda_1};\ast\Lambda_0)
 =\\
 \bigoplus_{\substack{J\subset\Lambda_1 \\ |J|=|\Lambda_1|+k}}
 \bigcap_{i=1}^p
 \lefttop{K_i}
 W_{m_i-|K_i\cap J|
 -|K_i\cap\Lambda_0\cap J|}
 \bigl(
 \Image f^{(\Lambda_0)}_J,W
 \bigr)
  \otimes\newTate(-|J|+|J\setminus\Lambda_0|)
 \otimes\cnum(J).
 \]
Thus, we obtain a complex
$\lefttop{\vecK}W_{\vecm}
\IntC^{\bullet}(\nbigt,W,\vecf_{\Lambda_1};\ast\Lambda_0)$
of mixed twistor structures.
If $\Lambda_0=\emptyset$,
we set
\[
\lefttop{\vecK}W_{\vecm}
\IntC^{\bullet}(\nbigt,W,\vecf_{\Lambda_1}):=
\lefttop{\vecK}W_{\vecm}
\IntC^{\bullet}(\nbigt,W,\vecf_{\Lambda_1};\ast\emptyset).
\]

We obtain the following proposition from
Proposition \ref{prop;22.2.7.2}.
\begin{prop}
\label{prop;22.2.7.10}
If $\Lambda_0\cap K_p=\emptyset$
and $K_p\subset\Lambda_1$,
the following natural morphism is a quasi-isomorphism:
\[
 \lefttop{\vecK}W_{\vecw(\vecK)}
 \IntC^{\bullet}(\nbigt,W,\vecf;\ast\Lambda_0)
 \lrarr
 \IntC^{\bullet}(\nbigt,W,\vecf;\ast\Lambda_0).
\] 
\hfill\qed
\end{prop}

Let $\nbigf^{\bullet}_{10}(\nbigt,W,\vecf)$,
$\nbigf^{\bullet}_{12}(\nbigt,W,\vecf)$
and $\nbigf^{\bullet_{13}}(\nbigt,W,\vecf)$
be functors from $\nbigsbar(\Lambda)$
to $\ttC(\MTS^{\adm}_{\reg}(X,H))$
defined as follows:
\[
 \nbigf^{\bullet}_{10}(\nbigt,W,\vecf)(\vecJ)
 =\lefttop{\vecJ}W_{\vecw(\vecJ)}
 \IntC^{\bullet}(\nbigt,W,\vecf;\ast\Lambda),
\]
\[
  \nbigf^{\bullet}_{12}(\nbigt,W,\vecf)(\vecJ)
 =\lefttop{\vecJ}W_{\vecw(\vecJ)}
 \IntC^{\bullet}(\nbigt,W,\vecf),
\]
\[
 \nbigf^{\bullet}_{13}(\nbigt,W,\vecf)(\vecJ)
=\IntC^{\bullet}(\nbigt,W,\vecf).
\]
There exist the following natural transforms:
\[
 \nbigf^{\bullet}_{10}(\nbigt,W,\vecf)
 \llarr
 \nbigf^{\bullet}_{12}(\nbigt,W,\vecf)
 \lrarr
 \nbigf^{\bullet}_{13}(\nbigt,W,\vecf).
\]
In general, for a functor
$\nbigf^{\bullet}:
\nbigsbar(\Lambda)\to
\ttC(\MTS^{\adm}_{\reg}(X,H))$,
we construct
$\Tot\nbigc^{\bullet}(\nbigf^{\bullet})\in
\ttC(\MTS^{\adm}_{\reg}(X,H))$
as in \S\ref{subsection;22.2.9.1}.
We obtain the following theorem from Theorem \ref{thm;22.2.8.7}.
\begin{thm}
\label{thm;22.2.9.2}
The following morphisms are quasi-isomorphisms:
\[
\Tot\nbigc^{\bullet}(\nbigf^{\bullet}_{10}(\nbigt,W,\vecf))
 \llarr
\Tot\nbigc^{\bullet}(\nbigf^{\bullet}_{12}(\nbigt,W,\vecf))
 \lrarr
\Tot\nbigc^{\bullet}(\nbigf^{\bullet}_{13}(\nbigt,W,\vecf)).
\] 
\hfill\qed
\end{thm}
 
\subsubsection{Complexes of the underlying $\nbigr_X$-modules}
\label{subsection;22.2.7.14}

We set $\nbigx=\cnum\times X$.
Let $p_{\lambda}:\nbigx\lrarr X$ denote the projection.
Let $\nbigd_{\nbigx}$ denote the sheaf of
holomorphic linear differential operators on $\nbigd_{\nbigx}$.
Recall that $\nbigr_X\subset\nbigd_{\cnum\times X}$ denotes
the sheaf of subalgebras generated by
$\lambda p_{\lambda}^{\ast}\Theta_X$
over $\nbigo_{\nbigx}$.
Let $\ttC(\nbigr_X)$ denote the category of bounded
complexes of $\nbigr_X$-modules.

Let $\nbigv$ be an underlying $\nbigr_X$-modules of $\nbigt$.
It is equipped with a commuting tuple of morphisms
$f_j:\nbigv\to \lambda^{-1}\nbigv$ $(j\in\Lambda)$.
We set $f_J:=\prod_{j\in J}f_j$.

Let $\Lambda_0\subset\Lambda$.
We set $f^{(\Lambda_0)}_J:=f_{J\setminus\Lambda_0}$.
For $-|\Lambda|\leq k\leq 0$,
we set
\[
 \IntC^k(\nbigv,\vecf;\ast\Lambda_0):=
 \bigoplus_{\substack{J\subset|\Lambda|\\ |J|=|\Lambda|+k }}
 \lambda^{-|J|+|J\setminus\Lambda_0|}
 \Image f^{(\Lambda_0)}_J
 \otimes\cnum(J).
\]
For $k<-|\Lambda|$ or $k>0$,
we set $\IntC^k(\nbigv,\vecf)=0$.
The differential is given as in \S\ref{subsection;22.2.6.5}.
Thus, we obtain 
$\IntC^{\bullet}(\nbigv,\vecf;\ast\Lambda_0)
\in\ttC(\nbigr_X)$.
If $\Lambda_0=\emptyset$,
it is denoted by 
$\IntC^{\bullet}(\nbigv,\vecf)$.

Let $\vecK=(K_1,\ldots,K_p)\in\nbigs(\Lambda)$.
Let $\vecm\in\seisuu^{p}$.
We set
\[
 \lefttop{\vecK}W_{\vecm}\IntC^k(\nbigv,\vecf;\ast\Lambda_0)
 =\\
 \bigoplus_{\substack{J\subset\Lambda \\ |J|=|\Lambda|+k}}
 \bigcap_{i=1}^p
 \lefttop{K_i}
 W_{m_i-|K_i\cap J|
 -|K_i\cap\Lambda_0\cap J|}
 \bigl(
\lambda^{-|J|+|J\setminus\Lambda_0|}
 \Image f^{(\Lambda_0)}_J
 \bigr)
 \otimes\cnum(J).
 \]
Thus, we obtain 
$\lefttop{\vecK}W_{\vecm}
\IntC^{\bullet}(\nbigv,\vecf;\ast\Lambda_0)
\in\ttC(\nbigr_X)$.
If $\Lambda_0=\emptyset$,
we set
\[
\lefttop{\vecK}W_{\vecm}\IntC^{\bullet}(\nbigv,\vecf):=
\lefttop{\vecK}W_{\vecm}\IntC^{\bullet}(\nbigv,\vecf;\ast\emptyset).
\]

We obtain the following proposition from
Proposition \ref{prop;22.2.7.10}.
\begin{prop}
\label{prop;22.2.7.12}
If $\Lambda_0\cap K_p=\emptyset$,
the following natural morphism is a quasi-isomorphism
in $\ttC(\nbigr_X)$:
\[
 \lefttop{\vecK}W_{\vecw(\vecK)}
 \IntC^{\bullet}(\nbigv,\vecf;\ast\Lambda_0)
 \lrarr
 \IntC^{\bullet}(\nbigv,\vecf;\ast\Lambda_0).
\] 
\hfill\qed
\end{prop}

Let $\nbigf^{\bullet}_{10}(\nbigv,\vecf)$,
$\nbigf^{\bullet}_{12}(\nbigv,\vecf)$
and $\nbigf^{\bullet}_{13}(\nbigv,\vecf)$
be functors from $\nbigsbar(\Lambda)$
to $\ttC(\nbigr_X)$
defined as follows:
\[
 \nbigf^{\bullet}_{10}(\nbigv,\vecf)(\vecJ)
 =\lefttop{\vecJ}W_{\vecw(\vecJ)}
 \IntC^{\bullet}(\nbigv,\vecf;\ast\Lambda),
\]
\[
  \nbigf^{\bullet}_{12}(\nbigv,\vecf)(\vecJ)
 =\lefttop{\vecJ}W_{\vecw(\vecJ)}
 \IntC^{\bullet}(\nbigv,\vecf),
\]
\[
 \nbigf^{\bullet}_{13}(\nbigv,\vecf)(\vecJ)
=\IntC^{\bullet}(\nbigv,\vecf).
\]
There exist the following natural transforms:
\[
 \nbigf^{\bullet}_{10}(\nbigv,\vecf)
 \llarr
 \nbigf^{\bullet}_{12}(\nbigv,\vecf)
 \lrarr
 \nbigf^{\bullet}_{13}(\nbigv,\vecf).
\]
In general, for a functor
$\nbigf^{\bullet}:
\nbigsbar(\Lambda)\to\ttC(\nbigr_X)$,
we construct
an $\nbigr_X$-complex $\Tot\nbigc^{\bullet}(\nbigf^{\bullet})$
as in \S\ref{subsection;22.2.9.1}.
We obtain the following theorem from Theorem \ref{thm;22.2.9.2}.
\begin{prop}
\label{prop;22.2.9.3}
The following morphisms
in $\ttC(\nbigr_X)$ are quasi-isomorphisms:
\[
\Tot\nbigc^{\bullet}(\nbigf^{\bullet}_{10}(\nbigv,\vecf))
 \llarr
\Tot\nbigc^{\bullet}(\nbigf^{\bullet}_{12}(\nbigv,\vecf))
 \lrarr
\Tot\nbigc^{\bullet}(\nbigf^{\bullet}_{13}(\nbigv,\vecf)).
\] 
\hfill\qed
\end{prop}

\subsubsection{Regular KMS-structure}

Let $V_H\nbigr_X\subset\nbigr_X$
denote the sheaf of subalgebras
generated by $\lambda p_{\lambda}^{\ast}\Theta_X(\log H)$
over $\nbigo_{\nbigx}$.
Let $\lambda_0\in\cnum$.
Let $\nbigx^{(\lambda_0)}$ denote a neighbourhood of
$\{\lambda_0\}\times X$ in $\cnum_{\lambda}\times X$.
Let $\ttC(V_H\nbigr_{X|\nbigxzero})$
denote the category of $\nbigr_{X|\nbigxzero}$-complexes.

We have the regular KMS-structure $\nbigp^{(\lambda_0)}_{\ast}\nbigv$
which is a filtration of $\nbigv_{|\nbigxzero}$
by $V_H\nbigr_{X|\nbigxzero}$-submodules.
(See \S\ref{subsection;22.3.25.120} for
regular-KMS-structure.)

\begin{lem}
We have
$\nbigp^{(\lambda_0)}_{\veca}(\Image f_{J})
=\lambda^{-1}\nbigp^{(\lambda_0)}_{\veca}(\nbigv)\cap\Image f_J$.
For the filtration $\lefttop{K}W$ on $\nbigpzero_{\veca}(\Image f_J)$,
we have
$\lefttop{K}W_{\ast}(\nbigp^{(\lambda_0)}_{\veca}\Image f_J)
 =\lefttop{K}W_{\ast}(\Image f_J)\cap
 \nbigp^{(\lambda_0)}_{\veca}(\Image f_J)$
for any $K\subset\Lambda$.
We also have
\[
 \bigcap_{i=1}^p
 \lefttop{K_i}W_{m_i}(\nbigp^{(\lambda_0)}_{\veca}\Image f_J)
=\bigcap_{i=1}^p
 \lefttop{K}W_{\ast}(\nbigv)\cap
 \nbigp^{(\lambda_0)}_{\veca}(\Image f_J).
\]
\end{lem}
\pf
If follows from Proposition \ref{prop;22.3.27.2} below.
(Note the inductive construction of
the filtrations $\lefttop{K}W$
as in \cite[(1.2.3)]{k3}.)
\hfill\qed

\vspace{.1in}

Let $\Lambda_0\subset\Lambda$.
For $-|\Lambda|\leq k\leq 0$,
we set
\[
 \IntC^k(\nbigp^{(\lambda_0)}_{\veca}\nbigv,\vecf;\ast\Lambda_0):=
 \bigoplus_{\substack{J\subset|\Lambda|\\ |J|=|\Lambda|+k }}
 \lambda^{-|J|+|J\setminus\Lambda_0|}
 \nbigp^{(\lambda_0)}_{\veca}(\Image f^{(\Lambda_0)}_J)
 \otimes\cnum(J).
\]
For $k<-|\Lambda|$ or $k>0$,
we set $\nbigpzero_{\veca}\IntC^k(\nbigv,\vecf)=0$.
The differential is given as in \S\ref{subsection;22.2.6.5}.
Thus, we obtain
$\IntC^{\bullet}(\nbigpzero_{\veca}\nbigv,\vecf;\ast\Lambda_0)
\in\ttC(V_H\nbigr_{X|\nbigxzero})$,
which is a $V_H\nbigr_{X|\nbigxzero}$-subcomplex of
$\IntC^{\bullet}(\nbigv,\vecf;\ast\Lambda_0)$.
If $\Lambda_0=\emptyset$,
it is denoted by 
$\IntC^{\bullet}(\nbigpzero_{\veca}\nbigv,\vecf)$.

Let $\vecK=(K_1,\ldots,K_p)\in\nbigs(\Lambda)$.
Let $\vecm\in\seisuu^{p}$.
We set
\[
 \lefttop{\vecK}W_{\vecm}
 \IntC^k(\nbigpzero_{\veca}\nbigv,\vecf;\ast\Lambda_0)
 =\\
 \bigoplus_{\substack{J\subset\Lambda \\ |J|=|\Lambda|+k}}
 \bigcap_{i=1}^p
 \lefttop{K_i}
 W_{m_i-|K_i\cap J|
 -|K_i\cap\Lambda_0\cap J|}
 \bigl(
\lambda^{-|J|+|J\setminus\Lambda_0|}
 \nbigpzero_{\veca}(
 \Image f^{(\Lambda_0)}_J)
 \bigr)
 \otimes\cnum(J).
 \]
We obtain
$\lefttop{\vecK}W_{\vecm}
\IntC^{\bullet}(\nbigpzero_{\veca}\nbigv,\vecf;\ast\Lambda_0)
\in\ttC(V_H\nbigr_{X|\nbigxzero})$,
which is an $\nbigr_{X|\nbigxzero}$-subcomplex of
$\IntC^{\bullet}(\nbigpzero_{\veca}\nbigv,\vecf;\ast\Lambda_0)$.
If $\Lambda_0=\emptyset$,
we set
\[
\lefttop{\vecK}W_{\vecm}
\IntC^{\bullet}(\nbigpzero_{\veca}\nbigv,\vecf):=
\lefttop{\vecK}W_{\vecm}
\IntC^{\bullet}(\nbigpzero_{\veca}\nbigv,\vecf;\ast\emptyset).
\]
We obtain the following proposition
from Proposition \ref{prop;22.2.7.12}.
\begin{prop}
\label{prop;22.3.14.12}
If $\Lambda_0\cap K_p=\emptyset$,
the following natural morphism is a quasi-isomorphism:
\[
 \lefttop{\vecK}W_{\vecw(\vecK)}
 \IntC^{\bullet}(\nbigpzero\nbigv,\vecf;\ast\Lambda_0)
 \lrarr
 \IntC^{\bullet}(\nbigpzero_{\veca}\nbigv,\vecf;\ast\Lambda_0).
\] 
\hfill\qed
\end{prop}

Let
$\nbigf^{\bullet}_{10}(\nbigp^{(\lambda_0)}_{\veca}\nbigv,\vecf)$,
$\nbigf^{\bullet}_{12}(\nbigp^{(\lambda_0)}_{\veca}\nbigv,\vecf)$
and
$\nbigf^{\bullet}_{13}(\nbigp^{(\lambda_0)}_{\veca}\nbigv,\vecf)$
be functors from $\nbigsbar(\Lambda)$
to $\ttC(V\nbigr_{X|\nbigx^{(\lambda_0)}})$
defined as follows:
\[
 \nbigf^{\bullet}_{10}( \nbigp^{(\lambda_0)}_{\veca}\nbigv,\vecf)(\vecJ)
 =\lefttop{\vecJ}W_{\vecw(\vecJ)}
 \IntC^{\bullet}( \nbigpzero_{\veca}\nbigv,\vecf;\ast\Lambda),
\]
\[
  \nbigf^{\bullet}_{12}(\nbigpzero_{\veca}\nbigv,\vecf)(\vecJ)
  =\lefttop{\vecJ}W_{\vecw(\vecJ)}
 \IntC^{\bullet}(\nbigpzero_{\veca}\nbigv,\vecf),
\]
\[
 \nbigf^{\bullet}_{13}(\nbigpzero_{\veca}\nbigv,\vecf)(\vecJ)
=\IntC^{\bullet}(\nbigpzero_{\veca}\nbigv,\vecf).
\]
There exist the following natural transforms:
\[
 \nbigf^{\bullet}_{10}(\nbigpzero_{\veca}\nbigv,\vecf)
 \llarr
 \nbigf^{\bullet}_{12}(\nbigpzero_{\veca}\nbigv,\vecf)
 \lrarr
 \nbigf^{\bullet}_{13}(\nbigpzero_{\veca}\nbigv,\vecf).
\]
In general, for a functor $\nbigf^{\bullet}:
\nbigsbar(\Lambda)\to\ttC(V\nbigr_{X|\nbigxzero})$,
we construct
$\Tot\nbigc^{\bullet}(\nbigf^{\bullet})\in
\ttC(V\nbigr_{X|\nbigxzero})$
as in \S\ref{subsection;22.2.9.1}.

\begin{prop}
\label{prop;22.4.25.20}
The following morphisms are quasi-isomorphisms:
\[
\Tot\nbigc^{\bullet}(\nbigf^{\bullet}_{10}(\nbigpzero_{\veca}\nbigv,\vecf))
 \llarr
\Tot\nbigc^{\bullet}(\nbigf^{\bullet}_{12}(\nbigpzero_{\veca}\nbigv,\vecf))
 \lrarr
\Tot\nbigc^{\bullet}(\nbigf^{\bullet}_{13}(\nbigpzero_{\veca}\nbigv,\vecf)).
\] 
\end{prop}
\pf
It follows from Proposition \ref{prop;22.2.9.3}
and Proposition \ref{prop;22.3.27.2}.
\hfill\qed

\subsubsection{Stalks}

Let $P$ be any point of
$\bigcap_{i\in\Lambda}\nbigh^{(\lambda_0)}_i$.
For any complex of sheaves $\nbigg^{\bullet}$ on $\nbigxzero$,
let $\nbigg^{\bullet}_P$ denote the stalk at $P$.
Let $\ttC((V_H\nbigr_X)_{P})$ denote category of
bounded $(V_H\nbigr_X)_{P}$-complexes.

Let
$\nbigf^{\bullet}_{10}(\nbigp^{(\lambda_0)}_{\veca}\nbigv_P,\vecf)$,
$\nbigf^{\bullet}_{12}(\nbigp^{(\lambda_0)}_{\veca}\nbigv_P,\vecf)$
and
$\nbigf^{\bullet}_{13}(\nbigp^{(\lambda_0)}_{\veca}\nbigv_P,\vecf)$
be functors from $\nbigsbar(\Lambda)$
to $\ttC((V_H\nbigr_X)_P)$
defined as follows:
\[
 \nbigf^{\bullet}_{10}(\nbigp^{(\lambda_0)}_{\veca}\nbigv_P,\vecf)(\vecJ)
 =\lefttop{\vecJ}W_{\vecw(\vecJ)}
 \IntC^{\bullet}( \nbigpzero_{\veca}\nbigv,\vecf;\ast\Lambda)_P,
\]
\[
  \nbigf^{\bullet}_{12}(\nbigpzero_{\veca}\nbigv_P,\vecf)(\vecJ)
  =\lefttop{\vecJ}W_{\vecw(\vecJ)}
 \IntC^{\bullet}(\nbigpzero_{\veca}\nbigv,\vecf)_P,
\]
\[
 \nbigf^{\bullet}_{13}(\nbigpzero_{\veca}\nbigv_P,\vecf)(\vecJ)
=\IntC^{\bullet}(\nbigpzero_{\veca}\nbigv,\vecf)_P.
\]
There exist the following natural transforms:
\[
 \nbigf^{\bullet}_{10}(\nbigpzero_{\veca}\nbigv_P,\vecf)
 \llarr
 \nbigf^{\bullet}_{12}(\nbigpzero_{\veca}\nbigv_P,\vecf)
 \lrarr
 \nbigf^{\bullet}_{13}(\nbigpzero_{\veca}\nbigv_P,\vecf).
\]
In general,
for a functor $\nbigf^{\bullet}:
\nbigsbar(\Lambda)\to \ttC((V_H\nbigr_X)_P)$we construct
$\Tot\nbigc^{\bullet}(\nbigf^{\bullet})
\in\ttC((V_H\nbigr_X)_P)$
as in \S\ref{subsection;22.2.9.1}.
\begin{cor}
\label{cor;22.4.25.21}
The following morphisms
in $\ttC((V_H\nbigr_X)_P)$
are quasi-isomorphisms:
\[
\Tot\nbigc^{\bullet}(\nbigf^{\bullet}_{10}(\nbigpzero_{\veca}\nbigv_P,\vecf))
 \llarr
\Tot\nbigc^{\bullet}(\nbigf^{\bullet}_{12}(\nbigpzero_{\veca}\nbigv_P,\vecf))
 \lrarr
\Tot\nbigc^{\bullet}(\nbigf^{\bullet}_{13}(\nbigpzero_{\veca}\nbigv_P,\vecf)).
\] 
\end{cor}
\pf
It follows from Proposition \ref{prop;22.4.25.20}.
\hfill\qed

\subsubsection{The associated weakly constructible sheaves}
\label{subsection;22.2.9.10}

We obtain complexes
$\ttF^{\bullet}_i(\nbigpzero_{\veca}\nbigv_P,\vecf)$
$(i=10,12,13)$
in $\ttC^{\wc}(\ttX(\Lambda)_{\geq 0};(V_H\nbigr_X)_P)$.
There exist the following natural morphisms:
\[
 \ttF^{\bullet}_{10}(\nbigpzero_{\veca}\nbigv_P,\vecf)
 \llarr
 \ttF^{\bullet}_{12}(\nbigpzero_{\veca}\nbigv_P,\vecf)
 \lrarr
 \ttF^{\bullet}_{13}(\nbigpzero_{\veca}\nbigv_P,\vecf).
\]
We obtain the following
from Corollary \ref{cor;22.4.25.21}
and (\ref{eq;22.4.25.22}).

\begin{prop}
\label{prop;22.2.17.50}
The induced morphisms
\[
 R\pi_{\ttXbar(\Lambda)\ast}
 \ttF^{\bullet}_{10}(\nbigpzero_{\veca}\nbigv_P,\vecf)
 \llarr
  R\pi_{\ttXbar(\Lambda)\ast}
 \ttF^{\bullet}_{12}(\nbigpzero_{\veca}\nbigv_P,\vecf)
 \lrarr
  R\pi_{\ttXbar(\Lambda)\ast}
 \ttF^{\bullet}_{13}(\nbigpzero_{\veca}\nbigv_P,\vecf)
\] 
in the derived category of $(V_H\nbigr_X)_P$-complexes
are isomorphisms
\hfill\qed
\end{prop}

\subsubsection{Complement}

Let $\Gamma$ be a finite subset.
We set
$\Lambdatilde:=\Lambda\sqcup\Gamma$.
For $\vecJ=(J_1,\ldots,J_{m-1},J_m)\in\nbigsbar(\Lambdatilde)$
with $J_m=\Lambdatilde$,
there exists $0\leq i(0)\leq m$ such that
$J_{i(0)}\subset\Lambda$
and $J_{i(0)+1}\not\subset\Lambda$.
We set
$\vecJ\cap\Lambda:=
(J_{1},\ldots,J_{i(0)})\in\nbigs(\Lambda)$.

We define the functors
$\nbigf_{10,\Lambdatilde}^{\bullet}(\nbigp_{\veca}\nbigv_P,\vecf)$,
$\nbigf_{12,\Lambdatilde}^{\bullet}(\nbigp_{\veca}\nbigv_P,\vecf)$,
and
$\nbigf_{13,\Lambdatilde}^{\bullet}(\nbigp_{\veca}\nbigv_P,\vecf)$
from $\nbigsbar(\Lambdatilde)$
to $\ttC((V_H\nbigr_X)_P)$
as follows:
\[
 \nbigf_{10,\Lambdatilde}^{\bullet}(\nbigpzero_{\veca}\nbigv_P,\vecf)
 (\vecJ)
 =\lefttop{\vecJ\cap\Lambda}
W_{\vecw(\vecJ\cap\Lambda)}
 \IntC^{\bullet}(\nbigpzero_{\veca}\nbigv,\vecf;\ast\Lambda)_P
\]
\[
 \nbigf_{12,\Lambdatilde}^{\bullet}(\nbigpzero_{\veca}\nbigv_P,\vecf)(\vecJ)
 =\lefttop{\vecJ\cap\Lambda}
W_{\vecw(\vecJ\cap\Lambda)}
 \IntC^{\bullet}(\nbigpzero_{\veca}\nbigv,\vecf)_P
\]
\[
 \nbigf_{13,\Lambdatilde}^{\bullet}(\nbigpzero_{\veca}\nbigv_P,\vecf)(\vecJ)
 =\IntC^{\bullet}(\nbigpzero_{\veca}\nbigv,\vecf)_P.
\]
There exist the following natural transforms:
\[
  \nbigf_{10,\Lambdatilde}^{\bullet}(\nbigpzero_{\veca}\nbigv_P,\vecf)
  \llarr
  \nbigf_{12,\Lambdatilde}^{\bullet}(\nbigpzero_{\veca}\nbigv_P,\vecf)
  \lrarr
  \nbigf_{13,\Lambdatilde}^{\bullet}(\nbigpzero_{\veca}\nbigv_P,\vecf).
\]
We obtain the complexes
$\ttF^{\bullet}_{i,\Lambdatilde}(\nbigpzero_{\veca}\nbigv_P,\vecf)$
in $\ttC^{\wc}(\ttX(\Lambdatilde)_{\geq 0};(V_H\nbigr_X)_P)$
corresponding to
$\nbigf_{i,\Lambdatilde}^{\bullet}(\nbigpzero_{\veca}\nbigv_P,\vecf)$.
There exists the following natural morphisms
\[
 \ttF^{\bullet}_{10,\Lambdatilde}(\nbigpzero_{\veca}\nbigv_P,\vecf)
 \llarr
 \ttF^{\bullet}_{12,\Lambdatilde}(\nbigpzero_{\veca}\nbigv_P,\vecf)
 \lrarr
 \ttF^{\bullet}_{13,\Lambdatilde}(\nbigpzero_{\veca}\nbigv_P,\vecf).
\]
For this construction,
$\ttF^{\bullet}_{i,\Lambda}(\nbigpzero_{\veca}\nbigv_P,\vecf)$
are equal to
$\ttF^{\bullet}_i(\nbigpzero_{\veca}\nbigv_P,\vecf)$
in \S\ref{subsection;22.2.9.10}.

We have
$\overline{\ttO(\Lambda,\Lambdatilde)}_{\geq 0}
\subset
\ttXbar(\Lambdatilde)_{\geq 0}$.
Let $\iota_{\Lambda,\Lambdatilde}$ denote the inclusion.
We also have the natural morphism
$\pi_{\Lambda,\Lambdatilde}:
\overline{\ttO(\Lambda,\Lambdatilde)}_{\geq 0}
\lrarr
 \ttXbar(\Lambda)_{\geq 0}$.
There exist the following natural quasi-isomorphisms
for $i=10,12,13$:
\[
R\pi_{\Lambda,\Lambdatilde \ast}
\pi_{\Lambda,\Lambdatilde}^{-1}
(\ttF^{\bullet}_{i,\Lambda}(\nbigpzero_{\veca}\nbigv_P,\vecf))
\simeq
\ttF^{\bullet}_{i,\Lambda}(\nbigpzero_{\veca}\nbigv_P,\vecf).
\]
By the construction,
there exist the following isomorphisms:
\[
\pi_{\Lambda,\Lambdatilde}^{-1}
(\ttF^{\bullet}_{i,\Lambda}(\nbigpzero_{\veca}\nbigv_P,\vecf))
\simeq
\iota_{\Lambda,\Lambdatilde}^{-1}
(\ttF^{\bullet}_{i,\Lambdatilde}(\nbigpzero_{\veca}\nbigv_P,\vecf)).
\]
Hence, we obtain the following commutative diagram:
\begin{equation}
\label{eq;22.2.9.11}
\begin{CD}
 R\pi_{\ttXbar(\Lambdatilde)\ast}
 \ttF^{\bullet}_{10,\Lambdatilde}(\nbigpzero_{\veca}\nbigv_P,\vecf)
 @<<<
 R\pi_{\ttXbar(\Lambdatilde)\ast}
 \ttF^{\bullet}_{12,\Lambdatilde}(\nbigpzero_{\veca}\nbigv_P,\vecf)
 @>>>
 R\pi_{\ttXbar(\Lambdatilde)\ast}
 \ttF^{\bullet}_{13,\Lambdatilde}(\nbigpzero_{\veca}\nbigv_P,\vecf)
 \\
 @V{\alpha_1}VV @V{\alpha_2}VV @V{\alpha_3}VV \\
 R\pi_{\ttXbar(\Lambda)\ast}
 \ttF^{\bullet}_{10,\Lambda}(\nbigpzero_{\veca}\nbigv_P,\vecf)
 @<{\simeq}<<
 R\pi_{\ttXbar(\Lambda)\ast}
 \ttF^{\bullet}_{12,\Lambda}(\nbigpzero_{\veca}\nbigv_P,\vecf)
 @>{\simeq}>>
 R\pi_{\ttXbar(\Lambda)\ast}
 \ttF^{\bullet}_{13,\Lambda}(\nbigpzero_{\veca}\nbigv_P,\vecf).
\end{CD}
\end{equation}

\begin{prop}
\label{prop;22.2.17.51}
The morphisms $\alpha_i$ in {\rm(\ref{eq;22.2.9.11})}
are quasi-isomorphisms.
\end{prop}
\pf
Let $\prec$ denote the order on $\nbigsbar(\Lambda)$
defined by
$\vecJ_1\prec\vecJ_2$ if $\vecJ_1$ is a refinement of $\vecJ_2$.
Let $S\subset\nbigsbar(\Lambda)$ be a subset
satisfying the following condition.
\begin{description}
 \item[(A)] Let $\vecJ\in S$.
       Then,
       $\{\vecJ'\in\nbigsbar(\Lambda)\,|\,\vecJ'\prec\vecJ\}
       \subset S$.
\end{description}
For any 
$\vecJ\in\nbigsbar(\Lambda)$,
we obtain
$\vecJ^{\Lambdatilde},\vecJ\cdot\Lambdatilde
\in \nbigsbar(\Lambdatilde)$
as in \S\ref{subsection;22.1.22.13}.
We set
\[
 \nbigz(\vecJ,\Lambda):=
 \overline{\ttO(\vecJ\cdot\Lambdatilde)}_{\geq 0}
 \cup
 \bigcup_{\vecJ_1\prec\vecJ}
 \overline{\ttO(\vecJ_1^{\Lambdatilde})}_{\geq 0}.
\]
We also set
\[
 \nbigy(\vecJ,\Lambda):=
 \overline{\ttO(\vecJ^{\Lambdatilde})}_{\geq 0}
 \setminus
 \nbigz(\vecJ,\Lambda).
\]
Let $k_{\vecJ}:\nbigy(\vecJ,\Lambda)\lrarr
 \ttXbar(\Lambdatilde)_{\geq 0}$
denote the inclusion.
By the construction,
$k_{\vecJ}^{-1}
\ttF_{j,\Lambdatilde}^{\bullet}(\nbigpzero_{\veca}\nbigv_P,\vecf)$
are complexes of constant sheaves on
$\nbigy(\vecJ,\Lambda)$.
Because
$\overline{\ttO(\vecJ^{\Lambdatilde})}_{\geq 0}$
and 
$\nbigz(\vecJ,\Lambda)$
are contractible,
we obtain the following vanishing:
\[
 R\pi_{\ttXbar(\Lambdatilde)\ast}
  k_{\vecJ!}k_{\vecJ}^{-1}\ttF_{j,\Lambdatilde}^{\bullet}
  (\nbigp_{\veca}\nbigv_P,\vecf)=0.
\]

We set
\[
 \nbigx(S):=
 \overline{\ttO(\Lambda)}_{\geq 0}
 \cup
 \bigcup_{\vecJ\in S}
 \overline{\ttO(\vecJ^{\Lambdatilde})}_{\geq 0}
 \subset \ttXbar(\Lambdatilde)_{\geq 0}.
\]
Then, $\nbigx(S)$ is a closed subset.
Let $\vecJ_0$ be a maximal element of $S$.
We set $S':=S\setminus\{\vecJ_0\}$.
Then, $S'$ also satisfies the above condition (A).
We have
\[
 \nbigx(S)
=\nbigx(S')
 \cup
 \overline{\ttO(\vecJ^{\Lambda}_0)}_{\geq 0},
 \quad
 \nbigx(S')\cap
 \overline{\ttO(\vecJ^{\Lambdatilde}_0)}_{\geq 0}
=\nbigz(\vecJ_0,\Lambda).
\]
Hence, we obtain the claim of the proposition
by an easy induction.
\hfill\qed

\section{Complexes associated with good mixed twistor $\nbigd$-modules}

\subsection{Multi-$V$-filtration}
\label{subsection;22.2.14.1}

For a complex manifold $X$,
we set $\nbigx=\cnum_{\lambda}\times X$.
Let $p_{\lambda}:\nbigx\lrarr X$ denote the projection.
Let $\nbigd_{\nbigx}$ denote the sheaf of
holomorphic linear differential operators on $\nbigx$.
Recall that $\nbigr_X\subset\nbigd_{\nbigx}$ denotes
the sheaf of subalgebras generated by
$\lambda p_{\lambda}^{\ast}\Theta_X$
over $\nbigo_{\nbigx}$.

Let $\nbigc(X)$ denote the full subcategory of $\nbigr_X$-modules
$\nbigm$ underlying a mixed twistor $\nbigd_X$-module.
(See \S\ref{subsection;22.4.25.30} for a brief explanation
of mixed twistor $\nbigd$-modules.)
For any $\lambda_0\in\cnum$,
let $\nbigxzero$ denote a small neighbourhood of
$\{\lambda_0\}\times X$ in $\nbigx:=\cnum\times X$.

Let $Y$ be a complex manifold.
We consider the case where
$X$ is a neighbourhood of
$\{0\}\times Y$ in $\cnum^{\Lambda}\times Y$.
We set $H_i=\{z_i=0\}\cap X=
\bigl(
\cnum^{\Lambda\setminus\{i\}}\times Y
\bigr)\cap X$
$(i\in\Lambda)$.
We set $H_I=\bigcap_{i\in I}H_i$
and $\del H_I:=H\cap\bigcup_{i\in \Lambda\setminus I}H_i$
for $I\subset\Lambda$.

For any $I\subset J\subset\Lambda$,
let $q_I:\real^J\lrarr\real^I$ denote the projection.
Let $\veciti_I\in\real^J$ denote the element
such that $q_i(\veciti_I)=1$ $(\in I)$
and $q_i(\veciti_I)=0$ $(i\not\in I)$.
For any element $\veca\in\real^J$,
$q_i(\veca)$ $(i\in J)$ is also denoted by $a_i$.

Let $\nbigm\in\nbigc(H_I)$.
It has a $V$-filtrations $\lefttop{i}\Vzero$ along $z_i$
$(i\in \Lambda\setminus I)$.
(See \S\ref{subsection;22.3.28.10}
for $V$-filtrations.)
For any subset $K\subset\Lambda\setminus I$
and $\veca\in\real^K$,
we set
$\lefttop{K}\Vzero_{\veca}(\nbigm):=
\bigcap_{i\in K}\lefttop{i}\Vzero_{a_i}(\nbigm)$.
The $\nbigr_{H_{I\cup\{i\}}}$-modules
$\psitilde_{z_i,u}(\nbigm)$
$(u\in\real\times\cnum\setminus(\seisuu_{\geq 0}\times\{0\}))$
$\psi_{z_i,-\vecdelta}(\nbigm)$
and $\psi_{z_i,0}(\nbigm)$
are denoted by 
$\lefttop{i}\psitilde_{u}(\nbigm)$,
$\lefttop{i}\psi_{-\vecdelta}(\nbigm)$
and $\lefttop{i}\psi_{0}(\nbigm)$,
respectively.
(See \S\ref{subsection;22.3.28.10} for the notation.)
We say that the $\nbigr_{H_I}$-module $\nbigm$
has a multi-$V$-filtration along $\del H_I$
if the following holds.
\begin{itemize}
 \item $\lefttop{i}\psitilde_u(\nbigm)$
       $(u\in\real\times\cnum)$,
       $\lefttop{i}\psi_{-\vecdelta}(\nbigm)$
       and
       $\lefttop{i}\psi_{0}(\nbigm)$
       on $H_{I\cup\{i\}}$
       have multi-$V$-filtrations along
       $\del H_{I\cup \{i\}}$.
 \item
      For any $K\subset\Lambda$, $i\in K$
      and $\veca\in\real^I$, the induced morphism
      $\lefttop{K}\Vzero_{\veca}(\nbigm)
      \lrarr
      \lefttop{K\setminus i}\Vzero_{q_{K\setminus i}(\veca)}
      \lefttop{i}\Gr^{\Vzero}_{q_i(\veca)}(\nbigm)$
      is an epimorphism.
      Here,
      $\lefttop{i}\Gr^{\Vzero}_{q_i(\veca)}(\nbigm)
      =\lefttop{i}\Vzero_{q_i(\veca)}(\nbigm)\big/
       \lefttop{i}\Vzero_{<q_i(\veca)}(\nbigm)$.
\end{itemize}

\begin{prop}
\label{prop;22.3.13.2}
Let $F:\nbigm_1\lrarr\nbigm_2$
be a morphism in $\nbigc(X)$.
Suppose that $\nbigm_i$ have multi-$V$-filtrations
along $H$.
We also assume that
$\Ker F$, $\Image F$ and $\Cok(F)$
are objects in $\nbigc(X)$.
Then, the following holds.
\begin{itemize}
 \item $\Ker(F)$, $\Image(F)$ and $\Cok(F)$
       also have multi-$V$-filtrations
       along $H$.
 \item $F$ is strict with respect to
       $\lefttop{I}V_{\bullet}$ for any $I\subset\Lambda$,
       i.e.,
$ F(\nbigm_1)\cap \lefttop{I}V_{\veca}(\nbigm_2)
=F(\lefttop{I}V_{\veca}\nbigm_1)$.
\end{itemize}
\end{prop}
\pf
We use an induction on $|\Lambda|$.
We also use an induction on $|I|$.
It is known that
$F$ is strict with respect to $\lefttop{i}\Vzero$.
Suppose we have already proved that
$F$ is strict with respect to $\lefttop{I}\Vzero$.
Let $j\in\Lambda\setminus I$.
We consider the induced morphism
$\lefttop{j}\Vzero_b\lefttop{I}\Vzero_{\veca}(\nbigm_1)
\lrarr
\lefttop{j}\Vzero_b\lefttop{I}\Vzero_{\veca}(\nbigm_2)$.
By the assumption of the induction,
the morphism
\begin{equation}
\label{eq;22.2.11.1}
\lefttop{j}\Gr^{\Vzero}_b(\nbigm_1)
\lrarr
\lefttop{j}\Gr^{\Vzero}_b(\nbigm_2)
\end{equation}
is strict with respect to $\lefttop{I}\Vzero$.
For any $n\in\seisuu_{\geq 0}$,
we consider
\[
 \nbign_n:=
 \bigl\{
 s\in \lefttop{I}\Vzero_{\veca}(\nbigm_2)\,\big|\,
 z_j^ns\in F(\lefttop{I}\Vzero_{\veca}(\nbigm_1))
 \bigr\}
 \subset \lefttop{I}\Vzero_{\veca}(\nbigm_2).
\]
They are coherent over $\lefttop{I}V\nbigr_X$.
We have
$\nbign_{n}\subset\nbign_{m}$ for $n\leq m$.
Hence, there exists $n_0$ such that
$\nbign_{n_0}=\nbign_{n_0+p}$ for any $p\geq 0$.
Let $b_0<0$.
We obtain
\[
 F(\lefttop{I}\Vzero_{\veca}\nbigm_1)
 \cap
 \lefttop{j}\Vzero_{b_0-n_0-p}
 \lefttop{I}\Vzero_{\veca}(\nbigm_2)
 =z_j^p\Bigl(
 F(\lefttop{I}\Vzero_{\veca}\nbigm_1)
 \cap
 \lefttop{j}\Vzero_{b_0-n_0}
 \lefttop{I}\Vzero_{\veca}(\nbigm_2)
 \Bigr).
\]
There exists $c$ such that
\[
 F(\lefttop{j}\Vzero_c\lefttop{I}\Vzero_{\veca}\nbigm_1)
\supset
 F(\lefttop{I}\Vzero_{\veca}\nbigm_1)
 \cap
 \lefttop{j}\Vzero_{b_0-n_0}(\nbigm_2).
\]
Take $m\in\seisuu_{>0}$ such that $c-m<b$.
We obtain
\begin{multline}
\label{eq;22.4.25.40}
 F(\lefttop{j}\Vzero_b\lefttop{I}\Vzero_{\veca}\nbigm_1)
  \supset
  z_j^m
   F(\lefttop{j}\Vzero_c\lefttop{I}\Vzero_{\veca}\nbigm_1)
   \supset
 z_j^m\Bigl(
 F(\lefttop{I}\Vzero_{\veca}\nbigm_1)
 \cap
 \lefttop{j}\Vzero_{b_0-n_0}
 \lefttop{I}\Vzero_{\veca}(\nbigm_2)
 \Bigr)
 \\
\supset
 F(\lefttop{I}\Vzero_{\veca}\nbigm_1)
 \cap
 \lefttop{j}\Vzero_{b_0-n_0-m}
 \lefttop{I}\Vzero_{\veca}(\nbigm_2).
\end{multline}

Let $s$ be a local section of
$\Image(F)\cap \lefttop{j}\Vzero_b\lefttop{I}\Vzero_{\veca}(\nbigm_2)$.
Let $[s]$ denote the induced section of
$\lefttop{j}\Gr^{\Vzero}_b\lefttop{I}\Vzero_{\veca}(\nbigm_2)$.
Because $F$ is strict with respect to
$\lefttop{j}\Vzero$,
$s$ is contained in
$F(\lefttop{j}\Vzero_b(\nbigm_1))$.
It implies that
$[s]$ is contained in
$F(\lefttop{j}\Gr^{\Vzero}_b(\nbigm_1))
\cap
 \lefttop{j}\Gr^{\Vzero}_b\lefttop{I}\Vzero_{\veca}(\nbigm_2)$.
There exists
$t_1\in \lefttop{j}\Vzero_b\lefttop{I}\Vzero_{\veca}(\nbigm_1)$
such that
$s-F(t_1)\in
\Image(F)\cap
 \lefttop{j}\Vzero_{<b}\lefttop{I}\Vzero_{\veca}(\nbigm_2)$.
By an easy induction,
we can prove that there exists
$t_2\in \lefttop{j}\Vzero_b\lefttop{I}\Vzero_{\veca}(\nbigm_1)$
such that
$s-F(t_2)\in
\Image(F)\cap
 \lefttop{j}\Vzero_{b_0-n_0-m}\lefttop{I}\Vzero_{\veca}(\nbigm_2)$.
By (\ref{eq;22.4.25.40}),
there exists
$t_3\in \lefttop{j}\Vzero_b\lefttop{I}\Vzero_{\veca}(\nbigm_1)$
such that
$s-F(t_2)=F(t_3)$.
Thus, we obtain the second claim.
The first claim follows from the second claim.
\hfill\qed

\vspace{.1in}
For any finite subset $S\subset \real^{I}$,
we set
$\lefttop{I}\Vzero_S(\nbigm)
=\sum_{\veca\in S}
\lefttop{I}\Vzero_{\veca}(\nbigm)$.
We have the following refinement.
\begin{prop}
\label{prop;22.2.16.9}
Let $F:\nbigm_1\lrarr\nbigm_2$
be a morphism in $\nbigc(X)$.
Suppose that $\nbigm_i$ have multi-$V$-filtrations
along $H$.
We also assume that
$\Ker F$, $\Image F$ and $\Cok(F)$
are objects in $\nbigc(X)$.
 Then, for any $I\subset\Lambda$
 and finite subset $S\subset \real^{I}$,
we have
 \[
 F(\nbigm_1)\cap \lefttop{I}V_{S}(\nbigm_2)
=F(\lefttop{I}V_{S}\nbigm_1).
\]
\end{prop}
\pf
We set
$\rho_i(S):=
\max\{q_i(\veca)\,|\,\veca\in S\}
-\min\{q_i(\veca)\,|\,\veca\in S\}$.
We set $\rho(S):=\max_i\rho_i(S)$.
We set $\kappa(S)=\bigl|\{i\,|\,\rho_i(S)=\rho(S)\}\bigr|$.
We use an induction on
the lexicographic order 
$(\dim(X),\rho(S),\kappa(S))\in
\seisuu_{>0}\times\real\times \seisuu_{>0}$.

We choose $i\in I$ such that $\rho_i(S)=\rho(S)$.
We set $I':=I\setminus\{i\}$.
We set
$a:=\max\{q_i(\veca)\,|\,\veca\in S\}$.
Let $s$ be a non-zero section of
$F(\nbigm_1)\cap \lefttop{I}V_{S}(\nbigm_2)$.
Let $S_1=\{\veca\in S\,|\,q_i(\veca)=a\}$.
We set $S_1':=q_{I'}(S_1)$.
We have the induced section
$[s]$ of
$\lefttop{i}\Gr^{\Vzero}_a(\nbigm_2)$.
It is contained in
\[
 F\bigl(\lefttop{i}\Gr^{\Vzero}_a(\nbigm_1)\bigr)
 \cap
 \lefttop{I'}
 \Vzero_{S'_1}
 \lefttop{i}\Gr^{\Vzero}_a(\nbigm_2).
\]
By the assumption of the induction,
there exists
$t_1\in \lefttop{I}\Vzero_{S_1}(\nbigm_1)$
such that
$F([t_1])=[s]$.
Let $S_2:=
\{\veca\in S\,|\,q_i(\veca)<a\}
\cup \{\veca-\epsilon\veciti_i\,|\,\veca\in S_1\}$
for small $\epsilon>0$.
Then,
we have
$s-F(t_1)\in
F(\nbigm_1)\cap
\lefttop{I}\Vzero_{S_2}(\nbigm_2)$.
Because
$(\rho(S_2),\kappa(S_2))<(\rho(S),\kappa(S))$,
there exists a section $t_2$
of $\lefttop{I}\Vzero_{S_2}(\nbigm_1)$
such that $F(t_2)=s-F(t_1)$.
Thus, we are done.
\hfill\qed

\begin{prop}
Suppose that $\nbigm\in\nbigc(X)$ has a multi-$V$-filtration
along $H$. 
We have
\[
\lefttop{\Lambda}\Vzero_S(\nbigm)
\cap
 \lefttop{I}\Vzero_{\vecb}(\nbigm)
 =\sum_{\veca\in S}
 \left(
 \lefttop{\Lambda}\Vzero_{\veca}(\nbigm)
 \cap
 \lefttop{I}\Vzero_{\vecb}(\nbigm)
 \right).
\]
\end{prop}
\pf
We use an induction on $\dim X$.
It is enough to prove the case $|I|=1$,
i.e., $I=\{i\}$.
Let $S$ be a finite subset of $\real^{\Lambda}$.
We shall prove
\begin{equation}
\label{eq;22.2.15.10}
 \lefttop{\Lambda}\Vzero_S(\nbigm)
\cap
 \lefttop{i}\Vzero_{b}(\nbigm)
 =\sum_{\veca\in S}
 \left(
 \lefttop{\Lambda}\Vzero_{\veca}(\nbigm)
 \cap
 \lefttop{i}\Vzero_{b}(\nbigm)
 \right).
\end{equation}
We set $a:=\max q_i(S)$.
If $a\leq b$,
the claim is obvious.
Suppose $a>b$.
We set $S_1:=\{\veca\in S\,|\,q_i(\veca)=a\}$.
We use an induction on $(a-b,|S_1|)$.
Let us study
\[
f\in
\lefttop{\Lambda}\Vzero_{S}(\nbigm)
 \cap
 \lefttop{i}\Vzero_b(\nbigm).
\]
The induced section
$[f]$ of $\lefttop{i}\Gr^{\Vzero}_a(\nbigm)$
is $0$.
We obtain
$\sum_{\veca\in S_1} [f_{\veca}]=0$.
If $|S_1|=1$,
we obtain $[f_{\veca}]=0$,
and we replace $S$
with $S_2=(S\setminus\{\veca\})\cup \{\veca-\epsilon\veciti_i\}$.
Because $q_1(S_2)<a$, we are done.
If $|S_1|>1$,
we choose $\veca_1\in S_1$.
We set $\Lambda'=\Lambda\setminus \{i\}$
and $S_1':=S_1\setminus\{\veca_1\}$.
We have
\[
 [f_{\veca_1}]
 =-\sum_{\veca\in S_1\setminus\{\veca_1\}} [f_{\veca}]
 \in
 \Bigl(
 \lefttop{\Lambda'}\Vzero_{q_{\Lambda'}(\veca_1)}
 \lefttop{i}\Gr^{\Vzero}_a(\nbigm)
 \Bigr)
 \cap
 \Bigl(
 \lefttop{\Lambda'}\Vzero_{q_{\Lambda'}(S_1')}
 \lefttop{i}\Gr^{\Vzero}_a(\nbigm)
 \Bigr).
\]
By assumption of the induction,
there exists
$g_{\veca}\in
\lefttop{\Lambda}\Vzero_{\veca}(\nbigm)
\cap
\lefttop{\Lambda}\Vzero_{\veca_1}(\nbigm)$
$(\veca\in S_1')$
such that
\[
 [f_{\veca_1}]
 =-\sum_{\veca\in S_1'} [f_{\veca}]
=\sum_{\veca\in S_1'}[g_{\veca}].
\]
We have
\[
f=\sum_{\veca\in S\setminus S_1}f_{\veca}
+\Bigl(
f_{\veca_1}-\sum_{\veca\in S_1'}g_{\veca}
\Bigr)
+\sum_{\veca\in S_1'}(f_{\veca}+g_{\veca}).
\]
Note that
$[f_{\veca_1}-\sum_{\veca\in S_1'}g_{\veca}]=0$
and $|S_1'|<|S_1|$.
Hence, the induction can proceed.
\hfill\qed

\subsection{Good mixed twistor $\nbigd$-modules}

\subsubsection{$\nbigr_X$-modules associated with
a regular-KMS smooth $\nbigr_{X(\ast H)}$-module}

Let $X$ and $H$ be as in \S\ref{subsection;22.2.14.1}.
For any $\vecn\in\seisuu_{\geq 0}^{\Lambda}$,
we set $\deldel^{\vecn}=\prod_{i\in\Lambda}\deldel_{z_i}^{n_i}$.
Let $\nbigv$ be a regular-KMS smooth $\nbigr_{X(\ast H)}$-module.
(See \S\ref{subsection;22.3.25.120}.)
Assume that it comes from a regular admissible mixed twistor structure
on $(X,H)$.
(See \S\ref{subsection;22.4.25.41}.)
Let $\Lambda=I\sqcup J$ be a decomposition.
Let $\veca(I,J)\in\real^{\Lambda}$ be determined by
$q_i(\veca(I,J))=1$ $(i\in I)$
and $q_i(\veca(I,J))=1-\epsilon$ $(i\in J)$
for any sufficiently small $\epsilon>0$.
We obtain the $\nbigr_X$-module
$\nbigv[\ast I!J]_{|\nbigxzero}=
\nbigr_{X|\nbigxzero}\otimes_{\lefttop{\Lambda}V\nbigr_{X|\nbigxzero}}
\nbigq^{(\lambda_0)}_{\veca(I,J)}\nbigv$.

For $\vecb\in\real_{\leq 0}^I\times\real_{<0}^J$,
we set
$\lefttop{\Lambda}\Vzero_{\vecb}(\nbigv[\ast I!J]):=
\nbigqzero_{\vecb+\veciti_{\Lambda}}(\nbigv)$.
For general $\vecc\in\real^{\Lambda}$,
we set
\[
 \lefttop{\Lambda}\Vzero_{\vecc}(\nbigv[\ast I!J])
 =\sum_{\substack{
 \vecn\in\seisuu_{\geq 0}^{\Lambda}\\
 \vecb\in \real^I_{\leq 0}\times\real^J_{<0}\\
 \vecn+\vecb\leq \vecc
 }}
 \deldel^{\vecn}
 \lefttop{\Lambda}\Vzero_{\vecb}(\nbigv[\ast I!J]).
\]
For a subset $K\subset\Lambda$ and $\vecd\in\real^K$,
we set
\[
 \lefttop{K}\Vzero_{\vecd}(\nbigv[\ast I!J])
 =\sum_{\substack{\vecc\in\real^{\Lambda}\\
 q_K(\vecc)=\vecd
  }}
 \lefttop{\Lambda}\Vzero_{\vecc}(\nbigv[\ast I!J]).
\]
The following is proved in \cite[Corollary 5.3.9]{Mochizuki-MTM}.
\begin{lem}
For each $i\in\Lambda$,
$\lefttop{i}\Vzero_{\bullet}(\nbigv[\ast I!J])$
is the $V$-filtration of
$\nbigv[\ast I!J]$ along $z_i$.
\hfill\qed
\end{lem}

\begin{prop}
For any $K\subset\Lambda$ and $\veca\in\real^K$,
the following holds:
\begin{equation}
\label{eq;22.2.16.1}
 \lefttop{K}\Vzero_{\veca}(\nbigv[\ast I!J])
=\bigcap_{i\in K}
 \lefttop{i}\Vzero_{a_i}(\nbigv[\ast I!J]).
\end{equation}
\end{prop}
\pf
We study the case $K=\Lambda$.
We use an induction on $\dim X$.
It is clear that the left hand side is contained
in the right hand side.
Choose $i(0)\in\Lambda$.
We set $\Lambda':=\Lambda\setminus\{i(0)\}$.
The image of 
\[
 \bigcap_{i\in\Lambda}
 \lefttop{i}\Vzero_{a_i}(\nbigv[\ast I!J])
 \lrarr
 \lefttop{i(0)}\Gr^{\Vzero}_{a_{i(0)}}
 (\nbigv[\ast I!J])
\]
is contained in
\begin{equation}
\label{eq;22.2.15.1}
 \bigcap_{i\in\Lambda'}
 \lefttop{i}\Vzero_{a_i}
 \bigl(
 \lefttop{i(0)}\Gr^{\Vzero}_{a_{i(0)}}
 (\nbigv[\ast I!J])
 \bigr)
 =\lefttop{\Lambda'}
 \Vzero_{q_{\Lambda'}(\veca)}
  \bigl(
 \lefttop{i(0)}\Gr^{\Vzero}_{a_{i(0)}}
 (\nbigv[\ast I!J])
 \bigr).
\end{equation}
Here, the equality follows from the assumption of the induction.
(Note Proposition \ref{prop;22.4.25.50}
and the description \cite[(5.16)]{Mochizuki-MTM}
of $\lefttop{i}\Gr^{\Vzero}_a(\nbigv[\ast I!J])$.)
The image of
\[
 \lefttop{\Lambda}\Vzero_{\veca}(\nbigv[\ast I!J])
 \lrarr
 \lefttop{i(0)}\Gr^{\Vzero}_{a_{i(0)}}
 (\nbigv[\ast I!J])
\]
is equal to
the right hand side of (\ref{eq;22.2.15.1}).
Hence, we obtain
\[
\bigcap_{i\in\Lambda}
\lefttop{i}\Vzero_{a_i}(\nbigv[\ast I!J])
=\lefttop{\Lambda}\Vzero_{\veca}(\nbigv[\ast I!J])
+\Bigl(
\lefttop{i(0)}\Vzero_{<a_{i(0)}}(\nbigv[\ast I!J])
\cap
\bigcap_{i\in\Lambda'}
\lefttop{i}\Vzero_{a_i}(\nbigv[\ast I!J])
\Bigr).
\]
By using an inductive argument, for any $b<0$,
we obtain
\[
\bigcap_{i\in\Lambda}
\lefttop{i}\Vzero_{a_i}(\nbigv[\ast I!J])
=\lefttop{\Lambda}\Vzero_{\veca}(\nbigv[\ast I!J])
+\Bigl(
\lefttop{i(0)}\Vzero_{b}(\nbigv[\ast I!J])
\cap
\bigcap_{i\in\Lambda'}
\lefttop{i}\Vzero_{a_i}(\nbigv[\ast I!J])
\Bigr).
\]
By apply the same argument to any $i\in\Lambda$,
for any $\vecb\in\real_{<0}^{\Lambda}$,
we obtain
\[
 \bigcap_{i\in\Lambda}
\lefttop{i}\Vzero_{a_i}(\nbigv[\ast I!J])
=\lefttop{\Lambda}\Vzero_{\veca}(\nbigv[\ast I!J])
+
\bigcap_{i\in\Lambda}
\lefttop{i}\Vzero_{b_i}(\nbigv[\ast I!J]).
\]

Let us consider the natural morphism
$\varphi:\nbigv[\ast I!J]\lrarr\nbigv$.
We set
\[
\lefttop{i}\Vzero_a(\nbigv)
=\sum_{\substack{\vecb\in\real^{\Lambda}\\ q_i(\vecb)=a+1}}
\nbigqzero_{\vecb}(\nbigv).
\]
We have 
$\varphi\bigl(
\lefttop{i}\Vzero_a(\nbigv[\ast I!J])
\bigr)
\subset
\lefttop{i}\Vzero_a(\nbigv)$.
For $\vecb\in\real^{\Lambda}_{<0}$,
we have
\[
\bigcap_{i\in\Lambda}
 \lefttop{i}\Vzero_{b_i}(\nbigv)
 =\varphi\Bigl(
\lefttop{\Lambda}\Vzero_{\vecb}(\nbigv[\ast I!J])
 \Bigr).
\]
Hence, we obtain
\[
 \bigcap_{i\in\Lambda}
\lefttop{i}\Vzero_{a_i}(\nbigv[\ast I!J])
=\lefttop{\Lambda}\Vzero_{\veca}(\nbigv[\ast I!J])
+
\Bigl(
\Ker\varphi\cap
\bigcap_{i\in\Lambda}
\lefttop{i}\Vzero_{b_i}(\nbigv[\ast I!J])
\Bigr).
\]

Let $s$ be any section of
$\Ker\varphi\cap
\bigcap_{i\in\Lambda}
\lefttop{i}\Vzero_{b_i}(\nbigv[\ast I!J])$.
Because $\varphi(s)=0$,
there exists $K\subset\Lambda$
and $N\in\seisuu_{>0}$
such that $\prod_{i\in K}z_i^N s=0$.
If $K$ is not empty,
we choose $i(1)\in K$
and we set $K'=K\setminus\{i(1)\}$
and $s'=\prod_{i\in K'}z_i^Ns$.
We note that
$s'\in \bigcap_{i\in \Lambda}\lefttop{i}\Vzero_{b_i}(\nbigv[\ast I!J])$.
We have $z_{i(1)}^Ns'=0$
and $s'\in \lefttop{i(1)}\Vzero_{<0}(\nbigv[\ast I!J])$.
It implies $s'=0$.
(For example, see Lemma \ref{lem;22.4.25.100}.)
Hence, by an easy induction, we obtain $s=0$.

Because
\[
 \bigcap_{i\in K}
 \lefttop{i}\Vzero_{a_i}(\nbigv[\ast I!J])
\subset\sum_{q_K(\vecb)=\veca}
  \bigcap_{i\in \Lambda}
 \lefttop{i}\Vzero_{b_i}(\nbigv[\ast I!J])
 =\sum_{q_K(\vecb)=\veca}
 \lefttop{\Lambda}
 \Vzero_{\vecb}(\nbigv[\ast I!J])
=\lefttop{K}\Vzero_{\veca}(\nbigv[\ast I!J]),
\]
we obtain the claim in the general case.
\hfill\qed

\begin{cor}
\label{cor;22.3.13.1}
The $\nbigr_X$-module $\nbigv[\ast I!J]$
has a multi-$V$-filtration along $H$. 
\hfill\qed
\end{cor}

\subsubsection{Good mixed twistor $\nbigd$-modules}

\begin{prop}
\label{prop;22.2.15.2}
Let $\nbigm$ be the $\nbigr_X$-module
underlying a good mixed twistor $\nbigd$-module
on $(X,H)$.
Then, $\nbigm$ has a multi-$V$-filtration along $H$.
(See {\rm\S\ref{subsection;22.4.25.101}} for the notion of
good mixed twistor $\nbigd$-module.)
\end{prop}
\pf
We use an induction on $\dim X$.
We assume that the claim is already proved
in the lower dimensional case.
Let $\nbigt$ be a good mixed twistor $\nbigd$-module on $(X,H)$.
Suppose that the support of $\nbigt$ is contained in $H$.
Let $\Lambda_0\subset\Lambda$ be a minimal subset
such that
the support of $\nbigt$ is contained in
$\bigcup_{i\in\Lambda_0}H_i$.
We shall prove that the underlying $\nbigr_X$-module of $\nbigt$
has a multi-$V$-filtration along $H$
by using an induction on $|\Lambda_0|$.
If $|\Lambda_0|=1$,
we may apply the assumption that
the claim is already proved in the lower dimensional case.
Let $i(0)\in\Lambda_0$.
We set $\Lambda_1=\Lambda_0\setminus\{i(0)\}$.
Suppose that $\nbigt(\ast \Lambda_1)\neq 0$.
Let $\iota_{i(0)}:H_{i(0)}\to X$ denote the inclusion.
There exists a good mixed twistor $\nbigd$-module
$\nbigt_{i(0)}$ on $(H_{i(0)},\del H_{i(0)})$
such that
$\iota_{i(0)\dagger}\nbigt_{i(0)}(\ast \Lambda_1)
\simeq
 \nbigt(\ast\Lambda_1)$.
We set $g:=\prod_{i\in\Lambda_1}z_i$.
Because
$\Xi^{(a)}_{g}\nbigt_{i(0)}$,
and 
$\psi^{(a)}_{g}\nbigt_{i(0)}$
are good mixed twistor $\nbigd$-module
on $(H_{i(0)},\del H_{i(0)})$,
they have multi-$V$-filtrations along $\del H_{i(0)}$.
Hence,
the underlying $\nbigr_X$-modules of
$\iota_{i(0)\dagger}
\Xi^{(a)}_{g}\nbigt_{i(0)}$
and 
$\iota_{i(0)\dagger}\psi^{(a)}_{g}\nbigt_{i(0)}$
have multi-$V$-filtration along $H$.
Because $\phi^{(0)}_{g}\nbigt$ is a good mixed twistor $\nbigd$-module
on $(X,H)$ whose support is contained in
$\bigcup_{i\in\Lambda_1}H_i$,
the underlying $\nbigr_X$-module has a multi-$V$-filtration.
Then, we obtain that the underlying $\nbigr_X$-module
$\nbigt$ also has a multi-$V$-filtration.

Let us study the case where the support of $\nbigt$
is not necessarily contained in $H$.
Let $f=\prod_{i\in\Lambda}z_i$.
Because
$\Pi^{a,b}_{f}\nbigt(\ast H)$ is an admissible mixed twistor structure
on $(X,H)$,
the underlying $\nbigr_X$-modules of
$\Pi^{a,b}_{f\star}\nbigt$ $(\star=!,\ast)$
have multi-$V$-filtrations along $H$
by Corollary \ref{cor;22.3.13.1}.
Hence,
the underlying $\nbigr_X$-module of
$\Xi^{(0)}_f(\nbigt)$ has a multi-$V$-filtration along $H$.
Because $\phi^{(0)}_f(\nbigt)$
and $\psi^{(a)}_f(\nbigt)$ are good mixed twistor
$\nbigd$-modules whose supports are contained in $H$,
their underlying $\nbigr_X$-modules
have multi-$V$-filtrations along $\del H$.
Then, we obtain that
the underlying $\nbigr_X$-module of $\nbigt$
also has a multi-$V$-filtration.
\hfill\qed

\subsubsection{A lifting}

Let $\veca\in\real^{\Lambda}$.
Let $I\subset\Lambda$ be any non-empty subset.
For an $\nbigr_X$-module
underlying a good mixed twistor $\nbigd$-module on $(X,H)$,
we define
\[
 \del^{(I)}_{\veca}\nbigm
 =\bigoplus_{i\in I}
 \lefttop{i}\Gr^{\Vzero}_{a_i}(\lefttop{\Lambda}\Vzero_{\veca}\nbigm).
\]
There exists the natural morphism
\[
 \kappa^{(I)}_{\veca,\nbigm}:
 \lefttop{\Lambda}\Vzero_{\veca}(\nbigm)
 \lrarr
 \del^{(I)}_{\veca}\nbigm.
\]

Let $f:\nbigm_1\lrarr\nbigm_2$
be a morphism of $\nbigr_X$-modules
underlying a morphism of
good mixed twistor $\nbigd$-modules on $(X,H)$.
We obtain the following induced morphism
\[
\del^{(I)}_{\veca} f:
\del^{(I)}_{\veca}\nbigm_1\lrarr\del^{(I)}_{\veca}\nbigm_2.
\]

\begin{lem}
\label{lem;22.2.16.23}
We have
\[
 \Image \kappa^{(I)}_{\veca,\nbigm_1}
 \cap
 \Ker (\del^{(I)}_{\veca} f)
 =\Image \kappa^{(I)}_{\veca,\Ker f},
 \quad
 \Image \kappa^{(I)}_{\veca,\nbigm_2}
 \cap
 \Image (\del^{(I)}_{\veca} f)
 =\Image \kappa^{(I)}_{\veca,\Image f}.
\] 
\end{lem}
\pf
We use an induction on $|I|$.
We choose $i\in I$.
If $|I|=1$, i.e., $I=\{i\}$,
it follows from the strictness of $f$
with respect to $\lefttop{\Lambda}\Vzero$
(Proposition \ref{prop;22.3.13.2}).
We clearly have
\[
  \Image \kappa^{(I)}_{\veca,\nbigm_1}
 \cap \Ker \del_{\veca}^{(I)}f
\supset\Image \kappa^{(I)}_{\veca,\Ker f},
 \quad
 \Image \kappa^{(I)}_{\veca,\nbigm_2}
 \cap \Image \del_{\veca}^{(I)}f
\supset\Image \kappa^{(I)}_{\veca,\Image f}.
\]
Let $s\in\Image\kappa^{(I)}_{\veca,\nbigm_1}
\cap\Ker\del_{\veca}^{(I)}f$.
We set $I':=I\setminus \{i\}$.
Let $p_{I',I}:\del^{(I)}\nbigm\lrarr\del^{(I')}\nbigm$
be the projection.
There exists $t_1\in\lefttop{\Lambda}\Vzero_{\veca}\Ker(f)$
such that
\[
 p_{I',I}(s)=\kappa^{(I')}_{\veca,\Ker f}(t_1).
\]
We obtain
\[
 s-\kappa^{(I)}_{\veca,\nbigm_1}(t_1)
 \in
 \lefttop{i}\Gr^{\Vzero}_{a_i}(\lefttop{\Lambda}\Vzero_{\veca}\nbigm_1)
 \cap \Ker\del_{\veca}^{(I)}f
 \subset
 \del^{(I)}\nbigm_1.
\]
We set $\veca'=\veca-\epsilon\veciti_{I'}$
for a sufficiently small $\epsilon>0$.
By the construction, we have
\[
 s-\kappa^{(i)}_{a_i,\nbigm_1}(t_1)
 \in \lefttop{I'}\Vzero_{<q_{I'}(\veca)}
\lefttop{i}\Gr^{\Vzero}_{a_i}(\lefttop{\Lambda}\Vzero_{\veca}\nbigm_1)
\cap\Ker\del_{a_i}^{(i)}f
=\lefttop{i}\Gr^{\Vzero}_{a_i}
 \lefttop{\Lambda}\Vzero_{\veca'}(\nbigm_1)
 \cap
 \Ker\del^{(i)}_{a_i}(f).
\]
There exists
$t_2\in
\lefttop{\Lambda}\Vzero_{\veca'}\Ker(f)$
which induces $s-\kappa^{(i)}_{a_i,\nbigm_1}(t_1)$.
We obtain
$s=\kappa^{(I)}_{\veca,\nbigm_1}(t_1+t_2)$.

Let $s$ be a section of
$\Image\kappa^{(I)}_{\veca,\nbigm_2}
\cap\Image\del^{(I)}_{\veca}f$.
There exists $t_1\in\nbigm_1$
such that
\[
 p_{I',I}(s)
=\kappa^{(I')}_{q_{I'}(\veca),\nbigm_2}(f(t_1)).
\]
We obtain
\[
 s-\kappa^{(I)}_{\veca,\nbigm_2}(f(t_1))
 \in
 \lefttop{i}\Gr^{\Vzero}_{a_i}
 \lefttop{\Lambda}\Vzero_{\veca}(\nbigm_2)
 \cap
 \Image\del_{\veca}^{(I)}(f)
 \subset
 \del_{\veca}^{(I)}(\nbigm_2).
\]
By the construction,
we have
\[
 s-\kappa^{(i)}_{a_i,\nbigm_2}(f(t_1))
 \in
 \lefttop{I'}\Vzero_{<q_{I'}(\veca)}
 \lefttop{i}\Gr^{\Vzero}_{a_i}
 \lefttop{\Lambda}(\nbigm_2)
 \cap
 \Image\del_{a_i}^{(i)}(f)
=\lefttop{i}\Gr^{\Vzero}_{a_i}
 \lefttop{\Lambda}\Vzero_{\veca'}(\nbigm_2)
 \cap
 \Image\del_{a_i}^{(i)}(f).
\]
Hence, there exists
$t_2\in
\lefttop{\Lambda}\Vzero_{\veca'}(\nbigm_1)$
such that
$\kappa^{(i)}_{a_i,\nbigm_2}(f(t_2))
=s-\kappa^{(i)}_{a_i,\nbigm_2}(f(t_1))$.
Hence, we obtain
$s=\kappa^{(I)}_{\veca,\nbigm_2}(f(t_1+t_2))$.
\hfill\qed

\subsubsection{The associated graded objects}

Let $\nbigm$ the $\nbigr_X$-module
underlying a good mixed twistor $\nbigd$-module
on $(X,H)$.
Let $J,I\subset\Lambda$.
For $\veca\in\real^{I\cup J}$,
we obtain
\[
 \lefttop{J}\Gr^{V^{(\lambda_0)}}_{q_J(\veca)}
 \lefttop{I}V_{q_I(\veca)}^{(\lambda_0)}(\nbigm)
 :=
 \frac{\lefttop{I\cup J}\Vzero_{\veca}(\nbigm)}
 {\sum_{\vecb\in\nbigs(\veca,J)}
 \lefttop{I\cup J}V_{\vecb}^{(\lambda_0)}(\nbigm)},
 \quad\quad
 \nbigs(\veca,J)=\bigl\{
 \vecb\in\real^{I\cup J}\,\big|\,
 q_{I\setminus J}(\vecb)=q_{I\setminus J}(\veca),\,\,
 q_J(\vecb)\lneq q_J(\veca)
 \bigr\}.
\]
They are strict in the sense that
the multiplication of $\lambda-\lambda_1$ is a monomorphism
for any $\lambda_1\in\cnum$.
There exist the induced endomorphisms $\Res_j(\DD)$ $(j\in J)$
on 
$\lefttop{J}\Gr^{V^{(\lambda_0)}}_{q_J(\veca)}
 \lefttop{I}\Vzero_{\veca}(\nbigm)$.
There exists the decomposition
\[
 \lefttop{J}\Gr^{V^{(\lambda_0)}}_{q_J(\veca)}
 \lefttop{I}\Vzero_{q_I(\veca)}(\nbigm)
 =\bigoplus_{\substack{
 \vecu\in(\real\times\cnum)^J\\
 \paramap(\lambda_0,\vecu)=q_J(\veca)
 }}
 \lefttop{J}\psi^{(\lambda_0)}_{\vecu}
 \lefttop{I}\Vzero_{q_I(\veca)}(\nbigm),
\]
where
$\nbign_j=-\Res_j(\DD)+(\eigenmap(\lambda,u_j)+\lambda)\id$
$(j\in J)$
are nilpotent on
$\lefttop{J}\psi^{(\lambda_0)}_{\vecu}
 \lefttop{I}V_{q_I(\veca)}^{(\lambda_0)}(\nbigm)$.
The natural morphism
\[
 \lefttop{J}\Gr^{\Vzero}_{q_J(\veca)}
 \lefttop{I}\Vzero_{q_I(\veca)}(\nbigm)
 \lrarr
 \lefttop{J}\Gr^{\Vzero}_{q_J(\veca)}(\nbigm)
\]
is a monomorphism.
It induces a monomorphism
\[
 \lefttop{J}\psizero_{\vecu}
 \lefttop{I}\Vzero_{q_I(\veca)}(\nbigm)
 \lrarr
 \lefttop{J}\psizero_{\vecu}(\nbigm).
\]
We also write
$\lefttop{J}\psizero_{\vecu}
\lefttop{I}\Vzero_{q_I(\veca)}(\nbigm)$
as
$\lefttop{I\setminus J}\Vzero_{q_{I\setminus J}(\veca)}
\lefttop{J}\psizero_{\vecu}(\nbigm)$.
We remark the following lemma.
\begin{lem}
\label{lem;22.4.26.2}
There exists the good $\nbigr_{H_J}$-module
$\lefttop{J}\psi_{\vecu}(\nbigm)$
such that
(i) $\lefttop{J}\psi_{\vecu}(\nbigm)$ underlies
a good mixed twistor $\nbigd$-module
$\psi_{\vecu}(\nbigt)$ on $(H_J,\del H_J)$,
(ii)
 $\lefttop{J}\psi_{\vecu}(\nbigm)_{|\nbighzero_{J}}
 =\lefttop{J}\psizero_{\vecu}(\nbigm)$.
Moreover,
each $\nbign_j$ underlies a morphism
 $\lefttop{J}\psi_{\vecu}(\nbigt)
 \to
 \lefttop{J}\psi_{\vecu}(\nbigt)\otimes\newTate(-1)$.
\end{lem}
\pf
Let $J=(j(1),\ldots,j(m))$.
If $q_i(\veca)>0$,
we have
$-\deldel_{z_i}:
\lefttop{J}\psizero_{\vecu-\vecdelta_i}(\nbigm)
\simeq
\lefttop{J}\psizero_{\vecu}(\nbigm)$.
If $q_i(\veca)<0$,
we have
$z_i:\lefttop{J}\psizero_{\vecu}(\nbigm)
\simeq
\lefttop{J}\psizero_{\vecu-\vecdelta_i}(\nbigm)$.
Hence, we may assume that $q_J(\veca)\in\real_{\leq 0}^J$.
We have
$\lefttop{J}\psi_{\vecu}(\nbigm)
=\lefttop{j_1}\psi_{q_{j(1)}(\vecu)}\circ\cdots\circ
 \lefttop{j_m}\psi_{q_{j(m)}(\vecu)}(\nbigm)$.
Then, the claim of the lemma follows.
\hfill\qed

\vspace{.1in}

For $K\subset\Lambda\setminus J$,
and for $\veca\in\real^{I\cup J\cup K}$,
we obtain
\[
 \lefttop{K}\Gr^{\Vzero}_{q_K(\veca)}
 \lefttop{J}\Gr^{\Vzero}_{q_J(\veca)}
 \lefttop{I}\Vzero_{q_I(\veca)}(\nbigm)
:=\frac{\lefttop{J}\Gr^{\Vzero}_{q_J(\veca)}
 \lefttop{I\cup K}\Vzero_{q_{I\cup K}(\veca)}(\nbigm)}
 {\sum_{\vecb\in\nbigs(q_{I\cup K}(\veca),K)}
 \lefttop{J}\Gr^{\Vzero}_{q_J(\veca)}
 \lefttop{I\cup K}\Vzero_{\vecb}(\nbigm)}.
\]
In this case, there exists the natural isomorphism:
\begin{equation}
\label{eq;22.2.10.1}
 \lefttop{K}\Gr^{\Vzero}_{q_K(\veca)}
 \lefttop{J}\Gr^{V^{(\lambda_0)}}_{q_J(\veca)}
 \lefttop{I}\Vzero_{\veca}(\nbigm)
 \simeq
  \lefttop{K\sqcup J}\Gr^{\Vzero}_{q_{K\sqcup J}(\veca)}
  \lefttop{I}\Vzero_{\veca}(\nbigm).
\end{equation}
For $\vecu\in(\real\times\cnum)^J$
and $\vecc\in \real^{I\cup K}$,
we obtain
\[
 \lefttop{K}\Gr^{\Vzero}_{q_K(\vecc)}
 \lefttop{J}\psizero_{\vecu}
 \lefttop{I}\Vzero_{q_I(\vecc)}(\nbigm)
:=\frac{\lefttop{J}\psizero_{\vecu}
 \lefttop{I\cup K}\Vzero_{\vecb}(\nbigm)}
 {\sum_{\vecb\in\nbigs(\vecc,K)}
 \lefttop{J}\psizero_{\vecu}
 \lefttop{I\cup K}\Vzero_{\vecb}(\nbigm)}.
\]
We obtain the decomposition
\[
 \lefttop{K}\Gr^{\Vzero}_{q_K(\vecc)}
 \lefttop{J}\psizero_{\vecu}
 \lefttop{I}\Vzero_{q_I(\vecc)}(\nbigm)
 =\bigoplus_{\substack{
 \vecv\in(\real\times\cnum)^K
 \\
 \paramap(\lambda_0,\vecv)=q_K(\vecc) }}
 \lefttop{K}\psizero_{\vecv}
 \lefttop{J}\psizero_{\vecu}
 \lefttop{I}\Vzero_{q_I(\vecc)}(\nbigm),
\]
where
$\nbign_j=-\Res_j(\DD)+(\eigenmap(\lambda,u_j)+\lambda)\id$
$(j\in K)$
are nilpotent on
$\lefttop{K}\psizero_{\vecv}
 \lefttop{J}\psizero_{\vecu}
 \lefttop{I}\Vzero_{q_I(\vecc)}(\nbigm)$.
The isomorphism (\ref{eq;22.2.10.1})
induces
\[
 \lefttop{K}\psizero_{\vecv}
 \lefttop{J}\psizero_{\vecu}
 \lefttop{I}\Vzero_{q_I(\vecc)}(\nbigm)
 \simeq
 \lefttop{J\cup K}\psizero_{(\vecu,\vecv)}
 \lefttop{I}\Vzero_{q_I(\vecc)}(\nbigm),
\]
where $(\vecu,\vecv)\in\real^{J\cup K}$
is the element
such that
$q_J(\vecu,\vecv)=\vecu$
and
$q_K(\vecu,\vecv)=\vecv$.

\subsubsection{A filtration}

For $\ell\in\seisuu_{\geq 0}$,
we have the following morphism:
\[
 \lefttop{I}\Gr^{\Vzero}_{q_I(\veca)}
 \lefttop{\Lambda}\Vzero_{\veca}(\nbigm)
 \lrarr
 \bigoplus_{\substack{
 K\supset I\\
 |K\setminus I|=\ell+1
 }}
 \lefttop{K}\Gr^{\Vzero}_{q_K(\veca)}
 \lefttop{\Lambda}\Vzero_{\veca}(\nbigm).
\]
Let
$F_{\ell}
\lefttop{I}\Gr^{\Vzero}_{q_I(\veca)}
\lefttop{\Lambda}\Vzero_{\veca}(\nbigm)$
denote the kernel.
Thus, we obtain an increasing filtration $F_{\bullet}$
on
$\lefttop{I}\Gr^{\Vzero}_{q_I(\veca)}
\lefttop{\Lambda}\Vzero_{\veca}(\nbigm)$.
We set
\[
 \lefttop{\Lambda\setminus I}\Vzero_{<q_{\Lambda\setminus I}(\veca)}
 \lefttop{I}\Gr^{\Vzero}_{q_I(\veca)}
 \lefttop{\Lambda}\Vzero_{\veca}(\nbigm):=
F_0
  \lefttop{I}\Gr^{\Vzero}_{q_I(\veca)}
 \lefttop{\Lambda}\Vzero_{\veca}(\nbigm).
\]

Because the kernel of the natural morphism
\[
F_{\ell}
\lefttop{I}\Gr^{\Vzero}_{q_I(\veca)}
\lefttop{\Lambda}\Vzero_{\veca}(\nbigm)
\lrarr
 \bigoplus_{\substack{
 K\supset I\\
 |K\setminus I|=\ell
 }}
 \lefttop{K}\Gr^{\Vzero}_{q_K(\veca)}
 \lefttop{\Lambda}\Vzero_{\veca}(\nbigm)
\]
is $F_{\ell-1}
\lefttop{I}\Gr^{\Vzero}_{q_I(\veca)}
\lefttop{\Lambda}\Vzero_{\veca}(\nbigm)$,
we obtain
\begin{equation}
\label{eq;22.2.16.2}
 \Gr^{F}_{\ell}
\lefttop{I}\Gr^{\Vzero}_{q_I(\veca)}
\lefttop{\Lambda}\Vzero_{\veca}(\nbigm)
\lrarr
 \bigoplus_{\substack{
 K\supset I\\
 |K\setminus I|=\ell
 }}
 \lefttop{K}\Gr^{\Vzero}_{q_K(\veca)}
 \lefttop{\Lambda}\Vzero_{\veca}(\nbigm).
\end{equation}

\begin{lem}
The morphism {\rm(\ref{eq;22.2.16.2})} induces
the following isomorphism:
\begin{equation}
\label{eq;22.2.16.3}
 \Gr^F_{\ell}
  \lefttop{I}\Gr^{\Vzero}_{q_I(\veca)}
  \lefttop{\Lambda}\Vzero_{\veca}(\nbigm)
  \simeq
  \bigoplus_{\substack{K\supset I\\
  |K\setminus I|=\ell
  }}
  \lefttop{\Lambda\setminus K}
  \Vzero_{<q_{\Lambda\setminus K}(\veca)}
  \lefttop{K}\Gr^{\Vzero}_{q_K(\veca)}
  \lefttop{\Lambda}\Vzero_{\veca}(\nbigm).
\end{equation}
\end{lem}
\pf
Because the image of {\rm(\ref{eq;22.2.16.2})}
is contained in the right hand side of (\ref{eq;22.2.16.3}),
we obtain the morphism (\ref{eq;22.2.16.3}).
It is easy to check the surjectivity by the construction.
\hfill\qed

\begin{lem}
\label{lem;22.2.16.30}
Let $f:\nbigm_1\to\nbigm_2$
be a morphism of $\nbigr_X$-modules
underlying a morphism of good mixed twistor $\nbigd$-modules
on $(X,H)$.
Then,
 for the induced morphism
$\lefttop{J}\psizero_{\vecu}
 \lefttop{I}\Vzero_{q_I(\veca)}(f):
\lefttop{J}\psizero_{\vecu}
 \lefttop{I}\Vzero_{q_I(\veca)}(\nbigm_1)
 \to
 \lefttop{J}\psizero_{\vecu}
 \lefttop{I}\Vzero_{q_I(\veca)}(\nbigm_2)$,
 we have 
\[
 \lefttop{J}\psizero_{\vecu}
 \lefttop{I}\Vzero_{q_I(\veca)}(\Ker f)
 \simeq
 \Ker \lefttop{J}\psizero_{\vecu}
 \lefttop{I}\Vzero_{q_I(\veca)}(f),
\]
\[
 \lefttop{J}\psizero_{\vecu}
 \lefttop{I}\Vzero_{q_I(\veca)}(\Image f)
 \simeq
 \Image \lefttop{J}\psizero_{\vecu}
 \lefttop{I}\Vzero_{q_I(\veca)}(f),
\]
\[
  \lefttop{J}\psizero_{\vecu}
 \lefttop{I}\Vzero_{q_I(\veca)}(\Cok f)
 \simeq
 \Cok \lefttop{J}\psizero_{\vecu}
 \lefttop{I}\Vzero_{q_I(\veca)}(f).
\]
Moreover,
$\lefttop{J}\psizero_{\vecu}
\lefttop{I}\Vzero_{q_I(\veca)}(f)$
is strict with respect to the filtration $F_{\bullet}$.
\end{lem}
\pf
It follows from Proposition \ref{prop;22.2.16.9}.
\hfill\qed

\subsubsection{Modification}
\label{subsection;22.2.17.20}

For any $I\subset J\subset\Lambda$,
the projection $(\real\times\cnum)^J\to(\real\times\cnum)^I$
is also denoted by $q_I$.
Let $\vecdelta_I$ denote the element
determined by
$q_i(\vecdelta_I)=\vecdelta$ $(i\in I)$
and 
$q_i(\vecdelta_I)=0$ $(i\not\in I)$.
For any $\vecu\in(\real\times\cnum)^J$,
$q_i(\vecu)$ $(i\in J)$ is also denoted by $u_i$.

Let $\veca\in\real^{\Lambda}$.
We fix $\vecu(0)\in(\real\times\cnum)^{\Lambda}$
such that $\paramap(\lambda_0,\vecu(0))=\veca$.
For $I\subset\Lambda$,
there exists the following morphism:
\begin{equation}
\label{eq;22.2.16.10}
 \lefttop{I}\psizero_{q_I(\vecu(0))}
 \lefttop{\Lambda}\Vzero_{\veca}(\nbigm)
 \lrarr
 \bigoplus_{j\in\Lambda\setminus I}
 \lefttop{j}\Gr^{\Vzero}_{a_j}
  \lefttop{I}\psizero_{q_I(\vecu(0))}
 \lefttop{\Lambda}\Vzero_{\veca}(\nbigm)
=
 \bigoplus_{j\in\Lambda\setminus I}
 \bigoplus_{\substack{u\in\real\times\cnum\\
 \paramap(\lambda_0,u)=a_j
 }}
 \lefttop{j}\psizero_{u}
 \lefttop{I}\psizero_{q_I(\vecu(0))}
 \lefttop{\Lambda}\Vzero_{\veca}(\nbigm).
\end{equation}
Let
$ \lefttop{I}\psizero_{q_I(\vecu(0))}
\lefttop{\Lambda}\Vtilde^{(\lambda_0)}_{\veca}(\nbigm)$
denote the inverse image of
$\bigoplus_{j\in\Lambda\setminus I}
 \lefttop{j}\psizero_{q_j(\vecu(0))}
  \lefttop{I}\psizero_{q_I(\vecu(0))}
\lefttop{\Lambda}\Vzero_{\veca}(\nbigm)$.

For $k\in\Lambda\setminus I$,
the morphism
$\lefttop{I}\psizero_{q_I(\vecu(0))}
\lefttop{\Lambda}\Vzero_{\veca}(\nbigm)
 \lrarr
 \lefttop{k}\Gr^{\Vzero}_{q_k(\veca)}
 \lefttop{I}\psizero_{q_I(\vecu(0))}
 \lefttop{\Lambda}\Vzero_{\veca}(\nbigm)$
induces
\begin{equation}
\label{eq;22.3.14.1}
  \lefttop{I}\psizero_{q_I(\vecu(0))}
\lefttop{\Lambda}\Vtilde^{(\lambda_0)}_{\veca}(\nbigm)
\lrarr
\lefttop{I\cup \{k\}}\psizero_{q_{I\sqcup \{k\}}(\vecu(0))}
 \lefttop{\Lambda}\Vtilde^{(\lambda_0)}_{\veca}(\nbigm).
\end{equation}
For $L\subset\Lambda\setminus I$,
we obtain the following morphism:
\begin{equation}
\label{eq;22.3.14.3}
\rho_{1,L}:
  \lefttop{I}\psizero_{q_I(\vecu(0))}
\lefttop{\Lambda}\Vtilde^{(\lambda_0)}_{\veca}(\nbigm)
\lrarr
\bigoplus_{k\in L}
\lefttop{I\cup \{k\}}\psizero_{q_{I\sqcup \{k\}}(\vecu(0))}
 \lefttop{\Lambda}\Vtilde^{(\lambda_0)}_{\veca}(\nbigm).
\end{equation}
There exists the following morphism of ambient sheaves:
\begin{equation}
\rho_{0,L}:
 \lefttop{I}\psizero_{q_I(\vecu(0))}
\lefttop{\Lambda}\Vzero_{\veca}(\nbigm)
\lrarr
\bigoplus_{k\in L}
 \lefttop{k}\Gr^{\Vzero}_{q_{k}(\veca)}
 \lefttop{I}\psizero_{q_I(\vecu(0))}
 \lefttop{\Lambda}\Vzero_{\veca}(\nbigm).
\end{equation}

\begin{prop}
\label{prop;22.3.14.2}
The following holds
in $\bigoplus_{k\in L}
 \lefttop{k}\Gr^{\Vzero}_{q_k(\veca)}
 \lefttop{I}\psizero_{q_I(\vecu(0))}
 \lefttop{\Lambda}\Vzero_{\veca}(\nbigm)$:
\[
 \Image\rho_{1,L}=
 \Image\rho_{0,L}
 \cap\left(
 \bigoplus_{k\in L}
\lefttop{I\cup \{k\}}\psizero_{q_{I\sqcup \{k\}}(\vecu(0))}
 \lefttop{\Lambda}\Vtilde^{(\lambda_0)}_{\veca}(\nbigm)
 \right).
\] 
In particular, the morphism {\rm(\ref{eq;22.3.14.1})}
is an epimorphism.  
\end{prop}
\pf
We use an induction on $\dim X$,
We assume that the claim of the proposition
is already proved in the lower dimensional case.
We also use an induction on $|L|$.

Let us consider the case $|L|=1$,
i.e., $L$ consists of one element $k$.
We set $\Itilde=I\sqcup\{k\}$.
Let $s$ be a section of
$\lefttop{\Itilde}\psizero_{q_{\Itilde}(\vecu(0))}
 \lefttop{\Lambda}\Vtilde^{(\lambda_0)}_{\veca}(\nbigm)$.
There exists a section $f$ of
$\lefttop{I}\psizero_{q_{I}(\vecu(0))}
\lefttop{\Lambda}V^{(\lambda_0)}_{\veca}(\nbigm)$
which induces $s$.
We obtain the induced sections
$[f]_j$ $(j\in \Lambda\setminus\Itilde)$
of
\[
\lefttop{j}\Gr^{\Vzero}_{a_j}
\lefttop{I}\psizero_{q_{I}(\vecu(0))}
\lefttop{\Lambda}V^{(\lambda_0)}_{\veca}(\nbigm)
=\bigoplus_{\substack{u\in\real\times\cnum\\
 \paramap(u)=a_j
 }}
\lefttop{j}\psizero_{u}
\lefttop{I}\psizero_{q_{I}(\vecu(0))}
\lefttop{\Lambda}V^{(\lambda_0)}_{\veca}(\nbigm).
\]
We have the decomposition
$[f]_j=\sum [f]_{j,u}$,
where $[f]_{j,u}$
is a section of 
$\lefttop{j}\psizero_{u}
\lefttop{I}\psizero_{q_{I}(\vecu(0))}
\lefttop{\Lambda}V^{(\lambda_0)}_{\veca}(\nbigm)$.
Let $T(f)$ denote the set of
$j\in \Lambda\setminus I$
such that the following holds.
\begin{itemize}
 \item $[f]_{j,q_j(\vecu(0))}\in\lefttop{j}\psizero_{q_j(\vecu(0))}
\lefttop{I}\psizero_{q_{I}(\vecu(0))}
\lefttop{\Lambda}\Vtilde^{(\lambda_0)}_{\veca}(\nbigm)$.
 \item $\sum_{u\neq q_j(\vecu(0))}[f]_{j,u}=0$.
\end{itemize}
Note that $k\in T(f)$.
If $T(f)=\Lambda\setminus I$,
$f$ is a section of 
$\lefttop{I}\psizero_{q_I(\vecu(0))}
\lefttop{\Lambda}\Vtildezero_{\veca}(\nbigm)$,
and we are done.
Suppose that $T(f)\neq\Lambda\setminus I$.
We take $j\in \Lambda\setminus (I\cup T(f))$.
If $u\neq q_j(\vecu(0))$,
$[f]_{j,u}$ is a section of
\[
 \lefttop{T(f)}\Vzero_{<q_{T(f)}(\veca)}
 \lefttop{j}\psizero_{u}
\lefttop{I}\psizero_{q_{I}(\vecu(0))}
\lefttop{\Lambda}V^{(\lambda_0)}_{\veca}(\nbigm).
\]
By the claim of the proposition in the lower dimensional case,
there exist
 $h\in
 \lefttop{j}\psizero_{q_j(\vecu(0))}
\lefttop{I}\psizero_{q_{I}(\vecu(0))}
\lefttop{\Lambda}\Vtilde^{(\lambda_0)}_{\veca}(\nbigm)$
such that the following holds
for any $i\in T(f)$.
\begin{itemize}
 \item $h$ and $[f]_{j,q_j(\vecu(0))}$
       induce the same local section of
       \[
       \lefttop{i}\psizero_{q_i(\vecu(0))}
       \lefttop{j}\psizero_{q_j(\vecu(0))}
       \lefttop{I}\psizero_{q_{I}(\vecu(0))}
       \lefttop{\Lambda}\Vtilde^{(\lambda_0)}_{\veca}(\nbigm)
       \subset
       \lefttop{i}\Gr^{\Vzero}_{a_i}
       \lefttop{j}\psizero_{q_j(\vecu(0))}
       \lefttop{I}\psizero_{q_{I}(\vecu(0))}
       \lefttop{\Lambda}V^{(\lambda_0)}_{\veca}(\nbigm).
       \]
\end{itemize}
We set $\vecu'(0)=\vecu(0)-\epsilon\vecdelta_{T(f)}$
and $\veca'=\paramap(\lambda_0,\vecu'(0))=\veca-\epsilon\veciti_{T(f)}$.
There exists a section $g_j$
of
$\lefttop{I}\psizero_{q_{I}(\vecu(0))}
\lefttop{\Lambda}V^{(\lambda_0)}_{\veca'}(\nbigm)$
which induces
$\sum_{u_j\neq q_j(\vecu(0))}[f]_{j,u_j}
+([f]_{j,u(0)_j}-h)$.
Then, $f-g_j$ induces $s$,
and $T(f-g_j)=T(f)\cup\{j\}$.
Hence, by an induction, we obtain the claim
of the proposition in the case $|L|=1$.

Let us consider the general case.
Let $f\in \lefttop{I}\psizero_{q_I(\vecu(0))}
\lefttop{\Lambda}\Vzero_{\veca}(\nbigm)$
such that
$\rho_{0,L}(f)$ is contained in
\[
\bigoplus_{k\in L}
\lefttop{I\cup \{k\}}\psizero_{q_{I\sqcup \{k\}}(\vecu(0))}
\lefttop{\Lambda}\Vtilde^{(\lambda_0)}_{\veca}(\nbigm).
\]
We have the decomposition
$\rho_{0,L}(f)=\sum_{k\in L}[f]_k$.
Take $k(1)\in L$.
There exists
$g_1\in \lefttop{I}\psizero_{q_I(\vecu(0))}
\lefttop{\Lambda}\Vtilde^{(\lambda_0)}_{\veca}(\nbigm)$
such that
$\rho_{1,k(1)}(g_1)=\rho_{1,k(1)}(f)$.
We set
$\vecu'=\vecu-\epsilon\vecdelta_{k(1)}$
and $\veca'=\veca-\epsilon\veciti_{k(1)}$.
We also set $L'=L\setminus\{k(1)\}$.
We obtain
$f-g_1\in\lefttop{I}\psizero_{q_I(\vecu(0))}
\lefttop{\Lambda}\Vzero_{\veca'}(\nbigm)$,
which satisfies
\[
 \rho_{1,L'}(f-g_1)
 \in
\bigoplus_{k\in L'}
\lefttop{I\cup \{k\}}\psizero_{q_{I\sqcup \{k\}}(\vecu(0))}
\lefttop{\Lambda}\Vtilde^{(\lambda_0)}_{\veca'}(\nbigm).
\]
By the assumption of the induction,
there exists
$g_2\in \lefttop{I}\psizero_{q_I(\vecu(0))}
\lefttop{\Lambda}\Vtilde^{(\lambda_0)}_{\veca'}(\nbigm)$
such that
$\rho_{1,L'}(g_2)=\rho_{1,L'}(f-g_1)$.
Then, we obtain
$f-(g_1+g_2)\in
\lefttop{I}\psizero_{q_I(\vecu(0))}
\lefttop{\Lambda}\Vtilde_{\veca}(\nbigm)$
which satisfies
$\rho_{1,L}(f-(g_1+g_2))=\rho_{0,L}(f)$.
\hfill\qed

\begin{cor}
For $K\subset \Lambda\setminus I$,
the image of the natural morphism
\[
 \lefttop{I}\psizero_{q_I(\vecu(0))}
\lefttop{\Lambda}\Vtilde^{(\lambda_0)}_{\veca}(\nbigm)
 \lrarr
 \lefttop{K}\Gr^{\Vzero}_{q_{K}(\veca)}
 \lefttop{I}\psizero_{q_I(\vecu(0))}
 \lefttop{\Lambda}\Vzero_{\veca}(\nbigm)
\]
 is
$\lefttop{I\cup K}\psizero_{q_{I\sqcup K}(\vecu(0))}
 \lefttop{\Lambda}\Vtilde^{(\lambda_0)}_{\veca}(\nbigm)$.
\hfill\qed
\end{cor}

For $\ell\in\seisuu_{\geq 0}$,
we have the following morphism:
\[
 \lefttop{I}\psizero_{q_I(\vecu(0))}
 \lefttop{\Lambda}\Vtildezero_{\veca}(\nbigm)
 \lrarr
 \bigoplus_{\substack{
 K\supset I\\
 |K\setminus I|=\ell+1
 }}
 \lefttop{K}\psizero_{q_K(\vecu(0))}
 \lefttop{\Lambda}\Vtildezero_{\veca}(\nbigm).
\]
Let
$F_{\ell}
\lefttop{I}\psizero_{q_I(\vecu(0))}
\lefttop{\Lambda}\Vtildezero_{\veca}(\nbigm)$
denote the kernel.
Namely,
\[
 F_{\ell}
 \lefttop{I}\psizero_{q_I(\vecu(0))}
 \lefttop{\Lambda}\Vtildezero_{\veca}(\nbigm)
 =
  \lefttop{I}\psizero_{q_I(\vecu(0))}
  \lefttop{\Lambda}\Vtildezero_{\veca}(\nbigm)
  \cap
 F_{\ell}
 \lefttop{I}\Gr^{\Vzero}_{q_I(\veca)}
 \lefttop{\Lambda}\Vzero_{\veca}(\nbigm).
\]
Thus, we obtain an increasing filtration $F_{\bullet}$.
We also set
\[
 \lefttop{\Lambda\setminus I}\Vzero_{<q_{\Lambda\setminus I}(\veca)}
 \lefttop{I}\psizero_{q_I(\vecu(0))}
 \lefttop{\Lambda}\Vtildezero_{\veca}(\nbigm):=
F_0
  \lefttop{I}\psizero_{q_I(\vecu(0))}
 \lefttop{\Lambda}\Vtildezero_{\veca}(\nbigm).
\]

Because the kernel of the morphism
\[
 F_{\ell}\lefttop{I}\psizero_{q_I(\vecu(0))}
 \lefttop{\Lambda}\Vtildezero_{\veca}(\nbigm)
 \lrarr
 \bigoplus_{\substack{
 K\supset I\\
 |K\setminus I|=\ell
 }}
 \lefttop{K}\psizero_{q_K(\vecu(0))}
 \lefttop{\Lambda}\Vtildezero_{\veca}(\nbigm)
\]
is
$F_{\ell-1}\lefttop{I}\psizero_{q_I(\vecu(0))}
 \lefttop{\Lambda}\Vtildezero_{\veca}(\nbigm)$,
we obtain the following monomorphism
 \begin{equation}
\label{eq;22.2.16.4}
 \Gr^F_{\ell}\lefttop{I}\psizero_{q_I(\vecu(0))}
 \lefttop{\Lambda}\Vtildezero_{\veca}(\nbigm)
 \lrarr
 \bigoplus_{\substack{
 K\supset I\\
 |K\setminus I|=\ell
 }}
 \lefttop{K}\psizero_{q_I(\vecu(0))}
 \lefttop{\Lambda}\Vtildezero_{\veca}(\nbigm).
\end{equation}

\begin{lem}
The morphism {\rm(\ref{eq;22.2.16.4})} induces
\begin{equation}
\label{eq;22.2.16.5}
 \Gr^F_{\ell}
  \lefttop{I}\psizero_{q_I(\vecu(0))}
  \lefttop{\Lambda}\Vtildezero_{\veca}(\nbigm)
  \simeq
  \bigoplus_{\substack{K\supset I\\
  |K\setminus I|=\ell
  }}
  \lefttop{\Lambda\setminus K}
  \Vzero_{<q_{\Lambda\setminus K}(\veca)}
  \lefttop{K}\psizero_{q_K(\vecu(0))}
  \lefttop{\Lambda}\Vtildezero_{\veca}(\nbigm).
\end{equation}
\end{lem}
\pf
Because the image of (\ref{eq;22.2.16.4})
is contained in the right hand side of (\ref{eq;22.2.16.5}),
we obtain the morphism (\ref{eq;22.2.16.5}).
It is enough to prove that (\ref{eq;22.2.16.5}) is an epimorphism.
We have only to prove that
the following morphism is an epimorphism
for any $K\supset I$:
\[
 \lefttop{\Lambda\setminus K}\Vzero_{<q_{\Lambda\setminus K}(\veca)}
  \lefttop{I}\psizero_{q_I(\vecu(0))}
 \lefttop{\Lambda}\Vtildezero_{\veca}(\nbigm)
 \lrarr
 \lefttop{\Lambda\setminus K}\Vzero_{<q_{\Lambda\setminus K}(\veca)}
 \lefttop{K}\psizero_{q_K(\vecu(0))}
 \lefttop{\Lambda}\Vtildezero_{\veca}(\nbigm).
\]
The case $|K\setminus I|\leq 1$ are easy.

For any $L\supset I$, there exists the following morphism
\[
\Phi_{L,I}:
 \lefttop{I}\psizero_{q_I(\vecu(0))}
 \lefttop{\Lambda}\Vzero_{\veca}(\nbigm)
 \lrarr
 \bigoplus_{\substack{
 \paramap(\lambda_0,\vecu)=q_L(\veca)
 \\
 q_I(\vecu)=q_I(\vecu(0))
 }}
 \lefttop{L}\psizero_{\vecu}
 \lefttop{\Lambda}\Vzero_{\veca}(\nbigm).
\]
The following induced morphism is an epimorphism:
\begin{equation}
\label{eq;22.2.16.8}
 \lefttop{\Lambda\setminus L}\Vzero_{<q_{\Lambda\setminus L}(\veca)}
  \lefttop{I}\psizero_{q_I(\vecu(0))}
 \lefttop{\Lambda}\Vzero_{\veca}(\nbigm)
 \lrarr
 \bigoplus_{\substack{
 \paramap(\lambda_0,\vecu)=q_L(\veca)
 \\
 q_I(\vecu)=q_I(\vecu(0))
 }}
  \lefttop{\Lambda\setminus L}\Vzero_{<q_{\Lambda\setminus L}(\veca)}
 \lefttop{L}\psizero_{\vecu}
 \lefttop{\Lambda}\Vzero_{\veca}(\nbigm).
\end{equation}

For a local section $s$ of
$\lefttop{\Lambda\setminus K}\Vzero_{<q_{\Lambda\setminus K}(\veca)}
 \lefttop{K}\psizero_{q_K(\vecu(0))}
 \lefttop{\Lambda}\Vtildezero_{\veca}(\nbigm)$,
there exists a section
$f_0$ of
\begin{equation}
\label{eq;22.2.16.7}
\lefttop{\Lambda\setminus K}\Vzero_{<q_{\Lambda\setminus K}(\veca)}
  \lefttop{I}\psizero_{q_I(\vecu(0))}
  \lefttop{\Lambda}\Vzero_{\veca}(\nbigm)
\end{equation}
such that $\Phi_{K,I}(f_0)=s$.
Let us construct sections $f_i$ $(i=1,2,\ldots,|K\setminus I|)$
of (\ref{eq;22.2.16.7})
satisfying the following conditions.
\begin{itemize}
 \item $\Phi_{K,I}(f_i)=s$.
 \item Let $I\subset L\subset K$ such that
       $|K\setminus L|\leq i$.
       Then, $\Phi_{L,I}(f_i)$
       is contained in
$\lefttop{\Lambda\setminus L}\Vzero_{<q_{\Lambda\setminus L}(\veca)}
 \lefttop{L}\psizero_{q_I(\vecu(0))}
 \lefttop{\Lambda}\Vzero_{\veca}(\nbigm)$.
\end{itemize}
Suppose that we have already obtained $f_{i-1}$.
For $I\subset L\subset K$ such that $|K\setminus L|=i$,
we have the decomposition
\[
 \Phi_{L,I}(f_{i-1})
 =\sum_{
 \substack{\paramap(\lambda_0,\vecu)=q_L(\veca)\\
 q_I(\vecu)=q_I(\vecu(0))
  }}
  \Phi_{L,I}(f_{i-1})_{\vecu}
  \in
  \bigoplus_{
 \substack{\paramap(\lambda_0,\vecu)=q_L(\veca)\\
 q_I(\vecu)=q_I(\vecu(0))
  }}
  \lefttop{L}\psizero_{\vecu}
  \lefttop{\Lambda}\Vzero_{\veca}(\nbigm).
\]
By the condition for $f_{i-1}$,
if $\vecu\neq q_L(\vecu(0))$,
$\Phi_{L,I}(f_{i-1})_{\vecu}$
is a section of
$\lefttop{\Lambda\setminus L}
 \Vzero_{<q_{\Lambda\setminus L}(\veca)}
 \lefttop{L}\psizero_{\vecu}
 \lefttop{\Lambda}\Vzero_{\veca}(\nbigm)$.
Because the morphisms (\ref{eq;22.2.16.8}) are epimorphisms,
there exists a section $h_L$ of
$\lefttop{\Lambda\setminus L}
 \Vzero_{<q_{\Lambda\setminus L}(\veca)}
 \lefttop{I}\psizero_{q_I(\vecu(0))}
 \lefttop{\Lambda}\Vzero_{\veca}(\nbigm)$
 such that
\[
 \Phi_{L,I}(h_L)
 =\sum_{\vecu\neq q_L(\vecu(0))}
 \Phi_{L,I}(f_{i-1})_{\vecu}.
\] 
By setting $f_i=f_{i-1}-\sum_Lh_L$,
the inductive construction can proceed.
Then, $f=f_{|K\setminus I|}$ is a section of
$\lefttop{\Lambda\setminus K}\Vzero_{<q_{\Lambda\setminus K}(\veca)}
  \lefttop{I}\psizero_{q_I(\vecu(0))}
  \lefttop{\Lambda}\Vtildezero_{\veca}(\nbigm)$
such that $\Phi_{K,I}(f)=s$.
\hfill\qed

\begin{prop}
\label{prop;22.2.16.31}
Let $f:\nbigm_1\to \nbigm_2$ be a morphism of $\nbigr_X$-modules
underlying a morphism of good mixed twistor $\nbigd$-modules
on $(X,H)$.
For the induced morphism
\begin{equation}
 \lefttop{I}\psizero_{q_I(\vecu(0))}
 \lefttop{\Lambda}\Vtildezero_{\veca}(f):
 \lefttop{I}\psizero_{q_I(\vecu(0))}
 \lefttop{\Lambda}\Vtildezero_{\veca}(\nbigm_1)
 \lrarr
 \lefttop{I}\psizero_{q_I(\vecu(0))}
 \lefttop{\Lambda}\Vtildezero_{\veca}(\nbigm_2),
\end{equation}
there exist the following natural isomorphisms:
\begin{equation}
\label{eq;22.2.16.20}
 \lefttop{I}\psizero_{q_I(\vecu(0))}
 \lefttop{\Lambda}\Vtildezero_{\veca}(\Ker f)
 \simeq
\Ker\lefttop{I}\psizero_{q_I(\vecu(0))}
 \lefttop{\Lambda}\Vtildezero_{\veca}(f),
\end{equation}
\begin{equation}
\label{eq;22.2.16.21}
 \lefttop{I}\psizero_{q_I(\vecu(0))}
 \lefttop{\Lambda}\Vtildezero_{\veca}(\Image f)
 \simeq
\Image\lefttop{I}\psizero_{q_I(\vecu(0))}
 \lefttop{\Lambda}\Vtildezero_{\veca}(f),
\end{equation}
\begin{equation}
\label{eq;22.2.16.22}
 \lefttop{I}\psizero_{q_I(\vecu(0))}
 \lefttop{\Lambda}\Vtildezero_{\veca}(\Cok f)
 \simeq
\Cok\lefttop{I}\psizero_{q_I(\vecu(0))}
 \lefttop{\Lambda}\Vtildezero_{\veca}(f).
 \end{equation}
Moreover,
$\lefttop{I}\psizero_{q_I(\vecu(0))}
 \lefttop{\Lambda}\Vtildezero_{\veca}(f)$
is strict with respect to the filtration $F_{\bullet}$.
\end{prop}
\pf
We use an induction on $\dim X$,
We assume that the claims are already proved
in the lower dimensional case.
By using Lemma \ref{lem;22.2.16.30},
we easily obtain (\ref{eq;22.2.16.20}).
We also obtain
\[
  \lefttop{I}\psizero_{q_I(\vecu(0))}
  \lefttop{\Lambda}\Vtildezero_{\veca}(\Image f)
  =\Image\Bigl(
  \lefttop{I}\psizero_{q_I(\vecu(0))}
  \lefttop{\Lambda}\Vzero_{\veca}(f)
  \Bigr)
\cap
\lefttop{I}\psizero_{q_I(\vecu(0))}
\lefttop{\Lambda}\Vtildezero_{\veca}(\nbigm_2).
\]
Hence, (\ref{eq;22.2.16.21}) is equal to
\begin{equation}
\label{eq;22.3.14.10}
 \Image\Bigl(
 \lefttop{I}\psizero_{q_I(\vecu(0))}
\lefttop{\Lambda}\Vzero_{\veca}(f)\Bigr)
\cap
\lefttop{I}\psizero_{q_I(\vecu(0))}
\lefttop{\Lambda}\Vtildezero_{\veca}(\nbigm_2)
=
 \Image
 \lefttop{I}\psizero_{q_I(\vecu(0))}
\lefttop{\Lambda}\Vtildezero_{\veca}(f).
\end{equation}
Let us prove (\ref{eq;22.3.14.10}).
The claim $\supset$ is clear.
Let $s$ be a section of a section of the left hand side of
(\ref{eq;22.3.14.10}).
For $j\in \Lambda\setminus I$,
we have the induced sections
$[s]_j$ of
\[
\Image\Bigl(
\lefttop{j}\psizero_{q_j(\vecu(0))}
\lefttop{I}\psizero_{q_I(\vecu(0))}
\Vzero_{\veca}(f)
\Bigr)
\cap
\lefttop{j}\psizero_{q_j(\vecu(0))}
\lefttop{I}\psizero_{q_I(\vecu(0))}
\Vtildezero_{\veca}(\nbigm_2).
\]
Let $T_1(s)$ denote the set of $j\in \Lambda\setminus I$
such that $[s]_j\neq 0$.
We prove that $s$ is contained in
the right hand side of (\ref{eq;22.3.14.10})
by using an induction on $|T_1(s)|$.
Let us consider the case $|T_1(s)|=0$.
Then, by using Lemma \ref{lem;22.2.16.23},
we can construct
$t\in \lefttop{I}\psizero_{q_I(\vecu(0))}
\Vzero_{\veca}(\nbigm_1)$
such that
(i) $\lefttop{I}\psizero_{q_I(\vecu(0))}
\Vzero_{\veca}(f)(t)=s$,
(ii)
the induced sections
$[t]_j$ of
$\lefttop{j}\Gr^{\Vzero}_{a_j}
\lefttop{I}\psizero_{q_I(\vecu(0))}
\Vzero_{\veca}(\nbigm_1)$
are $0$.
The second condition implies
$t\in \lefttop{I}\psizero_{q_I(\vecu(0))}
\Vtildezero_{\veca}(\nbigm_1)$,
and hence we are done in this case.
Let us consider the case $|T_1(s)|>0$.
Take $j(1)\in T_1(s)$.
By the assumption of the induction on the dimension of $X$,
$[s]_{j(1)}$ is a section of
\[
 \Image\Bigl(
 \lefttop{j(1)}\psizero_{q_{j(1)}(\vecu(0))}
\lefttop{I}\psizero_{q_I(\vecu(0))}
\Vtildezero_{\veca}(f)
\Bigr).
\]
There exists a section $t_1$ of 
$\lefttop{j(1)}\psizero_{q_{j(1)}(\vecu(0))}
\lefttop{I}\psizero_{q_I(\vecu(0))}
\Vtildezero_{\veca}(\nbigm_1)$
such that
\[
\lefttop{j(1)}\psizero_{q_{j(1)}(\vecu(0))}
\lefttop{I}\psizero_{q_I(\vecu(0))}
\Vtildezero_{\veca}(f)(t_1)=[s]_{j(1)}.
\]
For any $k\in\Lambda\setminus (I\cup\{j(1)\})$,
we obtain the induced sections
$[t_1]_k$
of 
$\lefttop{k}\psizero_{q_k(\vecu(0))}
\lefttop{j(1)}\psizero_{q_{j(1)}(\vecu(0))}
\lefttop{I}\psizero_{q_I(\vecu(0))}
\Vtildezero_{\veca}(\nbigm_1)$.
By using (\ref{eq;22.2.16.20}),
we may assume that $[t_1]_k=0$
for any $k\in \Lambda\setminus(I\cup T_1(s))$.
We set
$\vecu'(0)=\vecu(0)-\epsilon\vecdelta_{\Lambda\setminus(I\cup T_1(s))}$
and
$\veca'=\veca-\epsilon\veciti_{\Lambda\setminus(I\cup T_1(s))}$
for any sufficiently small $\epsilon>0$.
There exists
$\widetilde{t}_1$
of 
$\lefttop{I}\psizero_{q_I(\vecu'(0))}
\Vtildezero_{\veca'}(\nbigm_1)$
which induces $t_1$.
We obtain
\[
 T_1\Bigl(
s- 
\lefttop{I}\psizero_{q_I(\vecu'(0))}
\Vtildezero_{\veca}(f)(\widetilde{t}_1)
 \Bigr)
=T_1(s)\setminus\{j(1)\}
\]
By the assumption on $|T_1(s)|$,
we obtain that
$s-\lefttop{I}\psizero_{q_I\vecu'(0)}
\Vtildezero_{\veca'}(f)(\widetilde{t}_1)$
is contained in
$\Image
\lefttop{I}\psizero_{q_I(\vecu(0))}
\Vtildezero_{\veca}(f)$.
Hence, the induction can proceed,
and we obtain (\ref{eq;22.3.14.10})
and (\ref{eq;22.2.16.21}).

\vspace{.1in}
Let us study (\ref{eq;22.2.16.22}).
There exists a natural morphism
\begin{equation}
\label{eq;22.3.14.11}
 \Cok\Bigl(
 \lefttop{I}\psizero_{q_I(\vecu(0))}
 \lefttop{\Lambda}\Vtildezero_{\veca}(f)
 \Bigr)
\lrarr
  \lefttop{I}\psizero_{q_I(\vecu(0))}
 \lefttop{\Lambda}\Vtildezero_{\veca}(\Cok f).
\end{equation}
By (\ref{eq;22.3.14.10}), (\ref{eq;22.3.14.11}) is a monomorphism.
To prove that (\ref{eq;22.3.14.11}) is an epimorphism,
it is enough to apply (\ref{eq;22.2.16.20})
for the morphism $\nbigm_2\to \Cok(f)$.

By Lemma \ref{lem;22.2.16.30},
we obtain the strictness with respect to the filtration $F_{\bullet}$
from (\ref{eq;22.3.14.10}).
\hfill\qed

\subsubsection{Weight filtrations}

Let $I\subset\Lambda$.
We recall that $\nbigs(I)$ denotes
the set of increasing sequences in $I$
(see \S\ref{subsection;22.4.26.1}).
For $K\subset I$,
we have the weight filtration
$\lefttop{K}\nbigw$
of 
$\lefttop{I}\psizero_{q_I(\vecu(0))}
\lefttop{\Lambda}\Vtilde^{(\lambda_0)}_{\veca}(\nbigm)$
with respect to
$\sum_{j\in K}\nbign_j$.
For $\vecK=(K_1,\ldots,K_m)\in\nbigs(I)$
and $\vecp\in\seisuu^m$,
we set
$\lefttop{\vecK}
 \nbigw_{\vecp}:=
 \bigcap_{j=1}^m
 \lefttop{K_j}\nbigw_{p_j}$
on  
$\lefttop{I}\psizero_{q_I(\vecu(0))}
\lefttop{\Lambda}\Vtilde^{(\lambda_0)}_{\veca}(\nbigm)$.

\begin{lem}
For $I\subset L$,
the morphism
$\lefttop{I}\psizero_{q_I(\vecu(0))}
 \lefttop{\Lambda}\Vtilde^{(\lambda_0)}_{\veca}(\nbigm)
 \lrarr
 \lefttop{L}\psizero_{q_L(\vecu(0))}
 \lefttop{\Lambda}\Vtilde^{(\lambda_0)}_{\veca}(\nbigm)$
 is strict with respect to $\lefttop{\vecK}\nbigw$
 for any $\vecK\in\nbigs(I)$,
 i.e.,
 $\lefttop{\vecK}\nbigw_{\vecp}
 \lefttop{L}\psizero_{q_L(\vecu(0))}
\lefttop{\Lambda}\Vtilde^{(\lambda_0)}_{\veca}(\nbigm)$
 is equal to the image of
$\lefttop{\vecK}\nbigw_{\vecp}
\lefttop{I}\psizero_{q_I(\vecu(0))}
\lefttop{\Lambda}\Vtilde^{(\lambda_0)}_{\veca}(\nbigm)$.
\end{lem}
\pf
It follows from
Lemma \ref{lem;22.4.26.2},
Proposition \ref{prop;22.2.16.31}
and the formula of the weight filtration \cite[(1.2.3)]{k3}.
\hfill\qed

\vspace{.1in}
Let us define
$\lefttop{\vecK}\nbigw_{\vecp}$
on 
$\lefttop{I}\psizero_{q_I(\vecu(0))}
\lefttop{\Lambda}\Vtilde^{(\lambda_0)}_{\veca}(\nbigm)$
for $\vecK=(K_1,\ldots,K_m)\in\nbigs(\Lambda)$,
which is not necessarily contained in $\nbigs(I)$,
as the inverse image of
\[
 \bigoplus_{i=1}^{m}
 \lefttop{K_i}\nbigw_{p_i}
 \lefttop{I\cup K_i}\psizero_{q_{I\cup K_i}(\vecu(0))}
 \lefttop{\Lambda}\Vtildezero_{\veca}(\nbigm)
\]
by the natural morphism
\[
 \lefttop{I}\psizero_{q_I(\vecu(0))}
 \lefttop{\Lambda}\Vtildezero_{\veca}(\nbigm)
 \lrarr
 \bigoplus_{i=1}^{m}
 \lefttop{I\cup K_i}\psizero_{q_{I\cup K_i}(\vecu(0))}
 \lefttop{\Lambda}\Vtildezero_{\veca}(\nbigm).
\] 

If $K_m\not\subset I$,
there exists $0\leq i(I)<m$
such that
$K_{i(I)}\subset I$ and $K_{i(I)+1}\not\subset I$.
If $K_m\subset I$,
we set $i(I)=m$.
For any $1\leq i\leq m$,
we set
$\vecK_{\leq i}:=(K_1,\ldots,K_i)$
and $\vecp_{\leq i}=(p_1,\ldots,p_i)$.
Let $\vecK_{\leq 0}$ denote the sequence of the length $0$,
and let $\vecp_{\leq 0}$ denote the empty tuple.
We have the following morphism
\begin{equation}
\label{eq;22.2.17.20}
 \lefttop{\vecK_{\leq i(I)}}
 \nbigw_{\vecp_{\leq i(I)}}
 \lefttop{I}\psizero_{q_I(\vecu(0))}
 \lefttop{\Lambda}\Vtilde^{(\lambda_0)}_{\veca}
 (\nbigm)
 \lrarr
 \bigoplus_{I\subset L}
 \lefttop{\vecK_{\leq i(I)}}\nbigw_{\vecp_{\leq i(I)}}
 \lefttop{L}
 \psizero_{q_{L}(\vecu(0))}
 \lefttop{\Lambda}\Vtilde^{(\lambda_0)}_{\veca}(\nbigm).
\end{equation}
The following lemma is clear by the construction.
\begin{lem}
$\lefttop{\vecK}\nbigw_{\vecp}
\lefttop{I}\psizero_{q_I(\vecu(0))}
\lefttop{\Lambda}\Vtilde^{(\lambda_0)}_{\veca}
(\nbigm)$
is equal to the inverse image of
\[
 \bigoplus_{I\subset L}
 \lefttop{\vecK_{\leq i(L)}}\nbigw_{\vecp_{\leq i(L)}}
 \lefttop{L}
 \psizero_{q_{L}(\vecu(0))}
 \lefttop{\Lambda}\Vtilde^{(\lambda_0)}_{\veca}(\nbigm)
\]
by {\rm(\ref{eq;22.2.17.20})}.
\hfill\qed
\end{lem}

We also have the following recursive expression of
 $\lefttop{\vecK}\nbigw_{\vecp}
 \lefttop{I}\psizero_{q_I(\vecu(0))}
 \lefttop{\Lambda}\Vtildezero_{\veca}(\nbigm)$.
\begin{lem}
\label{lem;22.2.17.11}
 $\lefttop{\vecK}\nbigw_{\vecp}
 \lefttop{I}\psizero_{q_I(\vecu(0))}
 \lefttop{\Lambda}\Vtildezero_{\veca}(\nbigm)$
is equal to the inverse image of
\[
  \bigoplus_{\substack{I\subset L\\ |L\setminus I|=1}}
  \lefttop{\vecK}\nbigw_{\vecp}
 \lefttop{L}\psizero_{q_L(\vecu(0))}
 \lefttop{\Lambda}\Vtildezero_{\veca}(\nbigm)
\]
by the following natural morphism
\begin{equation}
\label{eq;22.2.17.10}
 \lefttop{\vecK_{\leq i(I)}}\nbigw_{\vecp_{\leq i(I)}}
 \lefttop{I}\psizero_{q_I(\vecu(0))}
 \lefttop{\Lambda}\Vtildezero_{\veca}(\nbigm)
 \lrarr
 \bigoplus_{\substack{I\subset L\\ |L\setminus I|=1}}
  \lefttop{\vecK_{\leq i(I)}}\nbigw_{\vecp_{\leq i(I)}}
 \lefttop{L}\psizero_{q_L(\vecu(0))}
 \lefttop{\Lambda}\Vtildezero_{\veca}(\nbigm).
\end{equation}
\hfill\qed
\end{lem}

\begin{lem}
For $I\subset L$,
the image of the morphism
\begin{equation}
\label{eq;22.2.17.1}
 \lefttop{\vecK}\nbigw_{\vecp}
 \lefttop{I}\psizero_{q_I(\vecu(0))}
 \Vtilde^{(\lambda_0)}_{\veca}
 (\nbigm)
\lrarr
 \lefttop{L}\psizero_{q_L(\vecu(0))}
 \Vtilde^{(\lambda_0)}_{\veca}
 (\nbigm)
\end{equation}
is equal to
$\lefttop{\vecK}\nbigw_{\vecp}
 \lefttop{L}\psizero_{q_L(\vecu(0))}
 \Vtilde^{(\lambda_0)}_{\veca}
 (\nbigm)$.
\end{lem}
\pf
The image of (\ref{eq;22.2.17.1}) is clearly contained in
$\lefttop{\vecK}\nbigw_{\vecp}
 \lefttop{L}\psizero_{q_L(\vecu(0))}
 \Vtilde^{(\lambda_0)}_{\veca}
 (\nbigm)$.
We have only to consider the case $m=|\Lambda|$.
It is enough to consider the case $|L\setminus I|=1$.
We use a descending induction on $|L|$.
If $L=\Lambda$,
we obtain the claim by using Lemma \ref{lem;22.2.17.11}.

We note that $i(I)\leq i(L)\leq i(I)+1$
under the assumption $m=|\Lambda|$.
Let us consider the case $i(L)=i(I)+1$.
Let
$\lefttop{\vecK}\nbigw'_{\vecp}
\lefttop{I}\psizero_{q_I(\vecu(0))}
\lefttop{\Lambda}\Vtildezero_{\veca}(\nbigm)$
denote the inverse image of
$\lefttop{\vecK}\nbigw_{\vecp}
 \lefttop{L}\psizero_{q_L(\vecu(0))}
 \lefttop{\Lambda}\Vtildezero_{\veca}(\nbigm)$
by the natural morphism
\[
\lefttop{\vecK_{\leq i(I)}}\nbigw_{\vecp_{\leq i(I)}}
 \lefttop{I}\psizero_{q_I(\vecu(0))}
 \lefttop{\Lambda}\Vtildezero_{\veca}(\nbigm)
 \lrarr
  \lefttop{\vecK_{\leq i(I)}}\nbigw_{\vecp_{\leq i(I)}}
 \lefttop{L}\psizero_{q_L(\vecu(0))}
 \lefttop{\Lambda}\Vtildezero_{\veca}(\nbigm).
\]
By the construction, we have
$\lefttop{\vecK}\nbigw_{\vecp}
\lefttop{I}\psizero_{q_I(\vecu(0))}
\lefttop{\Lambda}\Vtildezero_{\veca}(\nbigm)
\subset
\lefttop{\vecK}\nbigw'_{\vecp}
\lefttop{I}\psizero_{q_I(\vecu(0))}
\lefttop{\Lambda}\Vtildezero_{\veca}(\nbigm)$.
For any $L'$ such that $I\subset L'$ and $|L'\setminus I|=1$,
we have $i(L')=i(I)$.
We also have
$i(L\cup L')=i(L)=i(I)+1$.
Hence,
the image of 
$\lefttop{\vecK}\nbigw'_{\vecp}
\lefttop{I}\psizero_{q_I(\vecu(0))}
\lefttop{\Lambda}\Vtildezero_{\veca}(\nbigm)$
via the morphism
\[
 \lefttop{\vecK_{\leq i(I)}}\nbigw_{\vecp_{\leq i(I)}}
 \lefttop{I}\psizero_{q_I(\vecu(0))}
 \lefttop{\Lambda}\Vtildezero_{\veca}(\nbigm)
 \lrarr
  \lefttop{\vecK_{\leq i(I)}}\nbigw_{\vecp_{\leq i(I)}}
 \lefttop{L'}\psizero_{q_{L'}(\vecu(0))}
 \lefttop{\Lambda}\Vtildezero_{\veca}(\nbigm)
\]
is contained in
$\lefttop{\vecK}\nbigw_{\vecp}
\lefttop{L'}\psizero_{q_{L'}(\vecu(0))}
\lefttop{\Lambda}\Vtildezero_{\veca}(\nbigm)$.
Then, by Lemma \ref{lem;22.2.17.11},
we obtain that
$\lefttop{\vecK}\nbigw_{\vecp}
\lefttop{I}\psizero_{q_I(\vecu(0))}
\lefttop{\Lambda}\Vtildezero_{\veca}(\nbigm)
=
\lefttop{\vecK}\nbigw'_{\vecp}
\lefttop{I}\psizero_{q_I(\vecu(0))}
\lefttop{\Lambda}\Vtildezero_{\veca}(\nbigm)$.

Let us consider the case $i(L)=i(I)$.
Let $j\in\Lambda\setminus I$ be determined by
$K_{i(I)+1}=K_{i(I)}\sqcup\{j\}$.
We set $\Itilde:=I\sqcup\{j\}$
and $\Ltilde=L\sqcup\{j\}$.
We have $i(\Itilde)=i(\Ltilde)=i(I)+1$.
Let $s$ be a section of
$\lefttop{\vecK}\nbigw_{\vecp}
\lefttop{L}\psizero_{q_L(\vecu(0))}\Vtildezero_{\veca}(\nbigm)$.
By the assumption of the induction,
the claim of the lemma holds for $L\subset\Ltilde$
and $\Itilde\subset\Ltilde$.
By the previous consideration,
the claim of the lemma holds for $I\subset\Itilde$.
Therefore, there exists a section $f_1$ of 
$\lefttop{\vecK}\nbigw_{\vecp}
\lefttop{I}\psizero_{q_I(\vecu(0))}\Vtildezero_{\veca}(\nbigm)$
such that
$f_1$ and $s$ induce the same section of
$\lefttop{\vecK}\nbigw_{\vecp}
\lefttop{\Ltilde}\psizero_{q_{\Ltilde}(\vecu(0))}
\Vtildezero_{\veca}(\nbigm)$.
Let $[f_1]$ denote the induced section of
$\lefttop{\vecK}\nbigw_{\vecp}
\lefttop{L}\psizero_{q_{L}(\vecu(0))}
\Vtildezero_{\veca}(\nbigm)$.
Then, $s-[f_1]$ is a section of
\[
\Bigl(
\lefttop{\vecK}\nbigw_{\vecp}
\lefttop{L}\psizero_{q_{L}(\vecu(0))}
\Vtildezero_{\veca}(\nbigm)
\Bigr)
\cap
\lefttop{L}\psizero_{q_{L}(\vecu(0))}
\Vtildezero_{\veca-\epsilon\vecdelta_j}(\nbigm)
=
\lefttop{\vecK_{\leq i(L)}}\nbigw_{\vecp}
\lefttop{L}\psizero_{q_{L}(\vecu(0))}
\Vtildezero_{\veca-\epsilon\vecdelta_j}(\nbigm).
\]
Then, there exists a section $f_2$ of
$\lefttop{\vecK_{\leq i(I)}}\nbigw_{\vecp}
\lefttop{I}\psizero_{q_{I}(\vecu(0))}
\Vtildezero_{\veca-\epsilon\vecdelta_j}(\nbigm)$
which induces $s-[f_1]$.
Then, the section $f_1+f_2$ of
$\lefttop{\vecK}\nbigw_{\vecp}
\lefttop{I}\psizero_{q_{I}(\vecu(0))}
\Vtildezero_{\veca}(\nbigm)$
induces $s$.
\hfill\qed

\vspace{.1in}
For $\ell\in\seisuu_{\geq 0}$,
we obtain the following morphism:
\begin{equation}
\label{eq;22.2.17.2}
 \lefttop{\vecK}\nbigw_{\vecp}
\lefttop{I}\psizero_{q_{I}(\vecu(0))}
\Vtildezero_{\veca}(\nbigm)
\lrarr
\bigoplus_{\substack{I\subset L\\ |L\setminus I|=\ell+1}}
 \lefttop{\vecK}\nbigw_{\vecp}
\lefttop{L}\psizero_{q_{L}(\vecu(0))}
\Vtildezero_{\veca}(\nbigm).
\end{equation}
Let $F_{\ell} \lefttop{\vecK}\nbigw_{\vecp}
\lefttop{I}\psizero_{q_{I}(\vecu(0))}
\Vtildezero_{\veca}(\nbigm)$
denote the kernel of (\ref{eq;22.2.17.2}).
Namely,
\[
 F_{\ell} \lefttop{\vecK}\nbigw_{\vecp}
\lefttop{I}\psizero_{q_{I}(\vecu(0))}
\Vtildezero_{\veca}(\nbigm)
=
\lefttop{\vecK}\nbigw_{\vecp}
\lefttop{I}\psizero_{q_{I}(\vecu(0))}
\Vtildezero_{\veca}(\nbigm)
\cap
F_{\ell}
\lefttop{I}\psizero_{q_{I}(\vecu(0))}
\Vtildezero_{\veca}(\nbigm).
\]
Thus, we obtain an increasing filtration $F_{\bullet}$.
We also set
\[
\lefttop{\Lambda\setminus I}\Vzero_{<q_{\Lambda\setminus I}(\veca)}
\lefttop{\vecK}\nbigw_{\vecp}
\lefttop{I}\psizero_{q_{I}(\vecu(0))}
\Vtildezero_{\veca}(\nbigm)
:=
F_{0} \lefttop{\vecK}\nbigw_{\vecp}
\lefttop{I}\psizero_{q_{I}(\vecu(0))}
\Vtildezero_{\veca}(\nbigm).
\]
Because the kernel of
\[
  F_{\ell} \lefttop{\vecK}\nbigw_{\vecp}
\lefttop{I}\psizero_{q_{I}(\vecu(0))}
\Vtildezero_{\veca}(\nbigm)
\lrarr
\bigoplus_{\substack{I\subset L\\ |L\setminus I|=\ell}}
 \lefttop{\vecK}\nbigw_{\vecp}
\lefttop{L}\psizero_{q_{L}(\vecu(0))}
\Vtildezero_{\veca}(\nbigm)
\]
is $F_{\ell-1}$,
we obtain the following monomorphism:
\begin{equation}
 \label{eq;22.2.17.3}
\Gr^F_{\ell}\lefttop{\vecK}\nbigw_{\vecp}
\lefttop{I}\psizero_{q_{I}(\vecu(0))}
\Vtildezero_{\veca}(\nbigm)
\lrarr
\bigoplus_{\substack{I\subset L\\ |L\setminus I|=\ell}}
 \lefttop{\vecK}\nbigw_{\vecp}
\lefttop{L}\psizero_{q_{L}(\vecu(0))}
\Vtildezero_{\veca}(\nbigm).
\end{equation}

We shall use the following proposition later.
\begin{prop}
The morphism {\rm(\ref{eq;22.2.17.3})} induces
the following isomorphism
\begin{equation}
\label{eq;22.2.17.4}
 \Gr^F_{\ell}\lefttop{\vecK}\nbigw_{\vecp}
\lefttop{I}\psizero_{q_{I}(\vecu(0))}
\Vtildezero_{\veca}(\nbigm)
\simeq
\bigoplus_{\substack{I\subset L\\ |L\setminus I|=\ell}}
\lefttop{\Lambda\setminus L}\Vzero_{<q_{\Lambda\setminus L}(\veca)}
 \lefttop{\vecK}\nbigw_{\vecp}
\lefttop{L}\psizero_{q_{L}(\vecu(0))}
\Vtildezero_{\veca}(\nbigm).
\end{equation}
\end{prop}
\pf
The image of (\ref{eq;22.2.17.3}) is contained in
the right hand side of (\ref{eq;22.2.17.4}),
and hence we obtain the morphism (\ref{eq;22.2.17.4}).
It is enough to prove that
\begin{equation}
\label{eq;22.2.17.5}
\lefttop{\vecK}\nbigw_{\vecp}
\lefttop{I}\psizero_{q_{I}(\vecu(0))}
\Vtildezero_{\veca}(\nbigm)
\cap
\lefttop{I}\psizero_{q_{I}(\vecu(0))}
\Vtildezero_{\veca-\epsilon\veciti_{\Lambda\setminus L}}(\nbigm)
\lrarr
\lefttop{\Lambda\setminus L}\Vzero_{<q_{\Lambda\setminus L}(\veca)}
 \lefttop{\vecK}\nbigw_{\vecp}
\lefttop{L}\psizero_{q_{L}(\vecu(0))}
\Vtildezero_{\veca}(\nbigm)
\end{equation}
is an epimorphism
for any $I\subset L$.
We have only to consider the case where $m=|\Lambda|$.
It is enough to consider the case $|L\setminus I|=1$.
If $i(L)=i(I)$,
the morphism (\ref{eq;22.2.17.5})
is identified with
\begin{equation}
\label{eq;22.2.17.6}
\lefttop{\vecK_{\leq i(I)}}\nbigw_{\vecp}  
\lefttop{I}\psizero_{q_{I}(\vecu(0))}
\Vtildezero_{\veca-\epsilon\veciti_{\Lambda\setminus L}}(\nbigm)
\lrarr
\lefttop{\vecK_{\leq i(I)}}\nbigw_{\vecp}  
\lefttop{L}\psizero_{q_{L}(\vecu(0))}
\Vtildezero_{\veca-\epsilon\veciti_{\Lambda\setminus L}}(\nbigm).
\end{equation}
Hence, it is an epimorphism.
If $i(L)=i(I)+1$,
we can observe that
the left hand side of
(\ref{eq;22.2.17.5})
is identified with
the inverse image of
$\lefttop{\vecK_{\leq i(L)}}\nbigw_{\vecp}  
\lefttop{L}\psizero_{q_{L}(\vecu(0))}
\Vtildezero_{\veca-\epsilon\veciti_{\Lambda\setminus L}}(\nbigm)$
by the morphism (\ref{eq;22.2.17.6}).
Hence, (\ref{eq;22.2.17.5})
is an epimorphism.
\hfill\qed

\subsection{Sub-complexes of the de Rham complex}
\label{subsection;22.2.17.100}

Let $X$ and $H$ be as in \S\ref{subsection;22.2.14.1}.
On $\nbigx$,
we set
$\Omegatilde^{k}_{\nbigx/\cnum}:=
\lambda^{-k}p_{\lambda}^{\ast}\Omega^k_{X}$.
Let $(\nbigo_{\cnum_{\lambda}})_{\nbigx}$ denote
the pull back of $\nbigo_{\cnum_{\lambda}}$
by $\nbigx\lrarr X$.
For any open subset $U\subset\nbigx$,
let
$\Omegatilde^{k}_{U/\cnum}$
and
$(\nbigo_{\cnum_{\lambda}})_U$ denote the restrictions
$\Omegatilde^{k}_{\nbigx/\cnum}$
and $(\nbigo_{\cnum_{\lambda}})_{\nbigx}$ to $U$,
respectively.

Let $\nbigm$ be an $\nbigr_X$-module
underlying a good mixed twistor $\nbigd$-module on $(X,H)$.
The structure of $\nbigr_X$-module induces
a meromorphic family of the flat connections,
that is a morphism of sheaves
$\DD^f:\nbigm\lrarr \nbigm\otimes\Omegatilde^1_{\nbigx/\cnum}$.
Together with the exterior derivative on
$d:\Omega^k_X\lrarr \Omega^{k+1}_X$,
we obtain the de Rham complex
$\nbigm\otimes\Omegatilde^{\bullet}_{\nbigx/\cnum}$
of $\nbigm$.
It is an $(\nbigo_{\cnum_{\lambda}})_{\nbigx}$-complex.

We set
$\Omega^{\bullet}_{\nbigx/\cnum}:=
p_{\lambda}^{\ast}\Omega^{\bullet}_X$.
For any open subset $U\subset\nbigx$,
the restriction 
$\Omega^{\bullet}_{\nbigx/\cnum|U}$
is also denoted by
$\Omega^{\bullet}_{U/\cnum}$.
We obtain
$\DD:\nbigm\lrarr \nbigm\otimes \Omega^{1}_{\nbigx/\cnum}$
by setting $\DD=\lambda\DD^f$.
Together with
$\lambda d:\Omega^{k}_{\nbigx/\cnum}
\lrarr \Omega^{k+1}_{\nbigx/\cnum}$,
we obtain the $(\nbigo_{\cnum_{\lambda}})_{\nbigx}$-complex
$\nbigm\otimes \Omega^{\bullet}_{\nbigx/\cnum}$.

We have the isomorphisms
$\Omega^k_{\nbigx/\cnum}
\simeq
\Omegatilde^{k}_{\nbigx/\cnum}$
given by $\tau\longmapsto \lambda^{-k}\tau$.
They induce an isomorphism of
the $(\nbigo_{\cnum_{\lambda}})_{\nbigx}$-complexes
$\nbigm\otimes \Omega^{\bullet}_{\nbigx/\cnum}
\simeq
 \nbigm\otimes\Omegatilde^{\bullet}_{\nbigx/\cnum}$.

\subsubsection{Truncated de Rham complex}

For any $K\subset\Lambda$,
by applying the construction in \S\ref{subsection;22.2.17.20}
with $\veca=-\veciti_K\in\real^{\Lambda}$ and
$\vecu(0)=-\vecdelta_K\in(\real\times\cnum)^{\Lambda}$,
we obtain
$\lefttop{\Lambda}\Vtildezero_{-\veciti_K}(\nbigm)$.
For $J=(j_1,\ldots,j_k)\subset\Lambda$,
let $dz_J=dz_{j_1}\wedge\cdots\wedge dz_{j_k}$.
We set
\[
 \nbigctilde^k_{\tw}(\nbigm^{(\lambda_0)})
 =\bigoplus_{k_1+k_2=k}
 \bigoplus_{\substack{J\subset \Lambda \\ |J|=k_1}}
 \lefttop{\Lambda}\Vtildezero_{-\veciti_{\Lambda\setminus J}}(\nbigm)
 \cdot \lambda^{-k_1}dz_J\cdot
 \Omegatilde^{k_2}_{\nbigyzero/\cnum}
 \subset
 \bigl(
 \nbigm\otimes\Omegatilde^k_{\nbigx/\cnum}
 \bigr)_{|\nbigxzero}.
\]
Thus, we obtain an $(\nbigo_{\cnum_{\lambda}})_{\nbigxzero}$-subcomplex
$\nbigctilde^{\bullet}_{\tw}(\nbigm^{(\lambda_0)})
\subset
(\nbigm\otimes\Omegatilde^{\bullet}_{\nbigx/\cnum})
_{|\nbigxzero}$.
The following lemma is standard.
\begin{lem}
\label{lem;22.3.17.11}
The natural inclusion
$\nbigctilde^{\bullet}_{\tw}(\nbigm^{(\lambda_0)})
\lrarr
(\nbigm\otimes\Omegatilde^{\bullet}_{\nbigx/\cnum})
_{|\nbigxzero}$
is a quasi-isomorphism
of $(\nbigo_{\cnum_{\lambda}})_{\nbigxzero}$-complexes.
\hfill\qed
\end{lem}

\subsubsection{The induced complexes on the intersections}

On $\nbigxzero$, we set
\[
 \lefttop{I}\psizero
 \nbigctilde^k_{\tw}(\nbigm^{(\lambda_0)})
:=\bigoplus_{k_1+k_2=k}
 \bigoplus_{\substack{J\subset \Lambda \\ |J|=k_1}}
 \lefttop{I}\psizero_{-\vecdelta_{I\setminus J}}
 \lefttop{\Lambda}\Vtildezero_{-\veciti_{\Lambda\setminus J}}(\nbigm)
 \cdot \lambda^{-k_1}dz_J\cdot
 \Omegatilde^{k_2}_{\nbigyzero/\cnum}.
\]
Thus, we obtain an $(\nbigo_{\cnum_{\lambda}})_{\nbigxzero}$-complex 
$\lefttop{I}\psizero
\nbigctilde^{\bullet}_{\tw}(\nbigm^{(\lambda_0)})$.
For any $I\subset L$,
there exists a naturally defined morphism of 
$(\nbigo_{\cnum_{\lambda}})_{\nbigxzero}$-complexes:
\[
\lefttop{I}\psizero
\nbigctilde^k_{\tw}(\nbigm^{(\lambda_0)})
  \lrarr
\lefttop{L}\psizero
 \nbigctilde^k_{\tw}(\nbigm^{(\lambda_0)}).
\]

For any $\ell\in\seisuu_{\geq 0}$,
we set
\[
 F_{\ell}\lefttop{I}\psizero
 \nbigctilde^k_{\tw}(\nbigm^{(\lambda_0)})
:=\bigoplus_{k_1+k_2=k}
 \bigoplus_{\substack{J\subset \Lambda \\ |J|=k_1}}
 F_{\ell}\lefttop{I}\psizero_{-\vecdelta_{I\setminus J}}
 \lefttop{\Lambda}\Vtildezero_{-\veciti_{\Lambda\setminus J}}(\nbigm)
 \cdot \lambda^{-k_1}dz_J\cdot
 \Omegatilde^{k_2}_{\nbigyzero/\cnum}.
\]
Thus, we obtain
an $(\nbigo_{\cnum_{\lambda}})$-subcomplex
$F_{\ell}\lefttop{I}\psizero
\nbigctilde^{\bullet}_{\tw}(\nbigm^{(\lambda_0)})$
of $\lefttop{I}\psizero
\nbigctilde^{\bullet}_{\tw}(\nbigm^{(\lambda_0)})$.
We have
\[
 \Gr^{F}_{\ell}\lefttop{I}\psizero
 \nbigctilde^k_{\tw}(\nbigm^{(\lambda_0)})
 \simeq
 \bigoplus_{\substack{I\subset L\\ |L\setminus I|=\ell}}
 F_0\lefttop{L}\psizero
 \nbigctilde^k_{\tw}(\nbigm^{(\lambda_0)}).
\]

\subsection{Weight filtrations for some complexes associated with
tame harmonic bundles}
\label{subsection;22.2.17.101}

Let $(E,\delbar_E,\theta,h)$ be a tame harmonic bundle
on $(X,H)$.
Let $\gbige$ denote the $\nbigr_X$-module
underlying the associated pure twistor $\nbigd$-module $\gbigt(E)$.
(See \S\ref{subsection;22.4.26.10}.)
We also have the $\nbigr_X$-module $\gbige[\ast H]$
underlying $\gbigt(E)[\ast H]$.
(See \S\ref{subsection;22.4.26.11} and
\S\ref{subsection;22.4.25.41} for $\gbigt(E)[\ast H]$.)

\subsubsection{The case of $\gbige$}

For $K\subset\Lambda$,
we set
\[
 \lefttop{K}\nbigw_{|K|-1}
 \nbigctilde^{k}_{\tw}(\gbige^{(\lambda_0)})
=\bigoplus_{k_1+k_2=k}
 \bigoplus_{\substack{J\subset \Lambda\\ |J|=k_1
 }}
 \lefttop{K}\nbigw_{|K|-1-|K\cap J|}
 \lefttop{\Lambda}
 \Vtildezero_{-\veciti_{\Lambda\setminus J}}(\gbige)
 \cdot \lambda^{-k_1}dz_J\cdot
 \Omegatilde^{k_2}_{\nbigyzero/\cnum}.
\]
For $\vecK=(K_1,\ldots,K_m)\in\nbigs(\Lambda)$,
we set
\[
 \lefttop{\vecK}\nbigw_{\vecw(\vecK)}
 \nbigctilde^{k}_{\tw}(\gbige^{(\lambda_0)})
=\bigcap_{i=1}^m
 \lefttop{K_i}\nbigw_{|K_i|-1}
 \nbigctilde^{k}_{\tw}(\gbige^{(\lambda_0)}).
\]
Thus, we obtain an $(\nbigo_{\cnum_{\lambda}})_{\nbigxzero}$-subcomplex
$\lefttop{\vecK}\nbigw_{\vecw(\vecK)}
\nbigctilde^{\bullet}_{\tw}(\gbige^{(\lambda_0)})$
of
$\nbigctilde^{\bullet}_{\tw}(\gbige^{(\lambda_0)})$.

More generally,
for any $I\subset \Lambda$,
and for any $\ell\in\seisuu_{\geq 0}$,
we set
\[
 F_{\ell}\lefttop{K}\nbigw_{|K|-1}
 \lefttop{I}\psizero
 \nbigctilde^{k}_{\tw}(\gbige^{(\lambda_0)})
=\bigoplus_{k_1+k_2=k}
 \bigoplus_{\substack{J\subset \Lambda\\ |J|=k_1
 }}
 F_{\ell}
 \lefttop{K}\nbigw_{|K|-1-|K\cap J|}
 \lefttop{I}\psizero_{-\vecdelta_{I\setminus J}}
 \lefttop{\Lambda}
 \Vtildezero_{-\veciti_{\Lambda\setminus J}}(\gbige)
 \cdot \lambda^{-k_1}dz_J\cdot
 \Omegatilde^{k_2}_{\nbigyzero/\cnum}.
\]
For $\vecK=(K_1,\ldots,K_m)\in\nbigs(\Lambda)$,
we set
\[
 F_{\ell}
 \lefttop{\vecK}\nbigw_{\vecw(\vecK)}
  \lefttop{I}\psizero
 \nbigctilde^{k}_{\tw}(\gbige^{(\lambda_0)})
 =\bigcap_{i=1}^m
 F_{\ell}
 \lefttop{K_i}\nbigw_{|K_i|-1}
   \lefttop{I}\psizero
 \nbigctilde^{k}_{\tw}(\gbige^{(\lambda_0)}).
\]
We set
$\lefttop{\vecK}\nbigw_{\vecw(\vecK)}
  \lefttop{I}\psizero
 \nbigctilde^{k}_{\tw}(\gbige^{(\lambda_0)})
=F_{|\Lambda|+1}
\lefttop{\vecK}\nbigw_{\vecw(\vecK)}
\lefttop{I}\psizero
\nbigctilde^{k}_{\tw}(\gbige^{(\lambda_0)})$.
Thus, we obtain
$(\nbigo_{\cnum_{\lambda}})_{\nbigxzero}$-subcomplexes
\[
F_{\ell}\lefttop{\vecK}\nbigw_{\vecw(\vecK)}
\lefttop{I}\psizero
\nbigctilde^{\bullet}_{\tw}(\gbige^{(\lambda_0)})
\subset
\lefttop{\vecK}\nbigw_{\vecw(\vecK)}
\lefttop{I}\psizero
\nbigctilde^{\bullet}_{\tw}(\gbige^{(\lambda_0)})
\]
of
$\lefttop{I}\psizero
\nbigctilde^{\bullet}_{\tw}(\gbige^{(\lambda_0)})$.
We have the following natural isomorphism of
$(\nbigo_{\cnum_{\lambda}})_{\nbigxzero}$-complexes:
\begin{equation}
\label{eq;22.4.26.20}
 \Gr^F_{\ell}
 \lefttop{\vecK}\nbigw_{\vecw(\vecK)}
  \lefttop{I}\psizero
\nbigctilde^{\bullet}_{\tw}(\gbige^{(\lambda_0)})
\simeq
\bigoplus_{\substack{L\subset\Lambda \\ |L\setminus I|=\ell}}
 F_0
 \lefttop{\vecK}\nbigw_{\vecw(\vecK)}
 \lefttop{L}\psizero
\nbigctilde^{\bullet}_{\tw}(\gbige^{(\lambda_0)}).
\end{equation}

\begin{prop}
\label{prop;22.2.17.40}
The natural morphism
 $\lefttop{\vecK}\nbigw_{\vecw(\vecK)}
   \lefttop{I}\psizero
\nbigctilde^{\bullet}_{\tw}(\gbige^{(\lambda_0)})
 \lrarr
   \lefttop{I}\psizero
 \nbigctilde^{\bullet}_{\tw}(\gbige^{(\lambda_0)})$
 is a quasi-isomorphism of
$(\nbigo_{\cnum_{\lambda}})_{\nbigxzero}$-complexes.
\end{prop}
\pf
By Proposition \ref{prop;22.3.14.12},
the induced morphisms
\[
 F_0
 \lefttop{\vecK}\nbigw_{\vecw(\vecK)}
 \lefttop{L}\psizero
\nbigctilde^{\bullet}_{\tw}(\gbige^{(\lambda_0)})
\lrarr
 F_0
 \lefttop{L}\psizero
\nbigctilde^{\bullet}_{\tw}(\gbige^{(\lambda_0)})
\]
are isomorphisms
for any $L\supset I$.
Hence, we obtain the claim of Proposition \ref{prop;22.2.17.40}
from (\ref{eq;22.4.26.20}).
\hfill\qed

\subsubsection{The case of $\gbige[\ast H]$}

For $K\subset\Lambda$,
we set
\[
 \lefttop{K}\nbigw_{|K|-1}
 \nbigctilde^{k}_{\tw}(\gbige[\ast H]^{(\lambda_0)})
=\bigoplus_{k_1+k_2=k}
 \bigoplus_{\substack{J\subset \Lambda\\ |J|=k_1
 }}
 \lefttop{K}\nbigw_{|K|-1-2|K\cap J|}
 \lefttop{\Lambda}
 \Vtildezero_{-\veciti_{\Lambda\setminus J}}(\gbige[\ast H])
 \cdot \lambda^{-k_1}dz_J\cdot
 \Omegatilde^{k_2}_{\nbigyzero/\cnum}.
\]
For $\vecK=(K_1,\ldots,K_m)\in\nbigs(\Lambda)$,
we set
\[
 \lefttop{\vecK}\nbigw_{\vecw(\vecK)}
 \nbigctilde^{k}_{\tw}(\gbige[\ast H]^{(\lambda_0)})
=\bigcap_{i=1}^m
 \lefttop{K_i}\nbigw_{|K_i|-1}
 \nbigctilde^{k}_{\tw}(\gbige[\ast H]^{(\lambda_0)}).
\]
Thus, we obtain an $(\nbigo_{\cnum_{\lambda}})_{\nbigxzero}$-subcomplex
$\lefttop{\vecK}\nbigw_{\vecw(\vecK)}
\nbigctilde^{\bullet}_{\tw}(\gbige[\ast H]^{(\lambda_0)})$
of
$\nbigctilde^{\bullet}_{\tw}(\gbige[\ast H]^{(\lambda_0)})$.

More generally,
for any $I\subset \Lambda$,
and for any $\ell\in\seisuu_{\geq 0}$,
we set
\[
 F_{\ell}\lefttop{K}\nbigw_{|K|-1}
 \lefttop{I}\psizero
 \nbigctilde^{k}_{\tw}(\gbige[\ast H]^{(\lambda_0)})
=\bigoplus_{k_1+k_2=k}
 \bigoplus_{\substack{J\subset \Lambda\\ |J|=k_1
 }}
 F_{\ell}
 \lefttop{K}\nbigw_{|K|-1-2|K\cap J|}
 \lefttop{I}\psizero_{-\vecdelta_{I\setminus J}}
 \lefttop{\Lambda}
 \Vtildezero_{-\veciti_{\Lambda\setminus J}}(\gbige[\ast H])
 \cdot \lambda^{-k_1}dz_J\cdot
 \Omegatilde^{k_2}_{\nbigyzero/\cnum}.
\]
For $\vecK=(K_1,\ldots,K_m)\in\nbigs(\Lambda)$,
we set
\[
 F_{\ell}
 \lefttop{\vecK}\nbigw_{\vecw(\vecK)}
  \lefttop{I}\psizero
 \nbigctilde^{k}_{\tw}(\gbige[\ast H]^{(\lambda_0)})
 =\bigcap_{i=1}^m
 F_{\ell}
 \lefttop{K_i}\nbigw_{|K_i|-1}
   \lefttop{I}\psizero
 \nbigctilde^{k}_{\tw}(\gbige[\ast H]^{(\lambda_0)}).
\]
We also set
$\lefttop{\vecK}\nbigw_{\vecw(\vecK)}
  \lefttop{I}\psizero
 \nbigctilde^{k}_{\tw}(\gbige[\ast H]^{(\lambda_0)})
=
F_{|\Lambda|+1}
 \lefttop{\vecK}\nbigw_{\vecw(\vecK)}
  \lefttop{I}\psizero
 \nbigctilde^{k}_{\tw}(\gbige[\ast H]^{(\lambda_0)})$.
Thus, we obtain $(\nbigo_{\cnum_{\lambda}})_{\nbigxzero}$-subcomplexes
\[
F_{\ell}\lefttop{\vecK}\nbigw_{\vecw(\vecK)}
\lefttop{I}\psizero
\nbigctilde^{\bullet}_{\tw}(\gbige[\ast H]^{(\lambda_0)})
\subset
\lefttop{\vecK}\nbigw_{\vecw(\vecK)}
\lefttop{I}\psizero
\nbigctilde^{\bullet}_{\tw}(\gbige[\ast H]^{(\lambda_0)})
\] 
of
$\lefttop{\vecK}\nbigw_{\vecw(\vecK)}
\lefttop{I}\psizero
\nbigctilde^{\bullet}_{\tw}(\gbige[\ast H]^{(\lambda_0)})$.
We have the following natural isomorphism:
\[
 \Gr^F_{\ell}
 \lefttop{\vecK}\nbigw_{\vecw(\vecK)}
  \lefttop{I}\psizero
\nbigctilde^{\bullet}_{\tw}(\gbige[\ast H]^{(\lambda_0)})
\simeq
\bigoplus_{\substack{L\subset\Lambda \\ |L\setminus I|=\ell}}
 F_0
 \lefttop{\vecK}\nbigw_{\vecw(\vecK)}
 \lefttop{L}\psizero
\nbigctilde^{\bullet}_{\tw}(\gbige[\ast H]^{(\lambda_0)}).
\]

\begin{lem}
The morphism
$\lefttop{I}\psizero
 \nbigctilde^{\bullet}_{\tw}(\gbige^{(\lambda_0)})
 \lrarr
 \lefttop{I}\psizero
 \nbigctilde^{\bullet}_{\tw}(\gbige[\ast H]^{(\lambda_0)})$
 induces
 the morphisms of
$(\nbigo_{\cnum_{\lambda}})_{\nbigxzero}$-subcomplexes: 
 \[
 F_{\ell}\lefttop{\vecK}\nbigw_{\vecw(\vecK)}
 \lefttop{I}\psizero
 \nbigctilde^{\bullet}_{\tw}(\gbige^{(\lambda_0)})
 \lrarr
  F_{\ell}\lefttop{\vecK}\nbigw_{\vecw(\vecK)}
 \lefttop{I}\psizero
 \nbigctilde^{\bullet}_{\tw}(\gbige[\ast H]^{(\lambda_0)}).
\]
\end{lem}
\pf
We consider the following morphism:
\[
\lefttop{K}\psizero_{-\vecdelta_{K\setminus J}}
\lefttop{\Lambda}\Vtildezero_{-\veciti_{\Lambda\setminus J}}
(\gbige)
\lrarr
\lefttop{K}\psizero_{-\vecdelta_{K\setminus J}}
\lefttop{\Lambda}\Vtildezero_{-\veciti_{\Lambda\setminus J}}
(\gbige[\ast H]).
\]
By Proposition \ref{prop;22.2.17.30} below,
it induces
\begin{equation}
\label{eq;22.2.17.31}
\lefttop{K}\nbigw_{p}
\lefttop{K}\psizero_{-\vecdelta_{K\setminus J}}
\lefttop{\Lambda}\Vtildezero_{-\veciti_{\Lambda\setminus J}}
(\gbige)(\ast \del H_K)
\lrarr
\lefttop{K}\nbigw_{p-|J\cap K|}
\lefttop{K}\psizero_{-\vecdelta_{K\setminus J}}
\lefttop{\Lambda}\Vtildezero_{-\veciti_{\Lambda\setminus J}}
(\gbige[\ast H])(\ast \del H_K).
\end{equation}
Note that
\[
 \Gr^{\lefttop{K}\nbigw}_{p'}\Bigl(
\lefttop{K}\psizero_{-\vecdelta_{K\setminus J}}
\lefttop{\Lambda}\Vtildezero_{-\veciti_{\Lambda\setminus J}}
(\gbige[\ast H])
\Bigr)
=\lefttop{\Lambda\setminus K}
\Vtilde_{-\veciti_{\Lambda\setminus(K\cup J)}}
\Gr^{\lefttop{K}\nbigw}_{p'}
\Bigl(
\lefttop{K}\psizero_{-\vecdelta_{K\setminus J}}
(\gbige[\ast H])
\Bigr)
\]
is torsion-free sheaf on
$\nbigh_K^{(\lambda_0)}$ for any $p'\in\seisuu$.
Hence, (\ref{eq;22.2.17.31})
induces
\[
 \lefttop{K}\nbigw_{p}
\lefttop{K}\psizero_{-\vecdelta_{K\setminus J}}
\lefttop{\Lambda}\Vtildezero_{-\veciti_{\Lambda\setminus J}}
(\gbige)
\lrarr
\lefttop{K}\nbigw_{p-|J\cap K|}
\lefttop{K}\psizero_{-\vecdelta_{K\setminus J}}
\lefttop{\Lambda}\Vtildezero_{-\veciti_{\Lambda\setminus J}}
(\gbige[\ast H]).
\]
Then, we obtain the claim of the lemma
by the construction.
\hfill\qed

\subsubsection{Some quasi-isomorphisms}

Let $P$ be any point of $\nbigh_{\Lambda}^{(\lambda_0)}$.
Let $(\nbigo_{\cnum_{\lambda}})_P$ denote the stalk of
$(\nbigo_{\cnum_{\lambda}})_{\nbigxzero}$ at $P$.
Let $\ttC((\nbigo_{\cnum_{\lambda}})_P)$
denote the category of bounded
$(\nbigo_{\cnum_{\lambda}})_P$-complexes.

Let $\nbigf_{10}^{\bullet}(\gbige_P)$,
$\nbigf_{12}^{\bullet}(\gbige_P)$ and
$\nbigf_{13}^{\bullet}(\gbige_P)$
denote the functors from
$\nbigsbar(\Lambda)$
to $\ttC((\nbigo_{\cnum_{\lambda}})_P)$
defined as follows:
\[
 \nbigf^{\bullet}_{10}(\gbige_P)(\vecK):=
 \lefttop{\vecK}\nbigw_{\vecw(\vecK)}
 \nbigctilde^{\bullet}_{\tw}(\gbige^{(\lambda_0)}[\ast H])_P,
\]
\[
 \nbigf^{\bullet}_{12}(\gbige_P)(\vecK):=
 \lefttop{\vecK}\nbigw_{\vecw(\vecK)}
 \nbigctilde^{\bullet}_{\tw}(\gbige^{(\lambda_0)})_P,
\]
\[
 \nbigf^{\bullet}_{13}(\gbige_P)(\vecK):=
 \nbigctilde^{\bullet}_{\tw}(\gbige^{(\lambda_0)})_P.
\]
Let $\ttF^{\bullet}_i(\gbige_P)\in
\ttC^{\wc}(\ttX(\Lambda)_{\geq 0};
(\nbigo_{\cnum_{\lambda}})_P)$
denote the corresponding objects.

\begin{thm}
\label{thm;22.2.17.103}
The induced morphisms
\[
 \begin{CD}
  R\pi_{\ttXbar(\Lambda)\ast}
  \ttF^{\bullet}_{10}(\gbige_P)
  @<{a_1}<<
  R\pi_{\ttXbar(\Lambda)\ast}
  \ttF^{\bullet}_{12}(\gbige_P)
  @>{a_2}>>
  R\pi_{\ttXbar(\Lambda)\ast}
  \ttF^{\bullet}_{13}(\gbige_P)
 \end{CD}
\]
are quasi-isomorphisms.
\end{thm}
\pf
By Proposition \ref{prop;22.2.17.40},
$a_2$ is a quasi-isomorphism.
For any $\ell\in\seisuu_{\geq 0}$,
we define the functors
$F_{\ell} \nbigf^{\bullet}_{i}(\gbige_P)$ $(i=10,12)$
from $\nbigsbar(\Lambda)$
to $\ttC((\nbigo_{\cnum_{\lambda}})_P)$
as follows:
\[
F_{\ell} \nbigf^{\bullet}_{10}(\gbige_P)(\vecK):=
 F_{\ell}\lefttop{\vecK}\nbigw_{\vecw(\veck)}
 \nbigctilde^{\bullet}_{\tw}(\gbige^{(\lambda_0)}[\ast H])_P,
\]
\[
 F_{\ell}\nbigf^{\bullet}_{12}(\gbige_P)(\vecK):=
 F_{\ell}\lefttop{\vecK}\nbigw_{\vecw(\vecK)}
 \nbigctilde^{\bullet}_{\tw}(\gbige^{(\lambda_0)})_P.
\]
They induce the subcomplexes of sheaves
$F_{\ell}\ttF^{\bullet}_i(\gbige_P)
\subset\ttF^{\bullet}_i(\gbige_P)$.

For any $\ell\in\seisuu_{\geq 0}$,
we define the functors
$\Gr^F_{\ell} \nbigf^{\bullet}_{i}(\gbige_P)$ $(i=10,12)$
from $\nbigsbar(\Lambda)$
to $\ttC((\nbigo_{\cnum_{\lambda}})_P)$
as follows:
\[
\Gr^F_{\ell} \nbigf^{\bullet}_{10}(\gbige_P)(\vecK):=
\Gr^F_{\ell}\lefttop{\vecK}\nbigw_{\vecw(\vecK)}
 \nbigctilde^{\bullet}_{\tw}(\gbige^{(\lambda_0)}[\ast H])_P
 =\bigoplus_{|I|=\ell}
 F_0\lefttop{\vecK}\nbigw_{\vecw(\vecK)}
 \lefttop{I}\psizero
  \nbigctilde^{\bullet}_{\tw}(\gbige^{(\lambda_0)}[\ast H])_P,
\]
\[
 \Gr^F_{\ell}\nbigf^{\bullet}_{12}(\gbige_P)(\vecK):=
 \Gr^F_{\ell}\lefttop{\vecK}\nbigw_{\vecw(\vecK)}
 \nbigctilde^{\bullet}_{\tw}(\gbige^{(\lambda_0)})_P
 =\bigoplus_{|I|=\ell}
 F_0\lefttop{\vecK}\nbigw_{\vecw(\vecK)}
 \lefttop{I}\psizero
  \nbigctilde^{\bullet}_{\tw}(\gbige^{(\lambda_0)})_P.
\]
The complex of sheaves $\Gr^F_{\ell}\ttF^{\bullet}_i(\gbige_P)$
correspond to the functors
$\Gr^F_{\ell}\nbigf^{\bullet}_{i}(\gbige_P)$.

For $I\subset \Lambda$ with $|I|=\ell$
and $\veca\in\real^{\Lambda\setminus I}_{\leq 0}$,
we have the $(V_H\nbigr_X)_P$-module
$\lefttop{\Lambda\setminus I}\Vzero_{<\veca}
\lefttop{I}\psi_{-\vecdelta_I}(\gbige)_P$
with $\vecf_I=(f_i\,|\,i\in I)$.
For $\vecK\in\nbigsbar(\Lambda)$,
we have the following $(\nbigo_{\cnum_{\lambda}})_P$-complexes:
\[
\nbigf^{\bullet}_{10,\Lambda}
\bigl(
\lefttop{\Lambda\setminus I}\Vzero_{<\veca}
\lefttop{I}\psi_{-\vecdelta_I}(\gbige)_P,
\vecf_I
\bigr)(\vecK)
=
\lefttop{\vecK_{\leq i(I)}}
 W_{\vecw(\vecK_{\leq i(I)})}
 \IntC^{\bullet}\bigl(
\lefttop{\Lambda\setminus I}\Vzero_{<\veca}
\lefttop{I}\psi_{-\vecdelta_I}(\gbige)_P,
\vecf_I,\ast I
 \bigr),
\]
\[
\nbigf^{\bullet}_{12,\Lambda}
\bigl(
\lefttop{\Lambda\setminus I}\Vzero_{<\veca}
\lefttop{I}\psi_{-\vecdelta_I}(\gbige)_P,
\vecf_I
\bigr)(\vecK)
=
  \lefttop{\vecK_{\leq i(I)}}
 W_{\vecw(\vecK_{\leq i(I)})}
 \IntC^{\bullet}\bigl(
\lefttop{\Lambda\setminus I}\Vzero_{<\veca}
\lefttop{I}\psi_{-\vecdelta_I}(\gbige)_P,
\vecf_I
 \bigr).
\]
We obtain $(\nbigo_{\cnum_{\lambda}})_P$-complexes
$\nbigb^{\bullet}_{i}
 \bigl(
\lefttop{I}\psi_{-\vecdelta_I}(\gbige)_P,
\vecf_I
\bigr)(\vecK)$
$(i=10,12)$
by setting
\[
 \nbigb^k_{10}
 \bigl(
\lefttop{I}\psi_{-\vecdelta_I}(\gbige)_P,
\vecf_I
\bigr)(\vecK)
=\bigoplus_{k_1+k_2+k_3=k}
\bigoplus_{\substack{
J\subset\Lambda\setminus I\\
|J|=k_2
}}
\nbigf^{k_1}_{10,\Lambda}
\bigl(
\lefttop{\Lambda\setminus I}\Vzero_{<-\vecdelta_{\Lambda\setminus J}}
\lefttop{I}\psi_{-\vecdelta_I}(\gbige)_P,
\vecf_I
\bigr)(\vecK)\,\lambda^{-k_2}dz_J\,
\Omegatilde^{k_3}_{\nbigyzero/\cnum},
\]
\[
 \nbigb^k_{12}
 \bigl(
\lefttop{I}\psi_{-\vecdelta_I}(\gbige)_P,
\vecf_I
\bigr)(\vecK)
=\bigoplus_{k_1+k_2+k_3=k}
\bigoplus_{\substack{
J\subset\Lambda\setminus I\\
|J|=k_2
}}
\nbigf^{k_1}_{12,\Lambda}
\bigl(
\lefttop{\Lambda\setminus I}\Vzero_{<-\vecdelta_{\Lambda\setminus J}}
\lefttop{I}\psi_{-\vecdelta_I}(\gbige)_P,
\vecf_I
\bigr)(\vecK)\,\lambda^{-k_2}dz_J\,
\Omegatilde^{k_3}_{\nbigyzero/\cnum}.
\]
Thus, we obtain functors
$\nbigf^{\bullet}_{i,\Lambda}
 \bigl(
 \lefttop{\Lambda\setminus I}\Vzero_{<\veca}
\lefttop{I}\psi_{-\vecdelta_I}(\gbige)_P,
\vecf_I
\bigr)$
and 
$\nbigb^{\bullet}_i
 \bigl(
\lefttop{I}\psi_{-\vecdelta_I}(\gbige)_P,
\vecf_I
\bigr)$
from $\nbigsbar(\Lambda)$
to $\ttC((\nbigo_{\cnum_{\lambda}})_P)$.

There exists the following commutative diagram
of functors:
\[
 \begin{CD}
  \Gr^F_{\ell}\nbigf^{\bullet}_{12}(\gbige_P)
  @>{\simeq}>>
  \bigoplus_{|I|=\ell}
  \nbigb^{\bullet}_{12}(\lefttop{I}\psi_{-\vecdelta_I}(\gbige_P),\vecf_I)
  \\
  @VVV @VVV \\
  \Gr^F_{\ell}\nbigf^{\bullet}_{10}(\gbige_P)
  @>{\simeq}>>
  \bigoplus_{|I|=\ell}
  \nbigb^{\bullet}_{10}(\lefttop{I}\psi_{-\vecdelta_I}(\gbige_P),\vecf_I).
 \end{CD}
\]

We obtain the complex of sheaves
$\ttF^{\bullet}_{i,\Lambda}
 \bigl(
 \lefttop{\Lambda\setminus I}\Vzero_{<-\vecdelta_{\Lambda\setminus J}}
\lefttop{I}\psi_{-\vecdelta_I}(\gbige)_P,
\vecf_I
\bigr)\in
\ttC^{\wc}(\ttX(\Lambda)_{\geq 0};(\nbigo_{\cnum_{\lambda}})_P)$
corresponding to the functors
$\nbigf^{\bullet}_{i,\Lambda}
 \bigl(
 \lefttop{\Lambda\setminus I}\Vzero_{<-\vecdelta_{\Lambda\setminus J}}
 \lefttop{I}\psi_{-\vecdelta_I}(\gbige)_P,
\vecf_I
\bigr)$.
By Proposition \ref{prop;22.2.17.50}
and Proposition \ref{prop;22.2.17.51},
the induced morphisms
\[
 R\pi_{\ttXbar(\Lambda)\ast}
 \ttF^{\bullet}_{12,\Lambda}
 \bigl(
  \lefttop{\Lambda\setminus I}\Vzero_{<-\vecdelta_{\Lambda\setminus J}}
\lefttop{I}\psi_{-\vecdelta_I}(\gbige)_P,
\vecf_I
\bigr)
\lrarr
 R\pi_{\ttXbar(\Lambda)\ast}
 \ttF^{\bullet}_{10,\Lambda}
 \bigl(
  \lefttop{\Lambda\setminus I}\Vzero_{<-\vecdelta_{\Lambda\setminus J}}
\lefttop{I}\psi_{-\vecdelta_I}(\gbige)_P,
\vecf_I
\bigr)
\]
are quasi-isomorphisms for any $I$.
Hence, we obtain that
\[
R\pi_{\ttXbar(\Lambda)\ast}
\Gr^{F}_{\ell}\ttF_{12}^{\bullet}(\gbige)
\lrarr
R\pi_{\ttXbar(\Lambda)\ast}
\Gr^{F}_{\ell}\ttF_{10}^{\bullet}(\gbige)
\]
are quasi-isomorphisms for any $\ell$.
Thus, we obtain that $a_1$ is a quasi-isomorphism.
\hfill\qed

\subsection{Globalization
in the $\cnum_{\lambda}$-direction}
\label{subsection;22.3.17.30}

\subsubsection{Truncated de Rham complex}

Let $X$ and $H$ be as in \S\ref{subsection;22.2.14.1}.
Let $\nbigm\in\nbigc(X)$ be the $\nbigr_X$-module
underlying a good mixed twistor $\nbigd$-module on $(X,H)$.

\begin{lem}
\label{lem;22.3.17.21}
If $\nbigx^{(\lambda_1)}\subset\nbigx^{(\lambda_0)}$,
there exists the following natural inclusion:
\begin{equation}
\label{eq;22.3.17.10}
 \nbigctilde^{\bullet}_{\tw}(\nbigm^{(\lambda_0)})_{|\nbigx^{(\lambda_1)}}
 \lrarr
 \nbigctilde^{\bullet}_{\tw}(\nbigm^{(\lambda_1)}).
 \end{equation}
It is a quasi-isomorphism.
\end{lem}
\pf
There exists the inclusion (\ref{eq;22.3.17.10})
by the construction of the complexes.
By Lemma \ref{lem;22.3.17.11},
it is a quasi-isomorphism.
\hfill\qed

\vspace{.1in}
For any $(\lambda_0,P_1)\in\nbigx$,
we set
\begin{equation}
\label{eq;22.3.17.20}
 \nbigctilde^{\bullet}_{\tw}(\nbigm)_{(\lambda_0,P_1)}
 :=\nbigctilde^{\bullet}_{\tw}(\nbigm^{(\lambda_0)})_{(\lambda_0,P_1)}
 \subset
 (\nbigm\otimes\Omegatilde^{\bullet}_{\nbigx/\cnum})_{(\lambda_0,P_1)}.
\end{equation}
By Lemma \ref{lem;22.3.17.21},
there exists the unique subcomplex
$\nbigctilde^{\bullet}_{\tw}(\nbigm)
\subset
\nbigm\otimes\Omegatilde^{\bullet}_{\nbigx/\cnum}$
such that
the stalk of $\nbigctilde^{\bullet}_{\tw}(\nbigm)$
at $(\lambda_0,P_1)\in\nbigx$
is equal to (\ref{eq;22.3.17.20}).
The natural inclusion
$\nbigctilde^{\bullet}_{\tw}(\nbigm)
\lrarr
\nbigm\otimes\Omegatilde^{\bullet}_{\nbigx/\cnum}$
is a quasi-isomorphism.

\subsubsection{Weight filtrations}

Let $\gbige$ be as in \S\ref{subsection;22.2.17.101}.
Let $\vecK\in\nbigs(\Lambda)$.
For any $(\lambda_0,P_1)\in\nbigx$,
we set
\begin{equation}
\label{eq;22.3.17.22}
 \lefttop{\vecK}\nbigw_{\vecw(\vecK)}
 \nbigctilde^{\bullet}_{\tw}(\gbige)_{(\lambda_0,P_1)}
 :=
  \lefttop{\vecK}\nbigw_{\vecw(\vecK)}
 \nbigctilde^{\bullet}_{\tw}(\gbige^{(\lambda_0)})_{(\lambda_0,P_1)}
 \subset
 \nbigctilde^{\bullet}_{\tw}(\gbige)_{(\lambda_0,P_1)}.
\end{equation}
There exists the unique subcomplex
$\lefttop{\vecK}\nbigw_{\vecw(\vecK)}
\nbigctilde^{\bullet}_{\tw}(\gbige)
\subset
\nbigctilde^{\bullet}_{\tw}(\gbige)$
whose talk at each $(\lambda_0,P_1)$ is 
equal to (\ref{eq;22.3.17.22}).
By Proposition \ref{prop;22.2.17.40},
the natural inclusion
$\lefttop{\vecK}\nbigw_{\vecw(\vecK)}
\nbigctilde^{\bullet}_{\tw}(\gbige)
\lrarr
\nbigctilde^{\bullet}_{\tw}(\gbige)$
is a quasi-isomorphism.

For any $(\lambda_0,P_1)\in\nbigx$,
we set
\begin{equation}
\label{eq;22.3.17.23}
 \lefttop{\vecK}\nbigw_{\vecw(\vecK)}
 \nbigctilde^{\bullet}_{\tw}(\gbige[\ast H])_{(\lambda_0,P_1)}
 :=
  \lefttop{\vecK}\nbigw_{\vecw(\vecK)}
 \nbigctilde^{\bullet}_{\tw}(\gbige[\ast H]^{(\lambda_0)})_{(\lambda_0,P_1)}
 \subset
 \nbigctilde^{\bullet}_{\tw}(\gbige[\ast H])_{(\lambda_0,P_1)}.
\end{equation}
There exists the unique subcomplex
$\lefttop{\vecK}\nbigw_{\vecw(\vecK)}
\nbigctilde^{\bullet}_{\tw}(\gbige[\ast H])
\subset
\nbigctilde^{\bullet}_{\tw}(\gbige[\ast H])$
whose talk at each $(\lambda_0,P_1)$ is 
equal to (\ref{eq;22.3.17.23}).

\subsubsection{Some special cases}

Let $\nbigm\in\nbigc(X)$ be the $\nbigr_X$-module
underlying a good mixed twistor $\nbigd$-module.
Suppose that
$\lefttop{i}\psitilde_u(\nbigm)=0$
unless $u\in \real\times\{0\}$,
i.e., 
$\KMS(\nbigm,i)\subset\real\times\{0\}$ for any $i$.
Then, the $V$-filtrations $\lefttop{i}V(\nbigm)$ along $z_i$
are defined on $\nbigx$, not only on $\nbigxzero$.
The multi-$V$-filtrations $\lefttop{I}V(\nbigm)$
are also defined on $\nbigx$.
We obtain
$\lefttop{I}\psi_{q_I(\veca)}
\lefttop{\Lambda}V_{\veca}(\nbigm)$.

In this case,
if $\nbigx^{(\lambda_1)}\subset\nbigx^{(\lambda_0)}$,
we have
$\nbigctilde^{\bullet}_{\tw}(\nbigm^{(\lambda_0)})
 _{|\nbigx^{(\lambda_1)}}
=\nbigctilde^{\bullet}_{\tw}(\nbigm^{(\lambda_1)})$.
Hence, we have the following for any $\lambda_0$:
\[
 \nbigctilde^{\bullet}_{\tw}(\nbigm)_{|\nbigxzero}
=\nbigctilde^{\bullet}_{\tw}(\nbigm^{(\lambda_0)})
\]
If moreover $\nbigm=\gbige$
induced by a tame harmonic bundle on $(X,H)$,
we have
\[
 \lefttop{\vecK}
 \nbigw_{\vecw(\vecK)}
 \nbigctilde^{\bullet}_{\tw}(\gbige)_{|\nbigxzero}
=
 \lefttop{\vecK}
 \nbigw_{\vecw(\vecK)}
 \nbigctilde^{\bullet}_{\tw}(\gbige^{(\lambda_0)}),
\]
\[
 \lefttop{\vecK}
 \nbigw_{\vecw(\vecK)}
 \nbigctilde^{\bullet}_{\tw}(\gbige[\ast H])_{|\nbigxzero}
=
 \lefttop{\vecK}
 \nbigw_{\vecw(\vecK)}
 \nbigctilde^{\bullet}_{\tw}(\gbige[\ast H]^{(\lambda_0)}).
\]

Suppose moreover that $\nbigm$ is the $\nbigr_X$-module
associated with a Hodge module.
Then,
$\nbigm$ is $\cnum^{\ast}$-homogeneous
with respect to the $\cnum^{\ast}$-action on $\nbigx$
given by $a(\lambda,x)=(a\lambda,x)$ $(a\in\cnum^{\ast})$.
(See \S\ref{subsection;22.3.16.30} below.)
It induces $\cnum^{\ast}$-actions on
$\lefttop{I}V\nbigm$
and $\nbigctilde^{\bullet}_{\tw}(\nbigm)$.
For $\gbige$ induced by a polarized variation of
Hodge structure on $X\setminus H$,
we have the induced $\cnum^{\ast}$-actions
on 
$\lefttop{\vecK}
 \nbigw_{\vecw(\vecK)}
 \nbigctilde^{\bullet}_{\tw}(\gbige)$
 and 
$\lefttop{\vecK}
 \nbigw_{\vecw(\vecK)}
 \nbigctilde^{\bullet}_{\tw}(\gbige[\ast H])$.

\section{Twistor complexes and $L^2$-complexes}
\label{section;22.4.10.3}

\subsection{Some spaces}

Let $X$ be a complex manifold.
Let $H$ be a simple normal crossing hypersurface
with the irreducible decomposition $H=\bigcup_{i\in\Lambda}H_i$.
For simplicity,
we assume that $\Lambda$ is finite.
For $I\subset \Lambda$,
we set $H_I=\bigcap_{i\in I}H_i$
and $H(I):=\bigcup_{i\in I}H_i$.
For the empty set,
we put $H_{\emptyset}=X$.
We set $H_I^{\circ}:=H_I\setminus\bigcup_{j\not\in I}H_j$.
For each $p\geq 1$, we set
\[
 H_{[p]}:=\bigcup_{\substack{I\subset \Lambda\\ |I|=p}}H_I.
\]
If $p>\dim X$,
then $H_{[p]}=\emptyset$.

For $H_i$,
let $\nbigo_X(H_i)$ be the line bundle associated with $H_i$.
We take a Hermitian metric $h_i$ of $\nbigo_X(H_i)$.
Let $\sigma_i:\nbigo_X\lrarr\nbigo_X(H_i)$ denote the canonical section.
We obtain the $C^{\infty}$-function
$|\sigma_i|_{h_i}$ on $X$.
We assume that $|\sigma_i|_{h_i}\leq 1/2$ on $X$.
We set
\[
 \psi_i:=\bigl(-\log|\sigma_i|_{h_i}\bigr)^{-1}.
\]
We obtain the map
$\Psi:X\lrarr \real^{\Lambda}_{\geq 0}$.
Let $X(\ttX(\Lambda)_{\geq 0})$ denote the topological space
obtained as the fiber product
of $X$ and $\ttX(\Lambda)_{\geq 0}$
over $\real^{\Lambda}_{\geq 0}$.
There exist naturally defined continuous maps
$p_1:X(\ttX(\Lambda)_{\geq 0})\lrarr X$
and
$p_2:X(\ttX(\Lambda)_{\geq 0})
\lrarr\ttX(\Lambda)_{\geq 0}$.

There exists the orbit decomposition
(see \S\ref{subsection;22.4.26.30})
\[
 \ttX(\Lambda)_{\geq 0}=
 \bigsqcup_{\vecJ\in\nbigs(\Lambda)} \ttO(\vecJ)_{\geq 0}.
\]
We set
$X(\ttO(\vecJ)_{\geq 0}):=p_2^{-1}(\ttO(\vecJ)_{\geq 0})$.
We obtain the decomposition
\[
 X(\ttX(\Lambda)_{\geq 0})
 =\bigsqcup_{\vecJ\in\nbigs(\Lambda)}
 X(\ttO(\vecJ)_{\geq 0}).
\]
We may naturally regard
$X\setminus H=X(\ttO(\emptyset)_{\geq 0})$
as an open subset of
$X(\ttX(\Lambda)_{\geq 0})$.

Let $h_i'$ be other Hermitian metrics of $\nbigo_X(H_i)$
as above.
We obtain $\psi_i'$ and
$X'(\ttX(\Lambda)_{\geq 0})$.
\begin{lem}
There exists a natural homeomorphism
$X(\ttX(\Lambda)_{\geq 0})\simeq
X'(\ttX(\Lambda)_{\geq 0})$
whose restriction to $X\setminus H$
is the identity.
\hfill\qed 
\end{lem}

We set $\nbigx(\ttX(\Lambda)_{\geq 0})=
\cnum_{\lambda}\times X(\ttX(\Lambda)_{\geq 0})$.
A neighbourhood of $\{\lambda_0\}\times X(\ttX(\Lambda)_{\geq 0})$
is denoted by
$\nbigxzero(\ttX(\Lambda)_{\geq 0})$.
The projection
$\nbigx(\ttX(\Lambda)_{\geq 0})
\lrarr X(\ttX(\Lambda)_{\geq 0})$
is denoted by $p_{\lambda}$.
The projection
$\nbigxzero(\ttX(\Lambda)_{\geq 0})
\lrarr\nbigxzero$
is denoted by $p_1$.

If
$X=\{(z_1,\ldots,z_n)\,|\,|z_i|<1/2\}$
and $H=\bigcup_{i=1}^{\ell}\{z_i=0\}$,
we may assume
$\psi_i=(-\log|z_i|)^{-1}$,
and 
$X(\ttX(\ellsitabar)_{\geq 0})$
is identified with an open subset
of
$\ttX(\ellsitabar)_{\geq 0}
\times (S^1)^{\ell}
\times
\prod_{i=\ell+1}^n\{|z_i|<1/2\}$.
It is the inverse image of
$\closedopen{0}{(-\log 2)^{-1}}^{\ell}$
via the composition of the maps
$X(\ttX(\ellsitabar)_{\geq 0})\to
\ttX(\ellsitabar)_{\geq 0}
\to \real_{\geq 0}^{\ell}$.

In general,
for a continuous map
$f:Z\lrarr\cnum_{\lambda}$,
we set
$(\nbigo_{\cnum_{\lambda}})_Z:=
f^{-1}(\nbigo_{\cnum_{\lambda}})$.
Let $\ttD((\nbigo_{\cnum_{\lambda}})_Z)$
denote the derived category of
bounded $(\nbigo_{\cnum_{\lambda}})_Z$-complexes.

\subsection{Twistor complexes of tame harmonic bundles}

Let $(E,\delbar_E,\theta,h)$ be a tame harmonic bundle on $(X,H)$.
Let $\nbigx=\cnum_{\lambda}\times X$.
Let $p_{\lambda}:\nbigx\lrarr X$ denote the projection.
Let $\gbige$ denote the $\nbigr_X$-module
underlying the associated pure twistor $\nbigd$-module.
(See \S\ref{subsection;22.4.12.1}.)
We also have the $\nbigr_X$-module $\gbige[\ast H]$.
(See \S\ref{subsection;22.4.25.41}.)

\subsubsection{Truncated de Rham complexes}

For $\nbigm=\gbige,\gbige[\ast H]$,
we have the de Rham complexes
$\nbigm\otimes\Omegatilde^{\bullet}_{\nbigx/\cnum}$ on $\nbigx$
induced by
$\DD^f:\nbigm\lrarr\nbigm\otimes\Omegatilde^1_{\nbigx/\cnum}$
and $d:\Omegatilde^{\bullet}_{\nbigx/\cnum}
\lrarr \Omegatilde^{\bullet+1}_{\nbigx/\cnum}$.
It is a complex of $(\nbigo_{\cnum_{\lambda}})_{\nbigx}$-modules.
We obtain the truncated complexes
$\nbigctilde^{\bullet}_{\tw}(\nbigm)
\subset
\nbigm\otimes\Omegatilde^{\bullet}_{\nbigx/\cnum}$
by applying the local construction in 
\S\ref{subsection;22.2.17.100}
and \S\ref{subsection;22.3.17.30},
and gluing them.
They are $(\nbigo_{\cnum_{\lambda}})_{\nbigx}$-subcomplexes,
and the natural inclusions are quasi-isomorphisms.

As in \S\ref{subsection;22.2.17.100},
we also have the de Rham complexes
$\nbigm\otimes\Omega^{\bullet}_{\nbigx/\cnum}$ on $\nbigx$
induced by $\DD:\nbigm\lrarr\nbigm\otimes\Omega^1_{\nbigx/\cnum}$
and
$\lambda d:\Omega^{\bullet}_{\nbigx/\cnum}
\lrarr
 \Omega^{\bullet+1}_{\nbigx/\cnum}$.
We have the sub-complexes
$\nbigc^{\bullet}_{\tw}(\nbigm)
\subset
\nbigm\otimes\Omega^{\bullet}_{\nbigx/\cnum}$
corresponding to 
$\nbigctilde^{\bullet}_{\tw}(\nbigm)$
under the natural isomorphism
$\nbigm\otimes\Omega^{\bullet}_{\nbigx/\cnum}
\simeq
\nbigm\otimes\Omegatilde^{\bullet}_{\nbigx/\cnum}$.

\subsubsection{Complexes on $\nbigx(\ttX(\Lambda)_{\geq 0})$}

We obtain
the $(\nbigo_{\cnum_{\lambda}})_{\nbigx(\ttX(\Lambda)_{\geq 0})}$-complex
of sheaves
$p_1^{-1}\nbigctilde^{\bullet}_{\tw}(\nbigm)$
on $\nbigx(\ttX(\Lambda)_{\geq 0})$.
There exist natural isomorphisms
\[
 Rp_{1\ast}
 p_1^{-1}\nbigctilde^{\bullet}_{\tw}(\nbigm)
 \simeq
 p_{1\ast}
 p_1^{-1}\nbigctilde^{\bullet}_{\tw}(\nbigm)
 \simeq
 \nbigctilde^{\bullet}_{\tw}(\nbigm)
\]
in $\ttD((\nbigo_{\cnum_{\lambda}})_{\nbigx})$.
We define 
an $(\nbigo_{\cnum_{\lambda}})_{\nbigx(\ttX(\Lambda)_{\geq 0})}$-subcomplex
of 
$p^{-1}\nbigctilde^{\bullet}_{\tw}(\nbigm)$.
Let $P\in X$.
Let $\Lambda(P)=\{j\in\Lambda\,|\,P\in H_j\}$.
Note that
\[
 p_1^{-1}(P)=
 \bigsqcup_{\vecJ\in\nbigsbar(\Lambda(P))}
 \Bigl(
 p_1^{-1}(P)
 \cap
 X(\ttO(\vecJ)_{\geq 0})
 \Bigr).
\]
Let $Q\in X(\ttO(\vecJ)_{\geq 0})\cap p_1^{-1}(P)$.
There exists the natural isomorphism of the stalks
\[
 p_1^{-1}\nbigctilde^{\bullet}_{\tw}(\nbigm)_{(\lambda,Q)}
 \simeq
 \nbigctilde^{\bullet}_{\tw}(\nbigm)_{(\lambda,P)}.
\]
By using the constructions
in \S\ref{subsection;22.2.17.101} and \S\ref{subsection;22.3.17.30},
we set
\[
 \nbigw\nbigctilde^{\bullet}_{\tw}(\nbigm)_{(\lambda,Q)}
:=\lefttop{\vecJ}\nbigw_{\vecw(\vecJ)}
 \nbigctilde^{\bullet}_{\tw}(\nbigm)_{(\lambda,P)}
 \subset
  p_1^{-1}\nbigctilde^{\bullet}_{\tw}(\nbigm)_{(\lambda,Q)}.
\]
Then, we obtain the subcomplex of sheaves
\[
\nbigw\nbigctilde^{\bullet}_{\tw}(\nbigm)
\subset
p_1^{-1}\nbigctilde^{\bullet}_{\tw}(\nbigm).
\]
There exist the following naturally defined morphisms
of
$(\nbigo_{\cnum_{\lambda}})_{\nbigx(\ttX(\Lambda)_{\geq 0})}$-complexes:
\begin{equation}
\label{eq;22.2.17.102}
\begin{CD}
 \nbigw\nbigctilde^{\bullet}_{\tw}(\gbige[\ast H])
@<<<
 \nbigw\nbigctilde^{\bullet}_{\tw}(\gbige)
@>>>
 p_1^{-1}\nbigctilde^{\bullet}_{\tw}(\gbige). 
\end{CD}
\end{equation}

\begin{thm}
\label{thm;22.2.17.111}
The morphisms in {\rm(\ref{eq;22.2.17.102})}
induce the following isomorphisms
in $\ttD((\nbigo_{\cnum_{\lambda}})_{\nbigx})$:
\begin{equation}
\label{eq;22.3.17.40}
 \begin{CD}
  Rp_{1\ast} \nbigw\nbigctilde^{\bullet}_{\tw}
  (\gbige[\ast H])
@<<<
  Rp_{1\ast}\nbigw\nbigctilde^{\bullet}_{\tw}
  (\gbige)
@>>>
 \nbigctilde^{\bullet}_{\tw}(\gbige). 
 \end{CD}
\end{equation}
\end{thm}
\pf
It follows from Theorem \ref{thm;22.2.17.103}.
\hfill\qed

\vspace{.1in}
By the natural isomorphism
$p_1^{-1}\nbigc^{\bullet}_{\tw}(\nbigm)
\simeq
 p_1^{-1}\nbigctilde^{\bullet}_{\tw}(\nbigm)$,
we obtain the subcomplexes
\[
 \nbigw\nbigc^{\bullet}_{\tw}(\nbigm)
 \subset
 p_1^{-1}\nbigc^{\bullet}_{\tw}(\nbigm).
\]
There exist the naturally defined morphisms
of $(\nbigo_{\cnum_{\lambda}})_{\nbigx(\ttX(\Lambda)_{\geq 0})}$-complexes
corresponding to (\ref{eq;22.2.17.102}):
\begin{equation}
\label{eq;22.2.26.2}
\begin{CD}
 \nbigw\nbigc^{\bullet}_{\tw}(\gbige[\ast H])
@<<<
 \nbigw\nbigc^{\bullet}_{\tw}(\gbige)
@>>>
 p_1^{-1}\nbigc^{\bullet}_{\tw}(\gbige). 
\end{CD}
\end{equation}

\begin{cor}
\label{cor;22.2.26.1}
The morphisms in {\rm(\ref{eq;22.2.26.2})}
induce the following isomorphisms
in $\ttD((\nbigo_{\cnum_{\lambda}})_{\nbigx})$:
\begin{equation}
\label{eq;22.3.17.41}
 \begin{CD}
  Rp_{1\ast}\nbigw\nbigc^{\bullet}_{\tw}
  (\gbige[\ast H])
@<<<
  Rp_{1\ast}\nbigw\nbigc^{\bullet}_{\tw}
  (\gbige)
@>>>
 \nbigc^{\bullet}_{\tw}(\gbige).
 \end{CD}
\end{equation}
\hfill\qed
\end{cor}

\begin{rem}
If $(E,\delbar_E,\theta,h)$
underlies a polarized variation of Hodge structure,
$\gbige$ and $\gbige[\ast H]$ are $\cnum^{\ast}$-homogeneous.
For $\Omegatilde^{\bullet}_{\nbigx/\cnum}$,
we consider the $\cnum^{\ast}$-action 
 induced by
 $a^{\ast}(\lambda^{-p}p_{\lambda}^{-1}\tau^p)
 =a^{-p}\lambda^{-p}\tau^p$
for holomorphic $p$-forms $\tau$.
For $p_{\lambda}^{\ast}\Omega^{\bullet}_X$,
we consider the $\cnum^{\ast}$-action
induced by
 $a^{\ast}(p_{\lambda}^{-1}\tau^p)
 =a^{-p}p_{\lambda}^{-1}(\tau)$
for holomorphic $p$-forms $\tau$.
Then, the complexes of sheaves and 
the morphisms
{\rm(\ref{eq;22.2.17.102})},
{\rm(\ref{eq;22.3.17.40})},
{\rm(\ref{eq;22.2.26.2})}
and 
{\rm(\ref{eq;22.3.17.41})}
are $\cnum^{\ast}$-equivariant.
\hfill\qed
\end{rem}

\subsubsection{Specializations}
\label{subsection;22.4.2.12}

Let $\iota_{\lambda_0}:X\lrarr \nbigx$
denote the inclusion $\iota_{\lambda_0}(x)=(\lambda_0,x)$.
For an $(\nbigo_{\cnum_{\lambda}})_{\nbigx}$-module $\nbigf$,
we set
$\iota_{\lambda_0}^{\ast}\nbigf:=
\iota_{\lambda_0}^{-1}(\nbigf/(\lambda-\lambda_0)\nbigf)$.
We set
\[
 \nbigc^{\bullet}_{\tw}(\nbigm,\lambda_0):=
 \iota_{\lambda_0}^{\ast}
 \nbigc^{\bullet}_{\tw}(\nbigm).
\]
The inclusion
$X(\ttX(\Lambda)_{\geq 0})=
\{\lambda_0\}\times X(\ttX(\Lambda)_{\geq 0})
\to \nbigx(\ttX(\Lambda)_{\geq 0})$
is also denoted by $\iota_{\lambda_0}$.
We set
\[
 \nbigw\nbigc^{\bullet}_{\tw}(\nbigm,\lambda_0):=
 \iota_{\lambda_0}^{\ast}
 \nbigw\nbigc^{\bullet}_{\tw}(\nbigm).
\]
There exist the natural morphisms:
\begin{equation}
\label{eq;22.2.17.110}
 \begin{CD}
 \nbigw\nbigc^{\bullet}_{\tw}(\gbige[\ast H],\lambda_0)
  @<<<
  \nbigw\nbigc^{\bullet}_{\tw}(\gbige,\lambda_0)
@>>>
p_1^{-1}\nbigc^{\bullet}_{\tw}(\gbige,\lambda_0).
 \end{CD}
\end{equation}

\begin{cor}
\label{cor;22.2.18.10}
The morphisms {\rm(\ref{eq;22.2.17.110})}
induce the following quasi-isomorphisms:
\begin{equation}
\label{eq;22.3.6.2}
\begin{CD}
 Rp_{1\ast}
 \nbigw\nbigc^{\bullet}_{\tw}(\gbige[\ast H],\lambda_0)
  @<<<
 Rp_{1\ast}
 \nbigw\nbigc^{\bullet}_{\tw}(\gbige,\lambda_0)
@>>>
 \nbigc^{\bullet}_{\tw}(\gbige,\lambda_0).
 \end{CD}
\end{equation}
\end{cor}
\pf
Because the complexes in (\ref{eq;22.2.17.102})
are flat over
$(\nbigo_{\cnum_{\lambda}})_{\nbigx(\ttX(\Lambda)_{\geq 0})}$,
the claim follows from Theorem \ref{thm;22.2.17.111}.
\hfill\qed

\subsection{$L^2$-complexes}

\subsubsection{Preliminary}
\label{subsection;22.4.2.20}

Let $g_{X\setminus H}$ be a Poincar\'{e} like K\"ahler metric of
$X\setminus H$.
Let $j_{X\setminus H}:X\setminus H\lrarr X$
and
$j_{X\setminus H,\ttX(\Lambda)_{\geq 0}}:
X\setminus H\lrarr X(\ttX(\Lambda)_{\geq 0})$
denote the inclusions.
We have the double complex
$(\Omega^{\bullet,\bullet}_{X\setminus H},
\lambda\del_X,\delbar_X)$.
Together with $\DDlambda$ on $\nbigelambda$,
we obtain the complex
$(\Tot\nbigelambda\otimes\Omega^{\bullet,\bullet}_{X\setminus H},\DDlambda)$
on $X\setminus H$.
Let $\nbigc^{\bullet}_{L^2}(\nbigelambda,\DDlambda,h)$
denote the sheaf of sections $\tau$ of
$(j_{X\setminus H})_{\ast}
\Tot\nbigelambda\otimes\Omega^{\bullet,\bullet}_{X\setminus H}$
such that
$\tau$ and $\DDlambda\tau$ are $L^2$
with respect to $h$ and $g_{X\setminus H}$.
Let $\nbigc^{\bullet}_{L^2,\ttX(\Lambda)_{\geq 0}}(\nbigelambda,\DDlambda,h)$
denote the sheaf of sections $\tau$ of
$(j_{X\setminus H,\ttX(\Lambda)_{\geq 0}})_{\ast}
\Tot\nbigelambda\otimes\Omega^{\bullet,\bullet}_{X\setminus H}$
such that
$\tau$ and $\DDlambda\tau$ are $L^2$
with respect to $p_1^{-1}h$ and $p_1^{-1}g_{X\setminus H}$.

Let $j_{\nbigx\setminus\nbigh}:
\nbigx\setminus\nbigh\lrarr
\nbigx$
and 
$j_{\nbigx\setminus\nbigh,\ttX(\Lambda)_{\geq 0}}:
\nbigx\setminus\nbigh
\lrarr \nbigx(\ttX(\Lambda)_{\geq 0})$
denote the inclusions.
We obtain the complex
$(\Tot(\nbige\otimes p_{\lambda}^{-1}\Omega^{\bullet,\bullet}_{X\setminus H}),
\DD)$
on $\nbigx\setminus\nbigh$.
Let $\nbigc^{\bullet}_{L^2}(\nbige,\DD,h)$
denote the sheaf of sections $\tau$ of
$(j_{\nbigx\setminus\nbigh})_{\ast}
\Tot\nbige\otimes
 p_{\lambda}^{-1}\Omega^{\bullet,\bullet}_{X\setminus H}$
such that
(i) $\del_{\lambdabar}\tau=0$,
(ii) $\tau$ and $\DD\tau$ are $L^2$
with respect to
$p_{\lambda}^{-1}h$ and $g_{X\setminus H}+d\lambda\,d\lambdabar$.
Let $\nbigc^{\bullet}_{L^2,\ttX(\Lambda)_{\geq 0}}(\nbige,\DD,h)$
denote the sheaf of sections $\tau$ of
$(j_{\nbigx\setminus\nbigh,\ttX(\Lambda)_{\geq 0}})_{\ast}
\Tot\nbige\otimes p_{\lambda}^{-1}\Omega^{\bullet,\bullet}_{X\setminus H}$
such that
(i) $\del_{\lambdabar}\tau=0$,
(ii) $\tau$ and $\DD\tau$ are $L^2$
with respect to
$p_1^{-1}p_{\lambda}^{-1}h$
and $p_1^{-1}(g_{X\setminus H}+d\lambda\,d\lambdabar)$.

\begin{lem}
\label{lem;22.2.18.20}
There exist the following natural isomorphisms:
\begin{equation}
\label{eq;22.2.18.1}
 Rp_{1\ast}
 \nbigc^{\bullet}_{L^2,\ttX(\Lambda)_{\geq 0}}(\nbigelambda,\DDlambda,h)
\simeq p_{1\ast} 
 \nbigc^{\bullet}_{L^2,\ttX(\Lambda)_{\geq 0}}(\nbigelambda,\DDlambda,h)
\simeq \nbigc^{\bullet}_{L^2}(\nbigelambda,\DDlambda,h),
\end{equation}
\begin{equation}
\label{eq;22.2.18.2}
Rp_{1\ast}
\nbigc^{\bullet}_{L^2,\ttX(\Lambda)_{\geq 0}}(\nbige,\DD,h)
\simeq p_{1\ast} 
 \nbigc^{\bullet}_{L^2,\ttX(\Lambda)_{\geq 0}}(\nbige,\DD,h)
\simeq \nbigc^{\bullet}_{L^2}(\nbige,\DD,h).
\end{equation}
\end{lem}
\pf
The second isomorphisms in
(\ref{eq;22.2.18.1})
and (\ref{eq;22.2.18.2})
are clear.
Let us study the first isomorphisms.
There exist the following natural morphisms:
\begin{equation}
\label{eq;22.2.18.4}
 p_{1\ast} 
 \nbigc^{\bullet}_{L^2,\ttX(\Lambda)_{\geq 0}}(\nbigelambda,\DDlambda,h)
\lrarr
  Rp_{1\ast}
 \nbigc^{\bullet}_{L^2,\ttX(\Lambda)_{\geq 0}}(\nbigelambda,\DDlambda,h),
\end{equation}
\begin{equation}
\label{eq;22.2.18.5}
p_{1\ast}
\nbigc^{\bullet}_{L^2,\ttX(\Lambda)_{\geq 0}}(\nbige,\DD,h)
 \lrarr
 Rp_{1\ast}
\nbigc^{\bullet}_{L^2,\ttX(\Lambda)_{\geq 0}}(\nbige,\DD,h).
\end{equation}
Let us prove that they are quasi-isomorphisms.
It is enough to study the claim locally.
Hence, we may assume that
$X=\prod_{i\in\Lambda}\{|z_i|<1/2\}\times Y\subset
\cnum^{\Lambda}\times Y$ 
and $H=\bigcup_{i\in\Lambda}\{z_i=0\}$
for a complex manifold $Y$.
We may assume $|\sigma_i|_{h_i}=|z_i|$.

We choose a bijection $\Lambda\simeq \{1,\ldots,\ell\}$.
Let $(s_1,\ldots,s_{\ell})$ be the standard coordinate system of
$\real_{\geq 0}^{\ell}$.
Let $\ttU(\vecJ_{\st})_{\geq 0}$ denote the affine chart
as in \S\ref{subsection;22.2.18.3}.
The standard coordinate system
$(t_1,\ldots,t_{\ell})$
of $\ttU(\vecJ_{\st})_{\geq 0}$
is obtained as $t_i=s_i/s_{i+1}$.
For any point
$P\in \ttU(\vecJ_{\st})_{\geq 0}$
and for any relatively compact neighbourhood $A$ of $P$
in $\ttU(\vecJ_{\st})_{\geq 0}$,
there exists a $C^{\infty}$-function
$\chi_{A,P}:\ttU(\vecJ_{\st})_{\geq 0}\lrarr \real_{\geq 0}$
such that
(i) $\chi_{A,P}=1$ on a neighbourhood of $P$,
(ii) the support of $\chi_{A,P}$ is contained in $A$.

Let $Q$ be any point of $X(\ttX(\Lambda)_{\geq 0})$
such that $p_2(Q)\in \ttU(\vecJ_{\st})_{\geq 0}$.
Let $K$ be any relatively compact open neighbourhood of $Q$ in
$X(\ttX(\Lambda)_{\geq 0})$
such that
$p_2(K)\subset \ttU(\vecJ_{\st})_{\geq 0}$.
There exists a relatively compact neighbourhood $A$ of
$p_2(Q)$ in $\ttU(\vecJ_{\st})_{\geq 0}$
which is contained in $p_2(K)$.
Let $\chi_{A,p_2(Q)}$ as above.
Note that
\[
 dp_2^{\ast}(t_i)
=p_2^{\ast}(t_i)
 \left(
 \frac{1}{-\log|z_i|^2}
 \frac{d|z_i|^2}{|z_i|^2}
+\frac{1}{\log|z_{i+1}|^2}
 \frac{d|z_{i+1}|^2}{|z_{i+1}|^2}
 \right)
\]
is bounded on $K\setminus H$
with respect to the Poincar\'{e} metric.
Hence, 
$dp_2^{\ast}(\chi_{A,p_2(Q)})$ is bounded
with respect to the Poincar\'{e} metric.

For any compact subset $L_1\subset X(\ttX(\Lambda)_{\geq 0})$
and any relatively compact neighbourhood $L_2$ of $L_1$
in $X(\ttX(\Lambda)_{\geq 0})$,
there exists a function $\chi_{L_1,L_2}$
with the following property.
\begin{itemize}
 \item $\chi_{L_1,L_2}$ is continuous on $X(\ttX(\Lambda)_{\geq 0})$,
       and its restriction to $X\setminus H$ is $C^{\infty}$.
 \item $\chi_{L_1,L_2}=1$ on a neighbourhood of $L_1$.
       The support of $\chi_{L_1,L_2}$ is contained in $L_2$.
 \item $d\chi_{L_1,L_2|X\setminus H}$ is bounded
       with respect to the Poincar\'{e} metric.
\end{itemize}

Let $\tau$ be any section of
$\nbigc^{\bullet}_{L^2,\ttX(\Lambda)_{\geq 0}}(\nbigelambda,\DDlambda,h)$
on a neighbourhood $U$ of a compact subset $L$.
We choose a relatively compact neighbourhood $L'$ of $L$
in $U$.
Let $\chi_{L,L'}$ a function be as above.
Then, $\chi_{L,L'}\tau$ is a section of
$\nbigc^{\bullet}_{L^2,\ttX(\Lambda)_{\geq 0}}(\nbigelambda,\DDlambda,h)$
on $X(\ttX(\Lambda)_{\geq 0})$.
Then, it is standard that
(\ref{eq;22.2.18.4}) is a quasi-isomorphism.
We can prove that (\ref{eq;22.2.18.5}) is a quasi-isomorphism
similarly.
\hfill\qed

\vspace{.1in}
We set
$\Omegatilde^{p,q}_{(\nbigx\setminus\nbigh)/\cnum}:=
\lambda^{-p}
p_{\lambda}^{-1}(\Omega^{p,q}_{X\setminus H})$.
We obtain the double complex
$(\Omegatilde^{\bullet,\bullet}_{(\nbigx\setminus\nbigh)/\cnum},
\lambda\del_X,\delbar_X)$.
Together with $\DD^f$ on $\nbige$,
we obtain the complex
$\bigl(
\Tot(\nbige\otimes\Omegatilde^{\bullet,\bullet}_{(\nbigx\setminus\nbigh)/\cnum}),
\DD^f\bigr)$.
There exists the natural isomorphism
\begin{equation}
\label{eq;22.2.26.10}
(j_{\nbigx\setminus\nbigh})_{\ast}
(\nbige\otimes
p_{\lambda}^{-1}\Omega^{\bullet,\bullet}_{X\setminus H})
\simeq
 (j_{\nbigx\setminus\nbigh})_{\ast}\Bigl(
 \nbige\otimes\Omegatilde^{\bullet,\bullet}_{(\nbigx\setminus\nbigh)/\cnum}
 \Bigr)
\end{equation}
induced by $\tau^{p,q}\longmapsto \lambda^{-p}\tau^{p,q}$.
Let $\nbigctilde^{\bullet}_{L^2}(\nbige,\DD^f,h)
\subset
 (j_{\nbigx\setminus\nbigh})_{\ast}\Bigl(
 \Tot(\nbige\otimes\Omegatilde^{\bullet,\bullet}
 _{(\nbigx\setminus\nbigh)/\cnum})
 \Bigr)$
denote the subcomplex corresponding to 
$\nbigc^{\bullet}_{L^2}(\nbige,\DD,h)$
under (\ref{eq;22.2.26.10}).
Similarly,
we obtain the subcomplex
\[
 \nbigctilde^{\bullet}_{L^2,\ttX(\Lambda)_{\geq 0}}(\nbige,\DD^f,h)
\subset
 (j_{\nbigx\setminus\nbigh,\ttX(\Lambda)_{\geq 0}})_{\ast}\Bigl(
 \Tot(\nbige\otimes\Omegatilde^{\bullet,\bullet}
 _{(\nbigx\setminus\nbigh)/\cnum})
 \Bigr)
\]
corresponding to
$\nbigc^{\bullet}_{L^2,\ttX(\Lambda)_{\geq 0}}(\nbige,\DD,h)
\subset
 (j_{\nbigx\setminus\nbigh,\ttX(\Lambda)_{\geq 0}})_{\ast}\Bigl(
 \Tot(\nbige\otimes p_{\lambda}^{-1}\Omega^{\bullet,\bullet}
 _{X\setminus H}) \Bigr)$.

\subsubsection{Statements}

By the construction,
there exists the following natural quasi-isomorphism
of $\cnum_{X\setminus H}$-complexes:
\begin{equation}
\label{eq;22.1.24.2}
 \nbigw\nbigc_{\tw}^{\bullet}(\gbige[\ast H],\lambda)
  _{|X\setminus H}
 \lrarr
 \nbigc_{L^2,\ttX(\Lambda)_{\geq 0}}^{\bullet}(\nbigelambda,\DDlambda,h)
   _{|X\setminus H}.
\end{equation}

We shall prove the following theorem
in \S\ref{subsection;22.3.17.51}
after the preliminaries 
in \S\ref{subsection;22.3.17.50}--\S\ref{subsection;22.3.17.52}.

\begin{thm}
\label{thm;22.2.18.11}
The morphism {\rm(\ref{eq;22.1.24.2})}
extends to the following morphism
\begin{equation}
  \nbigw\nbigc_{\tw}^{\bullet}(\gbige[\ast H],\lambda)
 \lrarr
 \nbigc_{L^2,\ttX(\Lambda)_{\geq 0}}^{\bullet}(\nbigelambda,\DDlambda,h).
\end{equation}
Moreover, the induced morphism
\begin{equation}
 Rp_{1\ast}\nbigw\nbigc_{\tw}^{\bullet}(\gbige[\ast H],\lambda)
 \lrarr
 Rp_{1\ast}\nbigc_{L^2,\ttX(\Lambda)_{\geq 0}}^{\bullet}(\nbigelambda,\DDlambda,h)
\end{equation}
in $\ttD(\cnum_X)$
is an isomorphism.
\end{thm}

Similarly,
there exists the following natural quasi-isomorphism
of $(\nbigo_{\cnum_{\lambda}})_{\nbigx\setminus\nbigh}$-complexes
\begin{equation}
\label{eq;22.1.24.1}
 \nbigw\nbigc_{\tw}^{\bullet}(\gbige[\ast H])
  _{|\nbigx\setminus \nbigh}
 \lrarr
 \nbigc_{L^2,\ttX(\Lambda)_{\geq 0}}^{\bullet}(\nbige,\DD,h)
   _{|\nbigx\setminus\nbigh}.
\end{equation}

We shall prove the following theorem
in \S\ref{subsection;22.3.17.51}.

\begin{thm}
\label{thm;21.12.21.11}
The morphism {\rm(\ref{eq;22.1.24.1})}
extends to the following morphism of 
$(\nbigo_{\cnum_{\lambda}})_{\nbigx(\ttX(\Lambda)_{\geq 0})}$-complexes:
 \begin{equation}
\label{eq;22.3.17.60}
  \nbigw\nbigc_{\tw}^{\bullet}(\gbige[\ast H])
 \lrarr
 \nbigc_{L^2,\ttX(\Lambda)_{\geq 0}}^{\bullet}(\nbige,\DD,h).
 \end{equation}
Moreover, the induced morphism
\begin{equation}
\label{eq;22.3.17.61}
 Rp_{1\ast}\nbigw\nbigc_{\tw}^{\bullet}(\gbige[\ast H])
 \lrarr
 Rp_{1\ast}\Bigl(
 \nbigc_{L^2,\ttX(\Lambda)_{\geq 0}}^{\bullet}(\nbige,\DD,h)
 \Bigr)
 \end{equation}
in $\ttD((\nbigo_{\cnum_{\lambda}})_{\nbigx})$
is an isomorphism.
\end{thm}

\begin{cor}
\label{cor;22.3.17.100}
As a result
of Corollary {\rm\ref{cor;22.2.18.10}}
and Theorem {\rm\ref{thm;22.2.18.11}},
we obtain the following sequence of
isomorphisms in $\ttD(\cnum_X)$:
\begin{equation}
\label{eq;22.1.24.11}
 \nbigc_{\tw}^{\bullet}(\gbige,\lambda)
 \llarr
 Rp_{1\ast}
 \nbigw\nbigc_{\tw}^{\bullet}(\gbige,\lambda)
 \lrarr
 Rp_{1\ast}
 \nbigw\nbigc_{\tw}^{\bullet}(\gbige[\ast H],\lambda)
 \lrarr
 \nbigc^{\bullet}_{L^2}(\nbigelambda,\DDlambda,h).
\end{equation}
We also obtain the following sequence of isomorphisms
in $\ttD((\nbigo_{\cnum_{\lambda}})_{\nbigx})$
as a consequence of Theorem {\rm\ref{thm;22.2.17.111}}
and Theorem {\rm\ref{thm;21.12.21.11}}:
\begin{equation}
\label{eq;22.3.17.63}
 \nbigc_{\tw}^{\bullet}(\gbige)
 \llarr
 Rp_{1\ast}
 \nbigw\nbigc_{\tw}^{\bullet}(\gbige)
 \lrarr
 Rp_{1\ast}
 \nbigw\nbigc_{\tw}^{\bullet}(\gbige[\ast H])
 \lrarr
 \nbigc^{\bullet}_{L^2}(\nbige,\DD,h).
\end{equation}
By the natural twisting, we obtain the following isomorphisms
in $\ttD((\nbigo_{\cnum_{\lambda}})_{\nbigx})$:
\begin{equation}
\label{eq;22.3.17.64}
 \nbigctilde_{\tw}^{\bullet}(\gbige)
 \llarr
 Rp_{1\ast}
 \nbigw\nbigctilde_{\tw}^{\bullet}(\gbige)
 \lrarr
 Rp_{1\ast}
 \nbigw\nbigctilde_{\tw}^{\bullet}(\gbige[\ast H])
 \lrarr
 \nbigctilde^{\bullet}_{L^2}(\nbige,\DD^f,h).
\end{equation}
\hfill\qed
\end{cor}
 
\begin{rem}
We consider
the $\cnum^{\ast}$-action on $\Omegatilde^{p,q}$
induced by
$a^{\ast}(\lambda^{-p}p_{\lambda}^{-1}(\tau^{p,q}))
=a^{-p}\lambda^{-p}p_{\lambda}^{-1}(\tau^{p,q})$
for $(p,q)$-forms $\tau^{p,q}$.
We consider  
the $\cnum^{\ast}$-action on $p_{\lambda}^{-1}\Omega^{p,q}$
induced by
$a^{\ast}(p_{\lambda}^{-1}(\tau^{p,q}))
=a^{-p}p_{\lambda}^{-1}(\tau^{p,q})$
for $(p,q)$-forms $\tau^{p,q}$.
 If $(E,\delbar_E,\theta,h)$ underlies a polarized variation of
Hodge structure,
the complex of sheaves
 $\nbigctilde^{\bullet}_{L^2}(\nbige,\DD^f,h)$,
$\nbigctilde^{\bullet}_{\ttX(\Lambda)_{\geq 0},L^2}(\nbige,\DD^f,h)$
$\nbigc^{\bullet}_{\ttX(\Lambda)_{\geq 0},L^2}(\nbige,\DD,h)$
are naturally $\cnum^{\ast}$-equivariant,
and the morphisms
{\rm(\ref{eq;22.3.17.60})},
{\rm(\ref{eq;22.3.17.61})},
{\rm(\ref{eq;22.3.17.63})}
and
{\rm(\ref{eq;22.3.17.64})}
are $\cnum^{\ast}$-equivariant.
\hfill\qed
\end{rem}

\subsection{Dolbeault's lemma in the relatively one dimensional case}
\label{subsection;22.3.17.50}

\subsubsection{Setting}

Let $S$ be an open subset of $\real^N$.
We consider a Riemannian metric of $S$,
and let $\mu_S$ be the induced volume form.
Let $\mu_{S\times\Delta^{\ast}}$ denote the measure of
$S\times\Delta^{\ast}$ induced by
$\mu_S$ and $|dz\,d\zbar|$.

Let $a:S\lrarr \real$ be a $C^{\infty}$-function
such that either one of the following holds.
\begin{itemize}
 \item $a$ is a constant integer.
 \item There exists an integer $[a]$
       such that $[a]<a(x)<[a]+1$
       for any $x\in S$.
       In this case, we say $[a]<a<[a]+1$.
\end{itemize}
Let $k\in\real$ .
Let $\nbigl^2(a,k)$ denote the space of
measurable functions $f$ on $S\times\Delta^{\ast}$
such that
\[
 \int_{S\times\Delta^{\ast}}
 |f|^2 |z|^{2a}(-\log|z|+1)^{k}
 \,\frac{\mu_{S\times\Delta^{\ast}}}{|z|^2}
 <\infty.
\]
Let $z=re^{\sqrt{-1}\theta}$ be the polar decomposition.
For $f\in\nbigl^2(a,k)$,
there exists the Fourier expansion
\[
 f=\sum_{n\in\seisuu}f_n(x,r)e^{\sqrt{-1}n\theta}.
\]
If $a=0$, we set
\[
 \nbigl^2(0,k)^{\circ}:=
\bigl\{
 f\in\nbigl^2(0,k)\,\big|\,
 f_0=0
 \bigr\},
 \quad
 \nbigl^2(0,k)^{\bot}:=
\bigl\{
 f\in\nbigl^2(0,k)\,\big|\,
 f=f_0
\bigr\}.
\]

\subsubsection{Solvability of $\delbar$-equation}
\label{subsection;22.1.26.1}

For $f\in \nbigl^2(a,k)$,
we shall study the solvability of
the equation $\zbar\del_{\zbar}u=f$
in $f\in \nbigl^2(a,k-2)$
by applying the method of Zucker \cite{z}.

\begin{lem}
\label{lem;22.1.25.20}
Let $f\in\nbigl^2(a,k)$ $(0<a<1)$
 (resp. $f\in \nbigl^2(0,k)^{\bot}$).
We set
\[
 \nbiga(f)(x,z):=\frac{1}{\pi}\int_{|\zeta|<1}
 \frac{f(x,\zeta)}{(z-\zeta)\zetabar}
 \frac{\sqrt{-1}}{2}d\zeta\,d\zetabar.
\]
Then, $\nbiga(f)$ is an element of
$\nbigl^2(a,k-2)$
(resp. $\nbigl^2(0,k-2)^{\bot}$),
and it satisfies
$\zbar\del_{\zbar}\nbiga(f)=f$.
\end{lem}
\pf
For the Fourier expansion
$\nbiga(f)_n=\sum \nbiga(f)_n(x,r)e^{\sqrt{-1}n\theta}$,
we have
\[
 \nbiga(f)_n(x,r)
 =\left\{
\begin{array}{ll}
 2r^n\int_{0}^r\rho^{-n-1}f_{n}(x,\rho)\,d\rho &
  (\mbox{if $n<0$})
  \\
 -2r^n\int_r^{1}\rho^{-n-1}f_{n}(x,\rho)\,d\rho &
  (\mbox{if $n\geq 0$}).
\end{array}
 \right.
\]
By the argument in \cite{z},
we obtain the claim of the lemma.
\hfill\qed

\vspace{.1in}

Similarly, we obtain the following
by the argument in \cite{z}.
\begin{lem}
Let $f\in \nbigl^2(0,k)^{\circ}$,
i.e., $f=f_0$.
We set
\[
 \nbiga(f)=\left\{
\begin{array}{ll}
 2\int_{0}^rf(x,\rho)\,\rho^{-1}d\rho &
  (\mbox{if $k>1$})
  \\
 -2\int_r^{1}f(x,\rho)\,\rho^{-1}d\rho &
  (\mbox{if $k<1$}).
\end{array}
 \right.
\]
Then, $\nbiga(f)\in\nbigl^2(0,k-2)$,
and $\zbar\del_{\zbar}\nbiga(f)=f$. 
\hfill\qed
\end{lem}

Let $R:S\to \openopen{0}{1}$ be a continuous function.
We set
\[
 S(\leq R):=
 \bigl\{(x,z)\in S\times \Delta^{\ast}\,\big|\,|z|\leq R(x)\bigr\},
 \quad\quad
 S(\geq R):=
 \bigl\{(x,z)\in S\times \Delta^{\ast}\,\big|\,|z|\geq R(x)
 \bigr\}.
\]
The following lemma is obvious.
\begin{lem}
\label{lem;22.4.27.30}
Suppose $f\in\nbigl^{2}(a,k)$ $(0<a<1)$ or $f\in\nbigl^{2}(0,k)^{\bot}$.
\begin{itemize}
 \item If $f_n=0$ for $n\geq 0$,
       and if the support of $f$ is contained in
       $S(\geq R)$,
       then the support of $\nbiga(f)$ is contained in
       $S(\geq R)$.
 \item If $f_n=0$ for $n<0$,
       and if the support of $f$ is contained in
       $S(\leq R)$,
       then the support of $\nbiga(f)$
       is contained in
       $S(\leq R)$.
\end{itemize}
Suppose
$f\in\nbigl^2(0,k)^{\circ}$ $(k\neq 1)$.
\begin{itemize}
 \item If $k>1$ and
       if the support of $f$ is contained in
       $S(\geq R)$,
       then the support of $\nbiga(f)$ is contained in
       $S(\geq R)$.
 \item If $k<1$ and
       if the support of $f$ is contained in
       $S(\leq R)$,
       then the support of $\nbiga(f)$
       is contained in $S(\leq R)$.
 \hfill\qed
\end{itemize}
\end{lem}

We obtain the following lemma from Lemma \ref{lem;22.1.25.20}.
\begin{lem}
If $[a]<a<[a]+1$,
for $f\in \nbigl^2(a,k)$,
we set $\nbiga^{[a]}(f):=z^{-[a]}\nbiga(z^{[a]}f)$.
Then, 
$\nbiga^{[a]}(f)\in\nbigl^2(a,k-2)$
and
$\zbar\del_{\zbar}\nbiga^{[a]}(f)=f$.
\hfill\qed
\end{lem}

\subsubsection{Additional conditions}

Suppose that $S=S_1\times S_2\subset
\real^{N_1}\times\real^{N_2}=\real^N$.
Let $P_i$ $(i=1,\ldots,m)$ be linear differential operators on $S$
of the form
$P_i=P_{i,0}+\sum_{j=1}^{N_1} P_{i,j}\del_{x_j}$,
where $P_{i,j}$ are $C^{\infty}$-functions on $S_2$,
and $(x_1,\ldots,x_{N_1},x_{N_1+1},\ldots,x_N)$
denotes the standard coordinate system of $\real^N$.
Let
$\lambda$ and $\alpha$ be $C^{\infty}$-functions on $S_2$.

Suppose $0<a<1$ or $a=0$.
Let $k\in\real$.
Let $g_j\in \nbigl^2(a,k-2)$ $(j=1,2)$
and $h_i\in\nbigl^2(a,k)$ $(i=1,\ldots,m)$.
If $(a,k)=(0,1)$ or $(a,k)=(0,3)$,
we assume that 
$g_j\in \nbigl^2(0,k-2)^{\bot}$
and $h_i\in\nbigl^2(0,k)^{\bot}$.
We assume that the following equality holds
as distributions on $S\times\Delta^{\ast}$:
\begin{equation}
\label{eq;22.2.20.1}
 \zbar\del_{\zbar}g_1
 +(\lambda z\del_z+\alpha)g_2
 +\sum_{i=1}^m P_ih_i=0.
\end{equation}
We also assume that
there exists a continuous function $R:S\to \openopen{0}{1}$
such that
the supports of $g_j$ and $h_i$ are contained in
$S(\leq R)$.

\begin{prop}
\label{prop;22.1.25.21}
The following equality holds on $S\times\Delta^{\ast}$
as distributions:
\begin{equation}
\label{eq;22.2.19.1}
 g_1
 +(\lambda z\del_z+\alpha)\nbiga(g_2)
 +\sum_{i=1}^m P_i\nbiga(h_i)
 =0.
\end{equation} 
 In particular,
$(\lambda z\del_z+\alpha)\nbiga(g_2)
 +\sum_{i=1}^m P_i\nbiga(h_i)
 \in \nbigl^2(a,k-2)$.
\end{prop}
\pf
Let us begin with an easy case.
\begin{lem}
\label{lem;22.2.24.20}
Suppose that
(i) the supports of $g_j$ and $h_i$ are proper over $S$,
 (ii) $(\lambda z\del_z+\alpha)g_j$,
 $(\lambda z\del_z+\alpha)\nbiga(g_j)$,
$P_i(h_i)$ and $P_i(\nbiga(h_i))$ are also $L^2$.
Then, {\rm(\ref{eq;22.2.19.1})} holds. 
\end{lem}
\pf
We use the Fourier expansions
$g_j=\sum g_{j,n}e^{\sqrt{-1}n\theta}$
and
$h_{i}=\sum h_{i,n}e^{\sqrt{-1}n\theta}$.
We obtain the following equality for each $n$
from (\ref{eq;22.2.20.1}):
\[
 \frac{1}{2}(r\del_r-n)g_{1,n}
 +\frac{\lambda}{2}(r\del_r+n)g_{2,n}
 +\alpha g_{2,n}
 +\sum_{i=1}^m P_ih_{i,n}=0.
\]
Let $n<0$.
Because the supports of $g_{j,n}$ and $h_{i,n}$
are proper over $S$,
we obtain
\begin{equation}
 g_{1,n}(x,r)
 +\lambda g_{2,n}(x,r)
 +(\lambda n+\alpha)\nbiga(g_2)_n(x,r)
 +\sum_{i=1}^m \nbiga(P_ih_i)_n(x,r)=0.
\end{equation}
We also have
\[
\bigl(
 (\lambda z\del_z+\alpha)\nbiga(g_2)
 \bigr)_n(x,r)
 =\lambda g_{2,n}(x,r)+(\lambda n+\alpha)\nbiga(g_2)_n(x,r),
\quad\quad
 \bigl(
 P_i\nbiga(h_{i})
 \bigr)_n(x,r)
 =
  \nbiga(P_ih_i)_n(x,r).
\]
Hence, we obtain
\begin{equation}
 g_{1,n}(x,r)
+\bigl((\lambda z\del_z+\alpha)\nbiga(g_2)\bigr)_n(x,r)
 +\sum_{i=1}^m\bigl(P_i\nbiga(h_i)\bigr)_n(x,r)=0.
\end{equation}
Similarly, for $n\geq 0$,
we obtain
\begin{equation}
 g_{1,n}(x,r)
+\bigl((\lambda z\del_z+\alpha)\nbiga(g_2)\bigr)_n(x,r)
 +\sum_{i=1}^m\bigl(P_i\nbiga(h_i)\bigr)_n(x,r)=0.
\end{equation}
Thus, we obtain (\ref{eq;22.2.19.1}).
\hfill\qed

\begin{lem}
\label{lem;22.2.22.10}
Suppose that the supports of $g_j$ and $h_i$
are proper over $S$.
Then, {\rm(\ref{eq;22.2.19.1})}
holds on $S\times\Delta^{\ast}$
as distributions.
\end{lem}
\pf
Let $P$ be any point of $S$.
Let $S_{P}$ be a relatively compact open neighbourhood
of $P$ in $S$.
By using a convolution in the $(S_1\times\Delta^{\ast})$-direction,
we can construct sequences of functions
$\{g_j^{(p)}\}_{p\geq 1}$ $(j=1,2)$ and
$\{h_i^{(p)}\}_{p\geq 1}$ $(i=1,\ldots,m)$
on $S_{P}\times\Delta^{\ast}$
such that the following holds.
\begin{itemize}
 \item $g_j^{(p)}$ and $h_i^{(p)}$
       satisfy the assumption in Lemma \ref{lem;22.2.24.20}.
 \item For each $p$,
       the equation (\ref{eq;22.2.20.1})
       holds for
       $g_j^{(p)}$ and $h^{(p)}_i$.
 \item The sequences
       $\{g_j^{(p)}\}$
       are convergent to $g_j$
       in
       $\nbigl^2(a,k-2)$,       
       and $\{h^{(p)}_i\}$
       are convergent to
       $h_i$
       in $\nbigl^2(a,k)$ on $S_P\times\Delta^{\ast}$.
\end{itemize}
By Lemma \ref{lem;22.2.24.20}, for each $p$,
{\rm(\ref{eq;22.2.19.1})} holds
for $g_j^{(p)}$ and $h^{(p)}_i$
on $S_P\times\Delta^{\ast}$.
Because
$\{\nbiga(h^{(p)}_i)\}_{p\geq 1}$
are convergent to $\nbiga(h_i)$
in $\nbigl^2(a,k-2)$,
and 
$\{\nbiga(g^{(p)}_2)\}_{p\geq 1}$
are convergent to $\nbiga(g_2)$
in $\nbigl^2(a,k-4)$,
we obtain that
$\{P_i\nbiga(h^{(p)}_i)\}$
and
$\{(\lambda z\del_z+\alpha)\nbiga(g^{(p)}_2)\}$
are convergent to $P_i\nbiga(h_i)$
and
$(\lambda z\del_z+\alpha)\nbiga(g_2)$, respectively,
as distributions.
Hence, we obtain
(\ref{eq;22.2.19.1})
for $g$ and $h_i$ on $S_P\times\Delta^{\ast}$.
\hfill\qed

\vspace{.1in}
Let us complete the proof of Proposition \ref{prop;22.1.25.21}.
Let $\kappa$ be a $C^{\infty}$-function
$\real\lrarr\real$
such that
(i) $0\leq \kappa(t)\leq 1$ for any $t$,
(ii) $\kappa(t)=1$ $(t\leq 1/2)$ and $\kappa(t)=0$ $(t\geq 1)$.
For a large $N$,
we set
\[
\chi_{1,N}(z)=\kappa\Bigl(\frac{1}{N}\log(-\log|z|^2)\Bigr),
\quad
\chi_{2,N}(z)=
-\frac{1}{N}\frac{1}{\log|z|^2}
\kappa'\Bigl(\frac{1}{N}\log(-\log|z|^2)\Bigr).
\]
There exists $C>0$ such that
\[
|\chi_{2,N}(z)|
\leq C(-\log |z|^2)^{-1}(\log(-\log|z|^2))^{-1}
\]
We set
$g_{j}^{(N)}=\chi_{1,N}\cdot g_j$,
$h^{(N)}_{i}:=\chi_{1,N}\cdot h_i$
and
$k^{(N)}:=\chi_{2,N}(g_1+g_2)$.
We obtain
\[
\zbar\del_{\zbar}g^{(N)}_1
+(\lambda z\del_z+\alpha)g^{(N)}_2
+\sum_{i=1}^m
P_i h_i^{(N)}
+k^{(N)}=0.
\]
By the previous consideration,
we obtain
\[
 g^{(N)}_1
+(\lambda z\del_z+\alpha)\nbiga(g^{(N)}_2)
+\sum_{i=1}^m
 P_i\nbiga(h_i^{(N)})
+\nbiga(k^{(N)})=0.
\]
The sequences
$\{g_j^{(N)}\}$ are convergent to $g_j$
in $\nbigl^2(a,k-2)$.
The sequences
$\{h^{(N)}_{i}\}$ are convergent to $h_i$
in $\nbigl^2(a,k)$.
The sequence 
$\{k^{(N)}\}$ is convergent to $0$
in $\nbigl^2(a,k)$.
Hence, we obtain 
{\rm(\ref{eq;22.2.19.1})}
for $g_j$ and $h_i$.
\hfill\qed

\subsubsection{A variant}

Suppose $-1<a<0$.
Let $k,k',k''\in\real$.
Let $g_1\in\nbigl^2(0,k-2)$,
$g_2\in \nbigl^2(a,k')$
and $h_i\in\nbigl^2(a,k'')$ $(i=1,\ldots,m)$.
We assume that the equality (\ref{eq;22.2.20.1})
holds for $g_j$ and $h_i$
as distributions on $S\times\Delta^{\ast}$.
We also assume that
the supports of $g_j$ and $h_i$ are contained in
$S(\leq R)$
for a continuous function $R:S\to \openopen{0}{1}$.

\begin{prop}
\label{prop;22.2.22.12}
If $k\geq 1$,
the following equality holds on $S\times\Delta^{\ast}$
as distributions:
\begin{equation}
\label{eq;22.2.22.22}
 g_1
 +(\lambda z\del_z+\alpha)\nbiga^{[-1]}(g_2)
 +\sum_{i=1}^m P_i\nbiga^{[-1]}(h_i)
 =0.
\end{equation} 
 In particular,
$(\lambda z\del_z+\alpha)\nbiga^{[-1]}(g_2)
 +\sum_{i=1}^m P_i\nbiga^{[-1]}(h_i)
 \in \nbigl^2(0,k-2)$.

If $k<1$, {\rm(\ref{eq;22.2.19.1})} holds
on $S\times\Delta^{\ast}$
as distributions.
\end{prop}
\pf
Because $g_{j}\in \nbigl^2(0,k-2)$
and $h_{i}\in\nbigl^2(0,k)$,
the second claim directly follows from
Proposition \ref{prop;22.1.25.21}.
Let us indicate an outline of the proof of the first claim.
The following lemma is similar to
Lemma \ref{lem;22.2.24.20}.

\begin{lem}
Suppose that $g_j$ and $h_i$
satisfy the condition in Lemma {\rm\ref{lem;22.2.24.20}}.
Then, {\rm(\ref{eq;22.2.22.22})} holds.
\hfill\qed
\end{lem}

The following lemma is similar to
Lemma \ref{lem;22.2.22.10}.
\begin{lem}
If the supports of $g_j$ and $h_i$
are proper over $S$,
then {\rm(\ref{eq;22.2.22.22})}
holds on $S\times\Delta^{\ast}$
as distributions.
\hfill\qed
\end{lem}

Let $g_j^{(N)}$,
$h_i^{(N)}$ and $k^{(N)}$
be as in the proof of
Proposition \ref{prop;22.1.25.21}.
By the previous consideration,
we obtain
\[
 g^{(N)}_1
+(\lambda z\del_z+\alpha)\nbiga^{[-1]}(g^{(N)}_2)
+\sum_{i=1}^m
 P_i\nbiga^{[-1]}(h_i^{(N)})
+\nbiga^{[-1]}(k^{(N)})=0.
\]

The sequence
$\{g^{(N)}_1\}$
is convergent to $g_1$ in $\nbigl^2(0,k-2)$.
The sequence
$\{\nbiga^{[-1]}(g^{(N)}_2)\}$
is convergent to $g_2$ in $\nbigl^2(a,k'-2)$.
The sequences
$\{\nbiga^{[-1]}(h^{(N)}_i)\}$
are convergent to $h_i$ in $\nbigl^2(a,k''-2)$.
We obtain Proposition \ref{prop;22.2.22.12}
from the next lemma.
\begin{lem}
\label{lem;22.2.22.13}
The sequence
$\{\nbiga^{[-1]}(k^{(N)})\}$
is convergent to $0$
in $\nbigl^2(0,k-2)$. 
In particular, the sequence is convergent to $0$
as distributions.
\end{lem}
\pf
By using the Fourier expansion
$k^{(N)}=\sum k^{(N)}_ne^{\sqrt{-1}m\theta}$,
we decompose 
$k^{(N)}=k^{(N)}_0+k^{(N)}_{\neq 0}$,
where $k^{(N)}_{\neq 0}=\sum_{m\neq 0} k^{(N)}_me^{\sqrt{-1}m\theta}$.
It is standard and easy to see that
the sequence $\{\nbiga^{[-1]}(k^{(N)}_{\neq 0})\}$
is convergent in $\nbigl^2(0,k-2)$.
There exists $C>0$ such that the following holds:
\begin{equation}
\label{eq;22.4.26.31}
 |k^{(N)}_0(x,z)|
 \leq
C \bigl(|g_{1,0}(x,z)|+|g_{2,0}(x,z)|\bigr)
 \cdot
 (-\log|z|^2)^{-1}
 \bigl(
 \log(-\log |z|^{2})
 \bigr)^{-1}.
\end{equation}
We have
\[
 \nbiga^{[-1]}(k^{(N)}_0)
=\int_0^r k^{(N)}_0\,\frac{d\rho}{\rho}.
\]
We obtain
\begin{equation}
\label{eq;22.4.26.30}
 \int_{S\times\{|z|<1/2\}} \bigl|
 \nbiga^{[-1]}(k^{(N)}_0)
 \bigr|^2
 (-\log |z|)^{k-2}
 \frac{\mu_{S\times\Delta^{\ast}}}{|z|^2}
=\int_S\mu_S
\int_{0}^{1/2}(-\log r)^{k-2}\frac{dr}{r}
\left|
 \int_0^rk_0^{(N)}\frac{d\rho}{\rho}
\right|^2.
\end{equation}
If $k>1$,
there exists $C_1>0$ such that the following holds
for any $N$:
\[
\left|
 \int_0^rk_0^{(N)}\frac{d\rho}{\rho}
 \right|^2
 \leq
 C_1\int_0^r|k_0^{(N)}|^2
 (-\log \rho)^k\bigl(
  \log(-\log\rho)
  \bigr)^2
 \frac{d\rho}{\rho}
 \cdot
 \frac{1}{(-\log r)^{k-1}}
 \frac{1}{\bigl(\log(-\log r)\bigr)^2}.
\]
Hence, there exists $C_2>0$ such that
(\ref{eq;22.4.26.30})
is dominated by
\[
 \int_{S}\mu_S
 \int_{0}^{1/2}
 |k_0^{(N)}|^2(-\log\rho)^k\bigl(\log(-\log\rho)\bigr)^2
 \frac{d\rho}{\rho}.
\]
By (\ref{eq;22.4.26.31})
and Lebesgue's theorem,
we obtain that (\ref{eq;22.4.26.30}) goes to $0$
as $N\to\infty$.
If $k=1$, we take $0<\epsilon<1$,
and then there exists $C_3>0$
such that the following holds for any $N$:
\[
\left|
 \int_0^rk_0^{(N)}\frac{d\rho}{\rho}
 \right|^2
 \leq
 C_3\int_0^r|k_0^{(N)}|^2
 (-\log \rho)\bigl(
  \log(-\log\rho)
  \bigr)^{2-\epsilon}
 \frac{d\rho}{\rho}
 \cdot
 \frac{1}{\bigl(\log(-\log r)\bigr)^{1-\epsilon}}.
\]
There exists $C_4,C_5>0$ such that
(\ref{eq;22.4.26.30})
is dominated by
\begin{multline}
 C_4
 \int_S\mu_S
 \int_0^{1/2}
 (-\log r)^{-1}
 \bigl(
 \log(-\log r)
 \bigr)^{-1+\epsilon}
 \frac{dr}{r}
 \left(
 \int_{0}^{r}
 |k_0^{(N)}|^2(-\log\rho)
 \bigl(
 \log(-\log\rho)\bigr)^{2-\epsilon}\frac{d\rho}{\rho}
 \right)
 \leq
 \\
 C_5\int_S\mu_S
 \int_{0}^{1/2}
 |k_0^{(N)}|^2(-\log\rho)
 \bigl(
 \log(-\log\rho)\bigr)^{2}
 \frac{d\rho}{\rho}.
\end{multline}
By (\ref{eq;22.4.26.31})
and Lebesgue's theorem,
we obtain that (\ref{eq;22.4.26.30}) goes to $0$
as $N\to\infty$.
Thus, we obtain
Lemma \ref{lem;22.2.22.13}
and hence Proposition \ref{prop;22.2.22.12}.
\hfill\qed

\subsection{Estimates for the connection form}
\label{subsection;22.2.22.24}

\subsubsection{Some notation}

Let $X=\Delta^n$ and $H=\bigcup_{i=1}^{\ell}\{z_i=0\}$.
We set $H_i=\{z_i=0\}$.
For any subset $I\subset\nbar$,
we set $H_I:=\bigcap_{i\in I}H_i$.
We consider the $(S^1)^{\ell}$-action on $(X,H)$
induced by
$(a_1,\ldots,a_{\ell})\cdot (z_1,\ldots,z_n)
=(a_1z_1,\ldots,a_{\ell}z_{\ell},z_{\ell+1},\ldots,z_n)$.
We use the Poincar\'{e} like metric $g_{X\setminus H}$
of $X\setminus H$.

For any $J,K\subset\nbar=\{1,\ldots,n\}$,
we set
\[
 \omega_{J,K}:=
 \bigwedge_{\substack{j\in J\\ j\leq \ell}}
 \frac{dz_j}{z_j}
 \wedge
 \bigwedge_{\substack{j\in J \\ j>\ell}}
  dz_j
  \wedge
 \bigwedge_{\substack{j\in K\\ j\leq \ell}}
 \frac{d\zbar_j}{\zbar_j}
 \wedge
 \bigwedge_{\substack{j\in K \\ j>\ell}}
  d\zbar_j.
\]
Here, we use the naturally induced order on $J$ and $K$
for the exterior products.

Let $(r_j,\theta_j)$ denote the polar coordinate
induced by $z_j=r_je^{\sqrt{-1}\theta_j}$.
For $I\subset \ellsitabar$
and $\vecp\in\seisuu^I$,
we set
$\vecp\cdot\vectheta=\sum_{i\in I} p_i\theta_i$.
For an $L^2$-function $f$ on an $(S^1)^{\ell}$-invariant
open subset $U$ of $X\setminus H$,
let
\[
 f=\sum_{\vecp\in\seisuu^I}
 \gbigf_{I,\vecp}(f)e^{\sqrt{-1}\vecp\vectheta}
\]
denote the Fourier expansion with respect to
$(\theta_j)_{j\in I}$.

Let $\gbigs_{\ell}$
denote the $\ell$-th symmetric group.
For $\vecC=(C_i\,|\,i=1,\ldots,\ell)\in\real_{>0}^{\ell}$
and for $\sigma\in\gbigs_{\ell}$,
let us consider the region
\[
 Z(\vecC,\sigma)=
 \Bigl\{
 (z_1,\ldots,z_n)\,\Big|\,
 \frac{-\log|z_{\sigma(i+1)}|}{-\log|z_{\sigma(i)}|}<C_i
 \,\,\,(i=1,\ldots,\ell-1)
 \Bigr\}.
\]
Let $\Zbar(\vecC,\sigma)$
denote the closure in $X(\ttX(\ellsitabar)_{\geq 0})$.
If $C_i>1$, we obtain the covering
\[
 X\setminus H=\bigcup_{\sigma\in\gbigs_{\ell}}
 Z(\vecC,\sigma).
\]

We obtain the following lemma
by using an argument in the proof of Lemma \ref{lem;22.2.18.20}.
\begin{lem}
\label{lem;22.2.20.3}
Let $S\subset\gbigs_{\ell}$  be a subset.
We set
 $(X\setminus H)_S:=\bigcup_{\sigma\in S}
 Z(\vecC,\sigma)$.
There exists a partition of unity
$\chi_{\sigma}$  $(\sigma\in\gbigs_{\ell})$ of
$(X\setminus H)_S$
such that
(i) it subordinates to the open covering
 $\{Z(\vecC,\sigma)\}_{\sigma\in S}$,
(ii) $d\chi_{\sigma}$ is bounded with respect to $g_{X\setminus H}$,
(iii) each $\chi_{\sigma}$ is $(S^1)^{\ell}$-invariant.
\hfill\qed
\end{lem}

\subsubsection{KMS-structure of the $\lambda$-flat bundle
associated with a tame harmonic bundle}

Let $(E,\delbar_E,\theta,h)$ be a tame harmonic bundle
on $(X,H)$.
For $\lambda\in\cnum$,
let $(\nbigp_{\ast}\nbigelambda,\DDlambda)$
be the associated regular filtered $\lambda$-flat bundle
on $(X,H)$.
(See \S\ref{subsection;22.4.27.3}.)
Let $\veczero=(0,\ldots,0)\in\real^{\Lambda}$.
For $i\in\Lambda$,
we have the endomorphism
$\Res_i(\DDlambda)$
of $\nbigp_{\veczero}\nbigelambda_{|H_i}$
obtained as the residue of $\DDlambda$.
The eigenvalues are constant on $H_i$.
In particular,
there exists a decomposition
\begin{equation}
\label{eq;22.4.27.21}
\nbigp_{\veczero}\nbigelambda_{|H_i}
=\bigoplus_{\alpha\in\cnum}
\lefttop{i}\EE_{\alpha}
\nbigp_{\veczero}\nbigelambda_{|H_i}
\end{equation}
which is preserved by $\Res_i(\DDlambda)$,
and
$\Res_i(\DDlambda)-\alpha$
is nilpotent on
$\lefttop{i}\EE_{\alpha}$.       
For $-1<b\leq 0$,
let $\lefttop{i}F_b(\nbigp_{\veczero}\nbigelambda)$
denote the image of
$\nbigp_{b\veciti_i}\nbigelambda_{|H_i}
\lrarr
\nbigp_{\veczero}\nbigelambda_{|H_i}$,
which are naturally locally free $\nbigo_{H_i}$-modules.
They are preserved by
$\Res_i(\DDlambda)$.
We note that
$\lefttop{i}F_b=
\bigoplus_{\alpha}
\lefttop{i}F_b\cap
\lefttop{i}\EE_{\alpha}$.

For $I\subset\ellsitabar$,
we obtain the filtrations
$\lefttop{i}F$ $(i\in I)$
and
the decompositions
$\lefttop{i}\EE$ $(i\in I)$
of
$\nbigp_{\veczero}\nbigelambda_{|H_I}$.
There exists the decomposition
\[
 \nbigp_{\veczero}\nbigelambda_{|H_I}
 =\bigoplus_{\vecbeta\in\cnum^I}
 \lefttop{I}\EE_{\vecbeta} \nbigp_{\veczero}\nbigelambda_{|H_I}
\]
such that
$\bigl(
\lefttop{i}\EE_{\beta}
\bigr)_{|H_I}
=\bigoplus_{\beta_i=\beta}
\lefttop{I}\EE_{\vecbeta}$.
By the construction, we have
\[
 \bigl(
 \lefttop{i}F_b
 \bigr)_{|H_I}
 =\bigoplus_{\vecbeta}
 \Bigl(
 \lefttop{I}\EE_{\vecbeta}
 \cap
 \bigl(
 \lefttop{i}F_b
 \bigr)_{|H_I}
 \Bigr).
\]

Let $I\subset J$.
For $\vecc\in\openclosed{-1}{0}^I$,
we set
\[
\lefttop{I}F_{\vecc}(\nbigp_{\veczero}\nbigelambda_{|H_J})
=\bigcap_{i\in I}
 \lefttop{i}F_{c_i}(\nbigp_{\veczero}\nbigelambda_{|H_J}).
\]
We use the partial order $\leq$ on $\real^{I}$
defined by
$(c_i')\leq (c_i)$
$\stackrel{\rm}{\Longleftrightarrow}$
$c_i'\leq c_i$ for any $i$.
We say $\vecc'\lneq\vecc$ if $\vecc'\leq\vecc$
and $\vecc'\neq\vecc$.
Because $\nbigp_{\ast}\nbigelambda$ is a filtered bundle,
\[
 \lefttop{I}\Gr^F_{\vecc}
 (\nbigp_{\veczero}\nbigelambda_{|H_I}):=
 \frac{\lefttop{I}F_{\vecc}(\nbigp_{\veczero}\nbigelambda_{|H_I})}
 {\sum_{\vecc'\lneq\vecc}
 \lefttop{I}F_{\vecc'}(\nbigp_{\veczero}\nbigelambda_{|H_I})}
\]
is a locally free $\nbigo_{H_I}$-module.
It is equipped with
the induced endomorphisms
$\Res_i(\DDlambda)$ $(i\in I)$.
There exists the decomposition
\begin{equation}
\label{eq;22.2.24.20}
 \lefttop{I}\Gr^F_{\vecc}
  (\nbigp_{\veczero}\nbigelambda_{|H_I})=
  \bigoplus_{\vecalpha\in\cnum^I}
  \lefttop{I}\EE_{\vecalpha}
   \lefttop{I}\Gr^F_{\vecc}
  (\nbigp_{\veczero}\nbigelambda_{|H_I}),
\end{equation}
where $\Res_i(\DDlambda)-\alpha_i$
are nilpotent on 
$\lefttop{I}\EE_{\vecalpha}
   \lefttop{I}\Gr^F_{\vecc}
  (\nbigp_{\veczero}\nbigelambda_{|H_I})$.

\begin{prop}[\mbox{\cite[Corollary 4.42]{mochi2}}]
There exists a decomposition
of a locally free sheaf:
\begin{equation}
\label{eq;22.2.21.1}
 \nbigp_{\veczero}\nbigelambda
 =\bigoplus_{\veca\in \openclosed{-1}{0}^{\ell}}
 \bigoplus_{\vecalpha\in\cnum^{\ell}}
 \nbigv_{\veca,\vecalpha}
\end{equation}
which is compatible with 
the parabolic structure
and the generalized eigen decomposition
 i.e., the following holds
for any $I\subset\Lambda$
and
for any $\vecb\in\openclosed{-1}{0}^I$
and $\vecbeta\in\cnum^{I}$:
\[
\lefttop{I}\EE_{\vecbeta}
\lefttop{I}F_{\vecb}
(\nbigp_{\veczero}\nbigelambda)
=\bigoplus_{q_I(\veca)\leq q_I(\vecb)}
\bigoplus_{q_I(\vecalpha)=\vecbeta}
\nbigv_{\veca,\vecalpha|H_I}.
\]      
Such a decomposition is called a splitting
of the KMS-structure of $(\nbigp_{\ast}\nbigelambda,\DDlambda)$.
Note that a splitting is not unique. 

\hfill\qed
\end{prop}
 
We say that 
a frame $\vecv=(v_1,\ldots,v_{\rank E})$
of $\nbigp_{\veczero}\nbigelambda$
is compatible with the KMS-structure
of $(\nbigp_{\ast}\nbigelambda,\DDlambda)$
if there exists a splitting
$\nbigp_{\veczero}\nbigelambda
=\bigoplus\nbigv_{\veca,\vecalpha}$ of the KMS-structure
such that $v_i$ is a section of
$\nbigv_{\veca(v_i),\vecalpha(v_i)}$
for some $(\veca(v_i),\vecalpha(v_i))$.
  
\subsubsection{The weight filtrations}

Let $I\subset\ellsitabar$.
Let $\nbign_i$ denote the nilpotent part
of $-\Res_i(\DDlambda)$.
For each $J\subset I$,
let $W(J)$ denote the filtration of
$\lefttop{I}\Gr^F_{\vecc}
\nbigp_{\veczero}\nbigelambda_{|H_I}$
obtained as the weight filtration of
$\sum_{i\in J}\nbign_i$.
We recall the following lemmas
contained in \cite[Theorem 12.48]{mochi2}.
\begin{lem}
$W(J)$ is a filtration by subbundles of
$\lefttop{I}\Gr^F_{\vecc}
\nbigp_{\veczero}\nbigelambda_{|H_I}$ on $H_I$.
\hfill\qed
\end{lem}

\begin{lem}
Let $I_1\subset I_2$.
Let $\vecJ=(J_1,J_2,\ldots,J_m)\in\nbigs(I_1)$.
\begin{itemize}
 \item The decomposition {\rm(\ref{eq;22.2.24.20})},
       the filtrations
       $W(J_p)$ $(p=1,\ldots,m)$
       and $\lefttop{i}F$ $(i\in I_2\setminus I_1)$
	of $\lefttop{I_1}
	\Gr^F_{\vecc}(\nbigp_{\veczero}\nbigelambda_{|H_{I_1}})
	_{|H_{I_2}}$
	are compatible,
	i.e.,
	there exists a decomposition
\begin{equation}
\label{eq;22.2.21.2}
 \lefttop{I_1}\Gr^F_{\vecc}\bigl(
 \nbigp_{\veczero}\nbigelambda_{|H_{I_1}}
 \bigr)_{|H_{I_2}}
	=
	\bigoplus_{\vecbeta\in\cnum^{I_1}}
	\bigoplus_{\vecb\in \openclosed{-1}{0}^{I_2\setminus I_1}}
	\bigoplus_{\veck\in\seisuu^m}
	\nbigg_{\vecbeta,\vecb,\veck}
\end{equation}
	such that
\[
	\lefttop{I_2}\EE_{\vecbeta|H_{I_2}}
	\cap
       \bigcap_{i\in I_2\setminus I_1}
       \lefttop{i}F_{b_i|H_{I_2}}
	\cap
	\bigcap_{p=1}^m
	W_{k_p}(J_p)_{|H_{I_2}}
	=\bigoplus_{\vecb'\leq\vecb}
	\bigoplus_{\veck'\leq\veck}
	\nbigg_{\vecbeta,\vecb',\veck'}.
\]
 \item
      The induced filtrations
      $W(J_p)$ of
      $\lefttop{I_2}\Gr^F_{\vecb}(\nbigp_{\veczero}\nbigelambda_{|H_{I_2}})$
      are equal to the weight filtrations of
      $\sum_{i\in J_p}\lefttop{I_2}\Gr^F_{\vecb}(\nbign_i)$,
      where $\Gr^F_{\vecb}(\nbign_i)$
      denote the endomorphisms of
      $\lefttop{I_2}\Gr^F_{\vecb}(\nbigp_{\veczero}\nbigelambda_{|H_{I_2}})$
      induced by
      $\nbign_i$.
      \hfill\qed
\end{itemize}
 \end{lem}

\subsubsection{Some particular frames and the norm estimate}
\label{subsection;22.1.26.2}

Let $\sigma\in\gbigs_{\ell}$.
Let
$\vecw^{\sigma}=(w^{\sigma}_1,\ldots,w^{\sigma}_{\rank E})$
be a frame of $\nbigp_{\veczero}\nbigelambda$
such that the following holds.
\begin{itemize}
 \item It is compatible with the KMS-structure of
       $(\nbigp_{\ast}\nbigelambda,\DDlambda)$.
 \item For $j=1,\ldots,\ell$,
       the induced frame of
       $\lefttop{\sigma(\jbar)}
       \Gr^{F}_{\bullet}(\nbigp_{\veczero}\nbigelambda_{|H_{\sigma(\jbar)}})$
       is compatible with
       the weight filtration
       $W(\sigma(\jbar))$.
\end{itemize} 
Note that such a frame exists 
due to \cite[Corollary 4.47, Theorem 12.48]{mochi2}.
For each $w^{\sigma}_i$,
we have the associated data
$\veca(w^{\sigma}_i)\in \openclosed{-1}{0}^{\ell}$
and $\vecalpha(w^{\sigma}_i)\in\cnum^{\ell}$.
We also obtain
$\veck(w^{\sigma}_i)\in\seisuu^{\ell}$
by the condition
\[
 k_j(w^{\sigma}_i)
 =\max\bigl\{
 k\,\big|\,
 [w^{\sigma}_i]_{\ellsitabar}
 \in
 W_k(\sigma(\jbar))
 \bigr\}.
\]
Here, $[w^{\sigma}_i]_{\ellsitabar}$
denotes the induced section of
$\lefttop{\sigma(\ellsitabar)}
\Gr^F_{\veca(w_i)}\bigl(
 \nbigp_{\veczero}\nbigelambda_{|H_{\sigma(\ellsitabar)}}
 \bigr)$.
We note
\[
 k_j(w^{\sigma}_i)
 =\max\bigl\{
 k\,\big|\,
 [w^{\sigma}_i]_{\mbar}
 \in
 W_k(\sigma(\jbar))
 \bigr\}
\]
for any $j\leq m\leq \ell$,
where 
$[w^{\sigma}_i]_{\mbar}$
denotes the induced section of
$\lefttop{\sigma(\mbar)}
\Gr^F_{q_{\sigma(\mbar)}(\veca(w_i))}\bigl(
 \nbigp_{\veczero}\nbigelambda_{|H_{\sigma(\mbar)}}
 \bigr)$.

We recall the norm estimate
\cite[Theorem 13.29]{mochi2}.
The one dimensional case is proved in \cite{Simpson-non-compact}.
\begin{prop}
\label{prop;22.2.23.50}
Let $\vecw^{\sigma\prime}=(w^{\sigma\prime}_i)$
be the $C^{\infty}$-frame of $\nbigelambda$
determined as follows:
\[
 w^{\sigma\prime}_i=
 w^{\sigma}_i\cdot
 \prod_{j=1}^{\ell}
 |z_j|^{a_j(w^{\sigma}_i)}
 \prod_{p=1}^{\ell}
 (-\log|z_{\sigma(p)}|^2)^{-(k_p(w^{\sigma}_i)-k_{p-1}(w^{\sigma}_i))/2}.
\] 
Let $H(\vecw^{\sigma\prime})$
denote the Hermitian matrix valued function
determined by 
$H(\vecw^{\sigma\prime})_{i,i'}=h(w^{\sigma\prime}_i,w^{\sigma\prime}_{i'})$.
Then, $H(\vecw^{\sigma\prime})$
and $H(\vecw^{\sigma\prime})^{-1}$ are bounded
on $Z(\vecC,\sigma)$.
\hfill\qed
\end{prop}

\subsubsection{The estimate of the connection form}

Let $A$ be the diagonal matrix valued $1$-form
whose $(p,p)$-th entries are
$\sum_{j=1}^{\ell}\alpha_j(w_p^{\sigma})dz_j/z_j$.
Let $\DD^{\lambda(0)}$
denote the logarithmic $\lambda$-connection of
$\nbigp_{\veczero}\nbigelambda$
determined by
$\DD^{\lambda(0)}\vecw^{\sigma}=\vecw^{\sigma} A$.
We obtain the section
$f=\DD^{\lambda}-\DD^{\lambda(0)}$ of
$\End(\nbigp_{\veczero}\nbigelambda)\otimes\Omega^1(\log H)$.
We shall prove the following proposition
in \S\ref{subsection;22.2.22.1}
and \S\ref{subsection;22.2.22.2}.

\begin{prop}
\label{prop;22.2.21.3}
$f_{|Z(\vecC,\sigma)}$ is bounded 
with respect to $h$ and the Poincar\'{e} metric.
\end{prop}

For the proof of Proposition \ref{prop;22.2.21.3},
it is enough to study the case $\sigma=\id$.
We set $\vecw=\vecw^{\id}$
and $Z(\vecC)=Z(\vecC,\id)$.

\subsubsection{Preliminary}

Let $f_j$ be the endomorphisms determined by
$f=\sum_{j=1}^{\ell}f_j\,dz_j/z_j+\sum_{j=\ell+1}^n f_j\,dz_j$.
Let $B^j$ be the matrix valued functions
determined by $f_j\vecw=\vecw\cdot B^j$,
i.e.,
$f_jw_q=\sum B^j_{p,q}w_p$.
Let $C^j$ be determined by
$f_j(\vecw^{\prime})=\vecw^{\prime}C^j$.
Then, we have
\[
 C^j_{p,q}
 =B^j_{p,q}
 \prod_{m=1}^{\ell}
 |z_m|^{a_m(w_q)-a_m(w_p)}
 \prod_{m=1}^{\ell}
 (-\log |z_m|^2)^{-(k_m(w_q)-k_{m-1}(w_q)-k_m(w_p)+k_{m-1}(w_p))/2}.
\]
We note
\begin{multline}
  \prod_{m=1}^{\ell}
  (-\log |z_m|^2)^{-(k_m(w_q)-k_{m-1}(w_q)-k_m(w_p)+k_{m-1}(w_p))/2}
 =\\
 \prod_{m=1}^{\ell-1}
  \left(
  \frac{-\log|z_m|^2}{-\log|z_{m+1}|^2}
  \right)^{-(k_m(w_q)-k_m(w_p))/2}
  \times
  (-\log|z_{\ell}|^2)^{-(k_{\ell}(w_q)-k_{\ell}(w_p))/2}.
\end{multline}

For $w_p$ and $w_q$,
let $T_0(w_p,w_q)$ be the set of $1\leq i\leq \ell$
such that $\alpha_i(w_p)=\alpha_i(w_q)$
and $a_i(w_p)\leq a_i(w_q)$ hold,
and
let $T_1(w_p,w_q)$ be the set of $1\leq i\leq \ell$
such that
$(a_i(w_p),\alpha_i(w_p))
=(a_i(w_q),\alpha_i(w_q))$ holds.
Clearly, $T_1(w_p,w_q)\subset T_0(w_p,w_q)$.
There exists $\epsilon>0$ such that the following holds:
\[
 \prod_{m\in \ellsitabar\setminus T_0(w_p,w_q)}
 |z_m|^{1+a_m(w_q)-a_m(w_p)}
\times
  \prod_{m\in T_0(w_p,w_q)\setminus T_1(w_p,w_q)}
 |z_m|^{a_m(w_q)-a_m(w_p)}
 \leq
 \prod_{m\in \ellsitabar\setminus T_1(w_p,w_q)}
 |z_m|^{\epsilon}.
\]

If $T_1(w_p,w_q)=\ellsitabar$,
we set 
$i_0(w_p,w_q):=\ell$.
Otherwise,
we set
$i_0(w_p,w_q):=
-1+\min(\ellsitabar\setminus T_1(w_p,w_q))$.
If $T_1(w_p,w_q)\neq\emptyset$,
let $i_1(w_p,w_q)$ be determined by
the following condition.
\begin{itemize}
 \item If $k_i(w_p)\leq k_i(w_q)$ holds
       for any $i\leq i_0(w_p,w_q)$,
       we set
       $i_1(w_p,w_q):=i_0(w_p,w_q)+1$.
 \item Otherwise,
       $k_i(w_p)\leq k_i(w_q)$
       for any $i< i_1(w_p,w_q)$
       and
       $k_{i_1(w_p,w_q)}(w_p)>
        k_{i_1(w_p,w_q)}(w_q)$.
\end{itemize}

\subsubsection{The case $\ell<j\leq n$}
\label{subsection;22.2.22.1}

Let us study
$B^j$ and $C^j$ in the case $\ell<j\leq n$.

\begin{lem}
The following holds.
\begin{itemize}
 \item $B^j_{p,q|H_i}=0$
       if $i\in\ellsitabar\setminus T_0(w_p,w_q)$.
 \item Let $m\leq i_0(w_p,w_q)$.
       Then,
      $B^j_{p,q|H_{\mbar}}=0$
       unless $k_i(w_p)\leq k_i(w_q)$
       for any $i\leq m$.
\end{itemize}
\end{lem}
\pf
Let $\DDlambda_j$ be the differential operator
induced by $\DDlambda$ and $\del_{z_j}$.
It induces a differential operator $\DDlambda_{j|H_I}$
of $\nbigp_{\veczero}\nbigelambda_{|H_I}$
for $I\subset\ellsitabar$.
It is standard that
$\DDlambda_{j|H_i}(\Res_i(\DDlambda))=0$
for any $1\leq i\leq \ell$.
Moreover, because $\DDlambda_j$ preserves
$\nbigp_{\vecb}\nbigelambda$ for any $\vecb\in\real^{\ell}$,
$\DDlambda_{j|H_i}$ preserves
the parabolic filtrations $\lefttop{i}F$.
Because
$\DDlambda_{j|H_{\ibar}}(\Res_{p}(\DDlambda))=0$
for any $p\leq i$,
$\DDlambda_{j|H_{\ibar}}$
preserves the weight filtration
$W(\ibar)$.
Then, the claim of the lemma follows.
\hfill\qed

\vspace{.1in}

If $T_1(w_p,w_q)=\emptyset$,
there exists $\epsilon>0$
such that the following holds:
\[
 |C^j_{p,q}|
 =O\Bigl(
 \prod_{m\in\ellsitabar}
 |z_m|^{\epsilon}
 \Bigr).
\]
Let us consider the case where $T_1(w_p,w_q)\neq\emptyset$.
If $i_1(w_p,w_q)=i_0(w_p,w_q)+1$,
there exists $\epsilon>0$
such that the following holds on $Z(\vecC)$:
\[
 |C^j_{p,q}|
 =O\Bigl(
  \prod_{m\in\ellsitabar\setminus T_1(w_p,w_q)}
 |z_m|^{\epsilon}
 \Bigr).
\]
If $i_1(w_p,w_q)\leq i_0(w_p,w_q)$,
we obtain
$B^j_{p,q|H_{\underline{i_1(w_p,w_q)}}}=0$,
and hence
\[
B^j_{p,q}=
O\left(
\Bigl(\sum_{i\leq i_1(w_p,w_q)}|z_i|
\Bigr)
\times
 \prod_{m\in \ellsitabar\setminus T_0(w_p,w_q)}
 |z_m|
\right). 
\]
On $Z(\id,\vecC)$,
we obtain
\[
 |B^j_{p,q}|=
 O\Bigl(|z_{i_1(w_p,w_q)}|^{\epsilon}
 \prod_{m\in\ellsitabar\setminus T_0(w_p,w_q)}
 |z_m|
 \Bigr)
\]
for some $\epsilon>0$.
Hence, there exists $\epsilon'>0$
such that
\[
 |C^j_{p,q}|=
 O\Bigl(|z_{i_1(w_p,w_q)}|^{\epsilon'}
 \prod_{m\in\ellsitabar\setminus T_0(w_p,w_q)}
 |z_m|^{\epsilon'}
 \Bigr)
\]
Thus, we obtain that
$f_j$ is bounded on $Z(\id,\vecC)$
in the case $\ell+1\leq j\leq n$.

\vspace{.1in}
From the proof,
we shall also obtain the following:
\begin{lem}
Let $1\leq i\leq \ell$.
Assume that either one of the following holds.
\begin{itemize}
 \item $i_0(w_p,w_q)<i$.
 \item $i\leq i_0(w_p,w_q)$,
       and there exists $m\leq i$
       such that $k_m(w_p)>k_m(w_q)$.
\end{itemize}
Then, there exists $\epsilon>0$ such that
$|C^j_{p,q}|=O(|z_i|^{\epsilon})$
on $Z(\vecC)$.
\hfill\qed 
\end{lem}

\subsubsection{The case $1\leq j\leq \ell$}
\label{subsection;22.2.22.2}

Let us consider the case $1\leq j\leq \ell$.

\begin{lem}
The following holds.
\begin{itemize}
 \item $B^j_{p,q|H_i}=0$
       if $i\in \ellsitabar\setminus T_0(w_p,w_q)$.
 \item Suppose $p\neq q$.
       Let $m\leq i_0(w_p,w_q)$.
       Then,
       $B^j_{p,q|H_{\mbar}}=0$
       unless
       $k_i(w_p)\leq k_i(w_q)$
       for any $1\leq i\leq \min\{m,j-1\}$
       and
       $k_i(w_p)\leq k_i(w_q)-2$
       for any
       $\min\{m,j-1\}+1\leq i\leq m$.
 \item $B^j_{p,p|H_j}=0$.
\end{itemize} 
\end{lem}
\pf
Let $z_j\DDlambda_{j}$
be the differential operator
on $\nbigp_{\veczero}\nbigelambda$
induced by
$\DDlambda$ and $z_j\del_{z_j}$.
It induces differential operators
$z_j\DDlambda_{j|H_I}$
of $\nbigp_{\veczero}\nbigelambda_{|H_I}$
for $I\subset\ellsitabar\setminus\{j\}$,
which preserve the parabolic filtrations
$\lefttop{i}F$ $(i\in I)$.
We also have
$z_j\DDlambda_j(\Res_i(\DDlambda))=0$ for $i\neq j$.
The residue $\Res_j(\DDlambda)$
on $\nbigp_{\veczero}\nbigelambda_{|H_j}$
is also induced by $z_j\DDlambda_j$,
which preserves the parabolic filtration $\lefttop{j}F$.
For $i<j$,
$z_j\DDlambda_{j|H_{\ibar}}$ preserves
the weight filtrations $W(\ibar)$.
For $i\geq j$,
the nilpotent part $\nbign_j$ of $\Res_j(\DDlambda)$
satisfies
$\Res_j(\DDlambda)W_k(\ibar)
\subset
 W_{k-2}(\ibar)$.
Note that
$z_j\DDlambda_j\vecw=\vecw(A^j+B^j)$,
where $A^j$ denotes the diagonal matrix
whose $(p,p)$-th entries are $\alpha_j(w_p)$.
Then, we obtain the claim of the lemma.
 \hfill\qed

\vspace{.1in}

Because
$C^j_{p,p}=B^j_{p,p}$,
we have $|C^j_{p,p}|=O(|z_j|)$.
If $T_1(w_p,w_q)=\emptyset$,
there exists $\epsilon>0$
such that the following holds:
\[
 |C^j_{p,q}|
 =O\Bigl(
 \prod_{m\in\ellsitabar}
 |z_m|^{\epsilon}
 \Bigr)
=O(|z_j|^{\epsilon}).
\]
Suppose $T_1(w_p,w_q)\neq\emptyset$.
Let us consider the case where $i_0(w_p,w_q)<j$.
If $i_1(w_p,w_q)\leq i_0(w_p,w_q)$,
we obtain
\[
 |B^j_{p,q}|=
 O\Bigl(|z_{i_1(w_p,w_q)}|^{\epsilon_1}
 \prod_{m\in \ellsitabar\setminus T_0(w_p,w_q)}|z_m|
 \Bigr)
= O\Bigl(|z_{j}|^{\epsilon_2}
 \prod_{m\in \ellsitabar\setminus T_0(w_p,w_q)}|z_m|
 \Bigr)
\]
for some $\epsilon_1,\epsilon_2>0$ on $Z(\id,\vecC)$
by using the argument in \S\ref{subsection;22.2.22.1}.
Hence, there exists $\epsilon_3>0$
such that
\[
 |C^j_{p,q}|
=O\Bigl(|z_{j}|^{\epsilon_3}
 \prod_{m\in \ellsitabar\setminus T_1(w_p,w_q)}|z_m|^{\epsilon_3}
 \Bigr).
\]
If $i_1(w_p,w_q)=i_0(w_p,w_q)+1$,
we obtain
\[
|B^j_{p,q}|=
O\Bigl(
  \prod_{m\in \ellsitabar\setminus T_0(w_p,w_q)}|z_m|\Bigr).
\]
Hence, we obtain
\[
 |C^j_{p,q}|
 =O\Bigl(
 \prod_{m\in\ellsitabar\setminus T_1(w_p,w_q)}
 |z_m|^{\epsilon_0}
 \Bigr)
 =O\Bigl(|z_{i_0(w_p,w_q)+1}|^{\epsilon_1}\Bigr)
 =O\Bigl(|z_j|^{\epsilon_2}\Bigr)
\]
for some $\epsilon_0,\epsilon_1,\epsilon_2>0$
on $Z(\id,\vecC)$.
Let us consider the case where $j\leq i_0(w_p,w_q)$.
If $i_1(w_p,w_q)\leq j$,
then
we obtain
\[
 |B^j_{p,q}|=
 O\Bigl(|z_{i_1(w_p,w_q)}|^{\epsilon_1}
 \prod_{m\in \ellsitabar\setminus T_0(w_p,w_q)}|z_m|
 \Bigr)
= O\Bigl(|z_{j}|^{\epsilon_2}
 \prod_{m\in \ellsitabar\setminus T_0(w_p,w_q)}|z_m|
 \Bigr)
\]
for some $\epsilon_1,\epsilon_2>0$
on $Z(\id,\vecC)$
by using the argument in \S\ref{subsection;22.2.22.1}.
Hence, there exist $\epsilon_3,\epsilon_4>0$ such that
\[
 |C^j_{p,q}|
 =O\Bigl(|z_{i_1(w_p,w_q)}|^{\epsilon_3}
 \prod_{m\in \ellsitabar\setminus T_1(w_p,w_q)}|z_m|^{\epsilon_3}
 \Bigr)
= O\Bigl(|z_{j}|^{\epsilon_4}
 \prod_{m\in \ellsitabar\setminus T_1(w_p,w_q)}|z_m|^{\epsilon_4}
 \Bigr).
\]
If $j<i_1(w_p,w_q)$,
either one of the following holds:
(i) $k_i(w_p)\leq k_i(w_q)-2$ for any
$j\leq i\leq i_0(w_p,w_q)$,
(ii) there exists $j\leq i_2(w_p,w_q)<i_0(w_p,w_q)$
such that
$k_i(w_p)\leq k_i(w_q)-2$ for any
$j\leq i\leq i_2(w_p,w_q)$
and that
$k_{i_2(w_p,w_q)+1}(w_p)>
k_{i_2(w_p,w_q)+1}(w_q)-2$.
In the case (i),
we obtain
the following estimate on $Z(\id,\vecC)$:
\[
 |B^j_{p,q}|=O\Bigl(
 \prod_{m\in \ellsitabar\setminus T_0(w_p,w_q)}
 |z_m|
 \Bigr).
\]
Hence, there exists $\epsilon,\epsilon'>0$ such that
\begin{multline}
 |C^j_{p,q}|=O\Bigl(
 \prod_{m=j}^{i_0(w_p,w_q)}
 \frac{-\log|z_{m+1}|^2}{-\log|z_m|^2}
 \prod_{m\in \ellsitabar\setminus T_1(w_p,w_q)}
 |z_m|^{\epsilon}
 \Bigr)
 =O\Bigl(
 (-\log|z_m|^2)^{-1}
 \prod_{m\in \ellsitabar\setminus T_1(w_p,w_q)}
 |z_m|^{\epsilon'}
 \Bigr)
 \\
 =O\Bigl(
 (-\log|z_j|^2)^{-1}
 \Bigr).
\end{multline}
In the case (ii),
we obtain the following estimate on $Z(\id,\vecC)$:
\[
 |B^j_{p,q}|=O\Bigl(
 |z_{i_2(w_p,w_q)+1}|
 \prod_{m\in \ellsitabar\setminus T_0(w_p,w_q)}
 |z_m|
 \Bigr).
\]
Hence, there exists $\epsilon,\epsilon'>0$ such that
\begin{multline}
 |C^j_{p,q}|=O\Bigl(
 \prod_{m=j}^{i_2(w_p,w_q)}
 \Bigl(
 \frac{-\log|z_{m+1}|^2}{-\log|z_m|^2}
 \Bigr)
 |z_{i_2(w_p,w_q)+1}|^{1/2}
 \prod_{m\in \ellsitabar\setminus T_1(w_p,w_q)}
 |z_m|^{\epsilon}
 \Bigr)
 \\
 =O\Bigl(
 (-\log|z_j|^2)^{-1}
 |z_{i_2(w_p,w_q)+1}|^{\epsilon'}
  \prod_{m\in \ellsitabar\setminus T_1(w_p,w_q)}
 |z_m|^{\epsilon}
 \Bigr)
 =O\Bigl(
 (-\log|z_j|^2)^{-1}
 \Bigr).
\end{multline}
Thus, we obtain Proposition \ref{prop;22.2.21.3}.
\hfill\qed

\vspace{.1in}
From the proof,
we also obtain the following lemma.
\begin{lem}
Let $1\leq i<j$.
Assume one of the following holds.
\begin{itemize}
 \item $i_0(w_p,w_q)<i$.
 \item $i_0(w_p,w_q)\geq i$,
       and there exists
       $m\leq i$ such that
       $k_m(w_p)>k_m(w_q)$.
\end{itemize}
Then, there exists $\epsilon>0$ such that
$|C^j_{p,q}|=O(|z_i|^{\epsilon})$.
\hfill\qed
\end{lem}

\begin{lem}
Let $j\leq i\leq \ell$.
Assume one of the following holds. 
\begin{itemize}
 \item $i_0(w_p,w_q)<i$.
 \item $i_0(w_p,w_q)\geq i$,
       and there exists
       $m<j$ such that
       $k_m(w_p)>k_m(w_q)$.
 \item $i_0(w_p,w_q)\geq i$,
       and there exists
       $j\leq m\leq i$ such that
       $k_m(w_p)>k_m(w_q)-2$.
\end{itemize}
Then, there exists $\epsilon>0$ such that
$|C^j_{p,q}|=O(|z_i|^{\epsilon})$.
\hfill\qed
\end{lem}

Let us give some complements.
For $1\leq s\leq \ell$,
we set
 \[
 (C^j_{p,q})'_s:=
 z_s\del_{z_s}B^j_{p,q}
 \prod_{m=1}^{\ell}
 |z_m|^{a_m(w_q)-a_m(w_p)}
 \prod_{m=1}^{\ell}
 (-\log |z_m|^2)^{-(k_m(w_q)-k_{m-1}(w_q)-k_m(w_p)+k_{m-1}(w_p))/2}.
 \]
For $\ell+1\leq s\leq n$,
we set
 \[
 (C^j_{p,q})'_s:=
 \del_{z_s}B^j_{p,q}
 \prod_{m=1}^{\ell}
 |z_m|^{a_m(w_q)-a_m(w_p)}
 \prod_{m=1}^{\ell}
 (-\log |z_m|^2)^{-(k_m(w_q)-k_{m-1}(w_q)-k_m(w_p)+k_{m-1}(w_p))/2}.
\]

\begin{lem}
\label{lem;22.2.24.11}
For $1\leq s\leq \ell$,
there exists $\epsilon>0$ such that
$|(C^j_{p,q})'_s|=O\bigl(|z_s|^{\epsilon}(-\log|z_j|^2)^{-1}\bigr)$.
For $\ell+1\leq s\leq n$, we have
$|(C^j_{p,q})'_s|=O\bigl((-\log|z_j|^2)^{-1}\bigr)$. 
\end{lem}
\pf
For $1\leq s\leq \ell$,
the following holds.
\begin{itemize}
 \item $z_s\del_{z_s}B^j_{p,p|H_j}=0$ and  $z_s\del_{z_s}B^j_{p,q|H_s}=0$.
 \item
      $z_s\del_{z_s}B^j_{p,q|H_i}=0$
      if $i\in\ellsitabar\setminus T_0(w_p,w_q)$.
 \item Suppose $p\neq q$.
       Let $m\leq i_0(w_p,w_q)$.
       Then,
       $z_s\del_{z_s}B^j_{p,q|H_{\mbar}}=0$
       unless
       $k_i(w_p)\leq k_i(w_q)$
       for any $1\leq i\leq \min\{m,j-1\}$
       and
       $k_i(w_p)\leq k_i(w_q)-2$
       for any
       $\min\{m,j-1\}+1\leq i\leq m$.
\end{itemize}
For $\ell+1\leq s\leq n$,
the following holds.
\begin{itemize}
 \item $\del_{z_s}B^j_{p,q|H_i}=0$
       if $i\in\ellsitabar\setminus T_0(w_p,w_q)$.
 \item Let $m\leq i_0(w_p,w_q)$.
       Then,
      $\del_{z_s}B^j_{p,q|H_{\mbar}}=0$
       unless $k_i(w_p)\leq k_i(w_q)$
       for any $i\leq m$.
\end{itemize}
Then, we obtain the claim of Lemma \ref{lem;22.2.24.11}
by a similar argument.
\hfill\qed

\subsubsection{Family version}

\paragraph{KMS-structure}

Let $\lambda_0\in\cnum$.
Let $U(\lambda_0)$ be a small neighbourhood of $\lambda_0$
in $\cnum$.
(See Remark \ref{rem;22.4.27.2}.)
We set $\nbigxzero:=U(\lambda_0)\times X$.
We use the notation
$\nbighzero$,
$\nbighzero_I$
and $\nbigzzero(\vecC,\sigma)$
similarly.
We set $\nbigx^{\lambda}=\{\lambda\}\times X\subset\nbigx$.
We use the notation
$\nbigh^{\lambda}$,
$\nbigh^{\lambda}_I$
and $\nbigz^{\lambda}(\vecC,\sigma)$
similarly.

Let $(\nbigpzero_{\ast}\nbige,\DD)$
denote the associated
regular filtered family of $\lambda$-flat bundles
on $(\nbigxzero,\nbighzero)$.
(See \S\ref{subsection;22.4.26.10} and \S\ref{subsection;22.3.25.120}.)
The family of $\lambda$-connections $\DD$
is logarithmic with respect to
$\nbigpzero_{\veca}\nbige$ for any $\veca\in\real^{\Lambda}$.
Let $\Res_i(\DD)$ denote the endomorphism of
$\nbigpzero_{\veczero}(\nbige)_{|\nbighzero_i}$
obtained as the residue of $\DD$.
There exists the decomposition
\begin{equation}
\label{eq;22.4.27.20}
\nbigpzero_{\veczero}(\nbige)_{|\nbighzero_i}
=\bigoplus_{\alpha\in\cnum}
 \lefttop{i}\EEzero_{\alpha}
 (\nbigpzero_{\veczero}\nbige_{|\nbighzero_i})
\end{equation}
uniquely determined by the conditions
(i) the decomposition is preserved by $\Res(\DD)$,
(ii) the restriction of (\ref{eq;22.4.27.20})
to $\nbigh^{\lambda_0}_i$
is equal to the decomposition
(\ref{eq;22.4.27.21}) for $\lambda=\lambda_0$.
For $-1<b\leq 0$,
let $\lefttop{i}\Fzero_b(\nbigpzero_{\veczero}\nbige)$
denote the image of
$\nbigpzero_{b\veciti_i}\nbige_{|\nbighzero_i}
\lrarr
\nbigpzero_{\veczero}\nbige_{|\nbighzero_i}$.
It is preserved by $\Res_i(\DD)$.
By the construction,
we have
$\lefttop{i}\Fzero_b=
\bigoplus\bigl(
\lefttop{i}\Fzero_b
\cap
\lefttop{i}\EEzero_{\alpha}
\bigr)$.

For $I\subset\ellsitabar$,
we obtain the filtrations
$\lefttop{i}\Fzero$ $(i\in I)$
and
the decompositions
$\lefttop{i}\EEzero$ $(i\in I)$
of
$\nbigpzero_{\veczero}\nbige_{|\nbighzero_I}$.
By setting
$\lefttop{I}\EEzero_{\vecbeta}
=\bigcap_{i\in I}\lefttop{i}\EEzero_{\beta_i}$,
we obtain the decomposition
\[
 \nbigpzero_{\veczero}\nbige_{|\nbighzero_I}
 =\bigoplus_{\vecbeta\in\cnum^I}
 \lefttop{I}\EEzero_{\vecbeta}
 \nbigpzero_{\veczero}\nbige_{|\nbighzero_I}.
\]
By the construction, the following holds:
\[
 \bigl(
 \lefttop{i}\Fzero_b
 \bigr)_{|\nbighzero_I}
 =\bigoplus_{\vecbeta}
 \Bigl(
 \lefttop{I}\EEzero_{\vecbeta}
 \cap
 \bigl(
 \lefttop{i}\Fzero_b
 \bigr)_{|\nbighzero_I}
 \Bigr).
\]

Let $I\subset J$.
For $\vecc\in\openclosed{-1}{0}^I$,
we set
\[
\lefttop{I}\Fzero_{\vecc}(\nbigpzero_{\veczero}\nbige_{|\nbighzero_J})
=\bigcap_{i\in I}
 \lefttop{i}\Fzero_{c_i}(\nbigpzero_{\veczero}\nbige_{|\nbighzero_J}).
\]
Because $\nbigpzero_{\ast}\nbige$ is a filtered bundle
on $(\nbigxzero,\nbighzero)$,
\[
 \lefttop{I}\Gr^{\Fzero}_{\vecc}
 (\nbigpzero_{\veczero}\nbige_{|\nbighzero_I}):=
 \frac{\lefttop{I}\Fzero_{\vecc}
  (\nbigpzero_{\veczero}\nbige_{|\nbighzero_I})}
 {\sum_{\vecc'\lneq\vecc}
 \lefttop{I}\Fzero_{\vecc'}(\nbigpzero_{\veczero}\nbige_{|\nbighzero_I})}
\]
is a locally free $\nbigo_{\nbighzero_I}$-module.
There exists the induced decomposition
\begin{equation}
\label{eq;22.2.24.21}
 \lefttop{I}\Gr^{\Fzero}_{\vecc}
  (\nbigpzero_{\veczero}\nbige_{|\nbighzero_I})=
  \bigoplus_{\vecalpha\in\cnum^I}
  \lefttop{I}\EEzero_{\vecalpha}
   \lefttop{I}\Gr^{\Fzero}_{\vecc}
  (\nbigpzero_{\veczero}\nbige_{|\nbighzero_I})
\end{equation}
which is preserved by $\Res_i(\DD)$ $(i\in I)$.
Moreover,
$\Res_i(\DD)-\eigenmap(\lambda,u_i)\id$
are nilpotent on 
$\lefttop{I}\EEzero_{\vecalpha}
 \lefttop{I}\Gr^{\Fzero}_{\vecc}
 (\nbigpzero_{\veczero}\nbige_{|\nbighzero_I})$,
where
$u_j\in \real\times\cnum$
are determined by
$\kmsmap(\lambda_0,u_j)=(c_j,\alpha_j)$.
(See \S\ref{subsection;22.3.25.120}
for the map
$\kmsmap(\lambda_0):\real\times\cnum\simeq\real\times\cnum$.)

\begin{lem}[\mbox{\cite[Corollary 4.42]{mochi2}}]
There exists a decomposition of locally free sheaves
\begin{equation}
\label{eq;22.3.15.1}
 \nbigpzero_{\veczero}\nbige
 =\bigoplus_{\veca\in \openclosed{-1}{0}^{\ell}}
  \bigoplus_{\vecalpha\in\cnum^{\ell}}
  \nbigvzero_{\veca,\vecalpha}
\end{equation}
such that the following holds
for any $I\subset\Lambda$
and
for any $\vecb\in\openclosed{-1}{0}^I$
and $\vecbeta\in\cnum^{I}$:
\[
\lefttop{I}\EEzero_{\vecbeta}
\lefttop{I}\Fzero_{\vecb}
(\nbigpzero_{\veczero}\nbige)
=\bigoplus_{q_I(\veca)\leq q_I(\vecb)}
\bigoplus_{q_I(\vecalpha)=\vecbeta}
\nbigv_{\veca,\vecalpha|\nbighzero_I}.
\]
Such a decomposition {\rm(\ref{eq;22.3.15.1})} is called a splitting of
the KMS-structure of $(\nbigpzero_{\ast}\nbige,\DD)$.
\hfill\qed
\end{lem}

 We say that a frame $\vecv=(v_i)$
of $\nbigpzero_{\veczero}\nbige$
is compatible with the KMS-structure of
$(\nbigpzero_{\ast}\nbige,\DD)$
if there exists a splitting (\ref{eq;22.3.15.1})
such that $v_i$
is a section of $\nbigvzero_{\veca(v_i),\vecalpha(v_i)}$.

\paragraph{The weight filtrations}

Let $\nbign_j$ denote
the endomorphisms of
$\lefttop{I}\Gr^{\Fzero}_{\vecc}
\bigl(\nbigpzero_{\veczero}\nbige_{|\nbighzero_I}\bigr)$
obtained as the nilpotent part of $-\Res_j(\DD)$.
For each $J\subset I$,
we obtain the weight filtration
$W(J)$ of $\sum_{j\in J} \nbign_j$.
The following lemmas are contained
in \cite[Theorem 12.48]{mochi2}.
\begin{lem}
$W(J)$ is a filtration by subbundles of
$\lefttop{I}\Gr^{\Fzero}_{\vecc}
(\nbigpzero_{\veczero}\nbige_{|\nbighzero_I})$. 
\hfill\qed
\end{lem}

\begin{lem}
Let $I_1\subset I_2$.
Let $\vecJ=(J_1,J_2,\ldots,J_m)\in\nbigs(I_1)$.
\begin{itemize}
 \item The decomposition {\rm(\ref{eq;22.2.24.21})},
       the filtrations
       $W(J_p)$ $(p=1,\ldots,m)$
       and $\lefttop{i}\Fzero$ $(i\in I_2\setminus I_1)$
       of
       \[
	\lefttop{I_1}
	\Gr^{\Fzero}_{\vecc}(\nbigpzero_{\veczero}\nbige_{|\nbighzero_{I_1}})
       _{|\nbighzero_{I_2}}
       \]
	are compatible.
 \item
      The induced filtrations
      $W(J_p)$ of
      $\lefttop{I_2}\Gr^{\Fzero}_{\vecb}
      (\nbigpzero_{\veczero}\nbige_{|\nbighzero_{I_2}})$
      equal the weight filtration of
      $\sum_{i\in J_p}\lefttop{I_2}\Gr^{\Fzero}_{\vecb}(\nbign_i)$,
      where $\lefttop{I_2}\Gr^F_{\vecb}(\nbign_i)$
      denote the endomorphisms of
      $\lefttop{I_2}\Gr^{\Fzero}_{\vecb}
      (\nbigpzero_{\veczero}\nbige_{|\nbighzero_{I_2}})$
      induced by
      $\nbign_i$.
      \hfill\qed
\end{itemize}
 \end{lem}

\paragraph{Some particular frames and the norm estimate}

Let $\sigma\in\gbigs_{\ell}$.
Let
$\vecw^{\sigma}=(w^{\sigma}_1,\ldots,w^{\sigma}_{\rank E})$
be a frame of $\nbigpzero_{\veczero}\nbige$
such that the following holds.
\begin{itemize}
 \item It is compatible with the KMS-structure of
       $(\nbigpzero_{\ast}\nbige,\DD)$.
 \item For $j=1,\ldots,\ell$,
       the induced frame of
       $\lefttop{\sigma(\jbar)}
       \Gr^{\Fzero}_{\bullet}
       (\nbigpzero_{\veczero}\nbige_{|\nbighzero_{\sigma(\jbar)}})$
       is compatible with
       the weight filtration
       $W(\sigma(\jbar))$.
\end{itemize}
Note that such a frame exists,
according to \cite[Corollary 4.47, Theorem 12.48]{mochi2}.
For each $w^{\sigma}_i$,
we have the associated data
$\veca(w^{\sigma}_i)\in \openclosed{-1}{0}^{\ell}$
and $\vecalpha(w^{\sigma}_i)\in\cnum^{\ell}$.
We also obtain
$\veck(w^{\sigma}_i)\in\seisuu^{\ell}$
by the condition
\[
 k_j(w^{\sigma}_i)
 =\max\bigl\{
 k\,\big|\,
 [w^{\sigma}_i]_{\ellsitabar}
 \in
 W_k(\sigma(\jbar))
 \bigr\}.
\]
Here, $[w^{\sigma}_i]_{\ellsitabar}$
denotes the induced section of
$\lefttop{\sigma(\ellsitabar)}
\Gr^{\Fzero}_{\veca(w_i)}\bigl(
 \nbigpzero_{\veczero}\nbige_{|\nbighzero_{\sigma(\ellsitabar)}}
 \bigr)$.
We note
\[
 k_j(w^{\sigma}_i)
 =\max\bigl\{
 k\,\big|\,
 [w^{\sigma}_i]_{\mbar}
 \in
 W_k(\sigma(\jbar))
 \bigr\}
\]
for any $j\leq m\leq \ell$,
where 
$[w^{\sigma}_i]_{\mbar}$
denotes the induced section of
$\lefttop{\sigma(\mbar)}
\Gr^{\Fzero}_{q_{\sigma(\mbar)}(\veca(w_i))}\bigl(
 \nbigpzero_{\veczero}\nbige_{|\nbighzero_{\sigma(\mbar)}}
 \bigr)$.

We recall the norm estimate.
For each $w_i$,
we have
$u_p(w_i)\in\real\times\cnum$ $(p=1,\ldots,\ell)$
determined by
$\kmsmap(\lambda_0,u_p(w_i))=(a_p(w_i),\alpha_p(w_i))$.
We obtain
$\vecu(w_i)=(u_p(w_i))\in(\real\times\cnum)^{\ell}$.

\begin{prop}
\label{prop;22.2.24.100}
Let $\vecw^{\sigma\prime}=(w^{\sigma\prime}_i)$
be the $C^{\infty}$-frame of $\nbige_{|\nbigxzero\setminus\nbighzero}$
determined as follows:
\[
 w^{\sigma\prime}_i=
 w^{\sigma}_i\cdot
 \prod_{j=1}^{\ell}
 |z_j|^{\paramap(\lambda,u_j(w^{\sigma}_i))}
 \prod_{p=1}^{\ell}
 (-\log|z_{\sigma(p)}|^2)^{-(k_p(w^{\sigma}_i)-k_{p-1}(w^{\sigma}_i))/2}.
\] 
Let $H(\vecw^{\sigma\prime})$
denote the Hermitian matrix valued function
determined by 
$H(\vecw^{\sigma\prime})_{i,i'}=h(w^{\sigma\prime}_i,w^{\sigma\prime}_{i'})$.
Then, $H(\vecw^{\sigma\prime})$
and $H(\vecw^{\sigma\prime})^{-1}$ are bounded
on $\nbigzzero(\vecC,\sigma)$.
\hfill\qed
\end{prop}

\paragraph{The estimate of the connection form}

Let $A$ be the diagonal matrix valued $1$-form
whose $(p,p)$-th entries are
\[
 \sum_{j=1}^{\ell}\eigenmap(\lambda,u_j(w_p^{\sigma}))dz_j/z_j.
\]
Let $\DD^{(0)}$
denote the logarithmic $\lambda$-connection of
$\nbigpzero_{\veczero}\nbige$
determined by
$\DD^{(0)}\vecw^{\sigma}=\vecw^{\sigma} A$.
We obtain the section
$f=\DD-\DD^{(0)}$ of
$\End(\nbigpzero_{\veczero}\nbige)\otimes\Omega^1(\log H)$.
The following proposition is similar to
Proposition \ref{prop;22.2.21.3}.

\begin{prop}
\label{prop;22.2.24.30}
 $f_{|\nbigzzero(\vecC,\sigma)}$
 is bounded with respect to $h$ and the Poincar\'{e} metric.
\hfill\qed
\end{prop}

Let $f_j$ be the endomorphisms determined by
$f=\sum_{j=1}^{\ell}f_j\,dz_j/z_j+\sum_{j=\ell+1}^n f_j\,dz_j$.
Let $B^j$ be the matrix valued functions
determined by $f_j\vecw=\vecw\cdot B^j$.
We obtain Proposition \ref{prop;22.2.24.30}
from the next two lemma
as in the case of Proposition \ref{prop;22.2.21.3}.

\begin{lem}
If $\ell<j\leq n$,
the following holds.
\begin{itemize}
 \item $B^j_{p,q|\nbighzero_i}=0$
       if $i\in\ellsitabar\setminus T_0(w_p,w_q)$.
 \item Let $m\leq i_0(w_p,w_q)$.
       Then,
      $B^j_{p,q|\nbighzero_{\mbar}}=0$
       unless $k_i(w_p)\leq k_i(w_q)$
       for any $i\leq m$.
       \hfill\qed
\end{itemize}
\end{lem}

\begin{lem}
If $1\leq j\leq \ell$,
the following holds.
\begin{itemize}
 \item $B^j_{p,q|\nbighzero_i}=0$
       if $i\in \ellsitabar\setminus T_0(w_p,w_q)$.
 \item Suppose $p\neq q$.
       Let $m\leq i_0(w_p,w_q)$.
       Then,
       $B^j_{p,q|\nbighzero_{\mbar}}=0$
       unless
       $k_i(w_p)\leq k_i(w_q)$
       for any $1\leq i\leq \min\{m,j-1\}$
       and
       $k_i(w_p)\leq k_i(w_q)-2$
       for any
       $\min\{m,j-1\}+1\leq i\leq m$.
 \item $B^j_{p,p|\nbighzero_j=0}$.
       \hfill\qed
\end{itemize} 
\end{lem}

Let $1\leq j\leq \ell$.
For $1\leq s\leq \ell$,
we set
 \[
 (C^j_{p,q})'_s:=
 z_s\del_{z_s}B^j_{p,q}
 \prod_{m=1}^{\ell}
 |z_m|^{\paramap(\lambda,u_m(w_q))-\paramap(\lambda,u_m(w_p))}
 \prod_{m=1}^{\ell}
 (-\log |z_m|^2)^{-(k_m(w_q)-k_{m-1}(w_q)-k_m(w_p)+k_{m-1}(w_p))/2}.
 \]
For $\ell+1\leq s\leq n$,
we set
 \[
 (C^j_{p,q})'_s:=
 \del_{z_s}B^j_{p,q}
 \prod_{m=1}^{\ell}
 |z_m|^{\paramap(\lambda,u_m(w_q))-\paramap(\lambda,u_m(w_p))}
 \prod_{m=1}^{\ell}
 (-\log |z_m|^2)^{-(k_m(w_q)-k_{m-1}(w_q)-k_m(w_p)+k_{m-1}(w_p))/2}.
\]

\begin{lem}
For $1\leq s\leq \ell$,
there exists $\epsilon>0$ such that
$|(C^j_{p,q})'_s|=O\bigl(|z_s|^{\epsilon}(-\log|z_j|^2)^{-1}\bigr)$.
For $\ell+1\leq s\leq n$, we have
$|(C^j_{p,q})'_s|=O\bigl((-\log|z_j|^2)^{-1}\bigr)$. 
\hfill\qed
\end{lem}

\subsection{Reduction to the non-negative part}

\subsubsection{Space of $L^2$-sections}
\label{subsection;22.1.27.1}

To simplify the description,
we set
$\Omega^k=\bigoplus_{\ell_1+\ell_2=k}
\Omega^{\ell_1,\ell_2}$.
For any $(S^1)^{\ell}$-invariant open subset $U$ of $X$,
let $L^2(\nbigelambda\otimes\Omega^{k},h)(U)$
denote the space of
$L^2$-sections of
$(\nbigelambda\otimes\Omega^k)
_{|U\setminus H}$
with respect to $h$ and $g_{X\setminus H}$.
Let $\vecv$ be a holomorphic frame of
$\nbigp_{\veczero}(\nbigelambda)$
on $X$.
For any
$s\in L^2(\nbigelambda\otimes\Omega^{k},h)(U)$,
we have the expression
\[
s=\sum_{\substack{J,K\subset\nbar \\ |J|+|K|=k}}
\sum_{i=1}^{\rank E}
s^{\vecv}_{i,J,K}\cdot \omega_{J,K}\cdot v_{i}
=
\sum_{\substack{J,K\subset\nbar \\ |J|+|K|=k}}
\sum_{i=1}^{\rank E}
\sum_{\vecp\in\seisuu^{\ell}}
\gbigf_{\ellsitabar,\vecp}(s^{\vecv}_{i,J,K})
e^{\sqrt{-1}\vecp\vectheta}
\cdot \omega_{J,K}\cdot v_{i}.
\]
Here,
$f=\sum_{\vecp}\gbigf_{\ellsitabar,\vecp}(f)e^{\sqrt{-1}\vecp\vectheta}$
is the Fourier expansion.

For $\vecm\in\seisuu^{\ell}$, we set
\[
 L^2(\nbigelambda\otimes\Omega^k,h,\vecm,\vecv)(U)
 =\Bigl\{
 s\in L^2(\nbigelambda\otimes\Omega^k,h,\vecv)(U)
 \,\Big|\,
 \gbigf_{\ellsitabar,\vecp}(s^{\vecv}_{i,J,K})=0\,\,(\vecp\neq\vecm)
 \Bigr\}.
\]
We obtain the orthogonal decomposition of the Hilbert space:
\[
 L^2(\nbigelambda\otimes\Omega^k,h)(U)
=\bigoplus_{\vecm\in\seisuu^{\ell}}
 L^2(\nbigelambda\otimes\Omega^k,h,\vecm,\vecv)(U).
\]

For $J\subset\ellsitabar$,
we set
\[
 \lefttop{J}FL^2(\nbigelambda\otimes\Omega^k,h)(U):=
 \bigoplus_{\substack{\vecm\in\seisuu^{\ell}\\
  m_j\geq 0\,(j\in J)
  }}
  L^2(\nbigelambda\otimes\Omega^k,h,\vecm,\vecv)(U).
\]
It is easy to see that
$\lefttop{J}FL^2(\nbigelambda\otimes\Omega^k,h)(U)$ is independent
of the choice of $\vecv$.
We have
$\lefttop{J}F\cap \lefttop{K}F=\lefttop{J\cup K}F$.

\begin{rem}
$\lefttop{J}FL^2(\nbigelambda\otimes\Omega^k,h)(U)$
can depend on the choice of a coordinate system
$(z_1,\ldots,z_n)$.
\hfill\qed
\end{rem}

\subsubsection{Reduction to the non-negative part}

Let $U$ be an $(S^1)^{\ell}$-invariant neighbourhood of
the origin $O=(0,\ldots,0)\in X$.
We set
\[
 \lefttop{J}F\nbigc^k_{L^2}(\nbigelambda,\DDlambda,h)(U):=
 \lefttop{J}FL^2(\nbigelambda\otimes\Omega^k,h)(U)
 \cap
 \nbigc^k_{L^2}(\nbigelambda,\DDlambda,h)(U).
\]
We obtain the subcomplex
$\lefttop{J}F\nbigc^{\bullet}_{L^2}(\nbigelambda,\DDlambda,h)(U)$.

Let $(\gbigu(O),\prec)$ denote the directed set
of $(S^1)^{\ell}$-invariant open neighbourhoods of $O$,
where $U_1\prec U_2$  is defined to be $U_1\supset U_2$.
For any $U\in\gbigu(O)$,
we set $\gbigu(O,U):=\{U'\in\gbigu(O)\,|\,U\prec U'\}$.
We set
\[
\lefttop{J}F\nbigc^{\bullet}_{L^2}(\nbigelambda,\DDlambda,h)_O:=
 \varinjlim_{U\in\gbigu(O)}
 \lefttop{J}F\nbigc^{\bullet}_{L^2}(\nbigelambda,\DDlambda,h)(U).
\]
 
\begin{prop}
\label{prop;22.1.28.10}
The natural morphism
$\lefttop{J}F\nbigc^{\bullet}_{L^2}(\nbigelambda,\DDlambda,h)_O
 \lrarr
 \nbigc^{\bullet}_{L^2}(\nbigelambda,\DDlambda,h)_O$
is a quasi-isomorphism.  
\end{prop}
\pf
Let $I\subset\ellsitabar$ be any subset.
Let $j\in\ellsitabar\setminus I$.
We set $\Itilde:=I\sqcup\{j\}$.
It is enough to prove that the natural morphism
\begin{equation}
\label{eq;22.1.25.22}
\lefttop{\Itilde}F\nbigc^{\bullet}_{L^2}(\nbigelambda,\DDlambda,h)_O
\lrarr
\lefttop{I}F
\nbigc^{\bullet}_{L^2}(\nbigelambda,\DDlambda,h)_O
\end{equation}
is a quasi-isomorphism.
Let $U\in\gbigu(O)$.
For a section
$s\in \lefttop{I}F\nbigc^{k}_{L^2}(\nbigelambda,\DDlambda,h)(U)$,
by using a frame $\vecv$ of $\nbigp_{\veczero}\nbigelambda$,
we set
\[
 s^{(0)}= \sum_{\substack{J,K\\ j\not\in K\\ |J|+|K|=k}}
 \sum_{i=1}^{\rank E}
 s^{\vecv}_{i,J,K}\cdot
 \omega_{J,K}\cdot v_i,
\quad\quad
 s^{(1)}=\sum_{\substack{J,K \\ j\in K \\ |J|+|K|=k}}
  \sum_{i=1}^{\rank E}
 s^{\vecv}_{i,J,K}\cdot
 \omega_{J,K}
 \cdot v_i.
\]
The decomposition $s=s^{(0)}+s^{(1)}$ is independent of
the choice of a frame $\vecv$.
We also set
\[
 \Pi_{<0}(s^{(1)},\vecv):=
  \sum_{\substack{J,K \\ j\in K \\ |J|+|K|=k}}
  \sum_{i=1}^{\rank E}
  \sum_{m<0}
  \gbigf_{j,m}(s^{\vecv}_{i,J,K})
  e^{\sqrt{-1}m\theta_j}
  \cdot
 \omega_{J,K}
 \cdot v_i.
\]
The following lemma is easy to see.
\begin{lem}
Let $U'\subset U$ be an $(S^1)^{\ell}$-invariant subset.
If $\Pi_{<0}(s^{(1)},\vecv)=0$ on $U'$
for a frame $\vecv$ of $\nbigp_{\veczero}\nbigelambda$,
we obtain
$\Pi_{<0}(s^{(1)},\vecv')=0$ on $U'$
for any frame $\vecv'$ of $\nbigp_{\veczero}\nbigelambda$.
\hfill\qed
\end{lem}

We set $T_j(p):=\bigl\{\sigma\in\gbigs_{\ell}\,|\,
\sigma(p)=j
\bigr\}$.
\begin{lem}
\label{lem;22.2.20.4}
Assume that $U=\prod_{i=1}^{n}\{|z_i|<\epsilon_i\}$
for some $\epsilon_i>0$,
and that $s=0$ on $\{|z_j|>\epsilon/2\}$.
Suppose that
the support of
$\Pi_{<0}(s^{(1)},\vecv)$ is contained in
 $\bigcup_{p\geq c}\bigcup_{\sigma\in T_j(p)}
 Z(\vecC,\sigma)$
for some $1\leq c\leq \ell$.
Then, there exists
$t_c\in \lefttop{I}F\nbigc^{k-1}_{L^2}(\nbigelambda,\DDlambda,h)(U)$
such that
(i) $t_c=t_c^{(0)}$,
(ii) the support of
$\Pi_{<0}\bigl(
 (s-d\zbar_j\del_{\zbar_j}(t_c))^{(1)},\vecv
\bigr)$
is contained in 
 $\bigcup_{p>c}\bigcup_{\sigma\in T_j(p)}
 Z(\vecC,\sigma)$.
\end{lem}
\pf
Let $\{\chi_{\sigma}\}$ denote a partition of unity
of $\bigcup_{p\geq c}\bigcup_{\sigma\in T_j(p)}Z(\vecC,\sigma)$
subordinate to the open covering
in Lemma \ref{lem;22.2.20.3}
for $S=\bigcup_{p\geq c}T_j(p)$.
Note that $\chi_{\sigma}s\in
\lefttop{I}F\nbigc^{\bullet}_{L^2}(\nbigelambda,\DDlambda,h)(U)$
for any $\sigma\in T_j(c)$.
Clearly,
$(\chi_{\sigma}s)^{(i)}
=\chi_{\sigma}s^{(i)}$ hold for $i=0,1$.
We have the expression
\[
 \chi_{\sigma}s^{(1)}
 =\sum_{\substack{J,K \\ j\in K\\ |J|+|K|=k}}
 \sum_{i=1}^{\rank E}
 \sum_{m\in\seisuu}
 \gbigf_{j,m} (\chi_{\sigma}s_{i,J,K}^{\vecw^{\sigma}})
 e^{\sqrt{-1}m\theta_j}
 \omega_{J,K} w_i^{\sigma}.
\]
Let $\epsilon(J,K,j)\in\{\pm 1\}$
be determined by
$(d\zbar_j/\zbar_j)\omega_{J,K\setminus\{j\}}=
\epsilon(J,K,j)\omega_{J,K}$.
We set
\[
 t^{\sigma}=
 \sum_{\substack{J,K \\ j\in K\\ |J|+|K|=k}}
 \sum_{i=1}^{\rank E}
 \epsilon(J,K,j)
 \nbiga\left(
 \sum_{m<0}
 \gbigf_{j,m}(\chi_{\sigma}s^{\vecw^{\sigma}}_{i,J,K})
 e^{\sqrt{-1}m\theta_j}
 \right)
 \omega_{J,K\setminus\{j\}}w_i^{\sigma}.
\]
Here, $\nbiga$ is given as in Lemma \ref{lem;22.1.25.20}
for the variable $z_j$.
By the construction, we have
$t^{\sigma}_{|Z(\vecC,\sigma)\cap U}
\in
L^2(\nbigelambda\otimes\Omega^{k-1},h)(Z(\vecC,\sigma)\cap U)$,
and
\[
\Pi_{<0}\bigl(
\chi_{\sigma}s^{(1)},\vecw^{\sigma}
\bigr)_{|Z(\vecC,\sigma)\cap U}
=(d\zbar_j\wedge \del_{\zbar_j}t^{\sigma})_{|Z(\vecC,\sigma)\cap U}.
\]
We obtain
\[
 \Pi_{<0}\bigl(
 \chi_{\sigma}s^{(1)}
 -d\zbar_j\wedge \del_{\zbar_j}t^{\sigma},\,\,
 \vecv
\bigr)_{|Z(\vecC,\sigma)\cap U}=0.
\]

\begin{lem}
\label{lem;22.2.20.20}
$(\DDlambda t^{\sigma})_{|Z(\vecC,\sigma)\cap U}
\in
L^2(\nbigelambda\otimes\Omega^{k},h)(Z(\vecC,\sigma)\cap U)$.
\end{lem}
\pf
Let $A$ be the diagonal matrix valued $1$-form
whose $(p,p)$-th term is equal to
$\sum_{j=1}^{\ell} \alpha_j(w_p^{\sigma})\frac{dz_j}{z_j}$.
We define a $\lambda$-connection $\DD^{\lambda(0)}$
by
$\DD^{\lambda(0)}\vecw^{\sigma}=\vecw^{\sigma}A$.
We set $f:=\DDlambda-\DD^{\lambda(0)}$.
Note that $f$ is bounded with respect to $h$ and the Poincar\'{e} metric
on $Z(\vecC,\sigma)$
by Proposition \ref{prop;22.2.21.3}.
Note also that 
$(d\zbar_j/\zbar_j)\cdot \zbar_j\del_{\zbar_j}(s^{(0)})
+\DD^{\lambda(0)}(s^{(1)})+f(s^{(1)})$
is an $L^2$-section.
By applying Proposition \ref{prop;22.1.25.21}
to each coefficients of $w_i^{\sigma}$,
we obtain that
$(\DD^{\lambda(0)}t^{\sigma})_{|Z(\vecC,\sigma)\cap U}
\in
L^2(\nbigelambda\otimes\Omega^k,h)(Z(\vecC,\sigma)\cap U)$.
Then, we obtain the claim of Lemma \ref{lem;22.2.20.20}.
\hfill\qed

\vspace{.1in}
By Lemma \ref{lem;22.4.27.30}
and Lemma \ref{lem;22.2.20.20},
we have
\[
  \Bigl(
 1-\sum_{p>c}\sum_{\sigma'\in T_j(p)}
 \chi_{\sigma'}
 \Bigr)\cdot
 t^{\sigma}
 \in
 \lefttop{I}F\nbigc^{\bullet}_{L^2}(\nbigelambda,\DDlambda,h)(U),
\]
and the section
\[
t_c:=
 \sum_{\sigma\in T_j(c)}
 \Bigl(
 1-\sum_{p>c}\sum_{\sigma'\in T_j(p)}
 \chi_{\sigma'}
 \Bigr)\cdot
 t^{\sigma}
\]
satisfies the desired condition.
Thus, we obtain Lemma \ref{lem;22.2.20.4}.
\hfill\qed

\vspace{.1in}

By applying Lemma \ref{lem;22.2.20.4}
in an inductive way,
we obtain the following lemma.

\begin{lem}
\label{lem;22.2.21.10}
We assume that $U=\prod_{j=1}^{n}\{|z_j|<\epsilon_j\}$
for some $\epsilon_j>0$,
and that $s=0$ on $\{|z_j|>\epsilon/2\}$.
Then, there exists a section
$t\in\nbigc^{k-1}_{L^2}(\nbigelambda,\DDlambda,h)(U)$
such that 
$(s-\DDlambda t)^{(1)}\in \lefttop{\Itilde}F
L^2(\nbigelambda\otimes\Omega^k,h)(U)$.
\hfill\qed
\end{lem}

We obtain the following lemma as a consequence of
Lemma \ref{lem;22.2.21.10}.

\begin{lem}
\label{lem;22.2.21.11}
For any $U\in\gbigu(O)$
and $s\in\lefttop{I}\nbigc^k_{L^2}(\nbigelambda,\DDlambda,h)(U)$,
there exist $U_1\in\gbigu(O,U)$
and $t\in \lefttop{I}\nbigc^{k-1}_{L^2}(\nbigelambda,\DDlambda,h)(U_1)$
such that 
$(s-\DDlambda t)^{(1)}\in \lefttop{\Itilde}F
L^2(\nbigelambda\otimes\Omega^k,h)(U_1)$.
\hfill\qed
\end{lem}

Let $U\in\gbigu(O)$
and $s\in\lefttop{I}\nbigc^k_{L^2}(\nbigelambda,\DDlambda,h)(U)$.
Let $U_1$ and $t$ be as in Lemma \ref{lem;22.2.21.11}.
We set $\stilde=s-\DDlambda t$ on $U_1$.
We have the decomposition
$\stilde=\stilde^{(0)}+\stilde^{(1)}$.
There exists the expression:
\[
 \stilde^{(0)}
 =\sum_{\substack{J,K \\ j\not\in K \\ |J|+|K|=k}}
 \sum_m
 \gbigf_{j,m}(\stilde^{\vecv}_{i,J,K})
 e^{\sqrt{-1}m\theta_j}\cdot \omega_{J,K}\cdot
 v_i.
\]
If $\DDlambda s=0$,
we obtain $\DDlambda\stilde=0$,
and
$\del_{\zbar_j}\Bigl(
\gbigf_{j,m}(\stilde^{\vecv}_{i,J,K})
e^{\sqrt{-1}m\theta_j}
\Bigr)=0$ $(m<0)$ if $|J|+|K|=k$
with $j\not\in K$.
Because of the $L^2$-property,
we obtain that
$\gbigf_{j,m}(\stilde^{\vecv}_{i,J,K})=0$ $(m<0)$,
i.e.,
$\stilde\in \lefttop{\Itilde}F\nbigc^k_{L^2}(\nbigelambda,\DDlambda,h)(U_1)$.
Suppose moreover $\stilde=\DDlambda u$
for some $u\in \lefttop{I}F\nbigc^{k-1}_{L^2}(\nbigelambda,\DDlambda,h)(U_1)$.
For the decomposition $u=u^{(0)}+u^{(1)}$,
we may assume that
$u^{(1)}\in \lefttop{\Itilde}L^2(\nbigelambda\otimes\Omega^{k-1},h)(U_1)$.
Then, by a similar argument, we obtain that
$u^{(0)}\in \lefttop{\Itilde}L^2(\nbigelambda\otimes\Omega^{k-1},h)(U_1)$.
Hence,
the morphism (\ref{eq;22.1.25.22}) is a quasi-isomorphism.
\hfill\qed

\subsubsection{Family version}

Let $\lambda_0\in\cnum$.
We have the family of regular filtered $\lambda$-flat bundles
$(\nbigpzero_{\ast}\nbige,\DDlambda)$.
For any $(S^1)^{\ell}$-invariant open subset $U$ of
$\nbigxzero$,
let $L^2(\nbige\otimes\Omega^k,h)(U)$
denote the space of sections $s$ of
$(\nbige\otimes\Omega^k)_{|U\setminus\nbighzero}$
such that
(i) $|s|_h$ is $L^2$,
(ii) $s$ is holomorphic with respect to $\lambda$.

Let $\vecv$ be a holomorphic frame of
$\nbigpzero_{\veczero}\nbige$.
Any section $s\in L^2(\nbige\otimes \Omega^k,h)(U)$
is expressed as
\[
 s=
 \sum_{\substack{J,K\\ |J|+|K|=k}}\sum_{i=1}^{\rank E}
 \sum_{\vecp\in\seisuu^{\ell}}
 \gbigf_{\ellsitabar,\vecp}
 (s^{\vecv}_{i,J,K})e^{\sqrt{-1}\vecp\vectheta}\cdot
 \omega_{J,K}\cdot v_i.
\]
For any subset $I\subset\ellsitabar$,
let
$\lefttop{I}FL^2(\nbige\otimes \Omega^k,h)(U)$
denote the subspace of $s$
such that
$\gbigf_{\ellsitabar,\vecp}(s^{\vecv}_{i,J,K})=0$
unless $p_j\geq 0$ $(j\in I)$.
The subspace is independent of the choice of $\vecv$.
We set
\[
 \lefttop{I}F\nbigc^{k}_{L^2}\bigl(
 \nbige,\DD,h
 \bigr)(U):=
 \nbigc^{k}_{L^2}\bigl(
 \nbige,\DD,h
 \bigr)(U)
 \cap
 \lefttop{I}F
 L^2(\nbige\otimes \Omega^k,h)(U).
\]
Thus, we obtain a subcomplex
$\lefttop{I}F\nbigc^{\bullet}_{L^2}
\bigl(
\nbige,\DD,h
\bigr)(U)$
of
$\nbigc^{\bullet}_{L^2}
\bigl(
\nbige,\DD,h
\bigr)(U)$.

Let $(\gbigu(\lambda,O),\prec)$
be a directed set of
$(S^1)^{\ell}$-invariant neighbourhoods of $(\lambda,O)$,
where $U_1\prec U_2$ is defined to be
$U_1\supset U_2$.
We set
\[
\lefttop{I}F\nbigc^{\bullet}_{L^2}
\bigl(
\nbige,\DD,h
\bigr)_{(\lambda,O)}
=\varinjlim_{U\in\gbigu(\lambda,O)}
\lefttop{I}F\nbigc^{\bullet}_{L^2}
\bigl(
\nbige,\DD,h
\bigr)(U).
\]
The following proposition is similar to
Proposition \ref{prop;22.1.28.10}.
\begin{prop}
The natural morphism
$\lefttop{I}F\nbigc^{\bullet}_{L^2}
\bigl(
\nbige,\DD,h
\bigr)_{(\lambda,O)}
\lrarr
\nbigc^{\bullet}_{L^2}
\bigl(
\nbige,\DD,h
\bigr)_{(\lambda,O)}$
is a quasi-isomorphism.
\hfill\qed
\end{prop}

\subsection{Reduction to the complex of holomorphic sections}
\label{subsection;22.3.17.52}

\subsubsection{Statement for a fixed $\lambda$}

To simplify the description,
$\Omega^{\bullet,0}_X$
and
$\Omega^{0,\bullet}_X$
are denoted by $\Omega^{\bullet,0}$
and $\Omega^{0,\bullet}$,
respectively.
For $1\leq j\leq n$,
let $\Omega^{0,\bullet}_{\leq j}$
denote the subbundle of $\Omega^{0,\bullet}$
generated by $d\zbar_i$ $(i\leq j)$.

Let
$\nbigc^{\bullet}_{L^2,\leq j,\ttX(\Lambda)_{\geq 0}}(\nbigelambda,\DDlambda,h)
\subset 
\nbigc^{\bullet}_{L^2,\ttX(\Lambda)_{\geq 0}}(\nbigelambda,\DDlambda,h)$
denote the subsheaf of local sections $s$
satisfying the following condition.
\begin{itemize}
 \item $s$ is a local section of
       $(j_{X\setminus H,\ttX(\Lambda)_{\geq 0}})_{\ast}
       \Bigl(
       \nbigelambda\otimes
 \Omega^{\bullet,0}
 \otimes
       \Omega^{0,\bullet}_{\leq j}\Bigr)$,
       and satisfies $\del_{\zbar_i}s=0$ for $i>j$.
\end{itemize}
We obtain the subcomplex
$\nbigc^{\bullet}_{L^2,\leq j,\ttX(\Lambda)_{\geq 0}}
(\nbigelambda,\DDlambda,h)$
of 
$\nbigc^{\bullet}_{L^2,\ttX(\Lambda)_{\geq 0}}(\nbigelambda,\DDlambda,h)$.
By the construction,
\[
\nbigc^{\bullet}_{L^2,\leq n,\ttX(\Lambda)_{\geq 0}}
(\nbigelambda,\DDlambda,h)=
\nbigc^{\bullet}_{L^2,\ttX(\Lambda)_{\geq 0}}(\nbigelambda,\DDlambda,h). 
\]

For an $(S^1)^{\ell}$-invariant open subset
$U$ of $X(\ttX(\Lambda)_{\geq 0})$,
and for $J\subset\ellsitabar$,
we define
\[
\lefttop{J}F\nbigc^{\bullet}_{L^2,\leq j,\ttX(\Lambda)_{\geq 0}}
(\nbigelambda,\DDlambda,h)(U)
\subset
\nbigc^{\bullet}_{L^2,\leq j,\ttX(\Lambda)_{\geq 0}}
(\nbigelambda,\DDlambda,h)(U)
\]
as in \S\ref{subsection;22.1.27.1}.
In particular,
we obtain the subcomplex
\[
 \lefttop{\ellsitabar}F
 \nbigc^{\bullet}_{L^2,\leq j,\ttX(\Lambda)_{\geq 0}}
 (\nbigelambda,\DDlambda,h)(U)
\subset
\nbigc^{\bullet}_{L^2,\leq j,\ttX(\Lambda)_{\geq 0}}
(\nbigelambda,\DDlambda,h)(U).
\]
Let $Q$ be a point of
$p_1^{-1}(O)\subset X(\ttX(\Lambda)_{\geq 0})$.
Let $(\gbigu(Q),\prec)$ be the directed set of
$(S^1)^{\ell}$-invariant open neighbourhoods of $Q$,
where $U_1\prec U_2$ is defined to be $U_1\supset U_2$.
For any $U\in\gbigu(Q)$,
we set
$\gbigu(Q,U):=\bigl\{
 U'\in\gbigu(Q)\,\big|\,
 U\prec U'
\bigr\}$.
We set
\[
  \lefttop{\ellsitabar}F
  \nbigc^{\bullet}_{L^2,\leq j,\ttX(\Lambda)_{\geq 0}}
  (\nbigelambda,\DDlambda,h)_Q:=
  \varinjlim_{U\in\gbigu(Q)}
 \lefttop{\ellsitabar}F
  \nbigc^{\bullet}_{L^2,\leq j,\ttX(\Lambda)_{\geq 0}}
  (\nbigelambda,\DDlambda,h)(U).
\]
We shall prove the following proposition
in \S\ref{subsection;22.4.27.40}
after preliminaries.

\begin{prop}
\label{prop;22.2.23.43}
For any $0\leq j\leq n$,
the natural morphism
\begin{equation}
\label{22.2.21.12}
 \lefttop{\ellsitabar}F
 \nbigc^{\bullet}_{L^2,\leq j,\ttX(\Lambda)_{\geq 0}}
 (\nbigelambda,\DDlambda,h)_Q
 \lrarr
 \lefttop{\ellsitabar}F
 \nbigc^{\bullet}_{L^2,\leq n,\ttX(\Lambda)_{\geq 0}}
 (\nbigelambda,\DDlambda,h)_Q
\end{equation}
is a quasi-isomorphism.
\end{prop}

It is enough to consider the case
$Q\in \overline{Z(\vecC)}$,
which we shall assume.

\subsubsection{Preliminary (1)}

Let $\vecw=\vecw^{\id}$ be a frame of
$\nbigp_{\veczero}\nbigelambda$ as in \S\ref{subsection;22.1.26.2}.
Let $U$ be an $(S^1)^{\ell}$-invariant open subset of
$\Zbar(\vecC)$.
An $L^2$-section $s$ of
$\nbigelambda\otimes\Omega^{\bullet,0}\otimes\Omega^{0,\bullet}_{\leq j}$
on $U$ has the expression
\[
s=\sum_{i=1}^{\rank E} s^{\vecw}_i\cdot w_i
 =\sum_{J,K\subset\nbar}\sum_{i}
 s^{\vecw}_{J,K,i}\cdot \omega_{J,K}\cdot w_i
 =\sum_{J,K}\sum_i
 \sum_{m\in\seisuu}
 \gbigf_{j,m}(s^{\vecw}_{J,K,i})e^{\sqrt{-1}m\theta_j}
 \cdot \omega_{J,K}
 \cdot w_i. 
\]
We set
$\gbigf_{j,m}(s^{\vecw}_{i}):=
\sum_{J,K}\gbigf_{j,m}(s^{\vecw}_{J,K,i})\cdot \omega_{J,K}$.
The following holds:
\[
 s=\sum_{i=1}^{\rank E}
 \sum_{m\in\seisuu}
  \gbigf_{j,m}(s^{\vecw}_i)e^{\sqrt{-1}m\theta_j}\cdot w_i.
\]
We set
\[
 \nbigb_{\neq 0}(s)
=\sum_{a_j(w_i)\neq 0}s^{\vecw}_iw_i
+\sum_{a_j(w_i)=0}
 \sum_{m\neq 0}
\gbigf_{j,m}(s^{\vecw}_{i})e^{\sqrt{-1}m\theta_j}
 w_i,
 \quad
\nbigb_{0}(s)
=\sum_{a_j(w_i)=0}
\gbigf_{j,0}(s^{\vecw}_{i})w_i.
\]

Let 
$\lefttop{j}F_{<0}
\lefttop{\ellsitabar}F
L^2(\nbigelambda\otimes\Omega^k,h)(U)
\subset
\lefttop{\ellsitabar}F
L^2(\nbigelambda\otimes\Omega^k,h)(U)$
denote the subspace of sections $s$
such that
$\nbigb_0(s)=0$.

\subsubsection{Preliminary (2)}
\label{subsection;22.2.23.100}

Let $\Omega^{\bullet,0}_{\neq j}$
denote the subbundle of $\Omega^{\bullet,0}$
generated by $dz_i$ $(i\neq j)$.
By setting
$\Omega^{\bullet,\bullet}_{\neq j,\leq j-1}:=
\Omega^{\bullet,0}_{\neq j}
\otimes
\Omega^{0,\bullet}_{\leq j-1}$,
we obtain the decomposition:
\[
\Omega^{\bullet,0}\otimes
\Omega^{0,\bullet}_{\leq j}
=
\Omega^{\bullet,\bullet}_{\neq j,\leq j-1}
\oplus
(d\zbar_j/\zbar_j)
\Omega^{\bullet,\bullet}_{\neq j,\leq j-1}
\oplus
(dz_j/z_j)
\Omega^{\bullet,\bullet}_{\neq j,\leq j-1}
\oplus
(d\zbar_j/\zbar_j)
(dz_j/z_j)
\Omega^{\bullet,\bullet}_{\neq j,\leq j-1}.
\]

A section $s$ of
$\nbigelambda\otimes\Omega^{\bullet,0}\otimes \Omega^{0,\bullet}_{\leq j}$
has the decomposition
\[
 s=s^{(0)}+(d\zbar_j/\zbar_j)s^{(1)}+(dz_j/z_j)s^{(2)}
 +(d\zbar_j/\zbar_j)\,(dz_j/z_j)\,s^{(3)}.
\]
Here, $s^{(a)}$ has the expression
$s^{(a)}=\sum s^{(a)}_iw_i$
such that $s^{(a)}_i$ are sections of
$\Omega^{\bullet,\bullet}_{\neq j,\leq j-1}$.

On $\lefttop{J}F\nbigc^{k}_{L^2,\leq j,\ttX(\Lambda)_{\geq 0}}
(\nbigelambda,\DDlambda,h)$,
we have the decomposition 
\[
 \DDlambda
 =\DD^{\lambda(1,0)}
 +\DD^{\lambda(0,1)}
 =\DD^{\lambda(1,0)}_{\neq j}
 +\DDlambda_{z_j}\,dz_j
 +\DD^{\lambda(0,1)}_{\leq j-1}
 +\DD_{\zbar_j}^{\lambda}\,d\zbar_j.
\]
We set
$\DD^{\lambda}_{\neq j}
=\DD^{\lambda(1,0)}_{\neq j}
+\DD^{\lambda(0,1)}_{\leq j-1}$.
We obtain the decomposition
$\DDlambda
 =\DDlambda_{\neq j}
 +\DDlambda_{z_j}\,dz_j
 +\DD_{\zbar_j}^{\lambda}\,d\zbar_j$.

For a section $s$ of
$\lefttop{J}F\nbigc^{k}_{L^2,\leq j,\ttX(\Lambda)_{\geq 0}}
(\nbigelambda,\DDlambda,h)$
satisfying $\DDlambda(s)=0$,
we obtain the following relations:
\[
 \DDlambda_{\neq j}s^{(0)}=0,
 \quad
 d\zbar_j\DDlambda_{\zbar_j}(s^{(0)})
 -(d\zbar_j/\zbar_j)\,\DDlambda_{\neq j}s^{(1)}=0,
 \quad
 dz_j\DDlambda_{z_j}(s^{(0)})
-(dz_j/z_j)\DDlambda_{\neq j}(s^{(2)})=0,
\]
\[
 dz_j(d\zbar_j/\zbar_j)\DDlambda_{z_j}(s^{(1)})
+d\zbar_j(dz_j/z_j)\DDlambda_{\zbar_j}s^{(2)}
+\,(d\zbar_j/\zbar_j)\,(dz_j/z_j)\DDlambda_{\neq j}s^{(3)}=0.
\]

\subsubsection{Decomposition of the bundle}
\label{subsection;22.2.24.2}

For $(\veca,\vecalpha,\veck)\in \openclosed{0}{1}^j
\times\cnum^j\times\seisuu^j$,
let $\nbigv_{\veca,\vecalpha,\veck}$
denote the subbundle of $\nbigp_{\veczero}\nbigelambda$
generated by
$w_i$ such that
\[
 q_{\jbar}(\veca(w_i),\vecalpha(w_i),\veck(w_i))
 =(\veca,\vecalpha,\veck).
\]
We obtain the decomposition
\begin{equation}
\label{eq;22.2.22.23}
 \nbigp_{\veczero}\nbigelambda
 =\bigoplus_{\veca\in\openclosed{-1}{0}^j}
 \bigoplus_{\vecalpha\in\cnum^j}
 \bigoplus_{\veck\in\seisuu^j}
 \nbigv_{\veca,\vecalpha,\veck}.
\end{equation}
For a section $s$ of
$\nbigp_{\veczero}\nbigelambda\otimes \nbigk$,
where $\nbigk$ denotes any sheaf,
we obtain the decomposition
\[
 s=
 \sum_{\veca\in\openclosed{-1}{0}^j}
 \sum_{\vecalpha\in\cnum^j}
 \sum_{\veck\in\seisuu^j}
 s_{\veca,\vecalpha,\veck},
\]
according to (\ref{eq;22.2.22.23}).

We obtain the decomposition
\[
\DDlambda_{\neq j}
=\DD^{\prime\lambda}_{\neq j}
+\sum_{(\veca,\vecalpha,\veck)\neq (\veca',\vecalpha',\veck')}
 G_{\neq j,(\veca,\vecalpha,\veck),(\veca',\vecalpha',\veck')},
\]
where $\DD^{\prime\lambda}_{\neq j}$
preserves the decomposition (\ref{eq;22.2.22.23}),
and
$G_{\neq j,(\veca,\vecalpha,\veck),(\veca',\vecalpha',\veck')}$
are sections of
\[
 \Hom\Bigl(
 \nbigv_{\veca',\vecalpha',\veck'},
 \nbigv_{\veca,\vecalpha,\veck}
  \Bigr)
  \otimes\Omega^1_{\neq j}.
\]

We obtain the following lemma from
the estimates in \S\ref{subsection;22.2.22.24}.
\begin{lem}
Assume one of the following.
\begin{itemize}
 \item $(\veca,\vecalpha)\neq (\veca',\vecalpha')$.
 \item $(\veca,\vecalpha)=(\veca',\vecalpha')$
       and there exists $1\leq m\leq j$
       such that $k_m>k'_m.$
\end{itemize}
Then, there exists $\epsilon>0$ such that the following holds
on $Z(\vecC)$:
\[
 |G_{\neq j,(\veca,\vecalpha,\veck),
 (\veca',\vecalpha',\veck')}|_{h,g_{X\setminus H}}
 =O(|z_j|^{\epsilon}).
\]
\hfill\qed
\end{lem}

We define the differential operator $z_j\del_{z_j}$
on $\nbigp_{\veczero}\nbigelambda$ by
\[
 z_j\del_{z_j}\Bigl(\sum_i f_iw_i\Bigr)
=\sum_i z_j\del_{z_j}(f_i)w_i.
\]
The following lemma is easy to see.
\begin{lem}
If $\nbigb_0(f)=f$,
then we obtain
$z_j\del_{z_j}f=\zbar_j\del_{\zbar_j}f
=\frac{1}{2}\sum r\del_r(f_i)w_i$. 
\hfill\qed
\end{lem}
We have the decomposition
\[
 z_j\DDlambda_{z_j}
 =z_j\del_{z_j}
 +\Bigl(
 \bigoplus_{\veca,\vecalpha,\veck}
  \alpha_j\id_{\nbigv_{\veca,\vecalpha,\veck}}
  \Bigr)
  +\sum_{(\veca,\vecalpha,\veck),(\veca',\vecalpha',\veck')}
  G_{j,(\veca,\vecalpha,\veck),(\veca',\vecalpha',\veck')},
\]
where
$G_{j,(\veca,\vecalpha,\veck),(\veca',\vecalpha',\veck')}$
are sections of
$\Hom\Bigl(
 \nbigv_{\veca',\vecalpha',\veck'},
 \nbigv_{\veca,\vecalpha,\veck}
 \Bigr)$.
We note that
$G_{j,(\veca,\vecalpha,\veck),(\veca,\vecalpha,\veck)}$
are not necessarily $0$.
 
We obtain the following lemma from
the estimates in \S\ref{subsection;22.2.22.24}.

\begin{lem}
\label{lem;22.2.23.3}
Assume one of the following holds:
\begin{itemize}
 \item $(\veca,\vecalpha)\neq (\veca',\vecalpha')$.
 \item $(\veca,\vecalpha)=(\veca',\vecalpha')$,
       and there exists $1\leq m<j$
       such that $k_m>k_m'$.
 \item $(\veca,\vecalpha)=(\veca',\vecalpha')$,
       $k_m\leq k_m'$ $(1\leq m\leq j-1)$,
       and
       $k_j>k_{j}'-2$.
\end{itemize}
 Then, there exists $\epsilon>0$ such that
$\bigl|
 G_{j,(\veca,\vecalpha,\veck),(\veca',\vecalpha',\veck')}
 \bigr|_{h,g_{X\setminus H}}
=O\bigl(|z_j|^{\epsilon}\bigr)$
on $Z(\vecC)$.
\hfill\qed
\end{lem}

We also obtain the following lemma from
the estimates in \S\ref{subsection;22.2.22.24}.

\begin{lem}
\label{lem;22.2.24.1}
Assume  $k_m\leq k'_{m}$ for $m\leq j-1$,
and $k_j\leq k_j'-2$.
 Then,
\[
 \bigl|
 G_{j,(\veca,\vecalpha,\veck),(\veca',\vecalpha',\veck')}
 \bigr|_{h,g_{X\setminus H}}
 =O\bigl(
 (-\log|z_j|)^{-1}
 \bigr)
\]
 holds on $Z(\vecC)$.
\hfill\qed
\end{lem}

We note the following lemma.
\begin{lem}
Suppose that $\veck(0)$ and $\veck(1)$ satisfy
(i) $k(0)_{m}=k(1)_{m}$ $(m=1,\ldots,j-1)$,
(ii) $k(0)_{j}-k(0)_{j-1}=1$,
(iii) $k(1)_{j}-k(1)_{j-1}=-1$.
Then, the restriction of
\[
G_{j,(\veca,\vecalpha,\veck(1)),(\veca,\vecalpha,\veck(0))} 
\]
to $H_{\jbar}$  is an isomorphism.
Hence, by shrinking $X$,
we may assume that
$G_{j,(\veca,\vecalpha,\veck(1)),(\veca,\vecalpha,\veck(0))}$
is an isomorphism.
\hfill\qed 
\end{lem}

\subsubsection{Preliminary consideration (1)}

Let $U=\prod_{p=1}^{n}\{|z_p|<\epsilon_p\}$.
Let
$s\in \lefttop{\ellsitabar}F
\nbigc^{k}_{L^2,\leq j,\ttX(\Lambda)_{\geq 0}}
(\nbigelambda,\DDlambda,h)(U\cap \Zbar(\vecC))$
such that $s=0$ on $\{|z_j|>\epsilon_j/2\}\cap(U\cap\Zbar(\vecC))$.
We formally set $k_0(w_i)=0$ for any $i=1,\ldots,\rank E$.

\begin{lem}  
\label{lem;22.2.22.3}
Let $1\leq i\leq \rank E$.
If $(a_j(w_i),k_j(w_i)-k_{j-1}(w_i))\neq (0,-1),(0,1)$,
then there exists 
\[
 u\in\lefttop{\ellsitabar}F
 \nbigc^{k-1}_{L^2,\leq j,\ttX(\Lambda)_{\geq 0}}
 (\nbigelambda,\DDlambda,h)(U\cap\Zbar(\vecC))
\]
 such that
(i) $u_{i'}=0$ $(i'\neq i)$,
(ii) $u^{(b)}=0$ $(b=1,3)$, 
(iii) $(s-\DDlambda u)^{(b)}_{i}=0$ $(b=1,3)$.
\end{lem}
\pf
We regard $s^{(b)}_{J,K,i}$
as sections on $U$.
For $b=1,3$,
we set
$u^{(b-1)}_{J,K,i}=\nbiga(s^{(b)}_{J,K,i})$
for any $J,K$.
(See \S\ref{subsection;22.1.26.1} for $\nbiga$.)
We set
\[
 u^{(b-1)}= \sum_{J,K} u^{(b-1)}_{J,K,i}\omega_{J,K}w_i.
\]
Then, the restriction of
$u^{(0)}$ and $(dz_j/z_j)u^{(2)}$
to $U\cap\Zbar(\vecC)$
are $L^2$-sections.
By the construction, we have
\[
 (d\zbar_j/\zbar_j) \zbar_j\DDlambda_{\zbar_j}u^{(0)}
=(d\zbar_j/\zbar_j)s^{(1)}_iw_i,
\quad\quad
 (d\zbar_j/\zbar_j)\zbar_j\DDlambda_{\zbar_j}\bigl(
 -(dz_j/z_j) u^{(2)}\bigr)
=(d\zbar_j/\zbar_j)(dz_j/z_j)s^{(3)}_iw_i.
\]
By Proposition \ref{prop;22.1.25.21},
$s^{(0)}-\DDlambda_{\neq j}(u^{(0)})$
and
$s^{(2)}-dz_j\DDlambda_{z_j}(u^{(0)})
+\DDlambda_{\neq j}\bigl((dz_j/z_j)u^{(2)}\bigr)$
are $L^2$.
Hence, $u=u^{(0)}-(dz_j/z_j)u^{(2)}$
satisfies the desired condition.
\hfill\qed

\begin{lem}
\label{lem;22.2.22.30}
If $(a_j(w_i),k_j(w_i)-k_{j-1}(w_i))=(0,-1),(0,1)$,
then there exists
\[
 u\in\lefttop{\ellsitabar}F
 \nbigc^{k-1}_{L^2,\leq j,\ttX(\Lambda)_{\geq 0}}
 (\nbigelambda,\DDlambda,h)(U\cap\Zbar(\vecC))
\]
such that
(i) $u_{i'}=0$ $(i'\neq i)$,
(ii) $u^{(b)}=0$ $(b=1,3)$, 
(iii) 
$\gbigf_{j,m}\Bigl(
 (s-\DDlambda u)^{(b)}_{J,K,i}
 \Bigr)=0$ $(b=1,3)$
for any $J,K$ and $m\neq 0$.
\end{lem}
\pf
For $b=1,3$,
we set 
$u^{(b-1)}_{i,J,K}:=
\nbiga\Bigl(
\sum_{m>0}\gbigf_{j,m}(s^{(b)}_{i,J,K})e^{\sqrt{-1}m\theta_j}
\Bigr)$
for any $J,K$.
We set
\[
 u^{(b-1)}= \sum_{i,J,K} u^{(b-1)}_{i,J,K}\omega_{J,K}w_i.
\]
Then, as in the proof of Lemma \ref{lem;22.2.22.3},
$u=u^{(0)}-(dz_j/z_j)u^{(2)}$
has the desired property.
\hfill\qed

\begin{lem}
\label{lem;22.2.22.31}
Suppose
$(a_j(w_i),k_j(w_i)-k_{j-1}(w_i))=(0,-1),(0,1)$.
Moreover, we assume that
there exists $\epsilon>0$ such that the following holds.
\begin{itemize}
 \item $|z|^{-\epsilon}s^{(b)}_iw_i$ $(b=1,3)$ are $L^2$.
 \item $|z|^{-\epsilon}\bigl(
       \DDlambda(s_iw_i)\bigr)^{(b)}$ $(b=1,3)$ are $L^2$.
\end{itemize}
Then, there exists
$u\in\lefttop{\ellsitabar}F
 \nbigc^{k-1}_{L^2,\leq j,\ttX(\Lambda)_{\geq 0}}
 (\nbigelambda,\DDlambda,h)(U\cap\Zbar(\vecC))$
such that
(i) $u_{i'}=0$ $(i'\neq i)$,
(ii) $u^{(b)}=0$ $(b=1,3)$, 
(iii) $(s-\DDlambda u)^{(b)}_{J,K,i}=0$ $(b=1,3)$
for any $J,K$.
\end{lem}
\pf
We construct $u$ by using $\nbiga^{[-1]}$.
(See \S\ref{subsection;22.1.26.1} for $\nbiga^{[-1]}$.)
We can check the condition (iii)
by using Proposition \ref{prop;22.2.22.12}.
\hfill\qed

\subsubsection{Preliminary consideration (2)}

Let $U\in\gbigu(Q)$.
Let
$s\in \lefttop{\ellsitabar}F
\nbigc^{k}_{L^2,\leq j,\ttX(\Lambda)_{\geq 0}}
(\nbigelambda,\DDlambda,h)(U)$
such that $(\DDlambda s)^{(b)}=0$ $(b=1,3)$.
Let
$(\veca,\vecalpha)
\in\openclosed{-1}{0}^j\times\cnum^j$
such that $a_j=0$.
We set $\alpha:=\alpha_j$.
We shall prove the following proposition
in \S\ref{subsection;22.2.23.1}
and \S\ref{subsection;22.2.23.2}.
\begin{prop}
\label{prop;22.2.22.40}
There exists 
$U'\in\gbigu(Q,U)$ and
$u\in\lefttop{\ellsitabar}F
 \nbigc^{k-1}_{L^2,\leq j,\ttX(\Lambda)_{\geq 0}}
 (\nbigelambda,\DDlambda,h)(U')$
such that the following holds: 
\begin{itemize}
 \item $(s-\DDlambda u)^{(b)}_{\veca,\vecalpha,\veck}=0$
       for $b=1,3$
       for any $\veck\in\seisuu^{\ell}$.
 \item $(s-\DDlambda u)^{(b)}_{\veca',\vecalpha',\veck'}=
       s^{(b)}_{\veca',\vecalpha',\veck'}$
       for $b=1,3$,
       for any $\veck'\in\seisuu^{\ell}$
       and for $(\veca',\vecalpha')\neq (\veca,\vecalpha)$.
\end{itemize}
\end{prop}

\subsubsection{The case $\alpha\neq 0$}
\label{subsection;22.2.23.1}

Let us study the case $\alpha\neq 0$.

\begin{prop}
 \label{prop;22.2.22.32}
Assume that there exists
$\veck\in\seisuu^{\ell}$
such that $s^{(b)}_{\veca,\vecalpha,\veck'}=0$ for $b=1,3$
unless
$k_m'\leq k_m$ for any $1\leq m\leq j$.
Then, there exists
$U'\in\gbigu(Q,U)$ and
$t\in\lefttop{\ellsitabar}F
 \nbigc^{k-1}_{L^2,\leq j,\ttX(\Lambda)_{\geq 0}}
 (\nbigelambda,\DDlambda,h)(U')$
such that the following holds: 
\begin{itemize}
 \item $(s-\DDlambda t)^{(b)}_{\veca,\vecalpha,\veck}=0$
       for $b=1,3$.
 \item There exists $\epsilon>0$ such that
       $|z_j|^{-\epsilon}
       (\DDlambda t)^{(b)}_{\veca',\vecalpha',\veck'}$
       are $L^2$ for $b=1,3$
       and for $(\veca',\vecalpha',\veck')\neq (\veca,\vecalpha,\veck)$.
\end{itemize}
\end{prop}
\pf
By Lemma \ref{lem;22.2.22.3} and Lemma \ref{lem;22.2.22.30},
we may assume
$\nbigb_0(s^{(b)}_{\veca,\vecalpha,\veck})=
 s^{(b)}_{\veca,\vecalpha,\veck}$ for $b=1,3$.
We consider the following $L^2$-sections:
\[
t^{(3,1)}:=
-(d\zbar_j/\zbar_j)\alpha^{-1}
s^{(3)}_{\veca,\vecalpha,\veck},
\quad\quad
t^{(3,2)}:=
-(dz_j/z_j)\alpha^{-1}
\lambda
s^{(3)}_{\veca,\vecalpha,\veck},
\quad
 t^{(0)}:=\alpha^{-1}\nbigb_0(s^{(2)}_{\veca,\vecalpha,\veck})
 -\alpha^{-1}\lambda s^{(1)}_{\veca,\vecalpha,\veck}.
\]
 
\begin{lem}
\label{lem;22.3.15.20}
There exists $\epsilon>0$ such that
the following section is $L^2$:
\[
 |z_j|^{-\epsilon}
\Bigl(
(d\zbar_j/\zbar_j)(dz_j/z_j)
s^{(3)}_{\veca,\vecalpha,\veck}
 -
 \bigl(
 dz_j\DDlambda_{z_j}t^{(3,1)}
+d\zbar_j\del_{\zbar_j}t^{(3,2)}
\bigr)
\Bigr).
\] 
\end{lem}
\pf
The following holds:
\[
 (dz_j/z_j)\bigl(
 \lambda z_j\del_{z_j}+\alpha
 \bigr)t^{(3,1)}
 +(d\zbar_j/\zbar_j)
 (\zbar_j\del_{\zbar_j})
 t^{(3,2)}
 =(d\zbar_j/\zbar_j)(dz_j/z_j)
 s^{(3)}_{\veca,\vecalpha,\veck}.
\]
Then, we obtain the claim of Lemma \ref{lem;22.3.15.20}
from Lemma \ref{lem;22.2.23.3}.
\hfill\qed

\begin{lem}
There exists $\epsilon>0$ such that the following
section is $L^2$:
\begin{equation}
\label{eq;22.3.15.22}
|z_j|^{-\epsilon}
\Bigl(
 (d\zbar_j/\zbar_j)s^{(1)}_{\veca,\vecalpha,\veck}
 -\bigl(
 d\zbar_j\del_{\zbar_j} t^{(0)}
 +\DD^{\prime\lambda}_{\neq j}t^{(3,1)}
 \bigr)
 \Bigr).
\end{equation}
Moreover,
$\DD^{\prime\lambda}_{\neq j}t^{(0)}
-\nbigb_0(s^{(0)}_{\veca,\vecalpha,\veck})$
and 
\begin{equation}
\label{eq;22.3.15.21}
 (dz_j/z_j)
\nbigb_0(s^{(2)}_{\veca,\vecalpha,\veck})
-\bigl(
 (dz_j/z_j)(\lambda z_j\del_{z_j}+\alpha)t^{(0)}
 +\DD^{\prime\lambda}_{\neq j}t^{(3,2)}
 \bigr)
\end{equation}
are $L^2$.
In particular,
 $t^{(0)}+t^{(3,1)}+t^{(3,2)}\in
 \lefttop{J}F\nbigc^{k-1}_{L^2,\leq j,\ttX(\Lambda)_{\geq 0}}
 (\nbigelambda,\DDlambda,h)(U)$.
\end{lem}
\pf
There exists $\epsilon_1>0$ such that
the following sections are $L^2$:
\begin{equation}
\label{eq;22.3.15.11}
|z_j|^{-\epsilon_1}
 (d\zbar_j/\zbar_j)\Bigl(
 \zbar_j\del_{\zbar_j}
 \nbigb_0(s^{(0)}_{\veca,\vecalpha,\veck})
 -\DD^{\prime\lambda}_{\neq j}
 s^{(1)}_{\veca,\vecalpha,\veck}
 \Bigr),
\end{equation}
\begin{equation}
\label{eq;22.3.15.10}
|z_j|^{-\epsilon_1}
  (d\zbar_j/\zbar_j)(dz_j/z_j)
  \Bigl(
  \DD^{\prime\lambda}_{\neq j}
  s^{(3)}_{\veca,\vecalpha,\veck}
  -(\alpha+\lambda z_j\del_{z_j})
  s^{(1)}_{\veca,\vecalpha,\veck}
  +\zbar_j\del_{\zbar_j}
  \nbigb_0(s^{(2)}_{\veca,\vecalpha,\veck})
  \Bigr).
\end{equation}
The following sections are $L^2$:
\begin{equation}
\label{eq;22.3.15.12}
 \DD^{\prime\lambda}_{\neq j}
 \nbigb_0(s^{(0)}_{\veca,\vecalpha,\veck}),
 \quad\quad
 (dz_j/z_j)
 \Bigl(
 (\alpha+\lambda z_j\del_{z_j})\nbigb_0(s^{(0)}_{\veca,\vecalpha,\veck})
 -\DD^{\prime\lambda}_{\neq j}
 \nbigb_0(s^{(2)}_{\veca,\vecalpha,\veck})
 \Bigr).
\end{equation}
Because (\ref{eq;22.3.15.10}) is $L^2$,
there exists $\epsilon>0$
such that (\ref{eq;22.3.15.22})
is $L^2$.
We also obtain that
(\ref{eq;22.3.15.21}) is $L^2$.
Because (\ref{eq;22.3.15.11}) and (\ref{eq;22.3.15.12})
are $L^2$,
$\DD^{\prime\lambda}_{\neq j}t^{(0)}
-\nbigb_0(s^{(0)}_{\veca,\vecalpha,\veck})$ is $L^2$.
\hfill\qed

\vspace{.1in}

We set
\[
\stilde=s-\DDlambda(t^{(0)}+t^{(3,1)}+t^{(3,2)})
\in 
 \lefttop{J}F\nbigc^k_{L^2,\leq j,\ttX(\Lambda)_{\geq 0}}
 (\nbigelambda,\DDlambda,h)(U).
\]
Then, 
there exists $\epsilon>0$ such that
the following sections are $L^2$ for $b=1,3$:
\[
|z|^{-\epsilon}\stilde^{(b)}_{\veca,\vecalpha,\veck},\,\,\,
\quad
|z|^{-\epsilon}\bigl(
\stilde^{(b)}-s^{(b)}
\bigr)_{\veca',\vecalpha',\veck'}\,\,\,
\bigl(
(\veca',\vecalpha',\veck')
\neq
(\veca,\vecalpha,\veck)
\bigr).
\]
Then, we obtain Proposition \ref{prop;22.2.22.32}
from Lemma \ref{lem;22.2.22.30}
and Lemma \ref{lem;22.2.22.31}.
\hfill\qed

\vspace{.1in}

Let $\leq_{\seisuu^{j}}$
be the partial order defined by
$\veck\leq_{\seisuu^j}\veck'
\stackrel{\rm def}{\Longrightarrow}
k_m\leq k_m'$ for any $1\leq m\leq j$.
Let $\leq'_{\seisuu^j}$ be a total order on $\seisuu^{\ell}$
which is a refinement of $\leq_{\seisuu^j}$,
i.e.,
$\veck\leq_{\seisuu^j}\veck'$
implies
$\veck\leq'_{\seisuu^j}\veck'$.
Then,
we obtain the claim of Proposition \ref{prop;22.2.22.40}
in the case $\alpha\neq 0$,
from Proposition \ref{prop;22.2.22.32}
by using an easy induction
on $(\seisuu^j,\leq'_{\seisuu^j})$.
\hfill\qed

\subsubsection{The case $\alpha=0$}
\label{subsection;22.2.23.2}

Let us study the case $\alpha=0$.
By Lemma \ref{lem;22.2.22.3}
we may assume that
$s^{(b)}_{\veca,\vecalpha,\veck}=0$ $(b=1,3)$
unless $k_j-k_{j-1}=\pm 1$.
Moreover, by Lemma \ref{lem;22.2.22.30},
we may assume that
$s^{(b)}_{\veca,\vecalpha,\veck}=
\nbigb_0(s^{(b)}_{\veca,\vecalpha,\veck})$ $(b=1,3)$.

We set
$\nbigj_{\pm 1}:=
\bigl\{
\veck\in\seisuu\,\big|\,
k_j-k_{j-1}=\pm 1,\,\,\nbigv_{\veca,\vecalpha,\veck}\neq 0
\bigr\}$,
and
$\nbigj:=\nbigj_{1}
\cup\nbigj_{-1}$.
Let $\leq_{\seisuu^{j}}$
be the partial order defined by
$\veck\leq_{\seisuu^j}\veck'
\stackrel{\rm def}{\Longrightarrow}
k_m\leq k_m'$ for any $1\leq m\leq j$.
The induced partial order on $\nbigj$ is also denoted by
$\leq_{\seisuu^{j}}$.
We say
$\veck\lneq_{\seisuu^{j}}\veck'$
if $\veck\leq_{\seisuu^{j}}\veck'$
and $\veck\neq\veck'$.

\begin{prop}
\label{prop;22.2.23.42}
Let $\veck(0)\in\nbigj_1$ and $\veck(1)\in\nbigj_{-1}$
such that $k(0)_m=k(1)_m$ $(1\leq m\leq j-1)$.
We assume that
$s^{(b)}_{\veca,\vecalpha,\veck}=0$ $(b=1,3)$
unless $\veck\in\nbigj$ and $\veck\leq_{\seisuu^j}\veck(0)$.
Then, there exists
$U'\in\gbigu(Q,U)$ and
$t\in\lefttop{\ellsitabar}F
 \nbigc^{k-1}_{L^2,\leq j,\ttX(\Lambda)_{\geq 0}}
 (\nbigelambda,\DDlambda,h)(U')$
such that the following holds: 
\begin{itemize}
 \item $(s-\DDlambda t)^{(b)}_{\veca,\vecalpha,\veck}=0$
       $(b=1,3)$
       unless
       $\veck\in\nbigj$
       and $\veck\lneq_{\seisuu^j}\veck(1)$.
 \item There exists $\epsilon>0$ such that
       $|z_j|^{-\epsilon}
       (\DDlambda t)^{(b)}_{\veca',\vecalpha',\veck'}$
       are $L^2$
       for $b=1,3$
       and for $(\veca',\vecalpha')\neq (\veca,\vecalpha)$.
\end{itemize}
\end{prop}
\pf
By using the expression
\[
 s^{(3)}_{\veca,\vecalpha,\veck(0)}
 =\sum_{J,K}
 \sum_{(a_j(w_i),\alpha_j(w_i),\veck(w_i))=(a,\alpha,\veck(0))}
 s^{(3)}_{J,K,i}\omega_{J,K}w_i,
\]
we set
\[
 u^{(2)}_{\veca,\vecalpha,\veck(0)}
=\sum_{J,K}
 \sum_{(a_j(w_i),\alpha_j(w_i),\veck(w_i))=(a,\alpha,\veck(0))}
 \nbiga(s^{(3)}_{J,K,i})\omega_{J,K}w_i.
\]
By the construction, we have
\[
 d\zbar_j\del_{\zbar_j}
 \bigl(
 (dz_j/z_j)\cdot u^{(2)}_{\veca,\vecalpha,\veck(0)}
 \bigr)
 =(d\zbar_j/\zbar_j)(dz_j/z_j)
 s^{(3)}_{\veca,\vecalpha,\veck(0)}.
\]
\begin{lem}
 $\DD^{\prime\lambda}_{\neq j}
 \bigl((dz_j/z_j)\cdot u^{(2)}_{\veca,\vecalpha,\veck(0)}\bigr)$
 is $L^2$.
As a result, we may assume 
$s^{(3)}_{\veca,\vecalpha,\veck(0)}=0$.
\end{lem}
\pf
Because the section
\begin{multline}
(d\zbar_j/\zbar_j)(dz_j/z_j)
\Bigl(
 \DD^{\prime\lambda}_{\neq j}
 s^{(3)}_{\veca,\vecalpha,\veck(0)}
 +\zbar_j\del_{\zbar_j}
 \bigl(
 -\lambda s^{(1)}_{\veca,\vecalpha,\veck(0)}
+\nbigb_0(s^{(2)}_{\veca,\vecalpha,\veck(0)})
 \bigr)
 \Bigr)
 =\\
 (d\zbar_j/\zbar_j)(dz_j/z_j)
\Bigl(
 \DD^{\prime\lambda}_{\neq j}
 s^{(3)}_{\veca,\vecalpha,\veck(0)}
 -\lambda z_j\del_{z_j}
 s^{(1)}_{\veca,\vecalpha,\veck(0)}
 +
 \zbar_j\del_{\zbar_j}\nbigb_0(s^{(2)}_{\veca,\vecalpha,\veck(0)})
 \Bigr)
\end{multline}
is $L^2$,
the claim follows from Proposition \ref{prop;22.1.25.21}.
\hfill\qed

\vspace{.1in}
In the following, we assume
$s^{(3)}_{\veca,\vecalpha,\veck(0)}=0$.
To simplify the description, we set
\[
\nbigf_{\neq j,\veca,\vecalpha,\veck(1),\veck(0)}:=
G_{\neq j,(\veca,\vecalpha,\veck(1)),(\veca,\vecalpha,\veck(1))},
\]
\[
\nbigf_{j,\veca,\vecalpha,\veck(1),\veck(0)}:=
G_{j,(\veca,\vecalpha,\veck(1)),(\veca,\vecalpha,\veck(0))}.
\]
We may assume that
$\nbigf_{j,\veca,\vecalpha,\veck(1),\veck(0)}$ is an isomorphism.
We note that
$|\nbigf_{j,\veca,\vecalpha,\veck(1),\veck(0)}|_h=
O\bigl((-\log|z_j|)^{-1}\bigr)$
and
$|\nbigf^{-1}_{j,\veca,\vecalpha,\veck(1),\veck(0)}|_h=
O\bigl(-\log|z_j|\bigr)$.
We set
\[
 t^{(3,1)}_{\veca,\vecalpha,\veck(0)}
 =-(d\zbar_j/\zbar_j)
 \nbigf_{j,\veca,\vecalpha,\veck(1),\veck(0)}^{-1}
 \bigl(
 s^{(3)}_{\veca,\vecalpha,\veck(1)}
 \bigr),
\]
\[
 t^{(3,2)}_{\veca,\vecalpha,\veck(0)}
 =-\lambda (dz_j/z_j)
  \nbigf_{j,\veca,\vecalpha,\veck(1),\veck(0)}^{-1}
  (s^{(3)}_{\veca,\vecalpha,\veck(1)}),
\]
\[
 x^{(0)}_{\veck(0)}:=
 -\lambda\nbigf^{-1}_{j,\veca,\vecalpha,\veck(1),\veck(0)}
 \bigl(
 s^{(1)}_{\veca,\vecalpha,\veck(1)}
 \bigr)
+\nbigf^{-1}_{j,\veca,\vecalpha,\veck(1),\veck(0)}
 \nbigb_0(s^{(2)}_{\veca,\vecalpha,\veck(1)}).
\]

\begin{lem}
\label{lem;22.2.23.33}
$t^{(3,1)}_{\veca,\vecalpha,\veck(0)}
+t^{(3,2)}_{\veca,\vecalpha,\veck(0)}
+x^{(0)}_{\veck(0)}$
is a section of
 $\lefttop{\ellsitabar}
 F\nbigc^{k-1}_{L^2,\leq j,\ttX(\Lambda)_{\geq 0}}
 (\nbigelambda,\DDlambda,h)(U)$.
Moreover, there exists $\epsilon>0$ 
such that the following sections 
are $L^2$:
\[
 |z_j|^{-\epsilon}
 \Bigl(
 s-\DDlambda
 \bigl(
 t^{(3,1)}_{\veca,\vecalpha,\veck(0)}
 +t^{(3,2)}_{\veca,\vecalpha,\veck(0)}
 +x^{(0)}_{\veck(0)}
 \bigr)
 \Bigr)^{(3)}_{\veca,\vecalpha,\veck}
 \quad
 \mbox{(unless
 $\veck\lneq_{\seisuu^{j}}\veck(1)$)}
\]
\[
  |z_j|^{-\epsilon}
 \Bigl(
 s-\DDlambda
 \bigl(
 t^{(3,1)}_{\veca,\vecalpha,\veck(0)}
 +t^{(3,2)}_{\veca,\vecalpha,\veck(0)}
 +x^{(0)}_{\veck(0)}
 \bigr)
 \Bigr)^{(1)}_{\veca,\vecalpha,\veck}
\quad
  \mbox{(unless
 $\veck\leq_{\seisuu^{j}}\veck(1)$)}.
\]
\[
 |z_j|^{-\epsilon}
 \Bigl(
 \DDlambda
 \bigl(
 t^{(3,1)}_{\veca,\vecalpha,\veck(0)}
 +t^{(3,2)}_{\veca,\vecalpha,\veck(0)}
 +x^{(0)}_{\veck(0)}
 \bigr)
 \Bigr)^{(b)}_{\veca',\vecalpha',\veck'}
\quad
(b=1,3,\,\,\,(\veca',\vecalpha')\neq (\veca,\vecalpha)). 
\]
\end{lem}
\pf
There exists $\epsilon>0$ such that
the following sections are $L^2$:
\begin{equation}
 |z_j|^{-\epsilon}\cdot (d\zbar_j/\zbar_j)
 \cdot
 \Bigl(
 \zbar_j\del_{\zbar_j}
 \nbigb_0(s^{(0)}_{\veca,\vecalpha,\veck(0)})
 -\DD^{\prime\lambda}_{\neq j}
 s^{(1)}_{\veca,\vecalpha,\veck(0)}
 \Bigr),
\end{equation}
\begin{equation}
\label{eq;22.2.23.20}
 |z_j|^{-\epsilon}\cdot
 (d\zbar_j/\zbar_j)
 \cdot
 \Bigl(
 \zbar_j\del_{\zbar_j}
 \nbigb_0(s^{(0)}_{\veca,\vecalpha,\veck(1)})
 -\DD^{\prime\lambda}_{\neq j}
 s^{(1)}_{\veca,\vecalpha,\veck(1)}
 -\nbigf_{\neq j,\veca,\vecalpha,\veck(1),\veck(0)}\bigl(
 s^{(1)}_{\veca,\vecalpha,\veck(0)}
 \bigr)
 \Bigr),
\end{equation}
\begin{equation}
 |z_j|^{-\epsilon}\cdot
  (d\zbar_j/\zbar_j)
  (dz_j/z_j)
  \Bigl(
  -\lambda z_j\del_{z_j}
  s^{(1)}_{\veca,\vecalpha,\veck(0)}
  +\zbar_j\del_{\zbar_j}
  \nbigb_0(s^{(2)}_{\veca,\vecalpha,\veck(0)})
  \Bigr),
\end{equation}
\begin{equation}
\label{eq;22.2.23.4}
 |z_j|^{-\epsilon} (d\zbar_j/\zbar_j)
  (dz_j/z_j)
  \Bigl[
  \DD^{\prime\lambda}_{\neq j}
  s^{(3)}_{\veca,\vecalpha,\veck(1)}
 -\lambda z_j\del_{z_j}
  s^{(1)}_{\veca,\vecalpha,\veck(1)}
  -\nbigf_{j,\veca,\vecalpha,\veck(1),\veck(0)}\bigl(
  s^{(1)}_{\veca,\vecalpha,\veck(0)}
  \bigr)
 +\zbar_j\del_{\zbar_j}
 \nbigb_0(s^{(2)}_{\veca,\vecalpha,\veck(1)})
 \Bigr].
\end{equation}
We also note that
the following sections are $L^2$:
\begin{equation}
\label{eq;22.2.23.21}
 \DD^{\prime\lambda}_{\neq j}
 \nbigb_0(s^{(0)}_{\veca,\vecalpha,\veck(0)}),
 \quad
 \DD^{\prime\lambda}_{\neq j}
 \nbigb_0(s^{(0)}_{\veca,\vecalpha,\veck(1)}),
\end{equation}
\begin{equation}
\label{eq;22.2.23.22}
 \frac{dz_j}{z_j}
 \Bigl(
 \lambda z_j\del_{z_j}
  \nbigb_0(s^{(0)}_{\veca,\vecalpha,\veck(0)})
  -\DD^{\prime\lambda}_{\neq j}
  \nbigb_0(s^{(2)}_{\veca,\vecalpha,\veck(0)})
  \Bigr),
\quad\quad
 \frac{dz_j}{z_j}
 \Bigl(
 \lambda z_j\del_{z_j}
  \nbigb_0(s^{(0)}_{\veca,\vecalpha,\veck(1)})
  -\DD^{\prime\lambda}_{\neq j}
  \nbigb_0(s^{(2)}_{\veca,\vecalpha,\veck(1)})
 \Bigr).
\end{equation}

\begin{lem}
\label{lem;22.2.23.30}
There exists $\epsilon_1>0$
such that following section is $L^2$:
\[
 |z_j|^{-\epsilon_1}
 \Bigl(
 (d\zbar_j/\zbar_j)(dz_j/z_j)
 s^{(3)}_{\veca,\vecalpha,\veck(1)}
 -\bigl(
 dz_j\DDlambda_{z_j}
 t^{(3,1)}_{\veca,\vecalpha,\veck(0)}
 +d\zbar_j\del_{\zbar_j}t^{(3,2)}_{\veca,\vecalpha,\veck(0)}
\bigr)
 \Bigr).
\]
\end{lem}
\pf
By the construction,
we obtain
\[
 (dz_j/z_j)(\lambda z_j\del_{z_j}+\nbigf_{j,\veca,\vecalpha,\veck(1),\veck(0)})
 t^{(3,1)}_{\veca,\vecalpha,\veck(0)}
 +d\zbar_j\del_{\zbar_j}t^{(3,2)}_{\veca,\vecalpha,\veck(0)}
 =(d\zbar_j/\zbar_j)(dz_j/z_j)
 s^{(3)}_{\veca,\vecalpha,\veck(1)}.
\]
Then, the claim of Lemma \ref{lem;22.2.23.30} follows from
Lemma \ref{lem;22.2.23.3}.
\hfill\qed

\begin{lem}
\label{lem;22.2.23.5}
There exists $\epsilon_2>0$ such that
the following section is $L^2$:
\[
 |z_j|^{-\epsilon_2}
 \left(
 \frac{d\zbar_j}{\zbar_j}
 s^{(1)}_{\veca,\vecalpha,\veck(0)}
 -\bigl(
 d\zbar_j\del_{\zbar_j}x^{(0)}_{\veck(0)}
+\DD^{\prime\lambda}_{\neq j}t^{(3,1)}_{\veca,\vecalpha,\veck(0)}
 \bigr)
 \right).
\]
\end{lem}
\pf
We have
\[
  d\zbar_j\del_{\zbar_j}x^{(0)}_{\veck(0)}
= \frac{d\zbar_j}{\zbar_j}
 \nbigf^{-1}_{j,\veca,\vecalpha,\veck(1),\veck(0)}
 \Bigl(
 \zbar_j\del_{\zbar_j}
 \bigl(
 -\lambda s^{(1)}_{\veca,\vecalpha,\veck(1)}
 +\nbigb_0s^{(2)}_{\veca,\vecalpha,\veck(1)}
 \bigr)
 \Bigr).
\]
Because
$\DDlambda_{\neq j}\bigl(
\Res_j(\DDlambda)\bigr)=0$,
there exists $\epsilon_2'>0$
such that the following holds
on $Z(\vecC)$
with respect to $h$ and the Poincar\'e metric:
\[
 \bigl[
 \DD^{\prime\lambda}_{\neq j},
 \nbigf_{j,\veca,\vecalpha,\veck(1),\veck(0)}
 \bigr]
=O(|z_j|^{\epsilon_2'}).
\]
There exists $\epsilon_2''>0$
such that the following holds
on $Z(\vecC)$
with respect to $h$ and the Poincar\'e metric:
\[
 \bigl[
 \DD^{\prime\lambda}_{\neq j},
 \nbigf^{-1}_{j,\veca,\vecalpha,\veck(1),\veck(0)}
 \bigr]
=O(|z_j|^{\epsilon_2''}).
\]
There exists $\epsilon_2'''>0$
such that the following section is $L^2$:
\[
 |z_j|^{-\epsilon_2'''}
 \Bigl(
\DD^{\prime\lambda}_{\neq j}t^{(3,1)}_{\veck(0)}
-
 \frac{d\zbar_j}{\zbar_j}
 \nbigf^{-1}_{j,\veca,\vecalpha,\veck(1),\veck(0)}
 \Bigl(
 \DD^{\prime\lambda}_{\neq j}
 s^{(3)}_{\veca,\vecalpha,\veck(1)}
 \Bigr)
 \Bigr).
\]
Because (\ref{eq;22.2.23.4}) is $L^2$,
there exists $\epsilon_2''''>0$
such that the following section is $L^2$:
\begin{equation}
\label{eq;22.2.23.10}
 |z_j|^{-\epsilon_2''''}
  \frac{d\zbar_j}{\zbar_j}
 \Bigl(
 s^{(1)}_{\veca,\vecalpha,\veck(0)}
- \nbigf^{-1}_{j,\veca,\vecalpha,\veck(1),\veck(0)}
 \Bigl(
 \zbar_j\del_{\zbar_j}
 \bigl(
 -\lambda s^{(1)}_{\veca,\vecalpha,\veck(1)}
 +\nbigb_0s^{(2)}_{\veca,\vecalpha,\veck(1)}
 \bigr)
 +\DD^{\prime\lambda}_{\neq j}
 s^{(3)}_{\veca,\vecalpha,\veck(1)}
 \Bigr)
 \Bigr).
\end{equation}
Thus, we obtain the claim of Lemma \ref{lem;22.2.23.5}.
\hfill\qed

\begin{lem}
\label{lem;22.2.23.31}
There exists $\epsilon>0$ such that
the following estimate holds on $Z(\vecC)$:
\[
 z_j\del_{z_j}
 \nbigf_{j,\veca,\vecalpha,\veck(1),\veck(0)}
=O\bigl(|z_j|^{\epsilon_3'}\bigr).
\]
\end{lem}
\pf
Let $\vecw_{\veca,\vecalpha,\veck}$ denote the tuple of
$w_i$ such that
$(a_p(w_i),\alpha_p(w_i),k_p(w_i))
=(a_p,\alpha_p,k_p)$ $(p\leq j)$,
which is a frame of
$\nbigv_{\veca,\vecalpha,\veck}$.
Let $B_{j,\veca,\vecalpha,\veck(1),\veck(0)}$
be the matrix valued holomorphic function
determined by
\[
 \nbigf_{j,\veca,\vecalpha,\veck(1),\veck(0)}
 \vecw_{\veca,\vecalpha,\veck(0)}
 =
  \vecw_{\veca,\vecalpha,\veck(1)}
  B_{j,\veca,\vecalpha,\veck(1),\veck(0)}.
\]
For $j<m\leq \ell$,
we have
$(B_{j,\veca,\vecalpha,\veck(1),\veck(0)})_{p,q|H_m}=0$
unless the following holds.
\begin{itemize}
 \item $\alpha_m(w_p)=\alpha_m(w_q)$
       and
       $a_m(w_p)\leq a_m(w_q)$.
\end{itemize}
We have
$(z_j\del_{z_j}B_{j,\veca,\vecalpha,\veck(1),\veck(0)})_{p,q|H_m}=0$
unless the same condition holds.
We have
$(z_j\del_{z_j}B_{j,\veca,\vecalpha,\veck(1),\veck(0)})_{p,q|H_j}=0$.
Then, the claim of Lemma \ref{lem;22.2.23.31}
is easy to see.
\hfill\qed

\begin{lem}
\label{lem;22.2.23.11}
The section
$dz_j\DDlambda_{z_j}(x^{(0)}_{\veck(0)})
+\DDlambda_{\neq j}t^{(3,2)}_{\veck(0)}$
is $L^2$.
\end{lem}
\pf
As in the proof of Lemma \ref{lem;22.2.23.5},
there exists $\epsilon_3>0$
such that the following section is $L^2$:
\[
 |z_j|^{-\epsilon_3}
 \Bigl(
 \DD^{\prime\lambda}_{\neq j}
 t^{(3,2)}_{\veck(0)}
 -
  \frac{dz_j}{z_j}
 \lambda \nbigf^{-1}_{j,\veca,\vecalpha,\veck(1),\veck(0)}
 \Bigl(
 \DD^{\prime\lambda}_{\neq j}
 s^{(3)}_{\veca,\vecalpha,\veck(1)}
 \Bigr)
 \Bigr)
\]
There exists $\epsilon_3'>0$ such that
the following holds on $Z(\vecC)$:
\[
 \bigl[
 z_j\del_{z_j},
 \nbigf^{-1}_{j,\veca,\vecalpha,\veck(1),\veck(0)}
 \bigr]
=O\bigl(|z_j|^{\epsilon_3'}\bigr).
\]
There exists $\epsilon_3''>0$ such that
the following section is $L^2$:
\[
 |z_j|^{-\epsilon_3''}
 \Bigl(
 \frac{dz_j}{z_j}
 \lambda z_j\del_{z_j}x^{(0)}_{\veck(0)}
-\frac{dz_j}{z_j}
 \lambda \nbigf^{-1}_{j,\veca,\vecalpha,\veck(1),\veck(0)}
 \Bigl(
 z_j\del_{z_j}\bigl(
 -\lambda s^{(1)}_{\veca,\vecalpha,\veck(1)}
 +\nbigb_0(s^{(2)}_{\veca,\vecalpha,\veck(1)})
 \bigr)
 \Bigr)
 \Bigr).
\]
Because (\ref{eq;22.2.23.10}) is $L^2$,
and because
$(dz_j/z_j)s^{(1)}_{\veca,\vecalpha,\veck(0)}$ is $L^2$,
we obtain the claim of Lemma \ref{lem;22.2.23.11}.
\hfill\qed

\begin{lem}
\label{lem;22.2.23.32}
$\DD^{\lambda}_{\neq j}
 x^{(0)}_{\veck(0)}$ is $L^2$.
\end{lem}
\pf
There exists $\epsilon_5>0$ such that the following section is $L^2$:
\begin{equation}
|z_j|^{-\epsilon_5}
 \Bigl(
  \DD^{\prime\lambda}_{\neq j}
  x^{(0)}_{\veck(0)}
-
 \nbigf^{-1}_{j,\veca,\vecalpha,\veck(1),\veck(0)}
 \Bigl(
 -\lambda\DD^{\prime\lambda}_{\neq j} 
 s^{(1)}_{\veca,\vecalpha,\veck(1)}
 +\DD^{\prime\lambda}_{\neq j}
 \nbigb_0(s^{(2)}_{\veca,\vecalpha,\veck(1)})
 \Bigr)
\Bigr).
\end{equation}
Because (\ref{eq;22.2.23.20})
and (\ref{eq;22.2.23.22}) are $L^2$,
we can check that
\[
  \nbigf^{-1}_{j,\veca,\vecalpha,\veck(1),\veck(0)}
 \Bigl(
 -\lambda\DD^{\prime\lambda}_{\neq j} 
 s^{(1)}_{\veca,\vecalpha,\veck(1)}
 +\DD^{\prime\lambda}_{\neq j}
 \nbigb_0(s^{(2)}_{\veca,\vecalpha,\veck(1)})
 \Bigr)
\]
is $L^2$,
by using the relation
$\zbar_j\del_{\zbar_j}\nbigb_0(s^{(0)}_{\veca,\vecalpha,\veck(0)})
=z_j\del_{z_j}\nbigb_0(s^{(0)}_{\veca,\vecalpha,\veck(0)})$.
Thus, we obtain Lemma \ref{lem;22.2.23.32}.
\hfill\qed

\vspace{.1in}
We obtain the first claim of Lemma \ref{lem;22.2.23.33}
from Lemma \ref{lem;22.2.23.30},
Lemma \ref{lem;22.2.23.5},
Lemma \ref{lem;22.2.23.11},
and Lemma \ref{lem;22.2.23.32}.
We obtain the second claim of Lemma \ref{lem;22.2.23.33}
from Lemma \ref{lem;22.2.23.30} and
Lemma \ref{lem;22.2.23.5}.
Thus, the proof of Lemma \ref{lem;22.2.23.33}
is completed.
\hfill\qed

\begin{lem}
\label{lem;22.2.23.40}
 There exist $U'\in\gbigu(Q,U)$
 and a section $t\in
 \lefttop{\ellsitabar}F
 \nbigc^{k-1}_{L^2,\leq j,\ttX(\Lambda)_{\geq 0}}
 (\nbigelambda,\DDlambda,h)(U')$
such that the following holds.
\begin{itemize}
 \item
      $(s-\DDlambda t)_{\veca,\vecalpha,\veck}^{(3)}=0$
      unless $\veck\in\nbigj$
      and $\veck\lneq \veck(1)$.
 \item  $(s-\DDlambda t)_{\veca,\vecalpha,\veck}^{(1)}=0$
      unless $\veck\in\nbigj$
      and $\veck\leq \veck(1)$.
 \item There exists $\epsilon>0$
       such that
       $|z_j|^{-\epsilon}\bigl(
       \DDlambda t
       \bigr)^{(b)}_{\veca',\vecalpha',\veck'}$
       are $L^2$
       for $(\veca',\vecalpha')\neq(\veca,\vecalpha)$.
\end{itemize}
\end{lem}
\pf
It follows from Lemma \ref{lem;22.2.23.33},
Lemma \ref{lem;22.2.22.3},
Lemma \ref{lem;22.2.22.30}
and Lemma \ref{lem;22.2.22.31}.
\hfill\qed

\vspace{.1in}

By Lemma \ref{lem;22.2.23.40},
we may assume
$s^{(3)}_{\veca,\vecalpha,\veck(1)}=0$
and $s^{(1)}_{\veca,\vecalpha,\veck(0)}=0$.
By using the expression
\[
 s^{(1)}_{\veca,\vecalpha,\veck(1)}
 =\sum_{J,K}
 \sum_{(\veca(w_i),\vecalpha(w_i),\veck(w_i))=(\veca,\vecalpha,\veck(1))}
 s^{(1)}_{J,K,i}\cdot\omega_{J,K}\cdot w_i,
\]
We set
\[
 u^{(0)}_{\veca,\vecalpha,\veck(1)}
=\sum_{(\veca(w_i),\vecalpha(w_i),\veck(w_i))=(a,\alpha,\veck(1))}
 \nbiga(s^{(1)}_{J,K,i})\cdot \omega_{J,K}\cdot w_i.
\]
By the construction, we have
\[
 d\zbar_j\del_{\zbar_j}
 u^{(0)}_{\veca,\vecalpha,\veck(1)}
 =(d\zbar_j/\zbar_j)
 s^{(1)}_{\veca,\vecalpha,\veck(1)}.
\]
Proposition \ref{prop;22.2.23.42}
is reduced to the next lemma.
\begin{lem}
\label{lem;22.2.23.41}
$\DDlambda u^{(0)}_{\veca,\vecalpha,\veck(1)}$ is $L^2$.
\end{lem}
\pf
Because
$z_j\del_{z_j}u^{(0)}_{\veca,\vecalpha,\veck(1)}
=\zbar_j\del_{\zbar_j}u^{(0)}_{\veca,\vecalpha,\veck(1)}$,
we obtain that
$dz_j\del_{z_j}u^{(0)}_{0,0,\veck(1)}$ is $L^2$.
Because
\[
d\zbar_j\del_{\zbar_j}
\nbigb_0(s^{(0)}_{\veca,\vecalpha,\veck(1)})
+\DD^{\prime\lambda}_{\neq j}
\bigl(
(d\zbar_j/\zbar_j)
s^{(1)}_{\veca,\vecalpha,\veck(1)}
\bigr)
\]
is $L^2$,
we obtain that
$\DD^{\prime\lambda}_{\neq j}(u^{(0)}_{\veca,\vecalpha,\veck(1)})$
is $L^2$
by Proposition \ref{prop;22.1.25.21}.
Thus, we obtain Lemma \ref{lem;22.2.23.41},
and Proposition \ref{prop;22.2.23.42}.

\hfill\qed

\vspace{.1in}

As in \S\ref{subsection;22.2.23.1},
let $\leq'_{\seisuu^j}$ be a total order on $\seisuu^{\ell}$
which is a refinement of $\leq_{\seisuu^j}$.
Then,
we obtain the claim of Proposition \ref{prop;22.2.22.40}
in the case $\alpha=0$
from by using Proposition \ref{prop;22.2.23.42}.
and an easy induction
on $(\nbigj,\leq'_{\seisuu^j})$.
\hfill\qed

\subsubsection{Proof of Proposition \ref{prop;22.2.23.43}}
\label{subsection;22.4.27.40}

Let $\nbigk$ denote the set of
$(\veca,\vecalpha)\in\openclosed{-1}{0}^j\times\cnum^j$
such that
$\nbigv_{\veca,\vecalpha}\neq 0$.
Let $\leq_{\nbigk}$ be a total order on $\nbigk$.

Let $U\in\gbigu(Q)$.
Let
$s\in \lefttop{\ellsitabar}F
\nbigc^{k}_{L^2,\leq j,\ttX(\Lambda)_{\geq 0}}
(\nbigelambda,\DDlambda,h)(U)$
such that $(\DDlambda s)^{(b)}=0$ $(b=1,3)$.
We can prove the following claim
by using an induction on $(\nbigk,\leq_{\nbigk})$
with
Lemma \ref{lem;22.2.22.3},
Lemma \ref{lem;22.2.22.30},
Lemma \ref{lem;22.2.22.31}
and Proposition \ref{prop;22.2.22.40}.
\begin{description}
 \item[$\boldsymbol{P(\veca,\vecalpha)}$:]
	    There exist $U'\in\gbigu(Q,U)$
	    and $u\in\lefttop{\ellsitabar}F
	    \nbigc^{k-1}_{L^2,\leq j,\ttX(\Lambda)_{\geq 0}}
	    (\nbigelambda,\DDlambda,h)(U')$
	    such that the following holds: 
\begin{itemize}
 \item $(s-\DDlambda u)^{(b)}_{\veca',\vecalpha',\veck'}=0$
       for $b=1,3$
       unless $(\veca',\vecalpha')\geq_{\nbigk}(\veca,\vecalpha)$.
\end{itemize}
\end{description}
Hence, there exist
$U'\in\gbigu(Q,U)$
and $u\in\lefttop{\ellsitabar}F
\nbigc^{k-1}_{L^2,\leq j,\ttX(\Lambda)_{\geq 0}}
(\nbigelambda,\DDlambda,h)(U')$
such that
\[
s-\DDlambda u\in
\lefttop{\ellsitabar}F
\nbigc^{k}_{L^2,\leq j-1,\ttX(\Lambda)_{\geq 0}}
(\nbigelambda,\DDlambda,h)(U').
\]

Let $s\in \lefttop{\ellsitabar}F
\nbigc^{k}_{L^2,\leq j,\ttX(\Lambda)_{\geq 0}}
(\nbigelambda,\DDlambda,h)(U)$
such that $\DDlambda s=0$.
We take $U'$ and $u$ as above for $s$,
and we set $\stilde:=s-\DDlambda u$.
Then, we obtain $\DDlambda\stilde=0$
and $\stilde\in\lefttop{\ellsitabar}F
\nbigc^{k}_{L^2,\leq j-1,\ttX(\Lambda)_{\geq 0}}
(\nbigelambda,\DDlambda,h)(U')$.

Let
$s\in \lefttop{\ellsitabar}F
\nbigc^{k}_{L^2,\leq j-1,\ttX(\Lambda)_{\geq 0}}
(\nbigelambda,\DDlambda,h)(U)$
and
$t\in \lefttop{\ellsitabar}F
\nbigc^{k-1}_{L^2,\leq j,\ttX(\Lambda)_{\geq 0}}
(\nbigelambda,\DDlambda,h)(U)$
such that $s=\DDlambda t$.
There exist
$U'\in\gbigu(Q,U)$
and 
$u\in\lefttop{\ellsitabar}F
\nbigc^{k-2}_{L^2,\leq j,\ttX(\Lambda)_{\geq 0}}
(\nbigelambda,\DDlambda,h)(U')$
such that
$(t-\DDlambda u)^{(b)}=0$ $(b=1,3)$.
Then, $s=\DDlambda(t-\DDlambda u)$.
Thus, 
Proposition \ref{prop;22.2.23.43} is proved.
\hfill\qed

\subsubsection{Family version}

Let
$\nbigc^{\bullet}_{L^2,\leq j,\ttX(\Lambda)_{\geq 0}}
(\nbige,\DD,h)
\subset 
\nbigc^{\bullet}_{L^2,\ttX(\Lambda)_{\geq 0}}(\nbige,\DD,h)$
denote the subsheaf of local sections $s$
satisfying the following condition.
\begin{itemize}
 \item $s$ is a local section of
       $(j_{\nbigx\setminus\nbigh,\ttX(\Lambda)_{\geq 0}})_{\ast}
       \Bigl(
       \nbige\otimes
       p_{\lambda}^{-1}\bigl(
 \Omega^{\bullet,0}
 \otimes
       \Omega^{0,\bullet}_{\leq j}
       \bigr)
       \Bigr)$,
       and satisfies $\del_{\zbar_i}s=0$ for $i>j$.
\end{itemize}
We obtain the subcomplex
$\nbigc^{\bullet}_{L^2,\leq j,\ttX(\Lambda)_{\leq 0}}
(\nbige,\DD,h)$
of 
$\nbigc^{\bullet}_{L^2,\ttX(\Lambda)_{\geq 0}}(\nbige,\DD,h)$.
By the construction,
$\nbigc^{\bullet}_{L^2,\leq n,\ttX(\Lambda)_{\geq 0}}
(\nbige,\DD,h)=
\nbigc^{\bullet}_{L^2,\ttX(\Lambda)_{\geq 0}}(\nbige,\DD,h)$.

For an $(S^1)^{\ell}$-invariant open subset
$U$ of $\nbigxzero(\ttX(\Lambda)_{\geq 0})$,
we define
\[
 \lefttop{J}F\nbigc^{\bullet}_{L^2,\leq j,\ttX(\Lambda)_{\geq 0}}
(\nbige,\DD,h)(U)
\subset
\nbigc^{\bullet}_{L^2,\leq j,\ttX(\Lambda)_{\geq 0}}
(\nbige,\DD,h)(U)
\]
as in \S\ref{subsection;22.1.27.1}.
In particular,
we obtain the subcomplex
\[
 \lefttop{\ellsitabar}F
 \nbigc^{\bullet}_{L^2,\leq j,\ttX(\Lambda)_{\geq 0}}(\nbige,\DD,h)(U)
\subset
 \nbigc^{\bullet}_{L^2,\leq j,\ttX(\Lambda)_{\geq 0}}(\nbige,\DD,h)(U).
\]
Let $Q$ be a point of
$p_1^{-1}(O)\subset X(\ttX(\Lambda)_{\geq 0})$.
Let $(\gbigu(\lambda,Q),\prec)$ be the directed set of
$(S^1)^{\ell}$-invariant open neighbourhoods of $(\lambda,Q)$
in $\nbigxzero$,
where $U_1\prec U_2$ is defined to be $U_1\supset U_2$.
For any $U\in\gbigu(Q)$,
we set
$\gbigu((\lambda,Q),U):=\bigl\{
 U'\in\gbigu(\lambda,Q)\,\big|\,
 U\prec U'
\bigr\}$.
We set
\[
  \lefttop{\ellsitabar}F
  \nbigc^{\bullet}_{L^2,\leq j,\ttX(\Lambda)_{\geq 0}}
  (\nbige,\DD,h)_Q:=
  \varinjlim_{U\in\gbigu(\lambda,Q)}
    \lefttop{\ellsitabar}F
    \nbigc^{\bullet}_{L^2,\leq j,\ttX(\Lambda)_{\geq 0}}
    (\nbige,\DD,h)(U).
\]
The following proposition is similar to
Proposition \ref{prop;22.2.23.43}.
\begin{prop}
The natural morphism
\begin{equation}
 \lefttop{\ellsitabar}F
 \nbigc^{\bullet}_{L^2,\leq j,\ttX(\Lambda)_{\geq 0}}
 (\nbige,\DD,h)_{(\lambda,Q)}
 \lrarr
 \lefttop{\ellsitabar}F
 \nbigc^{\bullet}_{L^2,\leq n,\ttX(\Lambda)_{\geq 0}}
 (\nbige,\DD,h)_{(\lambda,Q)}
\end{equation}
is a quasi-isomorphism.
\hfill\qed
\end{prop}

\subsection{Proof of Theorem \ref{thm;22.2.18.11}
and Theorem \ref{thm;21.12.21.11}}
\label{subsection;22.3.17.51}

\subsubsection{$L^2$-holomorphic sections}
\label{subsection;22.2.23.53}

Let $0=p(0)<p(1)<\ldots<p(m)=\ell$ be an increasing sequence
of integers.
We set $K_i=\{1,\ldots,p(i)\}$.
We obtain $\vecK=(K_1,\ldots,K_m)\in\nbigs(\ellsitabar)$.

Let $\vecx=(x_1,\ldots,x_{\ell})\in\real^{\ell-1}_{\geq 0}\times\{0\}$
such that $\{p(1),\ldots,p(m)\}=\{1\leq i\leq \ell\,|\,x_i=0\}$.
Let $\vecepsilon=(\epsilon_1,\ldots,\epsilon_{\ell})
\in\real_{>0}^{n}$.
Let $U(\vecx,\vecepsilon)$
denote the open subset of $X\setminus H$
determined by the following condition:
\[
 \Bigl|
 \frac{-\log|z_{i+1}|}{-\log|z_{i}|}
 -x_i
 \Bigr|<\epsilon_i\,\,\,(i=1,\ldots,\ell-1),
 \quad\quad
 \frac{1}{-\log|z_{\ell}|}<\epsilon_{\ell},
 \quad\quad
 |z_i|<\epsilon_i\,\,\,(i=\ell+1,\ldots,n).
\]
It is an $(S^1)^{\ell}$-invariant open subset of
$X\setminus H$.

Let $s\in
\lefttop{\ellsitabar}F
\nbigc^k_{L^2,\leq 0,\ttX(\Lambda)_{\geq 0}}
(\nbigelambda,\DDlambda,h)(\overline{U(\vecx,\vecepsilon)})$.
The condition implies that $s$ is an $L^2$-holomorphic section of
$\nbigelambda\otimes\Omega^{k,0}$ on $U(\vecx,\vecepsilon)$.
We have the expression
\[
s=\sum_{J}
 \sum_i
 \sum_{\vecm\in\seisuu_{\geq 0}^{\ell}}
 \gbigf_{\ellsitabar,\vecm}
 (s_{J,\emptyset,i})
 e^{\sqrt{-1}\vecm\vectheta}\cdot
 \omega_{J,\emptyset}\cdot
 w_i.
\]
The functions
$s_{J,\emptyset,i}
=\sum_{\vecm}
 \gbigf_{\ellsitabar,\vecm}(s_{J,\emptyset,i})
e^{\sqrt{-1}\vecm\vectheta}$
are holomorphic on $U(\vecx,\vecepsilon)$.
Hence, there exists a holomorphic function
$\stilde_{J,i,\vecm}$ on $\prod_{j=\ell+1}^n\{|z_j|<\epsilon_j\}$
such that
\[
\gbigf_{\ellsitabar,\vecm}
(s_{J,\emptyset,i})\cdot
e^{\sqrt{-1}\vecm\vectheta}
=\stilde_{J,i,\vecm}(z_{\ell+1},\ldots,z_{n})
 \prod_{j=1}^{\ell}z_j^{m_j}.
\]
\begin{lem}
There exists a neighbourhood $N$ of
$H_{\ellsitabar}$ in $X$
such that
(i) 
$s_{J,\emptyset,i}$ induces a holomorphic function
 on $N$,
 (ii)
$U(\vecx,\vecepsilon)$
is contained in $N$. 
\end{lem}
\pf
We take any $\vect=(t_1,\ldots,t_{\ell})\in\real^{\ell}_{>0}$
satisfying
\[
 \left|\frac{-\log t_{i+1}}{-\log t_i}-x_i\right|
 <\epsilon_i\,\,\,(i=1,\ldots,\ell-1),
 \quad\quad
 \frac{1}{-\log t_{\ell}}<\epsilon_{\ell}.
\]
There exists $\kappa>0$ such that
$\prod_{i=1}^{\ell}
\{t_i-\kappa<|z_i|< t_i+\kappa\}
\times
\prod_{i=\ell+1}^n\{|z_i|<\epsilon_i\}$
is a relatively compact open subset of
$U(\vecx,\vecepsilon)$.
By using Cauchy's integral formula,
we obtain that 
$s_{J,\emptyset,i}$ induces
a holomorphic function on
$A(\vect)=\prod_{i=1}^{\ell}\{|z_i|<t_i\}
\times
\prod_{i=\ell+1}^n\{|z_i|<\epsilon_i\}$,
and it is convergent on the multi-disc.
By varying $\vect$,
we obtain the claim of the lemma.
\hfill\qed

\vspace{.1in}

Let $N$ be a neighbourhood of $H_{\ellsitabar}$ in $X$
of the form
$\prod_{i=1}^{n}\{|z_i|<\delta_i\}$
for some $(\delta_1,\ldots,\delta_n)\in\real_{>0}^n$.
Let $s$ be a holomorphic section of
$\nbigelambda\otimes\Omega^{k,0}$ on $N$.
Suppose that $U(\vecx,\vecepsilon)$
is a relatively compact open subset of $N$.
By the norm estimate as in Proposition \ref{prop;22.2.23.50},
$s$ is $L^2$ on $U(\vecx,\vecepsilon)$
with respect to $h$ and $g_{X\setminus H}$
if and only if the following integral is convergent:
\begin{equation}
\label{eq;22.1.27.20}
 \sum_{J,i,\vecm}
 \int_{U(\vecx,\vecepsilon)}
 \bigl|\stilde_{J,i,\vecm}(z_{\ell+1},\ldots,z_n)\bigr|^2
 \prod_{p=1}^{\ell}
 |z_p|^{-2a_p(w_i)+2m_p}
 (-\log|z_p|)^{k_p(w_i)-k_{p-1}(w_i)}
 \prod_{\substack{p\in J\\ p\leq \ell }}(-\log|z_p|)^{2}
 \dvol_{g_{X\setminus H}}.
\end{equation}
The integral (\ref{eq;22.1.27.20})
is convergent if and only if
the following is convergent for each $J$ and $i$:
\begin{multline}
\label{eq;22.1.28.1}
 \int_{U(\vecx,\vecepsilon)}
 |s_{J,\emptyset,i}\omega_{J,\emptyset}w_i|^2\dvol_{g_{X\setminus H}}
 =\\
 \sum_{\vecm}
\int_{U(\vecx,\vecepsilon)}
 \bigl|\stilde_{J,i,\vecm}(z_{\ell+1},\ldots,z_n)\bigr|^2
 \prod_{p=1}^{\ell}
 |z_p|^{-2a_p(w_i)+2m_p}
 (-\log|z_p|)^{k_p(w_i)-k_{p-1}(w_i)}
 \prod_{\substack{p\in J\\ p\leq \ell}}(-\log|z_p|)^{2}
 \dvol_{g_{X\setminus H}}.
\end{multline}

Let $\nbigt(s_{J,\emptyset,i})=
\bigl\{\vecm\in\seisuu^{\ell}_{\geq 0}\,\big|\,
\stilde_{J,i,\vecm}\neq 0
\bigr\}$.
We consider the partial order $\leq$ on $\seisuu^{\ell}$
defined by $(p_i)\leq (q_i)\Longleftrightarrow
p_i\leq q_i\,\,(\forall i)$.
Let $\min\nbigt(s_{J,\emptyset,i})$ denote the set of
the minimal elements of $\nbigt(s_{J,\emptyset,i})$,
which is finite.
Because $-1<a_p(w_i)\leq 0$,
$q_{K_s}(\vecm-\veca(w_i))=0$ holds
if and only if
$a_p(w_i)=m_p=0$ for $p\in K_s$.
We can check the following lemma
by using Lemma \ref{lem;22.1.27.11} below.
\begin{lem}
\label{lem;22.1.28.3}
The integral {\rm(\ref{eq;22.1.28.1})} is convergent
if and only if
for each $\vecm\in\min\nbigt(s_{J,\emptyset,i})$
one of the following holds: 
\begin{itemize}
 \item $q_{K_1}(\vecm-\veca(w_i))\neq 0$.
 \item There exists $1\leq j_0\leq m-1$ such that
       $q_{K_{j_0}}(\vecm-\veca(w_i))=0$,
       $q_{K_{j_0+1}}(\vecm-\veca(w_i))\neq 0$,
       and $k_{p(s)}(w_i)\leq p(s)-2|J\cap K_s|-1$
       for $s=1,\ldots,j_0$.
 \item $\vecm-\veca(w_i)=0$
       and $k_{p(s)}(w_i)\leq p(s)-2|J\cap K_s|-1$
       for any $s$.
\hfill\qed
\end{itemize}
\end{lem}

We set
\[
 s_{\veca,\vecalpha}
 =\sum_{J\subset\nbar}\sum_{(\veca(w_i),\vecalpha(w_i))=(\veca,\vecalpha)}
 s_{J,\emptyset,i}\omega_{J,\emptyset}w_i.
\]
The following lemma is easy to see.
\begin{lem}
$s$ and $\DDlambda s$ are $L^2$
if and only if
$s_{\veca,\vecalpha}$
and 
$\DDlambda s_{\veca,\vecalpha}$
are $L^2$ for any $(\veca,\vecalpha)$.
\hfill\qed
\end{lem}

\subsubsection{Quasi-isomorphisms}

Let $Q$ be any point of $X(\ttX(\ellsitabar)_{\geq 0})$.
By Lemma \ref{lem;22.1.28.3},
for any $U\in\gbigu(Q)$,
there exists $U'\in\gbigu(Q,U)$
such that 
\begin{equation}
\label{eq;22.2.23.50}
 \nbigw
 \nbigc^{\bullet}_{\tw}(\gbige[\ast H],\lambda)(U)
 \subset
 \lefttop{\ellsitabar}F
 \nbigc^{\bullet}_{L^2,\leq 0,\ttX(\Lambda)_{\geq 0}}
 (\nbigelambda,\DDlambda,h)(U').
\end{equation} 
By taking the inductive limit,
we obtain
\begin{equation}
\label{eq;22.2.23.51}
 \nbigw \nbigc^{\bullet}_{\tw}(\gbige[\ast H],\lambda)_Q
 \lrarr
 \lefttop{\ellsitabar}F
 \nbigc^{\bullet}_{L^2,\leq 0,\ttX(\Lambda)_{\geq 0}}
 (\nbigelambda,\DDlambda,h)_Q.
\end{equation}

\begin{prop}
\label{prop;22.2.23.52}
The natural inclusion
{\rm(\ref{eq;22.2.23.51})} is a quasi-isomorphism.
\end{prop}
\pf
It is enough to consider the case where
$Q\in X(\ttO(\vecJ)_{\geq 0})$
for $\vecJ$ as in \S\ref{subsection;22.2.23.53}.
Let $U=U(\vecx,\vecepsilon)$
be as in \S\ref{subsection;22.2.23.53}.
Let $s\in \lefttop{\ellsitabar}F
\nbigc^{k}_{L^2,\leq 0,\ttX(\Lambda)_{\geq 0}}
(\nbigelambda,\DDlambda,h)(U)$
such that $\DDlambda s=0$.
For $(\veca,\vecalpha)\in\openclosed{-1}{0}^{\ell}\times\cnum^{\ell}$,
and for $1\leq j\leq \ell$,
we have the expansion of $s_{\veca,\vecalpha}$
with respect to the variable $z_j$:
\[
 s_{\veca,\vecalpha}
 =\sum_{J\subset\nbar}
 \sum_{(\veca(w_i),\vecalpha(w_i))=(\veca,\vecalpha)}
 \sum_{m=0}^{\infty}
 s^{\langle j\rangle}_{\veca,\vecalpha,J,i;m}\cdot z_j^m\cdot
 \omega_{J,\emptyset}\cdot w_i.
\]
We set
\[
 s_{\veca,\vecalpha;0}^{\langle j\rangle}:=
 \sum_{J\subset\nbar}
 \sum_{(\veca(w_i),\vecalpha(w_i))=(\veca,\vecalpha)}
 s^{\langle j\rangle}_{\veca,\vecalpha,J,i;0}\cdot 
 \omega_{J,\emptyset}\cdot w_i.
\]

We shall prove the following claim
by using a descending induction on $j$:
\begin{description}
 \item[$\boldsymbol{P(j)}$:] There exist $U'\in\gbigu(Q,U)$
       and $t\in \lefttop{\ellsitabar}F
       \nbigc^{k-1}_{L^2,\leq 0,\ttX(\Lambda)_{\geq 0}}
       (\nbigelambda,\DDlambda,h)(U')$
   such that
    $(s-\DDlambda t)^{\langle k\rangle}_{\veca,\vecalpha;0}=0$
    holds	    
    for any $j\leq k\leq \ell$
    and $(\veca,\vecalpha)$
    satisfying $a_k=0$ and $\alpha_k\neq 0$.
\end{description}
Let us prove $\boldsymbol{P(j)}$
by assuming $\boldsymbol{P(j+1)}$.
We may assume
$s^{\langle k\rangle}_{\veca,\vecalpha;0}=0$
for any $j+1\leq k\leq \ell$
and any $(\veca,\vecalpha)$
satisfying $a_k=0$ and $\alpha_k\neq 0$.
Let $\alpha\neq 0$.
We set
\[
 s_{0,\alpha}:=
 \sum_{\substack{(\veca,\vecalpha)\in \openclosed{-1}{0}^j\times\cnum^j\\
  a_j=0,\alpha_j=\alpha}}
  s_{\veca,\vecalpha}.
\]
We decompose
\[
 s_{0,\alpha}
 =s^{(0)}_{0,\alpha}
 +\frac{dz_j}{z_j}s^{(2)}_{0,\alpha}
\]
as in \S\ref{subsection;22.2.23.100}.
Because
$\frac{dz_j}{z_j}s^{(2)}_{0,\alpha}$ is $L^2$ on $U$
with respect to $h$ and $g_{X\setminus H}$,
$s^{(2)}_{0,\alpha}$ is $L^2$ on $U$
with respect to $h$ and $g_{X\setminus H}$.

\begin{lem}
\label{lem;22.2.24.10}
$\DDlambda s^{(2)}_{0,\alpha}$ is $L^2$ on $U$
 with respect to $h$ and $g_{X\setminus H}$.
\end{lem}
\pf
Note that $(dz_j/z_j)s^{(2)}_{0,\alpha}$ is $L^2$ on $U$
with respect to $h$ and $g_{X\setminus H}$.
Hence, we obtain that
$dz_j\del_{z_j}(s^{(2)}_{0,\alpha})
=(dz_j/z_j)\bigl(
z_j\del_{z_j}
s^{(2)}_{0,\alpha}\bigr)$
is also $L^2$ on $U$ with respect to $h$ and $g_{X\setminus H}$
by using Lemma \ref{lem;22.1.28.3}.
Then, we obtain that
$dz_j\,\DDlambda_{z_j}s^{(2)}_{0,\alpha}
=(dz_j/z_j)(z_j\DDlambda_{z_j}s^{(2)}_{0,\alpha})$
is also $L^2$ on $U$
with respect to $h$ and $g_{X\setminus H}$
by using Lemma \ref{lem;22.2.23.3} and Lemma \ref{lem;22.2.24.1}.

Because $s^{(0)}_{0,\alpha}$ is $L^2$ on $U$,
we obtain that
$z_j\del_{z_j}s^{(0)}_{0,\alpha}$
is also $L^2$ on $U$
with respect to $h$ and $g_{X\setminus H}$
by using Lemma \ref{lem;22.1.28.3}.
By Lemma \ref{lem;22.2.23.3} and Lemma \ref{lem;22.2.24.1},
$z_j\DDlambda_{z_j}s^{(0)}_{0,\alpha}$
is also $L^2$ on $U$
with respect to $h$ and $g_{X\setminus H}$.
Because
$\frac{dz_j}{z_j}
 \Bigl(
 z_j\DDlambda_{z_j}s^{(0)}_{0,\alpha}
-\DDlambda_{\neq j}s^{(2)}_{0,\alpha}
 \Bigr)$
is $L^2$ on $U$
with respect to $h$ and $g_{X\setminus H}$,
we obtain that
$\DDlambda_{\neq j}s^{(2)}_{\veca,\vecalpha}$
is also $L^2$ on $U$
with respect to $h$ and $g_{X\setminus H}$.
\hfill\qed

\vspace{.1in}

We set
\[
\nbigv_{0,\alpha}=
\bigoplus_{\substack{
(\veca,\vecalpha)\in\openclosed{-1}{0}^j\times\cnum^j
 \\ (a_j,\alpha_j)=(0,\alpha)}}
 \nbigv_{\veca,\vecalpha}.
\]
Let $G_{j,(\veca,\vecalpha,\veck),(\veca',\vecalpha',\veck')}$
be as in \S\ref{subsection;22.2.24.2}.
We obtain the following automorphism of $\nbigv_{0,\alpha}$:
\[
 \nbigf_{j,0,\alpha}:=
\alpha\id_{\nbigv_{0,\alpha}}
+ \sum_{\substack{(\veca,\vecalpha)\in\openclosed{-1}{0}^j\times\cnum^j
 \\ (a_j,\alpha_j)=(0,\alpha)}}
 \sum_{\veck,\veck'\in\seisuu^j}
 G_{j,(\veca,\vecalpha,\veck),(\veca,\vecalpha,\veck')}.
\]
By Lemma \ref{lem;22.2.24.1},
$\nbigf_{j,0,\alpha}$
and
$\nbigf^{-1}_{j,0,\alpha}$
are bounded on $Z(\vecC)$.

We set
$t=\nbigf_{j,0,\alpha}^{-1}(s^{(2)}_{0,\alpha})$,
which is $L^2$ on $U$ with respect to $h$ and $g_{X\setminus H}$.

\begin{lem}
$\DDlambda t$ is $L^2$ on $U$
with respect to $h$ and $g_{X\setminus H}$.
\end{lem}
\pf
Note that $(dz_j/z_j) t$ is $L^2$ on $U$
with respect to $h$ and $g_{X\setminus H}$.
Hence, by the argument in Lemma \ref{lem;22.2.24.10},
we obtain that $dz_j\DDlambda_{z_j}t$ is $L^2$ on $U$
with respect to $h$ and $g_{X\setminus H}$.
By using
Proposition \ref{prop;22.2.23.50}
and Lemma \ref{lem;22.2.24.11},
we obtain that
$\DD^{\lambda}_{\neq j}(\nbigf^{-1}_{j,0,\alpha})$
is bounded on $U$ with respect to $h$ and $g_{X\setminus H}$.
Hence,
$\DD^{\lambda}_{\neq j}(t)$ is also $L^2$
on $U$ with respect to $h$ and $g_{X\setminus H}$.
\hfill\qed

\vspace{.1in}

Let $k>j$.
By the construction,
we have
$t^{\langle k\rangle}_{\veca,\vecalpha;0}=0$
for $(\veca,\vecalpha)$ satisfying
$a_k=0$ and $\alpha_k\neq 0$.
Hence, we obtain 
$(s-\DDlambda t)^{\langle k\rangle}_{\veca,\vecalpha;0}=0$
for $(\veca,\vecalpha)$
satisfying $a_k=0$ and $\alpha_k\neq 0$.
By the construction,
we have
$\bigl(
(s-\DDlambda t)^{(2)}
\bigr)^{\langle j\rangle}_{\veca,\vecalpha;0}=0$
for $(\veca,\vecalpha)$
satisfying $a_j=0$ and $\alpha_j\neq 0$.
Then, it is easy to see
$\bigl(
(s-\DDlambda t)^{(0)}
\bigr)^{\langle j\rangle}_{\veca,\vecalpha;0}=0$
for $(\veca,\vecalpha)$
satisfying $a_j=0$ and $\alpha_j\neq 0$.
Thus, we obtain ${\boldsymbol{P(j)}}$.

\vspace{.1in}
By ${\boldsymbol{P(1)}}$,
there exist $U'\in\gbigu(Q,U)$
and
$t\in \lefttop{\ellsitabar}F
\nbigc^{k-1}_{L^2,\leq 0,\ttX(\Lambda)_{\geq 0}}(\nbigelambda,\DDlambda,h)(U')$
such that
$(s-\DDlambda t)^{\langle j\rangle}_{\veca,\vecalpha;0}=0$
for any $1\leq j\leq \ell$
and $a_j=0$ and $\alpha_j\neq 0$,
i.e.,
$s-\DDlambda t$ is contained in
$\nbigw\nbigc^{\bullet}_{\tw}(\gbige[\ast H],\lambda)(U')$.
Then, by the standard argument,
we can prove that
{\rm(\ref{eq;22.2.23.51})} is a quasi-isomorphism.
\hfill\qed

\subsubsection{End of the proof of Theorem \ref{thm;22.2.18.11}}
\label{subsection;22.2.24.111}

Let $P=(0,\ldots,0)\in X$.
Let $\iota_P:p_1^{-1}(P)\lrarr X(\ttX(\ellsitabar)_{\geq 0})$
denote the inclusion.
We consider the following morphism of sheaves
\[
 \iota_P^{-1}
 \nbigw\nbigc_{\tw}^{\bullet}(\gbige[\ast H])
 \lrarr
 \iota_P^{-1}
 \nbigc_{L^2,\ttX(\Lambda)_{\geq 0}}^{\bullet}(\nbigelambda,\DDlambda,h).
\]
For each $Q\in p_1^{-1}(P)$,
there exists a decreasing sequence of
$(S^1)^{\ell}$-invariant open neighbourhoods
$U^{(j)}_Q$ of $Q$ in $X(\ttX(\ellsitabar)_{\geq 0})$
such that
$\bigcap U^{(i)}_Q=\{Q\}$.
We have
\[
 \nbigc_{L^2,\ttX(\Lambda)_{\geq 0}}^{\bullet}(\nbigelambda,\DDlambda,h)_Q
 =\varinjlim
 \nbigc_{L^2,\ttX(\Lambda)_{\geq 0}}^{\bullet}
 (\nbigelambda,\DDlambda,h)(U^{(i)}_Q).
\]
We define
\[
 \lefttop{\ellsitabar}F
 \nbigc_{L^2,\ttX(\Lambda)_{\geq 0}}^{\bullet}
 (\nbigelambda,\DDlambda,h)_Q
 =\varinjlim \lefttop{\ellsitabar}F
 \nbigc_{L^2,\ttX(\Lambda)_{\geq 0}}^{\bullet}
 (\nbigelambda,\DDlambda,h)(U^{(i)}_Q).
\]
It defines the subsheaf
$\lefttop{\ellsitabar}F
\iota_P^{-1}
\nbigc_{L^2,\ttX(\Lambda)_{\geq 0}}^{\bullet}
(\nbigelambda,\DDlambda,h)$
of
$\iota_P^{-1}
\nbigc_{L^2,\ttX(\Lambda)_{\geq 0}}^{\bullet}
(\nbigelambda,\DDlambda,h)$.

The complex of sheaves 
$\lefttop{\ellsitabar}F
\iota_P^{-1}
 \nbigctilde_{L^2}^{\bullet}(\nbigelambda,h)$
 and 
$\iota_P^{-1}
 \nbigctilde_{L^2}^{\bullet}(\nbigelambda,h)$
are $c$-soft with respect to $p_1$,
and 
the induced morphism
\[
 p_{1\ast}
\lefttop{\ellsitabar}F\bigl(
\iota_P^{-1}
 \nbigctilde_{L^2}^{\bullet}(\nbigelambda,h)
 \bigr)
 \lrarr
p_{1\ast}
 \iota_P^{-1}
 \nbigctilde_{L^2}^{\bullet}(\nbigelambda,h)
\]
is a quasi-isomorphism
by Proposition \ref{prop;22.1.28.10}.
The natural morphism
\[
 \iota_P^{-1}
 \nbigw\nbigctilde^{\bullet}_{\tw}
 (\gbige[\ast H])
\lrarr
 \lefttop{\ellsitabar}F\bigl(
\iota_P^{-1}
 \nbigctilde_{L^2}^{\bullet}(\nbigelambda,h)
 \bigr)
\]
is a quasi-isomorphism
by Proposition \ref{prop;22.2.23.43}
and Proposition \ref{prop;22.2.23.52}.
Thus, we obtain the claim of
Theorem \ref{thm;22.2.18.11}.
\hfill\qed

\subsubsection{Proof of Theorem \ref{thm;21.12.21.11}}

Let $U(\vecx,\vecepsilon)$ be as in \S\ref{subsection;22.2.23.53}.
Let $\nbiguzero(\vecx,\vecepsilon)$
denote the product of $U(\vecx,\vecepsilon)$
and a neighborhood $U(\lambda_0)$ of $\lambda_0$ in $\cnum$.

Let $s\in
\lefttop{\ellsitabar}F
\nbigc^k_{L^2,\leq 0,\ttX(\Lambda)_{\geq 0}}
(\nbige,\DD,h)(\overline{\nbiguzero(\vecx,\vecepsilon)})$.
We have the expression
$s=\sum s_{J,\emptyset,i}\omega_{J,\emptyset}w_i$.
By the argument in \S\ref{subsection;22.2.23.53},
we obtain that
$s_{J,\emptyset,i}$ are holomorphic functions
on a neighbourhood $N$
of $\nbigh_{\ellsitabar}$,
and hence
\[
 s_{J,\emptyset,i}
 =\sum_{\vecm}
 \stilde_{J,i,\vecm}(\lambda,z_{\ell+1},\ldots,z_n)
 \prod_{j=1}^{\ell}z_j^{m_j}.
\]

Let $s$ be a holomorphic section of
$(j_{\nbigx\setminus\nbigh})_{\ast}\Bigl(
\nbige\otimes\Omega^{k,0}\Bigr)$ on $N$.
Suppose that $\nbiguzero(\vecx,\vecepsilon)$
is a relatively compact open subset of $N$.
By the norm estimate as in Proposition \ref{prop;22.2.24.100},
$s$ is $L^2$ on $\nbiguzero(\vecx,\vecepsilon)$
with respect to $h$ and $g_{X\setminus H}$
if and only if the following integral is convergent:
\begin{equation}
\label{eq;22.2.24.101}
 \sum_{J,i,\vecm}
 \int_{\nbiguzero(\vecx,\vecepsilon)}
 \bigl|\stilde_{J,i,\vecm}(\lambda,z_{\ell+1},\ldots,z_n)\bigr|^2
 \prod_{p=1}^{\ell}
 |z_p|^{-2\paramap(\lambda,u_p(w_i)+2m_p}
 (-\log|z_p|)^{k_p(w_i)-k_{p-1}(w_i)}
 \prod_{\substack{p\in J\\ p\leq \ell }}(-\log|z_p|)^{2}
 \dvol_{g_{X\setminus H}}.
\end{equation}
The integral (\ref{eq;22.2.24.101})
is convergent if and only if
the following is convergent for each $J$ and $i$:
{\small
\begin{multline}
\label{eq;22.2.24.102}
 \int_{\nbiguzero(\vecx,\vecepsilon)}
 |s_{J,\emptyset,i}\omega_{J,\emptyset}w_i|^2\dvol_{g_{X\setminus H}}
 =\\
 \sum_{\vecm}
\int_{\nbiguzero(\vecx,\vecepsilon)}
 \bigl|\stilde_{J,i,\vecm}(\lambda,z_{\ell+1},\ldots,z_n)\bigr|^2
 \prod_{p=1}^{\ell}
 |z_p|^{-2\paramap(\lambda,u_p(w_i))+2m_p}
 (-\log|z_p|)^{k_p(w_i)-k_{p-1}(w_i)}
 \prod_{\substack{p\in J\\ p\leq \ell}}(-\log|z_p|)^{2}
 \dvol_{g_{X\setminus H}}.
\end{multline}
}

Let $\nbigt(s_{J,\emptyset,i})=
\bigl\{\vecm\in\seisuu^{\ell}_{\geq 0}\,\big|\,
\stilde_{J,i,\vecm}\neq 0
\bigr\}$.
For any $\vecm\in\nbigt(s_{J,\emptyset,i})$,
we have
$\vecm-\paramap(\lambda_0,\vecu(w_i))\in \real_{\geq 0}^{\ell}$.
We can check the following lemma.
\begin{lem}
\label{lem;22.2.24.110}
The integral {\rm(\ref{eq;22.2.24.102})} is convergent
if and only if
for each $\vecm\in\min\nbigt(s_{J,\emptyset,i})$
one of the following holds: 
\begin{itemize}
 \item $q_{K_1}(\vecm-\paramap(\lambda_0,\vecu(w_i)))\neq 0$.
       Moreover,
       if $m_j-\paramap(\lambda_0,u_j(w_i))=0$,
       then $m_j=0$ and $u_j(w_i)=(0,0)$ hold.

 \item There exists $1\leq j_0\leq m-1$ such that the following holds:
     \begin{itemize}
      \item
	   $q_{K_{j_0}}(\vecm)=0$,
	   $q_{K_{j_0}}(\vecu(w_i))=0\in(\real\times\cnum)^{K_{j_0}}$,
	   and
	   $q_{K_{j_0+1}}(\vecm-\veca(w_i))\neq 0$.
	   Moreover,
	   if $m_j-\paramap(\lambda_0,u_j(w_i))=0$,
       then $m_j=0$ and $u_j(w_i)=(0,0)$ hold.
      \item $k_{p(s)}(w_i)\leq p(s)-2|J\cap K_s|-1$
       for $s=1,\ldots,j_0$.
     \end{itemize}
 \item $\vecm=0$, $\vecu(w_i)=0\in(\real\times\cnum)^{\ellsitabar}$,
       and $k_{p(s)}(w_i)\leq p(s)-2|J\cap K_s|-1$
       for any $s$.
\hfill\qed
\end{itemize}
\end{lem}

By using Lemma \ref{lem;22.2.24.110},
we obtain that the natural morphism
\[
 \nbigw
 \nbigc^{\bullet}_{\tw}
 (\gbige[\ast H])_{(\lambda_0,Q)}
 \lrarr
  \lefttop{\ellsitabar}F
 \nbigc^{\bullet}_{L^2,\leq 0,\ttX(\Lambda)_{\geq 0}}
 (\nbige,\DD,h)_{(\lambda_0,Q)}.
\]
is an isomorphism,
not only a quasi-isomorphism.
Then, by using the argument in the proof of 
\S\ref{subsection;22.2.24.111},
we obtain Theorem \ref{thm;22.2.17.111}.
\hfill\qed

\subsubsection{Convergence of some integrals (Appendix)}

Let $\vecx=(x_1,\ldots,x_{\ell})\in\real^{\ell-1}_{\geq 0}\times\{0\}$.
We obtain the increasing sequence
$0=p(0)< p(1)\leq \cdots< p(m)=\ell$
determined by
$\{p(1),\ldots,p(m)\}=\{1\leq i\leq \ell\,|\,x_i=0\}$.
We set $K_i=\{1,\ldots,p(i)\}$.
Let $q_{K_i}:\real^{\ell}\lrarr\real^{K_i}$
denote the projection.
Let $\vecepsilon=(\epsilon_1,\ldots,\epsilon_{\ell})\in\real_{>0}^{\ell}$.
We set
\[
 Z(\vecx,\vecepsilon)=\Bigl\{
 (s_1,\ldots,s_{\ell})\in\real_{>0}^{\ell}\,\Big|\,
 |s_{j+1}s_j^{-1}-x_j|<\epsilon_j\,\,(j=1,\ldots,\ell-1),\,\,
 s_{\ell}^{-1}\leq \epsilon_{\ell}
 \Bigr\}.
\]

For $\vecb\in\real_{\geq 0}^{\ell}$,
$\veck\in\seisuu^{\ell}$
and $J\subset\ellsitabar$,
we consider the following integral:
\begin{equation}
\label{eq;22.1.27.10}
 \int_{Z(\vecx,\vecepsilon)}
 e^{-\sum_{j=1}^{\ell} b_js_j}
 \prod_{j=1}^{\ell} s_j^{(k_j-k_{j-1})-2}
 \prod_{j\in J}s_j^2
 \prod_{j=1}^{\ell} ds_j
\end{equation}
We can check the following lemma.
\begin{lem}
\label{lem;22.1.27.11}
The integral {\rm(\ref{eq;22.1.27.10})}
 is convergent if and only if
 either one of the following holds:
\begin{itemize}
 \item $q_{K_1}(\vecb)\neq 0$.
 \item There exists $j_0=1,\ldots,m-1$
       such that the following holds:
       $q_{K_{j_0}}(\vecb)=0$,
       $q_{K_{j_0+1}}(\vecb)\neq 0$,
       and
       $k_{p(i)}\leq p(i)-2|J\cap K_i|-1$ for $i=1,\ldots,j_0$.
 \item $\vecb=0$ and
       $k_{p(i)}\leq p(i)-2|J\cap K_i|-1$ for $i=1,\ldots,m$.
\hfill\qed
\end{itemize}
\end{lem}

\subsection{Appendix: The complex of the sheaves of $L^2$-holomorphic sections}
\label{subsection;22.3.18.30}

\subsubsection{A local example}
\label{subsection;22.4.13.2}

We set
$X=\bigl\{(z_1,z_2)\in\cnum^2\,\big|\,|z_i|<1\bigr\}$
and $H=X\cap \{z_1z_2=0\}$.
We use the Poincar\'{e} metric $g_{X\setminus H}$.

We set
$E=\nbigo_{X\setminus H}v_1\oplus\nbigo_{X\setminus H}v_2$,
and we consider the Higgs field $\theta$ defined as follows:
\[
\theta(v_1)=
\left(
\frac{dz_1}{z_1}+\frac{dz_2}{z_2}
\right)\cdot v_2,
\quad
\quad
\theta(v_2)=0.
\]
We consider the Hermitian metric of $h$ given as follows:
\[
 h(v_1,v_1)=-(\log|z_1|^2+\log|z_2|^2),
 \quad
 h(v_2,v_2)=-(\log|z_1|^2+\log|z_2|^2)^{-1},
 \quad
 h(v_1,v_2)=0. 
\]
Then, $(E,\theta,h)$ is a harmonic bundle.
This harmonic bundle underlies a polarized variation of complex Hodge structure.
(See \cite{s1}).
It is obtained as the pull back of
a basic nilpotent orbit in one variable on a punctured disc
by the map $(z_1,z_2)\to z_1z_2$.

We obtain the $\lambda$-flat bundle
$(\nbigelambda,\DDlambda)$ with the Hermitian metric $h$.
The $(0,1)$-part of $\DDlambda$ is denoted by
$\delbar_{\nbigelambda}$.

For any open subset $U\subset X$,
let $\nbigc^k_{L^2}(\nbigelambda,\DDlambda,h)(U)$
denote the space of measurable sections $\tau$ of
\[
 \bigoplus_{k_1+k_2=k} \Omega^{k_1,k_2}_{X\setminus H}\otimes E
\]
on $U\setminus H$
such that $\tau$ and $\DDlambda\tau$ are $L^2$
around any point of $U$ with respect to $h$ and $g_{X\setminus H}$.
Let $\nbigc^k_{L^2,\hol}(\nbigelambda,\DDlambda,h)(U)$
denote the space of holomorphic sections $\tau$
of $\Omega_{X\setminus H}^{k,0}\otimes\nbigelambda$ on $U\setminus H$
such that $\tau$ and $\DDlambda\tau$ are $L^2$
with respect to $h$ and $g_{X\setminus H}$
locally around any point of $U$.
We obtain the complexes of sheaves
$\nbigc^{\bullet}_{L^2}(\nbigelambda,\DDlambda,h)$
and
$\nbigc^{\bullet}_{L^2,\hol}(\nbigelambda,\DDlambda,h)$.
There exists the natural inclusion
\begin{equation}
\label{eq;22.3.18.1}
 \nbigc^{\bullet}_{L^2,\hol}(\nbigelambda,\DDlambda,h)
 \lrarr
\nbigc^{\bullet}_{L^2}(\nbigelambda,\DDlambda,h).
\end{equation}
Let us remark the following.
\begin{prop}
\label{prop;22.3.18.20}
The morphism {\rm(\ref{eq;22.3.18.1})} is not a quasi-isomorphism. 
\end{prop}
\pf
For $\lambda\neq 1$,
it is proved by Cattani and Kaplan in \cite{ck2}.
Their argument essentially works even in the case $\lambda=0$,
which we explain in the following.

We set $T=(S^1)^2$.
Let us consider the $T$-action on $X$
given by $(a_1,a_2)(z_1,z_2)=(a_1z_1,a_2z_2)$.
Note that $(v_1,v_2)$ is a frame of $E$ on $X\setminus H$.
We define the $T$-action on $E$ by $(a_1,a_2)^{\ast}(v_i)=v_i$.
Then, the harmonic bundle is naturally $T$-equivariant.
Let $U\subset X$ be any $T$-invariant open subset.
Any $\tau\in \nbigc^{1}_{L^2}(\nbige^0,\DD^0,h)(U)$
is expressed as
\[
\tau=
\sum_{i=1,2}
\Bigl(
 \tau^{(i)}_{1,\emptyset}\frac{dz_1}{z_1}
+\tau^{(i)}_{2,\emptyset}\frac{dz_2}{z_2}
+\tau^{(i)}_{\emptyset,1}\frac{d\zbar_1}{\zbar_1}
+\tau^{(i)}_{\emptyset,2}\frac{d\zbar_2}{\zbar_2}
\Bigr)
v_i.
\]
Let $z_j=r_je^{\sqrt{-1}\theta_j}$ be the polar decomposition.
There exists the Fourier expansions
\[
 \tau^{(i)}_{J,K}
=\sum_{\vecm=(m_1,m_2)\in\seisuu^2}
 \tau^{(i)}_{J,K;\vecm}(r_1,r_2)
 e^{\sqrt{-1}(m_1\theta_1+m_2\theta_2)}.
\]
We set
\[
 \tau^T:=
 \sum_{i=1,2}
 \Bigl(
 \tau^{(i)}_{1,\emptyset;(0,0)}\frac{dz_1}{z_1}
+\tau^{(i)}_{2,\emptyset;(0,0)}\frac{dz_2}{z_2}
+\tau^{(i)}_{\emptyset,1;(0,0)}\frac{d\zbar_1}{\zbar_1}
+\tau^{(i)}_{\emptyset,2;(0,0)}\frac{d\zbar_2}{\zbar_2}
\Bigr)
v_i,
\quad\quad
\tau^{\bot}:=\tau-\tau^{T}.
\]
We obtain the decomposition
$\tau=\tau^T+\tau^{\bot}$.
This induces a decomposition
\[
 \nbigc^{\bullet}_{L^2}(\nbige^0,\DD^0,h)(U)
=
 \nbigc^{\bullet}_{L^2}(\nbige^0,\DD^0,h)(U)^{T}
 \oplus
 \nbigc^{\bullet}_{L^2}(\nbige^0,\DD^0,h)(U)^{\bot}.
\]
Similarly,
we obtain a decomposition
\[
 \nbigc^{\bullet}_{L^2,\hol}(\nbige^0,\DD^0,h)(U)
=
 \nbigc^{\bullet}_{L^2,\hol}(\nbige^0,\DD^0,h)(U)^{T}
 \oplus
 \nbigc^{\bullet}_{L^2,\hol}(\nbige^0,\DD^0,h)(U)^{\bot}.
\]

Let 
$\nbigc^{\bullet}_{L^2}(\nbige^0,\DD^0,h)_O$
and
$\nbigc^{\bullet}_{L^2,\hol}(\nbige^0,\DD^0,h)_O$
denote the stalks of
$\nbigc^{\bullet}_{L^2}(\nbige^0,\DD^0,h)$
and
$\nbigc^{\bullet}_{L^2,\hol}(\nbige^0,\DD^0,h)$
at $O=(0,0)$, respectively.
Let $\gbigu(O)$ denote the set of
$T$-invariant open neighbourhoods of $O$.
We have
\[
 \nbigc^{\bullet}_{L^2}(\nbige^0,\DD^0,h)_O
=\varinjlim_{U\in\gbigu(O)}
  \nbigc^{\bullet}_{L^2}(\nbige^0,\DD^0,h)(U).
\]
We set
\[
\nbigc^{\bullet}_{L^2}(\nbige^0,\DD^0,h)_O^{T}
=\varinjlim_{U\in\gbigu(O)}
\nbigc^{\bullet}_{L^2}(\nbige^0,\DD^0,h)(U)^T,
\quad\quad
\nbigc^{\bullet}_{L^2}(\nbige^0,\DD^0,h)_O^{\bot}
=\varinjlim_{U\in\gbigu(O)}
\nbigc^{\bullet}_{L^2}(\nbige^0,\DD^0,h)(U)^{\bot}.
\]
We obtain a decomposition
\[ 
 \nbigc^{\bullet}_{L^2}(\nbige^0,\DD^0,h)_O
=\nbigc^{\bullet}_{L^2}(\nbige^0,\DD^0,h)_O^{T}
\oplus
\nbigc^{\bullet}_{L^2}(\nbige^0,\DD^0,h)_O^{\bot}.
\]
Similarly, by setting
\[
 \nbigc^{\bullet}_{L^2,\hol}(\nbige^0,\DD^0,h)_O^{T}
=\varinjlim_{U\in\gbigu(O)}
\nbigc^{\bullet}_{L^2,\hol}(\nbige^0,\DD^0,h)(U)^T,
\quad\quad
\nbigc^{\bullet}_{L^2,\hol}(\nbige^0,\DD^0,h)_O^{\bot}
=\varinjlim_{U\in\gbigu(O)}
\nbigc^{\bullet}_{L^2,\hol}(\nbige^0,\DD^0,h)(U)^{\bot},
\]
we obtain a decomposition
\[
  \nbigc^{\bullet}_{L^2,\hol}(\nbige^0,\DD^0,h)_O
=\nbigc^{\bullet}_{L^2,\hol}(\nbige^0,\DD^0,h)_O^{T}
\oplus
\nbigc^{\bullet}_{L^2,\hol}(\nbige^0,\DD^0,h)_O^{\bot}.
\]
It is enough to prove that
$\nbigc^{\bullet}_{L^2,\hol}(\nbige^0,\DD^0,h)_O^{T}
\lrarr
\nbigc^{\bullet}_{L^2}(\nbige^0,\DD^0,h)_O^{T}$
is not a quasi-isomorphism.

For any $U\in\gbigu(O)$,
we obtain
$\nbigc^{0}_{L^2,\hol}(\nbige^0,\DD^0,h)(U)^{T}
=\cnum\cdot v_2$
and
\[
\nbigc^{1}_{L^2,\hol}(\nbige^0,\DD^0,h)(U)^{T}
=
\nbigc^{2}_{L^2,\hol}(\nbige^0,\DD^0,h)(U)^{T}
=0.
\]
Hence, we have
\[
 \nbigc^{0}_{L^2,\hol}(\nbige^0,\DD^0,h)_O^{T}
=\cnum\cdot v_2,
\quad\quad
\nbigc^{1}_{L^2,\hol}(\nbige^0,\DD^0,h)_O^{T}
=
\nbigc^{2}_{L^2,\hol}(\nbige^0,\DD^0,h)_O^{T}
=0.
\]

We set $U_0=\bigl\{(z_1,z_2)\,\big|\,|z_i|<1/2\bigr\}$.
Let us consider the following function on $U_0\setminus H$:
\[
 f(z_1,z_2)=\frac{\log|z_2|-\log|z_1|}{\log|z_2|+\log|z_1|}.
\]
As in \cite{ck2},
we set
\[
 \tau:=
 \delbar(f)v_1+
 \left(
 (f+1)\frac{dz_1}{z_1}
+(f-1)\frac{dz_2}{z_2}
 \right)
 v_2
 =\delbar(f)v_1+
 f\cdot\left(
 \frac{dz_1}{z_1}
+\frac{dz_2}{z_2}
 \right) v_2
+
 \left(
 \frac{dz_1}{z_1}
-\frac{dz_2}{z_2}
 \right)
 v_2.
\]
Then, we can check that $\tau$ is $L^2$ on any $U\subset U_0$,
and that it satisfies $\DD^0\tau=0$.
Moreover,
for any $U\in \gbigu(O)$ such that $U\subset U_0$,
there does not exist $u\in\nbigc^0_{L^2}(\nbige^0,\DD^0)(U)^T$
such that $\delbar u=\tau_{|U}$.
Thus, we obtain Proposition \ref{prop;22.3.18.20}.
\hfill\qed

\subsubsection{A global example}

The following lemma is well known.
\begin{lem}
There exists a compact Riemann surface $C$
with a finite non-empty subset $D$
and a polarized variation of Hodge structure
$(V=V^{1,0}\oplus V^{0,1},\nabla,\langle\cdot,\cdot\rangle)$
of weight $1$ 
with an integral structure on $C\setminus D$
such that the local monodromy around each point of $D$
is unipotent but non-trivial. 
\end{lem}
\pf
For example, we can obtain such an example as follows.
We set $\hyperh:=\bigl\{\zeta\,\big|\,\Image(\zeta)>0\bigr\}$.
Let $m$ be an integer larger than $5$.
Let $\Gamma(m)\subset\SL(2,\seisuu)$
be the principal congruence subgroup,
i.e.,
$\Gamma(m)$ is the subgroup of
$\gamma\in\SL(2,\seisuu)$ satisfying
\[
 \gamma\equiv
 \left(
 \begin{array}{cc}
  1 & 0 \\ 0 & 1
 \end{array}
 \right)
 \mod m.
\]
The natural action of $\Gamma(m)$ on $\hyperh$ is proper.
If $\gamma\in \SL(2,\seisuu)$ has a fixed point,
$\gamma$ is an elliptic element and hence $|\tr(\gamma)|<2$ holds.
Because $\tr(\gamma)=2$ modulo $m$
for $\gamma\in\Gamma(m)\setminus\{1\}$,
the action of $\Gamma(m)$ is free.
The quotient space $C^{\circ}=\hyperh/\Gamma(m)$
is the complement of a non-empty finite subset $D$
in a compact Riemann surface $C$.
We have the universal elliptic curve
$\pi:\mathbb{E}\to \hyperh$ with a section.
By considering $R^1\pi_{\ast}\seisuu_{\mathbb{E}}$,
we obtain the polarized variation of Hodge structure
of weight $1$ with an integral structure on $\hyperh$,
which is equivariant with respect to $\Gamma(m)$.
There exists the induced polarized variation of Hodge structure
on $C^{\circ}$.
Because the trace of $\gamma\in\Gamma(m)$ cannot be $-2$,
the local monodromy is unipotent.
\hfill\qed

\vspace{.1in}
We fix a point $P\in D$.
The following lemma is easy to see.
\begin{lem}
There exists a morphism
$\rho:C\lrarr \proj^1$
such that
 (i) $\rho(P)=0$,
(ii)  $0$ and $\infty$ are regular values of $\rho$.

\hfill\qed
\end{lem}

Let us consider the morphism
$\varphi:\cnum^2\to \cnum$ given by
$(z_1,z_2)\longmapsto z_1z_2$.
Let $Y\lrarr (\proj^1)^2$
denote the blow up at $\{(\infty,0),(0,\infty)\}$.
There exists the morphism
$\varphi_Y:Y\lrarr \proj^1$
whose restriction to $\cnum^2$ is equal to $\varphi$.

Let $\Ytilde$ be the fiber product of
$\varphi_Y:Y\to \proj^1$
and $\rho:C\to\proj^1$.
Because $\rho$ is a ramified covering,
and because the critical value of $\rho$
is a regular value of $\varphi_Y$,
$\Ytilde$ is smooth.

Let $\pi_1:\Ytilde\to Y$
and $\pi_2:\Ytilde\to C$
denote the projections.
For any $Q\in C$,
we set $K_Q:=\pi_2^{-1}(Q)$.
If $\rho(Q)\neq 0,\infty$,
$K_Q$ is smooth.
If $\rho(Q)=0,\infty$,
then $K_Q$ is normal crossing.

We set $K_D:=\bigcup_{Q\in D}K_Q$.
Let $(E_1,\delbar_{E_1},\theta_1,h_1)$
denote the tame harmonic bundle
underlying the polarized variation of Hodge structure
$(V=V^{1,0}\oplus V^{0,1},\nabla,\langle\cdot,\cdot\rangle)$.
We obtain the tame harmonic bundle
$(E_2,\delbar_{E_2},\theta_2,h_2)
=(\rho\circ p_2)^{-1}(E_1,\delbar_{E_1},\theta_1,h_1)$
on $(\Ytilde,K_D)$.
It underlies a polarized variation of Hodge structure
of weight $1$ with an integral structure,
and the local monodromy around any irreducible component
of $K_D$ is unipotent.

Let $X$, $H$ and $(E,\theta,h)$ be as in \S\ref{subsection;22.4.13.2}.
Let $\Ptilde\in\Ytilde$ denote the point determined by
$(0,0)\in Y$ and $P\in C$.
We may naturally regard $X$ as a neighbourhood of $\Ptilde$
in $\Ytilde$.
There exists the natural isomorphism
$(E,\delbar_E,\theta)\simeq (E_2,\delbar_{E_2},\theta_2)$
of the Higgs bundles.
The pluri-harmonic metrics
$h$ and $h_{2}$ are mutually bounded,
which follows from the norm estimate of Simpson
in \cite{Simpson-non-compact}.
Hence, by Proposition \ref{prop;22.3.18.20},
the $L^2$-holomorphic Dolbeault complex of
$(E_2,\delbar_{E_2},\theta_2,h_2)$
is not quasi-isomorphic to
the $L^2$-Dolbeault complex of
$(E_2,\delbar_{E_2},\theta_2,h_2)$.

\begin{rem}
It is more appropriate to add some assumptions 
in {\rm\cite[Theorem A]{Jost-Yang-Zuo}}
because
{\rm\cite[Proposition 2]{Jost-Yang-Zuo}}
does not necessarily hold, in general.
J. Jost kindly informed that
X. M. Ye already discovered a counterexample
in the Hilbert modular surface case
around 2011. 
\hfill\qed
\end{rem}

\section{Hard Lefschetz Theorem}
\label{section;22.4.22.1}

\subsection{Main Theorem}

\subsubsection{The absolute case}

Let $X$ be a complex manifold with a K\"ahler form $\omega$.
Let $\pt$ be the one point set.
Let $a_X:X\lrarr \pt$ denote the canonical morphism.
Let $\nbigt$ be a regular pure twistor $\nbigd$-module of weight $w$
on $X$ with a polarization $\nbigs$.
Assume that the support of $\nbigt$ is compact.
We obtain the graded $\nbigr$-triple
$\bigoplus_ia_{X\dagger}^i(\nbigt)$
equipped with a Lefschetz operator $L_{\omega}$
and a sesqui-linear duality $a_{X\dagger}(\nbigs)$
of weight $w$.
(See \S\ref{subsection;22.4.28.2} and \S\ref{subsection;22.3.25.100}.)
We shall prove the following theorem
in \S\ref{subsection;22.3.10.2}--\S\ref{subsection;22.4.28.1}.

\begin{thm}
\label{thm;22.2.25.40}
$(\bigoplus a_{X\dagger}^j(\nbigt),L_{\omega},a_{X\dagger}\nbigs)$
is a polarized graded Lefschetz twistor structure of weight $w$ and type $1$
in the sense of {\rm\cite{sabbah2}}.
(See {\rm\S\ref{subsection;22.4.28.3}}.)
Namely, the Hard Lefschetz theorem holds for
$(\nbigt,\nbigs)$.
\end{thm}

\subsubsection{The relative case}
\label{subsection;22.3.10.10}

Let $f:X\lrarr Y$ be a morphism of complex manifolds.
Let $\nbigt$ be a pure twistor $\nbigd$-module of weight $w$ on $X$
with a polarization $\nbigs$
whose support is proper over $Y$.
Let $\omega$ be a K\"ahler class on $X$.

\begin{thm}
\label{thm;22.3.5.1}
$\bigl(
\bigoplus_jf_{\dagger}^j(\nbigt),L_{\omega},f_{\dagger}\nbigs
\bigr)$
is a polarized graded Lefschetz twistor $\nbigd$-module
of weight $w$ and type $1$ on $Y$ in the sense of {\rm \cite{sabbah2}}
(see {\rm\S\ref{subsection;22.4.28.4}}).
\end{thm}
\pf
We recall that there exists a unique decomposition
$(\nbigt,\nbigs)=
\bigoplus (\nbigt_Z,\nbigs_Z)$ by strict supports,
where each $(\nbigt_Z,\nbigs_Z)$ is a polarized
pure twistor $\nbigd$-module of weight $w$
whose strict support is a closed irreducible complex subvariety $Z$
of $X$.
Hence, it is enough to study the case where
$\nbigt$ has an irreducible strict support,
denoted by $\supp(\nbigt)$.

Let us consider the following claim for $n,m\in\seisuu_{\geq 0}$.
\begin{description}
 \item[$\boldsymbol{P(n,m)}$:]
	    If $\dim\supp(\nbigt)\leq n$
	    and $\dim f(\supp(\nbigt))\leq m$,
	    then
	    $(\bigoplus_jf^j_{\dagger}\nbigt,L_{\omega},f_{\dagger}\nbigs)$
	    is a polarized graded Lefschetz twistor $\nbigd$-module of
	    weight $w$ and type $1$ on $Y$.
\end{description}
By the argument
in the proof of Hard Lefschetz Theorem for projective morphisms
in \cite{sabbah2, mochi2, Mochizuki-wild},
which goes back to the Hodge case \cite{saito1},
we can prove that $P(n,m)$ implies $P(n+1,m+1)$
for any $n,m\in\seisuu_{\geq 0}$.

To deduce $P(n,0)$ from Theorem \ref{thm;22.2.25.40},
we give a preliminary consideration.
Let $Z_0$ be a compact irreducible complex analytic subset
of $X$ such that $f(Z_0)$ is a point.
Let $Z_1\subset Z_0$ be a closed complex analytic subset
$Z_0\setminus Z_1$ is a locally closed complex submanifold of $X$.
There exists a projective morphism of complex manifolds
$\varphi:\Xtilde\to X$ 
such that
(i) $\varphi$ induces an isomorphism
$\Xtilde\setminus \varphi^{-1}(Z_1)\to X\setminus Z_1$,
(ii) the proper transform $\Ztilde_0$ of $Z_0$
is a complex submanifold of $\Xtilde$,
(iii) $\varphi^{-1}(Z_1)$ and $\Ztilde_1=\Ztilde_0\cap\varphi^{-1}(Z_1)$
are simple normal crossing hypersurfaces of $\Xtilde$
and $\Ztilde_0$, respectively,
(iv) $\varphi$ is a composition of blowing-ups
along smooth centers.
(See \cite{Wlodarczyk}.)
Note that $\Xtilde$ and 
$\Ztilde_0$ are K\"ahler manifolds.
The induced morphism $\Ztilde_0\to X$ is also denoted by $\varphi$.
Let $(E,\delbar_E,\theta,h)$ be a tame harmonic bundle
on $(\Ztilde_0,\Ztilde_1)$.
We have the associated polarized pure twistor $\nbigd$-module
$(\gbigt(E),\gbigs)$ of weight $0$ on $\Ztilde_0$.
(See \S\ref{subsection;22.4.12.1}.)
By Theorem \ref{thm;22.2.25.40},
$a_{\Ztilde_0\dagger}^j(\gbigt(E))$
are pure twistor structures of weight $j$.
Let $i_0:\pt\to Y$ denote the inclusion
defined by $i_0(\pt)=f(Z_0)$.
Because $i_0\circ a_{\Ztilde_0}=f\circ\varphi$,
we obtain
$(f\circ\varphi)^j_{\dagger}(\gbigt(E))$
are polarizable pure twistor $\nbigd$-modules 
of weight $j$ on $Y$.
By the Hard Lefschetz Theorem for the projective morphisms
in \cite{sabbah2,mochi2,Mochizuki-wild},
$\varphi_{\dagger}(\gbigt(E))$
is isomorphic to
$\bigoplus \varphi_{\dagger}^j(\gbigt(E))[-j]$
in the derived category of $\nbigr_{X}$-triples.
It implies that
$(f\circ\varphi)^j_{\dagger}(\gbigt(E))
\simeq
 \bigoplus_{k+\ell=j}
 f_{\dagger}^k(\varphi_{\dagger}^{\ell}(\gbigt(E)))$,
and hence each 
$f_{\dagger}^k(\varphi_{\dagger}^{\ell}(\gbigt(E)))$
is a polarizable pure twistor $\nbigd$-module of weight $k+\ell$
on $Y$.
Let $\nbigt_1$ be any direct summand 
of $\varphi_{\dagger}^0(\gbigt(E))$.
Because
$f_{\dagger}^k(\nbigt_1)$ is
a direct summand of $f_{\dagger}^k\varphi_{\dagger}^0(\gbigt(E))$,
$f_{\dagger}^k(\nbigt_1)$ is
a polarizable pure twistor $\nbigd$-module
of weight $k$ on $Y$
whose support is contained in $f(Z_0)$.
We obtain
$a_{Y\dagger}^jf_{\dagger}^k(\nbigt_1)=0$
unless $j=0$,
and 
$a_{X\dagger}^k(\nbigt_1)=
(a_Y\circ f)_{\dagger}^k(\nbigt_1)
=a_{Y\dagger}^0f_{\dagger}^k(\nbigt_1)$.
We obtain
\begin{equation}
\label{eq;22.4.20.10}
f_{\dagger}^k(\nbigt_1)
\simeq
 i_{0\dagger}\circ a_{Y\dagger}^0
 \circ f_{\dagger}^k(\nbigt_1)
\simeq
 i_{0\dagger}a_{X\dagger}^k(\nbigt_1).
\end{equation}
Because $\nbigt_1$ is a direct summand of
a polarizable pure twistor $\nbigd$-module $\varphi_{\dagger}^0\gbigt(E)$
of weight $0$ on $X$,
$\nbigt_1$ is a polarizable pure twistor $\nbigd$-module
of weight $0$ on $X$.
For a polarization
$\nbigs_1:\nbigt_1\simeq \nbigt_1^{\ast}$ of weight $0$,
the induced morphisms
$f_{\dagger}(\nbigs_1):
f_{\dagger}^k(\nbigt_1)
\simeq
\bigl(
f_{\dagger}^{-k}(\nbigt_1)
\bigr)^{\ast}$
are equal to
$i_{0\dagger}a_{X\dagger}(\nbigs_1):
 i_{0\dagger}a_{X\dagger}^k(\nbigt_1)
 \simeq
 \bigl(
 i_{0\dagger}a_{X\dagger}^{-k}(\nbigt_1)
 \bigr)^{\ast}$
under the isomorphism (\ref{eq;22.4.20.10}).
The morphisms
$L_{\omega}:f_{\dagger}^j(\nbigt_1)\to
f_{\dagger}^{j+2}(\nbigt_1)\otimes\newTate(1)$
are equal to
$i_{0\dagger}(L_{\omega}):
 i_{0\dagger}a_{X\dagger}^j(\nbigt_1)\to
 i_{0\dagger}a_{X\dagger}^{j+2}(\nbigt_1)
 \otimes\newTate(1)$
under the isomorphism (\ref{eq;22.4.20.10}).
By Theorem \ref{thm;22.2.25.40},
$\bigl(
\bigoplus a_{X\dagger}^j(\nbigt_1),
L_{\omega},
a_{X\dagger}\nbigs_1
\bigr)$
is a polarized graded Lefschetz twistor structure of weight $0$
and type $1$,
we obtain that 
$\bigl(
\bigoplus f_{X\dagger}^j(\nbigt_1),
L_{\omega},
f_{X\dagger}\nbigs_1
\bigr)$
is a polarized graded Lefschetz twistor structure of weight $0$
and type $1$.

Let us prove $P(n,0)$.
It is enough to consider the case where $w=0$.
We set $Z_0:=\supp(\nbigt)$.
By Theorem \ref{thm;22.4.20.20} and Theorem \ref{thm;22.4.18.30},
there exists a projective morphism of complex manifolds
$\varphi:\Ztilde_0\to X$ 
and a polarizable pure twistor $\nbigd$-module
$\gbigt(E)$ of weight $0$ on $\Ztilde_0$ as above
such that
$\nbigt$ is a direct summand of
$\varphi_{\dagger}^0(\gbigt(E))$.
Hence, by the above preliminary consideration,
$(\bigoplus_jf^j_{\dagger}(\nbigt),L_{\omega},f_{\dagger}\nbigs)$
is a polarized graded Lefschetz $\nbigd$-module on $Y$ of weight $0$
and type $1$.
Thus, we obtain $P(n,0)$ and hence Theorem \ref{thm;22.3.5.1}.
\hfill\qed

\subsubsection{Mixed case}

Let $f:X\lrarr Y$ be as in \S\ref{subsection;22.3.10.10}.
Let $(\nbigt,W)$ be a mixed twistor $\nbigd$-module on $X$.
We always assume the graded polarizability,
i.e.,
$\Gr^W_j\nbigt$ are polarizable.
We assume that the support of $\nbigt$ is proper over $Y$.
We obtain the $\nbigr$-triple $f^j_{\dagger}(\nbigt)$.
Let $W_kf_{\dagger}^j(\nbigt)$
denote the image of
$f_{\dagger}^j(W_{k-j}\nbigt)\lrarr f_{\dagger}^j(\nbigt)$.
As in the case of projective morphisms in \cite[\S7.1,\S7.2]{Mochizuki-MTM},
we obtain the following corollary from Theorem \ref{thm;22.3.5.1}.
\begin{cor}
$(f_{\dagger}^j(\nbigt),W)$ $(j\in \seisuu)$
are mixed twistor $\nbigd$-modules on $Y$.
\hfill\qed
\end{cor}

\subsection{The normal crossing singular case}
\label{subsection;22.3.10.2}

Let $X$ be a compact complex manifold 
equipped with a K\"ahler form $\omega_0$.
Let $H$ be a simply normal crossing hypersurface of $X$.
Let $(E,\delbar_E,\theta,h)$ be a tame harmonic bundle
on $(X,H)$.
We obtain the pure twistor $\nbigd_X$-module
$\gbigt(E)=(\gbige,\gbige,C_h)$ of weight $0$
with the polarization $\nbigs=(\id,\id)$.
(See \S\ref{subsection;22.4.12.1}.)
The following theorem is a special but most essential
case of Theorem \ref{thm;22.2.25.40}.

\begin{thm}
\label{thm;22.3.15.20}
$(\bigoplus a_{X\dagger}^j(\gbigt(E)),L_{\omega},a_{X\dagger}\nbigs)$
is a polarized graded Lefschetz twistor structure of weight $0$ and type $1$.
\end{thm}

\subsubsection{K\"ahler form}
\label{subsection;22.4.2.21}

As explained in \cite{cg,k3,z},
let $\varphi$ be a $C^{\infty}$-function on $X\setminus H$
such that the following holds
around any point $P$ of $H$.
\begin{itemize}
 \item Let $(X_P;z_1,\ldots,z_n)$ be a holomorphic coordinate
       neighborhood of $P$
       such that $X_P\cap H=\bigcup_{i=1}^{\ell}\{z_i=0\}$.
       Then,
       $\varphi_{|X_P\setminus H}
       =\sum_{i=1}^{\ell}\log(a_i-\log|z_i|^2)$
       for some $C^{\infty}$-functions $a_i$.
\end{itemize}
If $c>0$ is sufficiently small,
$\omega=\omega_0+\sqrt{-1} c\del\delbar\varphi$
is a K\"ahler form on $X\setminus H$
which is Poincar\'{e} like around any point of $H$.
In the following,
we shall use the associated K\"ahler metric $g_{X\setminus H}$.

\subsubsection{Harmonic forms}
\label{subsection;22.4.10.1}

For any $\lambda$,
let $\DD^{\lambda\ast}$ denote the formal adjoint of $\DDlambda$.
We set $\Delta_E:=\DD^{0\ast}\DD^0+\DD^0\DD^{0\ast}$.
Recall that
$\DD^{\lambda\ast}\DD^{\lambda}
+\DD^{\lambda}\DD^{\lambda\ast}
=(1+|\lambda|^2)\Delta_E$.
This is essentially due to Simpson,
and explained in \cite{s5}.
(See also \cite{mochi5}.)

\begin{prop}
\label{prop;22.2.25.10}
Let $\tau$ be an $L^2$-section of
$\Tot^k(E\otimes\Omega^{\bullet,\bullet})=
\bigoplus_{k_1+k_2=k}E\otimes\Omega^{k_1,k_2}$
such that $\Delta_E\tau=0$.
Then, we obtain
$\DDlambda\tau=\DD^{\lambda\ast}\tau=0$
for any $\lambda$. 
\end{prop}
\pf
We explain only an outline.
Let $\rho:\real\to\closedclosed{0}{1}$ be a $C^{\infty}$-function
such that
(i) $\rho(t)=0$ $(t\geq 2/3)$,
(ii) $\rho(t)=1$ $(t\leq 1/2)$,
(iii) $\rho^{-1/2}\del_t\rho$ is $C^{\infty}$.
Let $\varphi$ be as in \S\ref{subsection;22.4.2.21}.
For $N>0$,
we set $\chi_N=\rho\bigl(N^{-1}\varphi\bigr)$.
Because $\Delta_E\tau=0$,
we obtain
\begin{multline}
 \int |\chi_N^{1/2}\DDlambda \tau|^2
 +\int|\chi_N^{1/2}\DD^{\lambda\ast}\tau|^2
 \leq 
 \left(
 \int|\chi_N^{1/2}\DDlambda\tau|^2
 \right)^{1/2}
 \left(
 \int |\chi_N^{-1/2}(\lambda\del\chi_N+\delbar\chi_N)|^2|\tau|^2
  \right)^{1/2}  
\\
 +\left(
  \int|\chi_N^{1/2}\DD^{\lambda\ast}\tau|^2
 \right)^{1/2}
 \left(
 \int |\chi_N^{-1/2}(\del_X\chi_N-\lambdabar\delbar_X\chi_N)|^2|\tau|^2
  \right)^{1/2}.
\end{multline}
There exists $C_1>0$, which is independent of $N$,
such that
\[
 \int |\chi_N^{1/2}\DDlambda \tau|^2
 +\int|\chi_N^{1/2}\DD^{\lambda\ast}\tau|^2
 \leq
 C_1\int
 \bigl|
 \chi_N^{-1/2}(\lambda\del_X\chi_N+\delbar_X\chi_N)
 \bigr|^2|\tau|^2
 +C_1\int
 \bigl|
 \chi_N^{-1/2}(\del_X\chi_N-\lambdabar\delbar_X\chi_N)
 \bigr|^2|\tau|^2.
\]
Note that
(i) $\rho^{-1/2}\del_t\rho$ is $C^{\infty}$,
(ii) $\del\varphi$ is bounded with respect to $g_{X\setminus H}$,
(iii) $|\tau|^2$ is $L^1$.
Hence, the right hand side is convergent to $0$ as $N\to \infty$.
Hence, we obtain the claim of Proposition \ref{prop;22.2.25.10}.
\hfill\qed

\vspace{.1in}

Let $\Harm^j(E)$ denote the space of
the $L^2$-sections $\tau$ of
$\Tot^j(E\otimes \Omega^{\bullet,\bullet})$
such that $\Delta_E\tau=0$.

\begin{cor}
\label{cor;22.3.17.120}
$\Harm^j(E)$ is equal to the space of
$L^2$-sections $\tau$
of $E\otimes\Tot^j\Omega^{\bullet,\bullet}$ such that
$\DDlambda\tau=\DD^{\lambda\ast}\tau=0$.
\hfill\qed 
\end{cor}

\subsubsection{The fiber-wise isomorphism and
the finite dimensional property}

For each $\lambda$,
there exists the following isomorphism
in $\ttD(\cnum_X)$
by Corollary \ref{cor;22.3.17.100}:
\begin{equation}
 \nbigc^{\bullet}_{\tw}(\gbige,\lambda)
 \lrarr
 \nbigc^{\bullet}_{L^2}(\nbigelambda,\DDlambda,h).
\end{equation}

\begin{prop}
\label{prop;22.3.17.110}
 $H^{j}(X,\nbigc^{\bullet}_{L^2}(\nbigelambda,\DDlambda,h))$
are finite dimensional.
\end{prop}
\pf
Because each term of 
$\nbigc^{\bullet}_{\tw}(\gbige,\lambda)$
is $\nbigo_X$-coherent,
we obtain that
$H^j(X,\nbigc^{\bullet}_{\tw}(\gbige,\lambda))$
are finite dimensional for any $j\in\seisuu$.
Then, we obtain the claim of 
Proposition \ref{prop;22.3.17.110}.
\hfill\qed

\vspace{.1in}
We obtain the following corollary
from Corollary \ref{cor;22.3.17.120}
and Proposition \ref{prop;22.3.17.110}.
(See \cite{k3} or \cite{z}.)

\begin{cor}
We obtain the isomorphism
$H^{j}(X,\nbigc^{\bullet}_{L^2}(\nbigelambda,\DDlambda,h))
\simeq
\Harm^{j}(E)$
for any $j\in\seisuu$.
 In particular,
$\Harm^{j}(E)$ are finite dimensional.
\hfill\qed
\end{cor}

The multiplication of $\sqrt{-1}\omega_0$
induces a morphism
$\sqrt{-1}\omega_0:
H^j\bigl(
X,\nbigc^{\bullet}_{\tw}(\gbige,\lambda)
\bigr)
\lrarr
H^{j+2}\bigl(
X,\nbigc^{\bullet}_{\tw}(\gbige,\lambda)
\bigr)$.
The multiplication of $\sqrt{-1}\omega$
induces a morphism
$\sqrt{-1}\omega:\Harm^{j}(E)\lrarr\Harm^{j+2}(E)$.
\begin{lem}
\label{lem;22.2.25.30}
The following diagram is commutative:
\[
 \begin{CD}
H^j\bigl(X,
  \nbigc^{\bullet}_{\tw}(\gbige,\lambda)
  \bigr)
  @>{\sqrt{-1}\omega_0}>>
H^{j+2}\bigl(X,
  \nbigc^{\bullet}_{\tw}(\gbige,\lambda)
  \bigr)
  \\
 @V{\simeq}VV @V{\simeq}VV \\
 \Harm^{j}(E)
 @>{\sqrt{-1}\omega}>>
 \Harm^{j+2}(E).
 \end{CD} 
\]
\end{lem}
\pf
The following diagram is commutative:
\[
\begin{CD}
H^j\bigl(X,
  \nbigc^{\bullet}_{\tw}(\gbige,\lambda)
  \bigr)
 @>{\sqrt{-1}\omega_0}>>
H^{j+2}\bigl(X,
  \nbigc^{\bullet}_{\tw}(\gbige,\lambda)
  \bigr)
 \\
 @V{\simeq}VV @V{\simeq}VV \\
 H^{j}(X,\nbigc^{\bullet}_{L^2}(\nbigelambda,\DDlambda,h))
 @>{\sqrt{-1}\omega_0}>>
 H^{j+2}(X,\nbigc^{\bullet}_{L^2}(\nbigelambda,\DDlambda,h)).
\end{CD}
\]
Because $\del\varphi$ is bounded with respect to $g_{X\setminus H}$,
the multiplications of $\sqrt{-1}\omega$ and
$\sqrt{-1}\omega_0$
induce the same morphism
\[
H^{j}(X,\nbigc^{\bullet}_{L^2}(\nbigelambda,\DDlambda,h))
\lrarr
 H^{j+2}(X,\nbigc^{\bullet}_{L^2}(\nbigelambda,\DDlambda,h)).
\]
Then, the claim of the lemma follows.
\hfill\qed

\subsubsection{Polarized graded Lefschetz twistor structure associated with
the space of harmonic sections}
\label{subsection;22.4.2.41}

We set $d_X:=\dim X$.
By Proposition \ref{prop;22.2.25.10},
we obtain the natural embedding
\[
 \nbigo_{\cnum}\otimes_{\cnum}
 \Harm^{d_X+j}(E)
\lrarr
 (\id_{\cnum}\times a_X)_{\ast}
 \nbigc^{d_X+j}_{L^2}(\nbige,\DD,h).
\]
The isomorphism
$\nbigc^{d_X+j}_{L^2}(\nbige,\DD,h)
\simeq
\nbigctilde^{d_X+j}_{L^2}(\nbige,\DD^f,h)$
induce a morphism
\[
 \nbigo_{\cnum}\otimes_{\cnum}
 \Harm^{d_X+j}(E)
 \lrarr
 (\id_{\cnum}\times a_X)_{\ast}
 \nbigctilde^{d_X+j}_{L^2}(\nbige,\DD^f,h).
\]
Let $\nbightilde^j(E)$ denote the image.
We obtain the following morphism of complexes:
\[
 \bigoplus_j\nbightilde^j(E)[-j]
 \lrarr
  (\id_{\cnum}\times a_X)_{\ast}
 \nbigctilde^{\bullet}_{L^2}(\nbige,\DD^f,h)[d_X].
\]

We set $\vecS=\bigl\{\lambda\in\cnum\,\big|\,|\lambda|=1\bigr\}$.
Let $\nbigc_{\vecS}$ denote the sheaf of continuous functions on $\vecS$.
We obtain the sesqui-linear pairing
\[
 a_{X\dagger}(C_h):
 (\id_{\cnum}\times a_X)_{\ast}\Bigl(
 \nbigctilde^{\bullet}_{L^2}(\nbige,\DD,h)[d_X]
 \Bigr)_{|\vecS}
 \times
 \overline{
 \sigma^{-1}\Bigl(
 (\id_{\cnum}\times a_X)_{\ast}\bigl(
 \nbigctilde^{\bullet}_{L^2}(\nbige,\DD,h)[d_X]
 \bigr)
  \Bigr)_{|\vecS}}
  \lrarr
  \nbigc_{\vecS}
\]
induced by the following pairing
for
$\eta_1^{d_X-j}\otimes m_1
\in \nbigctilde^{d_X-j}_{L^2}(\nbige,\DD,h)$
and 
$\eta_2^{d_X+j}\otimes m_2
\in \nbigctilde^{d_X+j}_{L^2}(\nbige,\DD,h)$:
\[
 a_{X\dagger}C_h(\tau_1,\tau_2)
 =\frac{(-1)^{(d_X+j)(d_X+j-1)/2}}{(2\pi\sqrt{-1})^{d_X}}
 \int_X h(m_1,\sigma^{\ast}(m_2))
 \eta_1^{d_X-j}
 \cdot
 \overline{\sigma^{\ast}(\eta_2^{d_X+j})}.
\]
Let $C_{\nbightilde}:
\nbightilde^{-j}(E)_{|\vecS}\times
\sigma^{-1}\bigl(
 \nbightilde^j(E)_{|\vecS}
 \bigr)
 \lrarr \nbigc_{\vecS}$
denote the induced sesqui-linear pairing.
We set
\[
\nbightilde^j(\gbigt(E)):=
(\nbightilde^{-j}(E),\nbightilde^j(E),C_{\nbightilde}).
\]
The multiplication of $\sqrt{-1}\omega$ induces
$\Harm^j(E)\lrarr
 \Harm^{j+2}(E)$,
which induces 
\[
\sqrt{-1}\omega:
\nbightilde^j(E)
 \lrarr
 \lambda\cdot \nbightilde^{j+2}(E).
\]
We obtain the morphism of
$\nbigr$-triples
$L_{\omega}:\nbightilde^j(\gbigt(E))
\lrarr
\nbightilde^{j+2}(\gbigt(E))\otimes\newTate(1)$.
The pair of the identity morphisms induces
the Hermitian sesqui-linear duality
$a_{X\dagger}\nbigs:
\nbightilde^j(\gbigt(E))
\lrarr
\nbightilde^{-j}(\gbigt(E))^{\ast}$.
We obtain the following proposition
as in the smooth case 
\cite[Hodge-Simpson Theorem 2.24]{sabbah2}.

\begin{prop}
\label{prop;22.3.1.2}
$(\bigoplus_{j}\nbightilde^j(\gbigt(E)),L_{\omega})$
is a graded Lefschetz twistor structure of weight $0$ and type $1$,
and $a_{X\dagger}\nbigs$ is a polarization.
\hfill\qed
\end{prop}

\begin{rem}
The signature of $L_{\omega}$ is opposite to that in
{\rm\cite{sabbah2, Mochizuki-MTM}}.
We impose the condition that
$(L_{\omega}^{\ast})^j\circ a_{X\dagger}(\nbigs)$
is a polarization on $P\nbightilde^{-j}(\gbigt(E))$
$(j\geq 0)$
instead of
$a_{X\dagger}(\nbigs)\circ L_{\omega}^j$
in {\rm\cite{sabbah2}}.
Hence, the conditions are equivalent.
\hfill\qed
\end{rem}

\subsubsection{The isomorphisms in family}

By Corollary \ref{cor;22.3.17.100},
there exists the following isomorphism
in $\ttD((\nbigo_{\cnum_{\lambda}})_{\nbigx})$:
\begin{equation}
\label{eq;22.2.25.3}
 \nbigctilde^{\bullet}_{\tw}(\gbige)
\simeq
 \nbigctilde^{\bullet}_{L^2}(\nbige,\DD^f,h).
\end{equation}
It induces the isomorphism of
$\nbigo_{\cnum_{\lambda}}$-modules:
\begin{equation}
\label{eq;22.2.25.1}
 a_{X\dagger}^{j}(\gbige)
 \simeq
 R^{j}\bigl(\id_{\cnum}\times a_X\bigr)_{\ast}
 \Bigl(
 \nbigctilde^{\bullet}_{L^2}(\nbige,\DD^f,h)[d_X]
 \Bigr).
\end{equation}

\begin{prop}
The induced morphisms
\begin{equation}
\label{eq;22.2.25.22}
 \nu^{j}:
 \nbightilde^{j}(E)
 \lrarr
  a_{X\dagger}^{j}(\gbige)
\end{equation}
are isomorphisms.
\end{prop}
\pf
We use a descending induction on $j$.
Suppose that we have already proved the claim for $j\geq j_0+1$.
Let $\iota_{\lambda}:\{\lambda\}\lrarr\cnum$
denote the inclusion.
Under the assumption,
we have
\[
 \iota_{\lambda}^{\ast}
 R^{j_0}(\id_{\cnum}\times a_X)_{\ast}
 \nbigc^{\bullet}_{\tw}(\gbige)[d_X]
 \simeq
 H^{j_0}(X,\nbigc^{\bullet}_{\tw}(\gbige,\lambda)[d_X]).
\]
Hence,
the following induced morphisms are isomorphisms for any $\lambda$:
\[
 \nbightilde^{j_0}(E)_{|\lambda}
 \lrarr
 \iota_{\lambda}^{\ast}
  a_{X\dagger}^{j_0}(\gbige).
\]
Then, it is easy to obtain that
(\ref{eq;22.2.25.22}) for $j_0$
is an isomorphism.
\hfill\qed

\begin{lem}
The following diagram is commutative.
\[
 \begin{CD}
 \nbightilde^{j}(E)
  @>{\sqrt{-1}\omega}>>
 \lambda\cdot \nbightilde^{j+2}(E)
  \\
  @V{\simeq}VV @V{\simeq}VV \\
  a_{X\dagger}^{j}(\gbige)
  @>{\sqrt{-1}\omega_0}>>
  \lambda \cdot a_{X\dagger}^{j+2}(\gbige).
 \end{CD}
\] 
\end{lem}
\pf
It follows from Lemma \ref{lem;22.2.25.30}.
\hfill\qed

\vspace{.1in}
Let $L_{\omega_0}:
a_{X\dagger}^j\gbigt(E)
\lrarr
a_{X\dagger}^{j+2}\gbigt(E)
\otimes\newTate(1)$
denote the morphism induced by
$(\sqrt{-1}\omega_0,\sqrt{-1}\omega_0)$.

\begin{lem}
\label{lem;22.3.1.1}
The pair of the morphisms
$(\nu^{-j})^{-1}:
 a_{X\dagger}^{-j}(\gbige)
 \simeq
 \nbigh^{-j}(E)$
and  
$\nu^j:
\nbigh^{j}(E)\simeq
a_{X\dagger}^{j}(\gbige)$
induces an isomorphism of
the $\nbigr$-triples
$\nu:\nbightilde^j(\gbigt(E))
\simeq
a_{X\dagger}^{j}(\gbigt(E))$
such that
$L_{\omega_0}\circ\nu
=\nu\circ L_{\omega}$. 
\end{lem}
\pf
Let us compare the sesqui-linear pairing 
of
$\nbightilde^{-j}(E)_{|\lambda}$
and
$\overline{\nbightilde^{j}(E)_{|-\lambda}}$,
and the sesqui-linear pairing of
$a_X^{-j}(\gbige)_{|\lambda}$
and
$\overline{a_X^{j}(\gbige)_{|-\lambda}}$
for general $\lambda$.
The complexes
$\nbigc^{\bullet}_{\tw}(\gbige,\lambda)$,
$\nbigc^{\bullet}_{\tw}(\gbige,-\lambda)$,
$\nbigc^{\bullet}_{L^2}(\nbigelambda,\DDlambda,h)$
and
$\nbigc^{\bullet}_{L^2}(\nbige^{-\lambda},\DDlambda,h)$
are the intersection complexes.
The restriction of the sesqui-linear pairings of
$\nbigc^{\bullet}_{\tw}(\gbige,\lambda)[d_X]$,
and $\overline{\nbigc^{\bullet}_{\tw}(\gbige,-\lambda)[d_X]}$
to $X\setminus H$
is equal to the restriction of
the sesqui-linear pairings of
$\nbigc^{\bullet}_{L^2}(\nbigelambda,\DDlambda,h)[d_X]$
and 
$\overline{\nbigc^{\bullet}_{L^2}(\nbige^{-\lambda},\DDlambda,h)[d_X]}$.
Hence, the sesqui-linear pairings are equal on $X$.
Then, the claim of Lemma \ref{lem;22.3.1.1} follows.
\hfill\qed

\vspace{.1in}
We obtain Theorem \ref{thm;22.3.15.20}
for $\gbigt(E)$ from Lemma \ref{lem;22.3.1.1}
and Proposition \ref{prop;22.3.1.2}.
\hfill\qed

\subsection{Reduction to the normal crossing singular case}
\label{subsection;22.4.28.1}

\subsubsection{Setting}

Let $X$ be a K\"ahler manifold
with a K\"ahler class $\omega_1$.
Let $f:\Xtilde\lrarr X$ be a projective birational morphism.
We assume that $f$ is the composition of
a finite sequence of blowing-ups along smooth complex submanifolds.
There exists an effective divisor $\Htilde$ such that
(i) $f$ induces an isomorphism
$\Xtilde\setminus|\Htilde|\simeq X\setminus f(|\Htilde|)$,
(ii) $\nbigo_{\Xtilde}(-\Htilde)$ is relatively ample
with respect to $f$.
Let $\omega_2$ denote the first Chern class of
$\nbigo_{\Xtilde}(-\Htilde)$.
The pull back $f^{\ast}\omega_1$ is also denoted by $\omega_1$.

Let $\nbigt$ be a pure twistor $\nbigd$-module
of weight $w$ on $\Xtilde$ with a polarization $\nbigs$.
We assume the following on the support $\supp(\nbigt)$
of $\nbigt$.
\begin{itemize}
 \item $\supp(\nbigt)$ is
       a connected complex closed submanifold of $\Xtilde$,
       and $\supp(\nbigt)\cap |\Htilde|$ is
       a simple normal crossing hypersurface of
       $\supp(\nbigt)$.
\end{itemize}

As the hypothesis of the induction,
we assume the following.
\begin{assumption}
\label{assumption;22.3.10.1}
 The claim of Theorem {\rm\ref{thm;22.2.25.40}} holds
for a polarized pure twistor $\nbigd$-module $\nbigt'$ on $X$ and $\Xtilde$
such that the dimension of the support of $\nbigt'$
is strictly smaller than $\dim\supp(\nbigt)$.
\end{assumption}

\subsubsection{Statement}

By the Hard Lefschetz Theorem for the push-forward of
a polarized pure twistor $\nbigd$-module via a projective morphism
in \cite{sabbah2,mochi2,Mochizuki-wild},
$f_{\dagger}^0\nbigt$ is a polarizable
pure twistor $\nbigd$-module of weight $w$ on $X$.
There exists a unique decomposition
\[
 f_{\dagger}^0\nbigt
 =\nbigt_X\oplus\nbigt_X',
\]
where the dimension of the strict support of $\nbigt_X$
is equal to $\dim\supp(\nbigt)$,
and the dimension of the support of $\nbigt_X'$
is strictly smaller than $\dim\supp(\nbigt)$.
Note that the dimension of the support of
$f_{\dagger}^j(\nbigt)$ $(j\neq 0)$ are strictly smaller
than $\dim\supp(\nbigt)$.
Hence, $\nbigt_X$ is contained
in the primitive part of $f_{\dagger}^0\nbigt$,
i.e., the kernel of
$L_{\omega_2}:f_{\dagger}^0\nbigt\lrarr f_{\dagger}^2\nbigt$.
Hence, $f_{\dagger}\nbigs$ induces a polarization
$\nbigs_X$ of $\nbigt_X$.
We shall prove the claim of Theorem \ref{thm;22.2.25.40}
for $(\nbigt_X,\nbigs_X)$
under the hypothesis of the induction.

\begin{prop}
\label{prop;22.3.4.30}
Under Assumption {\rm\ref{assumption;22.3.10.1}},
 $\Bigl(
 \bigoplus a_{X\dagger}^j\nbigt_X,
 L_{\omega_1},
 a_{X\dagger}(\nbigs_X)
 \Bigr)$
 is a polarized Lefschetz graded twistor structure of weight $w$
 and type $1$.
\end{prop}

\subsubsection{The associated polarized graded Lefschetz twistor structure}

As explained in \cite[\S2]{MSaito-proper},
there exist
$c_0>0$ and a relatively compact neighbourhood $U$ of
$f(\supp(\nbigt))$ in $X$
such that
$(\omega_1+c\omega_2)_{|f^{-1}(U)}$ is a K\"ahler class
for any $0<c\leq c_0$.
By replacing $X$ with $U$,
we may assume that $\omega_1+c\omega_2$ is a K\"ahler class
for any $0<c\leq c_0$ on $\Xtilde$.

In \S\ref{subsection;22.3.10.2},
we have already proved the following claim.
\begin{itemize}
 \item $\bigl(\bigoplus_ja_{\Xtilde\dagger}^j(\nbigt),
       L_{\omega_1+c\omega_2},
       a_{\Xtilde\dagger}\nbigs\bigr)$
       is a polarized graded Lefschetz twistor structure of weight $w$
       and type $1$
       for any $0<c\leq c_0$.
\end{itemize}

We set
$\abar_{\Xtilde\dagger}^j(\nbigt):=
 a_{\Xtilde\dagger}^j(\nbigt)\otimes\newTate(j)$.
We obtain the induced morphisms
$\Lbar_{\omega_i}:
\abar_{\Xtilde\dagger}^j(\nbigt)
\lrarr
\abar_{\Xtilde\dagger}^{j+2}(\nbigt)\otimes\newTate(-1)$
$(i=1,2)$.
Let $\abar_{X\dagger}\nbigs$
denote the induced Hermitian adjoint of weight $w$.
Then,
\[
\Bigl(
 \bigoplus_{j}\abar_{\Xtilde\dagger}^j(\nbigt),
 L_{\omega_1+c\omega_2},
 \abar_{X\dagger}\nbigs
\Bigr)
\]
is a polarized graded Lefschetz twistor structure of weight $w$
and type $-1$
for any $c\leq c_0$.
It implies that
\[
\Bigl(
 \bigoplus_{j}\abar_{\Xtilde\dagger}^j(\nbigt),
 L_{\omega_1},L_{\omega_2+c_0^{-1}\omega_1},
  \abar_{X\dagger}\nbigs
\Bigr)
\]
is a polarized mixed twistor structure of weight $w$.

\subsubsection{The induced graded Lefschetz twistor $\nbigd$-module on $X$}

By \cite{sabbah2, mochi2, Mochizuki-wild},
$\Bigl(
\bigoplus_j f^{j}_{\dagger}(\nbigt),
L_{\omega_2},f_{\dagger}\nbigs
\Bigr)$
is a polarized graded Lefschetz twistor $\nbigd$-module on $X$
of weight $w$ and type $1$.
We set
\[
 \fbar^j_{\dagger}(\nbigt):=
 f^j_{\dagger}(\nbigt)
 \otimes\newTate(j).
\]
Then,
$\Bigl(
 \bigoplus \fbar^j_{\dagger}(\nbigt),
 L_{\omega_2},
 \fbar_{\dagger}(\nbigs)
\Bigr)$
is a polarized graded Lefschetz twistor $\nbigd$-module
of weight $w$ and type $-1$.
We also set
\[
 \abar_{X\dagger}^k(\fbar^j_{\dagger}(\nbigt)):=
 a_{X\dagger}^k(\fbar^j_{\dagger}(\nbigt))
 \otimes\newTate(k)
=a_{X\dagger}^k(f^j_{\dagger}(\nbigt))
 \otimes\newTate(k+j).
\]

Because there exists an isomorphism
$f_{\dagger}(\nbigt)
 \simeq
 \bigoplus
 f^j_{\dagger}(\nbigt)[-j]$
in the derived category of $\nbigr_X$-triples,
we obtain the Leray filtration
$L(\abar_{\Xtilde\dagger}^i(\nbigt))$,
which is a decreasing filtration indexed by integers
with a canonical isomorphism
\[
 \Gr_L^k(\abar_{\Xtilde\dagger}^i(\nbigt))
 \simeq
 \abar_{X\dagger}^{k}
 \fbar_{\dagger}^{i-k}(\nbigt).
\]
We also note that 
$\abar_{X\dagger}^{k}
\fbar_{\dagger}^{i-k}(\nbigt)$
is a pure twistor structure of weight $w-i$
because it is isomorphic to a direct summand of
the pure twistor structure
$\abar_{\Xtilde\dagger}^i(\nbigt)$
of weight $w-i$. 

\subsubsection{The induced pure twistor $\nbigd$-modules}

We have the mixed twistor $\nbigd$-modules
$\nbigt[\ast \Htilde]$
and $\nbigt[!\Htilde]$.
We set
\[
\nbigt_1:=W_{w+1}(\nbigt[\ast \Htilde]),
\quad
\nbigt_2:=\nbigt[!\Htilde]\big/W_{w-2}(\nbigt[!\Htilde]).
\]
We set
\[
\nbigk:=\Gr^W_{w-1}(\nbigt[!\Htilde]),\quad
\nbigc:=\Gr^{W}_{w+1}(\nbigt[\ast\Htilde]).
\]
The polarization
$\nbigs:\nbigt\lrarr\nbigt^{\ast}\otimes\newTate(-w)$
induces the following commutative diagram:
\[
 \begin{CD}
  0 @>>>
  \nbigk @>>>\nbigt_2 @>>>
  \nbigt @>>> 0 \\
  @VVV @V{\nbigs_1}VV @VVV @V{\nbigs}VV @VVV \\
  0@>>>
  \nbigc^{\ast}\otimes\newTate(-w)
  @>>>
  \nbigt_1^{\ast}\otimes\newTate(-w)
  @>>>
  \nbigt^{\ast}\otimes\newTate(-w)
  @>>> 0.
 \end{CD}
\]
There also exists the following induced commutative diagram:
\[
  \begin{CD}
  0 @>>>
  \nbigt @>>>\nbigt_1 @>>>
  \nbigc @>>> 0 \\
  @VVV @V{\nbigs}VV @VVV @V{\nbigs_2}VV @VVV \\
  0@>>>
  \nbigt^{\ast}\otimes\newTate(-w)
   @>>>
  \nbigt_2^{\ast}\otimes\newTate(-w)   
  @>>>
   \nbigk^{\ast}\otimes\newTate(-w)
  @>>> 0.
 \end{CD}
\]

\begin{lem}
\label{lem;22.4.20.1}
There exists an isomorphism
$\upsilon:\nbigk\simeq \nbigc\otimes\newTate(1)$
such that
\[
\nbigs_{\nbigk}:=\widetilde{\upsilon}^{\ast}\circ\nbigs_1
 =-\nbigs_2\circ\upsilon:
\nbigk\lrarr\nbigk^{\ast}\otimes\newTate(-w+1)
\]
is a polarization of $\nbigk$.
Here, $\widetilde{\upsilon}^{\ast}$ is induced by
 $\upsilon^{\ast}:
 \nbigc^{\ast}
 \simeq
 (\nbigk\otimes\newTate(-1))^{\ast}$ 
 and
$\iota_{\newTate(-1)}=(-1,-1):\newTate(-1)^{\ast}\simeq\newTate(1)$. 
\end{lem}
\pf
Let $\pi:\Xhat\to X$ denote the holomorphic line bundle
such that the sheaf of holomorphic sections of $\Xhat$
is $\nbigo(\Htilde)$.
Let $\iota_0:\Xtilde\to \Xhat$ denote the $0$-section,
and let $\iota_1:\Xtilde\to \Xhat$ denote the section
induced by $1:\nbigo_{\Xtilde}\to \nbigo_{\Xtilde}(\Htilde)$.
We set
$(\nbigthat,\nbigshat):=
\iota_{1\dagger}(\nbigt,\nbigs)$,
which is a polarized pure twistor $\nbigd$-module on $\Xhat$
of weight $w$.
Let $\Hhat$ denote $\iota_0(\Xtilde)$.
We may regard it as a reduced effective divisor.
By \cite[Proposition 11.2.7]{Mochizuki-MTM},
we have
$\nbigthat[\star\Hhat]=
\iota_{1\dagger}(\nbigt[\star \Htilde])$
$(\star=!,\ast)$.
We set
$\nbigthat_1:=W_{w+1}(\nbigthat[\ast \Hhat])
\simeq
 \iota_{1\dagger}(\nbigt_1)$
and
$\nbigthat_2:=\nbigthat[!\Hhat]\big/
W_{w-2}(\nbigthat[! \Hhat])
\simeq
\iota_{1\dagger}(\nbigt_2)$.
Similarly, we set
$\nbigkhat:=
\iota_{1\dagger}(\nbigk)=\Gr^W_{w-1}(\nbigthat)$
and
$\nbigchat:=\iota_{1\dagger}(\nbigc)
=\Gr^W_{w+1}(\nbigthat)$.
We obtain
$\nbigshat_1:=
\iota_{1\dagger}(\nbigs_1):
\nbigkhat\to
\nbigchat^{\ast}\otimes\newTate(-w)$
and
$\nbigshat_2:=
\iota_{1\dagger}(\nbigs_2):
\nbigchat\to
\nbigkhat^{\ast}\otimes\newTate(-w)$.

Let $P$ be any point of $\Xtilde$.
Let $\Xtilde_P$ be an open neighbourhood of $P$ in $\Xtilde$
with a holomorphic trivialization
$\Xhat_P:=\pi^{-1}(\Xtilde_P)\simeq\cnum\times\Xtilde_P$.
Let $z$ be the holomorphic function on $\Xhat_P$
induced by the projection
$\Xhat_P\simeq \cnum\times\Xtilde_P\lrarr \cnum$.
We have the holomorphic vector field
$\del_z$ on $\Xhat_P\simeq\cnum\times\Xtilde_P$.
We set $\Hhat_P:=\Xhat_P\cap\Hhat$,
and let $\iota_{\Hhat_P}:\Hhat_P\to \Xhat_P$
denote the inclusion.
We set $\nbigthat_P:=\nbigthat_{|\Xhat_P}$.
We obtain the mixed twistor $\nbigd$-module
$\psitilde_{z,-\vecdelta}(\nbigthat_P)$
on $\Hhat_P$
with the induced morphism
$\nbign:\psitilde_{z,-\vecdelta}(\nbigthat_P)
\to\psitilde_{z,-\vecdelta}(\nbigthat_P)\otimes\newTate(-1)$.
Let $M(\nbign)$ denote the monodromy weight filtration
of $\nbign$ on $\psitilde_{z,-\vecdelta}(\nbigthat_P)$.
Let $P\Gr^{M(\nbign)}_0\psitilde_{z,-\vecdelta}(\nbigthat_P)$
denote the primitive part of
$\Gr^{M(\nbign)}_0\psitilde_{z,-\vecdelta}(\nbigthat_P)$.

In \cite[\S4.3.2, \S4.3.3]{Mochizuki-MTM}
we constructed isomorphisms
\begin{equation}
\label{eq;22.4.18.1}
 \iota_{\Hhat_P\dagger}\Ker\nbign\otimes\nbigu(-1,0)
 \simeq
 \Ker(\nbigthat[!\Hhat]\to\nbigthat[\ast\Hhat])_{|\Xhat_P},
\end{equation}
\begin{equation}
\label{eq;22.4.18.2}
 \iota_{\Hhat_P\dagger}\Cok\nbign\otimes\nbigu(-1,0)
 \simeq
\Cok(\nbigthat[!\Hhat]\to\nbigthat[\ast\Hhat])_{|\Xhat_P}.
\end{equation}
Let us recall the construction.
Let $\nbigmhat'_P$, $\nbigmhat''_P$
denote the $\nbigr_{\Xhat_P}$-module
underlying $\nbigthat_P$,
i.e., there exists a sesqui-linear pairing $C_P$ of
$\nbigmhat'_P$ and $\nbigmhat''_P$
such that 
$\nbigthat_P$ is the $\nbigr_{\Xhat_P}$-triple
$(\nbigmhat'_P,\nbigmhat''_P,C_P)$.
The $\nbigr_{\Xhat_P}$-modules underlying $\nbigthat_P[!\Hhat_P]$
are 
$\nbigmhat'_P[\ast \Hhat_P]$
and
$\nbigmhat''_P[!\Hhat_P]$.
The $\nbigr_{\Xhat_P}$-modules underlying $\nbigthat_P[\ast\Hhat_P]$
are 
$\nbigmhat'_P[!\Hhat_P]$
and
$\nbigmhat''_P[\ast\Hhat_P]$.
For $\nbigm=\nbigmhat'_P,\nbigmhat''_P$,
let $\Ker(\nbigm)$ and $\Cok(\nbigm)$
denote the kernel and the cokernel of the natural morphism
$\nbigm[!\Hhat_P]\to\nbigm[\ast\Hhat_P]$,
respectively.
Then,
the underlying $\nbigr_{\Xhat_P}$-modules of
$\Ker\bigl(
 \nbigthat[!\Hhat]\to
 \nbigthat[\ast \Hhat]
 \bigr)_{|\Xhat_P}$
are $\Cok(\nbigmhat'_P)$
and $\Ker(\nbigmhat''_P)$,
and 
the underlying $\nbigr_{\Xhat_P}$-modules of
$\Cok\bigl(
 \nbigthat[!\Hhat]\to
 \nbigthat[\ast \Hhat]
 \bigr)_{|\Xhat_P}$
are $\Ker(\nbigmhat'_P)$
and $\Cok(\nbigmhat''_P)$.

Let $\lambda_0\in\cnum$.
Let $\nbigm=\nbigmhat'_P,\nbigmhat_P''$.
There exists a neighbourhood 
$\nbigxhat^{(\lambda_0)}_P$
of $\{\lambda_0\}\times\Xhat_P$ in $\{\lambda_0\}\times\cnum$
such that
$\nbigm_{|\nbigxhat^{(\lambda_0)}_P}$,
$\nbigm[\star \Hhat_P]_{|\nbigxhat^{(\lambda_0)}_P}$
$(\star=!,\ast)$,
$\Ker(\nbigm)_{|\nbigxhat^{(\lambda_0)}_P}$
and $\Cok(\nbigm)_{|\nbigxhat^{(\lambda_0)}_P}$
have the $V$-filtrations along $z$.
(See \S\ref{subsection;22.3.28.10} for $V$-filtrations
of $\nbigr$-modules.)
We set $\nbigmzero=\nbigm_{|\nbigxhat^{(\lambda_0)}_P}$,
$\nbigmzero[\star \Hhat_P]=
\nbigm[\star \Hhat_P]_{|\nbigxhat^{(\lambda_0)}_P}$,
$\Ker(\nbigmzero)=\Ker(\nbigm)_{|\nbigxhat^{(\lambda_0)}_P}$
and
$\Cok(\nbigmzero)=\Cok(\nbigm)_{|\nbigxhat^{(\lambda_0)}_P}$.
For any $a<0$, we have
\[
 \Vzero_{a}(\nbigmzero)
= \Vzero_{a}\bigl(
  \nbigmzero[\star \Hhat_P]
  \bigr).
\]
We have
\[
 \nbigmzero[\ast\Hhat_P]
 =
 \nbigr_{\Xhat_P|\nbigxhat^{(\lambda_0)}_P}
 \otimes_{V\nbigr_{\Xhat_P|\nbigxhat^{(\lambda_0)}_P}}
 \Vzero_{0}(\nbigmzero),
\quad
 \nbigmzero[!\Hhat_P]
=
\nbigr_{\Xhat_P|\nbigxhat^{(\lambda_0)}_P}
 \otimes_{V\nbigr_{\Xhat_P|\nbigxhat^{(\lambda_0)}_P}}
 \Vzero_{<0}(\nbigmzero). 
\]
We have the following commutative diagram:
\[
  \begin{CD}
   @.
 \Gr^{\Vzero}_{-1}(\nbigmzero)
 @>{-z\del_z}>>
   \lambda^{-1}\Gr^{\Vzero}_{-1}(\nbigmzero)
   @.
   \\
   @. @V{-\lambda\del_z}V{\simeq}V @A{\lambda^{-1}z}A{\simeq}A @.\\
 \Vzero_0(\Ker(\nbigmzero))
@>>>   
  \Gr^{\Vzero}_{0}(\nbigmzero[!\Hhat_P])
@>>>  
  \Gr^{\Vzero}_{0}(\nbigmzero[\ast\Hhat_P])
   @>>>
 \Vzero_0(\Cok(\nbigmzero)).
  \end{CD}
\]
Here, the lower horizontal morphisms give an exact sequence.
We have the decomposition
\[
 \Gr^{\Vzero}_{-1}(\nbigmzero)
 =\bigoplus_{
 \substack{
 u\in\real\times\cnum\\
 \paramap(\lambda_0,u)=-1
 }}
 \psizero_{z,u}(\nbigmzero),
\]
where the decomposition is preserved by
the action of $-\lambda z\del_z$,
and $-\lambda z\del_z+\eigenmap(\lambda,u)$
are nilpotent on
$\psizero_{z,u}(\nbigmzero)$.
Let $\psitilde_{z,-\vecdelta}(\nbigmhat'_P)$
and $\psitilde_{z,-\vecdelta}(\nbigmhat''_P)$
denote the $\nbigr_{\Hhat_P}$-modules
underlying
$\psitilde_{z,-\vecdelta}(\nbigthat_P)$.
Then, we have
\[
\psitilde_{z,-\vecdelta}(\nbigm)_{|\nbighhat_P^{(\lambda_0)}}
=\psizero_{z,-\vecdelta}(\nbigmzero),
\]
and the morphism
$\nbign:
\psitilde_{z,-\vecdelta}(\nbigthat_P)
\to
\psitilde_{z,-\vecdelta}(\nbigthat_P)\otimes\newTate(-1)$
is induced by $-z\del_z$.
The multiplication of $-\lambda\del_z$
induces
\begin{equation}
\label{eq;22.4.19.10}
 \Ker\bigl(
 -z\del_z:
 \psizero_{z,-\vecdelta}(\nbigmzero)
 \to
 \lambda^{-1}\psizero_{z,-\vecdelta}(\nbigmzero)
 \bigr)
 \simeq
 \Vzero_0(\Ker(\nbigmzero)).
\end{equation}
The multiplication of $\lambda^{-1}z$
induces
\begin{equation}
\label{eq;22.4.19.11}
 \Vzero_0(\Cok(\nbigmzero))
 \simeq
 \Cok\bigl(
 -z\del_z:
 \psizero_{z,-\vecdelta}(\nbigmzero)
 \to
 \lambda^{-1}\psizero_{z,-\vecdelta}(\nbigmzero)
 \bigr). 
\end{equation}
We obtain
the isomorphisms
\begin{equation}
\label{eq;22.4.19.12}
 \iota_{\Hhat_P\dagger}
\Ker\bigl(
-z\del_z:
\psitilde_{z,-\vecdelta}(\nbigm)\to
\lambda^{-1}\psitilde_{z,-\vecdelta}(\nbigm)
\bigr)
\simeq
 \Ker(\nbigm),
\end{equation}
\begin{equation}
\label{eq;22.4.19.13}
\Cok(\nbigm)
\simeq
\lambda^{-1}\iota_{\Hhat_P\dagger}
\Cok\bigl(
-z\del_z:
\psitilde_{z,-\vecdelta}(\nbigm)
\to
\lambda^{-1}\psitilde_{z,-\vecdelta}(\nbigm)
\bigr)
\end{equation}
from (\ref{eq;22.4.19.10}) and (\ref{eq;22.4.19.11}),
respectively.
We obtain 
(\ref{eq;22.4.18.1}) and (\ref{eq;22.4.18.2})
from the isomorphisms 
(\ref{eq;22.4.19.12})
and
(\ref{eq;22.4.19.13})
for $\nbigm=\nbigmhat_P',\nbigmhat''_P$.

From (\ref{eq;22.4.18.1}) and (\ref{eq;22.4.18.2}),
we obtain the following isomorphisms,
respectively:
\begin{equation}
\label{eq;22.4.19.3}
 \iota_{\Hhat_P\dagger}P\Gr^{M(\nbign)}_0
  \psitilde_{z,-\vecdelta}(\nbigthat_P)
  \otimes\nbigu(-1,0)
 \simeq
 \nbigkhat_{|\Xhat_P},
\end{equation}
\begin{equation}
 \iota_{\Hhat_P\dagger}P\Gr^{M(\nbign)}_0
 \psitilde_{z,-\vecdelta}(\nbigthat_P)
 \otimes\nbigu(0,-1)
 \simeq
\nbigchat_{|\Xhat_P}.
\end{equation}
We obtain the isomorphism
\begin{equation}
 \label{eq;22.4.18.3}
\widehat{\upsilon}_P:
\nbigkhat_{|\Xhat_P}\simeq \nbigchat_{|\Xhat_P}\otimes\newTate(1).
\end{equation}
Note that $\psitilde_{z,-\vecdelta}(\nbigshat)$
induces a polarization
$\nbigs_0$ on
$P\Gr^{M(\nbign)}_0\psitilde_{z,-\vecdelta}(\nbigthat_{P})$.
We obtain the polarization
$\nbigs_0\otimes\iota_{\nbigu(-1,0)}$
of
$P\Gr^{M(\nbign)}_0\psitilde_{z,-\vecdelta}(\nbigthat_{P})
\otimes\nbigu(-1,0)$.
(See Example \ref{example;22.4.19.1} for $\iota_{\nbigu(-1,0)}$.)
By the construction of the morphisms,
the following diagram is commutative:
\[
 \begin{CD}
P\Gr^{M(\nbign)}_0\psitilde_{z,-\vecdelta}(\nbigthat_{P})
  \otimes\nbigu(-1,0)
  @>{\simeq}>> \nbigkhat_{|\Xhat_P}
  \\
  @V{\nbigs_0\otimes\iota_{\nbigu(-1,0)}}VV
  @V{-\nbigshat_{2|\Xhat_P}\circ\widehat{\upsilon}_P}VV \\
\bigl(
 P\Gr^{M(\nbign)}_0\psitilde_{z,-\vecdelta}(\nbigthat_{P})
  \otimes\nbigu(-1,0)
  \bigr)^{\ast}
  \otimes\newTate(-w+1)
  @>{\simeq}>>
  \bigl(
  \nbigkhat_{|\Xhat_P}
  \bigr)^{\ast}\otimes\newTate(-w+1).
 \end{CD}
\]
Here, the horizontal isomorphisms are induced by
(\ref{eq;22.4.19.3}).
Hence,
$-\nbigshat_{2|\Xhat_P}\circ\widehat{\upsilon}_P$
is a polarization of
$\nbigkhat_{|\Xhat_P}$.

Let us check that 
the induced isomorphism $\upsilonhat_P$
is independent of the choice of the trivialization
$\Xhat_P\simeq \cnum\times\Xtilde_P$.
Let $\nbigc(\nbigmhat'_P)$
and $\nbigk(\nbigmhat''_P)$
denote the $\nbigr_{\Xhat_P}$-modules
underlying $\nbigkhat_P=\nbigkhat_{|\Xhat_P}$.
Let $\nbigk(\nbigmhat'_P)$
and $\nbigc(\nbigmhat''_P)$
denote the $\nbigr_{\Xhat_P}$-modules
underlying $\nbigchat_P=\nbigchat_{|\Xhat_P}$.
For $\nbigm=\nbigmhat'_P,\nbigmhat''_P$,
there exist the natural inclusion
$\nbigc(\nbigm)\to \Cok(\nbigm)$
and the projection
$\Ker(\nbigm)\to\nbigk(\nbigm)$.
Let us describe
the induced isomorphism
\begin{equation}
\upsilonhat_P^{-1}:
 \nbigc(\nbigm)
  \simeq
 \lambda^{-1}\nbigk(\nbigm)
\end{equation}
more explicitly.
The restrictions
$\nbigc(\nbigm)_{|\nbigxhat^{(\lambda_0)}_P}$
and
$\nbigk(\nbigm)_{|\nbigxhat^{(\lambda_0)}_P}$
are denoted by
$\nbigc(\nbigmzero)$
and
$\nbigk(\nbigmzero)$,
respectively.
Let $M(\nbign)$ denote the monodromy weight filtration of
$-\lambda\del_zz$
on
$\psizero_{z,0}(\nbigmzero[\ast \Hhat_P])$.
The primitive part of 
$\Gr^{M(\nbign)}_{0}\psizero_{z,0}
(\nbigmzero[\ast \Hhat_P])$
is denoted by
$P\Gr^{M(\nbign)}_{0}\psizero_{z,0}(\nbigmzero[\ast \Hhat_P])$.
The natural morphism
\[
 \Gr^{\Vzero}_0(\nbigmzero[\ast \Hhat_P])
 \lrarr
 \Vzero_0(\Cok(\nbigmzero))
\]
induces the following isomorphism:
\[
P\Gr^{M(\nbign)}_{0}\psizero_{z,0}(\nbigmzero[\ast \Hhat_P])
\simeq
 \Vzero_0(\nbigc(\nbigmzero)).
\]
Let $s$ be a local section
of $\Vzero_0\nbigc(\nbigmzero)$.
We may assume that there exists
a section $\stilde$ of
$\Vzero_0(\nbigmzero[\ast \Hhat_P])$
such that
(i) $\stilde$ induces $s$,
(ii) $-\lambda\del_z(z\stilde)\in
\Vzero_{<0}(\nbigmzero[\ast \Hhat_P])$.
We may naturally regard
$\lambda^{-1}z\stilde$ and
$-\del_z(z\stilde)$
as sections of
$\lambda^{-1}\Vzero_{-1}(\nbigmzero)$
and
$\lambda^{-1}\Vzero_{<0}(\nbigmzero)$,
respectively.
We obtain the section
$\bigl(
 -\del_z\otimes z\stilde
 \bigr)
 +\bigl(
 1\otimes(\del_zz\stilde)
 \bigr)$
 of
\[
\lambda^{-1}\Vzero_0\Ker(\nbigmzero)
\subset
\lambda^{-1}\Vzero_0(\nbigmzero[!\Hhat_P]).
\]
By the construction,
the section
$\upsilonhat_P^{-1}(s)$
is the image of
$\bigl(-\del_z\otimes z\stilde\bigr)
+\bigl(1\otimes (\del_zz\stilde)\bigr)$
by the projection
\[
\lambda^{-1}\Vzero_0\Ker(\nbigmzero)
\to
\lambda^{-1}\Vzero_0\nbigk(\nbigmzero).
\]

Let $\Xhat_P\simeq \cnum\times \Xtilde_P$
be another holomorphic trivialization
from which we obtain a holomorphic function $w$
and the holomorphic vector field $\del_w$.
We note that the $V$-filtrations for $z$ and $w$ are the same.
There exists a nowhere vanishing holomorphic function
$f$ on $\Xtilde_P$ such that $z=fw$.
We have the relation
$\del_{w}=\del_{w}(fw)\del_{z}
=f\del_z+w\del_w(f)\cdot\del_z$.
We obtain
\begin{equation}
 -\del_{w}\otimes (w\stilde)
=-f\del_z\otimes (zf^{-1} \stilde)
-w\del_w(f)\del_z\otimes w\stilde
=
 -\del_z\otimes z\stilde
 -1\otimes f\del_z(f^{-1})z\stilde
 -1\otimes (w\del_w(f)\del_z(w\stilde)).
\end{equation}
We also have
\begin{equation}
 1\otimes\del_{w}(w\stilde)
=1\otimes f\del_z(zf^{-1} \stilde)
+1\otimes w\del_w(f)\del_z( w\stilde)
=
 1\otimes \del_z(z\stilde)
 +1\otimes f\del_z(f^{-1})z\stilde
 +1\otimes (w\del_w(f)\del_z(w\stilde)).
\end{equation}
We obtain
\[
  -\del_{w}\otimes (w\stilde)
  +1\otimes\del_{w}(w\stilde)
  =-\del_z\otimes z\stilde
  +1\otimes \del_z(z\stilde).
\]
Hence, the morphism $\upsilonhat_P$
is independent of the choice of
a holomorphic trivialization $\Xhat_P\simeq\cnum\times\Xtilde_P$.

Therefore,
we obtain a global isomorphism
$\upsilonhat:\nbigkhat\simeq\nbigchat\otimes\newTate(1)$
such that
$-\nbigshat_2\circ\upsilonhat$ is a polarization of
$\nbigkhat$.
Because
$\nbigkhat=\iota_{1\dagger}\nbigk$
and $\nbigchat=\iota_{1\dagger}\nbigc$,
there exists an isomorphism
$\upsilon:\nbigk\simeq\nbigc\otimes\newTate(1)$
such that $\upsilonhat=\iota_{1\dagger}\upsilon$
for which
$-\nbigs_2\circ\upsilon$ is a polarization
of $\nbigk$.
\hfill\qed

\vspace{.1in}
We note $\dim\supp(\nbigk)<\dim\supp(\nbigt)$.
The tuple
$\Bigl(
\bigoplus a_{\Xtilde\dagger}^j(\nbigk),
L_{\omega_1+c\omega_2},
a_{\Xtilde\dagger}(\nbigs_{\nbigk})
\Bigr)$
is a polarized graded Lefschetz twistor structure of
weight $w-1$ and type $1$
for any $0<c\leq c_0$.
We set
\[
 \abar_{\Xtilde\dagger}^j(\nbigk):=
 a_{\Xtilde\dagger}^j(\nbigk)
 \otimes\newTate(j),
 \quad
 \fbar_{\dagger}^j(\nbigk):=
 f_{\dagger}^j(\nbigk)
 \otimes\newTate(j),
 \quad
 \abar_{\dagger}^k\fbar_{\dagger}^j(\nbigk):=
 a_{\dagger}^k\fbar_{\dagger}^j(\nbigk)
 \otimes\newTate(k).
\]
Then, 
for any $0<c\leq c_0$,
$\Bigl(
\bigoplus \abar_{\Xtilde\dagger}^j(\nbigk),
L_{\omega_1+c\omega_2},
\abar_{\Xtilde\dagger}(\nbigs_{\nbigk})
\Bigr)$
is a polarized graded Lefschetz twistor structure of
weight $w-1$ and type $-1$,
and hence 
\[
\Bigl(
\bigoplus \abar_{\Xtilde\dagger}^j(\nbigk),
L_{\omega_1},
L_{\omega_2+c_0^{-1}\omega_1},
\abar_{\Xtilde\dagger}(\nbigs_{\nbigk})
\Bigr)
\]
is a polarized mixed twistor structure of weight $w-1$.
There exists the Leray filtration $L(\abar_{\Xtilde\dagger}^i(\nbigk))$,
and we have
$\Gr_L^k\abar^i_{\Xtilde\dagger}(\nbigk)
\simeq
 \abar_{X\dagger}^{k}
 \fbar^{i-k}_{\dagger}(\nbigk)$.

\subsubsection{Boundary morphisms}

There exist the following exact sequences of
mixed twistor $\nbigd$-modules:
\begin{equation}
 0\lrarr
  \nbigt\lrarr
  \nbigt_1\lrarr\nbigk\otimes\newTate(-1)
  \lrarr 0,
\end{equation}
\begin{equation}
 0\lrarr\nbigk\lrarr
  \nbigt_2\lrarr
  \nbigt\lrarr 0.
\end{equation}
We obtain the morphisms
$\del_1:a_{\Xtilde\dagger}^{j}\nbigk
\lrarr
a_{\Xtilde\dagger}^{j+1}(\nbigt)
\otimes\newTate(1)$
and
$\del_2:a_{\Xtilde\dagger}^j(\nbigt)
\lrarr
a_{\Xtilde\dagger}^{j+1}\nbigk$.
As the composition,
we obtain
\begin{equation}
 \label{eq;21.12.21.1}
 \del_1\circ\del_2:
 a_{\Xtilde\dagger}^j(\nbigt)
 \lrarr
 a_{\Xtilde\dagger}^{j+2}(\nbigt)\otimes
 \newTate(1),
 \quad
 \del_2\circ\del_1:
 a_{\Xtilde\dagger}^j(\nbigk)
 \lrarr
 a_{\Xtilde\dagger}^{j+2}(\nbigk)\otimes
 \newTate(1).
\end{equation}

\begin{lem}
The morphisms in {\rm(\ref{eq;21.12.21.1})}
are $-L_{\omega_2}$.
\end{lem}
\pf
By Proposition \ref{prop;22.3.2.30} below,
the specializations of the morphisms in
(\ref{eq;21.12.21.1}) along general $\lambda$
are equal to the specialization of $-L_{\omega_2}$
along $\lambda$.
Then, we obtain the claim of the lemma
by the pure twistor property.
\hfill\qed

\vspace{.1in}

Similarly, we have
$\del_1:f_{\dagger}^{j}\nbigk
\lrarr
f_{\dagger}^{j+1}(\nbigt)\otimes\newTate(1)$
and
$\del_2:f_{\dagger}^j(\nbigt)
 \lrarr
 f_{\dagger}^{j+1}\nbigk$.
We obtain
\begin{equation}
 \label{eq;21.12.21.2}
  \del_1\circ\del_2:
 f_{\dagger}^j(\nbigt)
 \lrarr
 f_{\dagger}^{j+2}(\nbigt)\otimes
 \newTate(1),
 \quad\quad
 \del_2\circ\del_1:
 f_{\dagger}^j(\nbigk)
 \lrarr
 f_{\dagger}^{j+2}(\nbigk)\otimes
 \newTate(1).
\end{equation}
We obtain the following lemma
from Proposition \ref{prop;22.3.2.30} below.
\begin{lem}
\label{lem;22.3.4.10}
The morphisms in {\rm(\ref{eq;21.12.21.2})}
are $-L_{\omega_2}$.
\hfill\qed
\end{lem}

Note that the following diagram is commutative:
\begin{equation}
\begin{CD}
a_{\Xtilde\dagger}^j(\nbigk)
@>{-a_{\Xtilde\dagger}\nbigs_{\nbigk}}>>
a_{\Xtilde\dagger}^{-j}(\nbigk)^{\ast}\otimes\newTate(-w+1)\\
@V{\del_1}VV @V{\del_2^{\ast}}VV \\
a_{\Xtilde\dagger}^{j+1}(\nbigt)\otimes\newTate(1)
@>{a_{\Xtilde\dagger}\nbigs}>>
a_{\Xtilde\dagger}^{-j-1}(\nbigt)^{\ast}\otimes\newTate(-w+1).
\end{CD}
\end{equation}

We have the induced morphisms:
\[
\del_1:
 \abar^{j}_{\Xtilde\dagger}(\nbigk)
\lrarr
 \abar^{j+1}_{\Xtilde\dagger}(\nbigt),
 \quad\quad
 \del_2:
\abar^j_{\Xtilde\dagger}(\nbigt)
 \lrarr
 \abar^{j+1}_{\Xtilde\dagger}(\nbigk)
 \otimes\newTate(-1).
\]
\[
\del_1:
 \fbar^{j}_{\dagger}(\nbigk)
\lrarr
 \fbar^{j+1}_{\dagger}(\nbigt),
\quad\quad
\del_2:
\fbar^j_{\dagger}(\nbigt)
 \lrarr
 \fbar^{j+1}_{\dagger}(\nbigk)
 \otimes\newTate(-1).
\]
The following diagram is commutative:
\begin{equation}
\label{eq;22.3.4.1}
\begin{CD}
\abar_{\Xtilde\dagger}^j(\nbigk)
@>{\abar_{\Xtilde\dagger}\nbigs_{\nbigk}}>>
a_{\Xtilde\dagger}^{-j}(\nbigk)^{\ast}\otimes\newTate(-w+1)\\
@V{\del_1}VV @V{\iota_{\newTate(1)}\circ\del_2^{\ast}}VV \\
\abar_{\Xtilde\dagger}^{j+1}(\nbigt)
@>{\abar_{\Xtilde\dagger}\nbigs}>>
\abar_{\Xtilde\dagger}^{-j-1}(\nbigt)^{\ast}\otimes\newTate(-w).
\end{CD}
\end{equation}

\subsubsection{The graded object associated with
the monodromy weight filtration of $L_{\omega_1}$}

Let $W(\omega_1)$ be the monodromy weight filtration
of $\bigoplus \abar^j_{\Xtilde\dagger}(\nbigt)$
with respect to $L_{\omega_1}$.
Because $L_{\omega_1}$ preserves the grading,
$W(\omega_1)$ is compatible
with the grading
$\bigoplus \abar^j_{\Xtilde\dagger}(\nbigt)$,
i.e.,
\[
 W_{k}(\omega_1)\Bigl(
\bigoplus_j \abar^j_{\Xtilde\dagger}(\nbigt)
 \Bigr)
 =\bigoplus_{j}
 \Bigl(
 W_k(\omega_1)\cap \abar^j_{\Xtilde\dagger}(\nbigt)
 \Bigr).
\] 
We set
\[
 \nbigt^{(0)}_{\ell_1,\ell_2}=
 \Gr^{W(\omega_1)}_{\ell_1}
 \bigl(
 \abar_{\Xtilde\dagger}^{-\ell_1-\ell_2}(\nbigt)
 \bigr).
\]
Each $\nbigt^{(0)}_{\ell_1,\ell_2}$
is pure of weight $w+\ell_1+\ell_2$.
We obtain the induced morphisms
\[
L^{(0)}_{\omega_1}:
\nbigt^{(0)}_{\ell_1,\ell_2}
 \lrarr
\nbigt^{(0)}_{\ell_1-2,\ell_2}
 \otimes\newTate(-1),
\quad\quad
 L^{(0)}_{\omega_2}:
 \nbigt^{(0)}_{\ell_1,\ell_2}
 \lrarr
\nbigt^{(0)}_{\ell_1,\ell_2-2}
 \otimes\newTate(-1),
\]
induced by $L_{\omega_1}$ and $L_{\omega_2}$,
respectively.
They define bi-graded morphisms
$L^{(0)}_{\omega_i}:
\nbigt^{(0)}\lrarr
\nbigt^{(0)}\otimes\newTate(-1)$.
The bi-degree of $L^{(0)}_{\omega_1}$
and $L^{(0)}_{\omega_2}$
are $(-2,0)$ and $(0,-2)$, respectively.
We have the isomorphisms
\[
 \nbigs^{(0)}_{\ell_1,\ell_2}:
 \nbigt^{(0)}_{\ell_1,\ell_2}
 \simeq
 (\nbigt^{(0)}_{-\ell_1,-\ell_2})^{\ast}\otimes\newTate(-w)
\]
induced by $a_{\Xtilde\dagger}\nbigs$.
They induce a bi-graded Hermitian sesqui-linear duality
\[
 \nbigs^{(0)}:
 \nbigt^{(0)}\simeq
 (\nbigt^{(0)})^{\ast}\otimes\newTate(-w).
\]

\begin{lem}
 $\Bigl(
 \bigoplus_{\ell_1,\ell_2}
 \nbigt^{(0)}_{\ell_1,\ell_2},
 L^{(0)}_{\omega_1},
 L^{(0)}_{\omega_2},
 \nbigs^{(0)}
 \Bigr)$
 is a polarized bi-graded Lefschetz twistor structure
 of weight $w$ and type $(-1,-1)$.
\end{lem}
\pf
It follows from \cite[Proposition 3.115]{mochi2}.
\hfill\qed

\vspace{.1in}

Similarly, there exists the weight filtration $W(\omega_1)$
of $\bigoplus_j \abar^j_{\Xtilde\dagger}(\nbigk)$
with respect to $L_{\omega_1}$.
Because $L_{\omega_1}$
preserves the grading $\bigoplus_j \abar^j_{\Xtilde\dagger}(\nbigk)$,
$W(\omega_1)$ is compatible with the grading.
In this case, by the hypothesis of the induction
(Assumption \ref{assumption;22.3.10.1}),
it is essentially the same as the Leray filtration.
\begin{lem}
\label{lem;22.3.15.30}
The weight filtration $W(\omega_1)$ of $L_{\omega_1}$
and the Leray filtration $L$
on 
$a_{\Xtilde\dagger}^j(\nbigk)$
are related as follows:
 \[
 W_k(\omega_1)\bigl(
 \abar^j_{\Xtilde\dagger}(\nbigk)
 \bigr)
 =L^{-k}\bigl(
 \abar^j_{\Xtilde\dagger}(\nbigk)
 \bigr).
\]
\end{lem}
\pf
Note that the dimension of the support of
$\bigoplus f_{\dagger}^j(\nbigk)$
is strictly smaller than
the dimension of the support of $\nbigt$.
Hence, the claim follows from
Assumption \ref{assumption;22.3.10.1}.
\hfill\qed

\vspace{.1in}

 We set
\[
 \nbigt^{(1)}_{\ell_1,\ell_2}:=
 \Gr^{W(\omega_1)}_{\ell_1}
 \bigl(
 \abar_{\Xtilde\dagger}^{-\ell_1-\ell_2}(\nbigk)
 \bigr)
 \simeq
 \abar_{X\dagger}^{-\ell_1}
 \fbar_{X\dagger}^{-\ell_2}(\nbigk).
\]
Each $\nbigt^{(1)}_{\ell_1,\ell_2}$
is pure of weight $w-1+\ell_1+\ell_2$.
We obtain the morphisms
\[
 L^{(1)}_{\omega_1}:
 \nbigt^{(1)}_{\ell_1,\ell_2}
 \lrarr
 \nbigt^{(1)}_{\ell_1-2,\ell_2}\otimes\newTate(-1),
 \quad\quad
 L^{(1)}_{\omega_2}:
 \nbigt^{(1)}_{\ell_1,\ell_2}
 \lrarr
 \nbigt^{(1)}_{\ell_1,\ell_2-2}\otimes\newTate(-1)
\]
induced by $L_{\omega_1}$ and $L_{\omega_2}$,
respectively.
They define
$L^{(1)}_{\omega_i}:\nbigt^{(1)}\lrarr
\nbigt^{(1)}\otimes\newTate(-1)$.
We have the induced isomorphisms
\[
 \nbigs^{(1)}_{\ell_1,\ell_2}:
 \nbigt^{(1)}_{\ell_1,\ell_2}
 \simeq
 (\nbigt^{(1)}_{\ell_1,\ell_2})^{\ast}\otimes\newTate(-w+1),
\]
induced by $\abar_{\Xtilde\dagger}\nbigs_{\nbigk}$.
They define a bi-graded Hermitian sesqui-linear duality
$\nbigs^{(1)}:\nbigt^{(1)}\simeq
 (\nbigt^{(1)})^{\ast}\otimes\newTate(-w+1)$.
\begin{lem}
 $\Bigl(
 \bigoplus\nbigt^{(1)}_{\ell_1,\ell_2},
 L^{(1)}_{\omega_1},
 L^{(1)}_{\omega_2},
 \nbigs^{(1)}
 \Bigr)$
 is a polarized bi-graded Lefschetz twistor structure of weight $w-1$
 and type $(-1,-1)$.
\hfill\qed
\end{lem}

There exist the induced morphisms
$g_1:\nbigt^{(1)}\lrarr\nbigt^{(0)}$
induced by $\del_1$,
and
$g_2:\nbigt^{(0)}\lrarr\nbigt^{(1)}\otimes\newTate(-1)$
induced by $\del_2$.
More precisely, we have
\[
g_1:
 \nbigt^{(1)}_{\ell_1,\ell_2}
 \lrarr
  \nbigt^{(0)}_{\ell_1,\ell_2-1},
\quad\quad
g_2:
 \nbigt^{(0)}_{\ell_1,\ell_2}
 \lrarr
  \nbigt^{(1)}_{\ell_1,\ell_2-1}\otimes\newTate(-1).
\]
By the construction,
we have
$g_2\circ g_1=-L^{(0)}_{\omega_2}$
and
$g_1\circ g_2=-L^{(1)}_{\omega_2}$.
We also have
\[
 g_1\circ L^{(0)}_{\omega_i}
 =L^{(1)}_{\omega_i}\circ g_1,
 \quad\quad
 g_2\circ L^{(1)}_{\omega_i}
 =L^{(0)}_{\omega_i}\circ g_2.
\]
\begin{lem}
\label{lem;22.3.4.2}
We have the decomposition
$\nbigt^{(0)}=\Image g_1\oplus \Ker g_2$.
\end{lem}
\pf
For $k\geq 0$,
we set
\[
P_{\omega_1}\nbigt^{(i)}_{k,\ell}=
\Ker\Bigl(
 (L^{(i)}_{\omega_1})^{k+1}:
 \nbigt^{(i)}_{k,\ell}
 \lrarr
 \nbigt^{(i)}_{-k-2,\ell}\otimes\newTate(-k-1)
\Bigr).
\]
We set
\[
 P_{\omega_1}\nbigt^{(i)}_{k,\bullet}:=
 \bigoplus_{\ell}
 P_{\omega_1}\nbigt^{(i)}_{k,\ell}.
\]
We have
$L^{(i)}_{\omega_2}:
P_{\omega_1}\nbigt^{(i)}_{k,\bullet}
\lrarr
P_{\omega_1}\nbigt^{(i)}_{k,\bullet}
\otimes\newTate(-1)$.
We also obtain the Hermitian sesqui-linear dualities:
\[
P_{\omega_1}\nbigs^{(0)}_k=
(-1)^k\nbigs^{(0)}\circ (L^{(0)}_{\omega_1})^k:
P_{\omega_1}\nbigt^{(0)}_{k,\bullet}
\lrarr
\Bigl(
P_{\omega_1}\nbigt^{(0)}_{k,\bullet}
\Bigr)^{\ast}\otimes\newTate(-w-k).
\]
\[
P_{\omega_1}\nbigs^{(1)}_k=
(-1)^k\nbigs^{(1)}\circ (L^{(1)}_{\omega_1})^k:
P_{\omega_1}\nbigt^{(1)}_{k,\bullet}
\lrarr
\Bigl(
P_{\omega_1}\nbigt^{(1)}_{k,\bullet}
\Bigr)^{\ast}\otimes\newTate(-w+1-k).
\]
Then,
$(P_{\omega_1}\nbigt^{(0)}_{k,\bullet},L^{(0)}_{\omega_2},
P_{\omega_1}\nbigs^{(0)}_k)$
is a polarized graded Lefschetz twistor structure
of weight $w+k$ and type $-1$,
and
$(P_{\omega_1}\nbigt^{(1)}_{k,\bullet},L^{(1)}_{\omega_2},
P_{\omega_1}\nbigs^{(1)}_k)$
is a polarized
graded Lefschetz twistor structure of weight $w+k-1$ and type $-1$.

We obtain the induced morphisms
\[
 P_{\omega_1}(g_1)_k:
 P_{\omega_1}\nbigt^{(1)}_{k,\bullet}
 \lrarr
 P_{\omega_1}\nbigt^{(0)}_{k,\bullet},
\quad\quad
P_{\omega_1}(g_2)_k:
 P_{\omega_1}\nbigt^{(0)}_{k,\bullet}
 \lrarr
 P_{\omega_1}\nbigt^{(1)}_{k,\bullet}
 \otimes\newTate(-1).
\]
More precisely, we have
\[
 P_{\omega_1}(g_1)_k:
 P_{\omega_1}\nbigt^{(1)}_{k,\ell}
 \lrarr
 P_{\omega_1}\nbigt^{(0)}_{k,\ell-1},
\quad\quad
P_{\omega_1}(g_2)_k:
 P_{\omega_1}\nbigt^{(0)}_{k,\ell}
 \lrarr
 P_{\omega_1}\nbigt^{(1)}_{k,\ell-1}
 \otimes\newTate(-1).
\]
We have
$P_{\omega_1}(g_2)_k\circ P_{\omega_1}(g_1)_k
=-L^{(1)}_{\omega_1}$
and
$P_{\omega_1}(g_1)_k\circ
 P_{\omega_1}(g_2)_k
=-L^{(0)}_{\omega_1}$.
The following diagram is commutative
because (\ref{eq;22.3.4.1}) is commutative:
\[
 \begin{CD}
  P_{\omega_1}\nbigt^{(1)}_{k,\bullet}
  @>{P_{\omega_1}\nbigs^{(1)}_k}>>
  (P_{\omega_1}\nbigt^{(1)}_{k,\bullet})^{\ast}
  \otimes\newTate(-w+1-k)
  \\
  @V{P_{\omega_1}(g_1)_k}VV @V{P_{\omega_1}(g_2)_k^{\ast}}VV \\
  P_{\omega_1}\nbigt^{(0)}_{k,\bullet}
  @>{P_{\omega_1}\nbigs^{(0)}_k}>>
  (P_{\omega_1}\nbigt^{(0)}_{k,\bullet})^{\ast}
  \otimes\newTate(-w-k)
 \end{CD}
\]
By \cite[Proposition 2.1.19]{sabbah2},
which goes back to \cite[Lemma 5.2.15]{saito1},
we obtain 
$P_{\omega_1}\nbigt^{(0)}_{k,\bullet}
=\Image P_{\omega_1}(g_1)_k
\oplus
\Ker P_{\omega_1}(g_2)_k$.
Then, the claim of Lemma \ref{lem;22.3.4.2} follows.
\hfill\qed

\subsubsection{The Leray filtrations on
$\nbigt^{(0)}_{\ell_1,\ell_2}$
and $\nbigt^{(1)}_{\ell_1,\ell_2}$}

We can reformulate Lemma \ref{lem;22.3.15.30} as follows.
\begin{lem}
The induced Leray filtration $L$ of 
$\nbigt^{(1)}_{k,\ell}$ is pure:
\[
 \Gr_L^p\nbigt^{(1)}_{k,\ell}
 =\left\{
 \begin{array}{ll}
  \abar_{X\dagger}^{-k}
   \fbar_{\dagger}^{-\ell}(\nbigk)
   & (p=-k) \\
  0 & (p\neq -k)
 \end{array}
 \right.
\] 
\hfill\qed
\end{lem}

Let us study the induced Leray filtration $L$ of 
$\nbigt^{(0)}_{k,\ell}$.
\begin{lem}
\label{lem;22.3.4.40}
If $\ell\neq 0$, the following holds:
\begin{equation}
\label{eq;22.3.15.40}
 \Gr_L^p\nbigt^{(0)}_{k,\ell}
 =\left\{
 \begin{array}{ll}
  \abar_{X\dagger}^{-k}
  \fbar_{\dagger}^{-\ell}(\nbigt)& (p=-k) \\
  \Gr^{W(\omega_1)}_{k}(\abar_{X\dagger}^{-k-\ell}\nbigt_X)
   & (p=-k-\ell)\\
  0 & (p\neq -k,-k-\ell)
 \end{array}
  \right.
\end{equation}
In the case $\ell=0$,
we have
$\Gr_L^{-k}(\nbigt^{(0)}_{k,0})
=\nbigt^{(0)}_{k,0}$,
and there exists the decomposition
\begin{equation}
\label{eq;22.3.15.60}
 \nbigt_{k,0}^{(0)}
 =\abar_{X\dagger}^{-k}(\nbigt_X')
 \oplus
 \Gr^{W(\omega_1)}_k(\abar_{X\dagger}^{-k}(\nbigt_X)).
\end{equation}
\end{lem}
\pf
Because there exists an isomorphism
$f_{\dagger}(\nbigt)\simeq
\bigoplus f_{\dagger}^j(\nbigt)[-j]$
in the derived category of $\nbigr_X$-triples,
there exists a decomposition
\begin{equation}
\label{eq;22.3.4.100}
 \bigoplus_j
 \abar_{\Xtilde\dagger}^j(\nbigt)
 \simeq
 \bigoplus_j
 \Bigl(
 \bigoplus_{k_1+k_2=j}
 \abar_{X\dagger}^{k_1}\fbar_{\dagger}^{k_2}(\nbigt)
 \Bigr)
\end{equation}
in a way compatible with the actions of $L_{\omega_1}$,
and hence the filtrations $W(\omega_1)$.
It gives a splitting of the Leray filtration,
i.e.,
\[
 L^p\abar_{\Xtilde\dagger}^j(\nbigt)
=\bigoplus_{\substack{k_1+k_2=j \\ k_1\geq p}}
 \abar_{X\dagger}^{k_1}\fbar_{\dagger}^{k_2}(\nbigt).
\]
The dimension of the supports of
$f^j_{\dagger}(\nbigt)$ $(j\neq 0)$
and $\nbigt_X'$
are strictly smaller than the dimension of the support of $\nbigt$.
Hence, we have
\[
 W_{k}(\omega_1)
 \Bigl(
 \bigoplus_{\ell}
 \abar_{X\dagger}^{\ell}
 \fbar_{\dagger}^j(\nbigt)
 \Bigr)
 \simeq
 \bigoplus_{\ell\geq -k}
 \abar_{X\dagger}^{\ell}
 \fbar_{\dagger}^j(\nbigt)
 \quad(j\neq 0),
\]
\[
  W_{k}(\omega_1)
 \Bigl(
 \bigoplus_{\ell}
 \abar_{X\dagger}^{\ell}
 \nbigt_X'
 \Bigr)
 \simeq
 \bigoplus_{\ell\geq -k}
 \abar_{X\dagger}^{\ell}
 \nbigt_X'.
\]
Hence, there exist the following isomorphisms:
\begin{equation}
\label{eq;22.3.15.50}
 \Gr_{\ell}^{W(\omega_1)}
 \Bigl(
 \abar_{\Xtilde\dagger}^j(\nbigt)
 \Bigr)
 \simeq
 \abar_{X\dagger}^{-\ell}
 \fbar^{j+\ell}_{\dagger}(\nbigt)
 \oplus
 \Gr^{W(\omega_1)}_{\ell}
 \Bigl(
 \abar_{X\dagger}^{j}
 \nbigt_X
 \Bigr)
 \quad
 (j\neq -\ell),
\end{equation}
\begin{equation}
\label{eq;22.3.15.51}
 \Gr_{-j}^{W(\omega_1)}
 \Bigl(
 \abar_{\Xtilde\dagger}^j(\nbigt)
 \Bigr)
 \simeq
 \abar_{X\dagger}^{j}
 \nbigt_X'
 \oplus
 \Gr^{W(\omega_1)}_{-j}
 \Bigl(
 \abar_{X\dagger}^{j}
 \nbigt_X
 \Bigr).
\end{equation}
By (\ref{eq;22.3.15.50}),
if $j+\ell\neq 0$,
we have
\[
 \Gr_L^{-\ell}
 \Gr_{\ell}^{W(\omega_1)}
 \Bigl(
 \abar_{\Xtilde\dagger}^j(\nbigt)
 \Bigr)
 \simeq
  \abar_{X\dagger}^{-\ell}
  \fbar^{j+\ell}_{\dagger}(\nbigt),
  \quad\quad
  \Gr_L^{j}
   \Gr_{\ell}^{W(\omega_1)}
 \Bigl(
 \abar_{\Xtilde\dagger}^j(\nbigt)
 \Bigr)
 \simeq
  \Gr^{W(\omega_1)}_{\ell}
 \Bigl(
 \abar_{X\dagger}^{j}
 \nbigt_X
 \Bigr).
\]
We also have
$\Gr_L^{k}
   \Gr_{\ell}^{W(\omega_1)}
 \Bigl(
 \abar_{\Xtilde\dagger}^j(\nbigt)
 \Bigr)=0$
 for $k\neq j,-\ell$.
It implies (\ref{eq;22.3.15.40}).
We obtain (\ref{eq;22.3.15.60})
from (\ref{eq;22.3.15.51}).
Then, the claim of Lemma \ref{lem;22.3.4.40} follows. 
\hfill\qed

\subsubsection{Factorization}

Let us look at the morphisms
$g_1:
 \nbigt^{(1)}_{k,\ell+1}
 \lrarr
 \nbigt^{(0)}_{k,\ell}$.

\begin{lem}
\label{lem;22.3.4.110}
If $\ell>0$, $g_1$ is factorized as follows:
\begin{equation}
\label{eq;22.3.4.50}
\nbigt^{(1)}_{k,\ell+1}
\lrarr
L^{-k}\nbigt^{(0)}_{k,\ell}
=\abar_{X\dagger}^{-k}\fbar_{\dagger}^{-\ell}\nbigt
\lrarr
\nbigt^{(0)}_{k,\ell}
 \simeq\abar_{X\dagger}^{-k}\fbar_{\dagger}^{-\ell}\nbigt
 \oplus \Gr^{W(\omega_1)}_{k}
 \abar_{X\dagger}^{-k-\ell}(\nbigt_X).
\end{equation}
If $\ell=0$, $g_1$ is factorized as follows:
\begin{equation}
\label{eq;22.3.4.51}
 \nbigt^{(1)}_{k,1}
\lrarr
\abar_{X\dagger}^{-k}\nbigt_X'
\lrarr
\nbigt^{(0)}_{k,0}
=\abar_{X\dagger}^{-k}\nbigt_X'
 \oplus
 \Gr^{W(\omega_1)}_{k}
 \abar_{X\dagger}^{-k}(\nbigt_X).
 \end{equation}
\end{lem}
\pf
Because the Leray filtrations are preserved,
we obtain the factorization (\ref{eq;22.3.4.50}).
Because the Leray filtrations on
$\nbigt^{(1)}_{k,1}$ and
$\nbigt^{(0)}_{k,0}$ are pure of weight $-k$,
the morphism
$\nbigt^{(1)}_{k,1}
\lrarr\nbigt^{(0)}_{k,0}$ is
identified with the following morphism:
\[
 \Gr_L^{-k}
 \Gr^{W(\omega_1)}_k
 \bigl(
 \abar_{\Xtilde\dagger}^{-k-1}(\nbigk)
 \bigr)
 \lrarr
 \Gr_L^{-k}
 \Gr^{W(\omega_1)}_k
 \bigl(
 \abar_{\Xtilde\dagger}^{-k}(\nbigt)
 \bigr).
\]
It is equal to the morphism
$\Gr^{W(\omega_1)}_k
 \Gr_L^{-k}
 \bigl(
 \abar_{\Xtilde\dagger}^{-k-1}(\nbigk)
 \bigr)
 \lrarr
 \Gr^{W(\omega_1)}_k
 \Gr_L^{-k}
 \bigl(
 \abar_{\Xtilde\dagger}^{-k}(\nbigt)
 \bigr)$,
which is induced by
$f_{\dagger}^{-1}(\nbigk)\lrarr f_{\dagger}^0(\nbigt)
=\nbigt_X\oplus\nbigt_X'$.
Because the dimension of the support of
$f_{\dagger}^{-1}(\nbigk)$
is strictly smaller than
the dimension of the support of $\nbigt_X$,
it is factorized as
$f_{\dagger}^{-1}(\nbigk)\lrarr
\nbigt_X'\lrarr \nbigt_X\oplus\nbigt_X'$.
Then, we obtain the factorization
(\ref{eq;22.3.4.51}).
\hfill\qed

\vspace{.1in}

Let us look at the morphisms
$g_2:
 \nbigt^{(0)}_{k,\ell}
 \lrarr
 \nbigt^{(1)}_{k,\ell-1}
 \otimes\newTate(-1)$.

\begin{lem}
\label{lem;22.3.4.111}
If $\ell<0$, $g_2$ is factorized as follows:
 \begin{equation}
\label{eq;22.3.4.52}
  \nbigt^{(0)}_{k,\ell}
   \simeq
   \abar_{X\dagger}^{-k}\fbar_{\dagger}^{-\ell}\nbigt
 \oplus \Gr^{W(\omega_1)}_{k}
 \abar_{X\dagger}^{-k-\ell}(\nbigt_X)
 \lrarr
 \nbigt^{(0)}_{k,\ell}/L^{-k-\ell}\nbigt^{(0)}_{k,\ell}
 =\abar_{X\dagger}^{-k}\fbar_{\dagger}^{-\ell}\nbigt
 \lrarr
 \nbigt^{(1)}_{k,\ell-1}
  \otimes\newTate(-1).
 \end{equation}
If $\ell=0$,
$g_2$ is factorized as follows:
\begin{equation}
\label{eq;22.3.4.53}
   \nbigt^{(0)}_{k,0}
   \simeq
   \abar_{X\dagger}^{-k}\nbigt_X'
 \oplus \Gr^{W(\omega_1)}_{k}
 \abar_{X\dagger}^{-k}(\nbigt_X)
 \lrarr
\abar_{X\dagger}^{-k}\nbigt'_X
 \lrarr
 \nbigt^{(1)}_{k,-1}
\otimes\newTate(-1).
\end{equation}
\end{lem}
\pf
Because the Leray filtrations are preserved,
we obtain the factorization (\ref{eq;22.3.4.52}).
By a similar argument to the proof of Lemma \ref{lem;22.3.4.110},
we can observe that
$\nbigt^{(0)}_{k,0}\lrarr
\nbigt^{(1)}_{k,-1}
\otimes\newTate(-1)$
is induced by
$\fbar_{\dagger}^0(\nbigt)
\lrarr
\fbar_{\dagger}^1(\nbigk)\otimes\newTate(-1)$.   
By considering the supports of
$\nbigt_X$,
we obtain the factorization (\ref{eq;22.3.4.53}).
\hfill\qed

\subsubsection{Proof of Proposition \ref{prop;22.3.4.30}}

\begin{lem}
\label{lem;22.3.4.120}
We have $\Gr^{W(\omega_1)}_k
\bigl(
  \abar_{X\dagger}^{-k-\ell}\nbigt_X
  \bigr)=0$
unless $\ell\leq 0$.
\end{lem}
\pf
We consider the following morphisms:
\[
 \begin{CD}
  \nbigt^{(1)}_{k,\bullet}
  @>{g_1}>>
  \nbigt^{(0)}_{k,\bullet}
  @>{g_2}>>
  \nbigt^{(1)}_{k,\bullet}\otimes\newTate(-1).
 \end{CD}
\]
By taking a decomposition (\ref{eq;22.3.4.100})
as in the proof of Lemma \ref{lem;22.3.4.40},
we obtain a decomposition
\[
 \nbigt^{(0)}_{k,\bullet}
 =\bigoplus_{\ell\neq 0}
 \abar_{X\dagger}^{-k}\fbar_{\dagger}^{-k-\ell}(\nbigt)
 \oplus
  \abar_{X\dagger}^{-k}(\nbigt'_X)
  \oplus
\bigoplus_{\ell}
  \Gr^{W(\omega_1)}_{k}\Bigl(
  \abar_{X\dagger}^{-k-\ell}\nbigt_X
  \Bigr).
\]
By Lemma \ref{lem;22.3.4.110},
we obtain
\[
 \Image g_1\cap
 \bigoplus_{\ell}
  \Gr^{W(\omega_1)}_{k}\Bigl(
  \abar_{X\dagger}^{-k-\ell}\nbigt_X
  \Bigr)
  \subset
   \bigoplus_{\ell<0}
  \Gr^{W(\omega_1)}_{k}\Bigl(
  \abar_{X\dagger}^{-k-\ell}\nbigt_X
  \Bigr).
\]
By Lemma \ref{lem;22.3.4.111},
we obtain
\[
 \Image g_1\cap
 \bigoplus_{\ell}
  \Gr^{W(\omega_1)}_{k}\Bigl(
  \abar_{X\dagger}^{-k-\ell}\nbigt_X
  \Bigr)
\subset\Ker g_2.
\]
By Lemma \ref{lem;22.3.4.2},
we obtain
\[
  \Image g_1\cap
 \bigoplus_{\ell}
  \Gr^{W(\omega_1)}_{k}\Bigl(
  \abar_{X\dagger}^{-k-\ell}\nbigt_X
  \Bigr)
  \subset
  \Image g_1\cap\Ker g_2=0.
\]
The projection
$\nbigt^{(0)}_{k,\bullet}=\Image g_1\oplus\Ker g_2
\lrarr\Ker g_2$
induces a monomorphism
\[
  \bigoplus_{\ell}
  \Gr^{W(\omega_1)}_{k}\Bigl(
  \abar_{X\dagger}^{-k-\ell}\nbigt_X
  \Bigr)
  \lrarr
  \Ker g_2.
\]
Because $\Ker g_2$
is a pure twistor structure of weight $w+k$,
the proof of Lemma \ref{lem;22.3.4.120}
is reduced to the next lemma.
\begin{lem}
\label{lem;22.3.4.121}
$\Gr^{W(\omega_1)}_k(\abar_{X\dagger}^{j}\nbigt_X)$
is pure of weight $w-j$.
\end{lem}
\pf
Because 
$\abar_{X\dagger}^{j}\nbigt_X$
is a direct summand of
$\abar_{\Xtilde\dagger}^{j}\nbigt$,
it is a pure of weight $w-j$.
Then, by the construction of $W(\omega_1)$,
we obtain that
$\Gr^{W(\omega_1)}_k(\abar_{X\dagger}^{j}\nbigt_X)$
is pure of weight $w-j$.
Thus, we obtain Lemma \ref{lem;22.3.4.121}
and hence Lemma \ref{lem;22.3.4.120}.
\hfill\qed

\vspace{.1in}

We consider the following filtration:
\[
 W'_k:=
 \bigoplus_{j\geq -k}
\abar_{X\dagger}^{j}\nbigt_X.
\]
\begin{lem}
\label{lem;22.3.4.130}
We obtain $W_k(\omega_1)=W_k'$
for any $k\in\seisuu$.
As a result,
$L_{\omega_1}^j:
\abar_{X\dagger}^{-j}\nbigt_X
\lrarr
\abar_{X\dagger}^{j}\nbigt_X
\otimes\newTate(-j)$
are isomorphisms for any $j\geq 0$.
\end{lem}
\pf
By Lemma \ref{lem;22.3.4.120},
we have 
$W_k(\omega_1)\subset W_k'$ for any $k\in\seisuu$.
It is enough to prove that
$\rank\Gr^{W(\omega_1)}_k=\rank\Gr^{W'}_k$
for any $k\in\seisuu$.
Suppose that
$\{k\in\seisuu\,|\,
\rank\Gr^{W(\omega_1)}_k\neq \rank\Gr^{W'}_k\}$
is not empty,
and let $k_0$ denote the minimum.
Note that $a_{X\dagger}^j(\nbigt_X)$ are pure twistor structures
of weight $w+j$,
and that $\nbigs_X$ induces an isomorphism
$a^j_{X\dagger}(\nbigt_X)\simeq
a^{-j}_{X\dagger}(\nbigt_X)^{\ast}
\otimes\newTate(-w)$.
Hence, we have $\rank \Gr^{W'}_k=\rank\Gr^{W'}_{-k}$
for any $k\in\seisuu$.
We also have
$\rank \Gr^{W(\omega_1)}_k=\rank\Gr^{W(\omega_1)}_{-k}$
for any $k\in\seisuu$.
Then, it is easy to see that $k_0<0$.
Because $W_{k_0-1}(\omega_1)=W'_{k_0-1}$
and $W_{k_0}(\omega_1)\subset W'_{k_0}$,
we obtain $\rank\Gr^{W(\omega_1)}_{k_0}<\rank\Gr^{W'}_{k_0}$.
We obtain $\rank W(\omega_0)_{-k_0-1}>\rank W'_{-k_0-1}$
which contradicts with $W(\omega_0)_{-k_0-1}\subset W'_{-k_0-1}$.
\hfill\qed

\vspace{.1in}

We set
$P_{\omega_2}\nbigt^{(0)}_{k,0}:=
\Ker\bigl(
L^{(0)}_{\omega_2}:
\nbigt^{(0)}_{k,0}\to
\nbigt^{(0)}_{k,-2}\otimes\newTate(-1)
\bigr)$.
Then,
$P_{\omega_2}\nbigt^{(0)}_{\bullet,0}:=
\bigoplus_kP_{\omega_2}\nbigt^{(0)}_{k,0}$
with the induced morphisms
$L^{(0)}_{\omega_1}$
and $P_{\omega_2}\nbigs^{(0)}_0$
is a graded Lefschetz twistor structure of weight $w$ and type $-1$.
By Lemma \ref{lem;22.3.4.130},
the morphism
\begin{equation}
 L^{(0)}_{\omega_2}:
 \Gr^{W(\omega_1)}_{\ell}
 \Bigl(
 \bigoplus_{k}
 \abar_{\Xtilde\dagger}^k(\nbigt)
 \Bigr)
 \lrarr
 \Gr^{W(\omega_1)}_{\ell-2}
 \Bigl(
 \bigoplus_{k}
 \abar_{\Xtilde\dagger}^k(\nbigt)
 \Bigr)
 \otimes\newTate(-1)
\end{equation}
is identified with
\[
\bigoplus_{k+j=-\ell}
 \abar_{X\dagger}^k
 \fbar_{\dagger}^j(\nbigt)
 \lrarr
  \bigoplus_{k+j=-\ell+2}
 \abar_{X\dagger}^k
 \fbar_{\dagger}^j(\nbigt)
 \otimes
 \newTate(-1)
\]
induced by
$L_{\omega_2}:
\fbar_{\dagger}^j\nbigt
\lrarr
\fbar_{\dagger}^{j+2}\nbigt
\otimes\newTate(-1)$.
Hence,
the graded Lefschetz twistor structure
$\Bigl(
\bigoplus_{j}
\abar_{X\dagger}^j(\nbigt_X),
L_{\omega_1}
\Bigr)$
of weight $w$ and type $-1$
is a direct summand of
the graded Lefschetz twistor structure
$\bigl(
P_{\omega_2}\nbigt^{(0)}_{\bullet,0},
L^{(0)}_{\omega_1}
\bigr)$
of weight $w$ and type $-1$,
and
$a_{X\dagger}\nbigs_X$
is induced by
$P_{\omega_2}\nbigs^{(0)}_0$.
We obtain that
$\Bigl(
\bigoplus_{j}
\abar_{X\dagger}^j(\nbigt_X),
L_{\omega_1},a_{X\dagger}\nbigs_X
\Bigr)$
is also a polarized graded Lefschetz twistor structure
of weight $w$ and type $-1$.
(For example, see \cite[Lemma 8.2.7]{Mochizuki-MTM}.)
Thus, the proof of Proposition \ref{prop;22.3.4.30}
is completed.
\hfill\qed

\subsubsection{Appendix}

We note the following lemma
though we do not use it in the proof.

\begin{lem}
The following morphisms are isomorphisms:
\begin{equation}
\label{eq;22.3.4.12}
  \fbar^{j-1}_{\dagger}(\nbigk)
\lrarr
\fbar^{j}_{\dagger}(\nbigt),
\,\,(j\leq -1).
\end{equation}
 \begin{equation}
\label{eq;22.3.4.11}
  \fbar^j_{\dagger}(\nbigt)
 \lrarr
 \fbar^{j+1}_{\dagger}(\nbigk)
 \otimes\newTate(-1),
 \,\,
 (j\geq 1).
 \end{equation}
The following morphism is a monomorphism:
\begin{equation}
\label{eq;22.3.4.13}
 \fbar^{-1}_{\dagger}(\nbigk)
\lrarr
\fbar^{0}_{\dagger}(\nbigt).
\end{equation}
 The following morphism is an epimorphism
  \begin{equation}
\label{eq;22.3.4.14}
  \fbar^0_{\dagger}(\nbigt)
 \lrarr
 \fbar^{1}_{\dagger}(\nbigk)
 \otimes\newTate(-1).
 \end{equation}
\end{lem}
\pf
By the Hard Lefschetz theorem for $\nbigk$ and $f$,
the morphisms
$L_{\omega_2}:f_{\dagger}^{j-1}(\nbigk)\lrarr f_{\dagger}^{j+1}(\nbigk)$
are epimorphisms for $j\geq 0$.
By Lemma \ref{lem;22.3.4.10},
$f^{j}_{\dagger}(\nbigt)
\lrarr
f^{j+1}_{\dagger}(\nbigk)$
are epimorphisms for $j\geq 0$.
In particular, (\ref{eq;22.3.4.14}) is an epimorphism.

Let us observe that
$f^j_{\dagger}(\nbigt[!\Htilde])=0$ $(j>0)$.
It is enough to check it locally around any point $P$ of $X$.
There exists an open neighbourhood $X_P$ of $P$
with a finite family of hypersurfaces
$H_{P,i}\subset X_P$ $(i\in\Lambda)$
such that
$X_P\cap f(|\Htilde|)=\bigcap_{i\in\Lambda}H_{P,i}$.
We set $\Xtilde_P:=f^{-1}(X_P)$
and $\Htilde_{P,i}:=f^{-1}(H_{P,i})\subset\Xtilde_P$.
We also set
$\nbigt_P:=\nbigt_{|\Xtilde_P}$
and
$\Htilde_P:=\Htilde\cap \Xtilde_P$.
We have
$|\Htilde_P|=
\bigcap_{i\in\Lambda}\Htilde_{P,i}$.
For any $J\subset \Lambda$,
we set
$\Htilde_P(J)=\bigcup_{i\in J}\Htilde_{P,i}$.
For $k\leq 0$,
we set
\[
 \nbigc^k=\bigoplus_{|J|=-k+1}
 \nbigt_P[!\Htilde_P(J)]
 \otimes
 \cnum(J)^{\lor}.
\]
(See \S\ref{subsection;22.2.4.30} for $\cnum(J)^{\lor}$.)
As in \S\ref{subsection;22.4.28.10},
we obtain the complex $\nbigc^{\bullet}$.
The morphisms
$\nbigt_P[!\Htilde_{P,i}]\to
 \nbigt_P[!\Htilde_P]$
induce a quasi-isomorphism
$\nbigc^{\bullet}\to \nbigt_P[!\Htilde_P]$.
Let $f_P:\Xtilde_P\to X_P$
denote the induced morphism.
Because $f_P$ induces isomorphisms
$\Xtilde_P\setminus\Htilde_{P,i}\simeq
 X_P\setminus H_{P,i}$,
we obtain
$f^j_{P\dagger}\nbigt_P[!\Htilde_{P,i}]=0$  $(j\neq 0)$.
Hence, we obtain
$f_{P\dagger}^j(\nbigt_P[!\Htilde_P])=0$ $(j>0)$.

Let $\nbigktilde$ denote the kernel of
$\nbigt[!\Htilde]\lrarr\nbigt$.
Then, we obtain the following isomorphisms:
\[
 f_{\dagger}^j\nbigt
 \simeq
 f_{\dagger}^{j+1}\nbigktilde
 \quad(j\geq 1).
\]
We have the following exact sequence of mixed twistor $\nbigd$-modules
on $X$:
\[
 f_{\dagger}^{j+1}
 W_{w-2}\nbigktilde
 \lrarr
  f_{\dagger}^{j+1}
  \nbigktilde
  \lrarr
  f_{\dagger}^{j+1}\nbigk
  \lrarr
   f_{\dagger}^{j+2}
 W_{w-2}\nbigktilde.
\]
By the Hard Lefschetz theorem for projective morphisms
in \cite{sabbah2,mochi2,Mochizuki-wild},
we have
\[
 f_{\dagger}^{j+1}
 W_{w-2}\nbigktilde
 =W_{w+j-1}\Bigl(
 f_{\dagger}^{j+1}
 W_{w-2}\nbigktilde
  \Bigr),
  \quad
 f_{\dagger}^{j+1}
  \nbigktilde
 =W_{w+j}\Bigl(
 f_{\dagger}^{j+1}
 \nbigktilde
  \Bigr),
\]
and $f_{\dagger}^{j+1}\nbigk$ is pure of weight $w+j$.
Hence, the following induced morphisms are monomorphisms:
\[
 f_{\dagger}^j\nbigt
 =\Gr^W_{w+j} f_{\dagger}^j\nbigt
 \lrarr
 \Gr^W_{w+j}
  f_{\dagger}^{j+1}
  \nbigktilde
  \lrarr
   \Gr^W_{w+j}
  f_{\dagger}^{j+1}
  \nbigk
 =  f_{\dagger}^{j+1}
  \nbigk.
\]
Thus, we obtain that
(\ref{eq;22.3.4.11}) are isomorphisms.
By using the Hermitian adjoint,
we obtain that
(\ref{eq;22.3.4.13}) is a monomorphism,
and that (\ref{eq;22.3.4.12}) are isomorphisms.
\hfill\qed

\subsection{End of proof of Theorem \ref{thm;22.2.25.40}}

Let $(\nbigt,\nbigs)$ be a polarized pure twistor $\nbigd$-module
of weight $w$ on $X$ with compact support.
It is enough to consider the case where
$\nbigt$ has a strict support denoted by $Z_0$.
We use the induction on $\dim Z_0$.
We assume that we have already proved in the lower dimensional case.
According to \cite{mochi2},
there exists a closed complex analytic subvariety
$Z_1\subset Z_0$
such that the following holds.
\begin{itemize}
 \item $Z_0\setminus Z_1$
       is a locally closed complex submanifold of $X$.
 \item There exists
       a tame polarized variation of pure twistor structure
       $(\nbigv,\nbigs)$
       of weight $w$ on $(Z_0,Z_1)$
       such that
       $\iota_{\supp(\nbigt)\setminus Z\,\dagger}(\nbigv,\nbigs)
       =\nbigt_{|X\setminus Z_1}$.
\end{itemize}
According to \cite{Wlodarczyk},
there exists a birational morphism $f:\Xtilde\lrarr X$
such that the following holds.
\begin{itemize}
 \item $f$ is the composition of a finite sequence of
       blow ups along smooth centers.
       Moreover, $f$ induces
       $\Xtilde\setminus f^{-1}(Z_1)\simeq X\setminus Z_1$.
 \item The strict transform $\Ztilde_0$ of $Z_0$
       is a complex closed submanifold of $\Xtilde$.
 \item $f^{-1}(Z_1)$ is a simple normal crossing hypersurface of $\Xtilde$,
       and $f^{-1}(Z_1)\cap \Ztilde_0$
       is a simple normal crossing hypersurface of $\Ztilde_0$.
\end{itemize}
We obtain the tame polarized variation of pure twistor structure
$(\nbigvtilde,\nbigstilde)$ on
$(\Ztilde_0,f^{-1}(Z_1)\cap\Ztilde_0)$
from $(\nbigv,\nbigs)$.
Let $(\nbigt_{\Xtilde},\nbigs_{\Xtilde})$
denote the regular pure twistor $\nbigd$-module of weight $w$
on $\Xtilde$
corresponding to $(\nbigvtilde,\nbigstilde)$.
Then,
$\nbigt$ is a direct summand of
$f_{\dagger}^0(\nbigt_{\Xtilde})$,
and $\nbigs$ is induced by $\nbigs_{\Xtilde}$.
(See Theorem \ref{thm;22.4.20.20} and Theorem \ref{thm;22.4.18.30}.)
By Theorem \ref{thm;22.3.15.20}
and Proposition \ref{prop;22.3.4.30},
we obtain the claim of Theorem \ref{thm;22.2.25.40}
for $(\nbigt,\nbigs)$.
\hfill\qed

\subsection{Comparison of the Hodge structures}
\label{subsection;22.3.16.30}

We revisit the theorem of Kashiwara and Kawai
\cite[Theorem 1]{Kashiwara-Kawai-Hodge-holonomic}
in the Hodge case.
Let $X$ be a compact K\"ahler manifold
with a simple normal crossing hypersurface $H$.
We fix a Poincar\'{e} like K\"ahler metric $g_{X\setminus H}$
of $X\setminus H$ (see \S\ref{subsection;22.4.2.21}).

Let $(V=\bigoplus_{p+q=w} V^{p,q},\nabla,\langle\cdot,\cdot\rangle)$
be a polarized variation of complex Hodge structure \cite{s1}
on $X\setminus H$:
\begin{itemize}
 \item $\nabla$ is a flat connection
       satisfying the Griffiths transversality,
       i.e.,
\[
       \nabla V^{p,q}
       \subset
       (V^{p-1,q+1}\oplus V^{p,q})\otimes\Omega^{1,0}
       \oplus
       (V^{p+1,q-1}\oplus V^{p,q})\otimes\Omega^{0,1}.
\]
 \item $\langle\cdot,\cdot\rangle$
       is a $\nabla$-flat sesqui-linear pairing.
      The decomposition
       $V=\bigoplus V^{p,q}$ is orthogonal with respect to
       $\langle\cdot,\cdot\rangle$.
       Moreover,
       $h=\bigoplus (\sqrt{-1})^{p-q}
       \langle\cdot,\cdot\rangle_{|V^{p,q}}$
       is positive definite.
\end{itemize}
By the Griffiths transversality,
we obtain the decomposition
$\nabla=\nabla_h+\theta+\theta^{\dagger}$,
where
$\nabla_h$ is a connection preserving the Hodge decomposition,
$\theta$ is a section of
$\bigoplus \Hom(V^{p,q},V^{p-1,q+1})\otimes\Omega^{1,0}$
and
$\theta^{\dagger}$ is a section of
$\bigoplus \Hom(V^{p,q},V^{p+1,q-1})\otimes\Omega^{0,1}$.
Note that $\nabla_h$ is unitary with respect to $h$,
and $\theta+\theta^{\dagger}$ is self adjoint with respect to $h$.
We have the decomposition
$\nabla_h=\delbar_V+\del_V$ into
the $(0,1)$-part and the $(1,0)$-part.
Then,
$(V,\delbar_V,\theta,h)$ is a harmonic bundle.
We obtain the holomorphic vector bundle
$\nbigv=(p_{\lambda}^{-1}(V),
p_{\lambda}^{\ast}(\delbar_V)+\lambda \theta^{\dagger})$,
which is equipped with the family of flat $\lambda$-connections
$\DD=\delbar_V+\lambda\theta^{\dagger}+\lambda\del_V+\theta$.
As the specialization to $\{\lambda\}\times (X\setminus H)$,
we obtain the $\lambda$-flat bundle
$(\nbigv^{\lambda},\DDlambda)$ with
the pluri-harmonic metric $h$ on $X\setminus H$.
The flat bundle 
$(\nbigv^{1},\DD^1)$
is equal to the original flat bundle
$(V,\nabla)$ on $X\setminus H$.

\subsubsection{Hodge filtration of Cattani-Kaplan-Schmid and Kashiwara-Kawai}
\label{subsection;22.4.11.21}

Let $V_{\min}$ be the regular singular $\nbigd_X$-module
obtained as the minimal extension of $(V,\nabla)$.
Let $L_V$ denote the local system on $X\setminus H$
obtained as the sheaf of flat sections of $(V,\nabla)$.
Let $\Omega^{\bullet}_X$ denote the complex of
the sheaves of holomorphic differential forms on $X$.
Note that $V_{\min}\otimes\Omega^{\bullet}_X$
is naturally isomorphic to the intersection complex of $L_V$.
As proved in \cite{cks2} and \cite{k3},
$\nbigc^{\bullet}_{L^2}(V,\nabla,h)$
satisfies the condition of intersection complex of $L_V$,
and hence there exists a unique isomorphism
$V_{\min}\otimes\Omega^{\bullet}_X\simeq
\nbigc^{\bullet}_{L^2}(V,\nabla,h)$
in the derived category of cohomologically constructible sheaves on $X$.
We obtain the isomorphism
\begin{equation}
\label{eq;22.4.11.1}
 H^k\bigl(X,V_{\min}\otimes\Omega^{\bullet}_X\bigr)
 \simeq
 H^k\bigl(
 X,\nbigc^{\bullet}_{L^2}(V,\nabla,h)
 \bigr).
\end{equation}

Let $\Harm^k(V)$ denote the space of
$L^2$-sections $\tau$ of $V\otimes \Tot^k\Omega^{\bullet,\bullet}$
such that $\Delta_V\tau=0$,
where $\Delta_V$ denote the Laplacian
induced by $\nabla$ and $g_{X\setminus H}$.
As explained in Proposition \ref{prop;22.2.25.10},
it is equal to
the space of $L^2$-sections $\tau$
of $V\otimes \Tot^k\Omega^{\bullet,\bullet}$
such that $\nabla\tau=0$ and $\nabla^{\ast}\tau=0$,
where $\nabla^{\ast}$ denotes the adjoint of $\nabla$
with respect to $h$ and $g_{X\setminus H}$.
Because $\bigoplus H^k(X,\nbigc^{\bullet}_{L^2}(V,\nabla,h))$
is finite dimensional,
we obtain that
\begin{equation}
\label{eq;22.4.11.2}
 H^k(X,\nbigc^{\bullet}_{L^2}(V,\nabla,h))
 \simeq
 \Harm^k(V).
\end{equation}

For $(p,q)\in\seisuu^2$ satisfying $p+q=w+k$,
we set
\[
\bigl(
V\otimes\Tot^k\Omega^{\bullet,\bullet}
\bigr)^{p,q}
:=\bigoplus_j
V^{p-j,q-k+j}
\otimes
\Omega^{j,k-j}.
\]
We obtain the decomposition
\[
V\otimes\Tot^k\Omega^{\bullet,\bullet}
=\bigoplus_{p+q=w+k}
\bigl(
V\otimes\Tot^k\Omega^{\bullet,\bullet}
\bigr)^{p,q}.
\]
We set
$\Harm^{p,q}(V)\subset\Harm^k(V)$
denote the subspace of $\tau\in\Harm^k(V)$
which is an $L^2$-section of
$\bigl(
V\otimes\Tot^k\Omega^{\bullet,\bullet}
\bigr)^{p,q}$.
As proved in \cite[\S6]{k3},
we obtain the decomposition
\begin{equation}
\label{eq;22.4.10.3}
 \Harm^k(V)=\bigoplus_{p+q=w+k}\Harm^{p,q}(V).
\end{equation}
By (\ref{eq;22.4.11.1}), (\ref{eq;22.4.11.2}) and (\ref{eq;22.4.10.3}),
we obtain
\[
 H^k(X,V_{\min}\otimes\Omega^{\bullet})
 \simeq
 \bigoplus_{p+q=k+w}
 \Harm^{p,q}(V).
\]
By using the isomorphism,
we set
\[
F^pH^k(X,V_{\min}\otimes\Omega^{\bullet})
=\bigoplus_{p'\geq p}
 \Harm^{p,q}(V),
 \quad\quad
 \Fbar^qH^k(X,V_{\min}\otimes\Omega^{\bullet})
=\bigoplus_{q'\geq q}
 \Harm^{p,q}(V).
\]
Kashiwara and Kawai \cite{k3}
proved that
the filtrations
$F^{\bullet}H^k(X,V_{\min}\otimes\Omega^{\bullet})$
and
$\Fbar^{\bullet}H^k(X,V_{\min}\otimes\Omega^{\bullet})$
are independent of the choice
of a polarization $\langle\cdot,\cdot\rangle$
and a Poincar\'{e} like K\"ahler metric $g_{X\setminus H}$.
As a result, the Hodge decomposition
\[
 H^k(X,V_{\min}\otimes\Omega^{\bullet})
=\bigoplus_{p+q=k+w}
 F^p H^k(X,V\otimes\Omega^{\bullet})
 \cap
 \Fbar^q H^k(X,V_{\min}\otimes\Omega^{\bullet})
\]
is also independent of the choice of
$\langle\cdot,\cdot\rangle$
and $g_{X\setminus H}$.

We set
$F_pH^k(X,V_{\min}\otimes\Omega^{\bullet})
=F^{-p}H^k(X,V_{\min}\otimes\Omega^{\bullet})$
and
$\Fbar_qH^k(X,V_{\min}\otimes\Omega^{\bullet})
=\Fbar^{-q}H^k(X,V_{\min}\otimes\Omega^{\bullet})$.

\subsubsection{Theorem of Kashiwara-Kawai}

Let $H=\bigcup_{i\in\Lambda}H_i$ be the irreducible decomposition.
Because $(V,\nabla)$ underlies a polarized variation of Hodge structure,
any eigenvalue $\alpha$
of the local monodromy of $(V,\nabla)$
around $H_i$ satisfies $|\alpha|=1$.
Let $(\nbigp V,\nabla)$ denote the regular singular meromorphic
flat bundle on $(X,H)$
obtained as the extension of $(V,\nabla)$.
Let $\nbigp_{<\veciti_{\Lambda}}V\subset\nbigp V$
be the lattice satisfying the following conditions.
\begin{itemize}
 \item $\nbigp_{<\veciti_{\Lambda}}V$ is logarithmic
       with respect to $\nabla$.
 \item For any $i\in \Lambda$,
       any eigenvalue $\beta$ of $\Res_i(\DD^1)$
       satisfies
       $-1<\beta\leq 0$.
\end{itemize}
Note that
$V_{\min}=
\nbigd_X\cdot \nbigp_{<\veciti_{\Lambda}}(V)
\subset\nbigp V$.

We set $F^p(V)=\bigoplus_{p'\geq p} V^{p',q'}$
$(p\in\seisuu)$.
They are holomorphic subbundle of $(V,\nabla^{0,1})$,
called the Hodge filtration
of the variation of Hodge structure.
For any holomorphic section $v$ of $F^p(V)$,
$\nabla(v)$ is a holomorphic section of
$F^{p-1}(V)\otimes\Omega^1$.
By the nilpotent orbit theorem \cite{sch},
$F^{p}(V)$
extend to holomorphic subbundles
$F^p(\nbigp_{<\veciti_{\Lambda}}V)$ of $\nbigp_{<\veciti_{\Lambda}}V$.

We set $F_j(V):=F^{-j}(V)$
and $F_j(\nbigp_{<\veciti_{\Lambda}}V)
:=F^{-j}(\nbigp_{<\veciti_{\Lambda}}V)$.
Let $F_j\nbigd_X\subset\nbigd_X$ denote
the sheaf of differential operators whose orders are less than $j$.
We set
\[
F_p(V_{\min}):=
\sum_{i+j=p}F_i(\nbigd_X)\cdot F_j(\nbigp_{<\veciti_{\Lambda}}V)
\subset
V_{\min}.
\]
In this way,
we obtain the filtered $\nbigd_X$-module $(V,F)$.
The differential of 
$V_{\min}\otimes\Omega^{\bullet}$
induces
$F_p(V_{\min})\otimes\Omega^{k}_X
\lrarr
F_{p+1}(V_{\min})\otimes\Omega^{k+1}_X$.
Hence, 
$\bigoplus_k
F_{p+k}(V_{\min})
\otimes\Omega^k_X$
is a subcomplex of
$V_{\min}\otimes\Omega^{\bullet}$,
which is denoted by
$F_p(V_{\min}\otimes\Omega^{\bullet}_X)$.
Kashiwara and Kawai announced the following theorem
in \cite[Theorem 1]{Kashiwara-Kawai-Hodge-holonomic},
which we shall revisit in \S\ref{subsection;22.4.10.2}.

\begin{thm}[Kashiwara-Kawai]
\label{thm;22.4.10.11}
The natural morphism
\[
 H^{k}(X,F_p(V_{\min}\otimes\Omega^{\bullet}_X))
 \lrarr
 H^{k}(X,V_{\min}\otimes\Omega^{\bullet}_X)
\]
is injective,
and the image is equal to
$F_pH^k(X,V_{\min}\otimes\Omega^{\bullet}_{X})$. 
\end{thm}

Let $X^{\dagger}$ and $H^{\dagger}$
denote the conjugate of $X$ and $H$, respectively.
We may naturally regard $(V,\nabla)$
as a flat bundle on $X^{\dagger}\setminus H^{\dagger}$.
We have the $\nbigd_{X^{\dagger}}$-module
$V^{\dagger}_{\min}$ on $X^{\dagger}$
obtained as the minimal extension of $(V,\nabla)$.
Because 
$V^{\dagger}_{\min}\otimes\Omega^{\bullet}_{X^{\dagger}}$
is the intersection complex of the local system $L_V$,
we have the natural isomorphism
\[
 H^{\ast}(X^{\dagger},
 V^{\dagger}_{\min}\otimes\Omega^{\bullet}_{X^{\dagger}})
 \simeq
 H^{\ast}(X,V_{\min}\otimes\Omega^{\bullet}_X).
\]
From the Hodge filtration
$\Fbar^q(V)=\bigoplus_{q'\geq q}V^{p,q}$ $(q\in\seisuu)$,
we obtain the filtration
$\Fbar_q(V^{\dagger}_{\min}\otimes\Omega^{\bullet}_{X^{\dagger}})$
$(q\in\seisuu)$.
Theorem {\rm\ref{thm;22.4.10.11}} gives us a similar description of
the Hodge filtration
$\Fbar_{\bullet}H^{k}(X,V_{\min}\otimes\Omega^{\bullet}_X)
 =\Fbar_{\bullet}H^{k}(X^{\dagger},
 V^{\dagger}_{\min}\otimes\Omega^{\bullet}_{X^{\dagger}})$
as the filtration associated with the spectral sequence of 
the filtered complex
$\Fbar_{\bullet}(V^{\dagger}_{\min}\otimes\Omega^{\bullet}_{X^{\dagger}})$.

\subsubsection{Comparison}
\label{subsection;22.4.10.2}

We consider the $\cnum^{\ast}$-action on
$\nbigx=\cnum_{\lambda}\times X$
defined by $a(\lambda, x)=(a\lambda,x)$ $(a\in\cnum^{\ast})$.
We consider the $\cnum^{\ast}$-action on
$p_{\lambda}^{-1}(V^{p,q}\otimes\Omega^{r,s})$
satisfying
$a^{\ast}\bigl(
p_{\lambda}^{-1}(v^{p,q}\otimes \tau^{r,s})
\bigr)
=a^{-p-r}p_{\lambda}^{-1}(v^{p,q}\otimes\tau^{r,s})$,
where $v^{p,q}$ and $\tau^{r,s}$ denote
$C^{\infty}$-sections of
$V^{p,q}$ and $\Omega^{r,s}$, respectively.
Then, $(\nbigv,\DD)$ is $\cnum^{\ast}$-equivariant.
We note that
$a^{\ast}(p_{\lambda}^{-1}h)$
and $h$ are mutually bounded.
Hence,
$\nbigc^{\bullet}_{L^2}(\nbigv,\DD,h)$
is naturally $\cnum^{\ast}$-equivariant.

Because of the decomposition (\ref{eq;22.4.10.3}),
there exists the $\cnum^{\ast}$-action on $\Harm^{k}(V)$
defined by
$a^{\ast}(\tau^{p,q})=a^{-p}\tau^{p,q}$.
It induces the $\cnum^{\ast}$-action
on $\Harm^k(V)\otimes\nbigo_{\cnum_{\lambda}}$.
The natural morphism
\[
 \Harm^k(V)\otimes\nbigo_{\cnum_{\lambda}}
 \lrarr
 (\id\times a_X)_{\ast}
 \nbigc^k_{L^2}(\nbigv,\DD,h)
\]
is $\cnum^{\ast}$-equivariant.
Hence, the following isomorphism is $\cnum^{\ast}$-equivariant:
\[
 \Harm^k(V)\otimes\nbigo_{\cnum_{\lambda}}
\simeq
 R^k(\id\times a_X)_{\ast}
 \nbigc^{\bullet}_{L^2}(\nbigv,\DD,h).
\]

Let $\nbigr_F(V)$
denote the Rees module of $(V,F)$,
i.e.,
the analytification of the Rees module
$\sum F_j(V)\lambda^{j}=\sum F^{-j}(V)\lambda^j$.
We have the natural $\cnum^{\ast}$-action
defined by the pull back,
i.e.,
$a^{\ast}(\lambda^{j}v^{p,q})=a^{j}\lambda^{j}v^{p,q}$.
There exists the $\cnum^{\ast}$-equivariant isomorphism
$\nbigr_F(V)\simeq \nbigv$
induced by
$\lambda^{-p}v^{p,q}\longmapsto p_{\lambda}^{-1}(v^{p,q})$,
under which $\nabla$ on $\nbigr_F(V)$
is identified with $\DD^f$ on $\nbigv$.

We have the $\nbigr_X$-module $\gbigv\subset\nbigp\nbigv$
underlying the pure twistor $\nbigd_X$-module
associated with the harmonic bundle $(V,\delbar_V,\theta,h)$.
(See \S\ref{subsection;22.4.12.1}.)
We also obtain the $\nbigr_X$-module
$\nbigr_F(V_{\min})$
as the analytification of the Rees module
$\sum F_j(V_{\min})\lambda^j
=\sum F^{-j}(V_{\min})\lambda^j$.
By the construction,
the $\cnum^{\ast}$-equivariant isomorphism
$\nbigr_F(V)\simeq \nbigv$
uniquely extends to a $\cnum^{\ast}$-equivariant isomorphism
$\nbigr_F(V_{\min})\simeq\gbigv$.

Let $\nbigr_F(V_{\min}\otimes\Omega^{\bullet}_X)$
denote the analytification of the Rees module
of the filtered complex
$(V_{\min}\otimes\Omega^{\bullet}_X,F)$.
It is naturally isomorphic to
$\nbigr_F(V_{\min})\otimes
\Omegatilde^{\bullet}_{\nbigx/\cnum}$,
where
$\Omegatilde^{\bullet}_{\nbigx/\cnum}
=\bigwedge^{\bullet}(\lambda^{-1}p_{\lambda}^{\ast}\Omega^1_X)$.
Hence, we obtain the $\cnum^{\ast}$-equivariant isomorphism
\[
 \nbigr_F(V_{\min}\otimes
 \Omega^{\bullet}_X)
 \simeq
 \gbigv\otimes\Omegatilde^{\bullet}_{\nbigx/\cnum}.
\]
Here, we consider the $\cnum^{\ast}$-action on
$\Omegatilde^{\bullet}_{\nbigx/\cnum}$
given by $a^{\ast}(\lambda^{j}p_{\lambda}^{\ast}\tau^k)
=a^{j}\lambda^{j}p_{\lambda}^{\ast}\tau^k$.
\begin{lem}
\label{lem;22.4.11.12}
The natural morphism
\begin{equation}
\label{eq;22.4.10.10}
 H^k(X,F_p(V_{\min}\otimes\Omega^{\bullet}_X))
 \lrarr
 H^k(X,V_{\min}\otimes\Omega^{\bullet}_X)
\end{equation}
is injective.
\end{lem}
\pf
For any $\lambda\in\cnum$,
let $\iota_{\lambda}$ denote the inclusion
$\{\lambda\}\to\cnum$,
and let $\iota_{\lambda,X}$ denote the inclusion
$\{\lambda\}\times X\to \nbigx$.
Because
$R^k(\id\times a_X)_{\ast}
\Bigl(
\nbigr_F(V_{\min}\otimes
\Omega^{\bullet}_X)
\Bigr)$ $(k\in\seisuu)$
are locally free,
we have
\[
 \iota_{\lambda}^{\ast}
R^k(\id\times a_X)_{\ast}
\Bigl(
\nbigr_F(V_{\min}\otimes
\Omega^{\bullet}_X)
\Bigr)
=H^k\Bigl(
 X,\iota_{\lambda,X}^{\ast}
 \nbigr_F(V_{\min}\otimes
\Omega^{\bullet}_X)
\Bigr).
\]
We have 
$\iota_{0,X}^{\ast}
 \nbigr_F(V_{\min}\otimes
\Omega^{\bullet}_X)
=\Gr^F(V_{\min}\otimes\Omega^{\bullet}_X)$
and
$\iota_{1,X}^{\ast}
 \nbigr_F(V_{\min}\otimes
\Omega^{\bullet}_X)
=V_{\min}\otimes\Omega^{\bullet}_X$.
Because
\[
 \dim H^k\bigl(
 X,\Gr^F(V_{\min}\otimes\Omega^{\bullet}_X)
 \bigr)
 =\dim H^k\bigl(
 X,V_{\min}\otimes\Omega^{\bullet}_X
 \bigr)
\]
for any $k\in\seisuu$, 
we obtain the claim of the lemma
by the standard argument using the spectral sequence
associated with the filtered complex.
\hfill\qed

\vspace{.1in}

Let
$F'_pH^k(X,V_{\min}\otimes\Omega^{\bullet}_X)$
denote the image of (\ref{eq;22.4.10.10}),
and let
$\nbigr_{F'}
H^k(X,V_{\min}\otimes\Omega^{\bullet}_X)$
denote the analytification of the Rees module
$\sum F'_j
H^k(X,V_{\min}\otimes\Omega^{\bullet}_X)\lambda^{j}$.
By Lemma \ref{lem;22.4.11.12},
we have
\begin{equation}
\label{eq;22.4.11.20}
 R^k(\id\times a_X)_{\ast}\Bigl(
 \nbigr_F(V_{\min}\otimes\Omega^{\bullet}_X)
 \Bigr)
 \simeq
 \nbigr_{F'}
 H^k(X,V_{\min}\otimes\Omega^{\bullet}_X). 
\end{equation}

We consider the $\cnum^{\ast}$-action
on $p_{\lambda}^{\ast}\Omega^j$ by $a^{\ast}(\tau^j)=a^{-j}\tau^j$.
By the construction of the morphisms in \S\ref{section;22.4.10.3},
the following isomorphism is $\cnum^{\ast}$-equivariant:
\begin{equation}
\label{eq;22.4.11.10}
R^k(\id\times a_X)_{\ast}\Bigl(
\gbigv\otimes p_{\lambda}^{\ast}(\Omega_X^{\bullet})
\Bigr)
\simeq
R^k(\id\times a_X)_{\ast}\Bigl(
\nbigc^{\bullet}_{\tw}(\gbige)
\Bigr)
\simeq
\Harm^k(V)\otimes\nbigo_{\cnum_{\lambda}}.
\end{equation}

We have the $\cnum^{\ast}$-equivariant isomorphism
$\Omegatilde^{\bullet}_{\nbigx/\cnum}
\simeq
 p_{\lambda}^{\ast}\Omega^{\bullet}_X$
induced by
$\lambda^{-j}p_{\lambda}^{\ast}\tau\longmapsto p_{\lambda}^{\ast}\tau$
for holomorphic $j$-forms $\tau$.
By replacing
$p_{\lambda}^{\ast}\Omega^{\bullet}_X$
with $\Omegatilde^{\bullet}_{\nbigx/\cnum}$,
we obtain
\begin{equation}
\label{eq;22.4.11.11}
R^k(\id\times a_X)_{\ast}\Bigl(
\gbigv\otimes \Omegatilde_{\nbigx/\cnum}^{\bullet}
\Bigr)
\simeq
R^k(\id\times a_X)_{\ast}\Bigl(
\nbigctilde^{\bullet}_{\tw}(\gbige)
\Bigr)
\simeq
\Harm^k(V)\otimes\nbigo_{\cnum_{\lambda}}.
\end{equation}

We set
$F_j\Harm^k(V)=
F^{-j}\Harm^k(V)=\bigoplus_{p'\geq -j}\Harm^{p',q'}(V)$.
Let $\nbigr_F\Harm^k(V)$ denote the analytification
of the Rees module
$\sum \lambda^jF_j\Harm^k(V)=
\sum \lambda^{j}F^{-j}\Harm^k(V)$.
We have the $\cnum^{\ast}$-equivariant isomorphism
\[
 \nbigr_F\Harm^k(V)
 \simeq
 \Harm^k(V)\otimes\nbigo_{\cnum}
\]
induced by
$\lambda^{-p}\tau^{p,q}
\longmapsto \tau^{p,q}$.
By (\ref{eq;22.4.11.20}) and (\ref{eq;22.4.11.11}),
we obtain the following $\cnum^{\ast}$-equivariant isomorphism:
\begin{equation}
\nbigr_{F'}
H^k(X,V_{\min}\otimes\Omega^{\bullet}_X)
\simeq
\nbigr_F\Harm^k(V).
\end{equation}
We note that
the isomorphism
$\iota_1^{\ast}(\gbigv\otimes\Omega^{\bullet}_X)
\simeq
 \nbigc^{\bullet}_{L^2}(\nbigv^1,\DD^1,h)
 =\nbigc^{\bullet}_{L^2}(V,\nabla,h)$
induced by Corollary \ref{cor;22.2.18.10}
and Theorem \ref{thm;22.2.18.11}
is equal to the isomorphism
$V_{\min}\otimes\Omega^{\bullet}_X
\simeq
 \nbigc^{\bullet}_{L^2}(V,\nabla,h)$
because they are isomorphic to the intersection complexes.
Therefore,
we obtain that
$F_pH^k(X,V_{\min}\otimes\Omega^{\bullet}_X)
=F'_pH^k(X,V_{\min}\otimes\Omega^{\bullet}_X)$,
and Theorem \ref{thm;22.4.10.11} is proved.
\hfill\qed

\subsection{First Chern class (Appendix)}

Let $F:X\to Y$ be a morphism of complex manifolds.
Let $M$ be a holonomic $\nbigd_X$-module
whose support is proper over $Y$.
We obtain the holonomic $\nbigd_Y$-modules $F_{\dagger}^j(M)$
$(j\in\seisuu)$.

\subsubsection{The action of a $2$-cohomology class}

Let us recall that
$c\in H^2(X,\cnum)$ induces morphisms
\begin{equation}
\label{eq;22.3.3.1}
 c:F_{\dagger}^j(M)\lrarr F_{\dagger}^{j+2}(M).
\end{equation}
Let
$\iota_F:X\lrarr X\times Y$ denote the graph embedding,
and let $\pi_F:X\times Y\lrarr Y$ denote the projection.
We have
\[
 F_{\dagger}(M)=
 \pi_{F\dagger}(\iota_{F\dagger}M)
 =\pi_{F\ast}
 \Tot
 \Bigl(
 \iota_{F\dagger}(M)
 \otimes
 \Omega_{X\times Y/Y}^{\bullet,\bullet}[\dim X]
 \Bigr).
\]
Let $\omega$ be a closed $2$-form representing $c$.
The multiplication of $\omega$
induces
\[
 \pi_{F\dagger}(\iota_{F\dagger}M)
 \lrarr
  \pi_{F\dagger}(\iota_{F\dagger}M)[2].
\]
Thus, we obtain (\ref{eq;22.3.3.1}).

\subsubsection{Weight filtrations}

Let $H$ be an effective divisor of $X$.
We have $M(\ast H)$ and $M(!H)$.
Let $\rho_0:M(!H)\lrarr M(\ast H)$
denote the natural morphism.
Let us observe that there exist the naturally defined filtrations
$W$ on $\Ker(\rho_0)$ and $\Cok(\rho_0)$.

First, let us consider the case where $H$ is smooth and reduced.
Let $P$ be any point of $H$.
There exists a neighbourhood $X_P$ with
a holomorphic coordinate system
$(z_1,\ldots,z_n)$
such that $H_P:=X_P\cap H=\{z_1=0\}$.
We set $M_P:=M_{|X_P}$.
We obtain
$\rho_{1,P}:\psi^{(1)}_{z_1}(M_P)\to \psi^{(0)}_{z_1}(M_P)$,
and there exist natural isomorphisms
$\Ker\rho_{0|X_P}\simeq \Ker \rho_{1,P}$
and $\Cok\rho_{0|X_P}\simeq\Cok\rho_{1,P}$.
Here, $\psi^{(i)}_{z_1}(M)$ denote the nearby cycle sheaf.
(See \cite{beilinson2}.
See \S\ref{subsection;22.4.4.20}
for the explanation in the context of $\nbigr$-modules,
and \cite[\S4.1]{Mochizuki-MTM} for more details.)

Under the isomorphism
$\psi^{(1)}_{z_1}(M_P)\simeq \psi^{(0)}_{z_1}(M_P)$,
we obtain the locally nilpotent endomorphism $N$
on $\psi^{(1)}_{z_1}(M_P)$.
We obtain the weight filtration $W$
on $\psi^{(1)}_{z_1}(M_P)$.
It induces a filtration $W$
on $\Ker\rho_0$ and $\Cok\rho_0$.
It is independent of the choice of
a coordinate system $(z_1,\ldots,z_n)$.
(For example, see \cite[Lemma 4.2.7]{Mochizuki-MTM}
for the dependence of Beilinson functors
on a coordinate system in the context of $\nbigr$-triples.)
Hence, the filtration $W$ on $\Ker(\rho_0)$ and $\Cok(\rho_0)$
are globally well defined.

Let us consider the general case.
Let $L(H)$ denote the holomorphic line bundle on $X$
such that the sheaf of holomorphic sections of $L(H)$ is $\nbigo_X(H)$.
Let $\iota_1:X\lrarr L(H)$ denote the embedding induced by
the canonical morphism $\nbigo_X\lrarr\nbigo_X(H)$.
Let $\iota_0:X\lrarr L(H)$ denote the $0$-section.
We obtain the $\nbigd_{L_H}$-module $\iota_{1\dagger}(M)$.
We have
$\iota_{1\dagger}(M(\star H))
 \simeq
 \iota_{1\dagger}(M)(\star \iota_0(X))$ for $\star=\ast,!$.
Let
$\rhotilde_0:\iota_{1\dagger}(M)(!\iota_0(X))
\lrarr \iota_{1\dagger}(M)(\ast\iota_0(X))$
denote the natural morphism,
which equals $\iota_{1\dagger}(\rho_0)$.
Then, we have
\[
 \Ker(\rhotilde_0)\simeq
 \iota_{1\dagger}(\Ker\rho_0),
 \quad
 \quad
 \Cok(\rhotilde_0)\simeq
 \iota_{1\dagger}(\Cok\rho_0).
\]
There exists the filtration $W$
on $\Ker(\rho_0)$ and $\Cok(\rho_0)$
such that
$W_j\Ker(\rhotilde_0)
=\iota_{1\dagger}W_j\Ker(\rho_0)$
and
$W_j\Cok(\rhotilde_0)
=\iota_{1\dagger}W_j\Cok(\rho_0)$.

\subsubsection{The induced $\nbigd$-modules and some isomorphisms}
\label{subsection;22.4.21.11}

We set
$K_0=\Ker(\rho_0)/W_{-1}\Ker(\rho_0)$
and
$C_0=W_0\Cok(\rho_0)$.

\begin{lem}
There exists an isomorphism $K_0\simeq C_0$.
\end{lem}
\pf
Let us consider the case where $H$ is smooth.
For any $P\in H$,
let $(X_P,z_1,\ldots,z_n)$ be as above.
Then, there exist natural isomorphisms
$K_{0|X_P}\simeq P\Gr^W_0\psi^{(1)}_{z_1}(M_P)$
and
$C_{0|X_P}\simeq P\Gr^W_0\psi^{(0)}_{z_1}(M_P)$.
Hence, a natural isomorphism
$(\psi_{z_1}^{(1)}(M_P),N)
\simeq
(\psi_{z_1}^{(0)}(M_P),N)$
induces
$K_{0|X_P}\simeq C_{0|X_P}$,
which is independent of the choice of
$(z_1,\ldots,z_n)$.
(For example, see the proof of Lemma \ref{lem;22.4.20.1}.)
Hence, we obtain $K_{0}\simeq C_0$.
The general case is reduced to the above case
by the construction of the filtrations $W$.
\hfill\qed

\vspace{.1in}

For later use,
let us specify an isomorphism $\kappa:K_0\simeq C_0$
more explicitly.
It is enough to study the case where $H$ is smooth and reduced.
We fix a splitting
$s:\cnum/\seisuu\lrarr \seisuu$,
and a total order on $\cnum/\seisuu$.
We obtain a bijection
$\seisuu\times (\cnum/\seisuu)\simeq\cnum$
by $(n,x)\longmapsto n+s(x)$.
Let $\leq_1$ denote the total order on $\cnum$
induced by the lexicographic order.
Let $V\nbigd_X\subset\nbigd_X$ denote
the sheaf of subalgebras of differential operators $\nbigp$
satisfying $\nbigp\cdot \nbigo_X(-H)\subset\nbigo_X(-H)$.
Let $P\in H$ be any point.
Let $(X_P,z_1,\ldots,z_n)$ be as above.
Let $V_{\bullet}(M_P)$ denote the $V$-filtration along $z_1$
indexed by $(\cnum,\leq_1)$.
For $a\in\real$, we set
\[
 \psi_{a}(M_P):=V_{a}(M_P)/V_{<a}(M_P).
\]
Note that $-\del_{z_1}z_1-a$ is nilpotent on $\psi_a(M_P)$.
We consider $N_P=-z_1\del_{z_1}$ on $\psi_{-1}(M_P)$,
which is nilpotent.
We obtain the weight filtration $W$ of $N_P$
on $\psi_{-1}(M_P)$.
It induces filtrations
on $\Ker N_P$ and $\Cok N_P$.

We have $V_{-1}(M_P)=V_{-1}(M(!H)_P)$.
The natural morphism
$-\del_{z_1}:V_{-1}(M_P)\lrarr V_0(M(!H)_P)$ 
induces the isomorphism
$\can:\psi_{-1}(M_P)\simeq \psi_0(M(!H)_P)$.
It induces
\[
\can:\Ker N_P
\simeq
 \psi_0(\Ker\rho_{0|X_P})
 =V_0(\Ker\rho_{0|X_P})
 \subset \Ker(\rho_{0|X_P}).
\]

We regard $\Ker(N_P)$ as the $\nbigd_{H_P}$-module.
Let $\iota_{H_P}:H_P\lrarr X_P$ denote the inclusion.
There exists the isomorphism
$\iota_{H_P\dagger}\bigl(\Ker N_P\bigr)
 \simeq
 \Ker(\rho_{0|X_P})$
induced by the inclusion
\[
 \iota_{H_P\ast}\bigl(
 \Ker N_P\otimes (dz_1)^{-1}
 \bigr)
 \lrarr
 \Ker(\rho_{0|X_P}),
 \quad
 v\cdot (dz_1)^{-1}
 \longmapsto
 \can(v).
\]
We obtain
\[
\iota_{H_P\dagger}\Bigl(
 \Ker N_P/W_{-1}\Ker N_P
 \Bigr)
 \simeq
 K_{0|X_P}.
\]

We have $V_{-1}(M(\ast H)_P)=V_{-1}(M_P)$.
The morphism
$z_1:V_{0}(M(\ast H)_P)\lrarr V_{-1}(M_P)$
induces
the isomorphism
$\var:\psi_0(M(\ast H)_P)\simeq \psi_{-1}(M_P)$.
It induces an isomorphism
\[
\var:V_0(\Cok(\rho_{0|X_P}))\simeq \Cok(N_P).
\]
We regard $\Cok(N_P)$ as a $\nbigd_{H_P}$-module.
There exists the isomorphism
$\iota_{H_P\dagger}\bigl(\Cok N_P\bigr)
 \simeq
 \Cok(\rho_{0|X_P})$
induced by the following inclusion:
\[
 \iota_{H_P\ast}\bigl(
 \Cok N_P\otimes (dz_1)^{-1}
 \bigr)
 \lrarr
 \Cok(\rho_{0|X_P}),
 \quad
 v\cdot (dz_1)^{-1}
 \longmapsto
 \var^{-1}(v).
\]
We obtain
\[
\iota_{H_P\dagger}\Bigl(
 W_0\Cok N_P
 \Bigr)
 \simeq
 C_{0|X_P}.
\]
The natural morphism
$\Ker N_P\to \Cok N_P$
induces an isomorphism
$\Ker N_P/W_{-1}\Ker N_P\simeq W_0\Cok N_P$.
Therefore, we obtain an isomorphism
$K_{0|X_P}\simeq C_{0|X_P}$.
It is independent of the choice of
$(z_1,\ldots,z_n)$.
(See the proof of Lemma \ref{lem;22.4.20.1}, for example.)
Hence, we obtain the globally defined isomorphism
$K_0\simeq C_0$.

\vspace{.1in}
In the case where $H$ is smooth and reduced,
it is convenient to introduce some $\nbigd$-modules on $H$.
In the above local construction,
the filtration $V(M_{P})$ is independent of
the choice of $(z_1,\ldots,z_n)$.
Hence, we obtain a global filtration $V(M)$
by $V\nbigr_X$-modules.
We obtain the $V\nbigr_X$-module
$\psi_{-1}(M):=V_{-1}(M)/V_{<-1}(M)$.
It is obtained as the gluing of $\psi_{-1}(M_P)$.
By gluing the operators $-z_1\del_{z_1}$ on $\psi_{-1}(M_P)$,
we obtain a globally defined operator
$N:\psi_{-1}(M)\lrarr\psi_{-1}(M)$.
It is nilpotent, and we obtain the weight filtration $W$.
The kernel $\Ker N$
and the cokernel $\Cok N$
are naturally $\nbigd_H$-modules.
Let $\iota_H:H\lrarr X$ denote the inclusion.
For $P\in H$,
we have
$(\iota_{H\dagger}\Ker N)_{|X_P}
=\iota_{H_P\dagger}(\Ker N_P)$
and
$(\iota_{H\dagger}\Cok N)_{|X_P}
=\iota_{H_P\dagger}(\Cok N_P)$.
We can check that
the isomorphisms
$\iota_{H_P\dagger}(\Ker N_P)
\simeq
 \Ker(\rho_0)_{|X_P}$
and 
$\iota_{H_P\dagger}(\Cok N_P)
\simeq
 \Cok(\rho_0)_{|X_P}$
are independent of the coordinate system.
Hence, we obtain the globally defined
isomorphisms
$\iota_{H\dagger}\Ker N\simeq\Ker\rho_0$
and
$\iota_{H\dagger}\Cok N\simeq\Cok\rho_0$.
We set $K_{0,H}:=\Ker N/W_{-1}\Ker N$
and $C_{0,H}:=W_0\Cok N$.
We have the naturally induced isomorphisms
$\iota_{H\dagger}K_{0,H}\simeq K_0$
and
$\iota_{H\dagger}C_{0,H}\simeq C_0$.
The natural morphism $\Ker N\lrarr \Cok N$
induces an isomorphism $K_{0,H}\simeq C_{0,H}$,
which induces $K_0\simeq C_0$.

\subsubsection{The first Chern class
in the case where $M$ is the minimal extension}

Let $M$ be a holonomic $\nbigd_X$-module.
We assume that $M$ is isomorphic to the image of
$\rho_0:M(!H)\to M(\ast H)$.
We set
\[
 \Mtilde_1:=\Ker\Bigl(
 M(\ast H)
 \lrarr \Cok(\rho_0)/W_0\Cok(\rho_0)
 \Bigr),
\quad\quad
 \Mtilde_2:= M(!H)/W_{-1}\Ker(\rho_0).
\]
There exist the following exact sequences:
\[
 0\lrarr M\lrarr \Mtilde_1
 \lrarr
 C_0\lrarr 0,
\quad\quad
0\lrarr K_0\lrarr \Mtilde_2
 \lrarr M\lrarr 0.
\]
We obtain the following morphisms:
\begin{equation}
 \del_1:F_{\dagger}^j(K_0)
  \simeq F_{\dagger}^j(C_0)
  \lrarr F_{\dagger}^{j+1}(M),
\end{equation}
\begin{equation}
 \del_2:F_{\dagger}^{j}(M)
  \lrarr 
  F_{\dagger}^{j+1}(K_0).
\end{equation}

We shall prove the following proposition in
\S\ref{subsection;22.3.16.10}--\S\ref{subsection;22.4.4.21}.

\begin{prop}
\label{prop;22.3.2.30}
The morphism $\del_1\circ\del_2$
is equal to
\[
 2\pi\sqrt{-1}c_1(\nbigo_X(H)):
 F_{\dagger}^j(M)\lrarr F_{\dagger}^{j+2}(M).
\]
The morphism $\del_2\circ\del_1$
is equal to  
\[
 2\pi\sqrt{-1}c_1(\nbigo_X(H)):
 F_{\dagger}^j(K_0)\lrarr F_{\dagger}^{j+2}(K_0).
\]
\end{prop}

\subsubsection{Truncated de Rham complexes}
\label{subsection;22.3.16.10}

Let $Y$ be a complex manifold.
Let $X_0$ be a complex manifold with a smooth reduced hypersurface $H_0$.
We set $X=X_0\times Y$ and $H=H_0\times Y$.
Let $M$ be a holonomic $\nbigd_X$-module
which is not necessarily
isomorphic to the image of $M(!H)\lrarr M(\ast H)$.

Let $P\in H_0$ be any point.
Let $X_{0,P}$ be a neighbourhood of $P$ in $X_{0}$
with a holomorphic coordinate neighbourhood
$(z_1,\ldots,z_n)$
such that $H_{0,P}:=H_0\cap X_{0,P}=\{z_1=0\}$.

We set $X_P:=X_{0,P}\times Y$.
We set $M_P:=M_{|X_{P}}$.
Let $V_{\bullet}(M_P)$ denote the $V$-filtration along $z_1$.
Let $\Omega^{\bullet}_{X_P/Y,\neq 1}$ denote the subcomplex of
$\Omega^{\bullet}_{X_P/Y}$
generated by $dz_i$ $(i\neq 1)$.
There exists a subcomplex
\[
\nbigc^{\bullet}(M_P):=
\Bigl(
 V_{-1}(M_P)\otimes\Omega^{\bullet}_{X_P/Y,\neq 1}
\Bigr)
 \oplus
\Bigl(
 V_{0}(M_P)\otimes\,\bigl(dz_1\cdot \Omega^{\bullet}_{X_P/Y,\neq 1}\bigr)
 \Bigr)
 \subset
 M_P\otimes \Omega^{\bullet}_{X_P/Y}.
\]
\begin{lem}
The inclusion is a quasi-isomorphism. 
\end{lem}
\pf
Because
$\del_{z_1}:\psi_a(M_P)\lrarr \psi_{a+1}(M_P)$
are isomorphisms for $a>-1$,
we obtain the claim of the lemma.
\hfill\qed

\vspace{.1in}
The subcomplex $\nbigc^{\bullet}(M_P)$
is independent of the choice of
$(z_1,\ldots,z_n)$.
Hence, we obtain a global subcomplex
$\nbigc^{\bullet}(M)\subset M\otimes\Omega^{\bullet}_{X/Y}$.

\vspace{.1in}
If the support of $M$ is contained in $H$,
there exists a $\nbigd_H$-module $M_H$
with an isomorphism $\rho_M:M\simeq \iota_{H\dagger}M_H$
by Kashiwara's equivalence.
Note that $\rho_M$ induces an isomorphism
$V_0(M)\simeq
 \iota_{H\ast}\Bigl(
 M_H\otimes \bigl(\iota_H^{\ast}\Omega^{\dim X}_X\bigr)^{-1}
 \Omega_H^{\dim H}
 \Bigr)$
and the following isomorphism
\begin{equation}
\label{eq;22.4.21.20}
 \nbigc^{\bullet}(M)[1]\simeq
 \iota_{H\ast}
 \bigl(
 M_H\otimes\Omega^{\bullet}_H
 \bigr).
\end{equation}
Locally,
for a section $\tau$ of $\Omega^{\bullet}_{X_P/Y,\neq 1}$
and a section $v$ of $V_0(M_P)$,
the isomorphism is described as
\[
v\otimes (dz_1\wedge \tau)\longmapsto
\iota_{H\ast}\bigl(
\bigl(
\rho_M(v)\otimes dz_1
\bigr)\otimes \tau_{|H}
\bigr)
\] 

\vspace{.1in}
Let us consider the case where
$M$ is isomorphic to the image of $M(!H)\to M(\ast H)$.
Let $V'_{-1}(M)$ denote the inverse image of
$\Ker N$ by the projection $V_{-1}(M)\to \psi_{-1}(M)$.
Around any $P\in H$,
there exists the following subcomplex:
\[
 \nbigc^{\prime\bullet}(M_P):=
 \Bigl(
 V'_{-1}(M_P)\otimes\Omega^{\bullet}_{X_P/Y,\neq 1}
 \Bigr)
 \oplus
 \Bigl(
 V_{<0}(M_P)\otimes\,\bigl(
 dz_1\cdot \Omega^{\bullet}_{X_P/Y,\neq 1}
 \bigr)
 \Bigr)
 \subset
 \nbigc^{\bullet}(M_P).
\]
\begin{lem}
The inclusion is a quasi-isomorphism.
\end{lem}
\pf
Note that
$\can:V_{-1}(M)/V_{<-1}(M)\to V_{0}(M)/V_{<0}(M)$ is an epimorphism,
and that
$\var:V_{0}(M)/V_{<0}(M)\to V_{-1}(M)/V_{<-1}(M)$ is a monomorphism.
Then, we obtain the claim of the lemma.
\hfill\qed

\vspace{.1in}
Because 
$\nbigc^{\prime\bullet}(M_P)$ is independent
of the choice of a coordinate system $(z_1,\ldots,z_n)$,
we obtain the global subcomplex
$\nbigc^{\prime\bullet}(M)
\subset
\nbigc^{\bullet}(M)$.
We may naturally regard
$\nbigc^{\prime\bullet}(M)$
as a subcomplex of
$\nbigc^{\bullet}(M(\star H))$ $(\star=!,\ast)$
and
$\nbigc^{\bullet}(\Mtilde_i)$ $(i=1,2)$.

\subsubsection{The boundary map (1)}
\label{subsection;22.3.3.5}

Let $Y$, $X$, $H$ and $M$ be as in \S\ref{subsection;22.3.16.10}.
Suppose moreover that $M$ is isomorphic to the image of
$M(!H)\lrarr M(\ast H)$.
There exists the following exact sequence of complexes:
\begin{equation}
\label{eq;22.4.21.10}
 0\lrarr \nbigc^{\bullet}(M)
 \lrarr\nbigc^{\bullet}(\Mtilde_1)
 \lrarr\nbigc^{\bullet}(C_0)\lrarr 0.
\end{equation}
Let $p_X:X\lrarr X_0$ denote the projection.
Let $\Omega^{0,\bullet}_{X/Y}$ denote
the sheaf of $C^{\infty}$-sections of $p_X^{-1}\Omega^{0,\bullet}_{X_0}$
which are holomorphic in the $Y$-direction.
We obtain the following exact sequence of complexes:
\begin{equation}
\label{eq;22.4.4.30}
 0\lrarr
 \Tot\Bigl(
 \nbigc^{\bullet}(M)
 \otimes
 \Omega^{0,\bullet}_{X/Y}
 \Bigr)
 \lrarr
 \Tot\Bigl(
 \nbigc^{\bullet}(\Mtilde_1)
  \otimes
 \Omega^{0,\bullet}_{X/Y}
 \Bigr)
 \lrarr
 \Tot\Bigl(
 \nbigc^{\bullet}(C_0)
  \otimes
 \Omega^{0,\bullet}_{X/Y}
 \Bigr)
 \lrarr 0.
\end{equation}
We note that the isomorphism
$C_0\simeq \iota_{H\dagger}(C_{0,H})$
in \S\ref{subsection;22.4.21.11}
induces the following isomorphisms
as in (\ref{eq;22.4.21.20}):
\begin{equation}
 \rho_0:
  \nbigc^{\bullet}(C_0)[1]\simeq
 \iota_{H\ast}(C_{0,H}\otimes\Omega^{\bullet}_{H/Y}),
\end{equation}
\begin{equation}
 \rho_1:
 \Tot\Bigl(
 \nbigc^{\bullet}(C_0)\otimes_{\nbigo_X}\Omega^{0,\bullet}_{X/Y}
 \Bigr)[1]
\simeq
 \Tot\Bigl(
\iota_{H\ast}(C_{0,H}\otimes\Omega^{\bullet}_{H/Y})
 \otimes_{\nbigo_X}
 \Omega^{0,\bullet}_{X/Y}
 \Bigr).
\end{equation}

Let $h$ be a $C^{\infty}$-metric of $\nbigo_{X_0}(H_0)$.
Let $1:\nbigo_{X_0}\lrarr \nbigo_{X_0}(H_0)$ denote the canonical section.
We obtain the $C^{\infty}$-function $\kappa=|1|_{h}^2$.
The $(1,0)$-form $\del_{X_0}\log\kappa$ is a $C^{\infty}$-section of
$\Omega_{X_0}^{1,0}(\log H_0)$.
Around any point $P$ of $H_0$
with a coordinate system $(z_1,\ldots,z_n)$ as above,
$\del_{X_0}\log\kappa-dz_1/z_1$ is a $C^{\infty}$-section of
$\Omega^{1,0}_{X_{0,P}}$.

Let us describe the image of
$c\in F_{\dagger}^j(K_0)$
by the boundary map
$\del_1:
F_{\dagger}^j(K_0)\simeq
F_{\dagger}^j(C_0)\lrarr F_{\dagger}^{j+1}(M)$.
There exists a section $u$
of
$\Tot^{j+\dim H_0}\Bigl(
\bigl(
K_{0,H}\otimes\Omega^{\bullet}_{H/Y}
\bigr)\otimes\Omega^{0,\bullet}_{X/Y}
\Bigr)$
such that $du=0$,
which represents $c$.
There exists a section $u_1$
of
$\Tot^{j+\dim H_0}\Bigl(
\bigl(
\Ker N\otimes\Omega^{\bullet}_{H/Y}
\bigr)\otimes\Omega^{0,\bullet}_{X/Y}
\Bigr)$
which induces $u$.
There exists a section $u_2$ of
$\Tot^{j+\dim H_0}\Bigl(
\nbigc^{\prime\bullet}(M)
\otimes\Omega^{0,\bullet}_{X/Y}
\Bigr)$
whose restriction to $H$
induces $u_1$.
We obtain the section
$\del(\log\kappa) u_2$
of
$\Tot^{j+\dim X_0}\Bigl(
 \nbigc^{\bullet}(M(\ast H))
 \otimes\Omega^{0,\bullet}_{X/Y}
 \Bigr)$.

\begin{lem}
$\del(\log\kappa) u_2$
is a section of
$\Tot^{j+\dim X_0}\Bigl(
 \nbigc^{\bullet}(\Mtilde_1)
 \otimes\Omega^{0,\bullet}_{X/Y}
 \Bigr)$.
Moreover, the induced section $u_3$ of
$\Tot^{j+\dim X_0}\Bigl(
 \nbigc^{\bullet}(C_0)
 \otimes\Omega^{0,\bullet}_{X/Y}
 \Bigr)$
satisfies $du_3=0$,
and it represents $c$. 
\end{lem}
\pf
Let $P$ be any point of $H_0$.
We use the notation in \S\ref{subsection;22.3.16.10}.
We have the expression
$u_{2|X_P}=u_{2}'+dz_1\cdot u_2''$,
where
$u_2'$
is a section of
$\Tot^{j+\dim H_0}\bigl(
V'_{-1}(M_P)\otimes \Omega^{\bullet}_{X_P/Y,\neq 1}
\otimes\Omega^{0,\bullet}_{X_P/Y}
\bigr)$,
and 
$u_2''$
is a section of
$\Tot^{j+\dim H_0-1}\bigl(
V_{<0}(M_P)\otimes \Omega^{\bullet}_{X_P/Y,\neq 1}
\otimes\Omega^{0,\bullet}_{X_P/Y}
\bigr)$.
Moreover,
$u_2'$ induces $u_1$.
We set $\nu:=\del(\log\kappa)-dz_1/z_1$.
We have
$\bigl(
\del(\log\kappa)u_2\bigr)_{|X_P}
=(dz_1/z_1)\cdot u'_2
+\nu (u_2'+dz_1\cdot u_2'')$.
Note that
$\nu(u_2'+dz_1\cdot u_2'')$
is a section of
$\Tot^{j+\dim X_0}\Bigl(
\nbigc^{\bullet}(M_P)\otimes\Omega^{0,\bullet}_{X_P/Y}
\Bigr)$.
Because $(dz_1/z_1)\cdot u_2'$ is a section of
$\Tot^{j+\dim X_0}\bigl(
\nbigc^{\bullet}(\Mtilde_{1,P})\otimes\Omega^{0,\bullet}_{X_P/Y}
\bigr)$,
we obtain the first claim.
By a local computation,
we can check $\rho_1(u_3)=u$.
Hence, we obtain the second claim.
\hfill\qed

\vspace{.1in}

We obtain a section
$d\bigl((\del\log\kappa) u_2\bigr)$
of
$\Tot^{j+1+\dim X_0}\Bigl(
\nbigc^{\bullet}(\Mtilde_1)
\otimes\Omega^{0,\bullet}_{X/Y}
\Bigr)$.
By the construction, we obtain the following lemma.
\begin{lem}
$d\bigl((\del\log\kappa)u_2\bigr)$
is a section of 
$\Tot^{j+1+\dim X_0}\Bigl(
\nbigc^{\bullet}(M)
\otimes\Omega^{0,\bullet}_{X/Y}
\Bigr)$
and  represents $\del_1c$. 
\hfill\qed
\end{lem}

\subsubsection{The boundary map (2)}
\label{subsection;22.3.3.4}

There exists the following exact sequence of complexes:
\[
 0\lrarr \nbigc^{\bullet}(K_0)
 \lrarr \nbigc^{\bullet}(\Mtilde_2)
 \lrarr \nbigc^{\bullet}(M)\lrarr 0.
\]
We obtain the following exact sequence of complexes:
\begin{equation}
\label{eq;22.3.3.3}
 0\lrarr
 \Tot\Bigl(
 \nbigc^{\bullet}(K_0)\otimes_{\nbigo_X}\Omega^{0,\bullet}_{X/Y}
 \Bigr)
 \lrarr
 \Tot\Bigl(
 \nbigc^{\bullet}(\Mtilde_2)
 \otimes_{\nbigo_X}\Omega^{0,\bullet}_{X/Y}
 \Bigr)
 \lrarr
 \Tot\Bigl(
 \nbigc^{\bullet}(M)
 \otimes_{\nbigo_X}\Omega^{0,\bullet}_{X/Y}
 \Bigr)\lrarr 0.
\end{equation}
The isomorphism
$K_0\simeq \iota_{H\dagger}(K_{0,H})$
in \S\ref{subsection;22.4.21.11}
induces the following isomorphisms of complexes
as in (\ref{eq;22.4.21.20}):
\begin{equation}
 \rho_3:
 \nbigc^{\bullet}(K_0)[1]
 \simeq
   \iota_{H\ast}
 \bigl(
 K_{0,H}\otimes
 \Omega^{\bullet}_{H/Y}
 \bigr),
\end{equation}
\begin{equation}
 \rho_4:
 \Tot\Bigl(
 \nbigc^{\bullet}(K_0)
  \otimes
 \Omega^{0,\bullet}_{X/Y}
 \Bigr)[1]
\simeq
 \Tot
 \Bigl(
 \iota_{H\ast}
 \bigl(
 K_{0,H}\otimes
 \Omega^{\bullet}_{H/Y}
 \bigr)\otimes
 \Omega^{0,\bullet}_{X/Y}
 \Bigr).
\end{equation}

Let $x$ be a section of $F_{\dagger}^j(M)$.
Let us describe the image of $x$ by
the boundary map
$\del_2:F_{\dagger}^j(M)\to F_{\dagger}^{j+1}(K_0)$.
There exists a section $s$ of
$\Tot^{j+\dim X_0}(\nbigc^{\prime\bullet}(M)\otimes\Omega^{0,\bullet}_{X/Y})$
such that $ds=0$,
which represents $x$.
We obtain the induced section $[s]$
of
$\Tot^{j+1+\dim H_0}\Bigl(
\iota_{H\ast}
\bigl(
\Ker(N)\otimes\Omega^{\bullet}_{H/Y}
\bigr)
\otimes
\Omega^{0,\bullet}_{X/Y}
\Bigr)$
such that $d[s]=0$.
It induces a section
$\langle s\rangle$
of 
$\Tot^{j+1+\dim H_0}\Bigl(
\iota_{H\ast}
\bigl(
K_{0,H}\otimes\Omega^{\bullet}_{H/Y}
\bigr)
\otimes
\Omega^{0,\bullet}_{X/Y}
\Bigr)$
such that $d\langle s\rangle=0$.
\begin{lem}
$-\langle s\rangle$
represents $\del_2 x$. 
\end{lem}
\pf
When we regard $s$ as a section of
$\Tot^{j+\dim X_0}\Bigl(
\nbigc^{\bullet}(\Mtilde_2)
\otimes\Omega^{0,\bullet}_{X/Y}
\Bigr)$,
it is denoted by $\stilde_2$.
We obtain the section
$d\stilde_2$ of
$\Tot^{j+1+\dim X_0}\Bigl(
\nbigc^{\bullet}(\Mtilde_2)
\otimes\Omega^{0,\bullet}_{X/Y}
 \Bigr)$.
It gives the cocycle in
$\Tot^{j+1+\dim X_0}\Bigl(
\nbigc^{\bullet}(K_0)
\otimes\Omega^{0,\bullet}_{X/Y}
 \Bigr)$,
and represents
$\del_2 x$.
By a local computation,
we can check
$\rho_4(d\stilde_2)=-\langle s\rangle$.
\hfill\qed

\subsubsection{Proof of Proposition \ref{prop;22.3.2.30}}
\label{subsection;22.4.4.21}

It is enough to study the case
$F$ is the projection $X=X_0\times Y\lrarr Y$
and $H=H_0\times Y$.
By the considerations in \S\ref{subsection;22.3.3.4}
and \S\ref{subsection;22.3.3.5},
$\del_1\circ\del_2$ and $\del_2\circ\del_1$
are represented by the multiplication of
$-\delbar\del\log\kappa$.
Because $\frac{\sqrt{-1}}{2\pi}\delbar\del\log\kappa$
represents $c_1(\nbigo(H))$,
we obtain the claim of Proposition \ref{prop;22.3.2.30}.
\hfill\qed

\section{Appendix}
\label{section;22.4.2.1}

\subsection{Mixed twistor structures}
\label{subsection;22.3.23.1}

\subsubsection{Some fundamental properties of mixed twistor structures}

We recall the notion of mixed twistor structure
with some operations
introduced by Simpson \cite{s3}.
A twistor structure
is a holomorphic vector bundle of finite rank
on $\proj^1$.
A morphism of twistor structures is
an $\nbigo_{\proj^1}$-homomorphism.
Let $\TS$ denote the category of twistor structures.
A mixed twistor structure is a twistor structure
$\nbigv$ with an increasing filtration $W$ indexed by integers
such that
(i) $W_j(\nbigv)=0$ for $j<\!<0$,
(ii) $W_j(\nbigv)=\nbigv$ for $j>\!>0$,
(iii) $\Gr^W_j(\nbigv)\simeq\nbigo_{\proj^1}(j)^{\oplus r(j)}$.
The filtration $W$ is called the weight filtration of
the mixed twistor structure.
A morphism $f:(\nbigv_1,W)\to (\nbigv_2,W)$
is defined to be an $\nbigo_{\proj^1}$-homomorphism
$f:\nbigv_1\to\nbigv_2$
such that
$f(W_j(\nbigv_1))\subset W_j(\nbigv_2)$.
Let $\MTS$ denote the category of mixed twistor structures.
A mixed twistor structure $(\nbigv,W)$ is called
pure of weight $n$
if $\Gr^W_j(\nbigv)=0$ unless $j=n$.
In that case, we also say that
$\nbigv$ is pure of weight $n$.
Let $\PTS_n\subset\MTS$ denote the full subcategory
of pure twistor structure of weight $n$.

\begin{prop}[\cite{s3}]
\label{prop;22.3.24.1}
The category $\MTS$ is abelian.
Any morphism $f:(\nbigv_1,W)\to (\nbigv_2,W)$
is strict with respect to the weight filtration $W$,
i.e.,
$\Image f\cap W_j(\nbigv_2)
=f(W_j(\nbigv_1))$.
\hfill\qed
\end{prop}

Let $(\nbigv^{\bullet},W)$ be a complex
with the differential $d^i:(\nbigv^i,W)\to(\nbigv^{i+1},W)$
in the abelian category $\MTS$.
We obtain
$H^i(\nbigv^{\bullet},W)=\Ker d^i/\Image d^{i-1}\in\MTS$
as usual.
By the strictness with respect to the weight filtration,
we obtain the following lemma.
\begin{lem}
\label{lem;22.3.24.50}
There exist natural isomorphisms
\[
 \Gr^W_jH^i(\nbigv^{\bullet},W)\simeq
 H^i\bigl(\Gr^W_j(\nbigv^{\bullet})\bigr),
 \quad\quad
 W_jH^i(\nbigv^{\bullet},W)\simeq
 H^i\bigl(W_j(\nbigv^{\bullet})\bigr).
\]
\hfill\qed
\end{lem}

For $(\nbigv,W)$,
we set
$\Gr^W(\nbigv)=\bigoplus\Gr^W_j(\nbigv)$,
which is equipped with the induced weight filtration
$W_k(\Gr^W(\nbigv))=\bigoplus_{j\leq k}\Gr^W_j(\nbigv)$.
It induces a functor $\Gr^W:\MTS\to\MTS$.
Clearly, we have
$\Gr^W\circ\Gr^W=\Gr^W$.
By the strictness with respect to the weight filtration,
we obtain the following.
\begin{lem}
The functor $\Gr^W$ is faithful. 
\hfill\qed
\end{lem}

For $(\nbigv,W)\in\MTS$,
the fiber of $\nbigv$ over $\lambda\in\proj^1$
is denoted by $\nbigv_{|\lambda}$.
We obtain the functor $\Xi_{\lambda}$ from
$\MTS$ to the category of finite dimensional $\cnum$-vector spaces.
The following easy lemma is often used implicitly.

\begin{lem}
The functor $\Xi_{\lambda}$ is faithful.
\end{lem}
\pf
Let $f:(\nbigv_1,W)\to(\nbigv_2,W)$ be a morphism.
If $\Xi_{\lambda}(f)=0$,
we obtain $\Xi_{\lambda}\Gr^W_j(f)=0$,
which implies $\Gr^W_j(f)=0$.
Hence, we obtain $f=0$.
\hfill\qed

\begin{lem}
Let $(\nbigv,W)\in\TS$.
Let $(\nbigv_i,W)$ $(i=1,2)$ be subobjects of $(\nbigv,W)$.
If $\nbigv_{1|\lambda}=\nbigv_{2|\lambda}$  for some $\lambda\in\proj^1$,
we obtain
$(\nbigv_1,W)=(\nbigv_2,W)$. 
\end{lem}
\pf
Let us consider the induced morphism
$a:(\nbigv_2,W)\to (\nbigv,W)/(\nbigv_1,W)$.
Because $\Xi_{\lambda}(a)=0$,
we obtain that $a=0$.
\hfill\qed

\begin{lem}
\label{lem;22.4.25.1}
Let $f:(\nbigv_1^{\bullet},W)\to (\nbigv_2^{\bullet},W)$
be a morphism of complexes in $\MTS$.
If $\Xi_{\lambda_0}(f)$ is a quasi-isomorphism
for some $\lambda_0\in\proj^1$,
then $f$ is a quasi-isomorphism. 
Conversely,
if $f$ is a quasi-isomorphism,
then  $\Xi_{\lambda}(f)$ is a quasi-isomorphism
for any $\lambda\in\proj^1$.
\end{lem}
\pf
We have the induced morphism
$\Gr^W(f):
\Gr^W(\nbigv_1^{\bullet})
\to
\Gr^W(\nbigv_2^{\bullet})$.
Because the category of pure twistor structures of weight $n$
is equivalent to the category of $\cnum$-vector spaces,
we can easily check the claims for 
$\Gr^W(f)$.
We obtain the claim of Lemma \ref{lem;22.4.25.1}
from Lemma \ref{lem;22.3.24.50}.
\hfill\qed

\subsubsection{The anti-holomorphic involution
and Tate twists}

Let $\sigma:\proj^1\lrarr\proj^1$ be
the anti-holomorphic involution
defined by $\sigma(\lambda)=\overline{(-\lambda)^{-1}}$.
We have the natural isomorphism of sheaves of algebras
$\sigma^{-1}(\nbigo_{\proj^1})\simeq
\nbigo_{\proj^1}$,
induced by $\sigma^{-1}(f)\longmapsto\overline{\sigma^{-1}(f)}$.
For any $\nbigo_{\proj^1}$-module $\nbign$,
we obtain the $\nbigo_{\proj^1}$-module
\[
\sigma^{\ast}(\nbign):=
\nbigo_{\proj^1}\otimes_{\sigma^{-1}\nbigo}
\sigma^{-1}(\nbign).
\]
In particular,
for any mixed twistor structure $(\nbigv,W)$,
we obtain the mixed twistor structure
$\sigma^{\ast}(\nbigv,W)$.

We set
$\nbigo(p,q):=\nbigo(p\{0\}+q\{\infty\})$
and $\Tate(j):=\nbigo(-j,-j)$.
The latter is called the $j$-th Tate object.
There exist the following isomorphisms:
\[
 \iota_{p,q}:
 \sigma^{\ast}(\nbigo(p,q))
 \simeq
 \nbigo(p,q),
  \quad
  \iota_{p,q}(\sigma^{-1}(f))
  =(-1)^p\overline{\sigma^{\ast}(f)}
\]
In particular, we obtain
\[
\iota_{\Tate(n)}:
\sigma^{\ast}(\Tate(n))
\simeq
\Tate(n),
\quad
\iota_{\Tate(n)}(\sigma^{-1}f)
=(-1)^n\overline{\sigma^{\ast}(f)}.
\]

\begin{rem}
See {\rm\S\ref{subsection;22.3.23.10}}
for the comparison with
 $\nbigo(p,q)$  and $\Tate(n)$
in {\rm\cite{mochi2,mochi8}}.
\hfill\qed
\end{rem}

\subsubsection{Polarizations}

For a twistor structure $\nbigv$,
a morphism $S:\nbigv\otimes\sigma^{\ast}(\nbigv)\to\Tate(-n)$
is called a pairing of $\nbigv$ of weight $n$.
It is called $(-1)^n$-symmetric
if
\[
(-1)^n\iota_{\Tate(-n)}\circ\sigma^{\ast}S
=S\circ\exchange,
\]
where $\exchange$ denotes
the natural isomorphism
$\nbigv\otimes\sigma^{\ast}(\nbigv)
\simeq\sigma^{\ast}\nbigv\otimes\nbigv$.
It means that
for an open subset $U\subset \proj^1$,
and sections $u\in H^0(U,\nbigv)$ and $v\in H^0(\sigma(U),\nbigv)$,
we have
\[
 \overline{
 \sigma^{\ast}\bigl(
 S(u\otimes\sigma^{\ast}v)
 \bigr)}
=(-1)^n\iota_{\Tate(-n)}\circ
  \sigma^{\ast}S(\sigma^{\ast}(u)\otimes v)
=S(v\otimes\sigma^{\ast}u). 
\]

\begin{example}
The isomorphism
 $\iota_{p,q}:\sigma^{\ast}\nbigo(p,q)\simeq
 \nbigo(q,p)
 =\nbigo(p,q)^{-1}\otimes\Tate(-(p+q))$
induce the following morphism:
\[
 S_{p,q}:
 \nbigo(p,q)
 \otimes
 \sigma^{\ast}(\nbigo(p,q))
 \lrarr
 \Tate(-(p+q)),
 \quad
 \quad
 S_{p,q}(f\otimes\sigma^{\ast}(g))
=(-1)^pf\overline{\sigma^{\ast}(g)}. 
\]
It is a $(-1)^{p+q}$-symmetric pairing of weight $p+q$.
\hfill\qed
\end{example}

We recall the definition of a polarization of
a pure twistor structure of weight $n$.
\begin{df}[\cite{s3}]
Let $\nbigv$ be a pure twistor structure of weight $n$.
Let $S$ be a $(-1)^n$-symmetric pairing of $\nbigv$
of weight $n$.
It is called a polarization of $\nbigv$
if the following condition is satisfied.
\begin{itemize}
 \item The Hermitian pairing
       $H^0(S\otimes S_{0,-n})$
       of $H^0(\proj^1,\nbigv\otimes\nbigo(0,-n))$
       is positive definite.      
\hfill\qed
\end{itemize} 
\end{df}

Note that
a polarized pure twistor structure of weight $0$
is naturally equivalent to
a finite dimensional $\cnum$-vector space
with a Hermitian metric.

\begin{example}
For any $(p,q)\in\seisuu^2$,
$S_{p,q}$ is a polarization of 
$\nbigo(p,q)$.
Indeed,
we have
$S_{p,q}\otimes S_{0,-p-q}=S_{p,-p}$,
and
$\lambda^{-p}$ is a base of
$H^0(\proj^1,\nbigo(p,q)\otimes \nbigo(0,-p-q))$.
Hence, we obtain
\[
 H^0(S_{p,q}\otimes S_{0,-p-q})(\lambda^{-p},\sigma^{\ast}\lambda^{-p})
 =(-1)^p
 \cdot \lambda^{-p}
 \overline{\bigl(
 \overline{(-\lambda^{-1})}
 \bigr)^{-p}}
=1.
\]
\hfill\qed
\end{example}

\subsubsection{Twistor structures and $\nbigr$-triples}

Let us recall the basic notions and operations for $\nbigr$-triples
in \cite{sabbah2}.
We set
$\vecS=\bigl\{\lambda\in\cnum\,\big|\,|\lambda|=1\bigr\}$.
For a sheaf $\nbigf$ on $\cnum$,
we set $\nbigf_{|\vecS}:=\iota_{\vecS}^{-1}(\nbigf)$,
where $\iota_{\vecS}$ denotes
the inclusion $\vecS\lrarr \cnum_{\lambda}$.
We set $\nbigr:=\nbigo_{\cnum_{\lambda}}$.
Let $\nbigc_{\vecS}$ denote the sheaf of
continuous functions on $\vecS$.
It is naturally
an $\nbigr_{|\vecS}\otimes_{\cnum}
 \sigma^{\ast}(\nbigr_{|\vecS})$-module.
 
Let $\nbigm_i$ $(i=1,2)$ be $\nbigr$-modules.
A sesqui-linear paring of
$\nbigm_1$ and $\nbigm_2$
is a morphism
\[
 C:
  \nbigm_{1|\vecS}
  \otimes_{\cnum}
  \sigma^{-1}(\nbigm_{2|\vecS})
  \lrarr
  \nbigc_{\vecS}
\]
of
$\nbigr_{|\vecS}\otimes_{\cnum}
\sigma^{\ast}(\nbigr_{|\vecS})$-modules.
A morphism
of $\nbigr$-triples
$(\nbigm_1,\nbigm_2,C)\to (\nbigm_1',\nbigm_2',C')$
is a pair of
$\nbigr$-homomorphisms
$\varphi_1:\nbigm'_1\to\nbigm_1$
and
$\varphi_2:\nbigm_2\to\nbigm'_2$
such that
$C(\varphi_1(u'_1)\otimes u_2)
=C'(u_1'\otimes\varphi_2(u_2))$.
Let $\rtriplecat$ denote the category of
$\nbigr$-triples.

Let $\nbigv$ be a twistor structure.
We set
$\nbigv_1:=\nbigv_{|\cnum_{\lambda}}$
and
$\nbigv_2:=\nbigv_{|\cnum_{\mu}}$.
We obtain $\nbigr$-modules
$\nbigv_{1|\cnum_{\lambda}}^{\lor}$
and
$\sigma^{\ast}(\nbigv_{2|\cnum_{\mu}})$.
The isomorphism
$\bigl((\nbigv_1^{\lor})^{\lor}
\bigr)_{|\cnum_{\lambda}\cap\cnum_{\mu}}
\simeq
\sigma^{\ast}(\sigma^{\ast}\nbigv_2)_{|\cnum_{\lambda}\cap\cnum_{\mu}}$
induces a sesqui-linear pairing 
\[
 C_{\nbigv}:
 \nbigv_{1|\vecS}^{\lor}
 \otimes
 \sigma^{-1}\bigl(
 \sigma^{\ast}(\nbigv_{2|\vecS})
 \bigr)
\lrarr \nbigc_{\vecS}.
\]
Thus, we obtain
$\Theta(\nbigv)=
(\nbigv^{\lor}_1,
 \sigma^{\ast}(\nbigv_{2}),
 C_{\nbigv})$.
In this way,
we obtain the fully faithful functor 
$\Theta:\TS\to\rtriplecat$.
The essential image is also called
the category of twistor structures.
Similarly, the essential image of $\PTS(n)$ via $\Theta$
is also called the category of
pure twistor structure of weight $n$.
Let $\rtriplecat^{\filt}$ denote
the category of objects in $\rtriplecat$
equipped with an increasing filtration.
We have the naturally defined functor
$\Theta:\MTS\to\rtriplecat^{\filt}$.
The essential image is also called
the category of mixed twistor structures.

We set
$\nbigu(p,q)=\bigl(
\lambda^{p}\nbigo_{\cnum},
\lambda^{q}\nbigo_{\cnum},
C_0
\bigr)$,
where $C_0$ is the sesqui-linear pairing
induced by the natural multiplication.
We set
$\newTate(n)=\nbigu(-n,n)$.
There exist the natural isomorphisms
$\Theta(\nbigo(-p,q))\simeq\nbigu(p,q)$
and
$\Theta(\Tate(n))\simeq\newTate(n)$.

For $\nbigt=(\nbigm_1,\nbigm_2,C)\in\rtriplecat$,
we set
$\nbigt^{\ast}=(\nbigm_2,\nbigm_1,C^{\ast})$
called the Hermitian adjoint of $\nbigt$,
where
$C^{\ast}(u,\sigma^{-1}(v))
:=\overline{\sigma^{-1}(C(v,\sigma^{-1}(u)))}$.
For a morphism $f=(f_1,f_2):\nbigt_1\to\nbigt_2$ in $\rtriplecat$,
we obtain
$f^{\ast}=(f_2,f_1):\nbigt_2^{\ast}\to\nbigt_1^{\ast}$.
There exists the isomorphism
$\iota_{\newTate(n)}:\newTate(n)^{\ast}\simeq \newTate(-n)$
given by $((-1)^n,(-1)^n)$.

A sesqui-linear duality of $\nbigt$ of weight $n$
is a morphism
$S:\nbigt\to\nbigt^{\ast}\otimes\newTate(-n)$.
By using $\iota_{\newTate(-n)}$,
we obtain
$S^{\ast}:\nbigt\to \nbigt^{\ast}\otimes\newTate(-n)$.
The sesqui-linear duality $S$ is called Hermitian
if $S^{\ast}=(-1)^nS$.

Note that for any $\nbigv\in\TS$
we have
$\Theta(\nbigv)^{\ast}=
\Theta(\sigma^{\ast}\nbigv^{\lor})$.
A $(-1)^n$-symmetric pairing
$S:\nbigv\otimes\sigma^{\ast}\nbigv\to \Tate(-n)$
induces
a Hermitian sesqui-linear duality
$\Theta(S):
\Theta(\nbigv)\to
\Theta(\nbigv)^{\ast}\otimes\newTate(-n)$.
A Hermitian sesqui-linear duality $S_{\nbigt}$ of
a pure twistor structure $\nbigt$ of weight $n$
is called a polarization of $\nbigt$
if there exist a pure twistor structure $\nbigv$,
of weight $n$,
a polarization $S_{\nbigv}$ of $\nbigv$
and an isomorphism
$\nbigt\simeq\Theta(\nbigv)$
under which $S_{\nbigt}$ is equal to $\Theta(S_{\nbigv})$.

\begin{example}
\label{example;22.4.19.1}
There exists the isomorphism
$\iota_{\nbigu(p,q)}:\nbigu(p,q)\simeq
\nbigu(p,q)^{\ast}\otimes\newTate(-(p+q))$
 given by $((-1)^p,(-1)^p)$.
It is a polarization
because it is induced by
the $(-1)^{-p+q}$-symmetric pairing
 $S_{-p,q}$ of $\nbigo(-p,q)$.
\hfill\qed
\end{example}

\subsubsection{Polarized graded Lefschetz twistor structures}
\label{subsection;22.4.28.3}

We recall the notion of polarized graded Lefschetz twistor structure
in \cite{sabbah2},
to which we refer for more detailed properties.
Let $\epsilon=\pm 1$.
Let $w\in\seisuu$.
Let $\nbigt=\bigoplus_{i\in\seisuu}\nbigt_i$
be a graded twistor structure
such that each $\nbigt_i$ is pure of weight $w-\epsilon i$.
Assume that $\{i\in\seisuu\,|\,\nbigt_i\neq 0\}$
is a finite set.
Let 
$\nbigl:\nbigt_i\to
\nbigt_{i-2}\otimes\newTate(\epsilon)$ $(i\in\seisuu)$
be morphisms,
which induce
$\nbigl:\nbigt\to\nbigt\otimes\newTate(\epsilon)$.
Such a pair $(\nbigt,\nbigl)$ is called
a graded Lefschetz twistor structure
of weight $w$ and type $\epsilon$
if the induced morphisms
$\nbigt_{j}\to \nbigt_{-j}\otimes\newTate(\epsilon j)$
$(j\in\seisuu_{\geq 0})$
are isomorphisms.

Let  $(\nbigt,\nbigl)$ be
a graded Lefschetz twistor structure of weight $w$ and type $\epsilon$.
For $j\geq 0$,
let $P\nbigt_j$ denote the kernel of
$\nbigl^{j+1}:\nbigt_j\to \nbigt_{-j-2}\otimes\newTate(\epsilon(j+1))$,
and called the primitive part.
Let $\nbigs:\nbigt_j\to \nbigt^{\ast}_{-j}\otimes\newTate(-w)$
be morphisms.
They induce a graded morphism
$\nbigs:\nbigt\to\nbigt^{\ast}\otimes\newTate(-w)$.
It is called a sesqui-linear duality of weight $w$
of $(\nbigt,\nbigl)$
if $\nbigl^{\ast}\circ \nbigs+\nbigs\circ\nbigl=0$.
A Hermitian sesqui-linear duality of weight $w$
of $(\nbigt,\nbigl)$ is called
a polarization of $(\nbigt,\nbigl)$
if $(\nbigl^{\ast})^j\circ \nbigs=(-1)^j\nbigs\circ\nbigl^j$
induces a polarization on $P\nbigt_j$.       
for any $j\in\seisuu_{\geq 0}$.
\begin{rem}
We choose the signature of the nilpotent morphism $\nbigl$
in the way opposite to that in {\rm\cite{sabbah2}}.
Hence, the above condition is equivalent to
the condition in {\rm\cite{sabbah2}}.
See Remark {\rm\ref{rem;22.3.24.3}} below.
\hfill\qed 
\end{rem}

We recall a useful proposition \cite[Proposition 2.1.19]{sabbah2},
which goes back to \cite[Lemma 5.2.15]{saito1}.
Let $(\nbigt,\nbigl,\nbigs)$ be a polarized
graded Lefschetz twistor structure of weight $w$ and type $\epsilon$.
Let $(\nbigt',\nbigl',\nbigs')$ be a polarized
graded Lefschetz twistor structure of weight $w-\epsilon$
and type $\epsilon$.
Let $c:\nbigt_{j+1}\to\nbigt_j'$
and $v:\nbigt'_{j+1}\to\nbigt_j\otimes\newTate(\epsilon)$
be morphisms.
They induce graded morphisms
$c:\nbigt\to\nbigt'$
and
$v:\nbigt'\to\nbigt\otimes\newTate(\epsilon)$
of degree $-1$.
We assume $v\circ c=-\nbigl$ and $c\circ v=-\nbigl'$.
We also assume the commutativity of
the following diagram:
\[
 \begin{CD}
  \nbigt @>>>
  \nbigt^{\ast}\otimes\newTate(-w)
  \\
  @V{c}VV @V{\iota_{\newTate(-\epsilon)}\circ v^{\ast}}VV \\
  \nbigt'
  @>>>
  (\nbigt')^{\ast}\otimes\newTate(-w+\epsilon).
 \end{CD}
\]

\begin{prop}[\mbox{\cite[Proposition 2.1.19]{sabbah2}},
\mbox{\cite[Lemma 5.2.15]{saito1}}]
There exists a decomposition
$\nbigt'=\Image c\oplus\Ker v$
as a graded Lefschetz twistor structure of weight $w$.
\hfill\qed
\end{prop}

Let $\epsilon_1,\epsilon_2\in\{\pm 1\}$.
Let $\nbigt=\bigoplus \nbigt_{i,j}$
be a bi-graded twistor structure
such that each $\nbigt_{i,j}$ is pure twistor structure
of weight $w-\epsilon_1i-\epsilon_2j$.
We assume that
$\{(i,j)\,|\,\nbigt_{i,j}\neq 0\}$ is finite.
Let $\nbigl_1:\nbigt_{i,j}\to \nbigt_{i-2,j}\otimes\newTate(\epsilon_1)$
and
$\nbigl_2:\nbigt_{i,j}\to \nbigt_{i,j-2}\otimes\newTate(\epsilon_2)$
be morphisms
such that $\nbigl_1\circ\nbigl_2=\nbigl_2\circ\nbigl_1$.
Such a tuple $(\nbigt,\nbigl_1,\nbigl_2)$
is called a bi-graded Lefschetz twistor structure
of weight $w$ and type $(\epsilon_1,\epsilon_2)$
if the induced morphisms
$\nbigl_1^i:\nbigt_{i,j}\lrarr
 \nbigt_{-i,j}\otimes\newTate(\epsilon_1 i)$
$(i\geq 0)$
and
$\nbigl_2^j:\nbigt_{i,j}\lrarr
 \nbigt_{i,-j}\otimes\newTate(\epsilon_2 j)$
$(j\geq 0)$
are isomorphisms.
For $i,j\geq 0$,
let $P\nbigt_{i,j}$ denote the intersection of
$\Ker(\nbigl_1^{i+1}:
\nbigt_{i,j}\to
\nbigt_{-i-2,j}\otimes\newTate(\epsilon_1(i+1)))$
and
$\Ker(\nbigl_2^{j+1}:
\nbigt_{i,j}\to
\nbigt_{i,-j-2}\otimes\newTate(\epsilon_2(j+1)))$.
Let $\nbigs:\nbigt_{i,j}\to\nbigt_{-i,-j}^{\ast}\otimes\newTate(-w)$
be morphisms.
They induce a graded morphism
$\nbigt\to\nbigt^{\ast}\otimes\newTate(-w)$.
It is called a sesqui-linear duality of 
$(\nbigt,\nbigl_1,\nbigl_2)$
if
$\nbigs\circ\nbigl_p+\nbigl_p^{\ast}\circ\nbigs=0$ $(p=1,2)$.
A Hermitian sesqui-linear duality $\nbigs$ of weight $w$
of $(\nbigt,\nbigl_1,\nbigl_2)$ is called a polarization
if $(\nbigl_1^{\ast})^i\circ(\nbigl_2^{\ast})^j\circ\nbigs$
induce a polarization on
$P\nbigt_{i,j}$ $(i,j\geq 0)$.

\subsubsection{Appendix: Comparison of Tate objects}
\label{subsection;22.3.23.10}

Let $\lambda$ be the standard coordinate of
$\cnum=\proj^1\setminus\{\infty\}$,
and we set $\mu=\lambda^{-1}$
on $\proj^1\setminus\{0\}$.
We set $\cnum_{\lambda}:=\proj^1\setminus\{\infty\}$
and $\cnum_{\mu}=\proj^1\setminus\{0\}$.
Let $\nbigo'(p,q)$ denote the twistor structure
obtained as the gluing of
$\nbigo_{\cnum_{\lambda}} f^{(p,q)}_0$
and $\nbigo_{\cnum_{\mu}}f^{(p,q)}_{\infty}$
by the relation
$(\sqrt{-1}\lambda)^{p+q}f^{(p,q)}_0
=f^{(p,q)}_{\infty}$.
We set
$\Tate'(n):=\nbigo'(-n,-n)$.
(In {\rm\cite{mochi2,mochi8}},
$\nbigo'(p,q)$ and $\Tate'(n)$
are denoted by $\nbigo(p,q)$
and $\Tate(n)$, respectively.)

Let
 $\iota_{p,q}':\sigma^{\ast}\nbigo'(p,q)
 \simeq
 \nbigo'(q,p)$
denote the isomorphism
induced by
\[
 \iota_{p,q}'(\sigma^{\ast}(f_0^{p,q}))
 =(\sqrt{-1})^{p+q}f_{\infty}^{q,p},
 \quad
  \iota_{p,q}'(\sigma^{\ast}(f_{\infty}^{p,q}))
=(-\sqrt{-1})^{p+q}f_{0}^{q,p}.
\]
We obtain
$\iota'_{\Tate(n)}=\iota'_{-n,-n}:
\sigma^{\ast}\Tate'(n)\simeq\Tate(n)$.

We choose $a\in\cnum$ such that $a^2=-\sqrt{-1}$.
We obtain the isomorphism
$\Psi_{p,q}:\nbigo'(p,q)\simeq\nbigo(p,q)$
by 
\[
 \Psi_{p,q}(f_0^{p,q})=(\sqrt{-1}\lambda)^{-p}a^{p+q},
 \quad
 \Psi_{p,q}(f_{\infty}^{p,q})
 =(-\sqrt{-1}\mu)^{-q}a^{p+q}.
\] 
We obtain $\Psi_{\Tate(n)}:\Tate'(n)\simeq \Tate(n)$.
Then, the following diagrams are commutative:
\[
 \begin{CD}
  \sigma^{\ast}\nbigo'(p,q)
  @>{\iota'_{p,q}}>>
  \nbigo'(q,p)\\
  @V{\sigma^{\ast}\Psi_{p,q}}VV @V{\Psi_{q,p}}VV \\
  \sigma^{\ast}\nbigo(p,q)
  @>{\iota_{p,q}}>>
  \nbigo(q,p),
 \end{CD}
\quad\quad\quad\quad\quad\quad
  \begin{CD}
  \sigma^{\ast}\Tate'(n)
  @>{\iota'_{\Tate(n)}}>>
  \Tate'(n)\\
   @V{\sigma^{\ast}\Psi_{\Tate(n)}}VV
   @V{\Psi_{\Tate(n)}}VV \\
  \sigma^{\ast}\Tate(n)
  @>{\iota_{\Tate(n)}}>>
  \Tate(n).
 \end{CD}
\]
Hence, we may replace
$\nbigo(p,q)$ and $\Tate(n)$
in \cite{mochi2,mochi8}
with
$\nbigo(p,q)$ and $\Tate(n)$
in this paper.

\subsection{Twistor structures on a complex manifold}
\label{subsection;22.3.25.10}

\subsubsection{Twistor structures}

Let $X$ be a complex manifold.
We recall the notion of variation of twistor structure on $X$,
introduced by Simpson in \cite{s3},
as a generalization of variation of Hodge structure.
To simplify the description,
we shall often say
``twistor structure on $X$''
instead of
``variation of twistor structure on $X$''.

We set $\nbigx:=\cnum_{\lambda}\times X$.
Let $p_{\lambda}:\nbigx\lrarr X$ denote the projection.
We set
$\Omegatilde^{p,q}_{\nbigx/\cnum_{\lambda}}:=
\lambda^{-p} p_{\lambda}^{-1}\Omega^{p,q}_X$.
We set
$\Omegatilde^{\bullet,\bullet}_{\nbigx/\cnum_{\lambda}}:=
\bigoplus_{p,q}\Omegatilde^{p,q}_{\nbigx/\cnum_{\lambda}}$.
Let $X^{\dagger}$ denote the conjugate of $X$.
Let $p_{\mu}:\nbigx^{\dagger}\lrarr X^{\dagger}$
denote the projection.
We set
$\nbigx^{\dagger}:=\cnum_{\mu}\times X^{\dagger}$.
We set
$\Omegatilde^{p,q}_{\nbigx^{\dagger}/\cnum_{\mu}}:=
\mu^{-p} p_{\mu}^{-1}\Omega^{p,q}_{X^{\dagger}}$.
We obtain
$\Omegatilde^{\bullet,\bullet}_{\nbigx^{\dagger}/\cnum_{\mu}}
=\bigoplus_{p,q}\Omegatilde^{p,q}_{\nbigx^{\dagger}/\cnum_{\mu}}$.

We regard  $\proj^1\times X$
as the $C^{\infty}$-manifold
obtained by gluing of
$\nbigx$ and $\nbigx^{\dagger}$
by $\lambda=\mu^{-1}$.
We regard
$\cnum_{\lambda}^{\ast}=\proj^1\setminus\{0,\infty\}
=\cnum_{\mu^{\ast}}$.
Because
$\Omega_{X^{\dagger}}^{p,q}=\Omega^{q,p}_X$,
we have the natural identifications
$\bigl(
\Omegatilde^{p,q}_{\nbigx/\cnum_{\lambda}}
\bigr)_{|\cnum_{\lambda}^{\ast}\times X}
=\bigl(
\Omegatilde^{q,p}_{\nbigx^{\dagger}/\cnum_{\mu}}
\bigr)_{|\cnum_{\lambda}^{\ast}\times X^{\dagger}}$.
Let $\xi\Omega^{p,q}_X$
denote the bundle obtained from
$\Omegatilde^{p,q}_{\nbigx/\cnum_{\lambda}}$
and
$\Omegatilde^{q,p}_{\nbigx^{\dagger}/\cnum_{\lambda}}$.
We obtain
$\xi\Omega^{\bullet,\bullet}_{X}:=
\bigoplus_{p,q} \xi\Omega^{p,q}_X$.
The exterior derivative of $X$ induces a differential
$d_X:\Tot\xi\Omega^{\bullet,\bullet}_X\lrarr
\Tot\xi\Omega^{\bullet,\bullet}_X[1]$.
We set
$\Tot^k\xi\Omega^{\bullet,\bullet}_X
=\bigoplus_{p+q=k}
\xi\Omega^{p,q}_X$.

Let $\pi_1:\proj^1\times X\lrarr \proj^1$
denote the projection.
We obtain the bundle
$\pi_1^{-1}\Omega^{0,1}_{\proj^1}$.
The differential operator
$\delbar_{\proj^1}:\Omega^{0,0}_{\proj^1}
\lrarr\Omega^{0,1}_{\proj^1}$
induces
$\delbar_{\proj^1}:
\pi_1^{-1}\Omega^{0,0}_{\proj^1}
\lrarr
\pi_1^{-1}\Omega^{0,1}_{\proj^1}$.

Let $V^{\sankaku}$ be a $C^{\infty}$-vector bundle
on $\proj^1\times X$.
A $\proj^1$-holomorphic structure of $V^{\sankaku}$
is a differential operator
$\delbar_{\proj^1,V^{\sankaku}}:
V^{\sankaku}
\lrarr
V^{\sankaku}\otimes\pi_1^{-1}\Omega^{0,1}_{\proj^1}$
such that
$\delbar_{\proj^1,V^{\sankaku}}(fs)
=\delbar_{\proj^1}(f)\cdot s
+f\delbar_{\proj^1,V^{\sankaku}}(s)$
for any $C^{\infty}$-function $f$
and any $C^{\infty}$-section $s$
of $V^{\sankaku}$.

When $V^{\sankaku}$ is equipped with
a $\proj^1$-holomorphic structure
$\delbar_{\proj^1,V^{\sankaku}}$,
a $\widetilde{TT}$-structure
of $(V^{\sankaku},\delbar_{\proj^1,V^{\sankaku}})$
is a differential operator
$\DD^{\sankaku}:
V^{\sankaku}
\lrarr
V^{\sankaku}
\otimes
\Tot^1\xi\Omega^{\bullet,\bullet}_X$
such that
(i) $\DD^{\sankaku}(fs)
=df\cdot s+f\DD^{\sankaku}(s)$
for any $C^{\infty}$-function $f$
and a $C^{\infty}$-section $s$ of $V^{\sankaku}$,
(ii) $\DD^{\sankaku}\circ\DD^{\sankaku}=0$,
where we consider the naturally induced
differential operator
$\DD^{\sankaku}:
V^{\sankaku}
\otimes\Tot\xi\Omega^{\bullet,\bullet}_X
\lrarr
V^{\sankaku}
\otimes\Tot\xi\Omega^{\bullet,\bullet}_X[1]$,
(iii) $[\DD^{\sankaku},\delbar_{\proj^1,V^{\sankaku}}]=0$.
Such a tuple
$(V^{\sankaku},\delbar_{\proj^1,V^{\sankaku}},
\DD^{\sankaku})$
is called a twistor structure on $X$.
A morphism
$\varphi:(V^{\sankaku}_1,\delbar_{\proj^1,V^{\sankaku}_1},\DD^{\sankaku}_1)
\to(V^{\sankaku}_2,\delbar_{\proj^1,V^{\sankaku}_2},\DD^{\sankaku}_2)$
is a morphism of the sheaves of $C^{\infty}$-sections
$\varphi:V_1\to V_2$
such that 
$\varphi\circ\delbar_{\proj^1,V^{\sankaku}_1}
=\delbar_{\proj^1,V^{\sankaku}_2}\circ\varphi$
and
$\varphi\circ\DD^{\sankaku}_1
=\DD^{\sankaku}_2\circ\varphi$.
Let $\TS(X)$ denote the category of twistor structures on $X$.

Let $f:Y\lrarr X$ be any morphism of complex manifolds.
Then, we can naturally obtain
the twistor structure
$f^{\ast}\nbigv\in\TS(Y)$ as the pull back
of $\nbigv\in\TS(X)$.
We set
$\nbigo(p,q)_X:=a_X^{\ast}\nbigo(p,q)$
and 
$\Tate(n)_X:=a_X^{\ast}\Tate(n)$,
where $a_X$ denote the natural map from $X$ to
the one point set $\pt$,
and we regard $\nbigo(p,q)$ and $\Tate(n)$
as objects in $\TS(\pt)$.

\subsubsection{The anti-holomorphic involution}

We obtain the involution
$\sigma_X:=\sigma\times\id_X$ on $\proj^1\times X$
from $\sigma$.
There exists the naturally defined isomorphism
\[
 \varphi_0:\sigma_X^{-1}
 \xi\Omega^{p,q}_X
 \simeq
 \xi\Omega^{q,p}_X,
 \quad\quad
 \sigma_X^{-1}(\tau)
 \longmapsto
 \overline{
 \sigma_X^{-1}(\tau)
 }.
\]
It induces the isomorphism of complexes
$\varphi_0:\sigma_X^{-1}(\Tot\xi\Omega^{\bullet,\bullet}_X,d_X)
\simeq
(\Tot\xi\Omega^{\bullet,\bullet}_X,d_X)$.
There also exists the following isomorphism of bundles
\[
 \sigma_X^{-1}
 \pi_1\Omega^{0,q}_{\proj^1}
 \simeq
 \pi_1\Omega^{0,q}_{\proj^1},
 \quad
 \sigma_X^{-1}\tau
 \longmapsto
  \overline{
  \sigma_X^{-1}(\tau)
  }.
\]
It induces the isomorphism of complexes
$\varphi_0:\Tot\sigma_X^{-1}\xi\Omega^{\bullet,\bullet}_X\simeq
\Tot\xi\Omega^{\bullet,\bullet}$.

Let $\nbigv=(V^{\sankaku},\delbar_{\proj^1,V^{\sankaku}},\DD^{\sankaku})$
be a twistor structure on $X$.
We obtain the $C^{\infty}$-vector bundle
$\sigma_X^{-1}(V^{\sankaku})$
as the pull back.
Let 
$\overline{\sigma_X^{-1}(V^{\sankaku})}$
denote the conjugate of 
$\sigma_X^{-1}(V^{\sankaku})$,
i.e.,
we have
$\overline{\sigma_X^{-1}(V^{\sankaku})}
=\sigma_X^{-1}(V^{\sankaku})$
as an $\real$-vector bundle,
and the multiplication of $a\in\cnum$
on $\overline{\sigma_X^{-1}(V^{\sankaku})}$
is given by the multiplication
of $\abar$ on $\sigma_X^{-1}(V^{\sankaku})$.
We obtain
\[
\delbar_{\proj^1,\overline{\sigma_X^{-1}V^{\sankaku}}}:=
\varphi_0\circ
\sigma_X^{-1}(\delbar_{\proj^1,V^{\sankaku}}):
\overline{\sigma_X^{-1}(V^{\sankaku})}
\lrarr
\overline{\sigma_X^{-1}(V^{\sankaku})}
 \otimes
 \pi_1^{-1}(\Omega^{0,1}_{\proj^1}),
\]
\[
 \DD^{\sankaku}_{\overline{\sigma_X^{-1}(V^{\sankaku})}}:=
\varphi_0\circ
\sigma_X^{-1}(\DD^{\sankaku}):
\overline{\sigma_X^{-1}(V^{\sankaku})}
\lrarr
\overline{\sigma_X^{-1}(V^{\sankaku})}
 \otimes
\Tot^1\xi\Omega^{\bullet,\bullet}_X.
\]
In this way,
we obtain a twistor structure
$\sigma_X^{\ast}(\nbigv)
=(\overline{\sigma_X^{-1}V^{\sankaku}},
\delbar_{\proj^1,\overline{\sigma_X^{-1}V^{\sankaku}}},
\DD^{\sankaku}_{\overline{\sigma_X^{-1}V^{\sankaku}}})$
on $X$.

The isomorphisms
$\iota_{p,q}:\sigma^{\ast}\nbigo(p,q)\simeq \nbigo(q,p)$
and
$\iota_{\Tate(n)}:
\sigma^{\ast}\Tate(n)\simeq \Tate(n)$
induce the isomorphisms
$\sigma_X^{\ast}\nbigo(p,q)_X\simeq\nbigo_X(q,p)$
and
$\sigma_X^{\ast}\Tate(n)_X\simeq\Tate(n)_X$,
which are also denoted by
$\iota_{p,q}$ and $\iota_{\Tate(n)}$.

\begin{rem}
Let $\pi_2:\proj^1\times X\to X$
denote the projection.
We may naturally regard
$\xi\Omega^{p,q}_X
=\pi_1^{-1}\nbigo(p,q)\otimes
\pi_2^{-1}\Omega_X^{p,q}$.
In {\rm\cite{mochi2}},
we use the isomorphism
$\overline{\sigma_X^{-1}\xi\Omega^{p,q}_X}
\simeq
\xi\Omega^{q,p}_X$
induced by 
the isomorphisms
$\overline{\sigma^{\ast}\nbigo(p,q)}
\simeq
\nbigo(q,p)$ and
the isomorphisms
 $\overline{\sigma_X^{-1}\pi_2^{-1}\Omega^{p,q}_X}
 \simeq
\pi_2^{-1}\Omega^{q,p}_X$
induced by
 $\sigma_X^{-1}\pi_2^{-1}(dz_i)\longmapsto
 -\pi_2^{-1}(d\zbar_i)$
and
 $\sigma_X^{-1}\pi_2^{-1}(d\zbar_i)
 \longmapsto -\pi_2^{-1}(dz_i)$.
It is equal to the isomorphism considered in this subsection.
\hfill\qed
\end{rem}

\subsubsection{Pure twistor structure and mixed twistor structure}

Let
$\nbigv=(V^{\sankaku},\delbar_{\proj^1,V^{\sankaku}},\DD^{\sankaku})$
be a twistor structure on $X$.
For any point $P\in X$,
let $\iota_P:\{P\}\lrarr X$ denote the inclusion.
We obtain the twistor structure
$\iota_P^{\ast}(\nbigv)\in\TS$.
If $\iota_P^{\ast}\nbigv$ is pure of weight $n$
for any $P\in X$,
then $\nbigv$ is called pure of weight $n$.

Let $W_{\bullet}(\nbigv)=\{W_m(\nbigv)\,|\,m\in\seisuu\}$
be an increasing filtration of $\nbigv$ by subbundles
indexed by integers which is preserved by
$\delbar_{\proj^1,V^{\sankaku}}$ and $\DD^{\sankaku}$.
A filtered twistor structure
$(\nbigv,W)$ is called a mixed twistor structure on $X$
if $\iota_P^{\ast}(\nbigv,W)$ are mixed twistor structures
for any $P\in X$.
A morphism of mixed twistor structures 
$\varphi:(\nbigv_1,W)\to(\nbigv_2,W)$ on $X$
is a morphism 
$\varphi:\nbigv_1\to\nbigv_2$ in $\TS(X)$
such that $\varphi(W_m(\nbigv_1))\subset W_m(\nbigv_2)$
for any $m$.
The following lemma is well known,
which is a direct consequence of Proposition \ref{prop;22.3.24.1}.
\begin{lem}  
\label{lem;22.3.25.21}
Let $f:(\nbigv_1,W)\to(\nbigv_2,W)$
be any morphism of mixed twistor structures on $X$.
\begin{itemize}
\item  $f$ is strict with respect to $W$.
\item $\Ker f$, $\Image(f)$ and $\Cok(f)$ with the induced filtration $W$
      are mixed twistor structures on $X$.
\end{itemize}
As a result, the category of mixed twistor structure on $X$
is abelian.
\end{lem}
\pf
We only remark that
$\rank (f_{|(\lambda,P)})$ are independent of
$(\lambda,P)\in \proj^1\times X$.
Indeed,
we obtain 
$\rank (f_{|(\lambda,P)})
=\rank(f_{|(\lambda',P)})$  for any $\lambda,\lambda'\in\proj^1$
by Proposition \ref{prop;22.3.24.1},
and 
$\rank (f_{|(1,P)})
=\rank(f_{|(1,P')})$  for any $P,P'\in X$
because
$f_{|\{1\}\times X}$
is a morphism of flat bundles
$(\nbigv_{1|\{1\}\times X},\DD^{\sankaku}_{1|\{1\}\times X})
\to
(\nbigv_{2|\{1\}\times X},\DD^{\sankaku}_{2|\{1\}\times X})$.
Therefore, $\Ker(f)$, $\Image(f)$ and $\Cok(f)$ are
naturally vector bundles.
The other claims easily follow from 
Proposition \ref{prop;22.3.24.1}.
\hfill\qed

\subsubsection{Polarizations}
\label{subsection;22.3.25.1}
For $\nbigv\in\TS(X)$,
a morphism
$S:\nbigv\otimes\sigma^{\ast}(\nbigv)\to \Tate(-n)_X$
is called a pairing of $\nbigv$ of weight $n$.
It is called $(-1)^n$-symmetric if
$(-1)^n\iota_{\Tate(-n)}\circ \sigma_X^{\ast}(S)
=S\circ\exchange$ holds,
where $\exchange$ denotes the natural isomorphism
$\nbigv\otimes\sigma^{\ast}\nbigv
\simeq
\sigma^{\ast}(\nbigv)\otimes\nbigv$.
\begin{df}[\cite{s3}]
Let $\nbigv$ be a pure twistor structure of weight $n$ on $X$.
Let $S:\nbigv\otimes\sigma_X^{\ast}\nbigv\to \Tate(-n)_X$ 
be a $(-1)^n$-symmetric pairing.
It is called a polarization of $\nbigv$
if $\iota_P^{\ast}(S)$
is a polarization of $\iota_P^{\ast}(\nbigv)$
for any $P\in X$.
\hfill\qed
\end{df}

\begin{df}
A pure twistor structure $\nbigv$ of weight $n$ is called polarizable
if there exists a polarization of $\nbigv$.
A mixed twistor structure $(\nbigv,W)$ is called graded polarizable
if each $\Gr^W_n(\nbigv)$ is polarizable.
\hfill\qed
\end{df}

We mention some instructive lemmas.
\begin{lem}
\label{lem;22.3.25.20}
Let $\nbigv$ be a pure and polarizable twistor structure of weight $n$ on $X$.
Let $\nbigv'\subset\nbigv$ be a subobject in $\TS(X)$,
which is pure of weight $n$.
Then, $\nbigv'$ is also pure and polarizable.
Moreover, there exists a decomposition
$\nbigv=\nbigv'\oplus\nbigv''$,
and $\nbigv''$ is also pure and polarizable of weight $n$.
\end{lem}
\pf
It is enough to consider the case $n=0$.
Let $S$ be a polarization of $\nbigv$.
Let $S'$ be a symmetric pairing of weight $0$
of $\nbigv'$ induced by $S$.
It is easy to check that $S'$ is also a polarization of $\nbigv'$.
We have the morphism
$\Psi_S:\nbigv\lrarr \sigma^{\ast}(\nbigv)^{\lor}$
induced by $S$.
Let $\nbigv''$ denote the kernel of the composite of
$\Psi_S$ and $\sigma^{\ast}(\nbigv)\to\sigma^{\ast}(\nbigv')$.
It is also pure of weight $0$,
and we can easily check that
$\nbigv=\nbigv'\oplus\nbigv''$.
\hfill\qed

\begin{lem}
\label{lem;22.3.25.22}
The category of pure and polarizable twistor structures on $X$
is abelian and semisimple.
The category of graded polarizable mixed twistor structure on $X$
is abelian. 
\end{lem}
\pf
The first claim follows from Lemma \ref{lem;22.3.25.20}.
The second claim follows from the first claim and
Lemma \ref{lem;22.3.25.21}.
\hfill\qed

\subsubsection{Polarizations and harmonic bundles}

Let us recall the notion of harmonic bundles.
Let $(E,\delbar_E)$ be a holomorphic vector bundle on $X$.
Let $\theta$ be a Higgs field of $(E,\delbar_E)$,
i.e.,
$\theta$ is a holomorphic section of
$\End(E)\otimes\Omega^{1,0}_X$
such that $\theta\wedge\theta=0$.
Let $h$ be a Hermitian metric of $E$.
We obtain the Chern connection
$\nabla_h=\delbar_E+\del_{E,h}$
from $h$ and $\delbar_E$,
and the adjoint $\theta^{\dagger}_h$
from $h$ and $\theta$.
Then, $h$ is called a pluri-harmonic metric
if $\DD^1=\nabla_h+\theta+\theta^{\dagger}_h$ is flat.
Such a tuple $(E,\delbar_E,\theta,h)$ is called a harmonic bundle.

From a harmonic bundle
$(E,\delbar_E,\theta,h)$,
we set
$V^{\sankaku}(E):=\pi_2^{-1}(E)$ on $\proj^1\times X$.
It is naturally equipped with
the $\proj^1$-holomorphic structure
$\delbar_{\proj^1,V^{\sankaku}(E)}$.
We set
$\DD^{\sankaku}=\nabla_h+\lambda^{-1}\theta+\lambda\theta^{\dagger}$.
Then,
$\nbigv(E)=(V^{\sankaku}(E),\delbar_{\proj^1,V^{\sankaku}(E)},\DD^{\sankaku})$
is a pure twistor structure of weight $0$ on $X$.
We obtain the pairing
$S_h:V^{\sankaku}(E)\otimes\sigma^{\ast}V^{\sankaku}(E)
\to \Tate(0)_X$
by
$S_h(u,\sigma^{-1}(v))=h(u,\sigma^{-1}(v))$.
Then,
$(\nbigv(E),S_h)$ is a polarized pure twistor structure of weight $0$ on $X$.
The following theorem due to Simpson \cite{s3} is fundamental.
It is also explained in \cite{sabbah2,mochi2}.
\begin{thm}[Simpson]
The above procedure induces an equivalence
between harmonic bundles and 
polarized pure twistor structure of weight $0$ on $X$.
\hfill\qed
\end{thm}

\subsubsection{Deformation associated with
a commuting tuple of nilpotent endomorphisms}
\label{subsection;22.3.24.4}

Let $X$ be a complex manifold,
and let $\Lambda$ be a finite set.
Let $\TS(X,\Lambda)$ denote the category of
twistor structure $\nbigv$
equipped with a tuple $\vecf=(f_i\,|\,i\in\Lambda)$
of morphisms $f_i:\nbigv\to\nbigv\otimes\Tate(-1)_X$
such that $[f_i,f_j]=0$
and that $f_i^N=0$ for some $N>0$.
A morphism $\varphi:(\nbigv_1,\vecf_1)\to (\nbigv_2,\vecf_2)$
in $\TS(X,\Lambda)$
is a morphism
$\varphi:\nbigv_1\to\nbigv_2$ in $\TS(X)$
such that
$f_{2,i}\circ\varphi=\varphi\circ f_{1,i}$.

Let $\Lambda_1\subset\Lambda$ be any finite subset.
We set $Y(\Lambda_1)=X\times\cnum^{\Lambda_1}$.
We set $H_i=X\times\cnum^{\Lambda_1\setminus\{i\}}$,
and $H(\Lambda_1)=\bigcup_{i\in\Lambda_1} H_i\subset Y(\Lambda_1)$.
We set $Y^{\circ}(\Lambda_1):=Y(\Lambda_1)\setminus H(\Lambda_1)$.
Let us construct a functor
$\TNIL:
\TS(X,\Lambda)\to\TS(Y(\Lambda_1)^{\circ},\Lambda)$
which is a minor refinement of the construction
in \cite[\S3.6]{mochi2} and \cite[\S2.2.3]{mochi8}.
(See also \cite[\S2.2]{Mochizuki-MTM}
for the construction in the context of $\nbigr$-triples.)

Let $(\nbigv,\vecf)\in\TS(X,\Lambda)$.
We have $\nbigv=(V^{\sankaku},\delbar_{\proj^1,V^{\sankaku}},
\DD^{\sankaku})\in\TS(X)$.
We set
$V_0:=V^{\sankaku}_{|\cnum_{\lambda}\times X}$
and
$V_{\infty}:=V^{\sankaku}_{|\cnum_{\mu}\times X}$.
We obtain the following differential operators
as the restriction of $\DD^{\sankaku}$:
\[
\DD^f_{0}:
V_0\to
V_0\otimes
\bigl(
\Omegatilde^{1,0}_{\nbigx/\cnum_{\lambda}}
\oplus
\Omegatilde^{0,1}_{\nbigx/\cnum_{\lambda}}
\bigr),
\]
\[
 \DD^{\dagger\,f}_{\infty}:
V_{\infty}\to
V_{\infty}\otimes
\bigl(
\Omegatilde^{1,0}_{\nbigx^{\dagger}/\cnum_{\mu}}
\oplus
\Omegatilde^{0,1}_{\nbigx^{\dagger}/\cnum_{\mu}}
\bigr).
\]
We obtain the morphisms
$V_0\to \lambda^{-1}V_0$
and
$V_{\infty}\to \mu^{-1}V_{\infty}$
as the restriction of $f_i$
to $\cnum_{\lambda}\times X$
and $\cnum_{\mu}\times X$, respectively.
They are also denoted by $f_i$.

Let $q_X:Y(\Lambda_1)^{\circ}\to X$ denote the projection.
We obtain the $C^{\infty}$-vector bundles
$q_X^{-1}V_0$ and $q_X^{-1}V_{\infty}$
on $\nbigy^{\circ}(\Lambda_1)$ and $\nbigy^{\circ\dagger}(\Lambda_1)$,
respectively.
We have the following differential operator
\[
 \DD^f_{q_X^{-1}(V_0)}:
 q_X^{-1}(V_0)
 \lrarr
 q_X^{-1}(V_0)
 \otimes
 \bigl(
\Omegatilde^{1,0}_{\nbigy^{\circ}(\Lambda_1)/\cnum_{\lambda}}
\oplus
\Omegatilde^{0,1}_{\nbigy^{\circ}(\Lambda_1)/\cnum_{\lambda}}
\bigr),
\]
\[
\DD^f_{q_X^{-1}(V_0)}:=
q_X^{\ast}\DD^f_{0}
-\sum_{i\in\Lambda_1}
\sqrt{-1}q_X^{\ast}(f_{i})
\frac{dz_i}{z_i}.
\]
We also have the following differential operator
\[
 \DD^f_{q_X^{-1}(V_{\infty})}:
 q_X^{-1}(V_{\infty})
 \lrarr
 q_X^{-1}(V_{\infty})
 \otimes
 \bigl(
\Omegatilde^{1,0}_{\nbigy^{\circ}(\Lambda_1)^{\dagger}/\cnum_{\mu}}
\oplus
\Omegatilde^{0,1}_{\nbigy^{\circ}(\Lambda_1)^{\dagger}/\cnum_{\mu}}
\bigr),
\]
\[
\DD^f_{q_X^{-1}(V_{\infty})}:=
q_X^{\ast}\DD^f_{\infty}
+\sum_{i\in\Lambda_1}
\sqrt{-1}q_X^{\ast}(f_{i})
\frac{d\zbar_i}{\zbar_i}.
\]

There exists the natural $C^{\infty}$-isomorphism
$\Psi:V_{0|\cnum_{\lambda}^{\ast}\times X}
\simeq
V_{\infty|\cnum_{\mu}^{\ast}\times X}$.
It induces the $C^{\infty}$-isomorphism
$\Psitilde:
q_{X}^{-1}(V_{0})_{|\cnum_{\lambda}^{\ast}\times Y^{\circ}(\Lambda_1)}
\simeq
q_{X}^{-1}(V_{\infty})_{|\cnum_{\mu}^{\ast}\times Y^{\circ}(\Lambda_1)}$
defined by
\[
 \Psitilde
 =\Psi\circ
 \exp\Bigl(
 -\sum_{i\in\Lambda_1}
 \log|z_i|^2
 \cdot
 \sqrt{-1}q_X^{\ast}(f_i)
 \Bigr).
\]
We obtain the $C^{\infty}$-vector bundle
$\Vtilde^{\sankaku}$
by gluing
$q_X^{-1}(V_0)$
and
$q_X^{-1}(V_{\infty})$
by $\Psitilde$.
It is equipped with the differential operator
$\DD^{\sankaku}_{\Vtilde^{\sankaku}}$
by gluing
$\DD^f_{q_X^{-1}(V_0)}$
and
$\DD^f_{q_X^{-1}(V_{\infty})}$.
Because $\Psitilde$ is holomorphic in the $\proj^1$-direction,
there exists the naturally defined
$\proj^1$-holomorphic structure
$\delbar_{\proj^1,\Vtilde^{\sankaku}}$.
In this way,
we obtain a twistor structure
$\TNIL(\nbigv,\vecf_{\Lambda_1})=
\bigl(
\Vtilde^{\sankaku},
\delbar_{\proj^1,\Vtilde^{\sankaku}},
\DD^{\sankaku}_{\Vtilde^{\sankaku}}
\bigr)$.
We obtain
$f_i:\TNIL(\nbigv,\vecf)
\to \TNIL(\nbigv)\otimes\Tate(-1)_{Y^{\circ}(\Lambda_1)}$ $(i\in\Lambda)$.
Thus, we obtain
$\TNIL(\nbigv,\vecf_{\Lambda_1})$ with $\vecf$
in $\TS(Y^{\circ}(\Lambda_1),\Lambda)$.

\begin{rem}
\label{rem;22.3.24.3}
The signature of nilpotent endomorphisms
are chosen to be opposite from those in {\rm\cite{mochi2,mochi8}},
where
$\TNIL(\nbigv,-\vecf)$ is considered
 for $(\nbigv,\vecf)$.
As remarked in {\rm\cite[Remark 3.123]{mochi2}},
the choice of the signature in {\rm\cite{mochi2,mochi8}}
is chosen to be
opposite from the standard choice in the Hodge theory.
See Remark {\rm\ref{rem;22.3.24.2}} below.
\hfill\qed
\end{rem}

\subsubsection{Twistor nilpotent orbits}

Let $(\nbigv,\vecf)\in \TS(X,\Lambda)$.
A $(-1)^n$-symmetric pairing
$S$ of $\nbigv$ of weight $n$ is called
a $(-1)^n$-symmetric pairing of
$(\nbigv,\vecf)$ of weight $n$
if $S\circ(\id\otimes \sigma^{\ast}(f_i))+S\circ(f_i\otimes\id)=0$
$(i\in\Lambda)$.
It induces a $(-1)^n$-symmetric pairing
of 
$\bigl(
\TNIL(\nbigv,\vecf),\vecf
\bigr)$
of weight $n$.
We recall the following definition
\cite[\S3.6]{mochi2}.

\begin{df}
If there exists a neighbourhood $U$ of
$X\times\{(0,\ldots,0)\}$ in $Y(\Lambda)$
such that
\[
  (\TNIL(\nbigv,\vecf),\TNIL(S))_{|U\setminus H}
\]
is a polarized pure twistor structure of weight $n$,
$(\TNIL(\nbigv,\vecf),\TNIL(S))_{|U\setminus H}$
is called a twistor nilpotent orbit of weight $n$.
\hfill\qed
\end{df}

\subsubsection{Polarized mixed twistor structures}
\label{subsection;22.3.24.30}

Let $(\nbigv,W)$ be a mixed twistor structure on $X$.
Let $f:(\nbigv,W)\to(\nbigv,W)\otimes\Tate(-1)_X$
be a morphism.
Note that it is automatically nilpotent.
We recall the following definition \cite[Definition 3.48]{mochi2}.
\begin{df}
$(\nbigv,W,f)$ is called a polarizable
mixed twistor structure of weight $w$
if the following conditions are satisfied.
\begin{itemize}
 \item The weight filtration $W(f)$ of $f$
       is equal to $W$ up to shift by $w$,
       i.e., $W_{w+m}=W(f)_{m}$.
       In other words,
       $f^k:\Gr^W_{w+k}(\nbigv)\simeq
        \Gr^W_{w-k}(\nbigv)\otimes\Tate(-k)$
       for any $k\geq 0$.
\item There exists a $(-1)^w$-symmetric pairing
$S:\nbigv\otimes\sigma^{\ast}\nbigv
 \lrarr \Tate(-w)_X$ such that
      $S\circ (f\otimes\id)
      +S(\id\otimes \sigma^{\ast}(f))=0$.
 \item Moreover, for $k\geq 0$,
       $S(\id\otimes \sigma^{\ast}(f)^k)=
       (-1)^kS(f^k\otimes\id)$
       is a polarization of
       the pure twistor structure
       $P\Gr^W_{w+k}(\nbigv)$ of weight $w+k$.
       Here, $P\Gr^W_{w+k}(\nbigv)$ denotes the primitive part,
       i.e.,
       the kernel of
       $f^{k+1}:\Gr^{W}_{w+k}(\nbigv)\lrarr
       \Gr^W_{w-k-2}(\nbigv)
       \otimes\Tate(-k-1)$.
\end{itemize} 
 Such $S$ is called a polarization of
 $(\nbigv,W,f)$,
and such a tuple $(\nbigv,W,f,S)$ is called 
a polarized mixed twistor structure of weight $w$.
\hfill\qed
\end{df}

\begin{rem}
\label{rem;22.3.24.2}
We changed the signature condition
from that in {\rm\cite{mochi2,mochi8}},
where we adopted the condition that
$(-1)^kS(\id\otimes \sigma^{\ast}(f)^k)=
S(f^k\otimes\id)$
is a polarization of
$P\Gr^W_{w+k}(\nbigv)$.
See Remark {\rm\ref{rem;22.3.24.3}}.
\hfill\qed 
\end{rem}

Let $\Lambda$ be a finite set.
Let $\vecf=(f_i\,|\,i\in\Lambda)$ be a commuting tuple of
morphisms
$f_i:(\nbigv,W)\to (\nbigv,W)\otimes\Tate(-1)$.
For any $\veca=(a_i\,|\,i\in\Lambda)\in\real^{\Lambda}$,
we set $f(\veca)=\sum a_if_i$.
We recall the following definition in \cite[Definition 3.50]{mochi2}.
\begin{df}
$(\nbigv,W,\vecf)$ is called
a $(w,\Lambda)$-polarizable mixed twistor structure
if there exists
a $(-1)^w$-symmetric pairing
$\nbigv\otimes\sigma^{\ast}\nbigv
\lrarr \Tate(-w)$ such that
(i) $S(f_i\otimes \id)+S(\id\otimes \sigma^{\ast}(f_i))=0$
$(i\in\Lambda)$,
(ii) $(\nbigv,W,f(\veca),S)$
is a polarized mixed twistor structure of weight $w$
for any $\veca\in\real_{>0}^{\Lambda}$.
Such $S$ is called a polarization of
$(\nbigv,W,\vecf)$,
and such a tuple $(\nbigv,W,\vecf,S)$
is called a $(w,\Lambda)$-polarized mixed twistor structure.
\hfill\qed 
\end{df}

We have the following variant.
\begin{df}
\label{df;22.3.24.11}
Let $(\nbigv,\vecf)\in\TS(X,\Lambda)$.
Let $S$ be a $(-1)^w$-symmetric pairing of
$(\nbigv,\vecf)$ of weight $w$.
We say that $(\nbigv,f,S)$ is
a $(w,\Lambda)$-polarized mixed twistor structure
if the following holds. 
\begin{itemize}
 \item We define the filtration $W$ of $V$
       by $W_{w+k}=W_k(\sum_{i\in\Lambda} f_i)$.
      Then, $(\nbigv,W,f,S)$ is a polarized mixed twistor structure
       of weight $w$.
       \hfill\qed
\end{itemize} 
\end{df}
Note that the condition implies that
the monodromy weight filtrations
$W(f(\veca))$ are independent of $\veca\in\real_{>0}^{\Lambda}$.
It is clear that the two concepts of
polarized mixed twistor structure are equivalent.

We recall the following proposition.

\begin{prop}[\mbox{\cite[Theorem 4.2]{mochi8}}]
\label{prop;22.3.24.10}
Let $S$ be a $(-1)^n$-symmetric pairing of
$(\nbigv,W,\vecf)\in\TS(X,\Lambda)$ of weight $n$.
Let $Y(\Lambda)$ and $H$ be as in 
{\rm\S\ref{subsection;22.3.24.4}}.
Then, 
$(\nbigv,\vecf,S)$ is a $(w,\Lambda)$-polarized mixed twistor
 structure on $X$
in the sense of Definition {\rm\ref{df;22.3.24.11}} 
if and only if there exists a neighbourhood $U$ of
$X\times\{(0,\ldots,0)\}$ in $Y$ such that
$\TNIL(\nbigv,\vecf,S)$
is a polarized pure twistor structure of weight $w$
on $U\setminus H$.
\hfill\qed
\end{prop}

Let $\Lambda=\Lambda_1\sqcup\Lambda_2$
be a decomposition.
We set $\vecf_{\Lambda_j}:=(f_i\,|\,i\in\Lambda_j)$.
We can regard
$Y(\Lambda_1)\subset Y(\Lambda)$
in a natural way.
The following corollary is useful.
\begin{cor}
\label{cor;22.3.24.100}
Let $(\nbigv,\vecf,S)$ be
a $(w,\Lambda)$-polarized mixed twistor structure on $X$.
Let $U$ be a neighbourhood of
$X\times\{(0,\ldots,0)\}$ in $Y(\Lambda)$
as in Proposition {\rm\ref{prop;22.3.24.10}}.
Then,
$(\TNIL(\nbigv,\vecf_{\Lambda_1}),\vecf_{\Lambda_2})
 _{|Y^{\circ}(\Lambda_1)\cap U}$
is a $(w,\Lambda_2)$-polarizable mixed twistor structure. 
In particular, the following holds.
\begin{itemize}
 \item We set $W_{w+k}=W(\sum_{i\in\Lambda_2}f_i)$ $(k\in\seisuu)$
       on
       $\TNIL(\nbigv,\vecf_{\Lambda_1})_{|Y^{\circ}(\Lambda_1)\cap U}$.
       Then,
       $(\TNIL(\nbigv,\vecf_{\Lambda_1})_{|Y^{\circ}(\Lambda_1)\cap U},W)$
       is a mixed twistor structure.
\hfill\qed
\end{itemize}
\end{cor}

\subsubsection{Vanishing cycle theorem}
\label{subsection;22.4.3.31}

Let $(\nbigv,W,\vecf,S)$
be a $(w,\Lambda)$-polarized mixed twistor structure on $X$.
For any $J\subset\Lambda$,
we set
$f_J=\prod_{j\in J}f_j:
(\nbigv,W)\to (\nbigv,W)\otimes\Tate(-|J|)$.
We obtain the mixed twistor structure
$(\Image (f_J),W)\subset (\nbigv,W)\otimes\Tate(-|J|)$.
It is equipped with a tuple of the induced morphisms
$f_i:(\Image(f_J),W)\to (\Image(f_J),W)\otimes\Tate(-1)$
$(i\in \Lambda)$.
We obtain the induced pairing
$S:\nbigv\otimes\sigma^{\ast}\Image (f_J)\to\Tate(-w-|J|)$,
which vanishes on
$\Ker(f_J)\otimes\sigma^{\ast}\Image(f_J)$.
Hence, we obtain
$S_J:\Image(f_J)\otimes\sigma^{\ast}\Image(f_J)
\to \Tate(-w-|J|)$.
We recall the following proposition
which is an easy consequence of
the vanishing cycle theorem
in the Hodge case due to
Cattani-Kaplan-Schmid \cite{cks2}
and Kashiwara-Kawai \cite{k3}.

\begin{prop}[\mbox{\cite[Proposition 3.126]{mochi2}}]
\label{prop;22.2.6.1}
$(\Image(f_J),W,\vecf,S_J)$
is a polarized mixed twistor structure of weight $w+|J|$.
\hfill\qed
\end{prop}

For any subset $K\subset\Lambda$,
let $\lefttop{K}W$ denote the monodromy weight filtration
of $\sum_{j\in K} f_j$.
By definition,
we have
$\lefttop{\Lambda}W_{m}(\nbigv)=W_{m+w}(\nbigv)$.
Proposition \ref{prop;22.2.6.1} implies
$\lefttop{\Lambda}W_m(\Image f_J)
=W_{m+w+|J|}(\Image f_J)$,
and hence
\begin{equation}
\label{eq;22.2.6.2}
 f_J\bigl(
 \lefttop{\Lambda}W_{m+|J|}(\nbigv)
 \bigr)
=\lefttop{\Lambda}W_{m}(\Image f_J)
=\Image(f_J)\cap\lefttop{\Lambda}W_{m-|J|}(\nbigv)
\otimes\Tate(-|J|).
\end{equation}

We obtain the following refinement.
\begin{prop}
\label{prop;22.2.17.30}
For $J\subset K\subset\Lambda$,
we have
\begin{equation}
\label{eq;22.3.24.20}
  f_J\bigl(
 \lefttop{K}W_{m+|J|}(\nbigv)
 \bigr)
=\lefttop{K}W_{m}(\Image f_J)
=\Image(f_J)\cap\bigl(
\lefttop{K}W_{m-|J|}(\nbigv)
\otimes\Tate(-|J|)\bigr).
\end{equation}
\end{prop}
\pf
Let us explain how to obtain this kind of claim
by using twistor nilpotent orbits.
We set $L:=\Lambda\setminus K$.
There exists a neighbourhood
$U$ of $X\times\{(0,\ldots,0)\}$ in $X\times\cnum^L$
such that
$(\TNIL(\nbigv,\vecf_L),\vecf_K)_{|U\setminus H(L)}$
is a $(w,K)$-polarized mixed twistor structure on
$U\setminus H(L)$.
By Proposition \ref{prop;22.2.6.1},
we obtain (\ref{eq;22.3.24.20})
for 
$(\TNIL(\nbigv,\vecf_L),\vecf_K)_{|U\setminus H(L)}$.
Then, by the construction of $\TNIL(\nbigv,\vecf_L)$,
we obtain (\ref{eq;22.3.24.20})
for $(\nbigv,\vecf_K)$.
\hfill\qed

\vspace{.1in}

Let $J\sqcup\{j\}\subset \Lambda$.
If $j\in K$,
we obtain
\[
 f_j:\lefttop{K}W_{k+1}(\Image f_J)
 \lrarr
 \lefttop{K}W_{k}(\Image f_{J\sqcup\{j\}})
 =\Bigl(
 \lefttop{K}W_{k-1}(\Image(f_J))\otimes\Tate(-1)
 \Bigr)
 \cap
 \Image (f_{J\sqcup\{j\}}).
\]
If $j\not\in K$,
we obtain
\[
 f_j:
 \lefttop{K}W_{k}(\Image f_J)
 \lrarr
 \lefttop{K}W_{k}(\Image f_{J\sqcup\{j\}}).
\]

\subsection{$\nbigr_{X(\ast H)}$-modules}
\label{subsection;22.4.28.2}

Let $X$ be a complex manifold.
We set $\nbigx:=\cnum_{\lambda}\times X$.
Let $p_{\lambda}:\nbigx\lrarr X$ denote the projection.
Let $\nbigd_{\nbigx}$ denote the sheaf of holomorphic
linear differential operators on $\nbigx$.
Let $\Theta_X$ denote the tangent sheaf of $X$.
Let $\nbigr_X\subset\nbigd_{\nbigx}$ denote
the sheaf of subalgebras generated by
$\lambda p_{\lambda}^{\ast}\Theta_X$
over $\nbigo_{\nbigx}$.
Let $H$ be a hypersurface of $X$.
We set $\nbigh:=\cnum_{\lambda}\times X$.
For any $\nbigo_{\nbigx}$-module $\nbigf$,
we set $\nbigf(\ast H):=\nbigf\otimes\nbigo_{\nbigx}(\ast\nbigh)$.
We put $\nbigr_{X(\ast H)}:=\nbigr_X(\ast H)$.
We recall some standard concepts for $\nbigr_X$-modules
by following
\cite[\S0,\S1]{sabbah2} and \cite[\S2]{Mochizuki-MTM}.

A left $\nbigr_{X(\ast H)}$-module $\nbigm$
is equivalent to
an $\nbigo_{\nbigx}$-module $\nbigm$
equipped with a $\cnum$-linear operator
\[
 \DD^f:\nbigm\lrarr
  \nbigm\otimes\bigl(
  \lambda^{-1}
  p_{\lambda}^{\ast}\Omega^1_X(\ast H)
  \bigr)
\]
such that
$\DD^f(gs)=(dg)\cdot s+g\DD^f(s)$
for local sections $g$ and $s$
of $\nbigo_{\nbigx}$ and $\nbigm$,
respectively.
By setting $\DD=\lambda\DD^f$,
we obtain
\[
 \DD:\nbigm\lrarr
 \nbigm\otimes p_{\lambda}^{\ast}\Omega^1_X(\ast H)
\]
such that 
$\DD(gs)=(\lambda dg)\cdot s+g\DD(s)$
for local sections $g$ and $s$
of $\nbigo_{\nbigx}$ and $\nbigm$,
respectively.

An $\nbigr_{X(\ast H)}$-module $\nbigm$ is called strict
if it is flat over $\nbigo_{\cnum_{\lambda}}$.
For a strict $\nbigr_{X(\ast H)}$-module $\nbigm$,
$\nbigm^{\lambda_0}:=
\iota_{\lambda_0}^{-1}(\nbigm/(\lambda-\lambda_0)\nbigm)$
for any $\lambda_0\in\cnum$,
where $\iota_{\lambda_0}:\{\lambda_0\}\times X\lrarr \nbigx$
denote the inclusion.

The pull back of $\nbigo_{\cnum_{\lambda}}$
by the projection $\nbigx\to \cnum$
is denoted by $(\nbigo_{\cnum_{\lambda}})_{\nbigx}$.
For any open subset $U\subset\nbigx$,
the restriction to $U$ is denoted by
$(\nbigo_{\cnum_{\lambda}})_U$.
A holonomic $\nbigr_{X}$-module on $U$
is called strict
if it is flat over $(\nbigo_{\cnum_{\lambda}})_U$.

The sheaf of algebras $\nbigr_{X(\ast H)}$
has the filtration by the order of differential operators.
We have the notion of coherent filtration for
$\nbigr_{X(\ast H)}$-modules as in the case of
$\nbigd_X$-modules.
(See \cite{kashiwara_text}. It is also called a good filtration.)
A coherent $\nbigr_{X(\ast H)}$-module
has a coherent filtration locally around any point of $X$.
The associated graded module induces a coherent sheaf
on $\cnum_{\lambda}\times T^{\ast}(X\setminus H)$.
The support is called the characteristic variety of $\nbigm$
denoted by $\Ch(\nbigm)$.
A strict coherent $\nbigr_{X(\ast H)}$-module
$\nbigm$ is called holonomic
if there exists a homogeneous Lagrangian subvariety $\Lambda$
of $T^{\ast}(X\setminus H)$
such that $\Ch(\nbigm)\subset \cnum\times\Lambda$.

An $\nbigr_{X(\ast H)}$-module $\nbigm$
on $U\subset\nbigx$
is called good if for any compact subset $K\subset U$
there exist a neighbourhood $U'$ of $K$ in $U$
and a finite filtration $F$ of $\nbigm_{|U'}$
such that $\Gr^F(\nbigm_{|U'})$
has a coherent filtration.
Good $\nbigr_{X(\ast H)}$-modules are
$\nbigr_{X(\ast H)}$-coherent.

Let $f:X'\to X$ be a morphism of complex manifolds.
We set $H':=f^{-1}(H)$.
The induced morphism $\nbigx'\to\nbigx$
is also denoted by $f$.
We set
$\nbigr_{X'\to X}:=
\nbigo_{\nbigx'}\otimes_{f^{-1}\nbigo_{\nbigx}} f^{-1}\nbigr_X$,
which is naturally an $(\nbigr_{X'},f^{-1}\nbigr_X)$-bimodule.
For any $\nbigr_{X(\ast H)}$-module $\nbign$,
we set
$f^{\dagger}(\nbign)
=\nbigr_{X'\to X}\otimes^L_{f^{-1}\nbigr_X}\nbign$
in
$D^b(\nbigr_{X'(\ast H')})$.
Let $\omega_{\nbigx}$ denote
the sheaf of holomorphic sections of
$\Omegatilde^{\dim X,0}_{\nbigx/\cnum_{\lambda}}$.
We use the notation $\omega_{\nbigx'}$ with a similar meaning.
We set
$\nbigr_{X\larr X'}:=
\omega_{\nbigx'}\otimes \nbigr_{X'\to X}\otimes f^{-1}(\omega_{\nbigx}^{-1})$.
For any $\nbigr_{X'(\ast H')}$-module $\nbigm$,
we set
$f_{\dagger}(\nbigm):=
 Rf_{\ast}\Bigl(
 \nbigr_{X\larr X'}\otimes^L_{\nbigr_{X'}}
 \nbigm
 \Bigr)$
in $D^b(\nbigr_X)$.
The $i$-th cohomology sheaves are denoted by
$f_{\dagger}^i(\nbigm)$.
If $\nbigm$ is good,
and if the support of $\nbigm$ is proper over $X$,
$f_{\dagger}^i(\nbigm)$ are also good.
We decompose $f$ into the composition of
the graph embedding
$\iota_f:X\lrarr X\times X'$
and the projection
$\pi_f:X\times X'\to X'$.
Let $p_1:\cnum_{\lambda}\times X\times X'\to \cnum_{\lambda}\times X$
denote the projection.
Then, we have
\[
 f_{\dagger}^i(\nbigm)
 =R^i(\id_{\cnum_{\lambda}}\times \pi)_{\ast}
 \Bigl(
 \iota_{f\dagger}\nbigm
 \otimes
 p_1^{-1}
 \Tot
 \Omegatilde^{\bullet,\bullet}_{\nbigx/\cnum}[\dim X]
 \Bigr).
\]
Let $c\in H^{2}_{\DR}(X)$ be
a cohomology class induced by
a closed $(1,1)$-form.
Taking a closed $(1,1)$-form $\omega$
representing $c$,
we define
$L_{c}:
f_{\dagger}^i(\nbigm)
\to
\lambda f_{\dagger}^{i+2}(\nbigm)$
by the multiplication of $2\pi\sqrt{-1}\omega$,
which is independent of the choice of $\omega$.
\begin{rem}
The signature of $L_c$ is chosen to be opposite to that
in {\rm\cite{sabbah2}}.
The induced conditions for polarizations
are equivalent.
\hfill\qed
\end{rem}

\subsubsection{Strict specializability, the canonical prolongations,
and strict $S$-decomposability}
\label{subsection;22.3.28.10}

For the convenience of readers,
we recall the notion of $V$-filtration and
the strict specializability condition for $\nbigr_{X(\ast H)}$-modules
by following \cite{sabbah2,Mochizuki-MTM}.
Let $\cnum_t$ be a complex line with the standard coordinate $t$.
Let $X_0$ be a complex manifold.
For simplicity, $X_0$ is assumed to be connected.
Let $X$ be a neighbourhood of $X_0\times\{0\}$ in $X_0\times\cnum_t$.
We regard $X_0\subset X$ by using the identification
$X_0=X_0\times\{0\}$.

Let $\Theta_X(\log X_0)$ denote the sheaf of vector fields on $X$
which are logarithmic along $X_0$.
Let $V\nbigr_X\subset \nbigr_X$
denote the sheaf of subalgebras generated by
$\lambda p_{\lambda}^{\ast}\Theta_X(\log X_0)$
over $\nbigo_{\nbigx}$.

\paragraph{The case $X_0\not\subset H$}
Let $H$ be a hypersurface of $X$
such that $X_0\not\subset H$.
We set $H_0:=X_0\cap H$.
We set $\nbigr_{X(\ast H)}=\nbigr_X(\ast H)$.
For any $\lambda_0$,
let $\nbigxzero$ denote a neighbourhood of
$\{\lambda_0\}\times X$ in $\nbigx$.
We use the notation $\nbigxzero_0$ in a similar meaning.
A coherent strict $\nbigr_{X(\ast H)}$-module $\nbigm$
is strictly specializable along $t$ at $\lambda_0$
if $\nbigm_{|\nbigxzero}$ is equipped with
an increasing and exhaustive filtration
$\Vzero\nbigm=(\Vzero_a(\nbigm)\,|\,a\in\real)$
by coherent $V\nbigr_{X(\ast H)}$-modules
satisfying the following conditions.
\begin{description}
 \item[(i)] For any $a\in\real$
       and $P\in\nbigxzero_0$,
       there exists $\epsilon>0$ such that
       $\Vzero_{a+\epsilon}(\nbigm)=\Vzero_a(\nbigm)$
       on a neighbourhood of $P$.
 \item[(ii)]
	    Set $\Vzero_{<a}(\nbigm):=\bigcup_{b<a}\Vzero_b(\nbigm)$,
	    and $\Gr^{\Vzero}_a(\nbigm):=\Vzero_a(\nbigm)/\Vzero_{<a}(\nbigm)$.
	    Then, 
	    $\Gr^{\Vzero}_a(\nbigm)$ are strict,
	    i.e., flat over $\nbigo_{\cnum_{\lambda}}$.

 \item[(iii)] $t\Vzero_{a}(\nbigm)\subset \Vzero_{a-1}(\nbigm)$
       for any $a$,
       and $t\Vzero_{a}(\nbigm)=\Vzero_{a-1}(\nbigm)$
       for any $a<0$.
 \item[(iv)] Set $\deldel_t:=\lambda\del_t$.
       Then,
       we have $\deldel_t\Vzero_a(\nbigm)\subset\Vzero_{a+1}(\nbigm)$
       for any $a\in\real$,
       and
       the induced morphisms
       $\deldel_t:\Gr^{\Vzero}_a(\nbigm)
       \to \Gr^{\Vzero}_{a+1}(\nbigm)$
       are isomorphisms for any $a>-1$.
 \item[(v)] For any $a\in\real$,
       and $P\in\nbigxzero_0$,
       there exists a finite subset
\[
       \nbigk(a,\lambda_0,P)\subset
       \bigl\{
       u\in\real\times\cnum\,\big|\,
       \paramap(\lambda_0,u)=a
       \bigr\}
\]
       such that
       $\prod_{u\in\nbigk(a,\lambda_0,P)}
       (-\deldel_tt+\eigenmap(\lambda,u))$
       are nilpotent on
       $\Gr^{\Vzero}_a(\nbigm)$
       on a neighbourhood of $P$.
       Here,
\[
       \paramap(\lambda,(b,\beta))
       =b+2\Re(\lambda\betabar),
       \quad\quad
       \eigenmap(\lambda,(b,\beta))
       =\beta-b\lambda-\betabar\lambda^2.
\]
       We implicitly assume that
       $\nbigk(a,\lambda_0,P)$ is minimal
       among such finite subsets.
\end{description}
If such a $V$-filtration exists, it is unique.
We say that $\nbigm$ is strictly specializable along $t$
if it is strictly specializable along $t$ at any $\lambda_0$.

Let $\nbigm$ be a strict coherent $\nbigr_{X(\ast H)}$-module
which is strictly specializable along $t$.
For any $\lambda_0\in\cnum$ and $u\in\real\times\cnum$
such that $\paramap(\lambda_0,u)=a$,
we set
\[
 \psizero_{t,u}(\nbigm):=
 \bigcup_N
 \Ker\Bigl(
 (-\deldel_tt+\eigenmap(\lambda,u))^N
 \Bigr)
 \subset
 \Gr^{\Vzero}_a(\nbigm).
\]
We obtain the decomposition
\[
\Gr^{\Vzero}_{a}(\nbigm)=
\bigoplus_{\substack{u\in\real\times\cnum\\
\paramap(\lambda_0,u)=a}}
\psizero_{t,u}(\nbigm).
\]
If $\nbigx_0^{(\lambda_0)}\cap\nbigx_0^{(\lambda_1)}\neq\emptyset$,
there exist natural isomorphisms
\[
\psizero_{t,u}(\nbigm)_{|\nbigx_0^{(\lambda_0)}\cap\nbigx_0^{(\lambda_1)}}
\simeq
\psi^{(\lambda_1)}_{t,u}(\nbigm)
_{|\nbigx_0^{(\lambda_0)}\cap\nbigx_0^{(\lambda_1)}}.
\]
By gluing $\psizero_{t,u}(\nbigm)$,
we obtain $\nbigr_{X_0(\ast H_0)}$-module $\psi_{t,u}(\nbigm)$.
We set
$\vecdelta=(1,0)\in\real\times\cnum$.
There exist natural morphisms
$\deldel_t:\psi_{t,u}(\nbigm)\to\psi_{t,u+\vecdelta}(\nbigm)$
and
$t:\psi_{t,u}(\nbigm)\to\psi_{t,u-\vecdelta}(\nbigm)$.
We set
$\psitilde_{t,u}(\nbigm):=
\varinjlim_{N} \psi_{t,u-N\vecdelta}(\nbigm)$
for $u\not\in\real\times\cnum\setminus (\seisuu_{\geq 0}\times \{0\})$.

Let $\nbigm_i$ be coherent $\nbigr_{X(\ast H)}$-modules
which are strictly specializable along $t$.
Let $f:\nbigm_1\lrarr\nbigm_2$ be a morphism of
coherent $\nbigr_{X(\ast H)}$-modules.
Then, $f$ is compatible with the $V$-filtrations,
i.e.,
$f(\Vzero_a(\nbigm_1))\subset \Vzero_a(\nbigm_2)$
for any $\lambda_0\in\cnum$ and $a\in\real$.
We say that $f$ is strictly specializable along $t$
if the cokernel of the induced morphisms
$\psi_{t,u}(f):
\psi_{t,u}(\nbigm_1)\to
\psi_{t,u}(\nbigm_2)$
are strict.
The following proposition is well known
(see \cite{mochi2,sabbah2}).
\begin{prop}
If $f$ is strictly specializable along $t$,
then $f$ is strictly compatible with the $V$-filtrations,
i.e.,
$f(\Vzero_a(\nbigm_1))=
\Image(f)_{|\nbigxzero}\cap \Vzero_a(\nbigm_2)$
for any $a\in\real$ and $\lambda_0\in\cnum$.
\hfill\qed
\end{prop}

The following lemma is also well known.
\begin{lem}
\label{lem;22.4.25.100}
Let $\nbigm$ be a strict coherent $\nbigr_X(\ast H)$-module
which is strictly specializable along $t$.
Let $s$ be a section of $\Vzero_{<0}(\nbigm)$.
If there exists $m\in\seisuu_{\geq 0}$
such that $t^ms=0$,
then $s=0$.
\end{lem}
\pf
There exits $a<0$ such that
$s\in \Vzero_{a}(\nbigm)$.
Because $t:\Gr^{\Vzero}_a(\nbigm)\simeq \Gr^{\Vzero}_{a-1}(\nbigm)$
for any $a<0$,
the induced section $[s]$ of
$\Gr^{\Vzero}_{a}(\nbigm)$ is $0$.
Hence, we obtain $s\in \bigcap_{a\in\real}\Vzero_a(\nbigm)$.
We set $s_0=s$.
There exists $s_1\in \Vzero_{<a}(\nbigm)$
such that $ts_1=s_0$.
By the previous argument,
we obtain $s_1\in\bigcap_{a\in\real}\Vzero_a(\nbigm)$.
Inductively, we take $s_j\in\Vzero_{<0}(\nbigm)$ $(j=1,2,\ldots)$
such that $ts_{j+1}=s_{j}$.
We consider
the $V\nbigr_X$-submodules
$\nbign_n$ of $\Vzero_{<0}(\nbigm)$
generated by $s_j$ $(0\leq j\leq n)$.
Because $V\nbigr_{X}$ is N\"{o}ther,
and $V_{0}(\nbigm)$ is coherent,
there exists $n_0$ such that
$\nbign_n=\nbigm_{n_0}$ for any $n\geq n_0$.
There exist $P_j\in V\nbigr_X$ such that
$s_{n_0+1}=\sum_{j=0}^{n_0} P_j\cdot s_j$.
If $m>0$,
we obtain $t^{n_0+m}s_{n_0+1}=t^{m-1}s=0$
because $t^{j+m}s_j=0$.
By an induction, we obtain that $s=0$.
\hfill\qed

\paragraph{The case $X_0\subset H$}
Let $H$ be a hypersurface such that $X_0\subset H$.
We decompose $H=X_0\cup H'$,
where $\dim(H'\cap X_0)<\dim X_0$.
We have $\nbigr_{X(\ast H)}=\nbigr_{X(\ast H')}(\ast t)$.
A coherent $\nbigr_{X(\ast H')}(\ast t)$-module $\nbigm$
is called strictly specializable along $t$ at $\lambda_0$
if $\nbigm_{|\nbigxzero}$ has a filtration
$V_{\bullet}(\nbigm)$
satisfying the above conditions (i), (ii) and (v)
and the following conditions:
\begin{description}
 \item[(iii')]
	    $t\Vzero_a(\nbigm)=\Vzero_{a-1}(\nbigm)$
	    for any $a\in\real$.
 \item[(iv')]$\deldel_tV_a(\nbigm)\subset V_{a+1}(\nbigm)$
	    for any $a\in\real$.
\end{description}
We set $H_0':=X_0\cap H'$.
As in the previous case,
we obtain the $\nbigr_{X_0(\ast H'_0)}$-modules
$\psi_{t,u}(\nbigm)$ for any $u\in\real$,
and the natural morphisms
$\deldel_t:\psi_{t,u}(\nbigm)\lrarr\psi_{t,u+\vecdelta}(\nbigm)$
and
$t:\psi_{t,u}(\nbigm)\simeq\psi_{t,u-\vecdelta}(\nbigm)$.
They are also denoted by $\psitilde_{t,u}(\nbigm)$.

\paragraph{Canonical prolongations}

For simplicity, we assume that $X_0\not\subset H$.
Let $\nbigm$ be a coherent $\nbigr_{X(\ast H)}(\ast t)$-module
which is strictly specializable along $t$.
For any $\lambda_0\in\cnum$,
we set
\[
 \nbigm_{|\nbigxzero}[\ast t]:=
 \nbigr_X\otimes_{V\nbigr_X} \Vzero_0(\nbigm),
 \quad\quad
 \nbigm_{|\nbigxzero}[!t]:=
 \nbigr_X\otimes_{V\nbigr_X} \Vzero_{<0}(\nbigm).
\]
By gluing $\nbigm_{|\nbigxzero}[\star t]$ $(\star=!,\ast)$
for varying $\lambda_0\in\cnum$,
we obtain
$\nbigm[\star t]$ on $\nbigx$.
They are also strictly specializable along $t$.
We note that
$t:\psi_{t,0}(\nbigm[\ast t])\simeq
\psi_{t,-\vecdelta}(\nbigm[\ast t])$
and
$\deldel_t:\psi_{t,-\vecdelta}(\nbigm[!t])\simeq
\psi_{t,0}(\nbigm[!t])$.
There exist the following canonical morphisms:
\[
\begin{CD}
 \nbigm[!t]
 @>{\iota_1}>>
 \nbigm
 @>{\iota_2}>>
 \nbigm[\ast t].
\end{CD}
\]
See \cite[\S3]{Mochizuki-MTM}
for basic properties of $\nbigm[\star t]$.

\paragraph{Strict $S$-decomposability}

Let $\nbigm$ be a strict coherent $\nbigr_{X(\ast H)}$-module
which is strictly specializable along $t$.
There exist the induced morphisms
$\can:\psi_{t,-\vecdelta}(\nbigm)
\to \psi_{t,0}(\nbigm)$
and
$\var:\psi_{t,0}(\nbigm)\to\psi_{t,-\vecdelta}(\nbigm)$
induced by
$-\deldel_t$ and $t$, respectively.
We say that $\nbigm$ is strictly $S$-decomposable
if
$\psi_{t,0}(\nbigm)
 =\Image\can
 \oplus
 \Ker\var$.
If $\nbigm$ is strictly $S$-decomposable,
there exists a decomposition
$\nbigm=\nbigm'\oplus\nbigm''$
where $\nbigm'$ does not have any non-zero submodule nor
any non-zero quotient whose supports are contained in $X_0$,
and the support of $\nbigm''$ is contained in $X_0$.
We have $\nbigm'=\Image(\nbigm[!t]\lrarr\nbigm[\ast t])$,
and $\psi_0(\nbigm')\simeq \Image\can$.
We also have $\psi_0(\nbigm'')=\Ker\var$
and $\nbigm''\simeq \iota_{\dagger}\Ker\var$,
where $\iota:X_0\to X$ denotes the inclusion.

\paragraph{Regularity}

We say that $\nbigm$ is regular along $t$
if each $V_a(\nbigm)$ is
coherent over $\pi^{\ast}\nbigr_{X_0(\ast H_0)}$,
where $\pi:X\to X_0$ denotes the projection.
(See \cite[\S3.1.d]{sabbah2}.)
Note that $V\nbigr_X=\pi^{\ast}\nbigr_{X_0}\langle t\deldel_t\rangle$.

\subsubsection{Beilinson $\nbigr(\ast t)$-modules}

We put
$A:=\nbigo_{\cnum_{\lambda}}[\lambda s,(\lambda s)^{-1}]
=\bigoplus_n\nbigo_{\cnum_{\lambda}}\lambda^ns^n$,
where $s$ is a formal variable.
We set
$A^a:=(\lambda s)^a\nbigo[\lambda s]\subset A$.
For $a\leq b$,
we put
\[
\II_1^{a,b}:=A^{-b+1}/A^{-a+1}
 \simeq
 \bigoplus_{a\leq i<b}
 \nbigo_{\cnum_{\lambda}}\cdot
 (\lambda s)^{-i},
\quad\quad
\II_2^{a,b}:=A^a/A^b
 \simeq
 \bigoplus_{a\leq i<b}
 \nbigo_{\cnum_{\lambda}}\cdot
 (\lambda s)^i.
\]
The multiplication of $-s$ induces
$\nbign_{\II,1}:\lambda\II^{a,b}_1
\lrarr \II^{a,b}_1$
and
$\nbign_{\II,2}:
\II^{a,b}_2\lrarr
\lambda^{-1}\II^{a,b}_2$.
We put
\[
\IItilde_1^{a,b}:=
 \bigoplus_{a\leq i<b}
 \nbigo_{\cnum^2_{\lambda,t}}(\ast t)\cdot
 (\lambda s)^{-i},
\quad\quad
\IItilde_2^{a,b}:=
 \bigoplus_{a\leq i<b}
 \nbigo_{\cnum^2_{\lambda,t}}(\ast t)\cdot
 (\lambda s)^{i},
\]
where we consider the actions of $\nbigr_{\cnum_t}$
given by 
$\deldel_t(\lambda s)^j=(\lambda s)^{j+1}
=-\lambda\nbign_{\II,i}(\lambda s)^j$.

\subsubsection{Beilinson functors}
\label{subsection;22.4.4.20}

Let $X$ and $X_0$ be as in \S\ref{subsection;22.3.28.10}.
Let $H$ be a hypersurface such that $X_0\not\subset H$.
The pull back of $\IItilde_{2}^{a,b}$ by
the projection $X\to\cnum_t$
is also denoted by $\IItilde_2^{a,b}$.

Let $\nbigm$ be a strict coherent $\nbigr_X(\ast H)$-module
which is strictly specializable along $t$.
We note that
$\Pi_t^{a,b}\nbigm:=\nbigm\otimes\IItilde_2^{a,b}$
is also strictly specializable along $t$.
We obtain
$\Pi^{a,b}_{t,\star}(\nbigm):=\Pi^{a,b}(\nbigm)[\star t]$.
We define
\[
 \Pi^{a,b}_{t,\ast !}\nbigm:=
 \varprojlim_{N\to\infty}
 \Cok\Bigl(
 \Pi^{b,N}_{t,!}(\nbigm)
 \lrarr
 \Pi^{a,N}_{t,\ast}(\nbigm)
 \Bigr)
 \simeq
 \varinjlim_{N\to\infty}
 \Ker\Bigl(
 \Pi^{-N,b}_{t,!}\nbigm
 \lrarr
 \Pi^{-N,a}_{t,\ast}\nbigm
 \Bigr).
\]
For $a\in\seisuu$,
we set
\[
 \psi^{(a)}_t(\nbigm):=
 \Pi^{a,a}_{t,\ast!}(\nbigm),
 \quad
 \Xi^{(a)}_t(\nbigm):=
\Pi^{a,a+1}_{t,\ast!}(\nbigm).
\]
There exist natural exact sequences:
\[
\begin{CD}
 0@>>>
 \nbigm[!t]
 @>{\alpha}>>
 \Xi^{(0)}_t(\nbigm)
 @>{\beta}>>
 \psi^{(0)}_t(\nbigm)
 @>>> 0,
\end{CD}
\]
\[
\begin{CD}
 0@>>>
 \psi^{(1)}_t(\nbigm)
 @>{\gamma}>>
 \Xi^{(0)}_t(\nbigm)
 @>{\delta}>>
 \nbigm[\ast t]
 @>>>
 0.
\end{CD}
\]
We define
$\phi^{(0)}_t(\nbigm)$ as the cohomology of
the complex
\[
\begin{CD}
 \nbigm[!t]
 @>{\alpha+\iota_1}>>
 \Xi^{(0)}_t(\nbigm)
 \oplus
 \nbigm
 @>{\delta-\iota_2}>>
 \nbigm[\ast t].
\end{CD}
\]
The morphisms
$\psi^{(1)}_t(\nbigm)\to\Xi_t^{(0)}(\nbigm)\to\psi_t^{(0)}(\nbigm)$
induce the following morphisms:
\[
\begin{CD}
 \psi_t^{(1)}(\nbigm)
 @>>>
 \phi_t^{(0)}(\nbigm)
 @>>>
 \psi_t^{(0)}(\nbigm)
\end{CD}
\]

We can reconstruct $\nbigm$
as the cohomology of the following complex,
as in \cite{beilinson2}:
\[
 \begin{CD}
  \psi_t^{(1)}(\nbigm)
  @>>>
  \Xi_t^{(0)}(\nbigm)
  \oplus
  \phi_t^{(0)}(\nbigm)
  @>>>
  \psi_t^{(0)}(\nbigm).
 \end{CD}
\]

\begin{lem}[\mbox{\cite[Lemma 4.1.9]{Mochizuki-MTM}}]
Let $\iota:X_0\to X$ denote the inclusion.
Then, there exist isomorphisms
$\iota_{\dagger}\psi_{t,0}(\nbigm)
\simeq
\phi^{(0)}_t(\nbigm)$
and 
$\iota_{\dagger}\psi_{t,-\vecdelta}(\nbigm)
 \simeq\psi^{(a)}_t(\nbigm)$.
\hfill\qed
\end{lem}

\begin{rem}
$\phi_t^{(0)}(\nbigm)$
and $\Xi_t^{(0)}(\nbigm)$
are also denoted by $\phi_t(\nbigm)$
and $\Xi_t(\nbigm)$, respectively.
\hfill\qed
\end{rem}

\subsubsection{Strict specializability and localizability
along holomorphic functions}

Let $X$ be a complex manifold with a hypersurface $H$.
Let $\nbigm$ be an $\nbigr_{X(\ast H)}$-module.
Let $f$ be a holomorphic function on $X$.
Let $\iota_f:X\lrarr X\times\cnum_t$ denote
the graph embedding
defined by $\iota_f(Q)=(Q,f(Q))$.
We say that $\nbigm$ is strictly specializable
(resp. strictly $S$-decomposable)
along $f$
if $\iota_{f\dagger}(\nbigm)$ is strictly specializable
(resp. strictly $S$-decomposable) along $t$.
If $\nbigm$ is strictly specializable along $t$,
there exist
$\nbigr_{X(\ast H)}$-modules
$\psi^{(a)}_f(\nbigm)$
and $\phi^{(0)}_f(\nbigm)$
with isomorphisms
\[
\iota_{f\dagger}\bigl(
 \psi^{(a)}_f(\nbigm)
 \bigr)
 \simeq
 \psi^{(a)}_t\bigl(\iota_{f\dagger}\nbigm\bigr),
 \quad\quad
 \iota_{f\dagger}\bigl(
 \phi^{(0)}_f(\nbigm)
 \bigr)
 \simeq
 \phi^{(0)}_t\bigl(\iota_{f\dagger}\nbigm\bigr).
\]
If there exists
an $\nbigr_X(\ast H)$-modules
$\nbigm_{\star}$ $(\star=!,\ast)$ such that
$\iota_{f\dagger}\nbigm_{\star}\simeq
 \iota_{f\dagger}(\nbigm)[\star t]$,
we say that $\nbigm$ is localizable along $f$,
and $\nbigm_{\star}$ is denoted by 
$\nbigm[\star f]$.
We use the notation
$\Pi^{a,b}_{f\star}(\nbigm)$
$\Pi^{a,b}_{f!\ast}(\nbigm)$
and
$\Xi^{(a)}_f(\nbigm)$
in similar meanings.
For example,
if there exists an $\nbigr_{X(\ast H)}$-module $\nbigm'$
with an isomorphism
$\iota_{f\dagger}\nbigm'\simeq \Xi^{(a)}_t(\iota_{f\dagger}\nbigm)$,
then $\nbigm'$ is denoted by $\Xi^{(a)}_f(\nbigm)$.

\subsubsection{Strict specializability and localizability
along effective divisors}
\label{subsection;22.3.30.1}

Let $D$ be an effective divisor.
For any $P\in X$,
there exists a neighbourhood $X_P$ of $P$ in $X$
with a holomorphic function $f_P$ on $X_P$
which is a generator of $\nbigo(-D)_{|X_P}$.
Let $\iota_{f_P}:X\lrarr X\times\cnum_t$ denote
the graph embedding $\iota_{f_P}(Q)=(Q,f_P(Q))$.
We say that $\nbigm$ is strictly specializable
along $D$ at $P$
if $\iota_{f_P\dagger}(\nbigm_{|X_P})$ is strictly specializable along $t$.
We say that $\nbigm$ is strictly specializable
along $D$
if it is strictly specializable along $D$ at any $P\in X$.
We similarly define the strictly $S$-decomposability condition along $D$.
We say that $\nbigm$ is localizable along $D$
if there exists $\nbigm[\star D]$ $(\star=!,\ast)$
such that
$\iota_{f_P\dagger}(\nbigm[\star D]_{|X_P})
\simeq
\iota_{f_P\dagger}(\nbigm_{|X_P})[\star t]$
at any $P\in X$.
If such $\nbigm[\star D]$ exist, they are unique.

\subsubsection{Admissible specializability}

Let $\nbigm$ be a strict coherent $\nbigr_{X(\ast H)}$-module
with a locally finite increasing filtration $L$ by
coherent $\nbigr_{X(\ast H)}$-submodules.
Assume that $\Gr^L(\nbigm)$ is also strict.
Let $D$ be an effective divisor of $X$.
We say that $(\nbigm,L)$ is
filtered strictly specializable along $D$
if the following holds for any $P\in X$.
\begin{itemize}
 \item Let $(X_P,f_P)$ be as above.
       Then, for any $j$,
       $L_j\nbigm$
       are strictly specializable along $D$.
       Moreover,
       the cokernel of the morphisms
       $\psitilde_{f_P,u}(L_j\nbigm)\to \psi_{f_P,u}(\nbigm)$
       and
       $\psi_{f_P,0}(L_j\nbigm)\to \psi_{f_P,0}(\nbigm)$
       are strict.
\end{itemize}
The $\nbigr_{X}$-modules
$\psitilde_{f_P,u}(\nbigm)$
and $\psi_{f_P,0}(\nbigm)$
equipped with the induced filtrations
are denoted by
$\psitilde_{f_P,u}(\nbigm,L)$
and
$\psi_{f_P,0}(\nbigm,L)$.
We have the morphisms
$\nbign:\psitilde_{f_P,u}(\nbigm)\lrarr
\lambda^{-1}\psitilde_{f_P,u}(\nbigm)$
and
$\nbign:\psi_{f_P,0}(\nbigm)\lrarr
\lambda^{-1}\psi_{f_P,0}(\nbigm)$
obtained as the nilpotent part of $-\del_tt$.
We have the morphisms
$\can:\psitilde_{f_P,-\vecdelta}(\nbigm)
\to \psi_{f_P,0}(\nbigm)$
and
$\var:\psi_{f_P,0}(\nbigm)\to
\psitilde_{f_P,-\vecdelta}(\nbigm)$
induced by $-\deldel_t$ and $t$,
respectively.
We say $(\nbigm,L)$ is admissibly specializable
if moreover the following conditions are satisfied
for any $P\in X$.
\begin{itemize}
 \item There exist a relative monodromy filtrations $M(\nbign;L)$
     of $\nbign$ on
     $\psitilde_{f_P,u}(\nbigm,L)$ and
    $\psi_{f_P,0}(\nbigm,L)$.
 \item We have
       $\can\bigl(
       M_k(\nbign;L)\bigr)
       \subset M_{k-1}(\nbign;L)$
       and
       $\var\bigl(
       M_k(\nbign;L)\bigr)
       \subset M_{k-1}(\nbign;L)$
\end{itemize}

\begin{rem}
In {\rm\cite{Mochizuki-MTM}},
we also imposed the conditions for ramified exponential twist.
Because we are interested in only regular mixed twistor $\nbigd$-modules
in this paper,
we omit the consideration for
ramified exponential twist.
\hfill\qed
\end{rem}

\subsection{$\nbigr_{X(\ast H)}$-triples}

\subsubsection{Some sheaves}

We set $\vecS=\bigl\{\lambda\in\cnum\,\big|\,|\lambda|=1\bigr\}$.
We set $n:=\dim X$.
Let $\pi:\vecS\times X\lrarr X$ denote the projection.
Let $\nbigc^{\infty}_{\vecS\times X}$
denote the sheaf of $C^{\infty}$-functions on $\vecS\times X$.
We recall the sheaves
$\distribution_{\vecS\times X/\vecS}$,
$\distribution^{\moderate H}_{\vecS\times X/\vecS}$
and
$\nbigc^{\infty\moderate H}_{\vecS\times X}$
on $\vecS\times X$
by following \cite{sabbah2} and \cite[\S2.1.3]{Mochizuki-MTM}.

For any open subset $V$ of $\vecS\times X$,
let $\nbige^{(n,n)}_{\vecS\times X/\vecS,c}(V)$
denote the space of
$C^{\infty}$-sections of $\pi^{-1}\Omega_X^{n,n}$ on $V$
whose supports are compact.
Let $\Diff_{\vecS\times X/\vecS}(V)$
denote the space of
linear $C^{\infty}$ differential operators on $V$
relative to $\vecS$,
where we consider only differentials in the $X$-direction.
For any compact subset $K\subset V$
and $P\in \Diff_{\vecS\times X/\vecS}(V)$,
we obtain the semi-norm
$\|\varphi\|_{P,K}:=\sup_K|P\varphi|$.
Let $\nbige^{<H\,(n,n)}_{\vecS\times X/\vecS,c}(V)$
denote the subspace of
$\nbige^{(n,n)}_{\vecS\times X/\vecS,c}(V)$
which consists of $\varphi$
such that
$(P\varphi)_{|(\vecS\times H)\cap V}=0$
for any $P\in\Diff_{\vecS\times X/\vecS}(V)$.
By the semi-norms,
we regard
$\nbige^{(n,n)}_{\vecS\times X/\vecS,c}(V)$
and
$\nbige^{<H\,(n,n)}_{\vecS\times X/\vecS,c}(V)$
as topological vector spaces.
Let $C^0(\vecS)$ denote the space of continuous functions
on $\vecS$,
equipped with the sup norm.
Let $\distribution_{\vecS\times X/\vecS}(V)$
denote the space of continuous
$C^{\infty}(\vecS)$-linear maps
$\nbige^{(n,n)}_{\vecS\times X/\vecS,c}(V)
\to C^0(\vecS)$.
Let $\distribution_{\vecS\times X/\vecS}^{\moderate H}(V)$
denote the space of continuous 
$C^{\infty}(\vecS)$-linear maps
$\nbige^{<H\,(n,n)}_{\vecS\times X/\vecS,c}(V)
\to C^0(\vecS)$.
Thus, we obtain the sheaves
$\distribution_{\vecS\times X/\vecS}$
and
$\distribution^{\moderate H}_{\vecS\times X/\vecS}$.

For any open subset $V$ of $\vecS\times X$,
a $C^{\infty}$-function $\varphi$ on $V\setminus (\vecS\times H)$
is called of moderate growth along $H$
if the following condition is satisfied
for any $(\lambda,Q)\in V\cap(\vecS\times H)$:
\begin{itemize}
 \item Let $g$ be a holomorphic function defined on a neighbourhood $X_Q$
       of $Q$ such that $H\cap X_Q=g^{-1}(0)$.
       Then, for any $P\in \Diff_{\vecS\times X/\vecS}(V)$,
       there exists $N>0$ such that
       $|P\varphi|=O(|g|^{-N})$ around $(\lambda,Q)$.       
\end{itemize}
Let $\nbigc^{\infty\,\moderate H}_{\vecS\times X/\vecS}(V)$
denote the space of
$C^{\infty}$-functions on $V\setminus (\vecS\times H)$
of moderate growth along $H$.
Thus, we obtain the sheaf $\nbigc^{\infty\moderate}_{\vecS\times X/\vecS}$.

Note that $\distribution^{\moderate H}_{\vecS\times X/\vecS}$
and $\nbigc^{\infty}_{\vecS\times X/\vecS}$
are naturally $\nbigr_{X(\ast H)|\vecS\times X}$-modules.
Let $\sigma:\vecS\times X\to\vecS\times X$ be given by
$\sigma(\lambda,P)=(-\lambda,P)$.
Then, $\distribution^{\moderate H}_{\vecS\times X/\vecS}$
and $\nbigc^{\infty}_{\vecS\times X/\vecS}$
are also $\sigma^{-1}\nbigr_{X(\ast H)|\vecS\times X}$-modules
by $\sigma^{-1}(P)\bullet m=\overline{\sigma^{-1}(P)}m$.
The actions of
$\nbigr_{X(\ast H)|\vecS\times X}$
and
$\sigma^{-1}\nbigr_{X(\ast H)|\vecS\times X}$
are commutative.
Thus, we can regard
$\distribution^{\moderate H}_{\vecS\times X/\vecS}$
and $\nbigc^{\infty}_{\vecS\times X/\vecS}$
as
$\nbigr_{X(\ast H)|\vecS\times X}
  \otimes_{\cnum}\sigma^{-1}\nbigr_{X(\ast H)|\vecS\times X}$-modules.
  
\subsubsection{$\nbigr_{X(\ast H)}$-triples}
\label{subsection;22.3.25.100}

A sesqui-linear pairing of
$\nbigr_{X(\ast H)}$-modules
$\nbigm'$ and $\nbigm''$ is 
an
$\nbigr_{X(\ast H)|\vecS\times X}\otimes_{\cnum}
\sigma^{-1}\nbigr_{X(\ast H)|\vecS\times X}$-homomorphism
\[
 C:\nbigm'_{|\vecS\times X}
 \otimes_{\cnum}
 \sigma^{-1}\nbigm''_{|\vecS\times X}
 \to
 \distribution^{\moderate H}_{\vecS\times X/\vecS}.
\]
Such a tuple
$(\nbigm',\nbigm'',C)$ is called
an $\nbigr_{X(\ast H)}$-triple.
An $\nbigr_{X(\ast H)}$-triple is called
coherent (resp. good)
if $\nbigm'$ and $\nbigm''$
are coherent (resp. good).
A morphism of $\nbigr_{X(\ast H)}$-triples
$\varphi:(\nbigm'_1,\nbigm''_1,C_1)
\to
(\nbigm'_2,\nbigm''_2,C_2)$
is defined to be a pair of $\nbigr_{X(\ast H)}$-modules
$\varphi':\nbigm_2'\to\nbigm_1'$
and
$\varphi'':\nbigm_1''\to\nbigm_2''$
such that
\[
 C_1\bigl(\varphi'(u'_2),\sigma^{-1}(u_1'')\bigr)
=C_2\bigl(u_2',\sigma^{-1}\varphi''(u_1)\bigr).
\]
The category of $\nbigr_{X(\ast H)}$-triples
is denoted by $\rtriplecat(X,H)$,
which is abelian.
A subobject $\nbigt_1=(\nbigm_1',\nbigm_1'',C_1)$
of $\nbigt=(\nbigm',\nbigm'',C)$ in $\rtriplecat(X,H)$
consists of a quotient $\nbigr_{X(\ast H)}$-module $\nbigm_1'$ of $\nbigm'$,
an $\nbigr_{X(\ast H)}$-submodule $\nbigm''_1$ of $\nbigm''$
such that
the projection $\nbigm'\lrarr\nbigm'_1$
and the inclusion $\nbigm''_1\lrarr\nbigm''$
induces $\nbigt_1\lrarr\nbigt$.

For $\nbigt=(\nbigm',\nbigm'',C)\in\rtriplecat(X,H)$,
we obtain the Hermitian adjoint
$\nbigt^{\ast}=(\nbigm'',\nbigm',C^{\ast})\in\rtriplecat(X,H)$,
where
$C^{\ast}$ is defined by
$C^{\ast}(u''\otimes \sigma^{-1}(u'))
=\overline{\sigma^{-1}(C(u',\sigma^{\ast}u''))}$.

We set
$\nbigu(p,q)_X:=
\bigl(
 \lambda^p\nbigo_{\nbigx},
 \lambda^q\nbigo_{\nbigx},C_0
\bigr)$,
where $C_0$ is the sesqui-linear pairing
induced by the natural multiplication.
We set
$\newTate(n)_X:=
\nbigu(-n,n)_X=
\bigl(
\lambda^{-n}\nbigo_{\nbigx},
\lambda^n\nbigo_{\nbigx},
C_0
\bigr)$,
called the $n$-th Tate triple.
We use the isomorphisms
$\iota_{\newTate(n)}:\newTate(n)_X^{\ast}\simeq\newTate(-n)_X$
given by $((-1)^n,(-1)^n)$,
and
$\nbigu(p,q)_X^{\ast}\simeq \nbigu(q,p)_X$
given by $((-1)^p,(-1)^p)$.
We shall often omit the subscript $X$.

For an $\nbigr_{X(\ast H)}$-triple
$\nbigt=(\nbigm',\nbigm'',C)$,
we obtain
$\nbigt\otimes\newTate(n)=
(\lambda^{-n}\nbigm',\lambda^n\nbigm'',C\otimes C_0)$,
called the $n$-th Tate twist.
A sesqui-linear duality of $\nbigt\in\rtriplecat(X,H)$ of weight $w$
is a morphism
$\nbigs:\nbigt\lrarr\nbigt^{\ast}\otimes\newTate(-w)$.
It is called Hermitian
if $\nbigs^{\ast}\circ\iota_{\newTate(-w)}=(-1)^w\nbigs$.

Let $f:X'\lrarr X$ be a morphism of complex manifolds.
We set $H'=f^{-1}(H)$.
Let $\nbigt=(\nbigm',\nbigm'',C)$ be
a good $\nbigr_{X'(\ast H')}$-triple
whose support is proper over $X$.
Then, we obtain the induced
good $\nbigr_{X(\ast H)}$-triple
$f_{\dagger}^j(\nbigt)
=(f_{\dagger}^{-j}(\nbigm'),f_{\dagger}^j(\nbigm''),
f_{\dagger}^j(C))$
as defined in \cite{sabbah2}.
See \cite[\S2.1.4]{Mochizuki-MTM} for the formula.
For $c\in H^2_{\DR}(X)$,
we obtain
$L_c:f_{\dagger}^i(\nbigt)\to
f_{\dagger}^{i+2}(\nbigt)\otimes\newTate(1)$.

\begin{rem}
In {\rm\cite{Mochizuki-MTM}},
we adopt a different but equivalent signature rule
for the construction of
the sesqui-linear pairing $f_{\dagger}^j(C)$.
\hfill\qed
\end{rem}

\subsubsection{Filtered $\nbigr_{X(\ast H)}$-triples}

An increasing filtration $W_{\bullet}(\nbigt)$ of $\nbigt$
is a tuple of subobjects $W_j(\nbigt)$ $(j\in\seisuu)$
such that $W_{j-1}(\nbigt)$ is a subobject of $W_{j}\nbigt$
for any $j$.
If $\nbigt=(\nbigm',\nbigm'',C)$
is equipped with an increasing filtration,
$\nbigt^{\ast}$ is equipped with
the induced increasing filtration
$W_j(\nbigt^{\ast})
=\Ker\bigl(
\nbigt^{\ast}\to W_{-j-1}(\nbigt)^{\ast}
\bigr)$
for which
$\Gr^W_j(\nbigt^{\ast})=\Gr^W_{-j}(\nbigt)^{\ast}$.
By using the descriptions
$W_j(\nbigt)=(\nbigm'_j,\nbigm''_j,C_j)$,
the $\nbigr_{X(\ast H)}$-module $\nbigm''$
is equipped with the increasing filtration
$W_j(\nbigm'')=\nbigm_j''$,
and $\nbigm'$ is equipped with induced filtration
$W_j(\nbigm'):=\Ker(\nbigm'\lrarr \nbigm'_{-j-1})$.

\subsubsection{Deformation associated with nilpotent morphisms}

Let $\nbigt=(\nbigm',\nbigm'',C)$
be an $\nbigr_{X(\ast H)}$-triple.
Let $\Lambda$ be a finite set.
Let
$\vecnbign=(\nbign_i\,|\,i\in\Lambda)$
be a commuting tuple of morphisms
$\nbign_i=(\nbign'_i,\nbign_i''):\nbigt\to\nbigt\otimes\newTate(-1)$
$(i\in\Lambda)$.
We assume that $\nbign_i$ are locally nilpotent.
For any $\Lambda_1\subset\Lambda$,
we set $Y(\Lambda_1):=X\times\cnum^{\Lambda_1}$.
We set $\Htilde(\Lambda_1):=
(H\times\cnum^{\Lambda_1})\cup
\bigcup_{i\in\Lambda_1}(X\times\cnum^{\Lambda_1\setminus\{i\}})$.
Let $\vecnbign_{\Lambda_1}:=(\nbign_i\,|\,i\in\Lambda_1)$.
Recall that we obtain
the $\nbigr_{Y(\Lambda_1)(\ast \Htilde(\Lambda_1))}$-triple
$\TNIL(\nbigt,\vecnbign_{\Lambda_1})$
equipped with a tuple of the induced commuting tuple
$\vecnbign=(\nbign_i)$ of morphisms
$\nbign_i:\TNIL(\nbigt,\vecnbign_{\Lambda_1})
\to \TNIL(\nbigt,\vecnbign_{\Lambda_1})
\otimes\newTate(-1)$.
We consider the action of
the $\nbigr_{Y(\Lambda_1)(\ast \Htilde(\Lambda_1))}$
on 
\[
\nbigmtilde'':=
\nbigo_{\nbigy(\Lambda_1)}
(\ast\nbightilde(\Lambda_1))
\otimes
 \nbigm''
\]
given by the natural action of $\nbigr_{X(\ast H)}$
and the following action of $z_i\deldel_{z_i}=\lambda\del_{z_i}$:
\[
 z_i\deldel_{z_i}(m'')
=-\lambda\nbign''_im''.
\]
Here, $m''$ denotes a local section of $\nbigm''$.
We consider a similar action of
$\nbigr_{Y(\Lambda_1)}(\ast \Htilde(\Lambda_1))$
on
\[
 \nbigmtilde':=
 \nbigo_{\nbigy(\Lambda_1)}
(\ast\nbightilde(\Lambda_1))
 \otimes
 \nbigm'.
\]
There exists the sesqui-linear pairing
$\Ctilde$ of $\nbigm'$ and $\nbigm''$
determined by
\[
 \Ctilde(m',\sigma^{-1}m'')
 =C\Bigl(
 \exp\Bigl(
 \sum_{i\in \Lambda_1}-\nbign_i'\log|z_i|^2
 \Bigr)m',
 \sigma^{-1}(m'')
 \Bigr).
\]
Then,
$\TNIL(\nbigt,\vecnbign_{\Lambda_1})
=(\nbigmtilde',\nbigmtilde'',\Ctilde)$.
It is equipped with the induced commuting tuple of
morphisms
$\nbign_i=(\nbign_i',\nbign_i''):
\TNIL(\nbigt,\vecnbign_{\Lambda_1})
\to
\TNIL(\nbigt,\vecnbign_{\Lambda_1})
\otimes\newTate(-1)$.

We set $(\nbigt,\vecnbign)^{\ast}:=(\nbigt^{\ast},-\vecnbign^{\ast})$.
There exists a natural isomorphism
$\TNIL(\nbigt,\vecnbign)^{\ast}
\simeq
\TNIL\bigl((\nbigt,\vecnbign)^{\ast}\bigr)$.

\subsubsection{Beilinson triple}

Let $\II^{a,b}=(\II^{a,b}_1,\II^{a,b}_2,C_0)$,
where $C_0$ is induced by $C_0(s^{-i},s^j)=1$ $(i=j)$
and $C_0(s^{-i},s^j)=0$ $(i\neq j)$.
The multiplication of $-s$ induces
$\nbign_{\II}:\II^{a,b}\to \II^{a,b}\otimes\newTate(-1)$.
There exists the isomorphism
$\nbigs^{a,b}:(\II^{a,b},\vecnbign_{\II})^{\ast}
\simeq (\II^{-b+1,-a+1},\vecnbign)$
given by $s^i\mapsto s^i$.

We obtain
the $\nbigr_{\cnum_t(\ast t)}$-triple
$\IItilde^{a,b}=\TNIL(\II^{a,b},\nbign_{\II})$.
We have
$(\IItilde^{a,b})^{\ast}
\simeq
\IItilde^{-b+1,-a+1}$.

\subsubsection{Strict specializability, canonical prolongations,
and Beilinson functors}

Let $X$ and $X_0$ be as in \S\ref{subsection;22.3.28.10}.
Let $H$ be a hypersurface such that $X_0\not\subset H$.
Let $\nbigt=(\nbigm',\nbigm'',C)$ be
a coherent $\nbigr_{X(\ast H)}$-triple.
We say that $\nbigt$ is strictly specializable along $t$
if $\nbigm'$ and $\nbigm''$ are strictly specializable.

If $\nbigt$ is strictly specializable along $t$,
we obtain
$\nbigt[\ast t]=(\nbigm'[!t],\nbigm''[\ast t],C[\ast t])$
and 
$\nbigt[!t]=(\nbigm'[\ast t],\nbigm''[!t],C[!t])$
as explained in
\cite[\S3.2]{Mochizuki-MTM}.
There exist the natural morphisms
\[
 \begin{CD}
  \nbigt[!t]@>{\iota_1}>>
  \nbigt @>{\iota_2}>>
  \nbigt[\ast t].
 \end{CD}
\]
We obtain the $\nbigr_{X(\ast H)}(\ast t)$-triple
$\Pi^{a,b}_t\nbigt
=\nbigt\otimes\IItilde^{a,b}$
which is also strictly specializable along $t$.
Hence, we obtain
$\Pi^{a,b}_{t,\star}(\nbigt)=
 \Pi^{a,b}_{t,\star}(\nbigt)[\star t]$.
We set
\[
\Pi^{a,b}_{t,!\ast}(\nbigt)
=\varprojlim_{N\to \infty}
\Cok\bigl(
 \Pi^{b,N}_{t!}\nbigt
 \to
 \Pi^{a,N}_{t\ast}\nbigt
 \bigr)
\simeq
 \varinjlim_{N\to \infty}
\Ker\bigl(
 \Pi^{-N,b}_{t!}\nbigt
 \to
 \Pi^{-N,a}_{t\ast}\nbigt
 \bigr).
\]
We obtain
\[
 \psi^{(a)}_t(\nbigt):=
 \Pi^{a,a}_{t,!\ast}(\nbigt),
 \quad
 \Xi^{(a)}_t(\nbigt):=
 \Pi^{a,a+1}_{t,!\ast}(\nbigt).
\]
There exist the isomorphisms
$\II^{i,i+1}\simeq \bigoplus \newTate(i)$
induced by $s^i\longmapsto 1$.
It induces
$\psi_t^{(a)}(\nbigt)\simeq \psi_t^{(a')}(\nbigt)\otimes\newTate(a-a')$
and
$\Xi_t^{(a)}(\nbigt)\simeq \Xi_t^{(a')}(\nbigt)\otimes\newTate(a-a')$.
There exist the following natural exact sequences
\[
\begin{CD}
 0@>>>
 \psi^{(1)}_t(\nbigt)
 @>{\alpha}>>
 \Xi^{(0)}_t(\nbigt)
 @>{\beta}>>
 \nbigm[\ast t]
 @>>>0,
\end{CD}
\]
\[
 \begin{CD}
 0@>>>
 \nbigm[!t]  
 @>{\gamma}>>
 \Xi^{(0)}_t(\nbigt)
 @>{\delta}>>
 \psi^{(0)}_t(\nbigt)
 @>>>0.
\end{CD}
\]
We define
$\phi_t^{(0)}(\nbigt)$
as the cohomology of the following complex:
\[
\begin{CD}
 \nbigt[!t]
 @>{\alpha+\iota_1}>>
 \Xi_t^{(0)}(\nbigt)
 \oplus
 \nbigt
 @>{\delta-\iota_2}>>
 \nbigt[\ast t].
\end{CD}
\]
We obtain the following morphisms
from $\beta$ and $\gamma$:
\[
\begin{CD}
 \psi^{(1)}_t(\nbigt)
 @>{\can}>>
 \phi^{(0)}_t(\nbigt)
 @>{\var}>>
 \psi^{(0)}_t(\nbigt).
\end{CD}
\]
As in \cite{beilinson2},
we reconstruct $\nbigt$
as the cohomology of the following complex:
\[
 \begin{CD}
  \psi_t^{(1)}(\nbigt)
  @>{\gamma+\can}>>
  \Xi_t^{(0)}(\nbigt)
  \oplus
  \phi_t^{(0)}(\nbigt)
  @>{\delta-\var}>>
  \psi_t^{(0)}(\nbigt).
 \end{CD}
\]
There exist natural isomorphisms
$\psi^{(1)}(\nbigt)^{\ast}\simeq \psi^{(0)}(\nbigt^{\ast})$,
$\Xi^{(0)}(\nbigt)^{\ast}
\simeq
\Xi^{(0)}(\nbigt)^{\ast}$,
and
$\phi^{(0)}(\nbigt)^{\ast}\simeq \phi^{(0)}(\nbigt^{\ast})$.

Recall that for each
$u\in(\real\times\cnum)\setminus (\seisuu_{\geq 0}\times \{0\})$,
we obtain the $\nbigr_{X_0}$-triple
\[
 \psitilde_{t,u}(\nbigt)
 =\bigl(
 \psitilde_{t,u}(\nbigm'),
 \psitilde_{t,u}(\nbigm''),
 \psitilde_{t,u}(C)
 \bigr).
\]
See \cite{sabbah2} or \cite[\S14.4.6]{mochi2}
for the construction of $\psitilde_{t,u}(C)$.
We have the morphism
$\nbign=(\nbign',\nbign''):\psitilde_{t,u}(\nbigt)
\to\psitilde_{t,u}(\nbigt)\otimes\newTate(-1)$,
where $\nbign'$ and $\nbign''$
are obtained as the nilpotent part of
$-t\del_t$.
In particular, we obtain
$(\psitilde_{t,-\vecdelta}(\nbigt),\nbign)$.
Let $\iota:X_0\lrarr X$ denote the inclusion.

\begin{prop}[\mbox{\cite[Proposition 4.3.1]{Mochizuki-MTM}}]
\label{prop;22.4.25.10}
There exist the natural isomorphism
\[
 \iota_{\dagger}
 \psi_{t,-\vecdelta}(\nbigt)
 \simeq
 \psi^{(1)}_t(\nbigt)\otimes\nbigu(1,0)
\]
with the following property:
\begin{itemize}
 \item Let
       $\iota_{\dagger}\psi_{t,-\vecdelta}(\nbigt)
       \simeq
       \psi^{(0)}_t(\nbigt)\otimes\nbigu(0,1)$
       be the induced morphism.
       Then, the following diagram is commutative:
\[
 \begin{CD}
  \psi^{(1)}_t(\nbigt)
  @>{\delta\circ\gamma}>>
  \psi^{(0)}_t(\nbigt)
  \\
  @V{\simeq}VV @V{\simeq}VV \\
  \iota_{\dagger}\psi_{t,-\vecdelta}(\nbigt)
  \otimes\nbigu(-1,0)
  @>{\nbign}>>
 \iota_{\dagger}\psi_{t,-\vecdelta}(\nbigt)
  \otimes\nbigu(0,-1).
 \end{CD}
\]       
\end{itemize}
(See {\rm\cite{Mochizuki-MTM}} for the construction of the isomorphism.)
\hfill\qed
\end{prop}

Let $\nbigs:\nbigt\to\nbigt^{\ast}\otimes\newTate(-w)$
be a Hermitian sesqui-linear duality of weight $w$.
We also obtain the following induced
Hermitian sesqui-linear duality of weight $w$:
\[
 \phi_t^{(0)}(\nbigs):
 \phi_t^{(0)}(\nbigt)\lrarr
 \phi_t^{(0)}(\nbigt^{\ast})\otimes\newTate(-w)\simeq
 \phi_t^{(0)}(\nbigt)^{\ast}\otimes\newTate(-w),
\]
\[
 \Xi_t^{(0)}(\nbigs):
 \Xi_t^{(0)}(\nbigt)\lrarr
 \Xi_t^{(0)}(\nbigt^{\ast})\otimes\newTate(-w)\simeq
 \Xi_t^{(0)}(\nbigt)^{\ast}\otimes\newTate(-w).
\]
We obtain the following Hermitian sesqui-linear duality
of weight $w-1$:
\[
 \psi_t^{(1)}(\nbigs):
 \psi_t^{(1)}(\nbigt)
 \lrarr
 \psi_t^{(1)}(\nbigt^{\ast})
 \otimes
 \newTate(-w)
 \simeq
 \psi_t^{(0)}(\nbigt^{\ast})
 \otimes
 \newTate(-w+1)
 \simeq
 \psi_t^{(1)}(\nbigt)^{\ast}
 \otimes
 \newTate(-w+1).
\]
We also obtain the induced morphism
$\psitilde_{t,u}(\nbigs):
\psitilde_{t,u}(\nbigt)
\lrarr
\psitilde_{t,u}(\nbigt)^{\ast}(\nbigt)^{\ast}
\otimes\newTate(-w)$
for $u\in(\real\times\cnum)\setminus(\seisuu_{\geq 0}\times\{0\})$.

\begin{prop}[\mbox{\cite[Proposition 4.3.2]{Mochizuki-MTM}}]
The following diagram is commutative:
\[
 \begin{CD}
  \psi_t^{(1)}(\nbigt)
  @>>>
  \psi_t^{(1)}(\nbigt)^{\ast}
  \otimes\newTate(-w+1)
  \\
  @V{\simeq}VV @V{\simeq}VV \\
  \iota_{\dagger}
  \psitilde_{t,-\vecdelta}(\nbigt)
  \otimes\nbigu(-1,0)
  @>>>
  \Bigl(
  \iota_{\dagger}
  \psitilde_{t,-\vecdelta}(\nbigt)
  \otimes\nbigu(-1,0)
  \Bigr)^{\ast}\otimes\newTate(-w+1).
 \end{CD}
\] 
\hfill\qed
\end{prop}

\begin{rem}
$\phi^{(0)}_t(\nbigt)$
and $\Xi^{(0)}_t(\nbigt)$
are  also denoted by
 $\phi_t(\nbigt)$
and $\Xi_t(\nbigt)$, respectively.
\hfill\qed
\end{rem}

\subsubsection{Strict specializability and localizability
along holomorphic functions}

Let $X$ be a complex manifold with a hypersurface $H$.
Let $g$ be a holomorphic function on $X$.
Let $\iota_g:X\lrarr X\times\cnum_t$ denote the graph embedding.
A coherent $\nbigr_{X(\ast H)}$-triple $\nbigt$ is called
strictly specializable (resp. strictly $S$-decomposable) along $g$
if $\iota_{g\dagger}(\nbigt)$ is strictly specializable
(resp. strictly $S$-decomposable) along $t$.
For $u\in(\real\times\cnum)\setminus(\seisuu_{\geq 0}\times\{0\})$
We obtain the $\nbigr_{X(\ast H)}$-triples
$\psitilde_{g,u}(\nbigt)
=\psitilde_{t,u}(\iota_{g\dagger}\nbigt)$.
There exist the $\nbigr_{X(\ast H)}$-triples
$\psi^{(a)}_g(\nbigt)$ equipped with an isomorphism
$\iota_{g\dagger}\psi_g^{(a)}(\nbigt)
\simeq
 \psi^{(a)}_t(\iota_{g\dagger}\nbigt)$.
There exists the $\nbigr_{X(\ast H)}$-triple
$\phi_g(\nbigt)$ equipped with an isomorphism
$\iota_{g\dagger}\phi_g(\nbigt)
\simeq
 \phi^{(0)}_t(\iota_{g\dagger}\nbigt)$.
We say that $\nbigt$ is localizable along $g$
if the underlying $\nbigr_{X(\ast H)}$-modules are localizable
along $g$.
In that case, there exist $\nbigr_{X(\ast H)}$-triples
$\nbigt[\star g]$ with an isomorphism
$\iota_{f\dagger}\nbigt[\star g]
\simeq
\bigl(
\iota_{f\dagger}(\nbigt)
\bigr)[\star t]$.
If there exists an $\nbigr_{X(\ast H)}$-triple
$\nbigt'$
with an isomorphism
$\iota_{g\dagger}\nbigt'\simeq
\Xi^{(a)}_t(\iota_{g\dagger}\nbigt)$,
$\nbigt'$ is denoted by
$\Xi^{(a)}_g(\nbigt)$.
We use the notation
$\Pi^{a,b}_{g\star}(\nbigt)$
and
$\Pi^{a,b}_{g!\ast}(\nbigt)$
in similar meanings.

\subsubsection{Strict specializability and localizability
along holomorphic functions}

Let $D$ be an effective divisor.
We define the strictly specializability condition
and the strictly $S$-decomposability condition along $D$
for coherent $\nbigr_{X(\ast H)}$-triples
as in the case of coherent $\nbigr_{X(\ast H)}$-modules.
We also define the localizability condition for
coherent $\nbigr_{X(\ast H)}$-triples
as in the case of coherent $\nbigr_{X(\ast H)}$-modules.
We use the notation $\nbigt[\star D]$
as in the case of $\nbigm[\star D]$.
(See \S\ref{subsection;22.3.30.1}.)

\subsubsection{Admissible specializability}

Let $(\nbigt,L)$ be a filtered $\nbigr_{X(\ast H)}$-triple.
Let $D$ be an effective divisor.
We say that $(\nbigt,L)$ is admissibly
(resp. filtered strictly) specializable
along $D$
if the underlying filtered $\nbigr_{X(\ast H)}$-modules
are admissibly (filtered strictly) specializable along $D$.

\subsubsection{Strict $S$-decomposable $\nbigr_{X}$-triples}

Let $\nbigt$ be a coherent $\nbigr_X$-triple.
We say that it is strictly $S$-decomposable
if the following holds for any 
open subset $U\subset X$ with a holomorphic function $f_U\in\nbigo(U)$.
\begin{itemize}
 \item $\nbigt_{|U}$ is strictly $S$-decomposable along $f_U$.
\end{itemize}
If $\nbigt$ is strictly $S$-decomposable,
there exist a locally finite family
$\{Z_i\,|\,i\in\Gamma\}$ of closed irreducible
complex analytic subsets $Z_i$ of $X$
and decomposition
\[
 \nbigt=\bigoplus_{i\in\Gamma}
 \nbigt_{Z_i},
\]
such that the following holds for each $\nbigt_{Z_i}$.
\begin{itemize}
 \item 
      The support of $\nbigt_{Z_i}$ is $Z_i$.
       Moreover, $\nbigt_{Z_i}$ has neither
       non-zero subobject nor quotient
       whose support is strictly smaller than $Z_i$.
\end{itemize}
The decomposition is called
the strict support decomposition.
See \cite{sabbah2} for more details.

We recall the following lemma.
\begin{lem}
\label{lem;22.4.21.1}
Let $\nbigt$ be strictly $S$-decomposable
coherent $\nbigr_X$-triple.
If the support of $\nbigt$ is contained in
a closed complex submanifold $Z\subset X$,
there exists a strictly $S$-decomposable coherent
$\nbigr_Z$-triple $\nbigt_Z$
with an isomorphism
$\nbigt\simeq \iota_{Z\dagger}\nbigt_Z$,
where $\iota_Z:Z\to X$ denote the inclusion.
Such $\nbigt_Z$ is unique up to a canonical isomorphism.
\end{lem}
\pf
By the uniqueness, it is enough to study
the case where $X$ is a neighbourhood of
$Z\times\{(0,\ldots,0)\}$
in $Z\times\cnum^{\ell}$.
By an induction, we have only to study the case $\ell=1$.
Let $\iota_Z:Z=Z\times\{0\}\to X$ denote the inclusion,
and let $\pi:X\to Z$ denote the projection.
Because $\nbigt$ is strictly $S$-decomposable along
the projection $X\to \cnum$,
there exists a coherent $\nbigr_Z$-triple $\nbigt_Z$
with an isomorphism
$\iota_{Z\dagger}\nbigt_Z\simeq\nbigt$.
Let $\nbigm'_Z,\nbigm''_Z$ be the $\nbigr_Z$-modules
underlying $\nbigt_Z$,
i.e.,
$\nbigt_Z=(\nbigm'_Z,\nbigm''_Z,C_Z)$
for a sesqui-linear pairing $C_Z$ of
$\nbigm'_Z$ and $\nbigm''_Z$.
Let $\nbigm_Z$ denote $\nbigm'_Z$ or $\nbigm''_Z$.
It is enough to prove that $\nbigm_Z$ is
strictly $S$-decomposable under the assumption
that $\nbigm:=\iota_{Z\dagger}\nbigm_Z$ is
strictly $S$-decomposable.
As remarked in \cite[Lemma 4.4.9]{Mochizuki-MTM},
for any open subset $U\subset Z$
with a holomorphic function $f\in\nbigo_Z(U)$,
$(\nbigm_Z)_{|U}$ is strictly specializable along $f$.
We set $\Utilde=\pi^{-1}(U)$
and $\ftilde:=f\circ\pi$ on $\Utilde$.
Let $\iota_{Z|U}:U=U\times\{0\}\to \Utilde$ denote the inclusion.
We have
$(\iota_{Z|})_{\dagger}\psi_{f,u}(\nbigm_{Z|U})
=\psi_{\ftilde,u}(\nbigm_{|\Utilde})$
for any $u\in\real\times\cnum$.
Then, because $\nbigm_{|\Utilde}$
is strictly $S$-decomposable along $\ftilde$,
we obtain that
$\nbigm_{Z|U}$ is strictly $S$-decomposable along $f$.
\hfill\qed

\subsection{Regular polarizable pure twistor $\nbigd$-modules
and mixed twistor $\nbigd$-modules}
\label{subsection;22.4.25.30}

\subsubsection{Polarizable pure twistor $\nbigd$-modules}
Let $X$ be any complex manifold.
We recall the definition of regular
pure twistor $\nbigd$-modules
by following \cite{sabbah2},
which is given in an inductive way on the dimension
of the strict supports.

\begin{df}
The category $\MT_{\leq d,\reg}(X,w)$
is the full subcategory of $\rtriplecat(X)$
whose objects are triples $\nbigt$
satisfying the following conditions.
\begin{description}
 \item[(HSD)] $\nbigt$ is strict, holonomic
	    and strictly $S$-decomposable,
	    and the dimension of the support is less than $d$.
	    Let $\nbigt=\bigoplus\nbigt_Z$
	    denote the strict support decomposition.
	    Note that $\dim Z\leq d$.
 \item[($\boldsymbol{MT_{0}}$)]
	    If $\dim Z=0$,
	    $\nbigt_Z$
	    is the push-forward of
	    pure twistor structure of weight $w$
	    by the inclusion $Z\to X$.
 \item[(REG, $\boldsymbol{MT_{>0}}$)]
	    If $\dim Z>0$,
	    for any open subset $U\subset X$
	    with a non-zero holomorphic function $f_U\in\nbigo(U)$
  	    such that $Z\cap U\not\subset f_U^{-1}(0)$,
	    $\nbigt_{Z|U}$ is regular along $f_U$.
	    Moreover,
	    $\Gr^W_j\psitilde_{f_U,u}(\nbigt_{Z|U})\in
	    \MT_{\leq \dim Z-1,\reg}(U,w+j)$
	    $(u\in(\real\times\cnum)\setminus(\seisuu_{\geq 0}\times\{0\}))$
	    and
	    $\Gr^W_j\phi_{f_U}(\nbigt_{|U})
	    \in \MT_{\leq \dim Z-1,\reg}(U,w+j)$.

\end{description}
An object of $\MT_{\reg}(X,w):=\MT_{\leq\dim X,\reg}(X,w)$ is
called a regular pure twistor $\nbigd_X$-module of weight $w$.
\hfill\qed
\end{df}

\begin{df}
Let $\nbigt\in\MT_{\leq d,\reg}(X,w)$.
Let $\nbigs$ be 
a Hermitian sesqui-linear duality of $\nbigt$
of weight $w$.
Note that
 there exists the strict support
 decomposition $(\nbigt,\nbigs)=\bigoplus(\nbigt_Z,\nbigs_Z)$. 
We say that $\nbigs$
is a polarization of $\nbigt$
if the following condition is satisfied. 
\begin{description}
 \item[(Case $\dim Z=0$)]
	    $(\nbigt_Z,\nbigs_Z)$ is the push-forward of
	    a polarized pure twistor structure of weight $w$
	    by the inclusion
	    $Z\to X$.	    
 \item[(Case $\dim Z>0$)]
	    Let $U\subset X$ be
	    any open subset 
	    with a non-zero holomorphic function $f_U\in\nbigo(U)$
	    such that
	    $U\cap Z\not\subset f_U^{-1}(0)$.
	    Then, for any $j\in\seisuu_{\geq 0}$,
	    $(\nbign^{\ast})^j\circ \psitilde_{g,u}(\nbigs)$
	    $(u\in(\real\times\cnum)\setminus(\seisuu_{\geq 0}\times\{0\}))$
	    induce a polarization of
	    $P\Gr^W_j\psitilde_{g,u}(\nbigt_Z)$,
	    and
	    $(\nbign^{\ast})^j\circ
	    \phi_g(\nbigs)$
	    induce a polarization of
	    $P\Gr^W_j\phi_{g}(\nbigt_Z)$.
\end{description}
A regular pure twistor $\nbigd$-module is called polarizable
if it has a polarization. 
 \hfill\qed
\end{df}

We note the following easy lemma.
\begin{lem}
Let $\nbigt$ be an $\nbigr_X$-triple
equipped with a Hermitian sesqui-linear duality
$\nbigs:\nbigt\to\nbigt^{\ast}\otimes\newTate(-w)$
of weight $w$.
Then, 
$(\nbigt,\nbigs)$  is a pure twistor $\nbigd$-module on $X$
if and only if there exists an open covering
 $X=\bigcup_{i\in\Lambda} U_i$
such that each $(\nbigt,\nbigs)_{|U_i}$
is a pure twistor $\nbigd$-module of weight $w$ on $U$.
\end{lem}
\pf
The ``only if'' part if the claim is obvious.
Let us prove the ``if'' part by using
the induction on the dimension of the support.
We may assume that $\nbigt$ has the strict support $Z$.
Let $U$ be any open subset of $X$
equipped with a holomorphic function $g$ on $U$
such that $Z\cap U\not\subset f^{-1}(0)$.
For each $j\geq 0$,
$(\nbign^{\ast})^j\circ \psitilde_{g,u}(\nbigs)_{|U\cap U_i}$
is a polarization of
$P\Gr^W_j\psitilde_{g,u}(\nbigt)_{|U\cap U_i}$ for any $i\in\Lambda$.
By the assumption of the induction,
we obtain that
$(\nbign^{\ast})^j\circ \psitilde_{g,u}(\nbigs)$
is a polarization
of $P\Gr^W_j\psitilde_{g,u}(\nbigt)$.
Similarly, we obtain that
$P\Gr^W_j\phi_{g}(\nbigt,\nbigs)$
are polarized pure twistor $\nbigd$-modules
of weight $w+j$.
\hfill\qed

\subsubsection{Polarized graded Lefschetz twistor $\nbigd$-modules}
\label{subsection;22.4.28.4}

Let $\epsilon\in\{\pm 1\}$.
Let $\nbigt=\bigoplus \nbigt_i$
be a graded $\nbigr_X$-triple
such that each $\nbigt_i$ is
polarizable pure twistor $\nbigd$-module of weight $w-\epsilon i$.
We assume that the set $\{i\,|\,\nbigt_i\neq 0\}$ is finite.
Let $\nbigl:\nbigt_i\to\nbigt_{i-2}\otimes\newTate(\epsilon)$
$(i\in\seisuu)$,
which gives a morphism
$\nbigl:\nbigt\to\nbigt\otimes\newTate(\epsilon)$.
Such $(\nbigt,\nbigl)$ is called a graded Lefschetz twistor
$\nbigd$-module of weight $w$ if
$\nbigl^j$ induces
$\nbigt_j\simeq \nbigt_{-j}\otimes\newTate(\epsilon j)$
for any $j\in\seisuu_{\geq 0}$.
For $j\geq 0$,
$P\nbigt_j:=\Ker(\nbigl^{j+1}:\nbigt_j
\to \nbigt_{-j-2}\otimes\newTate(-j-1))$
is called the primitive part.
Let $\nbigs:\nbigt\to\nbigt^{\ast}\otimes\newTate(-w)$
denote the sesqui-linear duality of weight $w$
such that $\nbigs(\nbigt_j)\subset\nbigt^{\ast}_{-j}\otimes\newTate(-w)$.
It is called a sesqui-linear duality of
$(\nbigt,\nbigl)$
if $\nbigl^{\ast}\circ\nbigs+\nbigs\circ\nbigl=0$.
A Hermitian sesqui-linear duality of weight $w$
of $(\nbigt,\nbigl)$ is called a polarization
if $(\nbigl^{\ast})^j\circ\nbigs$
induces a polarization of $P\nbigt_j$
for any $j\in\seisuu_{\geq 0}$.
Such a tuple $(\nbigt,\nbigl,\nbigs)$ is called
a regular polarized graded Lefschetz twistor $\nbigd$-module
of weight $w$.

\subsubsection{Fundamental properties of
polarized pure twistor $\nbigd$-modules}

Let $\MT_{\reg}^{\pol}(X,w)$
denote the full subcategory of polarizable
regular pure twistor $\nbigd$-module
of weight $w$ on $X$.

\begin{prop}
The category $\MT_{\reg}^{\pol}(X,w)$
is abelian and semisimple.
\hfill\qed
\end{prop}

Let $f:X\lrarr Y$ be a projective morphism.
Let $\nbigt\in\MT_{\reg}^{\pol}(X,w)$
with a polarization $\nbigs$.
Let $c$ denote the first Chern class
of a line bundle on $X$
which is relatively ample with respect to $f$.
We obtain
$L_c:f_{\dagger}^i\nbigt\to f^{i+2}_{\dagger}\nbigt\otimes\newTate(1)$.
We have the induced Hermitian sesqui-linear duality
$f_{\dagger}(\nbigs)$ of $\bigoplus f_{\dagger}^i(\nbigt)$.

\begin{thm}
\label{thm;22.4.18.20}
$(\bigoplus f^i_{\dagger}(\nbigt),
L_c,f_{\dagger}(\nbigs))$
is a regular polarized graded Lefschetz twistor $\nbigd$-module
of weight $w$.
\hfill\qed
\end{thm}

See \cite{sabbah2} and \cite{mochi2} for more details.
See also \cite{Sabbah-wild} and \cite{Mochizuki-wild}
for polarized wild pure twistor $\nbigd$-modules.

\subsubsection{Mixed twistor $\nbigd$-modules}
\label{subsection;22.4.26.11}

A filtered $\nbigr_X$-triple $(\nbigt,W)$ is called
a mixed twistor $\nbigd$-module
if the following conditions are satisfied.
\begin{itemize}
 \item Each $\Gr^W_j(\nbigt)$ is
       a polarizable pure twistor $\nbigd$-module of weight $j$.
 \item Let $U\subset X$ be an open subset
       with a non-zero holomorphic function $f_U\in\nbigo(U)$.
       Then,
       $(\nbigt,W)_{|U}$ is admissibly specializable
       along $f_U$.
       Moreover, if 
       $U\cap\Supp(\nbigt)\not\subset f_U^{-1}(0)$,
       then
       $(\psitilde_{f_U,u}(\nbigt),W)$
       $(u\in(\real\times\cnum)\setminus(\seisuu_{\geq 0}\times\{0\}))$
       and
       $(\phi_{f_U}(\nbigt),W)$
       are mixed twistor $\nbigd$-modules
       whose supports are strictly smaller than $U\cap\Supp(\nbigt)$.
       Here, $W$ denotes the relative monodromy filtrations
       of $\nbign$ with respect to the naively induced filtrations.
\end{itemize}

Let $\MTM(X)$ denote the category of mixed twistor $\nbigd$-modules
on $X$.
It is an abelian category.
Any mixed twistor $\nbigd$-module $(\nbigt,W)$
is localizable along an effective divisor $D$,
and $\nbigt[\star D]$ $(\star=!,\ast)$
are naturally equipped with a filtration $W$,
and $(\nbigt[\star D],W)$ are mixed twistor $\nbigd$-modules.
Moreover, they depend only on the support of $D$.
See \cite{Mochizuki-MTM} for other fundamental properties
of mixed twistor $\nbigd$-modules.

We note the following lemma.
\begin{lem}
Let $(\nbigt,W)$ be a filtered $\nbigr_X$-triple
such that $\Gr^W_j(\nbigt)$ $(j\in\seisuu)$
are polarizable pure twistor $\nbigd$-module of weight $j$ on $X$.
Then, $(\nbigt,W)$ is a mixed twistor $\nbigd$-module
if and only if  
there exists an open covering
$X=\bigcup_{i\in\Lambda} U_i$
such that each $(\nbigt,W)_{|U_i}$
is a mixed twistor $\nbigd$-module.
\end{lem}
\pf
The ``only if'' part of the claim is obvious.
Let us prove the ``if part'' of the claim.
We use a N\"{o}therian induction on the support.
Let $U$ be an open subset of $X$
equipped with a non-zero holomorphic function $f_U$ on $U$.
Because
$(\nbigt,W)_{|U\cap U_i}$ $(i\in\Lambda)$
are admissibly specializable along $f_{U|U\cap U_i}$,
$(\nbigt,W)_{|U}$ is admissibly specializable along $f_U$.
Let $L$ be the filtration of
$\psitilde_{f_U,u}(\nbigt)$
induced by $W$,
i.e.,
$L_j\psitilde_{f_U,u}(\nbigt_{|U})
=\psitilde_{f_U,u}(W_j\nbigt_{|U})$.
By the canonical decomposition of Kashiwara \cite{k2}
with respect to the relative monodromy filtration
and a generalization in \cite{saito2}
(see \cite[\S6.1.2]{Mochizuki-MTM}),
there exist natural isomorphisms
$\Gr^W_m\psitilde_{f_U,u}(\nbigt_{|U})\simeq
 \bigoplus_k
 \Gr^W_m\Gr^L_k\psitilde_{f_U,u}(\nbigt_{|U})$.
Hence, $\Gr^W_m\psitilde_{f_U,u}(\nbigt_{|U})$
are polarizable pure twistor $\nbigd$-modules of weight $m$.
Then, because 
$(\psitilde_{f_U,u}(\nbigt_{|U\cap U_i}),W)$ $(i\in\Lambda)$
are mixed twistor $\nbigd$-modules,
we obtain that 
$(\psitilde_{f_U,u}(\nbigt_{|U}),W)$
is a mixed twistor $\nbigd$-module,
by the assumption of the induction.
Similarly, we obtain that
$(\phi_{f_U}(\nbigt_{|U}),W)$
are polarizable pure twistor $\nbigd$-module of weight $m$.
\hfill\qed

\subsection{Mixed twistor structures on $(X,H)$}

\subsubsection{Smooth $\nbigr_{X(\ast H)}$-triples}

Let $\nbigc^{\infty\,\lambda{\textrm -}\hol}_{\cnum^{\ast}\times (X\setminus H)}$
denote the sheaf of $C^{\infty}$-functions $\varphi$
on $\cnum^{\ast}\times (X\setminus H)$
such that $\del_{\lambdabar}\varphi=0$.
Let $\sigma:\cnum^{\ast}\times (X\setminus H)
\lrarr\cnum^{\ast}\times (X\setminus H)$
be defined by $\sigma(\lambda,P)=(-\lambdabar^{-1},P)$.
The sheaf 
$\nbigc^{\infty\,\lambda{\textrm -}\hol}_{\cnum^{\ast}\times (X\setminus H)}$
is naturally
an $\nbigr_{X|\cnum^{\ast}\times X}\otimes_{\nbigo_{\cnum^{\ast}}}
 \sigma^{-1}(\nbigr_{X|\cnum^{\ast}\times X})$-module.

A coherent $\nbigr_{X(\ast H)}$-module $\nbigm$ is called smooth
if $\nbigm_{|\nbigx\setminus\nbigh}$
is a locally free $\nbigo_{\nbigx\setminus\nbigh}$-module
of finite rank.
An $\nbigr_{X(\ast H)}$-triple
$(\nbigm',\nbigm'',C)$ is called smooth
if the following conditions are satisfied.
\begin{itemize}
 \item $\nbigm'$ and $\nbigm''$ are smooth $\nbigr_{X(\ast H)}$-modules.
 \item The image of $C$ is contained in
       $\nbigc^{\infty\moderate H}_{\vecS\times X/\vecS}$.
       Moreover,
       the restriction
\[
       C_{|\vecS\times(X\setminus H)}:
       \nbigm'_{|\vecS\times (X\setminus H)}
       \otimes
       \sigma^{-1}(\nbigm''_{|\vecS\times(X\setminus H)})
       \lrarr
       \nbigc^{\infty}_{\vecS\times (X\setminus H)}
\]
       extends to
       an $\nbigr_{X\setminus H|\cnum^{\ast}\times (X\setminus H)}
       \otimes
       \sigma^{-1}(\nbigr_{X\setminus H|
        \cnum^{\ast}\times (X\setminus H)})$-homomorphism
\[
       \nbigm'_{|\cnum^{\ast}\times X}
       \otimes_{\nbigo_{\cnum^{\ast}}}
       \sigma^{-1}(\nbigm''_{|\cnum^{\ast}\times X})
       \lrarr
       \nbigc^{\infty\,\lambda{\textrm-}\hol}_{\cnum^{\ast}\times (X\setminus H)}.
\]       
\end{itemize}
Let $\rtriplecat_{\sm}(X,H)\subset\rtriplecat(X,H)$
denote the full subcategory of smooth $\nbigr_{X(\ast H)}$-triples.

Let $f:Y\lrarr X$ be a morphism of complex manifolds.
We set $H_Y:=f^{-1}(H)$.
For any smooth $\nbigr_{X(\ast H)}$-triple
$\nbigt=(\nbigm',\nbigm'',C)$,
the pull back $f^{\ast}(\nbigt)$ is naturally defined.
If $f$ is the inclusion of a submanifold,
$f^{\ast}(\nbigt)$ is denoted by $\nbigt_{|Y}$
if there is no risk of confusion.

The pull back of
$\nbigu(p,q)$ and $\newTate(n)$ by $a_X:X\to\pt$
are often denoted by the same notation
$\nbigu(p,q)$ and $\newTate(n)$, respectively.

For $\nbigt_1=(\nbigm_1',\nbigm_1'',C_1)\in\rtriplecat_{\sm}(X,H)$
and $\nbigt_2=(\nbigm_2',\nbigm_2'',C_2)\in\rtriplecat(X,H)$,
the tensor product
$\nbigt_1\otimes\nbigt_2
=(\nbigm_1'\otimes\nbigm_2',\nbigm_1''\otimes\nbigm_2'',C_1\otimes C_2)$
is naturally defined.
In particular, for any $\nbigt\in\rtriplecat(X,H)$,
$\nbigt\otimes\newTate(n)$
is the $n$-th Tate twist of $\nbigt$.

\subsubsection{Twistor structures and smooth $\nbigr_X$-triples}
\label{subsection;22.3.27.1}

Let us observe that smooth $\nbigr_X$-triples
are equivalent to a twistor structure on $X$.
Let $\nbigv=(V^{\sankaku},\delbar_{\proj^1,V^{\sankaku}},\DD^{\sankaku})
\in \TS(X)$.
We set
$V_0:=V^{\sankaku}_{|\nbigx}$
and
$V_{\infty}:=V^{\sankaku}_{|\nbigx^{\dagger}}$.
Let $\nbigv_0$ and $\nbigv_{\infty}$
denote the sheaf of holomorphic sections of $V_0$ on $\nbigx$
and $V_{\infty}$ on $\nbigx^{\dagger}$,
respectively.
Then, $\nbigv_0^{\lor}$
and $\sigma^{\ast}(\nbigv)$
are naturally $\nbigr_X$-modules.
The isomorphism
$V_{0|\cnum_{\lambda}^{\ast}\times X}
\simeq
V_{\infty|\cnum_{\mu}^{\ast}\times X^{\dagger}}$
naturally induces
a sesqui-linear pairing
\[
C_{\nbigv}:(\nbigv_0^{\lor})_{|\vecS\times X}
\otimes_{\cnum}
\sigma^{-1}(\sigma^{\ast}(\nbigv_{\infty})_{|\vecS\times X})
\lrarr
\nbigc^{\infty}_{\vecS\times X}.
\]
In this way,
we obtain a smooth $\nbigr_X$-triple
$\Theta(\nbigv)=(\nbigv_0^{\lor},\sigma^{\ast}\nbigv_{\infty},C_{\nbigv})$.
It is easy to see that this induces a fully faithful functor
$\Theta:\TS(X)\lrarr\rtriplecat_{\sm}(X)$.
The essential image is also called
the category of twistor structures on $X$.
By this equivalence,
a $(-1)^n$-symmetric pairing of weight $n$ of $\nbigv$
corresponds to
a Hermitian sesqui-linear duality of weight $n$ of
$\Theta(\nbigv)$.

Let $\nbigt$ be a smooth $\nbigr_{X}$-triple
equipped with an increasing filtration $W$
indexed by integers.
It is called a graded polarizable mixed twistor structure on $X$
if there exists a graded polarizable mixed twistor structure
$(\nbigv,W)$ on $X$ in the sense of \S\ref{subsection;22.3.25.10}
such that
$(\nbigt,W)\simeq\Theta(\nbigv,W)$.
Let $\MTS(X)$ denote the category of
graded polarizable mixed twistor structure on $X$.
We shall often omit ``graded polarizable''
because we always impose it.
Let $\PTS(X,n)\subset\MTS(X)$
denote the full subcategory of objects
$(\nbigt,W)$ such that
$\Gr^W_m(\nbigt)=0$ $(m\neq n)$,
and called the category of
polarizable pure twistor structure of weight $n$ on $X$.
We shall often omit ``polarizable'' because we always impose it.
An object $(\nbigt,W)$ in $\PTS(X,n)$ is often denoted
just by $\nbigt$.
For an object $\nbigt\in\PTS(X,n)$,
there exists a polarizable pure twistor structure
$\nbigv$ of weight $n$ on $X$ in the sense of \S\ref{subsection;22.3.25.10}
such that $\nbigt\simeq\Theta(\nbigv)$.
In that case,
a Hermitian sesqui-linear duality of $\nbigt$
is called a polarization of $\nbigt$
if it corresponds to a polarization of $\nbigv$.
The following lemma is a translation of
Lemma \ref{lem;22.3.25.22}.

\begin{lem}
$\PTS(X,n)$ is an abelian and semisimple category,
and $\MTS(X)$ is an abelian category.
\hfill\qed 
\end{lem}

The following theorem is proved in \cite{sabbah2}.
\begin{thm}
Any $(\nbigt,W)\in\PTS(X,n)$ is
a polarizable pure twistor $\nbigd$-module of weight $n$. 
\hfill\qed
\end{thm}

The following theorem is contained in
\cite[Theorem 10.3.1]{Mochizuki-MTM}.
\begin{thm}
Any $(\nbigt,W)\in\MTS(X)$ is a mixed twistor $\nbigd$-module.
\hfill\qed
\end{thm}

\subsubsection{Mixed twistor structure on $(X,H)$}

Let $\nbigt=(\nbigm',\nbigm'',C)$
be a smooth $\nbigr_{X(\ast H)}$-triple
equipped with an increasing filtration $W$
by smooth $\nbigr_{X(\ast H)}$-triples.
It is called graded polarizable mixed twistor structure
on $(X,H)$
if $\nbigt_{|X\setminus H}\in\MTS(X\setminus H)$.
Let $\MTS(X,H)$ denote the category of
mixed twistor structure on $(X,H)$.
Let $\PTS(X,H,n)\subset\MTS(X,H)$ denote the full subcategory of
objects $(\nbigt,W)$ such that $\Gr^W_m(\nbigt)=0$ $(m\neq n)$,
which is called the category of polarizable pure twistor
structure of weight $n$ on $(X,H)$.
We obtain the following lemma.
\begin{lem}
$\PTS(X,H,n)$ is an abelian and semisimple category,
and $\MTS(X,H)$ is an abelian category.
\hfill\qed 
\end{lem}

\subsection{Regular admissible mixed twistor structure on $(X,H)$}
\label{subsection;22.3.26.1}

\subsubsection{Regular-KMS-structure}
\label{subsection;22.3.25.120}

Suppose that $H$ is a simple normal crossing hypersurface of $X$.
Let $H=\bigcup_{i\in\Lambda}H_i$ be the irreducible decomposition.
For any $\lambda_0$,
let $\nbigxzero$ denote a neighbourhood of
$\{\lambda_0\}\times X$ in $\nbigx$.
We use the notation
$\nbighzero_i$ and $\nbighzero$ in similar meanings.
Let $V_H\nbigr_X\subset\nbigr_X$ denote
the sheaf of subalgebras generated by
$\lambda p_{\lambda}^{\ast}\Theta_X(\log H)$
over $\nbigo_{\nbigx}$.

Let $\nbigm$ be a smooth $\nbigr_{X(\ast H)}$-module.
We say that $\nbigm$ is a regular-KMS $\nbigr_{X(\ast H)}$-module
if for any $\lambda_0\in\cnum$
there exist a neighbourhood $\nbigxzero$
of $\{\lambda_0\}\times X$
and $V_H\nbigr_{X|\nbigxzero}$-submodules
$\nbigpzero_{\veca}\nbigm$ $(\veca\in\real^{\Lambda})$
of $\nbigm_{|\nbigxzero}$
such that the following holds.
\begin{itemize}
 \item $\bigcup_{\veca\in\real^{\Lambda}}\nbigpzero_{\veca}\nbigm
       =\nbigm_{|\nbigxzero}$.
 \item $\nbigpzero_{\ast}\nbigm=
       (\nbigpzero_{\veca}\nbigmzero\,|\,\veca\in\real^{\Lambda})$
       is a filtered bundle
       on $(\nbigxzero,\nbighzero)$
       (see \cite{i-s} or
       \cite[\S2.3]{Mochizuki-KH-Higgs}).
       Note that for each $i\in\Lambda$ and $\veca\in\real^{\Lambda}$,
       we obtain the locally free $\nbigo_{\nbighzero_i}$-module
\[
       \lefttop{i}\Gr^{\nbigpzero}_{a_i}(\nbigpzero_{\veca}\nbigm)
       =\frac{\nbigpzero_{\veca}(\nbigm)}
       {\sum_{\epsilon>0}\nbigpzero_{\veca-\epsilon\veciti_i}\nbigm},
\]
       on which the endomorphism $\Res_j(\DD)$ is obtained
       as the residue.
       Here, $\veciti_i\in\real^{\Lambda}$
       denote the element whose $j$-th element is $1$ $(j=i)$
       or $0$ $(j\neq i)$.
 \item For each $i\in\Lambda$ and $\veca\in\real^{\Lambda}$,
       there exists a decomposition
\[
       \lefttop{i}\Gr^{\nbigpzero}_{a_i}(\nbigpzero_{\veca}\nbigm)
       =\bigoplus_{\alpha\in\cnum}
       \lefttop{i}\EEzero_{\alpha}\bigl(
       \lefttop{i}\Gr^{\nbigpzero}_{a_i}(\nbigpzero_{\veca}\nbigm)
       \bigr)
\]       
       which is preserved by $\Res_i(\DD)$.
       Moreover,
       for $\alpha\in\cnum$,
       let $u\in\real\times\cnum$
       be determined by
       $\kmsmap(\lambda_0,u)=(a_i,\alpha)$,
       and then $\Res_i(\DD)-\eigenmap(\lambda,u)$
       is nilpotent on
       $\lefttop{i}\EEzero_{\alpha}\bigl(
       \lefttop{i}\Gr^{\nbigpzero}_{a_i}(\nbigpzero_{\veca}\nbigm)
       \bigr)$.
       Here,
       $\kmsmap(\lambda_0,u)
       :=(\paramap(\lambda_0,u),\eigenmap(\lambda_0,u))$.
\end{itemize}

For $i\in\Lambda$,
we choose $P\in H_i^{\circ}:=H_i\setminus \bigcup_{j\neq i}H_j$.
Let $X_P$ be a neighbourhood of $P$
such that $X_P\cap H=X_P\cap H_i=:H_P$.
On $\nbighzero_P$,
we obtain
\[
\Gr^{\nbigpzero}_{a}(\nbigm_{|X_P})
=\bigoplus_{\alpha}
\EEzero_{\alpha}\Gr^{\nbigpzero}_{a}(\nbigm_{|X_P})
\]
as above.
We set
\[
\nbiggzero_{u}(\nbigm_{|X_P}):=
\EEzero_{\eigenmap(\lambda_0,u)}
\Gr^{\nbigpzero}_{\paramap(\lambda_0,u)}
(\nbigm_{|X_P}).
\]
Then, the set
$\KMS(\nbigm,i):=\bigl\{
u\in\real\times\cnum\,\big|\,
\nbiggzero_u(\nbigm_{|X_P})\neq 0
\bigr\}$
is independent of $P\in H_i^{\circ}$
and $\lambda_0\in\cnum$.
We can also obtain $\nbigr_{|\nbigh_P}$-module
$\nbigg_u(\nbigm_{|X_P})$
by gluing $\nbiggzero_u(\nbigm_{|X_P})$
for varying $\lambda_0\in\cnum$.

We note the following proposition
which is due to the regularity condition.

\begin{lem}
\label{lem;22.3.25.110}
Let $\nbigm$ and $\nbigm'$ be smooth regular-KMS
$\nbigr_{X(\ast H)}$-module.
An isomorphism
$f:\nbigm_{|\nbigx\setminus\nbigh}\simeq\nbigm'_{|\nbigx\setminus\nbigh}$
extends to an isomorphism $\nbigm\simeq\nbigm'$
under which 
we have
$\nbigpzero_{\veca}(\nbigm)=\nbigpzero_{\veca}(\nbigm')$
for any $\lambda_0\in\cnum$.
\end{lem}
\pf
For $\lambda\in\cnum$,
let $i_{\lambda}:\{\lambda\}\times X\lrarr \nbigx$
denote the inclusion.
We obtain
$\nbigo_X(\ast H)$-modules
$i_{\lambda}^{\ast}\nbigm$
and $i_{\lambda}^{\ast}\nbigm'$.
If $\lambda\neq 0$,
they are naturally regular singular meromorphic flat bundles.
Hence, 
$i_{\lambda}^{\ast}(f)$ extends to an isomorphism
$i_{\lambda}^{\ast}\nbigm\simeq
i_{\lambda}^{\ast}\nbigm'$
by the regularity.
Then, we can easily observe that
$f$ extends to an isomorphism
$\nbigm\simeq\nbigm'$.
By \cite[Lemma 2.8.3]{Mochizuki-wild},
the KMS-structures are the same.
\hfill\qed

\begin{rem}
\label{rem;22.4.27.2}
If $\Lambda$ is finite,
any $\lambda_0\in\cnum$
has an open neighbourhood $U(\lambda_0)$
in $\cnum$
such that
$\paramap(\lambda):\KMS(\nbigm,i)\to\real$
is injective
for any $i\in\Lambda$
and $\lambda\in U(\lambda_0)\setminus\{\lambda_0\}$.
In this case, we may assume $\nbigxzero=U(\lambda_0)\times X$.
 If $\lambda_1\in U(\lambda_1)\subset U(\lambda_0)$,
$\nbigp^{(\lambda_1)}_{\ast}\nbigm$ 
is reconstructed from 
$\nbigpzero_{\ast}\nbigm_{|\nbigx^{(\lambda_1)}}$
as follows.
Let $\vecb\in\real^{\Lambda}$.
We set
\[
 a_i(\vecb):=
 \max\bigl\{
 \paramap(\lambda_0,u)\,\big|\,
 u\in\KMS(\nbigm,i),\,
 \paramap(\lambda_1,u)\leq b_i
 \bigr\}.
\]
We obtain $\veca(\vecb)=(a_i(\vecb))\in\real^{\Lambda}$.
There exists the natural morphism
\begin{equation}
\label{eq;22.4.27.1}
 \nbigpzero_{\veca(\vecb)}\nbigm_{|\nbigx^{(\lambda_1)}}
 \lrarr
 \bigoplus_{i\in\Lambda}
 \lefttop{i}\Gr^{\nbigpzero}_{a_i(\vecb)}\bigl(
 \nbigpzero_{\veca(\vecb)}\nbigm
 \bigr)_{|\nbigh^{(\lambda_1)}_i}
=\bigoplus_{i\in\Lambda}
\bigoplus_{\substack{u\in\KMS(\nbigm,i)\\
 \paramap(\lambda_0,u)=a_i(\vecb)
 }}
 \lefttop{i}\EEzero_{\eigenmap(\lambda_0,u)}
 \lefttop{i}\Gr^{\nbigpzero}_{a_i(\vecb)}\bigl(
 \nbigpzero_{\veca(\vecb)}\nbigm
 \bigr)_{|\nbigh^{(\lambda_1)}_i}
\end{equation}
Then,
$\nbigp^{(\lambda_1)}_{\vecb}\nbigm$
is the inverse image of
\[
 \bigoplus_{i\in\Lambda}
\bigoplus_{\substack{u\in\KMS(\nbigm,i)\\
 \paramap(\lambda_0,u)=a_i(\vecb)\\
 \paramap(\lambda_1,u)\leq b_i
 }}
 \lefttop{i}\EEzero_{\eigenmap(\lambda_0,u)}
 \lefttop{i}\Gr^{\nbigpzero}_{a_i(\vecb)}\bigl(
 \nbigpzero_{\veca(\vecb)}\nbigm
 \bigr)_{|\nbigh^{(\lambda_1)}_i}.
\]
by {\rm(\ref{eq;22.4.27.1})}.
Similarly, 
we can reconstruct
$\nbigpzero_{\ast}\nbigm_{|\nbigx^{(\lambda_1)}}$
from $\nbigp^{(\lambda_1)}_{\ast}\nbigm$.
\hfill\qed
\end{rem}

\subsubsection{Compatibility of a filtration and a KMS-structure}

Let $\nbigm$ be a smooth regular-KMS $\nbigr_{X(\ast H)}$-module.
Let $W$ be an increasing filtration of $\nbigm$
by $\nbigr_{X(\ast H)}$-submodules indexed by $\seisuu$.
We say that $W$ is compatible with the regular KMS-structure of
$\nbigm$
if the following conditions are satisfied.
\begin{itemize}
 \item Each $\Gr^W_j(\nbigm)$ is smooth regular-KMS $\nbigr_{X(\ast H)}$-module.
 \item For any $\lambda_0\in\cnum$
       and $\veca\in\real^{\Lambda}$,
       the morphism
      $W_j(\nbigm)_{|\nbigxzero}\cap
\nbigpzero_{\veca}(\nbigm)
      \lrarr  \Gr^W_j(\nbigm)_{|\nbigxzero}$
       induces an epimorphism
       $W_j(\nbigm)_{|\nbigxzero}\cap
       \nbigpzero_{\veca}(\nbigm)
       \lrarr
       \nbigpzero_{\veca}\Gr^W_j(\nbigm)$.
\end{itemize}
This definition is equivalent to
\cite[Definition 5.2.6]{Mochizuki-MTM}.

\subsubsection{Regular admissible mixed twistor structure}
\label{subsection;22.4.25.41}

Let $(\nbigt,W)\in\MTS(X,H)$.
It is called regular admissible
if the following conditions are satisfied.
\begin{description}
 \item[(Adm1)] $\nbigm'$ and $\nbigm''$ are
       smooth regular-KMS $\nbigr_{X(\ast H)}$-modules.
       The filtrations of $\nbigm'$  and $\nbigm''$
       are compatible with the KMS-structure.
       (See \S\ref{subsection;22.3.25.100}
       for the induced filtrations
       on $\nbigm'$ and $\nbigm''$.)
 \item[(Adm2)]
	    Let $i\in\Lambda$,
	    and let $P$ be any smooth point of $H_i^{\circ}$.
   	    We use the notation in \S\ref{subsection;22.3.25.120}.
	    Let $\nbign$ denote the nilpotent part of
	    $\Res(\DD)$ on $\nbigg_u(\nbigm'_{|X_P})$
	    and $\nbigg_u(\nbigm''_{|X_P})$ $(u\in\real\times\cnum)$.
	    Note that
	    $\nbigg_u(\nbigm'_{|X_P})$ and $\nbigg_u(\nbigm''_{|X_P})$
	    are equipped with the induced filtration $W$.
	    Then,
	    there exist relative monodromy filtrations
	    of $\nbign$
	    on 
	    $(\nbigg_u(\nbigm'_{|X_P}),W)$ and
	    $(\nbigg_u(\nbigm''_{|X_P}),W)$.
	    (For example,
	    see \cite[\S6.1]{Mochizuki-MTM} for relative monodromy filtrations.)
\end{description}

We remark that {\bf(Adm0)} in \cite[\S9.1.2]{Mochizuki-MTM}
is always satisfied in the regular case,
as remarked in Lemma \ref{lem;22.3.25.110}.

Let $\MTS^{\adm}_{\reg}(X,H)\subset\MTS(X,H)$
denote the full subcategory of
admissible mixed twistor structure on $(X,H)$.
\begin{prop}[\mbox{\cite[Proposition 9.1.7]{Mochizuki-MTM}}]
$\MTS^{\adm}_{\reg}(X,H)$ is an abelian category.
\hfill\qed
\end{prop}

We also recall the following proposition
\cite[Proposition 9.1.8]{Mochizuki-MTM}.
\begin{prop}
\label{prop;22.3.27.2}
Let $(\varphi',\varphi''):(\nbigt_1,W)\to(\nbigt_2,W)$
be a morphism in $\MTS^{\adm}_{\reg}(X,H)$.
Then, $\varphi'$ and $\varphi''$ are 
strictly compatible with the KMS-structure,
i.e.,
for any $\lambda_0$ and $\veca\in\real^{\Lambda}$,
we have
$\varphi'(\nbigpzero_{\veca}\nbigm'_2)
 =\Image \varphi'_{|\nbigxzero}\cap
 \nbigpzero_{\veca}(\nbigm'_1)$
and
$\varphi''(\nbigpzero_{\veca}\nbigm''_1)
=\Image \varphi''_{|\nbigxzero}\cap
\nbigpzero_{\veca}(\nbigm''_2)$.
We also have
\[
 \nbigp_{\veca}(\Ker\varphi'')=
 \Ker\varphi''\cap\nbigp_{\veca}\nbigm''_1,
 \quad\quad
\nbigp_{\veca}(\Cok\varphi'')
 =\Image\bigl(
 \nbigp_{\veca}\nbigm''_2
 \to\Cok(\varphi'')
 \bigr),
\]
and similar equalities for $\varphi'$.
\hfill\qed
\end{prop}

The following theorem is proved
in \cite[Theorem 10.3.1]{Mochizuki-MTM}.
\begin{thm}
Let $(\nbigt,W)$ be a regular mixed twistor $\nbigd$-module
on $X$.
If $\nbigt(\ast H)$
is a smooth $\nbigr_{X(\ast H)}$-triple,
$(\nbigt,W)(\ast H)$ is
a regular admissible mixed twistor structure on $(X,H)$.
Conversely,
for an regular admissible mixed twistor structure $(\nbigt,W)$
on $(X,H)$,
there exists a mixed twistor $\nbigd$-module
$(\nbigttilde,W)$ on $X$
such that $(\nbigttilde,W)(\ast H)=(\nbigt,W)$.
\hfill\qed
\end{thm}

Indeed,
for any regular admissible mixed twistor structure
$(\nbigt,W)$ on $(X,H)$,
we have the mixed twistor $\nbigd$-modules
$(\nbigt,W)[\ast H]$
and
$(\nbigt,W)[! H]$.
See \cite[\S5]{Mochizuki-MTM}
for the construction of the $\nbigr_X$-triple
$\nbigt[\star H]$,
and see \cite[\S10.2]{Mochizuki-MTM}
for the construction of the weight filtrations
on $\nbigt[\star H]$.

\subsubsection{Nearby cycle functors}

Let $Y$ be a complex manifold.
Let $\Lambda$ be a finite set.
Suppose that $X$ is a neighbourhood of
$Y\times\{(0,\ldots,0)\}$
in $Y\times\cnum^{\Lambda}$.
Let
$H_i=(Y\times\cnum^{\Lambda\setminus\{i\}})\cap X$
and $H=\bigcup H_i$.
Let $(\nbigt,W)\in \MTS^{\adm}_{\reg}(X,H)$.

Let $i\in\Lambda$.
For any $u\in\real\times\cnum$,
we obtain
the $\nbigr$-triple
$\psitilde_{z_i,u}(\nbigt)$.
It is equipped with the filtration
$L_j\psitilde_{z_i,u}(\nbigt)
=\psitilde_{z_i,u}(W_j\nbigt)$.
By the admissibility condition,
$\nbign:\psitilde_{z_i,u}(\nbigt)
\to\psitilde_{z_i,u}(\nbigt)\otimes\newTate(-1)$
has a relative monodromy weight filtration $W=M(\nbign;L)$
with respect to the filtration $L$.
We obtain the following proposition from
\cite[Proposition 9.4.1]{Mochizuki-MTM}.
\begin{prop}
\label{prop;22.4.25.50}
$(\psitilde_{z_i,u}(\nbigt),W)\in\MTS^{\adm}_{\reg}(H_i,\del H_i)$,
where $\del H_i=H_i\cap\bigl(\bigcup_{j\neq i}H_j\bigr)$. 
\hfill\qed
\end{prop}

\subsection{Tame harmonic bundles with normal crossing singularity}
\label{subsection;22.4.26.10}

Let $X$ be a complex manifold
with a simple normal crossing hypersurface $H$.
Let $(E,\delbar_E,\theta,h)$ be a harmonic bundle
on $X\setminus H$.
The Higgs field induces a coherent sheaf on
the cotangent bundle of $X\setminus H$.
Its support $\Sigma_{E,\theta}$ is called the spectral variety.
The harmonic bundle $(E,\delbar_E,\theta,h)$ is called tame
if the closure of $\Sigma_{E,\theta}$
in the logarithmic cotangent bundle $T^{\ast}X(\log H)$
is proper over $X$.

\subsubsection{The associated polarized twistor structure on $X\setminus H$
as smooth $\nbigr_{X\setminus H}$-triple}

The bundle $p_{\lambda}^{-1}(E)$
is equipped with the natural holomorphic structure
$p_{\lambda}^{\ast}(\delbar_E)$.
We obtain the holomorphic vector bundle
$\nbige=(p_{\lambda}^{-1}(E),
p_{\lambda}^{\ast}(\delbar_E)+\lambda\theta^{\dagger})$.
The sheaf of holomorphic sections is also denoted by $\nbige$.
It is equipped with the differential operator
$\DD^f=\nabla_h+\lambda^{-1}\theta+\lambda\theta^{\dagger}$,
and naturally an $\nbigr_X$-module.
We have the sesqui-linear pairing
$C_h:
\nbige_{|\vecS\times (X\setminus H)}
\otimes
\sigma^{-1}(\nbige_{|\vecS\times (X\setminus H)})
\lrarr
\distribution_{\vecS\times (X\setminus H)/\vecS}$
given by
$C_h(u,\sigma^{-1}v)=h(u,\sigma^{-1}(v))$.
Thus,
we obtain the smooth $\nbigr_{X\setminus H}$-triple
$\nbigt(E)=(\nbige,\nbige,C_h)$.
It is a pure twistor structure of weight $0$.
The pair of the identity morphism
$(\id,\id):\nbigt(E)\to\nbigt(E)^{\ast}$
is a polarization of weight $0$ of $\nbigt(E)$.

\subsubsection{The smooth regular-KMS $\nbigr_{X(\ast H)}$-modules
associated with tame harmonic bundles}

For any open subset $U$ of $\nbigx$,
let $\nbigp\nbige(U)$
denote the space of holomorphic sections $s$ of
$\nbige$ on $U\setminus\nbigh$
such that the following condition is satisfied
for any $(\lambda,P)\in U$.
\begin{itemize}
 \item Let $(X_P,z_1,\ldots,z_n)$
       be a relatively compact
       holomorphic coordinate neighbourhood around $P$ in $U$
       such that $X_P\cap H=\bigcup_{i=1}^{\ell}\{z_i=0\}$.
       Then, there exists $N>0$ such that
       $|s|_h=O\bigl(\prod_{i=1}^{\ell}|z_i|^{-N}\bigr)$.
\end{itemize}
The following theorem is proved in \cite{mochi2}.
\begin{thm}
$\nbigp\nbige$ is a smooth regular-KMS $\nbigr_{X(\ast H)}$-module.
\hfill\qed
\end{thm}

\subsubsection{The associated pure twistor $\nbigd$-modules}
\label{subsection;22.4.12.1}

Let $H=\bigcup_{i\in\Lambda}H_i$
be the irreducible decomposition.
Let $\veciti_{\Lambda}\in\real^{\Lambda}$
be the element whose $j$-th elements are $1$.
For $\lambda_0\in\cnum$,
we set
\[
\nbigpzero_{<\veciti_{\Lambda}}\nbige
=\bigcup_{\substack{\veca\in\real^{\Lambda}\\ a_i<1 }}
\subset\nbigp\nbige_{|\nbigxzero}.
\]
There exists the $\nbigr_X$-submodule
$\gbige\subset\nbigp\nbige$
determined as follows on $\nbigxzero$:
\[
\gbige_{|\nbigxzero}
=\nbigr_{X|\nbigxzero}\cdot
\nbigpzero_{<\veciti_{\Lambda}}\nbige
\subset
\nbigp\nbige_{|\nbigxzero}.
\]
As proved in \cite{mochi2},
the $\nbigr_{X\setminus H}$-triple
$\nbigt(E)$
extends to the $\nbigr_X$-triple
$\gbigt(E):=(\gbige,\gbige,\gbigc)$.
The following theorem is proved in \cite{mochi2}.
\begin{thm}
\label{thm;22.4.18.10}
$\gbigt(E):=(\gbige,\gbige,\gbigc)$
is a pure twistor $\nbigd$-module of weight $0$,
and $\gbigs=(\id,\id)$ is a polarization.
\hfill\qed
\end{thm}

We note that
$\gbigt(E)(\ast H)$ is
a regular admissible mixed twistor structure
on $(X,H)$.
We obtain the mixed twistor $\nbigd$-modules
$\gbigt(E)[\star H]=
(\gbigt(E)(\ast H))[\star H]$ $(\star=!,\ast)$.
By the construction of the $\nbigr_X$-module,
it is easy to see that
$\gbigt(E)$ is the image of
the canonical morphism
$\gbigt(E)[!H]\to \gbigt(E)[\ast H]$.
In particular, $\gbigt(E)$ is a mixed twistor $\nbigd$-module.

\subsubsection{The induced polarized mixed twistor structures}
\label{subsection;22.4.5.10}

Let $Y$ be a complex manifold.
Let $\Lambda$ be a finite set.
Let $X$ be a neighbourhood of
$Y\times\{(0,\ldots,0)\}$ in $Y\times\cnum^{\Lambda}$,
and $H=\bigcup_{i\in\Lambda}(Y\times\cnum^{\Lambda\setminus \{i\}})$.
We use the partial order $\leq$ on $\real^{\Lambda}$
defined by
$\vecb\leq\veca$
$\stackrel{\rm def}{\Longleftrightarrow}$
$b_i\leq a_i$ $(\forall i\in\Lambda)$.
We say $\vecb\lneq\veca$ if
$\vecb\leq\veca$ and $\vecb\neq\veca$.

Let $(E,\delbar_E,\theta,h)$ be a tame harmonic bundle
on $(X,H)$.
We obtain the filtered bundle
$(\nbigpzero_{\ast}\nbige,\DD)$
on $\nbigxzero$.
For any $\veca\in\real^{\Lambda}$, we obtain
the following locally free $\nbigo_{\nbigyzero}$-module:
\[
 \lefttop{\Lambda}\Gr^{\nbigpzero}_{\veca}(\nbigp\nbige):=
 \frac{\nbigpzero_{\veca}\nbige}
 {\sum_{\vecb\lneq\veca}\nbigpzero_{\vecb}\nbige}.
\]
It is equipped with the endomorphisms
$\Res_i(\DD)$ $(i\in\Lambda)$.
We obtain the decomposition
\begin{equation}
\label{eq;22.4.5.1}
 \lefttop{\Lambda}\Gr^{\nbigpzero}_{\veca}(\nbigp\nbige)
 =\bigoplus_{\vecalpha\in\cnum^{\Lambda}}
 \lefttop{\Lambda}\nbiggzero_{\vecu}(\nbigp\nbige)
\end{equation}
such that the following holds:
\begin{itemize}
 \item The decomposition (\ref{eq;22.4.5.1}) is preserved by
       $\Res_i(\DD)$.
 \item Let $\vecu\in(\real\times\cnum)^{\Lambda}$
       be determined by
       $\kmsmap(\lambda_0,u_i)=(a_i,\alpha_i)$.
       Then,
       $\Res_i(\DD)-\eigenmap(\lambda,u_i)$
       are nilpotent on
       $\lefttop{\Lambda}\nbiggzero_{\vecu}(\nbigp\nbige)$.
\end{itemize}
If $\nbigx^{(\lambda_1)}\subset\nbigxzero$,
there exist natural isomorphisms
$\lefttop{\Lambda}\nbigg^{(\lambda_1)}_{\vecu}(\nbigp\nbige)
\simeq
\lefttop{\Lambda}\nbigg^{(\lambda_0)}_{\vecu}(\nbigp\nbige)
_{|\nbigx^{(\lambda_1)}}$.
Hence, by varying $\lambda_0\in\cnum$
and by gluing
$\lefttop{\Lambda}\nbiggzero_{\vecu}(\nbigp\nbige)$,
we obtain
a locally free $\nbigo_{\nbigy}$-module
$\lefttop{\Lambda}\nbigg_{\vecu}(\nbigp\nbige)$.
It is naturally a smooth $\nbigr_{Y}$-module.
The nilpotent part of $-\lambda^{-1}\Res_i(\DD)$ induces
$\nbign_i:
\lefttop{\Lambda}\nbigg_{\vecu}\to
\lambda^{-1}\lefttop{\Lambda}\nbigg_{\vecu}$.
There exists a naturally induced sesqui-linear pairing
$C_{h,\vecu}:
\lefttop{\Lambda}\nbigg_{\vecu}(\nbigp\nbige)_{|\vecS\times Y}
\otimes
\sigma^{-1}(
\lefttop{\Lambda}\nbigg_{\vecu}(\nbigp\nbige)_{|\vecS\times Y})
\lrarr
\distribution_{\vecS\times Y/\vecS}$.
Thus, we obtain a smooth $\nbigr_Y$-triple
$\lefttop{\Lambda}\nbigg_{\vecu}(\nbigt(E)):=
\bigl(
 \lefttop{\Lambda}\nbigg_{\vecu}(\nbigp\nbige),
 \lefttop{\Lambda}\nbigg_{\vecu}(\nbigp\nbige),
 C_{h,\vecu}
\bigr)$.
The morphism
$\nbign_i:
\lefttop{\Lambda}\nbigg_{\vecu}(\nbigp\nbige)
\to \lambda^{-1}
\lefttop{\Lambda}\nbigg_{\vecu}(\nbigp\nbige)$
induces
$\nbign_i:\lefttop{\Lambda}\nbigg_{\vecu}(\nbigt(E))
\to
\lefttop{\Lambda}\nbigg_{\vecu}(E)\otimes\newTate(-1)$.
Let $\vecnbign=(\nbign_i\,|\,i\in\Lambda)$.
We also obtain a Hermitian sesqui-linear duality
$\nbigs_{\vecu}=(\id,\id)$ of
$(\lefttop{\Lambda}\nbigg_{\vecu}(\nbigt(E)),\vecnbign)$.
The following theorem is proved in \cite{mochi2}.

\begin{thm}
 $(\lefttop{\Lambda}\nbigg_{\vecu}(\nbigt(E)),
 \vecnbign,\nbigs_{\vecu})$
is a $(0,\Lambda)$-polarized mixed twistor structure
on $Y$.
\hfill\qed 
\end{thm}

\subsubsection{The induced admissible mixed twistor structure}

Let $X$, $Y$, $\Lambda$ be as in \S\ref{subsection;22.4.5.10}.
Let $\Lambda=I\sqcup J$ be a decomposition.
We set $H_I:=Y\times\cnum^{\Lambda\setminus I}$.
Let $q_I:\real^{\Lambda}\to\real^I$ denote the projection.
For $\vecc\in\real^I$,
we set
\[
 \lefttop{I}\nbigpzero_{\vecc}\nbige:=
 \sum_{
 \substack{\veca\in\real^{\Lambda}\\ q_I(\veca)=\vecc
 }}
 \nbigpzero_{\veca}\nbige.
\]
We also set
\[
 \lefttop{I}\Gr^{\nbigpzero}_{\vecc}(\nbigp\nbige):=
 \frac{\lefttop{I}\nbigpzero_{\vecc}\nbige}
 {\sum_{\vecb\lneq\vecc}\lefttop{I}\nbigpzero_{\vecd}\nbige}.
\]
Then, it is a locally free
$\nbigo_{\nbighzero_I}(\ast\del \nbighzero_I)$-module.
Here $\del H_I:=H_I\cap\bigl(\bigcup_{j\in J}H_j\bigr)$.
It is equipped with the endomorphisms $\Res_i(\DD)$ $(i\in I)$.
There exists the decomposition
\[
 \lefttop{I}\Gr^{\nbigpzero}_{\vecc}(\nbigp\nbige)=
 \bigoplus_{\vecu\in(\real\times\cnum)^I}
 \lefttop{I}\nbiggzero_{\vecu}(\nbigp\nbige)
\]
such that
(i) the decomposition is preserved by
$\Res_i(\DD)$ $(i\in I)$,
(ii) $\Res_i(\DD)-\eigenmap(\lambda,u_i)$ are nilpotent,
where $u_i$ $(i\in I)$ are determined by
$\kmsmap(\lambda_0,u_i)=(a_i,\alpha_i)$.
By varying $\lambda_0\in\cnum$,
and by gluing $\lefttop{I}\nbiggzero(\nbigp\nbige)$,
we obtain
the $\nbigo_{\nbigh_I}(\ast\del\nbigh_I)$-module
$\lefttop{I}\nbigg_{\vecu}(\nbigp\nbige)$.
The nilpotent part of $\Res_i(\DD)$ induce
$\lefttop{I}\nbigg_{\vecu}(\nbigp\nbige)
\to
\lambda^{-1}\lefttop{I}\nbigg_{\vecu}(\nbigp\nbige)$
$(i\in I)$.
There exist the sesqui-linear pairing
 $C_{h,\vecu}:
 \lefttop{I}\nbigg_{\vecu}(\nbigp\nbige)_{|\vecS\times H_I}
 \otimes
 \sigma^{-1}\bigl(
  \lefttop{I}\nbigg_{\vecu}(\nbigp\nbige)_{|\vecS\times H_I}
  \bigr)
  \lrarr
  \nbigc^{\infty,\mod \del H_I}_{\vecS\times H_I/\vecS}$
 induced by $C_h$,
 and 
 $\lefttop{I}\nbigg_{\vecu}(\nbigt(E))
 =
 \bigl(
 \lefttop{I}\nbigg_{\vecu}(\nbigp\nbige),
 \lefttop{I}\nbigg_{\vecu}(\nbigp\nbige),
 C_{h,\vecu}
 \bigr)$
is a twistor structure on $(H_I,\del H_I)$.
We obtain $\nbign_i:\lefttop{I}\nbigg_{\vecu}(\nbigt(E))
\to\lefttop{I}\nbigg_{\vecu}(\nbigt(E))
 \otimes\newTate(-1)$
obtained as $(\nbign_i,\nbign_i)$ .
Let $W$ be the monodromy weight filtration of
$\sum_{i\in I}\nbign_i$.
Then,
$\bigl(
\lefttop{I}\nbigg_{\vecu}(\nbigt(E)),
W
\bigr)$
is a regular admissible mixed twistor structure
on $(H_I,\del H_I)$.
Indeed,
$\bigl(
\lefttop{I}\nbigg_{\vecu}(\nbigt(E)),
W,(\nbign_i\,|\,i\in I)
\bigr)$
is a regular polarizable admissible mixed twistor structure
in the sense of \cite{Mochizuki-MTM}.

\subsubsection{Description of the filtrations}
\label{subsection;22.4.27.3}

Let $X$, $Y$, $\Lambda$ be as in \S\ref{subsection;22.4.5.10}.
Let $(E,\delbar_E,\theta,h)$ be a tame harmonic bundle on $(X,H)$.
Let us describe the filtrations $\nbigpzero_{\ast}\nbige$
more explicitly.

Let $\lambda\in\cnum$.
We obtain the holomorphic vector bundle
$\nbigelambda=(E,\delbar_E+\lambda\theta^{\dagger})$
on $X\setminus H$.
We set
$\DDlambda=
\delbar_E+\lambda\theta^{\dagger}
+\lambda\del_E+\theta$.
It is a flat $\lambda$-connection of $\nbigelambda$,
i.e.,
it satisfies
(i) $\DDlambda(fs)=f\DDlambda(s)+
(\lambda\del_X+\delbar_X)f\cdot s$
for any $f\in C^{\infty}(X)$
and $s\in C^{\infty}(X,\nbigelambda)$,
and (ii) $\DDlambda\circ\DDlambda=0$.

Let $\veca\in\real^{\Lambda}$.
For any open subset $U$ of $X$,
let $\nbigp_{\veca}\nbigelambda(U)$
denote the space of
holomorphic sections $s$ of
$\nbigelambda$ on $U\setminus H$
such that
\[
 |s|_h=O\Bigl(
  \prod_{i\in\Lambda}
  |z_i|^{-a_i-\epsilon}
  \Bigr)
\]
locally around any point of $U$ for any $\epsilon>0$.
Thus, we obtain an $\nbigo_X$-module $\nbigp_{\veca}\nbigelambda$.
We obtain an $\nbigo_X(\ast H)$-module 
$\nbigp\nbigelambda:=\sum \nbigp_{\veca}\nbigelambda$
in $\iota_{X\setminus H\ast}(\nbigelambda)$,
where $\iota_{X\setminus H}:X\setminus H\to X$
denotes the inclusion.
The following theorem is proved
in \cite{Simpson-non-compact} for the one dimensional case,
and in
\cite[Theorem 8.58, Theorem 8.59, Corollary 8.89]{mochi2}
for the higher dimensional case.
\begin{thm}
 $\nbigp\nbigelambda$ is a locally
 free $\nbigo_X(\ast H)$-module,
 and 
 $\nbigp_{\ast}\nbigelambda=
\bigl(\nbigp_{\veca}\nbigelambda\,
\big|\,
\veca\in\real^{\Lambda}
\bigr)$
is a filtered bundle over $\nbigp\nbigelambda$.
The flat $\lambda$-connection
$\DDlambda$ is logarithmic with respect to
each $\nbigp_{\veca}\nbigelambda$.
Namely,
$(\nbigp_{\ast}\nbigelambda,\DDlambda)$ 
is a regular filtered $\lambda$-flat bundle.
\hfill\qed
\end{thm}

For any $\lambda\in\cnum$,
let $\iota_{\lambda}$ denote the inclusion
$X=\{\lambda\}\times X\subset\nbigx$.
By the construction,
there exist the natural morphisms
$\iota_{\lambda}^{\ast}\nbigp\nbige
\lrarr \nbigp\nbigelambda$
for any $\lambda\in\cnum$,
and they are isomorphisms.
We identify 
$\iota_{\lambda}^{\ast}\nbigp\nbige$
and $\nbigp\nbigelambda$.
We have
$\iota_{\lambda_0}^{\ast}\nbigpzero_{\veca}\nbige
=\nbigp_{\veca}\nbige^{\lambda_0}$
for any $\veca\in\real^{\Lambda}$.
By Remark \ref{rem;22.4.27.2},
we can describe the relation between
$\iota_{\lambda}^{\ast}\nbigpzero_{\ast}\nbige$
and $\nbigp_{\ast}\nbigelambda$
for any $\lambda$ which is sufficiently close to $\lambda_0$.

\subsection{Tame harmonic bundles and regular pure twistor $\nbigd$-modules}

\subsubsection{Extension of tame harmonic bundles
to pure twistor $\nbigd$-modules}

Let $X$ be any complex manifold.
Let $Z_1\subset Z_0\subset X$ be closed complex analytic subsets
such that $Z_0\setminus Z_1$ is a locally closed complex submanifold of $X$.
Let $\iota_{Z_0\setminus Z_1}$ denote the inclusion
of $Z_0\setminus Z_1\to X\setminus Z_1$.
A resolution $\varphi:\Ztilde_0\to X$ is
a proper morphism of complex manifolds
satisfying
(i) $\varphi(\Ztilde_0)=Z_0$,
(ii) $\Ztilde_1=\varphi^{-1}(Z_1)$ is a simple normal crossing
hypersurface of $\Ztilde_0$,
(iii) $\varphi$ induces
$\Ztilde_0\setminus\Ztilde_1\simeq Z_0\setminus Z_1$,
(iv) for any relatively compact open subset $U\subset X$,
the induced morphism
$\varphi^{-1}(U)\to U$ is a composition of
blowings-up along smooth centers.
Note that such a resolution exists
(see \cite{Wlodarczyk}).

A harmonic bundle
$(E,\delbar_E,\theta,h)$ on $Z_1\setminus Z_0$
is called tame on $(Z_0,Z_1)$
if
$\varphi^{-1}(E,\delbar_E,\theta,h)$
is a tame on $(\Ztilde_0,\Ztilde_1)$
for a resolution
$\varphi:\Ztilde_0\to X$ of $(Z_0,Z_1)$.
Note that the condition is independent of the choice of
a resolution.

For any tame harmonic bundle $(E,\delbar_E,\theta,h)$
on $Z_0\setminus Z_1$,
we have the polarized pure twistor $\nbigd$-module
$(\nbigt(E),\nbigs)$ on $Z_0\setminus Z_1$.
The following theorem is proved in \cite[Theorem 19.6]{mochi2}.

\begin{thm}
\label{thm;22.4.20.20}
There exists a unique polarized regular pure twistor
$\nbigd$-module $(\gbigt(E),\gbigs)$ of weight $0$ on $X$
such that
$(\gbigt(E),\gbigs)_{|X\setminus Z_1}
 =\iota_{Z_0\setminus Z_1\,\dagger}(\nbigt(E),\nbigs)$.
\hfill\qed
\end{thm}

See \cite[\S19.4]{mochi2} for the uniqueness.
Let us recall the construction of $(\gbigt(E),\gbigs)$.
Let $\varphi:\Ztilde_0\to X$ be a resolution of $(Z_0,Z_1)$.
We obtain a tame harmonic bundle
$(\Etilde,\delbar_{\Etilde},\thetatilde,\htilde)
=\varphi^{-1}(E,\delbar_E,\theta,h)$
on $(\Ztilde_0,\Ztilde_1)$.
We obtain the polarized pure twistor $\nbigd$-module
$(\gbigt(\Etilde),\gbigstilde)$ of weight $0$
on $\Ztilde_0$
associated with 
$(\Etilde,\delbar_{\Etilde},\thetatilde,\htilde)$
as in Theorem \ref{thm;22.4.18.10}.
We obtain the pure twistor $\nbigd$-module
$\varphi_{\dagger}^0\gbigt(\Etilde)$ of weight $0$
with the induced Hermitian sesqui-linear duality
$\varphi_{\dagger}(\gbigstilde):
\varphi_{\dagger}^0\gbigt(\Etilde)
\to
\varphi_{\dagger}^0\gbigt(\Etilde)^{\ast}$.
There exists the decomposition
$(\varphi_{\dagger}^0\gbigt(\Etilde),\varphi_{\dagger}(\gbigstilde))
=(\gbigt_0,\gbigs_0)\oplus (\gbigt_1,\gbigs_1)$,
where the strict support of $\gbigt_0$ is $Z_0$,
and the support of $\gbigt_1$ is contained in $Z_1$.
Let $U\subset X$ be any relatively compact open subset.
Let $c(U)$ be the first Chern class
of a relatively ample line bundle
of $\varphi^{-1}(U)\to U$.
By the support condition,
$\gbigt_{0|U}$ is contained in the primitive part of
$\varphi_{\dagger}^0\gbigt(\Etilde)_{|U}$
with respect to $L_{c(U)}$.
By Theorem \ref{thm;22.4.18.20},
we obtain that $(\gbigt_0,\gbigs_0)_{|U}$ is
a polarized pure twistor $\nbigd$-module of weight $0$.
Because $U$ is any relatively compact open subset of $X$,
$(\gbigt_0,\gbigs_0)$ is a polarized pure twistor
$\nbigd$-module of weight $0$.
By the construction,
we have
$(\gbigt_0,\gbigs_0)_{|X\setminus Z_1}=
\iota_{Z_0\setminus Z_1\dagger}(\nbigt(E),\nbigs)$.
We note the following lemma.
\begin{lem}
\label{lem;22.4.18.31}
Let $W$ be the increasing filtration on $\gbigt(E)$
defined by
$W_j(\gbigt(E))=0$ $(j<0)$
and $W_0(\gbigt(E))=\gbigt(E)$.
Then,
$(\gbigt(E),W)\in\MTM(X)$.
\end{lem}
\pf
Let $\varphi:\Ztilde_0\to X$ and
$(\Etilde,\delbar_{\Etilde},\thetatilde,\htilde)$
be as above.
We have the mixed twistor $\nbigd$-modules
$(\gbigt(\Etilde)[\star \Ztilde_1],W)$ $(\star=!,\ast)$.
Let $W$ be the filtration of $\gbigt(\Etilde)$
defined by
$W_j\gbigt(\Etilde)=0$ $(j<0)$
and
$W_0\gbigt(\Etilde)=\gbigt(\Etilde)$.
Because 
$(\gbigt(\Etilde),W)$ is
the image of
$(\gbigt(\Etilde)[! \Ztilde_1],W)
\to
(\gbigt(\Etilde)[\ast \Ztilde_1],W)$,
it is a mixed twistor
$\nbigd$-module on $\Ztilde_0$.
By \cite[Proposition 7.2.7]{Mochizuki-MTM},
for any relatively compact open subset $U\subset X$,
$\varphi_{\dagger}^0(\gbigt(\Etilde),W)_{|U}$
is a mixed twistor $\nbigd$-module on $U$.
Because $(\gbigt(E),W)_{|U}$ is a direct summand of
$\varphi_{\dagger}^0(\gbigt(\Etilde),W)_{|U}$
as a filtered $\nbigr_U$-triple,
we obtain that
$(\gbigt(E),W)_{|U}\in\MTM(U)$.
Hence, we obtain that
$(\gbigt(E),W)\in\MTM(X)$.
\hfill\qed
 
\subsubsection{The restriction of regular polarized
pure twistor $\nbigd$-modules}

Let $X$ be a complex manifold.
Let $Z_0\subset X$ be a closed complex analytic subset of $X$.
A smooth Zariski open subset $U$ of $Z_0$
means an open subset such that
(i) $U$ is a locally closed complex submanifold of $X$,
(ii) $Z_0\setminus U$ is a closed complex analytic subset of $Z_0$.
The following theorem is proved
in \cite[Theorem 19.6]{mochi2}.

\begin{thm}
\label{thm;22.4.18.30}
For a regular pure twistor $\nbigd$-module $(\nbigt,\nbigs)$
of weight $0$,
there exist a smooth Zariski open subset $U$ of $Z_0$
and a harmonic bundle
$(E,\delbar_E,\theta,h)$ on $U$,
which is tame on $(Z_0,Z_0\setminus U)$,
and $(\nbigt,\nbigs)$ is isomorphic to
the polarized pure twistor $\nbigd$-module
associated with $(E,\delbar_E,\theta,h)$.
\hfill\qed
\end{thm}

By Theorem \ref{thm;22.4.18.30} and Lemma \ref{lem;22.4.18.31},
we obtain the following lemma.
\begin{lem}
Any polarized pure twistor $\nbigd$-module on $X$
is a mixed twistor $\nbigd$-module on $X$.
\hfill\qed
\end{lem}

\subsection{Good mixed twistor $\nbigd$-modules}
\label{subsection;22.4.25.101}

Let $X$ be a complex manifold with a simple normal crossing
hypersurface $H$.
Let $H=\bigcup_{i\in\Lambda}H_i$ be the irreducible decomposition.
For any non-empty $I\subset\Lambda$,
we set $H_I:=\bigcap_{i\in I}H_i$,
$\del H_I:=H_I\cap\bigl(
\bigcup_{i\in \Lambda\setminus I}H_i
\bigr)$,
and $H_I^{\circ}:=H_I\setminus \del H_I$.
We formally set $H_{\emptyset}=X$.
Let $T_{H_I}^{\ast}X$ denote the conormal bundle of
$H_I$ in $X$.

\begin{df}
\label{df;22.3.30.10}
A regular mixed twistor $\nbigd$-module $(\nbigt,W)$ on $X$
is good on $(X,H)$
if the characteristic varieties of $\nbigt$
is contained in
$\cnum\times\bigcup_{I\subset\Lambda}T_{H_I}^{\ast}X$.
\hfill\qed
\end{df}

\subsubsection{Local case}

Let $\Lambda$ be a finite set.
Let $X$ be a neighbourhood of $\{(0,\ldots,0)\}$
in $\cnum^{\Lambda}$.
We put $H_i=\{z_i=0\}\cap X$
and we set $H=\bigcup_{i\in\Lambda}H_i$.
We set $H_I=\bigcap_{i\in I}H_i$
and $H(I)=\bigcup_{i\in I}H_i$.

Let $\nbigt$ be a non-zero
good regular mixed twistor $\nbigd$-module on $(X,H)$.
There exists a subset $S\subset 2^{\Lambda}$
such that the support of $\nbigt$ is
$\bigcup_{I\in S}H_I$.
Let $I_0$ be a minimal element in $S$,
i.e.,
$I\in S$ satisfies $I\subset I_0$ if and only if $I=I_0$.
There exists a mixed twistor structure
$(\nbigt_I,W)$
on $(H_I,\del H_I)$
with an isomorphism
$\iota_{I\dagger}(\nbigt_I,W)
\simeq
(\nbigt,W)(\ast H(\Lambda\setminus I))$.
By the admissibility of $\nbigt$,
we obtain that $(\nbigt_I,W)$
is an admissible mixed twistor structure on $(H_I,\del H_I)$.
We set $f=z^{\Lambda\setminus I}$.
Note that the characteristic variety of
$\Pi^{a,b}_{f\star}\nbigt_I$ $(\star=!,\ast)$
is contained in
$\cnum\times\bigl(
\bigcup_{I\subset J}T_{H_J}^{\ast}X
\bigr)$.
Hence,
the characteristic varieties of
$\nbigt[\star f]$ and
$\Xi^{(0)}_{f}(\nbigt_I)$
are contained in 
$\cnum\times\bigl(
\bigcup_{I\subset J}T_{H_J}^{\ast}X
\bigr)$.
Moreover,
the characteristic varieties of
$\psi^{(a)}_f(\nbigt_I)$,
and
$\phi_f(\nbigt_I)$
are contained in
$\cnum\times
\bigl(
\bigcup_{I\subsetneq J}T_{H_J}^{\ast}X
\bigr)$.

\begin{rem}
We can reconstruct
$\nbigt$
as the gluing of
good regular mixed twistor $\nbigd$-module
$\phi_{f}(\nbigt_I)$
and 
an admissible mixed twistor structure
$\nbigt_I(\ast H(\Lambda\setminus I))$.
By an inductive argument,
it is easy to see that
$\nbigt$ is good on $(X,H)$
in the sense of {\rm\cite{Mochizuki-MTM}}.
Similarly, we can see that
a good regular mixed twistor $\nbigd$-module
on $(X,H)$  in the sense of {\rm\cite{Mochizuki-MTM}}
 is good on $(X,H)$
 in the sense of Definition 
{\rm\ref{df;22.3.30.10}}.
\hfill\qed
\end{rem}

\vspace{.1in}

\noindent
{\em Address\\
Research Institute for Mathematical Sciences,
Kyoto University,
Kyoto 606-8502, Japan\\
takuro@kurims.kyoto-u.ac.jp
}


\begin{thebibliography}{99}

\bibitem{beilinson1}
A. Beilinson,
{\em On the derived category of perverse sheaves.}
in; {\em $K$-theory, arithmetic and geometry 
(Moscow, 1984--1986)},
Lecture Notes in Math., {\bf 1289}, 
Springer, Berlin, (1987),
27--41.

\bibitem{beilinson2}
A. Beilinson,
{\em How to glue perverse sheaves},
in; {\em $K$-theory, arithmetic and geometry
 (Moscow, 1984--1986)}, 
Lecture Notes in Math., {\bf 1289},
Springer, Berlin, (1987),
42--51.

 
\bibitem{bbd}
A. Beilinson, J. Bernstein, P. Deligne,
{\em Faisceaux pervers},
Analysis and topology on singular spaces, I (Luminy, 1981), 
Ast\'{e}risque, {\bf 100},
(1982), 5--171.

\bibitem{borne1}
N. Borne, {\em Fibr\'{e}s paraboliques 
et champ des racines},
Int. Math. Res. Not. IMRN  (2007),
no 16, 38 pages.

\bibitem{borne2}
N. Borne,
Sur les repr\'{e}sentations du groupe fondamental
d'une vari\'{e}t\'{e} priv\'{e}e d'un diviseur 
\`{a} croisements normaux simples. (French)
Indiana Univ. Math. J.  {\bf 58}  (2009), 
137--180.

	
\bibitem{de-Cataldo-Migliorini-semismall}
	M. A. A. de Cataldo,
	L. Migliorini,
	{\em The hard Lefschetz theorem and the topology of semismall maps},
	Ann. Sci. \'{E}cole Norm. Sup. {\bf 35} (2002),
	759--772.
\bibitem{de-Cataldo-Migliorini-decomposition}
	M. A. A. de Cataldo,
	L. Migliorini,
	{\em The perverse filtration and
	the Lefschetz hyperplane theorem},
	Ann. of Math. (2) {\bf 171} (2010), 2089--2113.

\bibitem{ck}
E. Cattani, and A. Kaplan,
{\it Polarized mixed Hodge structures
and the local monodromy of variation of Hodge structure},
Invent. Math. {\bf 67}
(1982), 101--115.

\bibitem{ck2}
	E. Cattani,
	A. Kaplan,
	{\em Sur la cohomologie $L^2$ et la cohomologie d'intersection
	\'{a} coefficients dans une variation de structure de Hodge.}
	C. R. Acad. Sci. Paris S\'{e}r. I Math. {\bf 300} (1985),
	351--353.
	
\bibitem{cks1}
E. Cattani, A. Kaplan and W. Schmid,
{\it Degeneration of Hodge structures},
Ann. of Math. {\bf 123} (1986), 457--535.

\bibitem{cks2}
E. Cattani, A. Kaplan and W. Schmid,
{\it $L^2$ and intersection cohomologies for a polarized variation
of Hodge structure,}
Invent. Math. {\bf 87} (1987), 217--252.

\bibitem{cg}
M. Cornalba and P. Griffiths,
{\it Analytic cycles and vector bundles 
on noncompact algebraic
varieties},
Invent. Math. {\bf 28} (1975), 1--106.

 \bibitem{Toric-Varieties-Book}
	 D. A. Cox, J. B. Little, H. K. Schenck,
	 {\em Toric Varieties},
	Graduate Studies in Mathematics, {\bf 124}.
	American Mathematical Society, Providence, RI, 2011. xxiv+841 pp.

\bibitem{d5}
P. Deligne,
{\em Th\'eor\`eme de Lefschetz et crit\`eres
de d\'eg\'en\'erescence de suites spectrales},
Inst. Hautes \'Etudes 
Sci. Publ. Math. {\bf 35} (1968) 259--278. 

\bibitem{deligne-theorie-HodgeII}
P. Deligne,
{\em Th\'{e}orie de Hodge. II.}
Inst. Hautes \'{E}tudes Sci. Publ. Math. {\bf 40} (1971), 5--57. 
	 
 \bibitem{Fulton-toric-varieties}
W. Fulton,
{\em Introduction to toric varieties},
	 Annals of Mathematics Studies, {\bf 131.}
	 Princeton University Press, Princeton, NJ, 1993.

\bibitem{hotta-tanisaki}
R. Hotta, K. Takeuchi and T. Tanisaki,
{\em $D$-modules, perverse sheaves, 
 and representation theory},
Progress in Mathematics, {\bf 236}. 
Birkh\"auser Boston, Inc., Boston, MA, 2008.

\bibitem{i-s}
J. Iyer and  C. Simpson,
{\em A relation between the parabolic Chern characters
of the de Rham bundles},
 Math. Ann.  {\bf 338},  (2007), 
347--383.
	
 \bibitem{Jost-Yang-Zuo}
	 J. Jost,
	 Y.-H. Yang,
	 K. Zuo,
	 {\em The cohomology of a variation of
	 polarized Hodge structures over a quasi-compact K\"{a}hler manifold},
	 J. Algebraic Geom. {\bf 16} (2007), no. 3, 401--434.
\bibitem{Jost-Yang-Zuo2}
	J. Jost,
	Y.-H. Yang,
	K. Zuo,
	{\em Cohomologies of unipotent harmonic bundles
	over noncompact curves},
	J. Reine Angew. Math. {\bf 609} (2007), 137--159.
\bibitem{k}
M. Kashiwara,
{\it The asymptotic behavior of a variation 
of polarized Hodge str.}
Publ. Res. Inst. Math. Sci {\bf 21} (1985), 853--875.

\bibitem{k2}
M. Kashiwara,
{\em A study of variation of mixed Hodge structure}.
Publ. Res. Inst. Math. Sci. {\bf 22}
	(1986), 991--1024.

\bibitem{kashiwara_text}
M. Kashiwara,
{\em $D$-modules and microlocal calculus},
Translations of Mathematical Monographs, 217. 
Iwanami Series in Modern Mathematics,
American Mathematical Society, 
2003

\bibitem{Kashiwara-Kawai-Hodge-holonomic}
	M. Kashiwara,
	T. Kawai, 
	{\em Hodge structure and holonomic systems},
	Proc. Japan Acad. Ser. A Math. Sci. {\bf 62} (1986), 1--4.
	
\bibitem{k3}
M. Kashiwara and T. Kawai,
{\it The Poincar\'{e} lemma for variations of polarized Hodge
    structure},
Publ. Res. Inst. Math. Sci. {\bf 23} (1987), 345--407. 

\bibitem{Kashiwara-Kawai-partition}
	M. Kashiwara,
	T. Kawai,
	{\em A particular partition of unity:
	an auxiliary tool in Hodge theory},
	IN
	{\em Theta functions--Bowdoin 1987,
	Part 1 (Brunswick, ME, 1987)},
	19--26, Proc. Sympos. Pure Math., {\bf 49}, Part 1,
	Amer. Math. Soc., Providence, RI, 1989.

\bibitem{Kashiwara-Schapira}
M. Kashiwara and P. Schapira,
{\em Sheaves on Manifolds},
Springer-Verlag, (1990).


	
\bibitem{mochi2}
T. Mochizuki,
{\em Asymptotic behaviour of tame harmonic bundles
and an application to pure twistor $D$-modules I, II},
Mem. AMS. {\bf 185}, (2007).

\bibitem{mochi5}
T. Mochizuki,
{\em Kobayashi-Hitchin correspondence
for tame harmonic bundles II},
Geometry\&Topology {\bf 13}, (2009),
359--455.

	
\bibitem{mochi8}
T. Mochizuki,
{\em Asymptotic behaviour of variation of
pure polarized TERP structures},
Publ. Res. Inst. Math. Sci. {\bf 47} (2011), 419--534.
	
\bibitem{Mochizuki-wild}
T. Mochizuki,
{\em Wild harmonic bundles and 
 wild pure twistor $D$-modules},
Ast\'{e}risque {\bf 340}, (2011)

\bibitem{Mochizuki-holonomic-Betti}
	T. Mochizuki,
	{\em Holonomic D-modules with Betti structure},
	M\'{e}m. Soc. Math. Fr. (N.S.) {\bf 138-139} (2014).
	
\bibitem{Mochizuki-MTM}
T. Mochizuki,
{\em Mixed twistor $D$-modules},
Springer, 2015.

\bibitem{Mochizuki-KH-Higgs}
	T. Mochizuki,
	{\em Good wild harmonic bundles and good filtered Higgs bundles},
	SIGMA Symmetry Integrability Geom. Methods Appl. {\bf 17} (2021),
	Paper No. 068, 66 pp.
\bibitem{sabbah2}
C. Sabbah,
{\em Polarizable twistor $D$-modules}
Ast\'{e}risque, {\bf 300},
Soci\'{e}t\'{e} Math\'{e}matique
de France, Paris, (2005).

\bibitem{Sabbah-wild}
C. Sabbah,
{\em Wild twistor $D$-modules},
in
{\em Algebraic analysis and around},
Adv. Stud. Pure Math., {\bf 54},
Math. Soc. Japan, Tokyo, (2009),
293--353.

\bibitem{saito1}
M. Saito,
{\em Modules de Hodge polarisables},
Publ. RIMS., {\bf 24}
(1988), 849--995.

\bibitem{saito2}
M. Saito,
{\em Mixed Hodge modules},
Publ. RIMS., {\bf 26}, (1990),
221-333.

\bibitem{MSaito-proper}
	M. Saito,
	{\em Decomposition theorem for proper K\"{a}hler morphisms},
	Tohoku Math. J. (2) {\bf 42} (1990), 127--147.
	
\bibitem{MSaito-young}
M. Saito,
	{\em A young person's guide to mixed Hodge modules},
	in {\em Hodge theory and $L^2$-analysis}, 517--553,
	Adv. Lect. Math. (ALM), {\bf 39},
	Int. Press, Somerville, MA, 2017. 
\bibitem{MSaito-2022}
	M. Saito,
	{\em Some remarks on decomposition theorem
	for proper K\"ahler morphisms},
	preprint (several versions).
\bibitem{sch}
W. Schmid,
{\it Variation of Hodge structure: 
the singularities of the period mapping},
Invent. Math. {\bf 22} (1973), 211--319.

 \bibitem{Shentu-Zhao}
	 J. Shentu, C. Zhao,
	 {\em $L^2$-representation of Hodge modules},
	 arXiv:2103.04030

 \bibitem{s1}
C. Simpson,
{\it Constructing variations of Hodge structure
using Yang-Mills theory
and application to uniformization},
J. Amer. Math. Soc. {\bf 1} (1988), 867--918.

\bibitem{Simpson-non-compact}
C. Simpson,
{\em Harmonic bundles on noncompact curves},
J. Amer. Math. Soc. {\bf 3} (1990), no. 3, 713--770.
	
		 
\bibitem{s5}
C. Simpson,
{\em
 Higgs bundles and local systems},
Publ. I.H.E.S.,
{\bf 75} (1992),
 5--95.

\bibitem{Simpson-family}
C. Simpson,
{\em Some families of local systems over smooth projective varieties},
Ann. of Math. (2) {\bf 138} (1993), 337--425

	
\bibitem{s3}
C. Simpson,
{\it Mixed twistor structures},
math.AG/9705006.

\bibitem{Wei-Yang}
	C. Wei, R. Yang,
	{\em Cohomology of semisimple local systems
	and the Decomposition theorem},
	arXiv:2109.11578

\bibitem{Wlodarczyk}
J. W{\l}odarczyk, 
{\em Resolution of singularities of analytic spaces},
Proceedings of G\"{o}kova Geometry-Topology Conference 2008, 
31–63, G\"{o}kova Geometry/Topology Conference (GGT), G\"{o}kova, 2009. 

	
\bibitem{z} S. Zucker,
 {\em Hodge theory with degenerating coefficients:
 $L^2$ cohomology in the Poincar\'{e} metric},
 Ann of Math. (2) {\bf 109} (1979), 415--476.	

 \bibitem{Williamson-Bourbaki}
	G. Williamson,
	{\em The Hodge theory of the decomposition theorem},
	S\'{e}minaire Bourbaki. Vol. 2015/2016.
	 Expos\'{e}s 1104--1119.
	 Ast\'{e}risque {\bf 390} (2017),
	Exp. No. 1115, 335--367. 
	 
	
\end{thebibliography}
\end{document}